\documentclass[letterpaper,twoside]{article}
\usepackage[OT1,T1]{fontenc}
\usepackage[dvips]{graphicx,color}
\usepackage{citesort}
\usepackage{pst-3d}
\usepackage{pst-char}
\usepackage{pst-coil}
\usepackage{pst-eps}
\usepackage{pst-fill}
\usepackage{pst-grad}
\usepackage{pst-node}
\usepackage{pst-plot}
\usepackage{pst-text}
\usepackage{pst-tree}
\usepackage{paralist}

\usepackage{ccfonts}
\usepackage{beton}
\usepackage{fancybox}
\usepackage{amssymb}
\usepackage{amsmath,bm}
\usepackage[mathscr]{eucal}
\usepackage{amsfonts}
\usepackage{latexsym}
\usepackage{amsxtra}
\usepackage{amsbsy}
\usepackage{amsthm}
\usepackage{amscd}
\usepackage{amsopn}
\usepackage{amstext}
\usepackage{oldgerm}
\usepackage[varumlaut]{yfonts}
\newcounter{z0}
\newtheorem{ay}{Lemma}[section]

\newtheorem{bbbb}{Proposition}[subsection]
\newtheorem{cccc}{Lemma}[subsection]
\newtheorem{dddd}{Theorem}[subsection]
\newtheorem{ffff}{Corollary}[subsection]
\newtheorem{aaaaa}{Definition}[section]
\newtheorem{bbbbb}{Proposition}[section]
\newtheorem{ccccc}{Lemma}[section]

\newtheorem{fffff}{Corollary}[section]
\newtheorem*{rhp1}{RHP1}
\newtheorem*{rhp2}{RHP2}

\theoremstyle{definition}

\theoremstyle{definition}
\newtheorem{eeee}{Remark}[subsection]

\theoremstyle{definition}
\newtheorem{eeeee}{Remark}[section]

\theoremstyle{definition}
\newtheorem{notrem}{Notational Remark}[section]
\newcommand{\me}{\mathrm{e}}
\newcommand{\mi}{\mathrm{i}}
\newcommand{\md}{\mathrm{d}}
\renewcommand{\Im}{\mathrm{Im}}
\renewcommand{\Re}{\mathrm{Re}}
\newcommand{\id}{\pmb{\mi \md}}
\newcommand{\pvi}{\ensuremath{\int \hspace{-2.75mm} \rule[2.5pt]{2mm}{0.25mm}}}
\newcommand{\vip}{\ensuremath{\int \hspace{-3.35mm} \rule[2.5pt]{2mm}{0.25mm}}}
\newcommand{\norm}[1]{\lVert#1\rVert}
\numberwithin{equation}{section}
\allowdisplaybreaks[4]

\usepackage{relsize,exscale}
\usepackage{stmaryrd}
\usepackage{pxfonts}
\begin{document}
\fontsize{10}{12}\selectfont
\fontencoding{T1}\selectfont
\baselineskip=12pt
\frenchspacing
\title{Asymptotics of Laurent Polynomials of Odd Degree Orthogonal with
Respect to Varying Exponential Weights}
\author{K.~T.-R.~McLaughlin\thanks{\texttt{E-mail: mcl@math.arizona.edu}} \\
Department of Mathematics \\
The University of Arizona \\
617 N.~Santa Rita Ave. \\
P.~O.~Box~ 210089 \\
Tucson, Arizona 85721--0089 \\
U.~S.~A. \and
A.~H.~Vartanian\thanks{\texttt{E-mail: arthurv@math.ucf.edu}} \\
Department of Mathematics \\
University of Central Florida \\
P.~O.~Box~161364 \\
Orlando, Florida 32816--1364 \\
U.~S.~A. \and
X.~Zhou\thanks{\texttt{E-mail: zhou@math.duke.edu}} \\
Department of Mathematics \\
Duke University \\
Box 90320 \\
Durham, North Carolina 27708--0320 \\
U.~S.~A.}
\date{24 January 2006}
\maketitle
\begin{abstract}
\noindent
Let $\Lambda^{\mathbb{R}}$ denote the linear space over $\mathbb{R}$ spanned 
by $z^{k}$, $k \! \in \! \mathbb{Z}$. Define the real inner product (with 
varying exponential weights) $\langle \boldsymbol{\cdot},\boldsymbol{\cdot} 
\rangle_{\mathscr{L}} \colon \Lambda^{\mathbb{R}} \times \Lambda^{\mathbb{R}} 
\! \to \! \mathbb{R}$, $(f,g) \! \mapsto \! \int_{\mathbb{R}}f(s)g(s) \exp 
(-\mathscr{N} \, V(s)) \, \md s$, $\mathscr{N} \! \in \! \mathbb{N}$, where 
the external field $V$ satisfies: (i) $V$ is real analytic on $\mathbb{R} 
\setminus \{0\}$; (ii) $\lim_{\vert x \vert \to \infty}(V(x)/\ln (x^{2} \! + 
\! 1)) \! = \! +\infty$; and (iii) $\lim_{\vert x \vert \to 0}(V(x)/\ln (x^{
-2} \! + \! 1)) \! = \! +\infty$. Orthogonalisation of the (ordered) base 
$\lbrace 1,z^{-1},z,z^{-2},z^{2},\dotsc,z^{-k},z^{k},\dotsc \rbrace$ with 
respect to $\langle \boldsymbol{\cdot},\boldsymbol{\cdot} \rangle_{\mathscr{
L}}$ yields the even degree and odd degree orthonormal Laurent polynomials 
$\lbrace \phi_{m}(z) \rbrace_{m=0}^{\infty}$: $\phi_{2n}(z) \! = \! \xi^{
(2n)}_{-n}z^{-n} \! + \! \dotsb \! + \! \xi^{(2n)}_{n}z^{n}$, $\xi^{(2n)}_{
n} \! > \! 0$, and $\phi_{2n+1}(z) \! = \! \xi^{(2n+1)}_{-n-1}z^{-n-1} \! 
+ \! \dotsb \! + \! \xi^{(2n+1)}_{n}z^{n}$, $\xi^{(2n+1)}_{-n-1} \! > \! 0$.
Define the even degree and odd degree monic orthogonal Laurent polynomials: 
$\boldsymbol{\pi}_{2n}(z) \! := \! (\xi^{(2n)}_{n})^{-1} \phi_{2n}(z)$ and 
$\boldsymbol{\pi}_{2n+1}(z) \! := \! (\xi^{(2n+1)}_{-n-1})^{-1} \phi_{2n+1}
(z)$. Asymptotics in the double-scaling limit as $\mathscr{N},n \! \to \! 
\infty$ such that $\mathscr{N}/n \! = \! 1 \! + \! o(1)$ of $\boldsymbol{
\pi}_{2n+1}(z)$ (in the entire complex plane), $\xi^{(2n+1)}_{-n-1}$, 
$\phi_{2n+1}(z)$ (in the entire complex plane), and Hankel determinant 
ratios associated with the real-valued, bi-infinite, strong moment sequence 
$\left\lbrace c_{k} \! = \! \int_{\mathbb{R}}s^{k} \exp (-\mathscr{N} \, 
V(s)) \, \md s \right\rbrace_{k \in \mathbb{Z}}$ are obtained by formulating 
the odd degree monic orthogonal Laurent polynomial problem as a matrix 
Riemann-Hilbert problem on $\mathbb{R}$, and then extracting the large-$n$ 
behaviour by applying the non-linear steepest-descent method introduced 
in \cite{a1} and further developed in \cite{a2,a3}.

\vspace{0.65cm}

47B36: (Secondary) 30C15, 30C70, 30E05, 30E10, 31A99, 41A20, 41A21, 41A60

\vspace{0.50cm}
{\bf Abbreviated Title.} Asymptotics of Odd Degree Orthogonal Laurent
Polynomials

\vspace{0.50cm}

Laurent-Jacobi matrices, parametrices, Riemann-Hilbert problems, singular
integral equations,

strong moment problems, variational problems
\end{abstract}
\clearpage
\section{Introduction}
Consider the \emph{strong Stieltjes} (resp., \emph{strong Hamburger}) 
\emph{moment problem} (SSMP) (resp., SHMP): given a doubly- or bi-infinite 
(moment) sequence $\{c_{n}\}_{n \in \mathbb{Z}}$ of real numbers:
\begin{enumerate}
\item[(i)] find necessary and sufficient conditions for the existence of a 
non-negative Borel measure $\mu^{\text{SS}}_{\text{MP}}$ (resp., $\mu^{\text{
SH}}_{\text{MP}})$ on $[0,+\infty)$ (resp., $(-\infty,+\infty))$, and with 
infinite support, such that $c_{n} \! = \! \int_{0}^{+\infty}t^{n} \, \md 
\mu^{\text{SS}}_{\text{MP}}(t)$, $n \! \in \! \mathbb{Z}$ (resp., $c_{n} \! 
= \! \int_{-\infty}^{+\infty}t^{n} \, \md \mu^{\text{SH}}_{\text{MP}}(t)$, 
$n \! \in \! \mathbb{Z})$, where the---improper---integral is to be understood 
in the Riemann-Stieltjes sense;
\item[(ii)] when there is a solution, in which case the SSMP (resp., SHMP) is
\emph{determinate}, find conditions for the uniqueness of the solution; and
\item[(iii)] when there is more than one solution, in which case the SSMP
(resp., SHMP) is \emph{indeterminate}, describe the family of all solutions.
\end{enumerate}
The SSMP (resp., SHMP) was introduced in $1980$ (resp., $1981)$ by Jones 
\emph{et al.} \cite{a4} (resp., Jones \emph{et al.} \cite{a5}), and studied 
further in \cite{a6,a7,a8,a9,a10} (see, also, the review article \cite{a11}). 
Unlike the moment theory for the \emph{classical Stieltjes} (resp., 
\emph{classical Hamburger}) \emph{moment problem} (SMP) \cite{a12} (resp., 
(HMP) \cite{a13}), wherein the theory of orthogonal polynomials \cite{a14} 
(and the analytic theory of continued fractions; in particular, $S$- and real 
$J$-fractions) enjoyed a prominent r\^{o}le (see, for example, \cite{a15}), 
the extension of the moment theory to the SSMP and the SHMP introduced a 
`rational generalisation' of the orthogonal polynomials, namely, the 
\emph{orthogonal Laurent} (or $L$-) \emph{polynomials} (as well as the 
introduction of special kinds of continued fractions commonly referred to 
as positive-$T$ fractions), which are now introduced \cite{a6,a7,a8,a9,a10,%
a11,a16,a17,a18,a19,a20,a21}. The SHMP can also be solved using the spectral 
theory of unbounded self-adjoint operators in Hilbert space \cite{a22} (see, 
also, \cite{a23}).

For any pair $(p,q) \! \in \! \mathbb{Z} \times \mathbb{Z}$, with $p \! 
\leqslant \! q$, let $\Lambda^{\mathbb{R}}_{p,q} \! := \! \left\lbrace 
\mathstrut f \colon \mathbb{C}^{\ast} \! \to \! \mathbb{C}; \, f(z) \! = \! 
\sum_{k=p}^{q} \widetilde{\lambda}_{k}z^{k}, \, \widetilde{\lambda}_{k} \! 
\in \! \mathbb{R}, \, k \! = \! p,\dotsc,q \right\rbrace$, where $\mathbb{
C}^{\ast} \! := \! \mathbb{C} \setminus \{0\}$. For any $m \! \in \! \mathbb{
Z}_{0}^{+} \! := \! \{0\} \cup \mathbb{N}$, set $\Lambda^{\mathbb{R}}_{2m} 
\! := \! \Lambda^{\mathbb{R}}_{-m,m}$, $\Lambda^{\mathbb{R}}_{2m+1} \! := \! 
\Lambda^{\mathbb{R}}_{-m-1,m}$, and $\Lambda^{\mathbb{R}} \! := \! \cup_{m 
\in \mathbb{Z}_{0}^{+}}(\Lambda^{\mathbb{R}}_{2m} \cup \Lambda^{\mathbb{R}
}_{2m+1})$. (Note: the sets $\Lambda^{\mathbb{R}}_{p,q}$ and $\Lambda^{
\mathbb{R}}$ form linear spaces over the field $\mathbb{R}$ with respect to 
the operations of addition and multiplication by a scalar.) The ordered base 
for $\Lambda^{\mathbb{R}}$ is $\lbrace 1,z^{-1},z,z^{-2},z^{2},\dotsc,z^{-k},
z^{k},\dotsc \rbrace$, corresponding to the \emph{cyclically-repeated pole 
sequence} $\lbrace \text{no pole},0,\infty,0,\infty,\dotsc,0,\infty,\dotsc 
\rbrace$. For each $l \! \in \! \mathbb{Z}_{0}^{+}$ and $0 \! \not\equiv 
\! f \! \in \! \Lambda^{\mathbb{R}}_{l}$, the $L$-\emph{degree} of $f$, 
symbolically $LD(f)$, is defined as
\begin{equation*}
LD(f) \! := \! l,
\end{equation*}
and the \emph{leading coefficient} of $f$, symbolically $LC(f)$, and the
\emph{trailing coefficient} of $f$, symbolically $TC(f)$, are defined as
follows:
\begin{equation*}
LC(f) \! := \!
\begin{cases}
\widetilde{\lambda}_{m}, &\text{$l \! = \! 2m$,} \\
\widetilde{\lambda}_{-m-1}, &\text{$l \! = \! 2m \! + \! 1$,}
\end{cases}
\end{equation*}
and
\begin{equation*}
TC(f) \! := \!
\begin{cases}
\widetilde{\lambda}_{-m}, &\text{$l \! = \! 2m$,} \\
\widetilde{\lambda}_{m}, &\text{$l \! = \! 2m \! + \! 1$.}
\end{cases}
\end{equation*}

Consider the positive measure on $\mathbb{R}$ (oriented throughout this
work, unless stated otherwise, {}from $-\infty$ to $+\infty)$ given by
\begin{equation*}
\md \widetilde{\mu}(z) \! = \! \widetilde{w}(z) \md z,
\end{equation*}
with varying exponential weight function of the form
\begin{equation*}
\widetilde{w}(z) \! := \! \exp (-\mathscr{N} \, V(z)), \quad \mathscr{N} \!
\in \! \mathbb{N},
\end{equation*}
where the \emph{external field} $V \colon \mathbb{R} \setminus \{0\} \! \to
\! \mathbb{R}$ satisfies the following conditions:
\begin{gather}
V \, \, \text{is real analytic on} \, \, \mathbb{R} \setminus \{0\}; \tag{V1}
\\
\lim_{\vert x \vert \to \infty} \! \left(V(x)/\ln (x^{2} \! + \! 1) \right)
\! = \! +\infty; \tag{V2} \\
\lim_{\vert x \vert \to 0} \! \left(V(x)/\ln (x^{-2} \! + \! 1) \right) \! =
\! +\infty. \tag{V3}
\end{gather}
(For example, a rational function of the form $V(z) \! = \! \sum_{k=-2m_{1}}^{
2m_{2}} \! \varrho_{k}z^{k}$, with $\varrho_{k} \! \in \! \mathbb{R}$, $k \! 
= \! -2m_{1},\dotsc,2m_{2}$, $m_{1,2} \! \in \! \mathbb{N}$, and $\varrho_{-
2m_{1}},\varrho_{2m_{2}} \! > \! 0$ would suffice.) Define (uniquely) the 
\emph{strong moment linear functional} $\mathscr{L}$ by its action on the 
basis elements of $\Lambda^{\mathbb{R}}$: $\mathscr{L} \colon \Lambda^{
\mathbb{R}} \! \to \! \Lambda^{\mathbb{R}}$, $f \! = \! \sum_{k \in \mathbb{
Z}} \widetilde{\lambda}_{k}z^{k} \! \mapsto \! \mathscr{L}(f) \! := \! \sum_{
k \in \mathbb{Z}} \widetilde{\lambda}_{k}c_{k}$, where $c_{k} \! = \! 
\mathscr{L}(z^{k}) \! = \! \int_{\mathbb{R}}s^{k} \exp (-\mathscr{N} \, V(s)) 
\, \md s$, $(k,\mathscr{N}) \! \in \! \mathbb{Z} \times \mathbb{N}$. (Note 
that $\left\lbrace c_{k} \! = \! \int_{\mathbb{R}}s^{k} \exp (-\mathscr{N} \, 
V(s)) \, \md s, \, \mathscr{N} \! \in \! \mathbb{N} \right\rbrace_{k \in 
\mathbb{Z}}$ is a bi-infinite, real-valued, \emph{strong moment sequence}: 
$c_{k}$ is called the \emph{$k$th strong moment of $\mathscr{L}$}.) Associated 
with the above-defined bi-infinite, real-valued, strong moment sequence 
$\lbrace c_{k} \rbrace_{k \in \mathbb{Z}}$ are the \emph{Hankel determinants} 
$H^{(m)}_{k}$, $(m,k) \! \in \! \mathbb{Z} \times \mathbb{N}$ 
\cite{a6,a7,a11,a17}:
\begin{equation}
H^{(m)}_{0} \! := \! 1 \qquad \quad \text{and} \qquad \quad H^{(m)}_{k} \! 
:= \! 
\begin{vmatrix}
c_{m} & c_{m+1} & \cdots & c_{m+k-2} & c_{m+k-1} \\
c_{m+1} & c_{m+2} & \cdots & c_{m+k-1} & c_{m+k} \\
c_{m+2} & c_{m+3} & \cdots & c_{m+k} & c_{m+k+1} \\
\vdots & \vdots & \ddots & \vdots & \vdots \\
c_{m+k-1} & c_{m+k} & \cdots & c_{m+2k-3} & c_{m+2k-2}
\end{vmatrix}.
\end{equation}

Define the real bilinear form $\langle \boldsymbol{\cdot},\boldsymbol{\cdot} 
\rangle_{\mathscr{L}}$ as follows: $\langle \boldsymbol{\cdot},\boldsymbol{
\cdot} \rangle_{\mathscr{L}} \colon \Lambda^{\mathbb{R}} \times \Lambda^{
\mathbb{R}} \! \to \! \mathbb{R}$, $(f,g) \! \mapsto \! \langle f,g \rangle_{
\mathscr{L}} \! := \! \mathscr{L}(f(z)g(z)) \! = \! \int_{\mathbb{R}}f(s)g(s) 
\exp (-\mathscr{N} \, V(s)) \, \md s$, $\mathscr{N} \! \in \! \mathbb{N}$. It 
is a fact \cite{a6,a7,a11,a17} that the bilinear form $\langle \boldsymbol{
\cdot},\boldsymbol{\cdot} \rangle_{\mathscr{L}}$ thus defined is an inner 
product if and only if $H^{(-2m)}_{2m} \! > \! 0$ and $H^{(-2m)}_{2m+1} \! > 
\! 0$ for each $m \! \in \! \mathbb{Z}_{0}^{+}$ (see Equations~(1.8) below, 
and Subsection~2.2, the proof of Lemma~2.2.2); and this fact is used, with 
little or no further reference, throughout this work (see, also, \cite{a24}).
\begin{eeeee}
These latter two (Hankel determinant) inequalities also appear when the
question of the solvability of the SHMP is posed (in this case, the $c_{k}$,
$k \! \in \! \mathbb{Z}$, which appear in Equations~(1.1) should be replaced
by $c_{k}^{\text{\tiny SHMP}}$, $k \! \in \! \mathbb{Z})$: indeed, if these
two inequalities are true $\forall \, \, m \! \in \! \mathbb{Z}_{0}^{+}$, then
there is a non-negative measure $\mu^{\text{\tiny SH}}_{\text{\tiny MP}}$ (on
$\mathbb{R})$ with the given (real) moments. For the case of the SSMP, there
are four (Hankel determinant) inequalities (in this latter case, the $c_{k}$,
$k \! \in \! \mathbb{Z}$, which appear in Equations~(1.1) should be replaced
by $c_{k}^{\text{\tiny SSMP}}$, $k \! \in \! \mathbb{Z})$ which guarantee the
existence of a non-negative measure $\mu^{\text{\tiny SS}}_{\text{\tiny MP}}$
(on $[0,+\infty))$ with the given moments, namely \cite{a4} (see, also,
\cite{a6,a7}): for each $m \! \in \! \mathbb{Z}_{0}^{+}$, $H^{(-2m)}_{2m}
\! > \! 0$, $H^{(-2m)}_{2m+1} \! > \! 0$, $H^{(-2m+1)}_{2m} \! > \! 0$, and
$H^{(-2m-1)}_{2m+1} \! < \! 0$. It is interesting to note that the former
solvability conditions do not automatically imply that the positive (real)
moments $\{c_{k}^{\text{\tiny SHMP}}\}_{k \in \mathbb{Z}_{0}^{+}}$ determine
a measure via the HMP: a similar statement holds true for the SMP (see the
latter four solvability conditions). \hfill $\blacksquare$
\end{eeeee}

If $f \! \in \! \Lambda^{\mathbb{R}}$, then
\begin{equation*}
\norm{f(\cdot)}_{\mathscr{L}} \! := \! (\langle f,f \rangle_{\mathscr{L}}
)^{1/2}
\end{equation*}
is called the \emph{norm of $f$ with respect to $\mathscr{L}$}: note that
$\norm{f(\cdot)}_{\mathscr{L}} \! \geqslant \! 0 \, \, \forall \, \, f \! \in
\! \Lambda^{\mathbb{R}}$, and $\norm{f(\cdot)}_{\mathscr{L}} \! > \! 0$ if $0
\! \not\equiv \! f \! \in \! \Lambda^{\mathbb{R}}$. $\{\phi_{n}^{\flat}(z)
\}_{n \in \mathbb{Z}_{0}^{+}}$ is called a (real) orthonormal Laurent (or
$L$-) polynomial sequence (ONLPS) with respect to $\mathscr{L}$ if, $\forall
\, \, m,n \! \in \! \mathbb{Z}_{0}^{+}$:
\begin{enumerate}
\item[(i)] $\phi_{n}^{\flat} \! \in \! \Lambda^{\mathbb{R}}_{n}$, that is,
$LD(\phi_{n}^{\flat}) \! := \! n$;
\item[(ii)] $\langle \phi_{m}^{\flat},\phi_{n^{\prime}}^{\flat} \rangle_{
\mathscr{L}} \! = \! 0$, $m \! \not= \! n^{\prime}$, or, alternatively,
$\langle f,\phi_{n}^{\flat} \rangle_{\mathscr{L}} \! = \! 0 \, \, \forall \,
\, f \! \in \! \Lambda^{\mathbb{R}}_{n-1}$;
\item[(iii)] $\langle \phi_{m}^{\flat},\phi_{m}^{\flat} \rangle_{\mathscr{L}}
\! =: \! \norm{\phi_{m}^{\flat}(\cdot)}^{2}_{\mathscr{L}} \! = \! 1$.
\end{enumerate}
Orthonormalisation of $\lbrace 1,z^{-1},z,z^{-2},z^{2},\dotsc,z^{-n},z^{n},
\dotsc \rbrace$, corresponding to the cyclically-repeated pole sequence
$\lbrace \text{no pole},0,\infty,0,\infty,\dotsc,0,\infty,\dotsc \rbrace$,
with respect to $\langle \boldsymbol{\cdot},\boldsymbol{\cdot} \rangle_{
\mathscr{L}}$ via the Gram-Schmidt orthogonalisation method, leads to the
ONLPS, or, simply, orthonormal Laurent (or $L$-) polynomials (OLPs), $\{
\phi_{m}(z)\}_{m \in \mathbb{Z}_{0}^{+}}$, which, by suitable normalisation,
may be written as, for $m \! = \! 2n$,
\begin{equation}
\phi_{2n}(z) \! = \! \xi^{(2n)}_{-n}z^{-n} \! + \! \dotsb \! + \!
\xi^{(2n)}_{n}z^{n}, \qquad \xi^{(2n)}_{n} \! > \! 0,
\end{equation}
and, for $m \! = \! 2n \! + \! 1$,
\begin{equation}
\phi_{2n+1}(z) \! = \! \xi^{(2n+1)}_{-n-1}z^{-n-1} \! + \! \dotsb \! + \!
\xi^{(2n+1)}_{n}z^{n}, \qquad \xi^{(2n+1)}_{-n-1} \! > \! 0.
\end{equation}
The $\phi_{n}$'s are normalised so that they all have real coefficients; in
particular, the leading coefficients, $LC(\phi_{2n}) \! := \! \xi^{(2n)}_{n}$
and $LC(\phi_{2n+1}) \! := \! \xi^{(2n+1)}_{-n-1}$, $n \! \in \! \mathbb{Z}_{
0}^{+}$, are both positive, $\xi^{(0)}_{0} \! = \! 1$, and $\phi_{0}(z) \!
\equiv \! 1$. Even though the leading coefficients $\xi^{(2n)}_{n}$ and
$\xi^{(2n+1)}_{-n-1}$, $n \! \in \! \mathbb{Z}_{0}^{+}$, are non-zero (in
particular, they are positive), no such restriction applies to the trailing
coefficients, $TC(\phi_{2n}) \! := \! \xi^{(2n)}_{-n}$ and $TC(\phi_{2n+1}) \!
:= \! \xi^{(2n+1)}_{n}$, $n \! \in \! \mathbb{Z}_{0}^{+}$. Furthermore, note
that, by construction:
\begin{enumerate}
\item[(1)] $\langle \phi_{2n},z^{j} \rangle_{\mathscr{L}} \! = \! 0$, $j \! =
\! -n,\dotsc,n \! - \! 1$;
\item[(2)] $\langle \phi_{2n+1},z^{j} \rangle_{\mathscr{L}} \! = \! 0$, $j \!
= \! -n,\dotsc,n$;
\item[(3)] $\langle \phi_{j},\phi_{k} \rangle_{\mathscr{L}} \! = \! \delta_{j
k}$, $j,k \! \in \! \mathbb{Z}_{0}^{+}$, where $\delta_{jk}$ is the Kronecker
delta.
\end{enumerate}

It is convenient to introduce the monic orthogonal Laurent (or $L$-)
polynomials, $\boldsymbol{\pi}_{j}(z)$, $j \! \in \! \mathbb{Z}_{0}^{+}$:
(i) for $j \! = \! 2n$, $n \! \in \! \mathbb{Z}_{0}^{+}$, with $\boldsymbol{
\pi}_{0}(z) \! \equiv \! 1$,
\begin{equation}
\boldsymbol{\pi}_{2n}(z) \! := \! \phi_{2n}(z)(\xi^{(2n)}_{n})^{-1} \! = \!
\nu^{(2n)}_{-n}z^{-n} \! + \! \dotsb \! + \! z^{n}, \qquad \nu^{(2n)}_{-n} \!
:= \! \xi^{(2n)}_{-n}/\xi^{(2n)}_{n};
\end{equation}
and (ii) for $j \! = \! 2n \! + \! 1$, $n \! \in \! \mathbb{Z}_{0}^{+}$,
\begin{equation}
\boldsymbol{\pi}_{2n+1}(z) \! := \! \phi_{2n+1}(z)(\xi^{(2n+1)}_{-n-1})^{-1}
\! = \! z^{-n-1} \! + \! \dotsb \! + \! \nu^{(2n+1)}_{n}z^{n}, \qquad \nu^{(2
n+1)}_{n} \! := \! \xi^{(2n+1)}_{n}/\xi^{(2n+1)}_{-n-1}.
\end{equation}
The monic orthogonal $L$-polynomials, $\boldsymbol{\pi}_{j}(z)$, $j \! \in \!
\mathbb{Z}_{0}^{+}$, possess the following properties:
\begin{enumerate}
\item[(1)] $\langle \boldsymbol{\pi}_{2n},z^{j} \rangle_{\mathscr{L}} \! = \!
0$, $j \! = \! -n,\dotsc,n \! - \! 1$;
\item[(2)] $\langle \boldsymbol{\pi}_{2n+1},z^{j} \rangle_{\mathscr{L}} \! =
\! 0$, $j \! = \! -n,\dotsc,n$;
\item[(3)] $\langle \boldsymbol{\pi}_{2n},\boldsymbol{\pi}_{2n} \rangle_{
\mathscr{L}} \! =: \! \norm{\boldsymbol{\pi}_{2n}(\cdot)}^{2}_{\mathscr{L}}
\! = \! (\xi^{(2n)}_{n})^{-2}$, whence $\xi^{(2n)}_{n} \! = \! 1/\norm{
\boldsymbol{\pi}_{2n}(\cdot)}_{\mathscr{L}}$ $(> \! 0)$;
\item[(4)] $\langle \boldsymbol{\pi}_{2n+1},\boldsymbol{\pi}_{2n+1} \rangle_{
\mathscr{L}} \! =: \! \norm{\boldsymbol{\pi}_{2n+1}(\cdot)}^{2}_{\mathscr{L}}
\! = \! (\xi^{(2n+1)}_{-n-1})^{-2}$, whence $\xi^{(2n+1)}_{-n-1} \! = \!
1/\norm{\boldsymbol{\pi}_{2n+1}(\cdot)}_{\mathscr{L}}$ $(> \! 0)$.
\end{enumerate}
Furthermore, in terms of the Hankel determinants, $H^{(m)}_{k}$, $(m,k) \!
\in \! \mathbb{Z} \times \mathbb{N}$, associated with the real-valued,
bi-infinite, strong moment sequence $\left\lbrace c_{k} \! = \! \int_{
\mathbb{R}}s^{k} \exp (-\mathscr{N} \, V(s)) \, \md s, \, \mathscr{N} \!
\in \! \mathbb{N} \right\rbrace_{k \in \mathbb{Z}}$, the monic orthogonal
$L$-polynomials, $\boldsymbol{\pi}_{j}(z)$, $j \! \in \! \mathbb{Z}_{0}^{+}$,
are represented via the following determinantal formulae \cite{a6,a7,a11,a17}:
for $m \! \in \! \mathbb{Z}_{0}^{+}$,
\begin{gather}
\boldsymbol{\pi}_{2m}(z) \! = \! \dfrac{1}{H^{(-2m)}_{2m}}
\begin{vmatrix}
c_{-2m} & c_{-2m+1} & \cdots & c_{-1} & z^{-m} \\
c_{-2m+1} & c_{-2m+2} & \cdots & c_{0} & z^{-m+1} \\
\vdots & \vdots & \ddots & \vdots & \vdots \\
c_{-1} & c_{0} & \cdots & c_{2m-2} & z^{m-1} \\
c_{0} & c_{1} & \cdots & c_{2m-1} & z^{m}
\end{vmatrix}, \\
\intertext{and}
\boldsymbol{\pi}_{2m+1}(z) \! = \! -\dfrac{1}{H^{(-2m)}_{2m+1}}
\begin{vmatrix}
c_{-2m-1} & c_{-2m} & \cdots & c_{-1} & z^{-m-1} \\
c_{-2m} & c_{-2m+1} & \cdots & c_{0} & z^{-m} \\
\vdots & \vdots & \ddots & \vdots & \vdots \\
c_{-1} & c_{0} & \cdots & c_{2m-1} & z^{m-1} \\
c_{0} & c_{1} & \cdots & c_{2m} & z^{m}
\end{vmatrix};
\end{gather}
moreover, it can be shown that (see, for example, \cite{a11,a17}), for $n \!
\in \! \mathbb{Z}_{0}^{+}$,
\begin{gather}
\xi^{(2n)}_{n} \! \left(= \! \dfrac{1}{\norm{\boldsymbol{\pi}_{2n}(\cdot)}_{
\mathscr{L}}} \right) \! =\sqrt{\dfrac{H^{(-2n)}_{2n}}{H^{(-2n)}_{2n+1}}},
\qquad \xi^{(2n+1)}_{-n-1} \! \left(= \! \dfrac{1}{\norm{\boldsymbol{\pi}_{2n
+1}(\cdot)}_{\mathscr{L}}} \right) \! = \! \sqrt{\dfrac{H^{(-2n)}_{2n+1}}{H^{
(-2n-2)}_{2n+2}}}, \\
\nu^{(2n)}_{-n} \! \left(:= \! \dfrac{\xi^{(2n)}_{-n}}{\xi^{(2n)}_{n}} \right)
\! =\dfrac{H^{(-2n+1)}_{2n}}{H^{(-2n)}_{2n}}, \qquad \quad \nu^{(2n+1)}_{n}
\! \left(:= \! \dfrac{\xi^{(2n+1)}_{n}}{\xi^{(2n+1)}_{-n-1}} \right) \! = \!
-\dfrac{H^{(-2n-1)}_{2n+1}}{H^{(-2n)}_{2n+1}}.
\end{gather}

For each $m \! \in \! \mathbb{Z}_{0}^{+}$, the monic orthogonal $L$-polynomial
$\boldsymbol{\pi}_{m}(z)$ and the index $m$ are called \emph{non-singular} if
$0 \! \not= \! TC(\boldsymbol{\pi}_{m}) \! := \!
\begin{cases}
\nu^{(2n)}_{-n}, &\text{$m \! = \! 2n$,} \\
\nu^{(2n+1)}_{n}, &\text{$m \! = \! 2n \! + \! 1$;}
\end{cases}$ otherwise, $\boldsymbol{\pi}_{m}(z)$ and $m$ are 
\emph{singular}. {}From Equations~(1.9), it can be seen that, for each $m \! 
\in \! \mathbb{Z}_{0}^{+}$:
\begin{enumerate}
\item[(i)] $\boldsymbol{\pi}_{2m}(z)$ is non-singular (resp., singular) if 
$H^{(-2m+1)}_{2m} \! \not= \! 0$ (resp., $H^{(-2m+1)}_{2m} \! = \! 0)$;
\item[(ii)] $\boldsymbol{\pi}_{2m+1}(z)$ is non-singular (resp., singular) 
if $H^{(-2m-1)}_{2m+1} \! \not= \! 0$ (resp., $H^{(-2m-1)}_{2m+1} \! = \! 
0)$.
\end{enumerate}
For each $m \! \in \! \mathbb{Z}_{0}^{+}$, let $\mu_{2m} \! := \!
\operatorname{card} \lbrace \mathstrut z; \, \boldsymbol{\pi}_{2m}(z) \! = \!
0 \rbrace$, and $\mu_{2m+1} \! := \! \operatorname{card} \lbrace \mathstrut z;
\boldsymbol{\pi}_{2m+1}(z) \! = \! 0 \rbrace$. It is an established fact
\cite{a6,a7,a17} that, for $m \! \in \! \mathbb{Z}_{0}^{+}$:
\begin{enumerate}
\item[(1)] the zeros of $\boldsymbol{\pi}_{2m}(z)$ are real, simple, and 
non-zero, and $\mu_{2m} \! = \! 2m$ (resp., $2m \! - \! 1)$ if $\boldsymbol{
\pi}_{2m}(z)$ is non-singular (resp., singular);
\item[(2)] the zeros of $\boldsymbol{\pi}_{2m+1}(z)$ are real, simple, and 
non-zero, and $\mu_{2m+1} \! = \! 2m \! + \! 1$ (resp., $2m)$ if $\boldsymbol{
\pi}_{2m+1}(z)$ is non-singular (resp., singular).
\end{enumerate}
For each $m \! \in \! \mathbb{Z}_{0}^{+}$, it can be shown that, via a
straightforward factorisation argument and using Equations~(1.6) and~(1.7):
\begin{enumerate}
\item[(i)] if $\boldsymbol{\pi}_{2m}(z)$ is non-singular, upon setting
$\left\lbrace \mathstrut \alpha^{(2m)}_{k}, \, k \! = \! 1,\dotsc,2m
\right\rbrace \! := \! \left\lbrace \mathstrut z; \, \boldsymbol{\pi}_{2m}(z)
\! = \! 0 \right\rbrace$,
\begin{equation*}
\prod_{k=1}^{2m} \alpha^{(2m)}_{k} \! = \nu^{(2m)}_{-m};
\end{equation*}
\item[(ii)] if $\boldsymbol{\pi}_{2m+1}(z)$ is non-singular, upon setting
$\left\lbrace \mathstrut \alpha^{(2m+1)}_{k}, \, k \! = \! 1,\dotsc,2m \! + \!
1 \right\rbrace \! := \! \left\lbrace \mathstrut z; \, \boldsymbol{\pi}_{2m+1}
(z) \! = \! 0 \right\rbrace$,
\begin{equation*}
\prod_{k=1}^{2m+1} \alpha^{(2m+1)}_{k} \! = \! - \! \left(\nu^{(2m+1)}_{m}
\right)^{-1}.
\end{equation*}
\end{enumerate}
\begin{eeeee}
It is important to note \cite{a10} that the classical and strong moment 
problems (SMP, HMP, SSMP, and SHMP) are special cases of a more general 
theory, where moments corresponding to an arbitrary, countable sequence 
of (fixed) points are involved (in the classical and strong moment cases, 
respectively, the points are $\infty$ repeated and $0,\infty$ cyclically 
repeated), and where \emph{orthogonal rational functions} \cite{a25} play 
the r\^{o}le of orthogonal polynomials and orthogonal Laurent (or $L$-) 
polynomials; furthermore, since $L$-polynomials are rational functions 
with (fixed) poles at the origin and at the point at infinity, the step 
towards a more general theory where poles are at arbitrary, but fixed, 
positions/locations in $\mathbb{C} \cup \{\infty\}$ is natural. \hfill 
$\blacksquare$
\end{eeeee}

Unlike orthogonal polynomials, which satisfy a system of three-term recurrence 
relations, monic orthogonal, and orthonormal, $L$-polynomials may satisfy 
recurrence relations consisting of a pair of four-term recurrence relations 
\cite{a11}, a pair of systems of three- or five-term recurrence relations 
(which is guaranteed in the case when the corresponding monic orthogonal, 
and orthonormal, $L$-polynomials are non-singular) \cite{a11,a16,a17}, or 
a system consisting of four five-term recurrence relations \cite{a23}.
\begin{eeeee}
The non-vanishing of the leading and trailing coefficients of the OLPs
$\lbrace \phi_{m}(z) \rbrace_{m=0}^{\infty}$, that is,
\begin{equation*}
LC(\phi_{m}) \! := \!
\begin{cases}
\xi^{(2n)}_{n}, &\text{$m \! = \! 2n$,} \\
\xi^{(2n+1)}_{-n-1}, &\text{$m \! = \! 2n \! + \! 1$,}
\end{cases}
\end{equation*}
and
\begin{equation*}
TC(\phi_{m}) \! := \!
\begin{cases}
\xi^{(2n)}_{-n}, &\text{$m \! = \! 2n$,} \\
\xi^{(2n+1)}_{n}, &\text{$m \! = \! 2n \! + \! 1$,}
\end{cases}
\end{equation*}
respectively, is of paramount importance: if both these conditions are not
satisfied, then the `length' of the recurrence relations may be greater
than three \cite{a16} (see, also, \cite{a24}). \hfill $\blacksquare$
\end{eeeee}
It can be shown that (see, for example, \cite{a17}; see, also, Chapter~11 of 
\cite{a25}), if $\lbrace \boldsymbol{\pi}_{m}(z) \rbrace_{m \in \mathbb{Z}_{
0}^{+}}$, as defined above, is a non-singular, monic orthogonal $L$-polynomial
sequence, that is, $H^{(-2n+1)}_{2n} \! \not= \! 0$ $(m \! = \! 2n)$ and $H^{
(-2n-1)}_{2n+1} \! \not= \! 0$ $(m \! = \! 2n \! + \! 1)$, then $\lbrace
\boldsymbol{\pi}_{m}(z) \rbrace_{m \in \mathbb{Z}_{0}^{+}}$ satisfy the pair
of three-term recurrence relations
\begin{align*}
\boldsymbol{\pi}_{2m+1}(z) &= \left(\dfrac{z^{-1}}{\beta_{2m}^{\natural}} \!
+ \! \beta_{2m+1}^{\natural} \right) \! \boldsymbol{\pi}_{2m}(z) \! + \!
\lambda_{2m+1}^{\natural} \boldsymbol{\pi}_{2m-1}(z), \\
\boldsymbol{\pi}_{2m+2}(z) &= \left(\dfrac{z}{\beta_{2m+1}^{\natural}} \! + \!
\beta_{2m+2}^{\natural} \right) \! \boldsymbol{\pi}_{2m+1}(z) \! + \! \lambda_{
2m+2}^{\natural} \boldsymbol{\pi}_{2m}(z),
\end{align*}
where $\boldsymbol{\pi}_{-1}(z) \! \equiv \! 0$,
\begin{gather*}
\beta_{2m}^{\natural}= \nu^{(2m)}_{-m}, \qquad \qquad \quad \beta_{2m+1}^{
\natural} = \nu^{(2m+1)}_{m}, \\
\lambda_{2m+1}^{\natural} \! = \! -\dfrac{H^{(-2m-1)}_{2m+1} H^{(-2m+2)}_{2
m-1}}{H^{(-2m)}_{2m} H^{(-2m+1)}_{2m}} \quad (\not= \! 0), \qquad \qquad 
\lambda_{2m+2}^{\natural} \! = \! -\dfrac{H^{(-2m-1)}_{2m+2} H^{(-2m)}_{2m}
}{H^{(-2m)}_{2m+1} H^{(-2m-1)}_{2m+1}} \quad (\not= \! 0),
\end{gather*}
and $\lambda_{j} \beta_{j-1}/\beta_{j} \! > \! 0$, $j \! \in \! \mathbb{N}$,
with $\lambda_{1} \! := \! -c_{-1}$, leading to a \emph{tri-diagonal-type
Laurent-Jacobi matrix} $\mathcal{F}$ for the `mixed' mapping
\begin{equation*}
\mathscr{F} \colon \Lambda^{\mathbb{R}} \! \to \! \Lambda^{\mathbb{R}}, \, \,
f(z) \! \mapsto \! (z^{-1}(\oplus_{n=0}^{\infty} \operatorname{diag}(1,0)) \!
+ \! z(\oplus_{n=0}^{\infty} \operatorname{diag}(0,1)))f(z),
\end{equation*}
where $\oplus_{n=0}^{\infty} \operatorname{diag}(1,0) \! := \! \operatorname{
diag}(1,0,\dotsc,1,0,\dotsc)$ and $\oplus_{n=0}^{\infty} \operatorname{diag}
(0,1) \! := \! \operatorname{diag}(0,1,\dotsc,0,1,\dotsc)$,
\begin{align*}
\setcounter{MaxMatrixCols}{13}
\mathcal{F} =& \, \operatorname{diag} \! \left(\beta_{0}^{\natural},\beta_{
1}^{\natural},\beta_{2}^{\natural},\dotsc \right)
\left(
\begin{smallmatrix}
-\beta_{1}^{\natural} & 1 & & & & & & & & & & & \\
-\lambda_{2}^{\natural} & -\beta_{2}^{\natural} & 1 & & & & & & & & & & \\
 & -\lambda_{3}^{\natural} & -\beta_{3}^{\natural} & 1 & & & & & & & & & \\
 & & -\lambda_{4}^{\natural} & -\beta_{4}^{\natural} & 1 & & & & & & & & \\
 & & & -\lambda_{5}^{\natural} & -\beta_{5}^{\natural} & 1 & & & & & & & \\
 & & & & -\lambda_{6}^{\natural} & -\beta_{6}^{\natural} & 1 & & & & & & \\
 & & & & & \ddots & \ddots & \ddots & & & & & \\
 & & & & & & -\lambda_{2m+1}^{\natural} & -\beta_{2m+1}^{\natural} & 1 & & & &
& \\
 & & & & & & & -\lambda_{2m+2}^{\natural} & -\beta_{2m+2}^{\natural} & 1 & & &
\\
 & & & & & & & & \ddots & \ddots & \ddots & &
\end{smallmatrix}
\right),
\end{align*}
with zeros outside the indicated diagonals (in terms of $\lbrace \phi_{m}(z)
\rbrace_{m \in \mathbb{Z}_{0}^{+}}$, the pair of three-term recurrence
relations reads \cite{a16}:
\begin{gather*}
\phi_{2m+1}(z) \! = \! (z^{-1} \! + \! \mathfrak{g}_{2m+1}) \phi_{2m}(z) \! +
\! \mathfrak{f}_{2m+1} \phi_{2m-1}(z), \\
\phi_{2m+2}(z) \! = \! (1 \! + \! \mathfrak{g}_{2m+2}z) \phi_{2m+1}(z) \! +
\! \mathfrak{f}_{2m+2} \phi_{2m}(z),
\end{gather*}
where $\mathfrak{f}_{2m+1},\mathfrak{f}_{2m+2} \! \not= \! 0$, $m \! \in \!
\mathbb{Z}_{0}^{+}$, $\phi_{-1}(z) \! \equiv \! 0$, and $\phi_{0}(z) \! \equiv
\! 1)$; otherwise, $\lbrace \boldsymbol{\pi}_{m}(z) \rbrace_{m \in \mathbb{
Z}_{0}^{+}}$ satisfy the following pair of five-term recurrence relations
\cite{a17}, with $\boldsymbol{\pi}_{-j}(z) \! \equiv \! 0$, $j \! = \! 1,2$,
\begin{align*}
\boldsymbol{\pi}_{2m+2}(z) =& \, \gamma_{2m+2,2m-2}^{\flat} \boldsymbol{\pi}_{
2m-2}(z) \! + \! \gamma_{2m+2,2m-1}^{\flat} \boldsymbol{\pi}_{2m-1}(z) \! + \!
(z \! + \! \gamma_{2m+2,2m}^{\flat}) \boldsymbol{\pi}_{2m}(z) \\
+& \, \gamma_{2m+2,2m+1}^{\flat} \boldsymbol{\pi}_{2m+1}(z), \\
\boldsymbol{\pi}_{2m+3}(z) =& \, \gamma_{2m+3,2m-1}^{\flat} \boldsymbol{\pi}_{
2m-1}(z) \! + \! \gamma_{2m+3,2m}^{\flat} \boldsymbol{\pi}_{2m}(z) \! + \!
(z^{-1} \! + \! \gamma_{2m+3,2m+1}^{\flat}) \boldsymbol{\pi}_{2m+1}(z) \\
+& \, \gamma_{2m+3,2m+2}^{\flat} \boldsymbol{\pi}_{2m+2}(z),
\end{align*}
where $\gamma_{l,k} \! = \! 0$, $k \! < \! 0$, $l \! \geqslant \! 2$, leading
to a \emph{penta-diagonal-type Laurent-Jacobi matrix} $\mathcal{G}$ for the
`mixed' mapping
\begin{equation*}
\mathscr{G} \colon \Lambda^{\mathbb{R}} \! \to \! \Lambda^{\mathbb{R}}, \,
\, g(z) \! \mapsto \! (z(\oplus_{n=0}^{\infty} \operatorname{diag}(1,0))
\! + \! z^{-1}(\oplus_{n=0}^{\infty} \operatorname{diag}(0,1)))g(z),
\end{equation*}
\begin{align*}
\setcounter{MaxMatrixCols}{14}
&\mathcal{G} =
\left(
\begin{smallmatrix}
-\gamma_{2,0}^{\flat} & -\gamma_{2,1}^{\flat} & 1 & & & & & & & & & & & \\
-\gamma_{3,0}^{\flat} & -\gamma_{3,1}^{\flat} & -\gamma_{3,2}^{\flat} & 1 & &
& & & & & & & & \\
-\gamma_{4,0}^{\flat} & -\gamma_{4,1}^{\flat} & -\gamma_{4,2}^{\flat} &
-\gamma_{4,3}^{\flat} & 1 & & & & & & & & \\
 & -\gamma_{5,1}^{\flat} & -\gamma_{5,2}^{\flat} & -\gamma_{5,3}^{\flat} &
-\gamma_{5,4}^{\flat} & 1 & & & & & & & & \\
 & & -\gamma_{6,2}^{\flat} & -\gamma_{6,3}^{\flat} & -\gamma_{6,4}^{\flat} &
-\gamma_{6,5}^{\flat} & 1 & & & & & & & \\
 & & & \ddots & \ddots & \ddots & \ddots & \ddots & & & & & \\
 & & & & -\gamma_{2m+2,2m-2}^{\flat} & -\gamma_{2m+2,2m-1}^{\flat} & -\gamma_{
2m+2,2m}^{\flat} & -\gamma_{2m+2,2m+1}^{\flat} & 1 & & & & & \\
 & & & & & -\gamma_{2m+3,2m-1}^{\flat} & -\gamma_{2m+3,2m}^{\flat} & -\gamma_{
2m+3,2m+1}^{\flat} & -\gamma_{2m+3,2m+2}^{\flat} & 1 & & & & \\
 & & & & & & \ddots & \ddots & \ddots & \ddots & \ddots & & &
\end{smallmatrix}
\right),
\end{align*}
with zeros outside the indicated diagonals. The general form of these 
(system of) recurrence relations is a pair of three- and five-term recurrence 
relations \cite{a23}: for $n \! \in \! \mathbb{Z}_{0}^{+}$,
\begin{gather*}
z \phi_{2n+1}(z) \! = \! b_{2n+1}^{\sharp} \phi_{2n}(z) \! + \! a_{2n+1}^{
\sharp} \phi_{2n+1}(z) \! + \! b_{2n+2}^{\sharp} \phi_{2n+2}(z), \\
z \phi_{2n}(z) \! = \! c_{2n}^{\sharp} \phi_{2n-2}(z) \! + \! b_{2n}^{\sharp}
\phi_{2n-1}(z) \! + \! a_{2n}^{\sharp} \phi_{2n}(z) \! + \! b_{2n+1}^{\sharp}
\phi_{2n+1}(z) \! + \! c_{2n+2}^{\sharp} \phi_{2n+2}(z),
\end{gather*}
where all the coefficients are real, $c_{0}^{\sharp} \! = \! b_{0}^{\sharp} \!
= \! 0$, and $c_{2k}^{\sharp} \! > \! 0$, $k \! \in \! \mathbb{N}$, and
\begin{gather*}
z^{-1} \phi_{2n}(z) \! = \! \beta_{2n}^{\sharp} \phi_{2n-1}(z) \! + \!
\alpha_{2n}^{\sharp} \phi_{2n}(z) \! + \! \beta_{2n+1}^{\sharp} \phi_{2n+1}
(z), \\
z^{-1} \phi_{2n+1}(z) \! = \! \gamma_{2n+1}^{\sharp} \phi_{2n-1}(z) \! + \!
\beta_{2n+1}^{\sharp} \phi_{2n}(z) \! + \! \alpha_{2n+1}^{\sharp} \phi_{2n+1}
(z) \! + \! \beta_{2n+2}^{\sharp} \phi_{2n+2}(z) \! + \! \gamma_{2n+3}^{
\sharp} \phi_{2n+3}(z),
\end{gather*}
where all the coefficients are real, $\beta_{0}^{\sharp} \! = \! \gamma_{1}^{
\sharp} \! = \! 0$, $\beta_{1}^{\sharp} \! > \! 0$, and $\gamma_{2l+1}^{
\sharp} \! > \! 0$, $l \! \in \! \mathbb{N}$, leading, respectively, to the
real-symmetric, \emph{tri-penta-diagonal-type Laurent-Jacobi matrices},
$\mathcal{J}$ and $\mathcal{K}$, for the mappings
\begin{equation*}
\mathscr{J} \colon \Lambda^{\mathbb{R}} \! \to \! \Lambda^{\mathbb{R}}, \, \,
j(z) \! \mapsto \! zj(z) \qquad \text{and} \qquad \mathscr{K} \colon \Lambda^{
\mathbb{R}} \! \to \! \Lambda^{\mathbb{R}}, \, \, k(z) \! \mapsto \! z^{-1}
k(z),
\end{equation*}
\begin{align*}
\setcounter{MaxMatrixCols}{18}
\mathcal{J} & = \!
\begin{pmatrix}
a_{0}^{\sharp} & b_{1}^{\sharp} & c_{2}^{\sharp} & & & & & & & & & & & & & \\
b_{1}^{\sharp} & a_{1}^{\sharp} & b_{2}^{\sharp} & & & & & & & & & & & & & \\
c_{2}^{\sharp} & b_{2}^{\sharp} & a_{2}^{\sharp} & b_{3}^{\sharp} & c_{4}^{
\sharp} & & & & & & & & & & & \\
 & & b_{3}^{\sharp} & a_{3}^{\sharp} & b_{4}^{\sharp} & & & & & & & & & & & \\
 & & c_{4}^{\sharp} & b_{4}^{\sharp} & a_{4}^{\sharp} & b_{5}^{\sharp} & c_{
6}^{\sharp} & & & & & & & & & \\
 & & & & b_{5}^{\sharp} & a_{5}^{\sharp} & b_{6}^{\sharp} & & & & & & & & & \\
 & & & & c_{6}^{\sharp} & b_{6}^{\sharp} & a_{6}^{\sharp} & b_{7}^{\sharp} &
c_{8}^{\sharp} & & & & & & & \\
 & & & & & & b_{7}^{\sharp} & a_{7}^{\sharp} & b_{8}^{\sharp} & & & & & & & \\
 & & & & & & c_{8}^{\sharp} & b_{8}^{\sharp} & a_{8}^{\sharp} & b_{9}^{\sharp}
& c_{10}^{\sharp} & & & & & \\
 & & & & & & & & \ddots & \ddots & \ddots & & & & & \\
 & & & & & & & & b_{2k+1}^{\sharp} & a_{2k+1}^{\sharp} & b_{2k+2}^{\sharp} & &
& & & \\
 & & & & & & & & c_{2k+2}^{\sharp} & b_{2k+2}^{\sharp} & a_{2k+2}^{\sharp} &
b_{2k+3}^{\sharp} & c_{2k+4}^{\sharp} & & & \\
 & & & & & & & & & & \ddots & \ddots & \ddots & & &
\end{pmatrix},
\end{align*}
and
\begin{align*}
\setcounter{MaxMatrixCols}{18}
\mathcal{K} &= \!
\begin{pmatrix}
\alpha_{0}^{\sharp} & \beta_{1}^{\sharp} & & & & & & & & & & & & & & \\
\beta_{1}^{\sharp} & \alpha_{1}^{\sharp} & \beta_{2}^{\sharp} & \gamma_{3}^{
\sharp}  & & & & & & & & & & & & \\
 & \beta_{2}^{\sharp} & \alpha_{2}^{\sharp} & \beta_{3}^{\sharp} & & & & & & &
& & & & & \\
 & \gamma_{3}^{\sharp} & \beta_{3}^{\sharp} & \alpha_{3}^{\sharp} & \beta_{4}^{
\sharp}  & \gamma_{5}^{\sharp} & & & & & & & & & & \\
 & & & \beta_{4}^{\sharp} & \alpha_{4}^{\sharp} & \beta_{5}^{\sharp} & & & & &
& & & & & \\
 & & & \gamma_{5}^{\sharp} & \beta_{5}^{\sharp} & \alpha_{5}^{\sharp} & \beta_{
6}^{\sharp} & \gamma_{7}^{\sharp} & & & & & & & & \\
 & & & & & \beta_{6}^{\sharp} & \alpha_{6}^{\sharp} & \beta_{7}^{\sharp} & & &
& & & & & \\
 & & & & & \gamma_{7}^{\sharp} & \beta_{7}^{\sharp} & \alpha_{7}^{\sharp} &
\beta_{8}^{\sharp} & \gamma_{9}^{\sharp} & & & & & & \\
 & & & & & & \beta_{8}^{\sharp} & \alpha_{8}^{\sharp} & \beta_{9}^{\sharp} & &
& & & & & \\
 & & & & & & \ddots & \ddots & \ddots & & & & & & & \\
 & & & & & & & & \gamma_{2k+1}^{\sharp} & \beta_{2k+1}^{\sharp} & \alpha_{2k+
1}^{\sharp} & \beta_{2k+2}^{\sharp} & \gamma_{2k+3}^{\sharp} & & & \\
 & & & & & & & & & & \beta_{2k+2}^{\sharp} & \alpha_{2k+2}^{\sharp} & \beta_{2
k+3}^{\sharp} & & & \\
 & & & & & & & & & & \ddots & \ddots & \ddots & & &
\end{pmatrix},
\end{align*}
with zeros outside the indicated diagonals; moreover, as shown in \cite{a23}, 
$\mathcal{J}$ and $\mathcal{K}$ are formal inverses, that is, $\mathcal{J} 
\mathcal{K} \! = \! \mathcal{K} \mathcal{J} \! = \! \operatorname{diag}(1,
\dotsc,1,\dotsc)$ (see, also, \cite{a26,a27,a28,a29,a30}). (Note: $\mathcal{
J}$ (resp., $\mathcal{K})$ is the matrix representation of the multiplication 
(resp., inversion) operator in the real linear space of rational functions 
$\mathbb{R}[z,z^{-1}]$ when $L$-polynomials are chosen as basis.)

It is worth mentioning that a subset of the multitudinous applications 
$L$-polynomials and their associated Laurent-Jacobi matrices find in numerical 
analysis (quadrature formulae) and trigonometric moment problems \cite{a31}, 
the spectral theory of self-adjoint operators in infinite-dimensional 
(necessarily separable) Hilbert spaces \cite{a22,a23}, complex approximation 
theory (two-point Pad\'{e} approximants) \cite{a32,a33,a34}, the 
direct/inverse scattering theory for the (finite) relativistic Toda lattice 
\cite{a35} (see, also, \cite{a36}), and the (classical) Pick-Nevanlinna 
problem \cite{a37} are discussed in Section~1 of \cite{a38} (see, also, 
\cite{a25}, and the references therein). It turns out that, as a recurring 
theme, $n \! \to \! \infty$ asymptotics of $L$-polynomials are an essential 
calculational ingredient in analyses related to the above-mentioned, seemingly 
disparate, topics.

Now that the principal objects have been defined, namely, the monic OLPs, 
$\lbrace \boldsymbol{\pi}_{m}(z) \rbrace_{m \in \mathbb{Z}_{0}^{+}}$, and 
OLPs, $\lbrace \phi_{m}(z) \rbrace_{m \in \mathbb{Z}_{0}^{+}}$, it's time to 
state that the purpose of the present (three-fold) series of works, of which 
the present article constitutes Part~II, is to analyse the behaviour in the 
double-scaling limit as $\mathscr{N},n \! \to \! \infty$ such that $z_{o} \! 
:= \! \mathscr{N}/n \! = \! 1 \! + \! o(1)$ (the simplified `notation' $n \! 
\to \! \infty$ will be adopted) of the $L$-polynomials $\boldsymbol{\pi}_{n}
(z)$ and $\phi_{n}(z)$ in $\mathbb{C}$, orthogonal with respect to the varying 
exponential measure\footnote{Note that $LD(\boldsymbol{\pi}_{m}) \! = \! LD
(\phi_{m}) \! = \!
\begin{cases}
2n, &\text{$m \! = \! \text{even}$,} \\
2n \! + \! 1, &\text{$m \! = \! \text{odd}$,}
\end{cases}$ coincides with the parameter in the measure of orthogonality:
the large parameter, $n$, enters simultaneously into the $L$-degree of the
$L$-polynomials and the (varying exponential) weight; thus, asymptotics of
the $L$-polynomials are studied along a `diagonal strip' of a doubly-indexed
sequence.} $\md \mu (z) \! = \! \exp (-n \widetilde{V}(z)) \, \md z$, 
where $\widetilde{V}(z) \! := \! z_{o}V(z)$, and the `scaled' external 
field\footnote{For real non-analytic external fields, see the recent work 
\cite{a39}.} $\widetilde{V} \colon \mathbb{R} \setminus \{0\} \! \to \! 
\mathbb{R}$ satisfies conditions~(2.3)--(2.5) (see Subsection~2.2), as well 
as of the associated norming constants and coefficients of the (system of) 
recurrence relations; more precisely, then:
\begin{enumerate}
\item[\pmb{(i)}] in this work (Part~II), asymptotics (as $n \! \to \! \infty)$ 
of $\boldsymbol{\pi}_{2n+1}(z)$ (in the entire complex plane) and $\xi^{(2n
+1)}_{-n-1}$, thus $\phi_{2n+1}(z)$ (cf. Equation~(1.5)) and the Hankel 
determinant ratio $H^{(-2n)}_{2n+1}/H^{(-2n-2)}_{2n+2}$ (cf. Equations~(1.8)), 
are obtained;
\item[\pmb{(ii)}] in the previous work \cite{a38} (Part~I), asymptotics 
(as $n \! \to \! \infty)$ of $\boldsymbol{\pi}_{2n}(z)$ (in the entire 
complex plane) and $\xi^{(2n)}_{n}$, thus $\phi_{2n}(z)$ (cf. Equation~(1.4)) 
and the Hankel determinant ratio $H^{(-2n)}_{2n}/H^{(-2n)}_{2n+1}$ (cf. 
Equations~(1.8)), were obtained;
\item[\pmb{(iii)}] in Part~III \cite{a40}, asymptotics (as $n \! \to \! 
\infty)$ of $\nu^{(2n)}_{-n}$ $(= \! H^{(-2n+1)}_{2n}/H^{(-2n)}_{2n})$ and 
$\xi^{(2n)}_{-n}$, $\nu^{(2n+1)}_{n}$ $(= \! -H^{(-2n-1)}_{2n+1}/H^{(-2n)}_{2
n+1})$ and $\xi^{(2n+1)}_{n}$, $\prod_{k=1}^{2n} \alpha^{(2n)}_{k}$ $(=\nu^{
(2n)}_{-n})$, and $\prod_{k=1}^{2n+1} \alpha^{(2n+1)}_{k}$ $(= \! -(\nu^{(2n+
1)}_{n})^{-1})$, as well as of the (elements of the) Laurent-Jacobi matrices, 
$\mathcal{J}$ and $\mathcal{K}$, and other, related, quantities constructed 
{}from the coefficients of the three- and five-term recurrence relations, are 
obtained.
\end{enumerate}

The above-mentioned asymptotics (as $n \! \to \! \infty)$ are obtained by 
reformulating, \emph{\`{a} la} Fokas-Its-Kitaev \cite{a41,a42}, the 
corresponding `even degree' and `odd degree' monic $L$-polynomial problems as 
(matrix) Riemann-Hilbert problems (RHPs) on $\mathbb{R}$, and then studying 
the large-$n$ behaviour of the corresponding solutions. The paradigm for the 
asymptotic (as $n \! \to \! \infty)$ analysis of the respective (matrix) 
RHPs is a union of the Deift-Zhou (DZ) non-linear steepest-descent method 
\cite{a1,a2}, used for the asymptotic analysis of undulatory RHPs, and the 
extension of Deift-Venakides-Zhou \cite{a3}, incorporating into the DZ method 
a non-linear analogue of the WKB method, making the asymptotic analysis of 
fully non-linear problems tractable (it should be mentioned that, in this 
context, the equilibrium measure \cite{a43} plays an absolutely crucial 
r\^{o}le in the analysis \cite{a44}); see, also, the multitudinous extensions 
and applications of the DZ method \cite{a45,a46,a47,a48,a49,a50,a51,a52,a53,%
a54,a55,a56,a57,a58,a59,a60,a61,a62,a63,a64,a65,a66,a67,a68,a69}. It is worth 
mentioning that asymptotics for Laurent-type polynomials and their zeros have 
been obtained in \cite{a33,a70} (see, also, \cite{a71,a72,a73}).

Unlike the large-$n$ asymptotic analysis for the orthogonal polynomials case, 
which is related to one (matrix) RHP normalised at the point at infinity, the 
large-$n$ asymptotic analysis for the OLPs requires the consideration of two 
different families of (matrix) RHPs, one for even degree (see Subsection~2.2, 
\textbf{RHP1} and Lemma~2.2.1), and one for odd degree (see Subsection~2.2, 
\textbf{RHP2} and Lemma~2.2.2): \textbf{RHP1} is normalised at the point at 
infinity, whereas \textbf{RHP2} is normalised at $z \! = \! 0$. The technical 
details, therefore, related to the large-$n$ asymptotic analyses of 
\textbf{RHP1} and \textbf{RHP2} are different, and must be carried through 
independently. Further, important, albeit technical, distinctions are:
\begin{compactenum}
\item[\textbullet] the associated $g$-functions are different (see \cite{a38}, 
Equation~(2.12), for the even degree case, and Equation~(2.13) of the present 
article for the odd degree case);
\item[\textbullet] the respective variational problems are different, which 
means that the supports of the associated equilibrium measures are different 
(see \cite{a38}, Lemmas~3.1--3.3, for the even degree case, and Lemmas 
3.1--3.3 of the present article for the odd degree case);
\item[\textbullet] even though the supports consist of the union of a finite 
number of compact real intervals, the systems of transcendental equations 
(finite in number) which characterise the end-points of the supports of the 
respective equilibrium measures are different (see \cite{a38}, Lemmas~3.5  
and~3.6, for the even degree case, and Lemmas~3.5 and~3.6 of the present 
article for the odd degree case);
\item[\textbullet] the associated `small-norm' RHPs, due to the difference in 
normalisations, are different, which gives rise to the appearance of several 
new---residue---terms in the solution of the odd degree small-norm RHP, 
which do not arise in the solution of the even degree small-norm RHP (see, in 
particular, \cite{a38}, Lemmas~4.8 and~5.2, for the even degree case, and 
Lemmas~4.8 and~5.2 of the present article for the odd degree case); and
\item[\textbullet] certain error analyses are more complicated for the odd 
degree case, because the difference in normalisations requires that twice 
as many matrix-operator norms be estimated for the odd degree case, which
makes the asymptotic analysis more tedious and involved (see, in particular, 
\cite{a38}, Proposition~5.2 and Lemma~5.2, for the even degree case, and 
Proposition~5.2 and Lemma~5.2 of the present article for the odd degree case).
\end{compactenum}
For these reasons, as well other, technical ones, attempting to meld together 
the even degree and odd degree cases into one article would hamper the 
readability of an already lengthy analysis. So, despite the repeated, 
over-arching scheme of analysis, the large-$n$ asymptotic behaviour for the 
even degree OLPs (and related quantities) is studied in \cite{a38}, and the 
large-$n$ asymptotic behaviour for the odd degree OLPs (and related 
quantities) is the subject of the present article.

This article is organised as follows. In Section~2, necessary facts {}from 
the theory of (compact) Riemann surfaces are given, the respective `even 
degree' and `odd degree' RHPs on $\mathbb{R}$ are stated and the corresponding 
variational problems for the associated equilibrium measures are discussed, 
and the main results of this work, namely, asymptotics (as $n \! \to \! 
\infty)$ of $\boldsymbol{\pi}_{2n+1}(z)$ (in $\mathbb{C})$, and $\xi^{(2n+1)
}_{-n-1}$ and $\phi_{2n+1}(z)$ (in $\mathbb{C})$ are stated in Theorems~2.3.1 
and~2.3.2, respectively. In Section~3, the detailed analysis of the `odd 
degree' variational problem and the associated equilibrium measure is 
undertaken, including the construction of the so-called $g$-function, and 
the RHP formulated in Section~2 is reformulated as an equivalent, auxiliary 
RHP, which, in Sections~4 and~5, is augmented, by means of a sequence of 
contour deformations and transformations \emph{\`{a} la} Deift-Venakides-Zhou, 
into simpler, `model' (matrix) RHPs which, as $n \! \to \! \infty$, and in 
conjunction with the Beals-Coifman construction \cite{a74} (see, also, the 
extension of Zhou \cite{a75}) for the integral representation of the solution 
of a matrix RHP on an oriented contour, are solved explicitly (in closed 
form) in terms of Riemann theta functions (associated with the underlying 
finite-genus hyperelliptic Riemann surface) and Airy functions, {}from 
which the final asymptotic (as $n \! \to \! \infty)$ results stated in 
Theorems~2.3.1 and~2.3.2 are proved. The paper concludes with an Appendix.
\begin{eeeee}
The even degree OLPs, $\phi_{2n}(z)$, $n \! \in \! \mathbb{Z}_{0}^{+}$, 
are related, in a way, to the polynomials orthogonal with respect to the 
varying weight $\widehat{w}(z) \! := \! z^{-2n} \exp (-\mathscr{N} \, 
V(z))$, $\mathscr{N} \! \in \! \mathbb{N}$: this follows directly {}from the 
orthogonality relation satisfied by $\phi_{2n}(z)$. This does not help with 
any of the algebraic relations, such as the system of three- and five-term 
recurrence relations; however, this does provide for an alternative approach 
to computing large-$n$ asymptotics for $\phi_{2n}(z)$. The connection is not 
so clear for the odd degree OLPs, $\phi_{2n+1}(z)$, $n \! \in \! \mathbb{Z}_{
0}^{+}$. Indeed, in this latter case, the associated (density of the) measure 
for the orthogonal polynomials would take the form $\md \widehat{\mu}(z) \! 
:= \! z^{-2n-1} \exp (-\mathscr{N} \, V(z)) \, \md z$, and this measure 
changes signs, which causes a number of difficulties in the large-$n$ 
asymptotic analysis. In this paper, these connections are not used, and a 
complete asymptotic analysis of the odd degree OLPs is carried out, directly. 
\hfill $\blacksquare$
\end{eeeee}
\section{Hyperelliptic Riemann Surfaces, The Riemann-Hi\-l\-b\-e\-r\-t
P\-r\-o\-b\-l\-e\-m\-s, and Su\-m\-m\-a\-r\-y of Results}
In this section, necessary facts {}from the theory of hyperelliptic Riemann 
surfaces are given (see Subsection~2.1), the respective RHPs on $\mathbb{R}$ 
for the even degree and the odd degree monic orthogonal $L$-polynomials are 
formulated and the corresponding variational problems for the associated 
equilibrium measures are discussed (see Subsection~2.2), and the asymptotics 
(as $n \! \to \! \infty)$ for $\boldsymbol{\pi}_{2n+1}(z)$ (in the entire 
complex plane), and $\xi^{(2n+1)}_{-n-1}$ and $\phi_{2n+1}(z)$ (in the entire 
complex plane) are given in Theorems~2.3.1 and~2.3.2, respectively (see 
Subsection~2.3).

Before proceeding, however, the notation/nomenclature used throughout this 
work is summarised.
\begin{center}
\Ovalbox{\textsc{Notational Conventions}}
\end{center}
\begin{compactenum}
\item[(1)] $\mathrm{I} \! = \!
\left(
\begin{smallmatrix}
1 & 0 \\
0 & 1
\end{smallmatrix}
\right)$ is the $2 \times 2$ identity matrix, $\sigma_{1} \! = \!
\left(
\begin{smallmatrix}
0 & 1 \\
1 & 0
\end{smallmatrix}
\right)$, $\sigma_{2} \! = \!
\left(
\begin{smallmatrix}
0 & -\mi \\
\mi & 0
\end{smallmatrix}
\right)$ and $\sigma_{3} \! = \!
\left(
\begin{smallmatrix}
1 & 0 \\
0 & -1
\end{smallmatrix}
\right)$ are the Pauli matrices, $\sigma_{+} \! = \!
\left(
\begin{smallmatrix}
0 & 1 \\
0 & 0
\end{smallmatrix}
\right)$ and $\sigma_{-} \! = \!
\left(
\begin{smallmatrix}
0 & 0 \\
1 & 0
\end{smallmatrix}
\right)$ are, respectively, the raising and lowering matrices, $\mathbf{0} \!
= \!
\left(
\begin{smallmatrix}
0 & 0 \\
0 & 0
\end{smallmatrix}
\right)$, $\mathbb{R}_{\pm} \! := \! \lbrace \mathstrut x \! \in \! \mathbb{
R}; \, \pm x \! > \! 0 \rbrace$, $\mathbb{C}_{\pm} \! := \! \lbrace \mathstrut
z \! \in \! \mathbb{C}; \, \pm \Im (z) \! > \! 0 \rbrace$, and $\mathrm{sgn}
(x) \! := \! 0$ if $x \! = \! 0$ and $x \vert x \vert^{-1}$ if $x \! \not= \! 
0$;
\item[(2)] for a scalar $\omega$ and a $2 \! \times \! 2$ matrix $\Upsilon$,
$\omega^{\mathrm{ad}(\sigma_{3})} \Upsilon \! := \! \omega^{\sigma_{3}} 
\Upsilon \omega^{-\sigma_{3}}$;
\item[(3)] a contour $\mathcal{D}$ which is the finite union of 
piecewise-smooth, simple curves (as closed sets) is said to be 
\emph{orientable} if its complement $\mathbb{C} \setminus \mathcal{D}$ can 
always be divided into two, possibly disconnected, disjoint open sets $\mho^{
+}$ and $\mho^{-}$, either of which has finitely many components, such that 
$\mathcal{D}$ admits an orientation so that it can either be viewed as a 
positively oriented boundary $\mathcal{D}^{+}$ for $\mho^{+}$ or as a 
negatively oriented boundary $\mathcal{D}^{-}$ for $\mho^{-}$ \cite{a75}, 
that is, the (possibly disconnected) components of $\mathbb{C} \setminus 
\mathcal{D}$ can be coloured by $+$ or $-$ in such a way that the $+$ regions 
do not share boundary with the $-$ regions, except, possibly, at finitely 
many points \cite{a76};
\item[(4)] for each segment of an oriented contour $\mathcal{D}$, according to
the given orientation, the ``+'' side is to the left and the ``-'' side is to
the right as one traverses the contour in the direction of orientation, that
is, for a matrix $\mathcal{A}_{ij}(z)$, $i,j \! = \! 1,2$, $(\mathcal{A}_{
ij}(z))_{\pm}$ denote the non-tangential limits $(\mathcal{A}_{ij}(z))_{
\pm} \! := \! \lim_{\genfrac{}{}{0pt}{2}{z^{\prime} \, \to \, z}{z^{\prime} \,
\in \, \pm \, \mathrm{side} \, \mathrm{of} \, \mathcal{D}}} \mathcal{A}_{ij}
(z^{\prime})$;
\item[(5)] for $1 \! \leqslant \! p \! < \! \infty$ and $\mathcal{D}$ some
point set,
\begin{equation*}
\mathcal{L}^{p}_{\mathrm{M}_{2}(\mathbb{C})}(\mathcal{D}) \! := \! \left\{
\mathstrut f \colon \mathcal{D} \! \to \! \mathrm{M}_{2}(\mathbb{C}); \,
\vert \vert f(\cdot) \vert \vert_{\mathcal{L}^{p}_{\mathrm{M}_{2}(\mathbb{
C})}(\mathcal{D})} \! := \! \left(\int_{\mathcal{D}} \vert f(z) \vert^{p} \,
\vert \md z \vert \right)^{1/p} \! < \! \infty \right\},
\end{equation*}
where, for $\mathcal{A}(z) \! \in \! \operatorname{M}_{2}(\mathbb{C})$, $\vert
\mathcal{A}(z) \vert \! := \! (\sum_{i,j=1}^{2} \overline{\mathcal{A}_{ij}(z)}
\, \mathcal{A}_{ij}(z))^{1/2}$ is the Hilbert-Schmidt norm, with $\overline{
\bullet}$ denoting complex conjugation of $\bullet$, for $p \! = \! \infty$,
\begin{equation*}
\mathcal{L}^{\infty}_{\mathrm{M}_{2}(\mathbb{C})}(\mathcal{D}) \! := \! \left\{
\mathstrut g \colon \mathcal{D} \! \to \! \mathrm{M}_{2}(\mathbb{C}); \, \vert
\vert g(\cdot) \vert \vert_{\mathcal{L}^{\infty}_{\mathrm{M}_{2}(\mathbb{C})}
(\mathcal{D})} \! := \! \max_{i,j = 1,2} \, \, \sup_{z \in \mathcal{D}} \vert
g_{ij}(z) \vert \! < \! \infty \right\},
\end{equation*}
and, for $f \! \in \! \mathrm{I} \! + \! \mathcal{L}^{2}_{\mathrm{M}_{2}
(\mathbb{C})}(\mathcal{D}) \! := \! \left\lbrace \mathstrut \mathrm{I} \! + \!
h; \, h \! \in \! \mathcal{L}^{2}_{\mathrm{M}_{2}(\mathbb{C})}(\mathcal{D})
\right\rbrace$,
\begin{equation*}
\vert \vert f(\cdot) \vert \vert_{\mathrm{I}+\mathcal{L}^{2}_{\mathrm{M}_{2}
(\mathbb{C})}(\mathcal{D})} \! := \! \left(\vert \vert f(\infty) \vert \vert_{
\mathcal{L}^{\infty}_{\mathrm{M}_{2}(\mathbb{C})}(\mathcal{D})}^{2} \! + \!
\vert \vert f(\cdot) \! - \! f(\infty) \vert \vert_{\mathcal{L}^{2}_{\mathrm{
M}_{2}(\mathbb{C})}(\mathcal{D})}^{2} \right)^{1/2};
\end{equation*}
\item[(6)] for a matrix $\mathcal{A}_{ij}(z)$, $i,j \! = \! 1,2$, to have
boundary values in the $\mathcal{L}^{2}_{\mathrm{M}_{2}(\mathbb{C})}(\mathcal{
D})$ sense on an oriented contour $\mathcal{D}$, it is meant that $\lim_{
\genfrac{}{}{0pt}{2}{z^{\prime} \, \to \, z}{z^{\prime} \, \in \, \pm \,
\mathrm{side} \, \mathrm{of} \, \mathcal{D}}} \int_{\mathcal{D}} \vert
\mathcal{A}(z^{\prime}) \! - \! (\mathcal{A}(z))_{\pm} \vert^{2} \, \vert \md
z \vert \! = \! 0$ (e.g., if $\mathcal{D} \! = \! \mathbb{R}$ is oriented
{}from $+\infty$ to $-\infty$, then $\mathcal{A}(z)$ has $\mathcal{L}^{2}_{
\mathrm{M}_{2}(\mathbb{C})}(\mathcal{D})$ boundary values on $\mathcal{D}$
means that $\lim_{\varepsilon \downarrow 0} \int_{\mathbb{R}} \vert \mathcal{A}
(x \! \mp \! \mi \varepsilon) \! - \! (\mathcal{A}(x))_{\pm} \vert^{2} \, \md
x \! = \! 0)$;
\item[(7)] for a $2 \times 2$ matrix-valued function $\mathfrak{Y}(z)$, the
notation $\mathfrak{Y}(z) \! =_{z \to z_{0}} \! \mathcal{O}(\ast)$ means
$\mathfrak{Y}_{ij}(z) \! =_{z \to z_{0}} \! \mathcal{O}(\ast_{ij})$, $i,j \! =
\! 1,2$ (\emph{mutatis mutandis} for $o(1))$;
\item[(8)] $\vert \vert \mathscr{F}(\cdot) \vert \vert_{\cap_{p \in J}
\mathcal{L}^{p}_{\mathrm{M}_{2}(\mathbb{C})}(\ast)} \! := \! \sum_{p \in J}
\vert \vert \mathscr{F}(\cdot) \vert \vert_{\mathcal{L}^{p}_{\mathrm{M}_{2}
(\mathbb{C})}(\ast)}$, with $\mathrm{card}(J) \! < \! \infty$;
\item[(9)] $\mathcal{M}_{1}(\mathbb{R})$ denotes the set of all non-negative,
bounded, unit Borel measures on $\mathbb{R}$ for which all moments exist,
\begin{equation*}
\mathcal{M}_{1}(\mathbb{R}) \! := \! \left\{\mathstrut \mu; \, \int_{\mathbb{
R}} \md \mu (s) \! = \! 1, \, \int_{\mathbb{R}}s^{m} \md \mu (s) \! < \!
\infty, \, m \! \in \! \mathbb{Z} \setminus \{0\} \right\};
\end{equation*}
\item[(10)] for $(\mu,\nu) \! \in \! \mathbb{R} \times \mathbb{R}$, denote the
function $(\bullet \! - \! \mu)^{\mi \nu} \colon \mathbb{C} \setminus (-\infty,
\mu) \! \to \! \mathbb{C}$, $\bullet \! \mapsto \! \exp (\mi \nu \ln (\bullet
-\mu))$, where $\ln$ denotes the principal branch of the logarithm;
\item[(11)] for $\widetilde{\pmb{\gamma}}$ a nullhomologous path in a region
$\mathscr{D} \subset \mathbb{C}$, $\mathrm{int}(\widetilde{\pmb{\gamma}}) \!
:= \! \left\lbrace \mathstrut \zeta \! \in \! \mathscr{D} \setminus
\widetilde{\pmb{\gamma}}; \, \mathrm{ind}_{\widetilde{\pmb{\gamma}}}(\zeta) \!
:= \! \int_{\widetilde{\pmb{\gamma}}} \tfrac{1}{z-\zeta} \, \tfrac{\md z}{2
\pi \mi} \! \not=\! 0 \right\rbrace$;
\item[(12)] for some point set $\mathcal{D} \subset \mathcal{X}$, with
$\mathcal{X} \! = \! \mathbb{C}$ or $\mathbb{R}$, $\overline{\mathcal{D}}
:= \! \mathcal{D} \cup \partial \mathcal{D}$, and $\mathcal{D}^{c} \! := \!
\mathcal{X} \setminus \overline{\mathcal{D}}$.
\end{compactenum}
\subsection{Riemann Surfaces: Preliminaries}
In this subsection, the basic elements associated with the construction of 
hyperelliptic and finite genus (compact) Riemann surfaces are presented (for 
further details and proofs, see, for example, \cite{a77,a78}).
\begin{eeee}
The superscripts ${}^{\pm}$, and sometimes the subscripts ${}_{\pm}$, in 
this subsection should not be confused with the subscripts ${}_{\pm}$ 
appearing in the various RHPs (this is a general comment which applies, 
unless stated otherwise, throughout the entire text). Although $\overline{
\mathbb{C}}$ (or $\mathbb{C} \mathbb{P}^{1})$ $:= \! \mathbb{C} \cup 
\{\infty\}$ (resp., $\overline{\mathbb{R}} \! := \! \mathbb{R} \cup \{-
\infty\} \cup \{+\infty\})$ is the standard definition for the (closed) 
Riemann sphere (resp., closed real line), the simplified, and somewhat 
abusive, notation $\mathbb{C}$ (resp., $\mathbb{R})$ is used to denote both 
the (closed) Riemann sphere, $\overline{\mathbb{C}}$ (resp., closed real 
line, $\overline{\mathbb{R}})$, and the (open) complex field, $\mathbb{C}$ 
(resp., open real line, $\mathbb{R})$, and the context(s) should make clear 
which object(s) the notation $\mathbb{C}$ (resp., $\mathbb{R})$ represents. 
\hfill $\blacksquare$
\end{eeee}

Let $N \! \in \! \mathbb{N}$ (with $N \! < \! \infty$ assumed throughout) and
$\varsigma_{k} \! \in \! \mathbb{R} \setminus \lbrace 0,\pm \infty \rbrace$,
$k \! = \! 1,\dotsc,2N \! + \! 2$, be such that $\varsigma_{i} \! \not= \!
\varsigma_{j} \, \, \forall \, \, i \! \not= \! j \! = \! 1,\dotsc,2N \! + \!
2$, and enumerated/ordered according to $\varsigma_{1} \! < \! \varsigma_{2}
\! < \! \cdots \! < \! \varsigma_{2N+2}$. Let $R(z) \! := \! \prod_{j=1}^{N}
(z \! - \! \varsigma_{2j-1})(z \! - \! \varsigma_{2j})$ $\in \! \mathbb{R}[z]$
(the algebra of polynomials in $z$ with coefficients in $\mathbb{R})$ be the
(unital) polynomial of even degree $\mathrm{deg}(R) \! = \! 2N \! + \! 2$
$(\mathrm{deg}(R) \! = \! 0$ $(\mathrm{mod} 2))$ whose (simple) zeros/roots
are $\{\varsigma_{j}\}_{j=1}^{2N+2}$. Denote by $\mathscr{R}$ the
hyperelliptic Riemann surface of genus $N$ defined by the equation $y^{2} \!
= \! R(z)$ and realised as a two-sheeted branched (ramified) covering of the
Riemann sphere such that its two sheets are two identical copies of $\mathbb{
C}$ with branch cuts along the intervals $(\varsigma_{1},\varsigma_{2})$,
$(\varsigma_{3},\varsigma_{4})$, $\dotsc$, $(\varsigma_{2N+1},\varsigma_{2N+
2})$, and glued/pasted to each other `crosswise' along the opposite banks of
the corresponding cuts $(\varsigma_{2j-1},\varsigma_{2j})$, $j \! = \! 1,
\dotsc,N \! + \! 1$. Denote the two sheets of $\mathscr{R}$ by $\mathscr{R}^{
+}$ (the first/upper sheet) and  $\mathscr{R}^{-}$ (the second/lower sheet):
to indicate that $z$ lies on the first (resp., second) sheet, one writes $z^{
+}$ (resp., $z^{-})$; of course, as points in the plane $\mathbb{C}$, $z^{+}
\! = \! z^{-} \! = \! z$. For points $z$ on the first (resp., second) sheet
$\mathscr{R}^{+}$ (resp., $\mathscr{R}^{-})$, one has that $z^{+} \! = \!
(z,+(R(z))^{1/2})$ (resp., $z^{-} \! = \! (z,-(R(z))^{1/2}))$, where the
single-valued branch of the square root is chosen such that $z^{-(N+1)}(R(z))^{
1/2} \! \sim_{\underset{z \in \mathscr{R}^{\pm}}{z \to \infty}} \! \pm 1$.

Let $\mathscr{E}_{j} \! := \! (\varsigma_{2j-1},\varsigma_{2j})$, $j \! = \!
1,\dotsc,N \! + \! 1$, and set $\mathscr{E} \! = \! \cup_{j=1}^{N+1} \mathscr{
E}_{j}$ (note that $\mathscr{E}_{i} \cap \mathscr{E}_{j} \! = \! \varnothing$,
$i \! \not= \! j \! = \! 1,\dotsc,N \! + \! 1)$. Denote by $\mathscr{E}_{j}^{
+}$ $(\subset \mathscr{R}^{+})$ (resp., $\mathscr{E}_{j}^{-}$ $(\subset 
\mathscr{R}^{-}))$ the upper (resp., lower) bank of the interval $\mathscr{
E}_{j}$, $j \! = \! 1,\dotsc,N \! + \! 1$, forming $\mathscr{E}$, and oriented 
in accordance with the orientation of $\mathscr{E}$ as the boundary of 
$\mathbb{C} \setminus \mathscr{E}$, namely, the domain $\mathbb{C} \setminus 
\mathscr{E}$ is on the left as one proceeds along the upper bank of the $j$th 
interval {}from $\varsigma_{2j-1}$ to the point $\varsigma_{2j}$ and back 
along the lower bank from $\varsigma_{2j}$ to $\varsigma_{2j-1}$; thus, 
$\mathscr{E}_{j}^{\pm} \! := \! (\varsigma_{2j-1},\varsigma_{2j})^{\pm}$, $j 
\! = \! 1,\dotsc,N \! + \! 1$, are two (identical) copies of $(\varsigma_{2
j-1},\varsigma_{2j}) \subset \mathbb{R}$ `lifted' to $\mathscr{R}^{\pm}$. Set 
$\Gamma \! := \! \cup_{j=1}^{N+1} \Gamma_{j}$ $(\subset \mathscr{R})$, where 
$\Gamma_{j} \! := \! \mathscr{E}_{j}^{+} \cup \mathscr{E}_{j}^{-}$, $j \! = 
\! 1,\dotsc,N \! + \! 1$ $(\Gamma \! = \! \mathscr{E}^{+} \cup \mathscr{E}^{
-})$: note that $\Gamma$, as a curve on $\mathscr{R}$ (defined by the equation 
$y^{2} \! = \! R(z))$, consists of a finitely denumerable number of disjoint 
analytic closed Jordan curves, $\Gamma_{j}$, $j \! = \! 1,\dotsc,N \! + \! 1$, 
which are \emph{cycles} on $\mathscr{R}$, and that correspond to the intervals 
$\mathscr{E}_{j}$. {}From the above construction, it is clear that $\mathscr{
R} \! = \! \mathscr{R}^{+} \cup \mathscr{R}^{-} \cup \Gamma$; furthermore, the 
canonical projection of $\Gamma$ onto $\mathbb{C}$ $(\boldsymbol{\pi} \colon 
\mathscr{R} \! \to \! \mathbb{C})$ is $\mathscr{E}$, that is, $\boldsymbol{
\pi}(\Gamma) \! = \! \mathscr{E}$ (also, $\boldsymbol{\pi}(\mathscr{R}^{+}) 
\! = \! \boldsymbol{\pi}(\mathscr{R}^{-}) \! = \! \mathbb{C} \setminus 
\mathscr{E}$, or, alternately, $\boldsymbol{\pi}(z^{+}) \! = \! \boldsymbol{
\pi}(z^{-}) \! = \! z)$. One moves in the `positive $(+)$' (resp., `negative 
$(-)$') direction along the (closed) contour $\Gamma \subset \mathscr{R}$ if 
the domain $\mathscr{R}^{+}$ is on the left (resp., right) and the domain 
$\mathscr{R}^{-}$ is on the right (resp., left): the corresponding notation 
is (see above) $\Gamma^{+}$ (resp., $\Gamma^{-})$. For a function $f$ defined 
on the two-sheeted hyperelliptic Riemann surface $\mathscr{R}$, one defines 
the non-tangential boundary values, provided they exist, of $f(z)$ as $z \! 
\in \! \mathscr{R}^{+}$ (resp., $z \! \in \! \mathscr{R}^{-})$ approaches 
$\lambda \! \in \! \Gamma$, denoted $\lambda_{+}$ (resp., $\lambda_{-})$, by 
$f(\lambda_{\pm}) \! := \! f_{\pm}(\lambda) \! := \! \lim_{\underset{z \in 
\Gamma^{\pm}}{z \to \lambda}}f(z)$.

One takes the first $N$ contours among the (closed) contours $\Gamma_{j}$ for
basis $\boldsymbol{\alpha}$-cycles $\lbrace \boldsymbol{\alpha}_{j}, \, j \! =
\! 1,\dotsc,N \rbrace$ and then completes/supplements this in the standard way
with $\boldsymbol{\beta}$-cycles $\lbrace \boldsymbol{\beta}_{j}, \, j \! = \!
1,\dotsc,N \rbrace$ so that the \emph{intersection matrix} has the (canonical)
form $\boldsymbol{\alpha}_{k} \circ \boldsymbol{\alpha}_{j} \! = \!
\boldsymbol{\beta}_{k} \circ \boldsymbol{\beta}_{j} \! = \! 0 \, \, \forall \,
\, k \! \not= \! j \! = \! 1,\dotsc,N$, and $\boldsymbol{\alpha}_{k} \circ
\boldsymbol{\beta}_{j} \! = \! \delta_{kj}$: the cycles $\lbrace \boldsymbol{
\alpha}_{j},\boldsymbol{\beta}_{j} \rbrace$, $j \! = \! 1,\dotsc,N$, form
the \emph{canonical 1-homology basis} on $\mathscr{R}$, namely, any cycle
$\widehat{\boldsymbol{\gamma}} \subset \mathscr{R}$ is homologous to an
integral linear combination of $\lbrace \boldsymbol{\alpha}_{j},\boldsymbol{
\beta}_{j} \rbrace$, that is, $\widehat{\boldsymbol{\gamma}} \! = \! \sum_{j=
1}^{N}(n_{j} \boldsymbol{\alpha}_{j} \! + \! m_{j} \boldsymbol{\beta}_{j})$,
where $(n_{j},m_{j}) \! \in \! \mathbb{Z} \times \mathbb{Z}$, $j \! = \! 1,
\dotsc,N$. The $\boldsymbol{\alpha}$-cycles $\lbrace \boldsymbol{\alpha}_{j},
\, j \! = \! 1,\dotsc,N \rbrace$, in the present case, are the intervals
$(\varsigma_{2j-1},\varsigma_{2j})$, $j \! = \! 1,\dotsc,N$, `going twice',
that is, along the upper ({}from $\varsigma_{2j-1}$ to $\varsigma_{2j})$ and
lower ({}from $\varsigma_{2j}$ to $\varsigma_{2j-1})$ banks $(\alpha_{j} \! 
= \! \mathscr{E}_{j}^{+} \cup \mathscr{E}_{j}^{-}$, $j \! = \! 1,\dotsc,N)$, 
and the $\boldsymbol{\beta}$-cycles $\lbrace \boldsymbol{\beta}_{j}, \, j \! 
= \! 1,\dotsc,N \rbrace$ are as follows: the $j$th $\boldsymbol{\beta}$-cycle
consists of the $\boldsymbol{\alpha}$-cycles $\boldsymbol{\alpha}_{k}$, $k \!
= \! j \! + \! 1,\dotsc,N$, and the cycles `linked' with them and consisting
of (the gaps) $(\varsigma_{2k},\varsigma_{2k+1})$, $k \! = \! 1,\dotsc,N$,
`going twice', that is, {}from $\varsigma_{2k}$ to $\varsigma_{2k+1}$ on the
first sheet and in the reverse direction on the second sheet. For an arbitrary
holomorphic Abelian differential (one-form) $\boldsymbol{\omega}$ on
$\mathscr{R}$, the function $\int^{z} \boldsymbol{\omega}$ is defined uniquely
modulo its $\boldsymbol{\alpha}$- and $\boldsymbol{\beta}$-periods, $\oint_{
\boldsymbol{\alpha}_{j}} \boldsymbol{\omega}$ and $\oint_{\boldsymbol{\beta}_{
j}} \boldsymbol{\omega}$, $j \! = \! 1,\dotsc,N$, respectively. It is well
known that the canonical $1$-homology basis $\lbrace \boldsymbol{\alpha}_{j},
\boldsymbol{\beta}_{j} \rbrace$, $j \! = \! 1,\dotsc,N$, constructed above
`generates', on $\mathscr{R}$, the corresponding 
$\boldsymbol{\alpha}$-normalised basis of holomorphic Abelian differentials 
(one-forms) $\lbrace \omega_{1},\omega_{2},\dotsc,\omega_{N} \rbrace$, where 
$\omega_{j} \! := \! \sum_{k=1}^{N} \tfrac{c_{jk}z^{N-k}}{\sqrt{\smash[b]{R
(z)}}} \, \md z$, $c_{jk} \! \in \! \mathbb{C}$, $j \! = \! 1,\dotsc,N$, and 
$\oint_{\boldsymbol{\alpha}_{k}} \omega_{j} \! = \! \delta_{kj}$, $k, \, 
j \! = \! 1,\dotsc,N$: the associated $N \times N$ matrix of 
$\boldsymbol{\beta}$-periods, $\tau \! = \! (\tau_{ij})_{i,j=1,\dotsc,N} \! 
:= \! \left(\oint_{\boldsymbol{\beta}_{j}} \omega_{i} \right)_{i,j=1,\dotsc,
N}$, is a \emph{Riemann matrix}, that is, it is symmetric $(\tau_{ij} \! = 
\! \tau_{ji})$, pure imaginary, and $-\mi \tau$ is positive definite $(\Im 
(\tau_{ij}) \! > \! 0)$; moreover, $\tau$ is non-degenerate $(\det (\tau) 
\! \not= \! 0)$. {}From the condition that the basis of the differentials 
$\omega_{l}$, $l \! = \! 1,\dotsc,N$, is canonical, with respect to the 
given basis cycles $\lbrace \boldsymbol{\alpha}_{j},\boldsymbol{\beta}_{j} 
\rbrace$, it is seen that this implies that each $\omega_{l}$ is real valued 
on $\mathscr{E} \! = \! \cup_{j=1}^{N+1}(\varsigma_{2j-1},\varsigma_{2j})$ 
and has exactly one (real) root/zero in any interval (band) $(\varsigma_{2j
-1},\varsigma_{2j})$, $j \! = \! 1,\dotsc,N \! + \! 1$, $j \! \not= \! l$; 
moreover, in the `gaps' $(\varsigma_{2j},\varsigma_{2j+1})$, $j \! = \! 1,
\dotsc,N$, these differentials take non-zero, purely imaginary values.

Fix the `standard basis' $\boldsymbol{e}_{1},\boldsymbol{e}_{2},\dotsc,
\boldsymbol{e}_{N}$ in $\mathbb{R}^{N}$, that is, $(\boldsymbol{e}_{j})_{k}
\! = \! \delta_{jk}$, $j,k \! = \! 1,\dotsc,N$ (these standard basis vectors
should be viewed as column vectors): the vectors $\boldsymbol{e}_{1},
\boldsymbol{e}_{2},\dotsc,\boldsymbol{e}_{N},\tau \boldsymbol{e}_{1},\tau
\boldsymbol{e}_{2},\dotsc,\tau \boldsymbol{e}_{N}$ are linearly independent
over the real field $\mathbb{R}$, and form a `basis' in $\mathbb{C}^{N}$.
The quotient space $\mathbb{C}^{N}/\{N \! + \! \tau M\}$, $(\mathrm{N},
\mathrm{M}) \! \in \! \mathbb{Z}^{N} \times \mathbb{Z}^{N}$, where $\mathbb{
Z}^{N} \! := \! \lbrace \mathstrut (m_{1},m_{2},\dotsc,m_{N}); \, m_{j} \!
\in \! \mathbb{Z}, \, j \! = \! 1,\dotsc,N \rbrace$, is a $2N$-dimensional
real torus $\mathbb{T}^{2N}$, and is referred to as the \emph{Jacobi variety},
symbolically $\mathrm{Jac}(\mathscr{R})$, of the two-sheeted (hyperelliptic)
Riemann surface $\mathscr{R}$ of genus $N$. Let $z_{0}$ be a fixed point in
$\mathscr{R}$. A vector-valued function $\boldsymbol{\mathscr{A}}(z) \! = \!
(\mathscr{A}_{1}(z),\mathscr{A}_{2}(z),\dotsc,\mathscr{A}_{N}(z)) \! \in \!
\mathrm{Jac}(\mathscr{R})$ with co-ordinates $\mathscr{A}_{k}(z) \! \equiv \!
\int_{z_{0}}^{z} \omega_{k}$, $k \! = \! 1,\dotsc,N$, where, hereafter, unless
stated otherwise and/or where confusion may arise, $\equiv$ denotes
`congruence modulo the period lattice', defines the \emph{Abel map}
$\boldsymbol{\mathscr{A}} \colon \mathscr{R} \! \to \! \mathrm{Jac}(\mathscr{
R})$. The unordered set of points $z_{1},z_{2},\dotsc,z_{N}$, with $z_{k} \!
\in \! \mathscr{R}$, form the $N$th symmetric power of $\mathscr{R}$,
symbolically $\mathscr{R}^{N}_{\mathrm{symm}}$ (or $\mathscr{S}^{N} \mathscr{
R})$. The vector function $\boldsymbol{\mathfrak{U}} \! = \! (\mathfrak{U}_{
1},\mathfrak{U}_{2},\dotsc,\mathfrak{U}_{N})$ with co-ordinates $\mathfrak{
U}_{j} \! = \! \mathfrak{U}_{j}(z_{1},z_{2},\dotsc,z_{N}) \! \equiv \! \sum_{
k=1}^{N} \mathscr{A}_{j}(z_{k}) \! \equiv \! \sum_{k=1}^{N} \int_{z_{0}}^{z_{
k}} \omega_{j}$, $j \! = \! 1,\dotsc,N$, that is, $(z_{1},z_{2},\dotsc,z_{N})
\! \to \! (\sum_{k=1}^{N} \int_{z_{0}}^{z_{k}} \omega_{1},\sum_{k=1}^{N}
\int_{z_{0}}^{z_{k}} \omega_{2},\dotsc,\sum_{k=1}^{N} \int_{z_{0}}^{z_{k}}
\omega_{N})$, is also referred to as the \emph{Abel map}, $\boldsymbol{
\mathfrak{U}} \colon \mathscr{R}^{N}_{\mathrm{symm}} \! \to \! \mathrm{Jac}
(\mathscr{R})$ (or $\boldsymbol{\mathfrak{U}} \colon \mathscr{S}^{N} \mathscr{
R} \! \to \! \mathrm{Jac}(\mathscr{R}))$. It is known that the Abel map
$\boldsymbol{\mathfrak{U}} \colon \mathscr{R}^{N}_{\mathrm{symm}} \! \to \!
\mathrm{Jac}(\mathscr{R})$ is surjective and locally biholomorphic, but not
injective globally. The \emph{dissected} Riemann surface, symbolically
$\widetilde{\mathscr{R}}$, is obtained {}from $\mathscr{R}$ by `cutting'
(canonical dissection) along the cycles of the canonical $1$-homology basis
$\boldsymbol{\alpha}_{k},\boldsymbol{\beta}_{k}$, $k \! = \! 1,\dotsc,N$, of
the original surface, namely, $\widetilde{\mathscr{R}} \! = \! \mathscr{R}
\setminus (\cup_{j=1}^{N}(\boldsymbol{\alpha}_{j} \cup \boldsymbol{\beta}_{
j}))$; the surface $\widetilde{\mathscr{R}}$ is not only connected, as one
can `pass' {}from one sheet to the other `across' $\Gamma_{N+1}$, but also
simply connected (a $4N$-sided polygon ($4N$-gon) of a canonical dissection
of $\mathscr{R}$ associated with the given canonical $1$-homology basis for
$\mathscr{R})$. For a given vector $\vec{\boldsymbol{v}} \! = \! (\upsilon_{
1},\upsilon_{2},\dotsc,\upsilon_{N}) \! \in \! \mathrm{Jac}(\mathscr{R})$, the
problem of finding an unordered collection of points $z_{1},z_{2},\dotsc,
z_{N}$, $z_{j} \! \in \! \mathscr{R}$, $j \! = \! 1,\dotsc,N$, for which
$\mathfrak{U}_{k}(z_{1},z_{2},\dotsc,z_{N}) \! \equiv \! \upsilon_{k}$, $k \!
= \! 1,\dotsc,N$, is called the \emph{Jacobi inversion problem} for Abelian
integrals: as is well known, the Jacobi inversion problem is always solvable;
but not, in general, uniquely.

By a \emph{divisor} on the Riemann surface $\mathscr{R}$ is meant a formal
`symbol' $\boldsymbol{d} \! = \! z_{1}^{n_{f}(z_{1})}z_{2}^{n_{f}(z_{2})}
\cdots z_{m}^{n_{f}(z_{m})}$, where $z_{j} \! \in \! \mathscr{R}$ and $n_{f}
(z_{j}) \! \in \! \mathbb{Z}$, $j \! = \! 1,\dotsc,m$: the number $\vert
\boldsymbol{d} \vert \! := \! \sum_{j=1}^{m}n_{f}(z_{j})$ is called the
\emph{degree} of the divisor $\boldsymbol{d}$: if $z_{i} \! \not= \! z_{j} \,
\, \forall \, \, i \! \not= \! j \! = \! 1,\dotsc,m$, and if $n_{f}(z_{j}) \!
\geqslant \! 0$, $j \! = \! 1,\dotsc,m$, then the divisor $\boldsymbol{d}$ is
said to be \emph{integral}. Let $g$ be a meromorphic function defined on
$\mathscr{R}$: for an arbitrary point $a \! \in \! \mathscr{R}$, one denotes
by $n_{g}(a)$ (resp., $p_{g}(a))$ the multiplicity of the zero (resp., pole)
of the function $g$ at this point if $a$ is a zero (resp., pole), and sets
$n_{g}(a) \! = \! 0$ (resp., $p_{g}(a) \! = \! 0)$ otherwise; thus, $n_{g}(a),
\, p_{g}(a) \! \geqslant \! 0$. To a meromorphic function $g$ on $\mathscr{
R}$, one assigns the divisor $(g)$ of zeros and poles of this function as $(g)
\! = \! z_{1}^{n_{g}(z_{1})}z_{2}^{n_{g}(z_{2})} \cdots z_{l_{1}}^{n_{g}(z_{
l_{1}})} \lambda_{1}^{-p_{g}(\lambda_{1})} \lambda_{2}^{-p_{g}(\lambda_{2})}
\cdots \lambda_{l_{2}}^{-p_{g}(\lambda_{l_{2}})}$, where $z_{i}, \, \lambda_{
j} \! \in \! \mathscr{R}$, $i \! = \! 1,\dotsc,l_{1}$, $j \! = \! 1,\dotsc,
l_{2}$, are the zeros and poles of $g$ on $\mathscr{R}$, and $n_{g}(z_{i}), \,
p_{g}(\lambda_{j}) \! \geqslant \! 0$ are their multiplicities (one can also
write $\lbrace \mathstrut (a,n_{g}(a),-p_{g}(a)); \, a \! \in \! \mathscr{R},
\, n_{g}(a), \, p_{g}(a) \! \geqslant \! 0 \rbrace$ for the divisor $(g)$ of
$g)$: these divisors are said to be \emph{principal}.

Associated with the Riemann matrix of $\boldsymbol{\beta}$-periods, $\tau$, 
is the \emph{Riemann theta function}, defined by
\begin{equation*}
\boldsymbol{\theta}(z;\tau) \! =: \! \boldsymbol{\theta}(z) \! = \! \sum_{m
\in \mathbb{Z}^{N}} \me^{2 \pi \mi (m,z)+\pi \mi (m,\tau m)}, \quad z \! \in
\! \mathbb{C}^{N},
\end{equation*}
where $(\boldsymbol{\cdot},\boldsymbol{\cdot})$ denotes the---real---Euclidean
inner/scalar product (for $\mathbf{A} \! = \! (A_{1},A_{2},\dotsc,A_{N}) \!
\in \! \mathbb{E}^{N}$ and $\mathbf{B} \! = \! (B_{1},B_{2},\dotsc,B_{N}) \!
\in \! \mathbb{E}^{N}$, $(A,B) \! := \! \sum_{k=1}^{N}A_{k}B_{k})$, with the
following evenness and (quasi-) periodicity properties,
\begin{equation*}
\boldsymbol{\theta}(-z) \! = \! \boldsymbol{\theta}(z), \qquad \boldsymbol{
\theta}(z \! + \! e_{j}) \! = \! \boldsymbol{\theta}(z), \qquad \mathrm{and}
\qquad \boldsymbol{\theta}(z \! \pm \! \tau_{j}) \! = \! \me^{\mp 2 \pi \mi
z_{j}-\mi \pi \tau_{jj}} \boldsymbol{\theta}(z),
\end{equation*}
where $e_{j}$ is the standard (basis) column vector in $\mathbb{C}^{N}$ with
$1$ in the $j$th entry and $0$ elsewhere (see above), and $\tau_{j} \! := \!
\tau e_{j}$ $(\in \! \mathbb{C}^{N})$, $j \! = \! 1,\dotsc,N$.

It turns out that, for the analysis of this work, the following multi-valued
functions are essential:
\begin{enumerate}
\item[$\boldsymbol{\bullet}$] $(R_{e}(z))^{1/2} \! := \! (\prod_{k=0}^{N}(z \!
- \! b_{k}^{e})(z \! - \! a_{k+1}^{e}))^{1/2}$, where, with the identification
$a_{N+1}^{e} \! \equiv \! a_{0}^{e}$ (as points on the complex sphere,
$\overline{\mathbb{C}})$ and with the point at infinity lying on the (open)
interval $(a_{0}^{e},b_{0}^{e})$, $-\infty \! < \! a_{0}^{e} \! < \! b_{0}^{e}
\! < \! a_{1}^{e} \! < \! b_{1}^{e} \! < \! \cdots \! < \! a_{N}^{e} \! < \!
b_{N}^{e} \! < \! +\infty$, $a_{0}^{e}$ $(\equiv \! a_{N+1}^{e})$ $\not= \!
-\infty,0$, and $b_{N}^{e} \! \not= \! 0,+\infty$ (see Figure~1);
\begin{figure}[tbh]
\begin{center}
\vspace{0.35cm}
\begin{pspicture}(0,0)(12,3)
\psset{xunit=1cm,yunit=1cm}
\psline[linewidth=0.9pt,linestyle=solid,linecolor=black]{o-o}(0.5,2)(2,2)
\psline[linewidth=0.9pt,linestyle=solid,linecolor=black]{o-o}(3,2)(4.5,2)
\psline[linewidth=0.9pt,linestyle=solid,linecolor=black]{o-o}(6.5,2)(8,2)
\psline[linewidth=0.9pt,linestyle=solid,linecolor=black]{o-o}(10,2)(11.5,2)
\psline[linewidth=0.7pt,linestyle=dotted,linecolor=darkgray](4.65,2)(6.35,2)
\psline[linewidth=0.7pt,linestyle=dotted,linecolor=darkgray](8.15,2)(9.85,2)
\rput(0.5,1.7){\makebox(0,0){$a_{0}^{e}$}}
\rput(0.5,0.9){\makebox(0,0){$a_{N+1}^{e}$}}
\rput{90}(0.5,1.3){\makebox(0,0){$\equiv$}}
\rput(1.25,2){\makebox(0,0){$\pmb{\times}$}}
\rput(1.25,2.3){\makebox(0,0){$\infty$}}
\rput(2,1.7){\makebox(0,0){$b_{0}^{e}$}}
\rput(3,1.7){\makebox(0,0){$a_{1}^{e}$}}
\rput(4.5,1.7){\makebox(0,0){$b_{1}^{e}$}}
\rput(6.5,1.7){\makebox(0,0){$a_{j}^{e}$}}
\rput(8,1.7){\makebox(0,0){$b_{j}^{e}$}}
\rput(10,1.7){\makebox(0,0){$a_{N}^{e}$}}
\rput(11.5,1.7){\makebox(0,0){$b_{N}^{e}$}}
\end{pspicture}
\end{center}
\vspace{-0.95cm}
\caption{Union of (open) intervals in the complex $z$-plane}
\end{figure}
\item[$\boldsymbol{\bullet}$] $(R_{o}(z))^{1/2} \! := \! (\prod_{k=0}^{N}(z \!
- \! b_{k}^{o})(z \! - \! a_{k+1}^{o}))^{1/2}$, where, with the identification
$a_{N+1}^{o} \! \equiv \! a_{0}^{o}$ (as points on the complex sphere,
$\overline{\mathbb{C}})$ and with the point at infinity lying on the (open)
interval $(a_{0}^{o},b_{0}^{o})$, $-\infty \! < \! a_{0}^{o} \! < \! b_{0}^{o}
\! < \! a_{1}^{o} \! < \! b_{1}^{o} \! < \! \cdots \! < \! a_{N}^{o} \! < \!
b_{N}^{o} \! < \! +\infty$, $a_{0}^{o}$ $(\equiv \! a_{N+1}^{o})$ $\not= \!
-\infty,0$, and $b_{N}^{o} \! \not= \! 0,+\infty$ (see Figure~2).
\begin{figure}[tbh]
\begin{center}
\vspace{0.35cm}
\begin{pspicture}(0,0)(12,3)
\psset{xunit=1cm,yunit=1cm}
\psline[linewidth=0.9pt,linestyle=solid,linecolor=black]{o-o}(0.5,2)(2,2)
\psline[linewidth=0.9pt,linestyle=solid,linecolor=black]{o-o}(3,2)(4.5,2)
\psline[linewidth=0.9pt,linestyle=solid,linecolor=black]{o-o}(6.5,2)(8,2)
\psline[linewidth=0.9pt,linestyle=solid,linecolor=black]{o-o}(10,2)(11.5,2)
\psline[linewidth=0.7pt,linestyle=dotted,linecolor=darkgray](4.65,2)(6.35,2)
\psline[linewidth=0.7pt,linestyle=dotted,linecolor=darkgray](8.15,2)(9.85,2)
\rput(0.5,1.7){\makebox(0,0){$a_{0}^{o}$}}
\rput(0.5,0.9){\makebox(0,0){$a_{N+1}^{o}$}}
\rput{90}(0.5,1.3){\makebox(0,0){$\equiv$}}
\rput(1.25,2){\makebox(0,0){$\pmb{\times}$}}
\rput(1.25,2.3){\makebox(0,0){$\infty$}}
\rput(2,1.7){\makebox(0,0){$b_{0}^{o}$}}
\rput(3,1.7){\makebox(0,0){$a_{1}^{o}$}}
\rput(4.5,1.7){\makebox(0,0){$b_{1}^{o}$}}
\rput(6.5,1.7){\makebox(0,0){$a_{j}^{o}$}}
\rput(8,1.7){\makebox(0,0){$b_{j}^{o}$}}
\rput(10,1.7){\makebox(0,0){$a_{N}^{o}$}}
\rput(11.5,1.7){\makebox(0,0){$b_{N}^{o}$}}
\end{pspicture}
\end{center}
\vspace{-0.95cm}
\caption{Union of (open) intervals in the complex $z$-plane}
\end{figure}
\end{enumerate}

The functions $R_{e}(z)$ and $R_{o}(z)$, respectively, are unital polynomials
$(\in \! \mathbb{R}[z])$ of even degree $(\deg (R_{e}(z)) \! = \! \deg (R_{o}
(z)) \! = \! 2(N \! + \! 1))$ whose (simple) roots/zeros are $\lbrace b_{j-
1}^{e},a_{j}^{e} \rbrace_{j=1}^{N+1}$ $(a_{N+1}^{e} \! \equiv \! a_{0}^{e})$
and $\lbrace b_{j-1}^{o},a_{j}^{o} \rbrace_{j=1}^{N+1}$ $(a_{N+1}^{o} \!
\equiv \! a_{0}^{o})$. The basic ingredients associated with the construction
of the hyperelliptic Riemann surfaces of genus $N$ corresponding,
respectively, to the multi-valued functions $y^{2} \! = \! R_{e}(z)$ and
$y^{2} \! = \! R_{o}(z)$ was given above. One now uses the above construction;
but particularised to the cases of the polynomials $R_{e}(z)$ and $R_{o}(z)$,
to arrive at the following:
\begin{enumerate}
\item[\shadowbox{$\sqrt{\smash[b]{R_{e}(z)}}$}] Let $\mathcal{Y}_{e}$ denote
the two-sheeted Riemann surface of genus $N$ associated with $y^{2} \! = \!
R_{e}(z)$, with $R_{e}(z)$ as characterised above: the first/upper (resp.,
second/lower) sheet of $\mathcal{Y}_{e}$ is denoted by $\mathcal{Y}_{e}^{+}$
(resp., $\mathcal{Y}_{e}^{-})$, points on the first/upper (resp.,
second/lower) sheet are represented as $z^{+} \! := \! (z,+(R_{e}(z))^{1/2})$
(resp., $z^{-} \! := \! (z,-(R_{e}(z))^{1/2}))$, where, as points on the plane
$\mathbb{C}$, $z^{+} \! = \! z^{-} \! = \! z$, and the single-valued branch
for the square root of the (multi-valued) function $(R_{e}(z))^{1/2}$ is
chosen such that $z^{-(N+1)}(R_{e}(z))^{1/2} \! \sim_{\underset{z \in
\mathcal{Y}_{e}^{\pm}}{z \to \infty}} \! \pm 1$. $\mathcal{Y}_{e}$ is realised
as a (two-sheeted) branched/ramified covering of the Riemann sphere such that
its two sheets are two identical copies of $\mathbb{C}$ with branch cuts
(slits) along the intervals $(a_{0}^{e},b_{0}^{e}),(a_{1}^{e},b_{1}^{e}),
\dotsc,(a_{N}^{e},b_{N}^{e})$ and pasted/glued together along $\cup_{j=1}^{N+1}
(a_{j-1}^{e},b_{j-1}^{e})$ $(a_{0}^{e} \! \equiv \! a_{N+1}^{e})$ in such a
way that the cycles $\boldsymbol{\alpha}_{0}^{e}$ and $\lbrace \boldsymbol{
\alpha}_{j}^{e},\boldsymbol{\beta}_{j}^{e} \rbrace$, $j \! = \! 1,\dotsc,
N$, where the latter forms the canonical $\mathbf{1}$-homology basis for
$\mathcal{Y}_{e}$, are characterised by the fact that (the closed contours)
$\boldsymbol{\alpha}_{j}^{e}$, $j \! = \! 0,\dotsc,N$, lie on $\mathcal{Y}_{
e}^{+}$, and (the closed contours) $\boldsymbol{\beta}_{j}^{e}$, $j \! = \!
1,\dotsc,N$, pass {}from $\mathcal{Y}_{e}^{+}$ (starting {}from the slit
$(a_{j}^{e},b_{j}^{e}))$, through the slit $(a_{0}^{e},b_{0}^{e})$ to
$\mathcal{Y}_{e}^{-}$, and back again to $\mathcal{Y}_{e}^{+}$ through the
slit $(a_{j}^{e},b_{j}^{e})$ (see Figure~3).
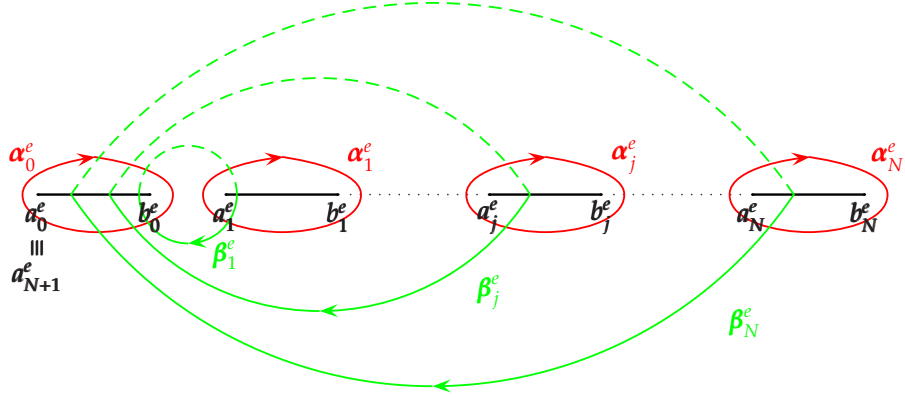
\begin{figure}[tbh]
\begin{center}
\vspace{0.45cm}
\begin{pspicture}(-2,0)(13,6)
\psset{xunit=1cm,yunit=1cm}
\psline[linewidth=0.9pt,linestyle=solid,linecolor=black]{c-c}(0.5,3)(2,3)
\psline[linewidth=0.9pt,linestyle=solid,linecolor=black]{c-c}(3,3)(4.5,3)
\psline[linewidth=0.9pt,linestyle=solid,linecolor=black]{c-c}(6.5,3)(8,3)
\psline[linewidth=0.9pt,linestyle=solid,linecolor=black]{c-c}(10,3)(11.5,3)
\psline[linewidth=0.7pt,linestyle=dotted,linecolor=darkgray](4.65,3)(6.35,3)
\psline[linewidth=0.7pt,linestyle=dotted,linecolor=darkgray](8.15,3)(9.85,3)
\pscurve[linewidth=0.7pt,linestyle=solid,linecolor=red,arrowsize=1.5pt 5]{->}%
(1.25,3.5)(2.3,3)(1.25,2.5)(0.3,3)(1.27,3.5)
\pscurve[linewidth=0.7pt,linestyle=solid,linecolor=red,arrowsize=1.5pt 5]{->}%
(3.75,3.5)(4.8,3)(3.75,2.5)(2.7,3)(3.77,3.5)
\pscurve[linewidth=0.7pt,linestyle=solid,linecolor=red,arrowsize=1.5pt 5]{->}%
(7.25,3.5)(8.3,3)(7.25,2.5)(6.2,3)(7.27,3.5)
\pscurve[linewidth=0.7pt,linestyle=solid,linecolor=red,arrowsize=1.5pt 5]{->}%
(10.75,3.5)(11.8,3)(10.75,2.5)(9.7,3)(10.77,3.5)
\rput(0.3,3.5){\makebox(0,0){{\color{red} $\boldsymbol{\alpha}_{0}^{e}$}}}
\rput(4.8,3.5){\makebox(0,0){{\color{red} $\boldsymbol{\alpha}_{1}^{e}$}}}
\rput(8.3,3.5){\makebox(0,0){{\color{red} $\boldsymbol{\alpha}_{j}^{e}$}}}
\rput(11.8,3.5){\makebox(0,0){{\color{red} $\boldsymbol{\alpha}_{N}^{e}$}}}
\psarc[linewidth=0.7pt,linestyle=solid,linecolor=green,arrowsize=1.5pt 5]{<-}%
(2.5,3){0.65}{270}{360}
\psarc[linewidth=0.7pt,linestyle=solid,linecolor=green](2.5,3){0.65}{180}{270}
\psarc[linewidth=0.7pt,linestyle=dashed,linecolor=green](2.5,3){0.65}{0}{180}
\psarc[linewidth=0.7pt,linestyle=solid,linecolor=green,arrowsize=1.5pt 5]%
{<-}(4.25,4.75){3.3}{270}{328}
\psarc[linewidth=0.7pt,linestyle=solid,linecolor=green](4.25,4.75){3.3}{212}%
{270}
\psarc[linewidth=0.7pt,linestyle=dashed,linecolor=green](4.25,1.25){3.3}{32}%
{148}
\psarc[linewidth=0.7pt,linestyle=solid,linecolor=green,arrowsize=1.5pt 5]%
{<-}(5.75,6.25){5.8}{270}{326}
\psarc[linewidth=0.7pt,linestyle=solid,linecolor=green](5.75,6.25){5.8}%
{214}{270}
\psarc[linewidth=0.7pt,linestyle=dashed,linecolor=green](5.75,-0.25){5.8}%
{34}{146}
\rput(3,2.2){\makebox(0,0){{\color{green} $\boldsymbol{\beta}_{1}^{e}$}}}
\rput(6.5,1.7){\makebox(0,0){{\color{green} $\boldsymbol{\beta}_{j}^{e}$}}}
\rput(9.9,1.3){\makebox(0,0){{\color{green} $\boldsymbol{\beta}_{N}^{e}$}}}
\rput(0.5,2.7){\makebox(0,0){$\pmb{a_{0}^{e}}$}}
\rput{90}(0.5,2.3){\makebox(0,0){$\pmb{\equiv}$}}
\rput(0.5,1.9){\makebox(0,0){$\pmb{a_{N+1}^{e}}$}}
\rput(2,2.7){\makebox(0,0){$\pmb{b_{0}^{e}}$}}
\rput(3,2.7){\makebox(0,0){$\pmb{a_{1}^{e}}$}}
\rput(4.5,2.7){\makebox(0,0){$\pmb{b_{1}^{e}}$}}
\rput(6.5,2.7){\makebox(0,0){$\pmb{a_{j}^{e}}$}}
\rput(8,2.7){\makebox(0,0){$\pmb{b_{j}^{e}}$}}
\rput(10,2.7){\makebox(0,0){$\pmb{a_{N}^{e}}$}}
\rput(11.5,2.7){\makebox(0,0){$\pmb{b_{N}^{e}}$}}
\end{pspicture}
\end{center}
\vspace{-0.55cm}
\caption{The Riemann surface $\mathcal{Y}_{e}$ of $y^{2} \! = \! \prod_{k=0
}^{N}(z \! - \! b_{k}^{e})(z \! - \! a_{k+1}^{e})$, $a_{N+1}^{e} \! \equiv
\! a_{0}^{e}$. The solid (resp., dashed) lines are on the first/upper (resp.,
second/lower) sheet of $\mathcal{Y}_{e}$, denoted $\mathcal{Y}_{e}^{+}$
(resp., $\mathcal{Y}_{e}^{-})$.}
\end{figure}

\hspace*{0.50cm}
The canonical $\mathbf{1}$-homology basis $\lbrace \boldsymbol{\alpha}_{j}^{e},
\boldsymbol{\beta}_{j}^{e} \rbrace$, $j \! = \! 1,\dotsc,N$, generates, on
$\mathcal{Y}_{e}$, the (corresponding) $\boldsymbol{\alpha}^{e}$-normalised
basis of holomorphic Abelian differentials (one-forms) $\lbrace \omega_{1}^{
e},\omega_{2}^{e},\dotsc,\omega_{N}^{e} \rbrace$, where $\omega_{j}^{e} \! :=
\! \sum_{k=1}^{N} \tfrac{c_{jk}^{e}z^{N-k}}{\sqrt{\smash[b]{R_{e}(z)}}} \, \md
z$, $c_{jk}^{e} \! \in \! \mathbb{C}$, $j \! = \! 1,\dotsc,N$, and $\oint_{
\boldsymbol{\alpha}_{k}^{e}} \omega_{j}^{e} \! = \! \delta_{kj}$, $k,j \! =
\! 1,\dotsc,N$: $\omega_{l}^{e}$, $l \! = \! 1,\dotsc,N$, is real valued on
$\cup_{j=1}^{N+1}(a_{j-1}^{e},b_{j-1}^{e})$, and has exactly one (real) root
in any (open) interval $(a_{j-1}^{e},b_{j-1}^{e})$, $j \! = \! 1,\dotsc,N \!
+ \! 1$; furthermore, in the intervals $(b_{j-1}^{e},a_{j}^{e})$, $j \! = \!
1,\dotsc,N$, $\omega_{l}^{e}$, $l \! = \! 1,\dotsc,N$, take non-zero, pure
imaginary values. Let $\boldsymbol{\omega}^{e} \! := \! (\omega_{1}^{e},
\omega_{2}^{e},\dotsc,\omega_{N}^{e})$ denote the basis of holomorphic
one-forms on $\mathcal{Y}_{e}$ as normalised above with the associated $N \!
\times \! N$ Riemann matrix of $\boldsymbol{\beta}^{e}$-periods, $\tau^{e} \!
= \! (\tau^{e})_{i,j=1,\dotsc,N} \! := \! (\oint_{\boldsymbol{\beta}_{j}^{e}}
\omega_{i}^{e})_{i,j=1,\dotsc,N}$: the Riemann matrix, $\tau^{e}$, is
symmetric $(\tau_{ij}^{e} \! = \! \tau_{ji}^{e})$ and pure imaginary, $-\mi
\tau^{e}$ is positive definite $(\Im (\tau^{e}_{ij}) \! > \! 0)$, and $\det
(\tau^{e}) \! \not= \! 0$ (non-degenerate). For the holomorphic Abelian
differential (one-form) $\boldsymbol{\omega}^{e}$ defined above, choose
$a_{N+1}^{e}$ as the \emph{base point}, and set $\boldsymbol{u}^{e} \colon
\mathcal{Y}_{e} \! \to \! \operatorname{Jac}(\mathcal{Y}_{e})$ $(:= \!
\mathbb{C}^{N}/\lbrace N \! + \! \tau^{e}M \rbrace$, $(N,M) \! \in \!
\mathbb{Z}^{N} \! \times \! \mathbb{Z}^{N})$, $z \! \mapsto \! \boldsymbol{
u}^{e}(z) \! := \! \int_{a_{N+1}^{e}}^{z} \boldsymbol{\omega}^{e}$, where the
integration {}from $a_{N+1}^{e}$ to $z$ $(\in \mathcal{Y}_{e})$ is taken along
any path on $\mathcal{Y}_{e}^{+}$.
\begin{eeee}
{}From the representation $\omega_{j}^{e} \! = \! \sum_{k=1}^{N} \tfrac{c_{j
k}^{e}z^{N-k}}{\sqrt{\smash[b]{R_{e}(z)}}} \, \md z$, $j \! = \! 1,\dotsc,N$,
and the normalisation condition $\oint_{\boldsymbol{\alpha}_{k}^{e}} \omega_{
j}^{e} \! = \! \delta_{kj}$, $k,j \! = \! 1,\dotsc,N$, one shows that $c_{j
k}^{e}$, $k,j \! = \! 1,\dotsc,N$, are obtained from
\begin{equation}
\begin{pmatrix}
c_{11}^{e} & c_{12}^{e} & \dotsb & c_{1N}^{e} \\
c_{21}^{e} & c_{22}^{e} & \dotsb & c_{2N}^{e} \\
\vdots     & \vdots     & \ddots & \vdots     \\
c_{N1}^{e} & c_{N2}^{e} & \dotsb & c_{NN}^{e}
\end{pmatrix} \! = \! \widetilde{\mathfrak{S}}_{e}^{-1} \tag{E1},
\end{equation}
where
\begin{equation}
\widetilde{\mathfrak{S}}_{e} \! := \!
\begin{pmatrix}
\oint_{\boldsymbol{\alpha}_{1}^{e}} \frac{\md s_{1}}{\sqrt{\smash[b]{R_{e}
(s_{1})}}} & \oint_{\boldsymbol{\alpha}_{2}^{e}} \frac{\md s_{2}}{\sqrt{
\smash[b]{R_{e}(s_{2})}}} & \dotsb & \oint_{\boldsymbol{\alpha}_{N}^{e}}
\frac{\md s_{N}}{\sqrt{\smash[b]{R_{e}(s_{N})}}} \\
\oint_{\boldsymbol{\alpha}_{1}^{e}} \frac{s_{1} \md s_{1}}{\sqrt{\smash[b]{
R_{e}(s_{1})}}} & \oint_{\boldsymbol{\alpha}_{2}^{e}} \frac{s_{2} \md s_{2}}{
\sqrt{\smash[b]{R_{e}(s_{2})}}} & \dotsb & \oint_{\boldsymbol{\alpha}_{N}^{e}}
\frac{s_{N} \md s_{N}}{\sqrt{\smash[b]{R_{e}(s_{N})}}} \\
\vdots & \vdots & \ddots & \vdots \\
\oint_{\boldsymbol{\alpha}_{1}^{e}} \frac{s_{1}^{N-1} \md s_{1}}{\sqrt{
\smash[b]{R_{e}(s_{1})}}} & \oint_{\boldsymbol{\alpha}_{2}^{e}} \frac{s_{
2}^{N-1} \md s_{2}}{\sqrt{\smash[b]{R_{e}(s_{2})}}} & \dotsb & \oint_{
\boldsymbol{\alpha}_{N}^{e}} \frac{s_{N}^{N-1} \md s_{N}}{\sqrt{\smash[b]{
R_{e}(s_{N})}}}
\end{pmatrix} \tag{E2}.
\end{equation}
For a (representation-independent) proof of the fact that $\det (\widetilde{
\mathfrak{S}}_{e}) \! \not= \! 0$, see, for example, Chapter~10, 
Section~10--2, of \cite{a77}. \hfill $\blacksquare$
\end{eeee}

\hspace*{0.50cm}
Set (see \cite{a38}), for $z \! \in \! \mathbb{C}_{+}$, $\gamma^{e}(z) \! 
:= \! (\prod_{k=1}^{N+1}(z \! - \! b_{k-1}^{e})(z \! - \! a_{k}^{e})^{-1})^{
1/4}$, and, for $z \! \in \! \mathbb{C}_{-}$, $\gamma^{e}(z) \! := \! -\mi 
(\prod_{k=1}^{N+1}(z \! - \! b_{k-1}^{e})(z \! - \! a_{k}^{e})^{-1})^{1/4}$. 
It is shown in \cite{a38} that $\gamma^{e}(z) \! =_{\underset{z \in \mathcal{
Y}_{e}^{\pm}}{z \to \infty}} \! (-\mi)^{(1 \mp 1)/2}(1 \! + \! \mathcal{O}
(z^{-1}))$, and
\begin{equation*}
\left\{z_{j}^{e,\pm} \right\}_{j=1}^{N} \! = \! \left\{\mathstrut z^{\pm} \!
\in \! \mathcal{Y}_{e}^{\pm}; \, (\gamma^{e}(z) \! \mp \! (\gamma^{e}(z))^{-
1}) \vert_{z=z^{\pm}} \! = \! 0 \right\},
\end{equation*}
with $z_{j}^{e,\pm} \! \in \! (a_{j}^{e},b_{j}^{e})^{\pm}$ $(\subset \!
\mathcal{Y}_{e}^{\pm})$, $j \! = \! 1,\dotsc,N$, where, as points on the
plane, $z_{j}^{e,+} \! = \! z_{j}^{e,-} \! := \! z_{j}^{e}$, $j \! = \! 1,
\dotsc,N$ (of course, on the plane, $z_{j}^{e} \! \in \! (a_{j}^{e},b_{j}^{
e})$, $j \! = \! 1,\dotsc,N)$.

\hspace*{0.50cm}
Corresponding to $\mathcal{Y}_{e}$, define $\boldsymbol{d}_{e} \! := \!
-\boldsymbol{K}_{e} \! - \! \sum_{j=1}^{N} \int_{a_{N+1}^{e}}^{z_{j}^{e,-}}
\boldsymbol{\omega}^{e}$ $(\in \! \mathbb{C}^{N})$, where $\boldsymbol{K}_{e}$
is the associated (`even') vector of Riemann constants, and the integration
{}from $a_{N+1}^{e}$ to $z_{j}^{e,-}$, $j \! = \! 1,\dotsc,N$, is taken along
a fixed path in $\mathcal{Y}_{e}^{-}$. It is shown in Chapter~VII of
\cite{a78} that $\boldsymbol{K}_{e} \! = \! \sum_{j=1}^{N} \int_{a_{j}^{e}
}^{a_{N+1}^{e}} \boldsymbol{\omega}^{e}$; furthermore, $\boldsymbol{K}_{e}$
is a point of order $2$, that is, $2 \boldsymbol{K}_{e} \! = \! 0$ and $s
\boldsymbol{K}_{e} \! \not= \! 0$ for $0 \! < \! s \! < \! 2$. Recalling the
definition of $\boldsymbol{\omega}^{e}$ and that $z^{-(N+1)}(R_{e}(z))^{1/2}
\! \sim_{\underset{z \in \mathbb{C}_{\pm}}{z \to \infty}} \! \pm 1$, using the
fact that $\boldsymbol{K}_{e}$ is a point of order $2$, one arrives at
\begin{align*}
\boldsymbol{d}_{e} =& \, -\boldsymbol{K}_{e} \! - \! \sum_{j=1}^{N} \int_{a_{N
+1}^{e}}^{z_{j}^{e,-}} \boldsymbol{\omega}^{e} \! = \! \boldsymbol{K}_{e} \! -
\! \sum_{j=1}^{N} \int_{a_{N+1}^{e}}^{z_{j}^{e,-}} \boldsymbol{\omega}^{e} \!
= \! -\boldsymbol{K}_{e} \! + \! \sum_{j=1}^{N} \int_{a_{N+1}^{e}}^{z_{j}^{e,
+}} \boldsymbol{\omega}^{e} \! = \! \boldsymbol{K}_{e} \! + \! \sum_{j=1}^{N}
\int_{a_{N+1}^{e}}^{z_{j}^{e,+}} \boldsymbol{\omega}^{e} \\
=& \, -\sum_{j=1}^{N} \int_{a_{j}^{e}}^{z_{j}^{e,-}} \boldsymbol{\omega}^{e}
\! = \! \sum_{j=1}^{N} \int_{a_{j}^{e}}^{z_{j}^{e,+}} \boldsymbol{\omega}^{e}.
\end{align*}

\hspace*{0.50cm}
Associated with the Riemann matrix of $\boldsymbol{\beta}^{e}$-periods,
$\tau^{e}$, is the (`even') Riemann theta function
\begin{equation*}
\boldsymbol{\theta}(z;\tau^{e}) \! =: \! \boldsymbol{\theta}^{e}(z) \! = \!
\sum_{m \in \mathbb{Z}^{N}} \me^{2 \pi \mi (m,z)+\pi \mi (m,\tau^{e}m)}, \quad
z \! \in \! \mathbb{C}^{N};
\end{equation*}
$\boldsymbol{\theta}^{e}(z)$ has the following evenness and (quasi-)
periodicity properties,
\begin{equation*}
\boldsymbol{\theta}^{e}(-z) \! = \! \boldsymbol{\theta}^{e}(z), \qquad
\boldsymbol{\theta}^{e}(z \! + \! e_{j}) \! = \! \boldsymbol{\theta}^{e}(z),
\qquad \mathrm{and} \qquad \boldsymbol{\theta}^{e}(z \! \pm \! \tau_{j}^{e})
\! = \! \me^{\mp 2 \pi \mi z_{j}-\mi \pi \tau_{jj}^{e}} \boldsymbol{\theta}^{e}
(z),
\end{equation*}
where $\tau_{j}^{e} \! := \! \tau^{e} e_{j}$ $(\in \! \mathbb{C}^{N})$, $j
\! = \! 1,\dotsc,N$. This entire latter apparatus is used extensively in
\cite{a38}.

\item[\shadowbox{$\sqrt{\smash[b]{R_{o}(z)}}$}] Let $\mathcal{Y}_{o}$ denote
the two-sheeted Riemann surface of genus $N$ associated with $y^{2} \! = \!
R_{o}(z)$, with $R_{o}(z)$ as characterised above: the first/upper (resp.,
second/lower) sheet of $\mathcal{Y}_{o}$ is denoted by $\mathcal{Y}_{o}^{+}$
(resp., $\mathcal{Y}_{o}^{-})$, points on the first/upper (resp.,
second/lower) sheet are represented as $z^{+} \! := \! (z,+(R_{o}(z))^{1/2})$
(resp., $z^{-} \! := \! (z,-(R_{o}(z))^{1/2}))$, where, as points on the plane
$\mathbb{C}$, $z^{+} \! = \! z^{-} \! = \! z$, and the single-valued branch
for the square root of the (multi-valued) function $(R_{o}(z))^{1/2}$ is
chosen such that $z^{-(N+1)}(R_{o}(z))^{1/2} \! \sim_{\underset{z \in
\mathcal{Y}_{o}^{\pm}}{z \to \infty}} \! \pm 1$. $\mathcal{Y}_{o}$ is realised
as a (two-sheeted) branched/ramified covering of the Riemann sphere such that
its two sheets are two identical copies of $\mathbb{C}$ with branch cuts
(slits) along the intervals $(a_{0}^{o},b_{0}^{o}),(a_{1}^{o},b_{1}^{o}),
\dotsc,(a_{N}^{o},b_{N}^{o})$ and pasted/glued together along $\cup_{j=1}^{N+1}
(a_{j-1}^{o},b_{j-1}^{o})$ $(a_{0}^{o} \! \equiv \! a_{N+1}^{o})$ in such a
way that the cycles $\boldsymbol{\alpha}_{0}^{o}$ and $\lbrace \boldsymbol{
\alpha}_{j}^{o},\boldsymbol{\beta}_{j}^{o} \rbrace$, $j \! = \! 1,\dotsc,
N$, where the latter forms the canonical $\mathbf{1}$-homology basis for
$\mathcal{Y}_{o}$, are characterised by the fact that (the closed contours)
$\boldsymbol{\alpha}_{j}^{o}$, $j \! = \! 0,\dotsc,N$, lie on $\mathcal{Y}_{
o}^{+}$, and (the closed contours) $\boldsymbol{\beta}_{j}^{o}$, $j \! = \!
1,\dotsc,N$, pass {}from $\mathcal{Y}_{o}^{+}$ (starting {}from the slit
$(a_{j}^{o},b_{j}^{o}))$, through the slit $(a_{0}^{o},b_{0}^{o})$ to
$\mathcal{Y}_{o}^{-}$, and back again to $\mathcal{Y}_{o}^{+}$ through the
slit $(a_{j}^{o},b_{j}^{o})$ (see Figure~4).
\begin{figure}[htb]
\begin{center}
\vspace{0.45cm}
\begin{pspicture}(-2,0)(13,6)
\psset{xunit=1cm,yunit=1cm}
\psline[linewidth=0.9pt,linestyle=solid,linecolor=black]{c-c}(0.5,3)(2,3)
\psline[linewidth=0.9pt,linestyle=solid,linecolor=black]{c-c}(3,3)(4.5,3)
\psline[linewidth=0.9pt,linestyle=solid,linecolor=black]{c-c}(6.5,3)(8,3)
\psline[linewidth=0.9pt,linestyle=solid,linecolor=black]{c-c}(10,3)(11.5,3)
\psline[linewidth=0.7pt,linestyle=dotted,linecolor=darkgray](4.65,3)(6.35,3)
\psline[linewidth=0.7pt,linestyle=dotted,linecolor=darkgray](8.15,3)(9.85,3)
\pscurve[linewidth=0.7pt,linestyle=solid,linecolor=red,arrowsize=1.5pt 5]{->}%
(1.25,3.5)(2.3,3)(1.25,2.5)(0.3,3)(1.27,3.5)
\pscurve[linewidth=0.7pt,linestyle=solid,linecolor=red,arrowsize=1.5pt 5]{->}%
(3.75,3.5)(4.8,3)(3.75,2.5)(2.7,3)(3.77,3.5)
\pscurve[linewidth=0.7pt,linestyle=solid,linecolor=red,arrowsize=1.5pt 5]{->}%
(7.25,3.5)(8.3,3)(7.25,2.5)(6.2,3)(7.27,3.5)
\pscurve[linewidth=0.7pt,linestyle=solid,linecolor=red,arrowsize=1.5pt 5]{->}%
(10.75,3.5)(11.8,3)(10.75,2.5)(9.7,3)(10.77,3.5)
\rput(0.3,3.5){\makebox(0,0){{\color{red} $\boldsymbol{\alpha}_{0}^{o}$}}}
\rput(4.8,3.5){\makebox(0,0){{\color{red} $\boldsymbol{\alpha}_{1}^{o}$}}}
\rput(8.3,3.5){\makebox(0,0){{\color{red} $\boldsymbol{\alpha}_{j}^{o}$}}}
\rput(11.8,3.5){\makebox(0,0){{\color{red} $\boldsymbol{\alpha}_{N}^{o}$}}}
\psarc[linewidth=0.7pt,linestyle=solid,linecolor=green,arrowsize=1.5pt 5]{<-}%
(2.5,3){0.65}{270}{360}
\psarc[linewidth=0.7pt,linestyle=solid,linecolor=green](2.5,3){0.65}{180}{270}
\psarc[linewidth=0.7pt,linestyle=dashed,linecolor=green](2.5,3){0.65}{0}{180}
\psarc[linewidth=0.7pt,linestyle=solid,linecolor=green,arrowsize=1.5pt 5]{<-}%
(4.25,4.75){3.3}{270}{328}
\psarc[linewidth=0.7pt,linestyle=solid,linecolor=green](4.25,4.75){3.3}{212}%
{270}
\psarc[linewidth=0.7pt,linestyle=dashed,linecolor=green](4.25,1.25){3.3}{32}%
{148}
\psarc[linewidth=0.7pt,linestyle=solid,linecolor=green,arrowsize=1.5pt 5]{<-}%
(5.75,6.25){5.8}{270}{326}
\psarc[linewidth=0.7pt,linestyle=solid,linecolor=green](5.75,6.25){5.8}%
{214}{270}
\psarc[linewidth=0.7pt,linestyle=dashed,linecolor=green](5.75,-0.25){5.8}%
{34}{146}
\rput(3,2.2){\makebox(0,0){{\color{green} $\boldsymbol{\beta}_{1}^{o}$}}}
\rput(6.5,1.7){\makebox(0,0){{\color{green} $\boldsymbol{\beta}_{j}^{o}$}}}
\rput(9.9,1.3){\makebox(0,0){{\color{green} $\boldsymbol{\beta}_{N}^{o}$}}}
\rput(0.5,2.7){\makebox(0,0){$\pmb{a_{0}^{o}}$}}
\rput{90}(0.5,2.3){\makebox(0,0){$\pmb{\equiv}$}}
\rput(0.5,1.9){\makebox(0,0){$\pmb{a_{N+1}^{o}}$}}
\rput(2,2.7){\makebox(0,0){$\pmb{b_{0}^{o}}$}}
\rput(3,2.7){\makebox(0,0){$\pmb{a_{1}^{o}}$}}
\rput(4.5,2.7){\makebox(0,0){$\pmb{b_{1}^{o}}$}}
\rput(6.5,2.7){\makebox(0,0){$\pmb{a_{j}^{o}}$}}
\rput(8,2.7){\makebox(0,0){$\pmb{b_{j}^{o}}$}}
\rput(10,2.7){\makebox(0,0){$\pmb{a_{N}^{o}}$}}
\rput(11.5,2.7){\makebox(0,0){$\pmb{b_{N}^{o}}$}}
\end{pspicture}
\end{center}
\vspace{-0.55cm}
\caption{The Riemann surface $\mathcal{Y}_{o}$ of $y^{2} \! = \! \prod_{k=0
}^{N+1}(z \! - \! b_{k}^{o})(z \! - \! a_{k+1}^{o})$, $a_{N+1}^{o} \! \equiv
\! a_{0}^{o}$. The solid (resp., dashed) lines are on the first/upper (resp.,
second/lower) sheet of $\mathcal{Y}_{o}$, denoted $\mathcal{Y}_{o}^{+}$
(resp., $\mathcal{Y}_{o}^{-})$.}
\end{figure}
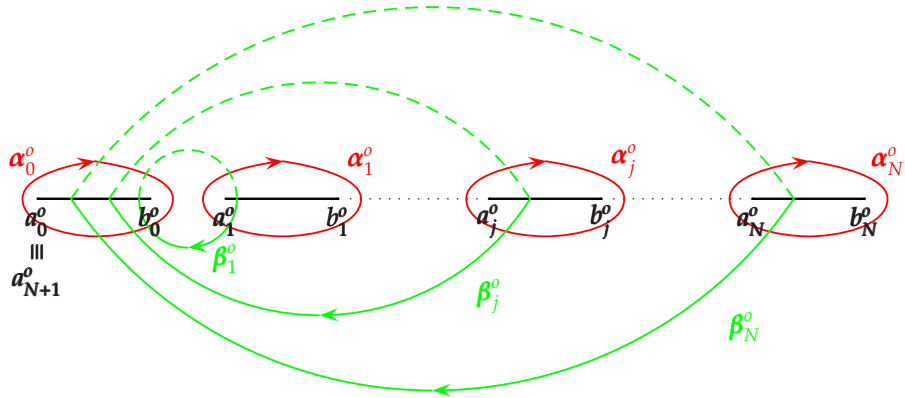

\hspace*{0.50cm}
The canonical $\mathbf{1}$-homology basis $\lbrace \boldsymbol{\alpha}_{j}^{o},
\boldsymbol{\beta}_{j}^{o} \rbrace$, $j \! = \! 1,\dotsc,N$, generates, on
$\mathcal{Y}_{o}$, the (corresponding) $\boldsymbol{\alpha}^{o}$-normalised
basis of holomorphic Abelian differentials (one-forms) $\lbrace \omega_{1}^{
o},\omega_{2}^{o},\dotsc,\omega_{N}^{o} \rbrace$, where $\omega_{j}^{o} \! :=
\! \sum_{k=1}^{N} \tfrac{c_{jk}^{o}z^{N-k}}{\sqrt{\smash[b]{R_{o}(z)}}} \, \md
z$, $c_{jk}^{o} \! \in \! \mathbb{C}$, $j \! = \! 1,\dotsc,N$, and $\oint_{
\boldsymbol{\alpha}_{k}^{o}} \omega_{j}^{o} \! = \! \delta_{kj}$, $k,j \! =
\! 1,\dotsc,N$: $\omega_{l}^{o}$, $l \! = \! 1,\dotsc,N$, is real valued on
$\cup_{j=1}^{N+1}(a_{j-1}^{o},b_{j-1}^{o})$, and has exactly one (real) root
in any (open) interval $(a_{j-1}^{o},b_{j-1}^{o})$, $j \! = \! 1,\dotsc,N \!
+ \! 1$; furthermore, in the intervals $(b_{j-1}^{o},a_{j}^{o})$, $j \! = \!
1,\dotsc,N$, $\omega_{l}^{o}$, $l \! = \! 1,\dotsc,N$, take non-zero, pure
imaginary values. Let $\boldsymbol{\omega}^{o} \! := \! (\omega_{1}^{o},
\omega_{2}^{o},\dotsc,\omega_{N}^{o})$ denote the basis of holomorphic
one-forms on $\mathcal{Y}_{o}$ as normalised above with the associated $N \!
\times \! N$ Riemann matrix of $\boldsymbol{\beta}^{o}$-periods, $\tau^{o} \!
= \! (\tau^{o})_{i,j=1,\dotsc,N} \! := \! (\oint_{\boldsymbol{\beta}_{j}^{o}}
\omega_{i}^{o})_{i,j=1,\dotsc,N}$: the Riemann matrix, $\tau^{o}$, is
symmetric $(\tau_{ij}^{o} \! = \! \tau_{ji}^{o})$ and pure imaginary, $-\mi
\tau^{o}$ is positive definite $(\Im (\tau^{o}_{ij}) \! > \! 0)$, and $\det
(\tau^{o}) \! \not= \! 0$ (non-degenerate). For the holomorphic Abelian
differential (one-form) $\boldsymbol{\omega}^{o}$ defined above, choose
$a_{N+1}^{o}$ as the base point, and set $\boldsymbol{u}^{o} \colon \mathcal{
Y}_{o} \! \to \! \operatorname{Jac}(\mathcal{Y}_{o})$ $(:= \! \mathbb{C}^{N}/
\lbrace N \! + \! \tau^{o}M \rbrace$, $(N,M) \! \in \! \mathbb{Z}^{N} \!
\times \! \mathbb{Z}^{N})$, $z \! \mapsto \! \boldsymbol{u}^{o}(z) \! := \!
\int_{a_{N+1}^{o}}^{z} \boldsymbol{\omega}^{o}$, where the integration {}from
$a_{N+1}^{o}$ to $z$ $(\in \mathcal{Y}_{o})$ is taken along any path on
$\mathcal{Y}_{o}^{+}$.
\begin{eeee}
{}From the representation $\omega_{j}^{o} \! = \! \sum_{k=1}^{N} \tfrac{c_{j
k}^{o}z^{N-k}}{\sqrt{\smash[b]{R_{o}(z)}}} \, \md z$, $j \! = \! 1,\dotsc,N$,
and the normalisation condition $\oint_{\boldsymbol{\alpha}_{k}^{o}} \omega_{
j}^{o} \! = \! \delta_{kj}$, $k,j \! = \! 1,\dotsc,N$, one shows that $c_{j
k}^{o}$, $k,j \! = \! 1,\dotsc,N$, are obtained from
\begin{equation}
\begin{pmatrix}
c_{11}^{o} & c_{12}^{o} & \dotsb & c_{1N}^{o} \\
c_{21}^{o} & c_{22}^{o} & \dotsb & c_{2N}^{o} \\
\vdots     & \vdots     & \ddots & \vdots     \\
c_{N1}^{o} & c_{N2}^{o} & \dotsb & c_{NN}^{o}
\end{pmatrix} \! = \! \widetilde{\mathfrak{S}}_{o}^{-1} \tag{O1},
\end{equation}
where
\begin{equation}
\widetilde{\mathfrak{S}}_{o} \! := \!
\begin{pmatrix}
\oint_{\boldsymbol{\alpha}_{1}^{o}} \frac{\md s_{1}}{\sqrt{\smash[b]{R_{o}
(s_{1})}}} & \oint_{\boldsymbol{\alpha}_{2}^{o}} \frac{\md s_{2}}{\sqrt{
\smash[b]{R_{o}(s_{2})}}} & \dotsb & \oint_{\boldsymbol{\alpha}_{N}^{o}}
\frac{\md s_{N}}{\sqrt{\smash[b]{R_{o}(s_{N})}}} \\
\oint_{\boldsymbol{\alpha}_{1}^{o}} \frac{s_{1} \md s_{1}}{\sqrt{\smash[b]{
R_{o}(s_{1})}}} & \oint_{\boldsymbol{\alpha}_{2}^{o}} \frac{s_{2} \md s_{2}}{
\sqrt{\smash[b]{R_{o}(s_{2})}}} & \dotsb & \oint_{\boldsymbol{\alpha}_{N}^{o}}
\frac{s_{N} \md s_{N}}{\sqrt{\smash[b]{R_{o}(s_{N})}}} \\
\vdots & \vdots & \ddots & \vdots \\
\oint_{\boldsymbol{\alpha}_{1}^{o}} \frac{s_{1}^{N-1} \md s_{1}}{\sqrt{
\smash[b]{R_{o}(s_{1})}}} & \oint_{\boldsymbol{\alpha}_{2}^{o}} \frac{s_{
2}^{N-1} \md s_{2}}{\sqrt{\smash[b]{R_{o}(s_{2})}}} & \dotsb & \oint_{
\boldsymbol{\alpha}_{N}^{o}} \frac{s_{N}^{N-1} \md s_{N}}{\sqrt{\smash[b]{
R_{o}(s_{N})}}}
\end{pmatrix} \tag{O2}.
\end{equation}
For a (representation-independent) proof of the fact that $\det (\widetilde{
\mathfrak{S}}_{o}) \! \not= \! 0$, see, for example, Chapter~10, 
Section~10--2, of \cite{a77}. \hfill $\blacksquare$
\end{eeee}

\hspace*{0.50cm}
Set (see Section~4), for $z \! \in \! \mathbb{C}_{+}$, $\gamma^{o}(z) \! := \! 
(\prod_{k=1}^{N+1}(z \! - \! b_{k-1}^{o})(z \! - \! a_{k}^{o})^{-1})^{1/4}$, 
and, for $z \! \in \! \mathbb{C}_{-}$, $\gamma^{o}(z) \! := \! -\mi (\prod_{
k=1}^{N+1}(z \! - \! b_{k-1}^{o})(z \! - \! a_{k}^{o})^{-1})^{1/4}$. It is 
shown in Section~4 that $\gamma^{o}(z) \! =_{\underset{z \in \mathcal{Y}_{
o}^{\pm}}{z \to 0}} \! (-\mi)^{(1 \mp 1)/2} \gamma^{o}(0) \linebreak[4]
\cdot (1 \! + \! \mathcal{O}(z))$, where $\gamma^{o}(0) \! := \! (\prod_{k=
1}^{N+1}b_{k-1}^{o}(a_{k}^{o})^{-1})^{1/4}$ $(> \! 0)$, and a set of $N$ 
upper-edge and lower-edge finite-length-gap roots/zeros are
\begin{equation*}
\left\{z_{j}^{o,\pm} \right\}_{j=1}^{N} \! = \! \left\{\mathstrut z^{\pm} \!
\in \! \mathcal{Y}_{o}^{\pm}; \, ((\gamma^{o}(0))^{-1} \gamma^{o}(z) \! \mp
\! \gamma^{o}(0)(\gamma^{o}(z))^{-1}) \vert_{z=z^{\pm}} \! = \! 0 \right\},
\end{equation*}
with $z_{j}^{o,\pm} \! \in \! (a_{j}^{o},b_{j}^{o})^{\pm}$ $(\subset \!
\mathcal{Y}_{o}^{\pm})$, $j \! = \! 1,\dotsc,N$, where, as points on the
plane, $z_{j}^{o,+} \! = \! z_{j}^{o,-} \! := \! z_{j}^{o}$, $j \! = \! 1,
\dotsc,N$ (of course, on the plane, $z_{j}^{o} \! \in \! (a_{j}^{o},b_{j}^{o}
)$, $j \! = \! 1,\dotsc,N)$.

\hspace*{0.50cm}
Corresponding to $\mathcal{Y}_{o}$, define $\boldsymbol{d}_{o} \! := \! -
\boldsymbol{K}_{o} \! - \! \sum_{j=1}^{N} \int_{a_{N+1}^{o}}^{z_{j}^{o,-}}
\boldsymbol{\omega}^{o}$ $(\in \! \mathbb{C}^{N})$, where $\boldsymbol{K}_{o}$
is the associated (`odd') vector of Riemann constants, and the integration
{}from $a_{N+1}^{o}$ to $z_{j}^{o,-}$, $j \! = \! 1,\dotsc,N$, is taken along
a fixed path in $\mathcal{Y}_{o}^{-}$. It is shown in Chapter~VII of
\cite{a78} that $\boldsymbol{K}_{o} \! = \! \sum_{j=1}^{N} \int_{a_{j}^{o}}^{
a_{N+1}^{o}} \boldsymbol{\omega}^{o}$; furthermore, $\boldsymbol{K}_{o}$ is a
point of order $2$. Recalling the definition of $\boldsymbol{\omega}^{o}$ and
that $z^{-(N+1)}(R_{o}(z))^{1/2} \! \sim_{\underset{z \in \mathbb{C}_{\pm}}{z
\to \infty}} \! \pm 1$, using the fact that $\boldsymbol{K}_{o}$ is a point of
order $2$, one arrives at
\begin{align*}
\boldsymbol{d}_{o} =& \, -\boldsymbol{K}_{o} \! - \! \sum_{j=1}^{N} \int_{a_{N
+1}^{o}}^{z_{j}^{o,-}} \boldsymbol{\omega}^{o} \! = \! \boldsymbol{K}_{o} \! -
\! \sum_{j=1}^{N} \int_{a_{N+1}^{o}}^{z_{j}^{o,-}} \boldsymbol{\omega}^{o} \!
= \! -\boldsymbol{K}_{o} \! + \! \sum_{j=1}^{N} \int_{a_{N+1}^{o}}^{z_{j}^{o,
+}} \boldsymbol{\omega}^{o} \! = \! \boldsymbol{K}_{o} \! + \! \sum_{j=1}^{N}
\int_{a_{N+1}^{o}}^{z_{j}^{o,+}} \boldsymbol{\omega}^{o} \\
=& \, -\sum_{j=1}^{N} \int_{a_{j}^{o}}^{z_{j}^{o,-}} \boldsymbol{\omega}^{o}
\! = \! \sum_{j=1}^{N} \int_{a_{j}^{o}}^{z_{j}^{o,+}} \boldsymbol{\omega}^{o}.
\end{align*}

\hspace*{0.50cm}
Associated with the Riemann matrix of $\boldsymbol{\beta}^{o}$-periods,
$\tau^{o}$, is the (`odd') Riemann theta function
\begin{equation}
\boldsymbol{\theta}(z;\tau^{o}) \! =: \! \boldsymbol{\theta}^{o}(z) \! = \!
\sum_{m \in \mathbb{Z}^{N}} \me^{2 \pi \mi (m,z)+\pi \mi (m,\tau^{o}m)}, \quad
z \! \in \! \mathbb{C}^{N};
\end{equation}
$\boldsymbol{\theta}^{o}(z)$ has the following evenness and (quasi-)
periodicity properties,
\begin{equation*}
\boldsymbol{\theta}^{o}(-z) \! = \! \boldsymbol{\theta}^{o}(z), \qquad
\boldsymbol{\theta}^{o}(z \! + \! e_{j}) \! = \! \boldsymbol{\theta}^{o}(z),
\qquad \mathrm{and} \qquad \boldsymbol{\theta}^{o}(z \! \pm \! \tau_{j}^{o})
\! = \! \me^{\mp 2 \pi \mi z_{j}-\mi \pi \tau_{jj}^{o}} \boldsymbol{\theta}^{
o}(z),
\end{equation*}
where $\tau_{j}^{o} \! := \! \tau^{o} e_{j}$ $(\in \! \mathbb{C}^{N})$, $j \!
= \! 1,\dotsc,N$. Extensive use of this apparatus will be made in Section~4.
\end{enumerate}
\subsection{The Riemann-Hilbert Problems for the Monic OLPs}
In this subsection, the RHPs corresponding to the even degree and odd 
degree monic OLPs $\boldsymbol{\pi}_{2n}(z)$ and $\boldsymbol{\pi}_{2n+1}
(z)$ defined, respectively, in Equations~(1.4) and~(1.5), are formulated 
\emph{\`{a} la} Fokas-Its-Kitaev \cite{a41,a42}. Furthermore, integral 
representations for the even degree and odd degree monic OLPs are also 
obtained.

Consider the varying exponential measure $\widetilde{\mu}$ $(\in \! \mathcal{
M}_{1}(\mathbb{R}))$ given by $\md \widetilde{\mu}(z) \! = \! \me^{-\mathscr{
N} \, V(z)} \, \md z$, $\mathscr{N} \! \in \! \mathbb{N}$, where (the external
field) $V \colon \mathbb{R} \setminus \{0\} \! \to \! \mathbb{R}$ satisfies
conditions~(V1)--(V3). The RHPs which characterise the even degree and odd
degree monic OLPs are now stated.
\begin{rhp1}
Let $V \colon \mathbb{R} \setminus \{0\} \! \to \! \mathbb{R}$ satisfy
conditions~{\rm (V1)--(V3)}. Find $\overset{e}{\mathrm{Y}} \colon \mathbb{C}
\setminus \mathbb{R} \! \to \! \mathrm{SL}_{2}(\mathbb{C})$ solving: {\rm (i)}
$\overset{e}{\mathrm{Y}}(z)$ is holomorphic for $z \! \in \! \mathbb{C}
\setminus \mathbb{R};$ {\rm (ii)} the boundary values $\overset{e}{\mathrm{
Y}}_{\pm}(z) \! := \! \lim_{\underset{\pm \Im (z^{\prime})>0}{z^{\prime} 
\to z}} \overset{e}{\mathrm{Y}}(z^{\prime})$ satisfy the jump condition
\begin{equation*}
\overset{e}{\mathrm{Y}}_{+}(z)= \overset{e}{\mathrm{Y}}_{-}(z) \! \left(
\mathrm{I} \! + \! \me^{-\mathscr{N} \, V(z)} \sigma_{+} \right), \quad z \!
\in \! \mathbb{R};
\end{equation*}
{\rm (iii)} $\overset{e}{\mathrm{Y}}(z)z^{-n \sigma_{3}} \! =_{\underset{z \in
\mathbb{C} \setminus \mathbb{R}}{z \to \infty}} \! \mathrm{I} \! + \! \mathcal{
O}(z^{-1});$ and {\rm (iv)} $\overset{e}{\mathrm{Y}}(z)z^{n \sigma_{3}} \! =_{
\underset{z \in \mathbb{C} \setminus \mathbb{R}}{z \to 0}} \! \mathcal{O}(1)$.
\end{rhp1}
\begin{rhp2}
Let $V \colon \mathbb{R} \setminus \{0\} \! \to \! \mathbb{R}$ satisfy
conditions~{\rm (V1)--(V3)}. Find $\overset{o}{\mathrm{Y}} \colon \mathbb{C}
\setminus \mathbb{R} \! \to \! \mathrm{SL}_{2}(\mathbb{C})$ solving: {\rm (i)}
$\overset{o}{\mathrm{Y}}(z)$ is holomorphic for $z \! \in \! \mathbb{C}
\setminus \mathbb{R};$ {\rm (ii)} the boundary values $\overset{o}{\mathrm{
Y}}_{\pm}(z) \! := \! \lim_{\underset{\pm \Im (z^{\prime})>0}{z^{\prime} 
\to z}} \overset{o}{\mathrm{Y}}(z^{\prime})$ satisfy the jump condition
\begin{equation*}
\overset{o}{\mathrm{Y}}_{+}(z)=\overset{o}{\mathrm{Y}}_{-}(z) \! \left(
\mathrm{I} \! + \! \me^{-\mathscr{N} \, V(z)} \sigma_{+} \right), \quad z 
\! \in \! \mathbb{R};
\end{equation*}
{\rm (iii)} $\overset{o}{\mathrm{Y}}(z)z^{n \sigma_{3}} \! =_{\underset{z \in
\mathbb{C} \setminus \mathbb{R}}{z \to 0}} \! \mathrm{I} \! + \! \mathcal{O}
(z);$ and {\rm (iv)} $\overset{o}{\mathrm{Y}}(z)z^{-(n+1) \sigma_{3}} \! =_{
\underset{z \in \mathbb{C} \setminus \mathbb{R}}{z \to \infty}} \! \mathcal{O}
(1)$.
\end{rhp2}
\begin{cccc}
Let $\overset{e}{\mathrm{Y}} \colon \mathbb{C} \setminus \mathbb{R} \! \to 
\! \mathrm{SL}_{2}(\mathbb{C})$ solve {\rm \pmb{RHP1}}. {\rm \pmb{RHP1}} 
possesses a unique solution given by: {\rm (i)} for $n \! = \! 0$,
\begin{equation*}
\overset{e}{\mathrm{Y}}(z) \! = \!
\begin{pmatrix}
1 & \int_{\mathbb{R}} \frac{\exp (-\mathscr{N} \, V(s))}{s-z} \, \frac{\md
s}{2 \pi \mi} \\
0 & 1
\end{pmatrix}, \quad z \! \in \! \mathbb{C} \setminus \mathbb{R},
\end{equation*}
where $\bm{\pi}_{0}(z) \! := \! \overset{e}{\mathrm{Y}}_{11}(z) \! \equiv \!
1$, with $\overset{e}{\mathrm{Y}}_{11}(z)$ the $(1 \, 1)$-element of
$\overset{e}{\mathrm{Y}}(z);$ and {\rm (ii)} for $n \! \in \! \mathbb{N}$,
\begin{equation*}
\overset{e}{\mathrm{Y}}(z) \! = \!
\begin{pmatrix}
\boldsymbol{\pi}_{2n}(z) & \int_{\mathbb{R}} \frac{\boldsymbol{\pi}_{2n}(s)
\exp (-\mathscr{N} \, V(s))}{s-z} \, \frac{\md s}{2 \pi \mi} \\
\overset{e}{\mathrm{Y}}_{21}(z) & \int_{\mathbb{R}} \frac{\overset{e}{\mathrm{
Y}}_{21}(s) \exp (-\mathscr{N} \, V(s))}{s-z} \, \frac{\md s}{2 \pi \mi}
\end{pmatrix}, \quad z \! \in \! \mathbb{C} \setminus \mathbb{R},
\end{equation*}
where $\overset{e}{\mathrm{Y}}_{21} \colon \mathbb{C}^{\ast} \! \to \!
\mathbb{C}$ denotes the $(2 \, 1)$-element of $\overset{e}{\mathrm{Y}}(z)$,
and $\boldsymbol{\pi}_{2n}(z)$ is the even degree monic {\rm OLP} defined in
Equation~{\rm (1.4)}.
\end{cccc}

\emph{Proof.} See \cite{a38}, the proof of Lemma~2.2.1. \hfill $\qed$

\begin{cccc}
Let $\overset{o}{\mathrm{Y}} \colon \mathbb{C} \setminus \mathbb{R} \! \to
\! \mathrm{SL}_{2}(\mathbb{C})$ solve {\rm \pmb{RHP2}}. {\rm \pmb{RHP2}}
possesses a unique solution given by: {\rm (i)} for $n \! = \! 0$,
\begin{equation*}
\overset{o}{\mathrm{Y}}(z) \! = \!
\begin{pmatrix}
z \bm{\pi}_{1}(z) & z \int_{\mathbb{R}} \frac{(s \bm{\pi}_{1}(s)) \exp (-
\mathscr{N} \, V(s))}{s(s-z)} \, \frac{\md s}{2 \pi \mi} \\
2 \pi \mi z & 1 \! + \! z \int_{\mathbb{R}} \frac{\exp (-\mathscr{N} \, 
V(s))}{s-z} \, \md s
\end{pmatrix}, \quad z \! \in \! \mathbb{C} \setminus \mathbb{R},
\end{equation*}
where $\bm{\pi}_{1}(z) \! = \! \tfrac{1}{z} \! + \! \tfrac{\xi^{(1)}_{0}}{
\xi^{(1)}_{-1}}$, with $\tfrac{\xi^{(1)}_{0}}{\xi^{(1)}_{-1}} \! = \! -\int_{
\mathbb{R}}s^{-1} \exp (-\mathscr{N} \, V(s)) \, \md s$, $\mathscr{N} \! \in
\! \mathbb{N};$ and {\rm (ii)} for $n \! \in \! \mathbb{N}$,
\begin{equation}
\overset{o}{\mathrm{Y}}(z) \! = \!
\begin{pmatrix}
z \boldsymbol{\pi}_{2n+1}(z) & z \int_{\mathbb{R}} \frac{(s \boldsymbol{\pi}_{
2n+1}(s)) \exp (-\mathscr{N} \, V(s))}{s(s-z)} \, \frac{\md s}{2 \pi \mi} \\
\overset{o}{\mathrm{Y}}_{21}(z) & z \int_{\mathbb{R}} \frac{\overset{o}{
\mathrm{Y}}_{21}(s) \exp (-\mathscr{N} \, V(s))}{s(s-z)} \, \frac{\md s}{2
\pi \mi}
\end{pmatrix}, \quad z \! \in \! \mathbb{C} \setminus \mathbb{R},
\end{equation}
where $\overset{o}{\mathrm{Y}}_{21} \colon \mathbb{C}^{\ast} \! \to \!
\mathbb{C}$ denotes the $(2 \, 1)$-element of $\overset{o}{\mathrm{Y}}(z)$,
and $\boldsymbol{\pi}_{2n+1}(z)$ is the odd degree monic {\rm OLP} defined in
Equation~{\rm (1.5)}.
\end{cccc}

\emph{Proof.} Set $\widetilde{w}(z) \! := \! \exp (-\mathscr{N} \, V(z))$,
$\mathscr{N} \! \in \! \mathbb{N}$, where $V \colon \mathbb{R} \setminus
\{0\} \! \to \! \mathbb{R}$ satisfies conditions~(V1)--(V3). Since $\int_{
\mathbb{R}}s^{j} \widetilde{w}(s) \, \md s \! < \! \infty$, $j \! \in \!
\mathbb{Z}$, and $\langle \bm{\pi}_{1},\bm{\pi}_{0} \rangle_{\mathscr{L}}
\! = \! \langle \bm{\pi}_{1},1 \rangle_{\mathscr{L}} \! = \! 0$, it follows
via an application of the Sokhotski-Plemelj formula (with the Cauchy kernel
normalised at zero) that, for $n \! = \! 0$, \textbf{RHP2} has the (unique)
solution
\begin{equation*}
\overset{o}{\mathrm{Y}}(z) \! = \!
\begin{pmatrix}
z \bm{\pi}_{1}(z) & z \int_{\mathbb{R}} \frac{(s \bm{\pi}_{1}(s)) \widetilde{
w}(s)}{s(s-z)} \, \frac{\md s}{2 \pi \mi} \\
2 \pi \mi z & 1 \! + \! z \int_{\mathbb{R}} \frac{\widetilde{w}(s)}{s-z} \,
\md s
\end{pmatrix}, \quad z \! \in \! \mathbb{C} \setminus \mathbb{R},
\end{equation*}
where $\bm{\pi}_{1}(z) \! = \! \tfrac{1}{z} \! + \! \tfrac{\xi^{(1)}_{0}}{
\xi^{(1)}_{-1}}$, with $\tfrac{\xi^{(1)}_{0}}{\xi^{(1)}_{-1}} \! = \! -\int_{
\mathbb{R}}s^{-1} \widetilde{w}(s) \, \md s$. Hereafter, $n \! \in \! \mathbb{
N}$ will be considered.

If $\overset{o}{\mathrm{Y}} \colon \mathbb{C} \setminus \mathbb{R} \! \to \! 
\operatorname{SL}_{2}(\mathbb{C})$ solves \textbf{RHP2}, then it follows 
{}from the jump condition~(ii) of \textbf{RHP2} that, for the elements of the 
first column of $\overset{o}{\mathrm{Y}}(z)$,
\begin{equation*}
\left(\overset{o}{\mathrm{Y}}_{j1}(z) \right)_{+} \! = \! \left(\overset{o}{
\mathrm{Y}}_{j1}(z) \right)_{-} \! := \! \overset{o}{\mathrm{Y}}_{j1}(z),
\quad j \! = \! 1,2,
\end{equation*}
and, for the elements of the second row,
\begin{equation*}
\left(\overset{o}{\mathrm{Y}}_{j2}(z) \right)_{+} \! - \! \left(\overset{o}{
\mathrm{Y}}_{j2}(z) \right)_{-} \! = \! \overset{o}{\mathrm{Y}}_{j1}(z)
\widetilde{w}(z), \quad j \! = \! 1,2.
\end{equation*}
{}From condition~(i), the normalisation condition~(iii), and the boundedness
condition~(iv) of \textbf{RHP2}, in particular, $\overset{o}{\mathrm{Y}}_{11}
(z)z^{n} \! =_{\underset{z \in \mathbb{C} \setminus \mathbb{R}}{z \to 0}} \!
1 \! + \! \mathcal{O}(z)$, $\overset{o}{\mathrm{Y}}_{11}(z)z^{-(n+1)} \! =_{
\underset{z \in \mathbb{C} \setminus \mathbb{R}}{z \to \infty}} \! \mathcal{O}
(1)$, $\overset{o}{\mathrm{Y}}_{21}(z)z^{n} \! =_{\underset{z \in \mathbb{C}
\setminus \mathbb{R}}{z \to 0}} \! \mathcal{O}(z)$, and $z^{-(n+1)} \overset{
o}{\mathrm{Y}}_{21} \linebreak[4]
(z) \! =_{\underset{z \in \mathbb{C} \setminus \mathbb{R}}{
z \to \infty}} \! \mathcal{O}(1)$, and the fact that $\overset{o}{\mathrm{Y}
}_{11}(z)$ and $\overset{o}{\mathrm{Y}}_{21}(z)$ have no jumps throughout
the $z$-plane, it follows that $\overset{o}{\mathrm{Y}}_{11}(z)$ is a monic
rational function with a pole at the origin and at the point at infinity, with
representation $z^{-1} \overset{o}{\mathrm{Y}}_{11}(z) \! = \! \sum_{l=-n-1}^{
n} \widetilde{\nu}_{l}z^{l}$, where $\widetilde{\nu}_{-n-1} \! = \! 1$, and
$\overset{o}{\mathrm{Y}}_{21}(z)$ is a rational function with a pole at the
origin and at the point at infinity, with representation $z^{-1} \overset{o}{
\mathrm{Y}}_{21}(z) \! = \! \sum_{l=-n}^{n} \nu_{l}^{\flat}z^{l}$. Application
of the Sokhotski-Plemelj formula to the jump relations for $\overset{o}{
\mathrm{Y}}_{j2}(z)$, $j \! = \! 1,2$, gives rise to the following Cauchy-type
integral representations (the Cauchy kernel is normalised at $z \! = \! 0)$:
\begin{equation}
\overset{o}{\mathrm{Y}}_{j2}(z) \! = \! z \int_{\mathbb{R}} \dfrac{\overset{
o}{\mathrm{Y}}_{j1}(s)w(s)}{s(s \! - \! z)} \dfrac{\md s}{2 \pi \mi}, \quad j
\! = 1,2, \quad z \! \in \! \mathbb{C} \setminus \mathbb{R}. \tag{CA1}
\end{equation}

One now studies $\overset{o}{\mathrm{Y}}_{j1}(z)$, $j \! = \! 1,2$, in more
detail. {}From the normalisation condition~(iii) of \textbf{RHP2}, in
particular, $\overset{o}{\mathrm{Y}}_{12}(z)z^{-n} \! =_{\underset{z \in
\mathbb{C} \setminus \mathbb{R}}{z \to 0}} \! \mathcal{O}(z)$ and $\overset{
o}{\mathrm{Y}}_{22}(z)z^{-n} \! =_{\underset{z \in \mathbb{C} \setminus
\mathbb{R}}{z \to 0}} \! 1 \! + \! \mathcal{O}(z)$, the formulae~(CA1), the
fact that $\int_{\mathbb{R}}s^{j} \widetilde{w}(s) \, \md s \! < \! \infty$,
$j \! \in \! \mathbb{Z}$, and the expansion (for $\vert z/s \vert \! \ll \!
1)$ $\tfrac{1}{s-z} \! = \! \sum_{k=0}^{l} \tfrac{z^{k}}{s^{k+1}} \! + \!
\tfrac{z^{l+1}}{s^{l+1}(s-z)}$, $l \! \in \! \mathbb{Z}_{0}^{+}$, it follows
that
\begin{equation*}
\int_{\mathbb{R}}(s^{-1} \overset{o}{\mathrm{Y}}_{11}(s))s^{-k} \widetilde{w}
(s) \, \md s \! = \! 0, \quad k \! = \! 1,2,\dotsc,n, \qquad \text{and} \qquad
\int_{\mathbb{R}}(s^{-1} \overset{o}{\mathrm{Y}}_{11}(s))s^{-(n+1)}
\widetilde{w}(s) \, \md s \! = \! -2 \pi \mi \mathfrak{p}^{o},
\end{equation*}
for some (pure imaginary) $\mathfrak{p}^{o}$ of the form $\mathfrak{p}^{o} \!
= \! \mi \mathfrak{q}^{o}$, with $\mathfrak{q}^{o} \! > \! 0$ (see below), and
\begin{equation*}
\int_{\mathbb{R}}(s^{-1} \overset{o}{\mathrm{Y}}_{21}(s))s^{-j} \widetilde{w}
(s) \, \md s \! = \! 0, \quad j \! = \! 1,2,\dotsc,n \! - \! 1, \qquad
\text{and} \qquad \int_{\mathbb{R}}(s^{-1} \overset{o}{\mathrm{Y}}_{21}(s))
s^{-n} \widetilde{w}(s) \, \md s \! = \! 2 \pi \mi;
\end{equation*}
and, {}from the boundedness condition~(iv) of \textbf{RHP2}, in particular,
$\overset{o}{\mathrm{Y}}_{12}(z)z^{n+1} \! =_{\underset{z \in \mathbb{C}
\setminus \mathbb{R}}{z \to \infty}} \! \mathcal{O}(1)$ and $\overset{o}{
\mathrm{Y}}_{22}(z) \linebreak[4]
\cdot z^{n+1} \! =_{\underset{z \in \mathbb{C} \setminus \mathbb{R}}{z \to
\infty}} \! \mathcal{O}(1)$, the formulae~(CA1), the fact that $\int_{\mathbb{
R}}s^{j} \widetilde{w}(s) \, \md s \! < \! \infty$, $j \! \in \! \mathbb{Z}$,
and the expansion (for $\vert s/z \vert \! \ll \! 1)$ $\tfrac{1}{s-z} \! = \!
-\sum_{k=0}^{l} \tfrac{s^{k}}{z^{k+1}} \! + \! \tfrac{s^{l+1}}{z^{l+1}(s-z)}$,
$l \! \in \! \mathbb{Z}_{0}^{+}$, it follows that
\begin{equation*}
\int_{\mathbb{R}}(s^{-1} \overset{o}{\mathrm{Y}}_{11}(s))s^{k} \widetilde{w}
(s) \, \md s \! = \! 0, \quad k \! = \! 0,1,\dotsc,n, \quad \text{and} \quad
\int_{\mathbb{R}}(s^{-1} \overset{o}{\mathrm{Y}}_{21}(s))s^{j} \widetilde{w}
(s) \, \md s \! = \! 0, \quad j \! = \! 0,1,\dotsc,n:
\end{equation*}
these give rise to $2n \! + \! 2$ conditions for $z^{-1} \overset{o}{\mathrm{
Y}}_{11}(z)$, and $2n \! + \! 1$ conditions for $z^{-1} \overset{o}{\mathrm{
Y}}_{21}(z)$. Consider, first, the $2n \! + \! 1$ conditions for $z^{-1}
\overset{o}{\mathrm{Y}}_{21}(z)$. Recalling that the strong moments are
defined by $c_{j} \! := \! \int_{\mathbb{R}}s^{j} \widetilde{w}(s) \, \md s$,
$j \! \in \! \mathbb{Z}$, it follows {}from the representation (established
above) $z^{-1} \overset{o}{\mathrm{Y}}_{21}(z) \! = \! \sum_{l=-n}^{n} \nu_{
l}^{\flat}z^{l}$ and the $2n \! + \! 1$ conditions for $z^{-1} \overset{o}{
\mathrm{Y}}_{21}(z)$ that
\begin{equation*}
\sum_{l=-n}^{n} \nu_{l}^{\flat}c_{l+k} \! = \! 0, \quad k \! = \! -(n \! - \!
1),-(n \! - \! 2),\dotsc,n, \qquad \quad \, \, \, \text{and} \, \, \, \qquad 
\quad \sum_{l=-n}^{n} \nu_{l}^{\flat}c_{l-n} \! = \! 2 \pi \mi,
\end{equation*}
that is,
\begin{equation*}
\begin{pmatrix}
c_{-2n} & c_{-2n+1} & \dotsb & c_{-1} & c_{0} \\
c_{-2n+1} & c_{-2n+2} & \dotsb & c_{0} & c_{1} \\
\vdots & \vdots & \ddots & \vdots & \vdots \\
c_{1} & c_{2} & \dotsb & c_{2n-2} & c_{2n-1} \\
c_{0} & c_{1} & \dotsb & c_{2n-1} & c_{2n}
\end{pmatrix} \!
\begin{pmatrix}
\nu_{-n}^{\flat} \\
\nu_{-n+1}^{\flat} \\
\vdots \\
\nu_{n-1}^{\flat} \\
\nu_{n}^{\flat}
\end{pmatrix} \! = \!
\begin{pmatrix}
2 \pi \mi \\
0 \\
\vdots \\
0 \\
0
\end{pmatrix}.
\end{equation*}
This linear system of $2n \! + \! 1$ equations for the $2n \! + \! 1$ unknowns
$\nu_{l}^{\flat}$, $l \! = \! -n,-(n \! - \! 1),\dotsc,n$, admits a unique
solution if, and only if, the determinant of the coefficient matrix, in this
case $H^{(-2n)}_{2n+1}$ (cf.~Equations~(1.1)), is non-zero; in fact, it will
be shown that $H^{(-2n)}_{2n+1} \! > \! 0$. An integral representation for the
Hankel determinants $H^{(m)}_{k}$, $(m,k) \! \in \! \mathbb{Z} \times \mathbb{
N}$, is now obtained; then the substitutions $m \! = \! -2n$ and $k \! = \!
2n \! + \! 1$ are made. In the calculations that follow, $\mathfrak{S}_{k}$
denotes the $k!$ permutations $\boldsymbol{\sigma}$ of $\lbrace 1,2,\dotsc,
k \rbrace$. Recalling that $c_{j} \! := \! \int_{\mathbb{R}}s^{j} \, \md
\widetilde{\mu}(s)$, $j \! \in \! \mathbb{Z}$, where $\md \widetilde{\mu}(z)
\! = \! \widetilde{w}(z) \, \md z \! = \! \exp (-\mathscr{N} \, V(z)) \,
\md z$, and using the multi-linearity property of the determinant, via
Equations~(1.1), one proceeds thus (recall
that $H^{(m)}_{0} \! := \! 1)$:
\begin{align*}
H^{(m)}_{k} :=& \,
\begin{vmatrix}
c_{m} & c_{m+1} & \dotsb & c_{m+k-1} \\
c_{m+1} & c_{m+2} & \dotsb & c_{m+k} \\
\vdots & \vdots & \ddots & \vdots \\
c_{m+k-2} & c_{m+k-1} & \dotsb & c_{m+2k-3} \\
c_{m+k-1} & c_{m+k} & \dotsb & c_{m+2k-2}
\end{vmatrix} \\
=& \,
\begin{vmatrix}
\int_{\mathbb{R}}s_{1}^{m} \, \md \widetilde{\mu}(s_{1}) & \int_{\mathbb{R}}
s_{2}^{m+1} \, \md \widetilde{\mu}(s_{2}) & \dotsb & \int_{\mathbb{R}}s_{k}^{
m+k-1} \, \md \widetilde{\mu}(s_{k}) \\
\int_{\mathbb{R}}s_{1}^{m+1} \, \md \widetilde{\mu}(s_{1}) & \int_{\mathbb{R}}
s_{2}^{m+2} \, \md \widetilde{\mu}(s_{2}) & \dotsb & \int_{\mathbb{R}}s_{k}^{
m+k} \, \md \widetilde{\mu}(s_{k}) \\
\vdots & \vdots & \ddots & \vdots \\
\int_{\mathbb{R}}s_{1}^{m+k-2} \, \md \widetilde{\mu}(s_{1}) & \int_{\mathbb{
R}}s_{2}^{m+k-1} \, \md \widetilde{\mu}(s_{2}) & \dotsb & \int_{\mathbb{R}}
s_{k}^{m+2k-3} \, \md \widetilde{\mu}(s_{k}) \\
\int_{\mathbb{R}}s_{1}^{m+k-1} \, \md \widetilde{\mu}(s_{1}) & \int_{\mathbb{
R}}s_{2}^{m+k} \, \md \widetilde{\mu}(s_{2}) & \dotsb & \int_{\mathbb{R}}s_{
k}^{m+2k-2} \, \md \widetilde{\mu}(s_{k})
\end{vmatrix} \\
=& \, \underbrace{\int_{\mathbb{R}} \int_{\mathbb{R}} \dotsi \int_{\mathbb{R}}
}_{k} \, \, \md \widetilde{\mu}(s_{1}) \, \md \widetilde{\mu}(s_{2}) \, \dotsb
\, \md \widetilde{\mu}(s_{k})
\begin{vmatrix}
s_{1}^{m} & s_{2}^{m+1} & \dotsb & s_{k}^{m+k-1} \\
s_{1}^{m+1} & s_{2}^{m+2} & \dotsb & s_{k}^{m+k} \\
\vdots & \vdots & \ddots & \vdots \\
s_{1}^{m+k-2} & s_{2}^{m+k-1} & \dotsb & s_{k}^{m+2k-3} \\
s_{1}^{m+k-1} & s_{2}^{m+k} & \dotsb & s_{k}^{m+2k-2}
\end{vmatrix} \\
=& \, \underbrace{\int_{\mathbb{R}} \int_{\mathbb{R}} \dotsi \int_{\mathbb{R}}
}_{k} \, \, \md \widetilde{\mu}(s_{1}) \, \md \widetilde{\mu}(s_{2}) \, \dotsb
\, \md \widetilde{\mu}(s_{k}) \,
s_{1}^{m}s_{2}^{m+1} \dotsb s_{k}^{m+k-1}
\underbrace{\begin{vmatrix}
1 & 1 & \dotsb & 1 \\
s_{1} & s_{2} & \dotsb & s_{k} \\
\vdots & \vdots & \ddots & \vdots \\
s_{1}^{k-2} & s_{2}^{k-2} & \dotsb & s_{k}^{k-2} \\
s_{1}^{k-1} & s_{2}^{k-1} & \dotsb & s_{k}^{k-1}
\end{vmatrix}}_{=: \, \mathrm{V}(s_{1},s_{2},\dotsc,s_{k})} \\
=& \, \dfrac{1}{k!} \sum_{\boldsymbol{\sigma} \in \mathfrak{S}_{k}} \,
\underbrace{\int_{\mathbb{R}} \int_{\mathbb{R}} \dotsi \int_{\mathbb{R}}}_{k}
\, \, \md \widetilde{\mu}(s_{\sigma (1)}) \, \md \widetilde{\mu}(s_{\sigma (2)
}) \, \dotsb \, \md \widetilde{\mu}(s_{\sigma (k)}) \prod_{j=1}^{k}s_{\sigma
(j)}^{m}s_{\sigma (j)}^{j-1} \, \mathrm{V} \! \left(s_{\sigma (1)},s_{\sigma
(2)},\dotsc,s_{\sigma (k)} \right) \\
=& \, \dfrac{1}{k!} \sum_{\boldsymbol{\sigma} \in \mathfrak{S}_{k}} \,
\underbrace{\int_{\mathbb{R}} \int_{\mathbb{R}} \dotsi \int_{\mathbb{R}}}_{k}
\, \, \md \widetilde{\mu}(s_{1}) \, \md \widetilde{\mu}(s_{2}) \, \dotsb \,
\md \widetilde{\mu}(s_{k}) \, s_{1}^{m}s_{2}^{m} \dotsb s_{k}^{m} \! \left(
\operatorname{sgn}(\boldsymbol{\sigma}) \prod_{j=1}^{k}s_{\sigma (j)}^{j-1}
\right) \\
\times& \, \mathrm{V} \! \left(s_{1},s_{2},\dotsc,s_{k} \right) \\
=& \, \dfrac{1}{k!} \underbrace{\int_{\mathbb{R}} \int_{\mathbb{R}} \dotsi
\int_{\mathbb{R}}}_{k} \, \, \md \widetilde{\mu}(s_{1}) \, \md \widetilde{\mu}
(s_{2}) \, \dotsb \, \md \widetilde{\mu}(s_{k}) \, s_{1}^{m}s_{2}^{m} \dotsb
s_{k}^{m} \, \mathrm{V} \! \left(s_{1},s_{2},\dotsc,s_{k} \right) \\
\times& \, \underbrace{\sum_{\boldsymbol{\sigma} \in \mathfrak{S}_{k}}
\operatorname{sgn}(\boldsymbol{\sigma})s_{\sigma (1)}^{0}s_{\sigma (2)}^{1}
\dotsb s_{\sigma (k)}^{k-1}}_{= \, \mathrm{V}(s_{1},s_{2},\dotsc,s_{k})} 
\quad \Rightarrow \\
H^{(m)}_{k} =& \, \dfrac{1}{k!} \underbrace{\int_{\mathbb{R}} \int_{\mathbb{
R}} \dotsi \int_{\mathbb{R}}}_{k} \, \, \md \widetilde{\mu}(s_{1}) \, \md
\widetilde{\mu}(s_{2}) \, \dotsb \, \md \widetilde{\mu}(s_{k}) \, s_{1}^{m}
s_{2}^{m} \dotsb s_{k}^{m} \! \left(\mathrm{V}(s_{1},s_{2},\dotsc,s_{k})
\right)^{2};
\end{align*}
using the well-known determinantal formula $\mathrm{V}(s_{1},s_{2},\dotsc,
s_{k}) \! = \! \prod_{\underset{j<i}{i,j=1}}^{k}(s_{i} \! - \! s_{j})$, one
arrives at
\begin{equation}
H^{(m)}_{k} \! = \! \dfrac{1}{k!} \underbrace{\int_{\mathbb{R}} \int_{\mathbb{
R}} \dotsi \int_{\mathbb{R}}}_{k} \, \, s_{1}^{m}s_{2}^{m} \dotsb s_{k}^{m}
\prod_{\underset{l<i}{i,l=1}}^{k}(s_{i} \! - \! s_{l})^{2} \, \md \widetilde{
\mu}(s_{1}) \, \md \widetilde{\mu}(s_{2}) \, \dotsb \, \md \widetilde{\mu}
(s_{k}), \quad (m,k) \! \in \! \mathbb{Z} \times \mathbb{N} \tag{HA1}.
\end{equation}
Letting $m \! = \! -2n$ and $k \! = \! 2n \! + \! 1$, it follows from the
formula~(HA1) that
\begin{equation*}
H^{(-2n)}_{2n+1} \! = \! \dfrac{1}{(2n \! + \! 1)!} \underbrace{\int_{\mathbb{
R}} \int_{\mathbb{R}} \dotsi \int_{\mathbb{R}}}_{2n+1} \, \, s_{1}^{-2n}
s_{2}^{-2n} \dotsb s_{2n+1}^{-2n} \prod_{\underset{l<i}{i,l=1}}^{2n+1}(s_{i}
\! - \! s_{l})^{2} \, \md \widetilde{\mu}(s_{1}) \, \md \widetilde{\mu}(s_{2})
\, \dotsb \, \md \widetilde{\mu}(s_{2n+1}) \! > \! 0,
\end{equation*}
whence the existence (and uniqueness) of $z^{-1} \overset{o}{\mathrm{Y}}_{21}
(z)$ (thus $\overset{o}{\mathrm{Y}}_{21}(z))$.

Similarly, it follows, {}from the representation (established above) $z^{-1}
\overset{o}{\mathrm{Y}}_{11}(z) \! = \! \sum_{l=-n-1}^{n} \widetilde{\nu}_{
l}z^{l}$, with $\widetilde{\nu}_{-n-1} \! = \! 1$, and the $2n \! + \! 2$
conditions for $z^{-1} \overset{o}{\mathrm{Y}}_{11}(z)$, that
\begin{equation*}
\sum_{l=-n-1}^{n} \widetilde{\nu}_{l}c_{l+k} \! = \! 0, \quad k \! = \! -n,
-(n \! - \! 1),\dotsc,n, \qquad \quad \, \, \text{and} \qquad \quad \, \, 
\sum_{l=-n-1}^{n} \widetilde{\nu}_{l}c_{l-n-1} \! = \! -2 \pi \mi 
\mathfrak{p}^{o},
\end{equation*}
that is,
\begin{equation*}
\begin{pmatrix}
2 \pi \mi & c_{-2n-1} & \dotsb  & c_{-n-1} & \dotsb & c_{-1} \\
0         & c_{-2n}   & \dotsb  & c_{-n}   & \dotsb & c_{0} \\
\vdots    & \vdots    & \ddots  & \vdots   & \ddots & \vdots \\
0         & c_{-n}    & \dotsb  & c_{0}    & \dotsb & c_{n} \\
\vdots    & \vdots    & \ddots  & \vdots   & \ddots & \vdots \\
0         & c_{0}     & \dotsb  & c_{n}    & \dotsb & c_{2n}
\end{pmatrix} \!
\begin{pmatrix}
\mathfrak{p}^{o} \\
\widetilde{\nu}_{-n} \\
\vdots \\
\widetilde{\nu}_{0} \\
\vdots \\
\widetilde{\nu}_{n}
\end{pmatrix} \! = \!
\begin{pmatrix}
-c_{-2(n+1)} \\
-c_{-2n-1} \\
\vdots \\
-c_{-n-1} \\
\vdots \\
-c_{-1}
\end{pmatrix}.
\end{equation*}
This linear system of $2n \! + \! 2$ equations for the $2n \! + \! 2$ unknowns
$\widetilde{\nu}_{l}$, $l \! = \! -n,-(n \! - \! 1),\dotsb,n$, and $\mathfrak{
p}^{o}$ admits a unique solution if, and only if, the determinant of the
coefficient matrix, in this case $2 \pi \mi H^{(-2n)}_{2n+1}$, is non-zero;
but, it was shown above that $H^{(-2n)}_{2n+1} \! > \! 0$. Furthermore, via
Cramer's Rule:
\begin{equation*}
\mathfrak{p}^{o} \! = \! \dfrac{
\begin{vmatrix}
-c_{-2n-2}   & c_{-2n-1} & \dotsb  & c_{-n-1}    & \dotsb & c_{-1} \\
-c_{-2n-1}   & c_{-2n}   & \dotsb  & c_{-n}      & \dotsb & c_{0} \\
\vdots       & \vdots    & \ddots  & \vdots      & \ddots & \vdots \\
-c_{-n-1}    & c_{-n}    & \dotsb  & c_{0}       & \dotsb & c_{n} \\
\vdots       & \vdots    & \ddots  & \vdots      & \ddots & \vdots \\
-c_{-1}      & c_{0}     & \dotsb  & c_{n}       & \dotsb & c_{2n}
\end{vmatrix}}{2 \pi \mi H^{(-2n)}_{2n+1}} \! = \! -\dfrac{1}{2 \pi \mi}
\dfrac{H^{(-2n-2)}_{2n+2}}{H^{(-2n)}_{2n+1}}.
\end{equation*}
Using the Hankel determinant formula~(HA1) with the substitutions $m \! = \!
-2(n \! + \! 1)$ and $k \! = \! 2(n \! + \! 1)$, one arrives at
\begin{align*}
H^{(-2n-2)}_{2n+2}=& \, \dfrac{1}{(2n \! + \! 2)!} \underbrace{\int_{\mathbb{
R}} \int_{\mathbb{R}} \dotsi \int_{\mathbb{R}}}_{2n+2} \, \, s_{1}^{-2n-2}
s_{2}^{-2n-2} \dotsb s_{2n+2}^{-2n-2} \prod_{\underset{l<i}{i,l=1}}^{2n+2}
(s_{i} \! - \! s_{l})^{2} \\
\times& \, \md \widetilde{\mu}(s_{1}) \, \md \widetilde{\mu}(s_{2}) \, \dotsb
\, \md \widetilde{\mu}(s_{2n+2}) \! > \! 0;
\end{align*}
hence, $H^{(-2n-2)}_{2n+2}/H^{(-2n)}_{2n+1} \! > \! 0$. Using, now, the fact
that $\int_{\mathbb{R}}(s^{-1} \overset{o}{\mathrm{Y}}_{11}(s))s^{k}
\widetilde{w}(s) \, \md s \! = \! 0$, $k \! = \! -n,-(n \! - \! 1),\dotsc,n$,
and the relation $\int_{\mathbb{R}}(s^{-1} \overset{o}{\mathrm{Y}}_{11}(s))
s^{-(n+1)} \widetilde{w}(s) \, \md s \! = \! -2 \pi \mi \mathfrak{p}^{o}$, one
notes, via the above formula for $\mathfrak{p}^{o}$, that
\begin{align*}
\int_{\mathbb{R}}(s^{-1} \overset{o}{\mathrm{Y}}_{11}(s))s^{-(n+1)} 
\widetilde{w}(s) \, \md s =& \, \int_{\mathbb{R}}(s^{-1} \overset{o}{\mathrm{
Y}}_{11}(s)) \underbrace{\left(s^{-(n+1)} \! + \! \widetilde{\nu}_{-n}s^{-
n} \! + \! \dotsb \! + \! \widetilde{\nu}_{n}s^{n} \right)}_{= \, s^{-1}
\overset{o}{\mathrm{Y}}_{11}(s)} \, \widetilde{w}(s) \, \md s \\
=& \, \int_{\mathbb{R}}(s^{-1} \overset{o}{\mathrm{Y}}_{11}(s))^{2}
\widetilde{w}(s) \, \md s \! = \! -2 \pi \mi \mathfrak{p}^{o} \! = \!
H^{(-2n-2)}_{2n+2}/H^{(-2n)}_{2n+1} \quad (> \! 0);
\end{align*}
but the right-hand side of the latter expression is equal to $(\xi^{(2n+1)}_{-
n-1})^{-2} \! = \! \lvert \lvert \pmb{\ast}^{-1} \overset{o}{\mathrm{Y}}_{11}
(\pmb{\ast}) \lvert \lvert_{\mathscr{L}}^{2}$ $(> \! 0)$ (cf.
E\-q\-u\-a\-t\-i\-o\-n\-s (1.8)): the existence and uniqueness of $z^{-1}
\overset{o}{\mathrm{Y}}_{11}(z) \! =: \! \boldsymbol{\pi}_{2n+1}(z)$, the 
odd degree monic OLP with respect to the inner product $\langle \pmb{\cdot},
\pmb{\cdot} \rangle_{\mathscr{L}}$, is thus established. \hfill $\qed$
\begin{ffff}
Let $V \colon \mathbb{R} \setminus \{0\} \! \to \! \mathbb{R}$ satisfy
conditions~{\rm (V1)--(V3)}. Let $\boldsymbol{\pi}_{2n}(z)$ and $\boldsymbol{
\pi}_{2n+1}(z)$ be the even degree and odd degree monic {\rm OLPs} with
respect to the inner product $\langle \pmb{\cdot},\pmb{\cdot} \rangle_{
\mathscr{L}}$ defined, respectively, in Equations~{\rm (1.4)} and~{\rm (1.5)},
and let $\xi^{(2n)}_{n}$ and $\xi^{(2n+1)}_{-n-1}$ be the corresponding `even'
and `odd' norming constants, respectively. Then, $\xi^{(2n)}_{n}$ and $\xi^{
(2n+1)}_{-n-1}$ have the following representations:
\begin{equation*}
\dfrac{\xi^{(2n)}_{n}}{\sqrt{\smash[b]{2n \! + \! 1}}} = \sqrt{\dfrac{
\underbrace{\int_{\mathbb{R}} \int_{\mathbb{R}} \dotsi \int_{\mathbb{R}}}_{2n}
\, \, s_{1}^{-2n}s_{2}^{-2n} \dotsb s_{2n}^{-2n} \mathlarger{\prod_{\underset{
l<i}{i,l=1}}^{2n}}(s_{i} \! - \! s_{l})^{2} \, \md \widetilde{\mu}(s_{1}) \,
\md \widetilde{\mu}(s_{2}) \, \dotsb \, \md \widetilde{\mu}(s_{2n})}{
\underbrace{\int_{\mathbb{R}} \int_{\mathbb{R}} \dotsi \int_{\mathbb{R}}}_{2
n+1} \, \, \lambda_{1}^{-2n} \lambda_{2}^{-2n} \dotsb \lambda_{2n+1}^{-2n}
\mathlarger{\prod_{\underset{l<i}{i,l=1}}^{2n+1}}(\lambda_{i} \! - \!
\lambda_{l})^{2} \, \md \widetilde{\mu}(\lambda_{1}) \, \md \widetilde{\mu}
(\lambda_{2}) \, \dotsb \, \md \widetilde{\mu}(\lambda_{2n+1})}},
\end{equation*}
\begin{equation*}
\dfrac{\xi^{(2n+1)}_{-n-1}}{\sqrt{\smash[b]{2(n \! + \! 1)}}} = \sqrt{
\dfrac{\underbrace{\int_{\mathbb{R}} \int_{\mathbb{R}} \dotsi \int_{\mathbb{
R}}}_{2n+1} \, \, \varpi_{1}^{-2n} \varpi_{2}^{-2n} \dotsb \varpi_{2n+1}^{-2n}
\mathlarger{\prod_{\underset{l<i}{i,l=1}}^{2n+1}}(\varpi_{i} \! - \! \varpi_{
l})^{2} \, \md \widetilde{\mu}(\varpi_{1}) \, \md \widetilde{\mu}(\varpi_{2})
\, \dotsb \, \md \widetilde{\mu}(\varpi_{2n+1})}{\underbrace{\int_{\mathbb{R}}
\int_{\mathbb{R}} \dotsi \int_{\mathbb{R}}}_{2n+2} \, \, \varsigma_{1}^{-2n-2}
\varsigma_{2}^{-2n-2} \dotsb \varsigma_{2n+2}^{-2n-2} \mathlarger{\prod_{
\underset{l<i}{i,l=1}}^{2n+2}}(\varsigma_{i} \! - \! \varsigma_{l})^{2} \, \md
\widetilde{\mu}(\varsigma_{1}) \, \md \widetilde{\mu}(\varsigma_{2}) \, \dotsb
\, \md \widetilde{\mu}(\varsigma_{2n+2})}},
\end{equation*}
where $\md \widetilde{\mu}(z) \! := \! \exp (-\mathscr{N} \, V(z)) \, \md 
z$, $\mathscr{N} \! \in \! \mathbb{N}$.
\end{ffff}

\emph{Proof.} Consider, without loss of generality, the representation for
$\xi^{(2n+1)}_{-n-1}$. Recall that (cf. Equations (1.8)) $(\xi^{(2n+1)}_{-n-
1})^{2} \! = \! H^{(-2n)}_{2n+1}/H^{(-2n-2)}_{2n+2}$ $(> \! 0)$: using the
integral representations for $H^{(-2n)}_{2n+1}$ and $H^{(-2n-2)}_{2n+2}$
derived in (the course of) the proof of Lemma~2.2.2, and taking positive
square roots of both sides of the resulting equality, one arrives at the
representation for $\xi^{(2n+1)}_{-n-1}$. See \cite{a38}, Corollary~2.2.1,
for the proof of the representation for $\xi^{(2n)}_{n}$. \hfill $\qed$
\begin{bbbb}
Let $V \colon \mathbb{R} \setminus \{0\} \! \to \! \mathbb{R}$ satisfy
conditions~{\rm (V1)--(V3)}. Let $\boldsymbol{\pi}_{2n}(z)$ and $\boldsymbol{
\pi}_{2n+1}(z)$ be the even degree and odd degree monic {\rm OLPs} with
respect to the inner product $\langle \pmb{\cdot},\pmb{\cdot} \rangle_{
\mathscr{L}}$ defined, respectively, in Equations~{\rm (1.4)} and~{\rm (1.5)}.
Then, $\boldsymbol{\pi}_{2n}(z)$ and $\boldsymbol{\pi}_{2n+1}(z)$ have,
respectively, the following integral representations:
\begin{align*}
\boldsymbol{\pi}_{2n}(z) =& \, \dfrac{z^{-n}}{(2n)!H^{(-2n)}_{2n}} \,
\underbrace{\int_{\mathbb{R}} \int_{\mathbb{R}} \dotsi \int_{\mathbb{R}}}_{2n}
\, \, s_{0}^{-2n}s_{1}^{-2n} \dotsb s_{2n-1}^{-2n} \, \prod_{\underset{l<i}{i,
l=0}}^{2n-1}(s_{i} \! - \! s_{l})^{2} \, \prod_{j=0}^{2n-1}(z \! - \! s_{j}) \\
\times& \, \md \widetilde{\mu}(s_{0}) \, \md \widetilde{\mu}(s_{1}) \, \dotsb
\, \md \widetilde{\mu}(s_{2n-1}),
\end{align*}
\begin{align*}
\boldsymbol{\pi}_{2n+1}(z) =& \, -\dfrac{z^{-n-1}}{(2n \! + \! 1)!H^{(-2n)}_{2
n+1}} \, \underbrace{\int_{\mathbb{R}} \int_{\mathbb{R}} \dotsi \int_{\mathbb{
R}}}_{2n+1} \, \, s_{0}^{-2n-1}s_{1}^{-2n-1} \dotsb s_{2n}^{-2n-1} \, \prod_{
\underset{l<i}{i,l=0}}^{2n}(s_{i} \! - \! s_{l})^{2} \, \prod_{j=0}^{2n}(z \!
- \! s_{j}) \\
\times& \, \md \widetilde{\mu}(s_{0}) \, \md \widetilde{\mu}(s_{1}) \, \dotsb
\, \md \widetilde{\mu}(s_{2n}),
\end{align*}
where
\begin{gather*}
H^{(-2n)}_{2n} \! = \! \dfrac{1}{(2n)!} \underbrace{\int_{\mathbb{R}} \int_{
\mathbb{R}} \dotsi \int_{\mathbb{R}}}_{2n} \, \, \lambda_{1}^{-2n} \lambda_{
2}^{-2n} \dotsb \lambda_{2n}^{-2n} \prod_{\underset{l<i}{i,l=1}}^{2n}
(\lambda_{i} \! - \! \lambda_{l})^{2} \, \md \widetilde{\mu}(\lambda_{1}) \,
\md \widetilde{\mu}(\lambda_{2}) \, \dotsb \, \md \widetilde{\mu}(\lambda_{2
n}), \\
H^{(-2n)}_{2n+1} \! = \! \dfrac{1}{(2n \! + \! 1)!} \underbrace{\int_{\mathbb{
R}} \int_{\mathbb{R}} \dotsi \int_{\mathbb{R}}}_{2n+1} \, \, \lambda_{1}^{-2n}
\lambda_{2}^{-2n} \dotsb \lambda_{2n+1}^{-2n} \prod_{\underset{l<i}{i,l=1}}^{2
n+1}(\lambda_{i} \! - \! \lambda_{l})^{2} \, \md \widetilde{\mu}(\lambda_{1})
\, \md \widetilde{\mu}(\lambda_{2}) \, \dotsb \, \md \widetilde{\mu}(\lambda_{
2n+1}),
\end{gather*}
with $\md \widetilde{\mu}(z) \! := \! \exp (-\mathscr{N} \, V(z)) \, \md z$,
$\mathscr{N} \! \in \! \mathbb{N}$.
\end{bbbb}

\emph{Proof.} Consider, without loss of generality, the integral
representation for the odd degree monic OLP $\boldsymbol{\pi}_{2n+1}(z)$.
Let $\mathfrak{S}_{k}$ denote the $k!$ permutations $\boldsymbol{\sigma}$
of $\lbrace 0,1,\dotsc,k \! - \! 1 \rbrace$. Recalling that $c_{j} \! := \!
\int_{\mathbb{R}}s^{j} \, \md \widetilde{\mu}(s)$, $j \! \in \! \mathbb{Z}$,
where $\md \widetilde{\mu}(z) \! := \! \widetilde{w}(z) \, \md z \! = \! \exp
(-\mathscr{N} \, V(z)) \, \md z$, $\mathscr{N} \! \in \! \mathbb{N}$, with
$V \colon \mathbb{R} \setminus \{0\} \! \to \! \mathbb{R}$ satisfying
conditons~(V1)--(V3), and using the multi-linearity property of the
determinant, via the determinantal representation for $\boldsymbol{\pi}_{2n+1}
(z)$ given in Equation~(1.7), one proceeds thus:
\begin{align*}
\boldsymbol{\pi}_{2n+1}(z) =& \, -\dfrac{1}{H^{(-2n)}_{2n+1}}
\begin{vmatrix}
c_{-2n-1} & c_{-2n} & \dotsb & c_{-1}  & z^{-n-1} \\
c_{-2n} & c_{-2n+1} & \dotsb & c_{0} & z^{-n} \\
\vdots & \vdots & \ddots & \vdots & \vdots \\
c_{-1} & c_{0} & \dotsb & c_{2n-1} & z^{n-1} \\
c_{0} & c_{1} & \dotsb & c_{2n} & z^{n}
\end{vmatrix} \\
=& \, -\dfrac{z^{-n-1}}{H^{(-2n)}_{2n+1}}
\begin{vmatrix}
c_{-2n-1} & c_{-2n} & \dotsb & c_{-1}  & c_{0} \\
c_{-2n} & c_{-2n+1} & \dotsb & c_{0} & c_{1} \\
\vdots & \vdots & \ddots & \vdots & \vdots \\
c_{-1} & c_{0} & \dotsb & c_{2n-1} & c_{2n} \\
z^{0} & z^{1} & \dotsb & z^{2n} & z^{2n+1}
\end{vmatrix} \\
=& \, -\dfrac{z^{-n-1}}{H^{(-2n)}_{2n+1}}
\begin{vmatrix}
\int_{\mathbb{R}}s_{0}^{-2n-1} \, \md \widetilde{\mu}(s_{0}) & \int_{\mathbb{
R}}s_{0}^{-2n} \, \md \widetilde{\mu}(s_{0}) & \dotsb & \int_{\mathbb{R}}s_{
0}^{0} \, \md \widetilde{\mu}(s_{0}) \\
\int_{\mathbb{R}}s_{1}^{-2n} \, \md \widetilde{\mu}(s_{1}) & \int_{\mathbb{R}}
s_{1}^{-2n+1} \, \md \widetilde{\mu}(s_{1}) & \dotsb & \int_{\mathbb{R}}s_{
1}^{1} \, \md \widetilde{\mu}(s_{1}) \\
\vdots & \vdots & \ddots & \vdots \\
\int_{\mathbb{R}}s_{2n}^{-1} \, \md \widetilde{\mu}(s_{2n}) & \int_{\mathbb{
R}}s_{2n}^{0} \, \md \widetilde{\mu}(s_{2n}) & \dotsb & \int_{\mathbb{R}}s_{2
n}^{2n} \, \md \widetilde{\mu}(s_{2n}) \\
z^{0} & z^{1} & \dotsb & z^{2n+1}
\end{vmatrix} \\
=& \, -\dfrac{z^{-n-1}}{H^{(-2n)}_{2n+1}} \, \underbrace{\int_{\mathbb{R}}
\int_{\mathbb{R}} \dotsi \int_{\mathbb{R}}}_{2n+1} \, \, \md \widetilde{\mu}
(s_{0}) \md \widetilde{\mu}(s_{1}) \, \dotsb \, \md \widetilde{\mu}(s_{2n})
\begin{vmatrix}
s_{0}^{-2n-1} & s_{0}^{-2n} & \dotsb & s_{0}^{-1} & s_{0}^{0} \\
s_{1}^{-2n} & s_{1}^{-2n+1} & \dotsb & s_{1}^{0} & s_{1}^{1} \\
\vdots & \vdots & \ddots & \vdots & \vdots \\
s_{2n}^{-1} & s_{2n}^{0} & \dotsb & s_{2n}^{2n-1} & s_{2n}^{2n} \\
z^{0} & z^{1} & \dotsb & z^{2n} & z^{2n+1}
\end{vmatrix} \\
=& \, -\dfrac{z^{-n-1}}{H^{(-2n)}_{2n+1}} \, \underbrace{\int_{\mathbb{R}}
\int_{\mathbb{R}} \dotsi \int_{\mathbb{R}}}_{2n+1} \, \, \md \widetilde{\mu}
(s_{0}) \md \widetilde{\mu}(s_{1}) \, \dotsb \, \md \widetilde{\mu}(s_{2n})
s_{0}^{-2n-1}s_{1}^{-2n} \dotsb s_{2n}^{-1} \\
\times& \,
\begin{vmatrix}
s_{0}^{0} & s_{0}^{1} & \dotsb & s_{0}^{2n} & s_{0}^{2n+1} \\
s_{1}^{0} & s_{1}^{1} & \dotsb & s_{1}^{2n} & s_{1}^{2n+1} \\
\vdots & \vdots & \ddots & \vdots & \vdots \\
s_{2n}^{0} & s_{2n}^{1} & \dotsb & s_{2n}^{2n} & s_{2n}^{2n+1} \\
z^{0} & z^{1} & \dotsb & z^{2n} & z^{2n+1}
\end{vmatrix} \\
=& \, -\dfrac{z^{-n-1}}{H^{(-2n)}_{2n+1}(2n \! + \! 1)!} \, \sum_{\boldsymbol{
\sigma} \in \mathfrak{S}_{2n+1}} \, \underbrace{\int_{\mathbb{R}} \int_{
\mathbb{R}} \dotsi \int_{\mathbb{R}}}_{2n+1} \, \, \md \widetilde{\mu}(s_{
\sigma (0)}) \md \widetilde{\mu}(s_{\sigma (1)}) \, \dotsb \, \md \widetilde{
\mu}(s_{\sigma (2n)}) \\
\times& \, s_{\sigma (0)}^{-2n-1}s_{\sigma (1)}^{-2n-1} \dotsb s_{\sigma (2n)
}^{-2n-1}s_{\sigma (0)}^{0}s_{\sigma (1)}^{1} \dotsb s_{\sigma (2n)}^{2n}
\begin{vmatrix}
s_{\sigma (0)}^{0} & s_{\sigma (0)}^{1} & \dotsb & s_{\sigma (0)}^{2n} &
s_{\sigma (0)}^{2n+1} \\
s_{\sigma (1)}^{0} & s_{\sigma (1)}^{1} & \dotsb & s_{\sigma (1)}^{2n} &
s_{\sigma (1)}^{2n+1} \\
\vdots & \vdots & \ddots & \vdots & \vdots \\
s_{\sigma (2n)}^{0} & s_{\sigma (2n)}^{1} & \dotsb & s_{\sigma (2n)}^{2n} &
s_{\sigma (2n)}^{2n+1} \\
z^{0} & z^{1} & \dotsb & z^{2n} & z^{2n+1}
\end{vmatrix} \\
=& \, -\dfrac{z^{-n-1}}{H^{(-2n)}_{2n+1}(2n \! + \! 1)!} \, \underbrace{\int_{
\mathbb{R}} \int_{\mathbb{R}} \dotsi \int_{\mathbb{R}}}_{2n+1} \, \, \md
\widetilde{\mu}(s_{0}) \md \widetilde{\mu}(s_{1}) \, \dotsb \, \md \widetilde{
\mu}(s_{2n})s_{0}^{-2n-1}s_{1}^{-2n-1} \dotsb s_{2n}^{-2n-1} \\
\times& \, \left(\sum_{\boldsymbol{\sigma} \in \mathfrak{S}_{2n+1}}
\operatorname{sgn}(\boldsymbol{\sigma})s_{\sigma (0)}^{0}s_{\sigma (1)}^{1}
\dotsb s_{\sigma (2n)}^{2n} \right) \!
\begin{vmatrix}
s_{0}^{0} & s_{0}^{1} & \dotsb & s_{0}^{2n} & s_{0}^{2n+1} \\
s_{1}^{0} & s_{1}^{1} & \dotsb & s_{1}^{2n} & s_{1}^{2n+1} \\
\vdots & \vdots & \ddots & \vdots & \vdots \\
s_{2n}^{0} & s_{2n}^{1} & \dotsb & s_{2n}^{2n} & s_{2n}^{2n+1} \\
z^{0} & z^{1} & \dotsb & z^{2n} & z^{2n+1}
\end{vmatrix} \\
=& \, -\dfrac{z^{-n-1}}{H^{(-2n)}_{2n+1}(2n \! + \! 1)!} \, \underbrace{\int_{
\mathbb{R}} \int_{\mathbb{R}} \dotsi \int_{\mathbb{R}}}_{2n+1} \, \, \md
\widetilde{\mu}(s_{0}) \md \widetilde{\mu}(s_{1}) \, \dotsb \, \md \widetilde{
\mu}(s_{2n})s_{0}^{-2n-1}s_{1}^{-2n-1} \dotsb s_{2n}^{-2n-1} \\
\times& \,
\begin{vmatrix}
s_{0}^{0} & s_{0}^{1} & \dotsb & s_{0}^{2n} \\
s_{1}^{0} & s_{1}^{1} & \dotsb & s_{1}^{2n} \\
\vdots & \vdots & \ddots & \vdots \\
s_{2n}^{0} & s_{2n}^{1} & \dotsb & s_{2n}^{2n}
\end{vmatrix}
\begin{vmatrix}
s_{0}^{0} & s_{0}^{1} & \dotsb & s_{0}^{2n} & s_{0}^{2n+1} \\
s_{1}^{0} & s_{1}^{1} & \dotsb & s_{1}^{2n} & s_{1}^{2n+1} \\
\vdots & \vdots & \ddots & \vdots & \vdots \\
s_{2n}^{0} & s_{2n}^{1} & \dotsb & s_{2n}^{2n} & s_{2n}^{2n+1} \\
z^{0} & z^{1} & \dotsb & z^{2n} & z^{2n+1}
\end{vmatrix};
\end{align*}
but a straightforward calculation shows that
\begin{equation*}
\begin{vmatrix}
s_{0}^{0} & s_{0}^{1} & \dotsb & s_{0}^{2n} & s_{0}^{2n+1} \\
s_{1}^{0} & s_{1}^{1} & \dotsb & s_{1}^{2n} & s_{1}^{2n+1} \\
\vdots & \vdots & \ddots & \vdots & \vdots \\
s_{2n}^{0} & s_{2n}^{1} & \dotsb & s_{2n}^{2n} & s_{2n}^{2n+1} \\
z^{0} & z^{1} & \dotsb & z^{2n} & z^{2n+1}
\end{vmatrix} \! = \! 
\begin{vmatrix}
s_{0}^{0} & s_{0}^{1} & \dotsb & s_{0}^{2n} \\
s_{1}^{0} & s_{1}^{1} & \dotsb & s_{1}^{2n} \\
\vdots & \vdots & \ddots & \vdots \\
s_{2n}^{0} & s_{2n}^{1} & \dotsb & s_{2n}^{2n}
\end{vmatrix} \prod_{j=0}^{2n}(z \! - \! s_{j}),
\end{equation*}
whence
\begin{align*}
\boldsymbol{\pi}_{2n}(z) =& \, -\dfrac{z^{-n-1}}{H^{(-2n)}_{2n+1}(2n \! + \!
1)!} \, \underbrace{\int_{\mathbb{R}} \int_{\mathbb{R}} \dotsi \int_{\mathbb{
R}}}_{2n+1} \, \, \md \widetilde{\mu}(s_{0}) \md \widetilde{\mu}(s_{1}) \,
\dotsb \, \md \widetilde{\mu}(s_{2n})s_{0}^{-2n-1}s_{1}^{-2n-1} \dotsb
s_{2n}^{-2n-1} \\
\times& \, \prod_{j=0}^{2n}(z \! - \! s_{j})
\begin{vmatrix}
s_{0}^{0} & s_{0}^{1} & \dotsb & s_{0}^{2n} \\
s_{1}^{0} & s_{1}^{1} & \dotsb & s_{1}^{2n} \\
\vdots & \vdots & \ddots & \vdots \\
s_{2n}^{0} & s_{2n}^{1} & \dotsb & s_{2n}^{2n}
\end{vmatrix}^{2} \\
=& \, -\dfrac{z^{-n-1}}{H^{(-2n)}_{2n+1}(2n \! + \! 1)!} \, \underbrace{\int_{
\mathbb{R}} \int_{\mathbb{R}} \dotsi \int_{\mathbb{R}}}_{2n+1} \, \, \md
\widetilde{\mu}(s_{0}) \md \widetilde{\mu}(s_{1}) \, \dotsb \, \md \widetilde{
\mu}(s_{2n})s_{0}^{-2n-1}s_{1}^{-2n-1} \dotsb s_{2n}^{-2n-1} \\
\times& \, \prod_{j=0}^{2n}(z \! - \! s_{j})
\underbrace{\begin{vmatrix}
1 & 1 & \dotsb & 1 \\
s_{0}^{1} & s_{1}^{1} & \dotsb & s_{2n}^{1} \\
\vdots & \vdots & \ddots & \vdots \\
s_{0}^{2n} & s_{1}^{2n} & \dotsb & s_{2n}^{2n}
\end{vmatrix}^{2}}_{= \, \prod_{\underset{l<i}{i,l=0}}^{2n}(s_{i}-s_{l})^{2}};
\end{align*}
hence the integral representation for $\boldsymbol{\pi}_{2n+1}(z)$ stated 
in the Proposition, with the integral representation for $H^{(-2n)}_{2n+1}$ 
derived in the proof of Lemma~2.2.2. See \cite{a38}, Proposition~2.2.1, for 
the proof of the integral representation for the even degree monic OLP 
$\boldsymbol{\pi}_{2n}(z)$. \hfill $\qed$
\begin{eeee}
For the purposes of the ensuing asymptotic analysis, it is convenient to 
re-write $\md \widetilde{\mu}(z) \! = \! \exp (-\mathscr{N} \, V(z)) \, \md 
z \! = \! \exp (-n \widetilde{V}(z)) \, \md z \! =: \! \md \mu (z)$, $n \! 
\in \! \mathbb{N}$, where
\begin{equation*}
\widetilde{V}(z) \! = \! z_{o}V(z),
\end{equation*}
with
\begin{equation*}
z_{o} \colon \mathbb{N} \times \mathbb{N} \! \to \! \mathbb{R}_{+}, \,
(\mathscr{N},n) \! \mapsto \! z_{o} \! := \mathscr{N}/n,
\end{equation*}
and where the `scaled' external field $\widetilde{V} \colon \mathbb{R}
\setminus \{0\} \! \to \! \mathbb{R}$ satisfies the following conditions:
\begin{gather}
\widetilde{V} \, \, \text{is real analytic on} \, \, \mathbb{R} \setminus
\{0\}; \\
\lim_{\vert x \vert \to \infty} \! \left(\widetilde{V}(x)/\ln (x^{2} \! + \!
1) \right) \! = \! +\infty; \\
\lim_{\vert x \vert \to 0} \! \left(\widetilde{V}(x)/\ln (x^{-2} \! + \! 1)
\right) \! = \! +\infty.
\end{gather}
(For example, a rational function of the form $\widetilde{V}(z) \! = \! \sum_{
k=-2m_{1}}^{2m_{2}} \widetilde{\varrho}_{k}z^{k}$, with $\widetilde{\varrho}_{
k} \! \in \! \mathbb{R}$, $k \! = \! -2m_{1},\dotsc,2m_{2}$, $m_{1,2} \! \in
\! \mathbb{N}$, and $\widetilde{\varrho}_{-2m_{1}},\widetilde{\varrho}_{2
m_{2}} \! > \! 0$ would satisfy conditions~(2.3)--(2.5).)

Hereafter, the double-scaling limit as $\mathscr{N},n \! \to \! \infty$ such 
that $z_{o} \! = \! 1 \! + \! o(1)$ is studied (the simplified `notation' $n 
\! \to \! \infty$ will be adopted). \hfill $\blacksquare$
\end{eeee}

It is, by now, a well-known, if not established, mathematical fact that 
variational conditions for minimisation problems in logarithmic potential 
theory, via the \emph{equilibrium measure} \cite{a43,a44,a79,a80,a81}, play 
a crucial r\^{o}le in the asymptotic analysis of (matrix) RHPs associated 
with (continuous and discrete) orthogonal polynomials, their roots, and 
corresponding recurrence relation coefficients (see, for example, 
\cite{a46,a47,a49,a53,a62}). The situation with respect to the large-$n$ 
asymptotic analysis for the monic OLPs, $\boldsymbol{\pi}_{n}(z)$, is 
analogous; but, unlike the asymptotic analysis for the orthogonal polynomials 
case, the asymptotic analysis for $\boldsymbol{\pi}_{n}(z)$ requires the 
consideration of two different families of RHPs, one for even degree 
(\textbf{RHP1}) and one for odd degree (\textbf{RHP2}). Thus, one must 
consider two sets of variational conditions for two (suitably posed) 
minimisation problems.

The following discussion is decomposed into two parts: one part corresponding 
to the RHP for $\overset{e}{\mathrm{Y}} \colon \mathbb{C} \setminus \mathbb{R}
\! \to \! \operatorname{SL}_{2}(\mathbb{C})$ formulated as \textbf{RHP1},
denoted by $\pmb{\mathrm{P}_{1}}$, and the other part corresponding to the RHP
for $\overset{o}{\mathrm{Y}} \colon \mathbb{C} \setminus \mathbb{R} \! \to \!
\operatorname{SL}_{2}(\mathbb{C})$ formulated as \textbf{RHP2}, denoted by
$\pmb{\mathrm{P}_{2}}$.
\begin{compactenum}
\item[\shadowbox{$\pmb{\mathrm{P}_{1}}$}] Let $\widetilde{V} \colon \mathbb{R}
\setminus \{0\} \! \to \! \mathbb{R}$ satisfy conditions~(2.3)--(2.5). Let
$\mathrm{I}_{V}^{e}[\mu^{e}] \colon \mathcal{M}_{1}(\mathbb{R}) \! \to \!
\mathbb{R}$ denote the functional
\begin{equation*}
\mathrm{I}_{V}^{e}[\mu^{e}] \! = \! \iint_{\mathbb{R}^{2}} \ln \! \left(
\dfrac{\lvert st \rvert}{\lvert s \! - \! t \rvert^{2}} \right) \md \mu^{e}(s)
\, \md \mu^{e}(t) \! + \! 2 \int_{\mathbb{R}} \widetilde{V}(s) \, \md \mu^{e}
(s),
\end{equation*}
and consider the associated minimisation problem,
\begin{equation*}
E_{V}^{e} \! = \! \inf \! \left\lbrace \mathstrut \mathrm{I}_{V}^{e}[\mu^{e}]; 
\, \mu^{e} \! \in \! \mathcal{M}_{1}(\mathbb{R}) \right\rbrace.
\end{equation*}
The infimum is finite, and there exists a unique measure $\mu_{V}^{e}$, 
referred to as the `even' equilibrium measure, achieving the infimum (that 
is, $\mathcal{M}_{1}(\mathbb{R}) \! \ni \! \mu_{V}^{e} \! = \! \inf \lbrace 
\mathstrut \mathrm{I}_{V}^{e}[\mu^{e}]; \, \mu^{e} \! \in \! \mathcal{M}_{1}
(\mathbb{R}) \rbrace)$. Furthermore, $\mu_{V}^{e}$ has the following 
`regularity' properties (see \cite{a38} for complete details and proofs):
\begin{compactenum}
\item[\textbullet] the `even' equilibrium measure has compact support which
consists of the disjoint union of a finite number of bounded real intervals;
in fact, as shown in \cite{a38}, $\mathrm{supp}(\mu_{V}^{e}) \! =: \!
J_{e}$\footnote{It would be more usual, {}from the outset, for the bounded
(and closed) set $\overline{J_{e}} \! := \! \cup_{j=1}^{N+1}[b_{j-1}^{e},
a_{j}^{e}]$ to denote the support of $\mu_{V}^{e}$; however, the open (and
bounded) set $J_{e}$ provides an effective description of (the interior of)
the support of $\mu_{V}^{e}$: for this reason, $J_{e}$ (and at other times
$\overline{J_{e}})$ is used to denote $\operatorname{supp}(\mu_{V}^{e})$;
\emph{mutatis mutandis} for $J_{o}$ and $\overline{J_{o}}$ (see
\shadowbox{$\mathrm{P}_{2}$} below). This should not cause confusion for 
the reader.} $\! = \! \cup_{j=1}^{N+1}(b_{j-1}^{e},a_{j}^{e})$ $(\subset \!
\mathbb{R} \setminus \{0\})$, where $\lbrace b_{j-1}^{e},a_{j}^{e} \rbrace_{j
=1}^{N+1}$, with $b_{0}^{e} \! := \! \min \lbrace \mathrm{supp}(\mu_{V}^{e})
\rbrace \! \notin \! \lbrace -\infty,0 \rbrace$, $a_{N+1}^{e} \! := \! \max
\lbrace \mathrm{supp}(\mu_{V}^{e}) \rbrace \! \notin \! \lbrace 0,+\infty
\rbrace$, and $-\infty \! < \! b_{0}^{e} \! < \! a_{1}^{e} \! < \! b_{1}^{e}
\! < \! a_{2}^{e} \! < \! \cdots \! < \! b_{N}^{e} \! < \! a_{N+1}^{e} \! <
\! +\infty$, constitute the end-points of the support of $\mu_{V}^{e}$;
\item[\textbullet] the end-points $\lbrace b_{j-1}^{e},a_{j}^{e} \rbrace_{j
=1}^{N+1}$ are not arbitrary; rather, they satisfy an $n$-dependent and
(locally) solvable system of $2(N \! + \! 1)$ \emph{moment conditions}
(transcendental equations; see \cite{a38}, Lemma~3.5);
\item[\textbullet] the `even' equilibrium measure is absolutely continuous
with respect to Lebesgue measure. The \emph{density} is given by
\begin{equation*}
\md \mu_{V}^{e}(x) \! := \! \psi_{V}^{e}(x) \, \md x \! = \! \dfrac{1}{2 \pi
\mi}(R_{e}(x))^{1/2}_{+}h_{V}^{e}(x) \pmb{1}_{J_{e}}(x) \, \md x,
\end{equation*}
where
\begin{equation*}
(R_{e}(z))^{1/2} \! := \! \left(\prod_{j=1}^{N+1}(z \! - \! b_{j-1}^{e})(z \!
- \! a_{j}^{e}) \right)^{1/2},
\end{equation*}
with $(R_{e}(x))^{1/2}_{\pm} \! := \! \lim_{\varepsilon \downarrow 0}(R_{e}
(x \! \pm \! \mi \varepsilon))^{1/2}$ and the branch of the square root is
chosen, as per the discussion in Subsection~2.1, such that $z^{-(N+1)}(R_{e}
(z))^{1/2} \! \sim_{\underset{z \in \mathbb{C}_{\pm}}{z \to \infty}} \! \pm
1$, $h_{V}^{e}(z) \! := \! \tfrac{1}{2} \oint_{C_{\mathrm{R}}^{e}}(R_{e}(s)
)^{-1/2}(\tfrac{\mi}{\pi s} \! + \! \tfrac{\mi \widetilde{V}^{\prime}(s)}{2
\pi})(s \! - \! z)^{-1} \, \md s$ (real analytic for $z \! \in \! \mathbb{R}
\setminus \{0\})$, where ${}^{\prime}$ denotes differentiation with respect
to the argument, $C_{\mathrm{R}}^{e}$ $(\subset \mathbb{C}^{\ast})$ is the
union of two circular contours, one outer one of large radius $R^{\natural}$
traversed clockwise and one inner one of small radius $r^{\natural}$ traversed
counter-clockwise, with the numbers $0 \! < \! r^{\natural} \! < \! R^{
\natural} \! < \! +\infty$ chosen such that, for (any) non-real $z$ in the
domain of analyticity of $\widetilde{V}$ (that is, $\mathbb{C}^{\ast})$,
$\mathrm{int}(C_{\mathrm{R}}^{e}) \! \supset \! J_{e} \cup \{z\}$, and $\pmb{
1}_{J_{e}}(x)$ denotes the indicator (characteristic) function of the set
$J_{e}$. (Note that $\psi_{V}^{e}(x) \! \geqslant \! 0 \, \, \forall \, \, x
\! \in \! \overline{J_{e}} \! := \! \cup_{j=1}^{N+1}[b_{j-1}^{e},a_{j}^{e}]$:
it vanishes like a square root at the end-points of the support of the `even'
equilibrium measure, that is, $\psi_{V}^{e}(s) \! =_{s \downarrow b_{j-1}^{e}}
\! \mathcal{O}((s \! - \! b_{j-1}^{e})^{1/2})$ and $\psi_{V}^{e}(s) \! =_{s
\uparrow a_{j}^{e}} \! \mathcal{O}((a_{j}^{e} \! - \! s)^{1/2})$, $j \! = \!
1,\dotsc,N \! + \! 1$.);
\item[\textbullet] the `even' equilibrium measure and its (compact) support
are uniquely characterised by the following Euler-Lagrange variational
equations: there exists $\ell_{e} \! \in \! \mathbb{R}$, the `even' Lagrange
multiplier, and $\mu^{e} \! \in \! \mathcal{M}_{1}(\mathbb{R})$ such that
\begin{gather*}
4 \int_{J_{e}} \ln (\vert x \! - \! s \vert) \, \md \mu^{e} (s) \! - \! 2 \ln
\vert x \vert \! - \! \widetilde{V}(x) \! - \! \ell_{e} \! = \! 0, \quad x \!
\in \! \overline{J_{e}}, \tag{$\mathrm{P}_{1}^{(a)}$} \\
4 \int_{J_{e}} \ln (\vert x \! - \! s \vert) \, \md \mu^{e} (s) \! - \! 2 \ln
\vert x \vert \! - \! \widetilde{V}(x) \! - \! \ell_{e} \! \leqslant \! 0,
\quad x \! \in \! \mathbb{R} \setminus \overline{J_{e}}; \tag{$\mathrm{P}_{
1}^{(b)}$}
\end{gather*}
\item[\textbullet] the Euler-Lagrange variational equations can be 
conveniently recast in terms of the complex potential $g^{e}(z)$ of $\mu_{
V}^{e}$:
\begin{equation*}
g^{e}(z) \! := \! \int_{J_{e}} \! \ln \! \left((z \! - \! s)^{2}(zs)^{-1}
\right) \md \mu_{V}^{e}(s), \quad z \! \in \! \mathbb{C} \setminus (-\infty,
\max \{0,a_{N+1}^{e}\}).
\end{equation*}
The function $g^{e} \colon \mathbb{C} \setminus (-\infty,\max \{0,a_{N+1}^{e}
\}) \! \to \! \mathbb{C}$ so defined satisfies:
\begin{compactenum}
\item[$(\mathrm{P}_{1}^{(1)})$] $g^{e}(z)$ is analytic for $z \! \in \!
\mathbb{C} \setminus (-\infty,\max \{0,a_{N+1}^{e}\})$;
\item[$(\mathrm{P}_{1}^{(2)})$] $g^{e}(z) \! =_{\underset{z \in \mathbb{C}
\setminus \mathbb{R}}{z \to \infty}} \! \ln (z) \! + \! \mathcal{O}(1)$;
\item[$(\mathrm{P}_{1}^{(3)})$] $g^{e}_{+}(z) \! + \! g^{e}_{-}(z) \! -
\! \widetilde{V}(z) \! - \! \ell_{e} \! + \! 2Q_{e} \! = \! 0$, $z \! \in \!
\overline{J_{e}}$, where $g^{e}_{\pm}(z) \! := \! \lim_{\varepsilon \downarrow
0}g^{e}(z \! \pm \! \mi \varepsilon)$, and $Q_{e} \! := \! \int_{J_{e}} \ln
(s) \, \md \mu_{V}^{e}(s) \! = \! \int_{J_{e}} \ln (\lvert s \rvert) \, \md
\mu_{V}^{e}(s) \! + \! \mi \pi \int_{J_{e} \cap \mathbb{R}_{-}} \md \mu_{V}^{
e}(s)$;
\item[$(\mathrm{P}_{1}^{(4)})$] $g^{e}_{+}(z) \! + \! g^{e}_{-}(z) \! - \!
\widetilde{V}(z) \! - \! \ell_{e} \! + \! 2Q_{e} \! \leqslant \! 0$, $z \! \in
\! \mathbb{R} \setminus \overline{J_{e}}$, where equality holds for at most a
finite number of points;
\item[$(\mathrm{P}_{1}^{(5)})$] $g^{e}_{+}(z) \! - \! g^{e}_{-}(z) \! = \!
\mi f_{g^{e}}^{\mathbb{R}}(z)$, $z \! \in \! \mathbb{R}$, where $f_{g^{e}}^{
\mathbb{R}} \colon \mathbb{R} \! \to \! \mathbb{R}$, and, in particular, $g^{
e}_{+}(z) \! - \! g^{e}_{-}(z) \! = \! \mi \operatorname{const.}$, $z \! \in
\! \mathbb{R} \setminus \overline{J_{e}}$, with $\operatorname{const.} \! \in
\! \mathbb{R}$;
\item[$(\mathrm{P}_{1}^{(6)})$] $\mi (g^{e}_{+}(z) \! - \! g^{e}_{-}(z))^{
\prime} \! \geqslant \! 0$, $z \! \in \! J_{e}$, where equality holds for 
at most a finite number of points.
\end{compactenum}
\end{compactenum}
\item[\shadowbox{$\pmb{\mathrm{P}_{2}}$}] Let $\widetilde{V} \colon \mathbb{
R} \setminus \{0\} \! \to \! \mathbb{R}$ satisfy conditions~(2.3)--(2.5). Let
$\mathrm{I}_{V}^{o}[\mu^{o}] \colon \mathcal{M}_{1}(\mathbb{R}) \! \to \!
\mathbb{R}$ denote the functional
\begin{equation*}
\mathrm{I}_{V}^{o}[\mu^{o}] \! = \! \iint_{\mathbb{R}^{2}} \ln \! \left(\dfrac{
\lvert st \rvert}{\lvert s \! - \! t \rvert^{2+\frac{1}{n}}} \right) \md \mu^{
o}(s) \, \md \mu^{o}(t) \! + \! 2 \int_{\mathbb{R}} \widetilde{V}(s) \, \md
\mu^{o}(s), \quad n \! \in \! \mathbb{N},
\end{equation*}
and consider the associated minimisation problem,
\begin{equation*}
E_{V}^{o} \! = \! \inf \! \left\lbrace \mathstrut \mathrm{I}_{V}^{o}[\mu^{o}]; 
\, \mu^{o} \! \in \! \mathcal{M}_{1}(\mathbb{R}) \right\rbrace.
\end{equation*}
The infimum is finite, and there exists a unique measure $\mu_{V}^{o}$, 
referred to as the `odd' equilibrium measure, achieving the infimum (that 
is, $\mathcal{M}_{1}(\mathbb{R}) \! \ni \! \mu_{V}^{o} \! = \! \inf \lbrace 
\mathstrut \mathrm{I}_{V}^{o}[\mu^{o}]; \, \mu^{o} \! \in \! \mathcal{M}_{
1}(\mathbb{R}) \rbrace)$. Furthermore, $\mu_{V}^{o}$ has the following 
`regularity' properties (all of these properties are proven in this work):
\begin{compactenum}
\item[\textbullet] the `odd' equilibrium measure has compact support which
consists of the disjoint union of a finite number of bounded real intervals;
in fact, as shown in Section~3 (see Lemma~3.5), $\mathrm{supp}(\mu_{V}^{o})
\! =: \! J_{o} \! = \! \cup_{j=1}^{N+1}(b_{j-1}^{o},a_{j}^{o})$ $(\subset \!
\mathbb{R} \setminus \{0\})$, where $\lbrace b_{j-1}^{o},a_{j}^{o} \rbrace_{j
=1}^{N+1}$, with $b_{0}^{o} \! := \! \min \lbrace \mathrm{supp}(\mu_{V}^{o}) 
\rbrace \! \notin \! \lbrace -\infty,0 \rbrace$, $a_{N+1}^{o} \! := \! \max 
\lbrace \mathrm{supp}(\mu_{V}^{o}) \rbrace \! \notin \! \lbrace 0,+\infty 
\rbrace$, and $-\infty \! < \! b_{0}^{o} \! < \! a_{1}^{o} \! < \! b_{1}^{o} 
\! < \! a_{2}^{o} \! < \! \cdots \! < \! b_{N}^{o} \! < \! a_{N+1}^{o} \! < 
\! +\infty$, constitute the end-points of the support of $\mu_{V}^{o}$; (The 
number of intervals, $N \! + \! 1$, is the same in the `odd' case as in the 
`even' case, which can be established by a lengthy analysis similar to that 
contained in \cite{a81}.)
\item[\textbullet] the end-points $\lbrace b_{j-1}^{o},a_{j}^{o} \rbrace_{j
=1}^{N+1}$ are not arbitrary; rather, they satisfy the $n$-dependent and 
(locally) solvable system of $2(N \! + \! 1)$ moment conditions 
(transcendental equations) given in Lemma~3.5;
\item[\textbullet] the `odd' equilibrium measure is absolutely continuous 
with respect to Lebesgue measure. The density is given by
\begin{equation*}
\md \mu_{V}^{o}(x) \! := \! \psi_{V}^{o}(x) \, \md x \! = \! \dfrac{1}{2 \pi
\mi}(R_{o}(x))^{1/2}_{+}h_{V}^{o}(x) \pmb{1}_{J_{o}}(x) \, \md x,
\end{equation*}
where
\begin{equation*}
(R_{o}(z))^{1/2} \! := \! \left(\prod_{j=1}^{N+1}(z \! - \! b_{j-1}^{o})(z \!
- \! a_{j}^{o}) \right)^{1/2},
\end{equation*}
with $(R_{o}(x))^{1/2}_{\pm} \! := \! \lim_{\varepsilon \downarrow 0}(R_{o}
(x \! \pm \! \mi \varepsilon))^{1/2}$ and the branch of the square root is
chosen, as per the discussion in Subsection~2.1, such that $z^{-(N+1)}(R_{o}
(z))^{1/2} \! \sim_{\underset{z \in \mathbb{C}_{\pm}}{z \to \infty}} \! \pm
1$, $h_{V}^{o}(z) \! := \! (2 \! + \! \tfrac{1}{n})^{-1} \oint_{C_{\mathrm{
R}}^{o}}(R_{o}(s))^{-1/2}(\tfrac{\mi}{\pi s} \! + \! \frac{\mi \widetilde{
V}^{\prime}(s)}{2 \pi})(s \! - \! z)^{-1} \, \md s$ (real analytic for $z \!
\in \! \mathbb{R} \setminus \{0\})$, where $C_{\mathrm{R}}^{o}$ $(\subset
\mathbb{C}^{\ast})$ is the union of two circular contours, one outer one of
large radius $R^{\flat}$ traversed clockwise and one inner one of small radius
$r^{\flat}$ traversed counter-clockwise, with the numbers $0 \! < \! r^{\flat}
\! < \! R^{\flat} \! < \! +\infty$ chosen such that, for (any) non-real $z$ in
the domain of analyticity of $\widetilde{V}$ (that is, $\mathbb{C}^{\ast})$,
$\mathrm{int}(C_{\mathrm{R}}^{o}) \! \supset \! J_{o} \cup \{z\}$, and $\pmb{
1}_{J_{o}}(x)$ denotes the indicator (characteristic) function of the set
$J_{o}$. (Note that $\psi_{V}^{o}(x) \! \geqslant \! 0 \, \, \forall \, \, x
\! \in \! \overline{J_{o}} \! := \! \cup_{j=1}^{N+1}[b_{j-1}^{o},a_{j}^{o}]$:
it vanishes like a square root at the end-points of the support of the `odd'
equilibrium measure, that is, $\psi_{V}^{o}(s) \! =_{s \downarrow b_{j-1}^{
o}} \! \mathcal{O}((s \! - \! b_{j-1}^{o})^{1/2})$ and $\psi_{V}^{o}(s) \!
=_{s \uparrow a_{j}^{o}} \! \mathcal{O}((a_{j}^{o} \! - \! s)^{1/2})$, $j \!
= \! 1,\dotsc,N \! + \! 1$.);
\item[\textbullet] the `odd' equilibrium measure and its (compact) support are
uniquely characterised by the following Euler-Lagrange variational equations:
there exists $\ell_{o} \! \in \! \mathbb{R}$, the `odd' Lagrange multiplier,
and $\mu^{o} \! \in \! \mathcal{M}_{1}(\mathbb{R})$ such that
\begin{gather*}
2 \! \left(2 \! + \! \dfrac{1}{n} \right) \! \int_{J_{o}} \ln (\vert x \! - \!
s \vert) \, \md \mu^{o} (s) \! - \! 2 \ln \vert x \vert \! - \! \widetilde{V}
(x) \! - \! \ell_{o} \! - \! 2 \left(2 \! + \! \dfrac{1}{n} \right) \!
\widetilde{Q}_{o} = \! 0, \quad x \! \in \! \overline{J_{o}},
\tag{$\mathrm{P}_{2}^{(a)}$} \\
2 \! \left(2 \! + \! \dfrac{1}{n} \right) \! \int_{J_{o}} \ln (\vert x \! - \!
s \vert) \, \md \mu^{o} (s) \! - \! 2 \ln \vert x \vert \! - \! \widetilde{V}
(x) \! - \! \ell_{o} \! - \! 2 \left(2 \! + \! \dfrac{1}{n} \right) \!
\widetilde{Q}_{o} \! \leqslant \! 0, \quad x \! \in \! \mathbb{R} \setminus
\overline{J_{o}}, \tag{$\mathrm{P}_{2}^{(b)}$}
\end{gather*}
where $\widetilde{Q}_{o} \! := \! \int_{J_{o}} \ln (\lvert s \rvert) \, \md
\mu^{o}(s)$;
\item[\textbullet] the Euler-Lagrange variational equations can be 
conveniently recast in terms of the complex potential $g^{o}(z)$ of $\mu_{
V}^{o}$:
\begin{equation*}
g^{o}(z) \! := \! \int_{J_{o}} \ln \! \left((z \! - \! s)^{2+\frac{1}{n}}
(zs)^{-1} \right) \md \mu_{V}^{o}(s), \quad z \! \in \! \mathbb{C} \setminus
(-\infty,\max \{0,a_{N+1}^{o}\}).
\end{equation*}
The function $g^{o} \colon \mathbb{C} \setminus (-\infty,\max \{0,a_{N+1}^{o}
\}) \! \to \! \mathbb{C}$ so defined satisfies:
\begin{compactenum}
\item[$(\mathrm{P}_{2}^{(1)})$] $g^{o}(z)$ is analytic for $z \! \in \! 
\mathbb{C} \setminus (-\infty,\max \{0,a_{N+1}^{o}\})$;
\item[$(\mathrm{P}_{2}^{(2)})$] $g^{o}(z) \! =_{\underset{z \in \mathbb{C}
\setminus \mathbb{R}}{z \to 0}} \! -\ln (z) \! + \! \mathcal{O}(1)$;
\item[$(\mathrm{P}_{2}^{(3)})$] $g^{o}_{+}(z) \! + \! g^{o}_{-}(z) \! - \!
\widetilde{V}(z) \! - \! \ell_{o} \! - \! \mathfrak{Q}_{\mathscr{A}}^{+} \! -
\! \mathfrak{Q}_{\mathscr{A}}^{-} \! = \! 0$, $z \! \in \! \overline{J_{o}}$,
where $g^{o}_{\pm}(z) \! := \! \lim_{\varepsilon \downarrow 0}g^{o}(z \! \pm
\! \mi \varepsilon)$, and $\mathfrak{Q}_{\mathscr{A}}^{\pm} \! := \! (1 \! +
\! \tfrac{1}{n}) \int_{J_{o}} \ln (\lvert s \rvert) \, \md \mu_{V}^{o}(s) \!
- \! \mi \pi \int_{J_{o} \cap \mathbb{R}_{-}} \md \mu_{V}^{o}(s) \! \pm \! \mi
\pi (2 \! + \! \tfrac{1}{n}) \int_{J_{o} \cap \mathbb{R}_{+}} \md \mu_{V}^{o}
(s)$;
\item[$(\mathrm{P}_{2}^{(4)})$] $g^{o}_{+}(z) \! + \! g^{o}_{-}(z) \! - \!
\widetilde{V}(z) \! - \! \ell_{o} \! - \! \mathfrak{Q}_{\mathscr{A}}^{+} \! -
\! \mathfrak{Q}_{\mathscr{A}}^{-} \! \leqslant \! 0$, $z \! \in \! \mathbb{R}
\setminus \overline{J_{o}}$, where equality holds for at most a finite number
of points;
\item[$(\mathrm{P}_{2}^{(5)})$] $g^{o}_{+}(z) \! - \! g^{o}_{-}(z) \! - \!
\mathfrak{Q}_{\mathscr{A}}^{+} \! + \! \mathfrak{Q}_{\mathscr{A}}^{-} \! = \!
\mi f_{g^{o}}^{\mathbb{R}}(z)$, $z \! \in \! \mathbb{R}$, where $f_{g^{o}}^{
\mathbb{R}} \colon \mathbb{R} \! \to \! \mathbb{R}$, and, in particular, $g^{
o}_{+}(z) \! - \! g^{o}_{-}(z) \! - \! \mathfrak{Q}_{\mathscr{A}}^{+} \! + \!
\mathfrak{Q}_{\mathscr{A}}^{-} \! = \! \mi \operatorname{const.}$, $z \! \in
\! \mathbb{R} \setminus \overline{J_{o}}$, with $\operatorname{const.} \! \in
\! \mathbb{R}$;
\item[$(\mathrm{P}_{2}^{(6)})$] $\mi (g^{o}_{+}(z) \! - \! g^{o}_{-}(z) \! -
\! \mathfrak{Q}_{\mathscr{A}}^{+} \! + \! \mathfrak{Q}_{\mathscr{A}}^{-} )^{
\prime} \! \geqslant \! 0$, $z \! \in \! J_{o}$, where equality holds for at
most a finite number of points.
\end{compactenum}
\end{compactenum}
\end{compactenum}

In this three-fold series of works on asymptotics of OLPs and related 
quantities, the so-called `regular case' is studied, namely:
\begin{compactenum}
\item[\textbullet] $\md \mu_{V}^{e}$, or $\widetilde{V} \colon \mathbb{R} 
\setminus \{0\} \! \to \! \mathbb{R}$ satisfying conditions~(2.3)--(2.5), 
is \emph{regular} if: (i) $h_{V}^{e}(x) \! \not\equiv \! 0$ on $\overline{
J_{e}}$; (ii) $4 \int_{J_{e}} \ln (\vert x \! - \! s \vert) \, \md \mu_{V}^{e}
(s) \! - \! 2 \ln \vert x \vert \! - \! \widetilde{V}(x) \! - \! \ell_{e} 
\! < \! 0$, $x \! \in \! \mathbb{R} \setminus \overline{J_{e}}$; and (iii) 
inequalities~$(\mathrm{P}_{1}^{(4)})$ and~$(\mathrm{P}_{1}^{(6)})$ in $\pmb{
\mathrm{P}_{1}}$ are strict, that is, $\leqslant$ (resp., $\geqslant)$ is 
replaced by $<$ (resp., $>)$;
\item[\textbullet] $\md \mu_{V}^{o}$, or $\widetilde{V} \colon \mathbb{R} 
\setminus \{0\} \! \to \! \mathbb{R}$ satisfying conditions~(2.3)--(2.5), 
is regular if: (i) $h_{V}^{o}(x) \! \not\equiv \! 0$ on $\overline{J_{o}}$; 
(ii) $2(2 \! + \! \tfrac{1}{n}) \int_{J_{o}} \ln (\vert x \! - \! s \vert) \, 
\md \mu_{V}^{o}(s) \! - \! 2 \ln \vert x \vert \! - \! \widetilde{V}(x) \! - 
\! \ell_{o} \! - \! 2(2 \! + \! \tfrac{1}{n})Q_{o} \! < \! 0$, $x \! \in \! 
\mathbb{R} \setminus \overline{J_{o}}$, where $Q_{o} \! := \! \int_{J_{o}} 
\ln (\lvert s \rvert) \, \md \mu_{V}^{o}(s)$; and (iii) 
inequalities~$(\mathrm{P}_{2}^{(4)})$ and~$(\mathrm{P}_{2}^{(6)})$ in $\pmb{
\mathrm{P}_{2}}$ are strict, that is, $\leqslant$ (resp., $\geqslant)$ is 
replaced by $<$ (resp., $>)$\footnote{There are three distinct situations in 
which these conditions may fail: (i) for at least one $\widetilde{x}_{e} \! 
\in \! \mathbb{R} \setminus \overline{J_{e}}$ (resp., $\widetilde{x}_{o} \!
\in \! \mathbb{R} \setminus \widetilde{J_{o}})$, $4 \int_{J_{e}} \ln (\lvert
\widetilde{x}_{e} \! - \! s \rvert) \, \md \mu_{V}^{e}(s) \! - \! 2 \ln \lvert
\widetilde{x}_{e} \rvert \! - \! \widetilde{V}(\widetilde{x}_{e}) \! - \!
\ell_{e} \! = \! 0$ (resp., $2(2 \! + \! \tfrac{1}{n}) \int_{J_{o}} \ln
(\lvert \widetilde{x}_{o} \! - \! s \rvert) \, \md \mu_{V}^{o}(s) \! - \! 2
\ln \lvert \widetilde{x}_{o} \rvert \! - \! \widetilde{V}(\widetilde{x}_{o})
\! - \! \ell_{o} \! - \! 2(2 \! + \! \tfrac{1}{n})Q_{o} \! = \! 0)$, that is,
for $n$ even (resp., $n$ odd) equality is attained for at least one point
$\widetilde{x}_{e}$ (resp., $\widetilde{x}_{o})$ in the complement of the
closure of the support of the `even' (resp., `odd') equilibrium measure $\mu_{
V}^{e}$ (resp., $\mu_{V}^{o})$, which corresponds to the situation in which a
`band' has just closed, or is about to open, about $\widetilde{x}_{e}$ (resp.,
$\widetilde{x}_{o})$; (ii) for at least one $\widehat{x}_{e}$ (resp.,
$\widehat{x}_{o})$, $h_{V}^{e}(\widehat{x}_{e}) \! = \! 0$ (resp., $h_{V}^{o}
(\widehat{x}_{o}) \! = \! 0)$, that is, for $n$ even (resp., $n$ odd) the
function $h_{V}^{e}$ (resp., $h_{V}^{o})$ vanishes for at least one point
$\widehat{x}_{e}$ (resp., $\widehat{x}_{o})$ within the support of the `even'
(resp., `odd') equilibrium measure $\mu_{V}^{e}$ (resp., $\mu_{V}^{o})$, which
corresponds to the situation in which a `gap' is about to open, or close,
about $\widehat{x}_{e}$ (resp., $\widehat{x}_{o})$; and (iii) there exists at
least one $j \! \in \! \lbrace 1,\dotsc,N \! + \! 1 \rbrace$, denoted $j_{e}$
(resp., $j_{o})$, such that $h_{V}^{e}(b_{j_{e}-1}^{e}) \! = \! 0$ and/or $h_{
V}^{e}(a_{j_{e}}^{e}) \! = \! 0$ (resp., $h_{V}^{o}(b_{j_{o}-1}^{o}) \! = \!
0$ and/or $h_{V}^{o}(a_{j_{o}}^{o}) \! = \! 0)$. Each of these three cases
can occur only a finite number of times due to the fact that $\widetilde{V}
\colon \mathbb{R} \setminus \lbrace 0 \rbrace \! \to \! \mathbb{R}$ satisfies
conditions~(2.3)--(2.5) \cite{a46,a81}.}.
\end{compactenum}
The (density of the) `even' and `odd' equilibrium measures $\md \mu_{V}^{e}$
and $\md \mu_{V}^{o}$, respectively, together with the corresponding
variational problems, emerge naturally in the asymptotic analyses of
\textbf{RHP1} and \textbf{RHP2}.
\begin{eeee}
The following correspondences should also be noted:
\begin{compactenum}
\item[\textbullet] $g^{e} \colon \mathbb{C} \setminus (-\infty, \max \{0,
a_{N+1}^{e}\}) \! \to \! \mathbb{C}$ solves the 
\emph{phase conditions}~$(\mathrm{P}_{1}^{(1)})$--$(\mathrm{P}_{1}^{(6)}) \! 
\Leftrightarrow \! \mathcal{M}_{1}(\mathbb{R}) \! \ni \! \mu_{V}^{e}$ solves 
the variational conditions~$(\mathrm{P}_{1}^{(a)})$ and~$(\mathrm{P}_{1}^{
(b)})$;
\item[\textbullet] $g^{o} \colon \mathbb{C} \setminus (-\infty, \max \{0,
a_{N+1}^{o}\}) \! \to \! \mathbb{C}$ solves the phase 
conditions~$(\mathrm{P}_{2}^{(1)})$--$(\mathrm{P}_{2}^{(6)}) \! 
\Leftrightarrow \! \mathcal{M}_{1}(\mathbb{R}) \! \ni \! \mu_{V}^{o}$ solves 
the variational conditions~$(\mathrm{P}_{2}^{(a)})$ and~$(\mathrm{P}_{2}^{
(b)})$. \hfill $\blacksquare$
\end{compactenum}
\end{eeee}

Since the main results of this paper are asymptotics (as $n \! \to \! \infty)$ 
for $\boldsymbol{\pi}_{2n+1}(z)$ $(z \! \in \! \mathbb{C})$, $\xi^{(2n+1)}_{-
n-1}$ and $\phi_{2n+1}(z)$ $(z \! \in \! \mathbb{C})$, which are, via 
Lemma~2.2.2, Equation~(2.2), and Equations~(1.3) and~(1.5), related to 
\textbf{RHP2} for $\overset{o}{\mathrm{Y}} \colon \mathbb{C} \setminus 
\mathbb{R} \! \to \! \operatorname{SL}_{2}(\mathbb{C})$, no further 
reference, henceforth, to \textbf{RHP1} (and Lemma~2.2.1) for $\overset{e}{
\mathrm{Y}} \colon \mathbb{C} \setminus \mathbb{R} \! \to \! 
\operatorname{SL}_{2}(\mathbb{C})$ will be made (see \cite{a38} for the 
complete details of the asymptotic analysis of \textbf{RHP1}). In the ensuing 
analysis, the large-$n$ behaviour of the solution of \textbf{RHP2} (see 
Lemma~2.2.2, Equation~(2.2)), hence asymptotics for $\boldsymbol{\pi}_{2n+1}
(z)$ (in the entire complex plane), $\xi^{(2n+1)}_{-n-1}$ and $\phi_{2n+1}
(z)$ (in the entire complex plane), are extracted.
\subsection{Summary of Results}
In this subsection, the final results of this work are presented (see
Sections~3--5 for the detailed analyses and proofs). Before doing so, however,
some notational preamble is necessary. For $j \! = \! 1,\dotsc,N \! + \! 1$,
let
\begin{equation*}
\Phi_{a_{j}}^{o}(z) \! := \! \left(\dfrac{3}{2} \! \left(n \! + \! \dfrac{1}{
2} \right) \! \int_{a_{j}^{o}}^{z}(R_{o}(s))^{1/2}h_{V}^{o}(s) \, \md s
\right)^{2/3},
\end{equation*}
and
\begin{equation*}
\Phi_{b_{j-1}}^{o}(z) \! := \! \left(-\dfrac{3}{2} \! \left(n \! + \! \dfrac{
1}{2} \right) \! \int_{z}^{b_{j-1}^{o}}(R_{o}(s))^{1/2}h_{V}^{o}(s) \, \md s
\right)^{2/3},
\end{equation*}
where $(R_{o}(z))^{1/2}$ and $h_{V}^{o}(z)$ are defined in Theorem~2.3.1,
Equations~(2.8) and~(2.9). Define the `small', mutually disjoint open discs
about the end-points of the support of the `odd' equilibrium measure,
$\lbrace b_{j-1}^{o},a_{j}^{o} \rbrace_{j=1}^{N+1}$, as follows: for $j \! =
\! 1,\dotsc,N \! + \! 1$,
\begin{equation*}
\mathbb{U}^{o}_{\delta_{a_{j}}} \! := \! \left\{ \mathstrut z \! \in \!
\mathbb{C}; \, \vert z \! - \! a_{j}^{o} \vert \! < \! \delta_{a_{j}}^{o}
\right\} \qquad \text{and} \qquad \mathbb{U}^{o}_{\delta_{b_{j-1}}} \! :=
\! \left\{ \mathstrut z \! \in \! \mathbb{C}; \, \vert z \! - \! b_{j-1}^{o}
\vert \! < \! \delta^{o}_{b_{j-1}} \right\},
\end{equation*}
where $(0,1) \! \ni \! \delta^{o}_{a_{j}}$ (resp., $(0,1) \! \ni \! \delta^{
o}_{b_{j-1}})$ are chosen `sufficiently small' so that $\Phi^{o}_{a_{j}}(z)$
(resp., $\Phi^{o}_{b_{j-1}}(z))$, which are bi-holomorphic, conformal, and
orientation preserving (resp., bi-holomorphic, conformal, and non-orientation
preserving), map $\mathbb{U}^{o}_{\delta_{a_{j}}}$ (resp., $\mathbb{U}^{o}_{
\delta_{b_{j-1}}})$, as well as the oriented skeletons (see Figure~5) $\cup_{
l=1}^{4} \Sigma^{o,l}_{a_{j}}$ (resp., $\cup_{l=1}^{4} \Sigma^{o,l}_{b_{j-1}
}$ (see Figure~6)), injectively onto open (and convex), $n$-dependent
neighbourhoods of $0$ such that:
\begin{figure}[tbh]
\begin{center}
\vspace{0.55cm}
\begin{pspicture}(0,0)(14,8)
\psset{xunit=1cm,yunit=1cm}
\pscircle[linewidth=0.7pt,linestyle=solid,linecolor=red](3.5,4){2.5}
\pscircle[linewidth=0.7pt,linestyle=solid,linecolor=red](10.5,4){2.5}
\psarcn[linewidth=0.6pt,linestyle=solid,linecolor=blue,arrowsize=1.5pt 5]{->}%
(0.5,4){3}{65}{36}
\pstextpath[c]{\psarcn[linewidth=0.6pt,linestyle=solid,linecolor=blue](0.5,4)%
{3}{36}{0}}{\makebox(0,0){$\pmb{\Sigma^{o,1}_{a_{j}}}$}}
\psarc[linewidth=0.6pt,linestyle=solid,linecolor=green,arrowsize=1.5pt 5]{->}%
(0.5,4){3}{295}{324}
\pstextpath[c]{\psarc[linewidth=0.6pt,linestyle=solid,linecolor=green](0.5,4)%
{3}{324}{360}}{\makebox(0,0){$\pmb{\Sigma^{o,3}_{a_{j}}}$}}
\psline[linewidth=0.6pt,linestyle=solid,linecolor=cyan,arrowsize=1.5pt 5]{->}%
(0.5,4)(1.6,4)
\pstextpath[c]{\psline[linewidth=0.6pt,linestyle=solid,linecolor=cyan](1.6,4)%
(3.5,4)}{\makebox(0,0){$\pmb{\Sigma^{o,2}_{a_{j}}}$}}
\psline[linewidth=0.6pt,linestyle=solid,linecolor=magenta](5.5,4)(6.3,4)
\pstextpath[c]{\psline[linewidth=0.6pt,linestyle=solid,linecolor=magenta,%
arrowsize=1.5pt 5]{->}(3.5,4)(5.5,4)}{\makebox(0,0){%
$\pmb{\Sigma^{o,4}_{a_{j}}}$}}
\rput(3.5,7.5){\makebox(0,0){$z-\text{plane}$}}
\rput(3.5,0.5){\makebox(0,0){$\mathbb{U}^{o}_{\delta_{a_{j}}}$}}
\pszigzag[coilwidth=0.3cm,coilarm=0.25cm,coilaspect=45]{->}(3.7,0.7)(4.1,2.2)
\rput(4.65,5.1){\makebox(0,0){$\pmb{\Omega^{o,1}_{a_{j}}}$}}
\rput(2.35,5.1){\makebox(0,0){$\pmb{\Omega^{o,2}_{a_{j}}}$}}
\rput(2.35,2.9){\makebox(0,0){$\pmb{\Omega^{o,3}_{a_{j}}}$}}
\rput(4.65,2.9){\makebox(0,0){$\pmb{\Omega^{o,4}_{a_{j}}}$}}
\rput(7,7.2){\makebox(0,0){$\zeta \! = \! \Phi^{o}_{a_{j}}(z)$}}
\rput(7,0.8){\makebox(0,0){$z \! = \! (\Phi^{o}_{a_{j}})^{-1}(\zeta)$}}
\psarcn[linewidth=0.8pt,linestyle=solid,linecolor=black,arrowsize=1.5pt 5]%
{->}(7,5.3){1.5}{135}{45}
\psarc[linewidth=0.8pt,linestyle=solid,linecolor=black,arrowsize=1.5pt 5]%
{<-}(7,2.7){1.5}{225}{315}
\rput(10.5,7.5){\makebox(0,0){$\zeta -\text{plane}$}}
\rput(10.5,0.5){\makebox(0,0){$\widehat{\mathbb{U}}^{o}_{\delta_{a_{j}}} \!
:= \! \Phi^{o}_{a_{j}}(\mathbb{U}^{o}_{\delta_{a_{j}}})$}}
\pszigzag[coilwidth=0.3cm,coilarm=0.25cm,coilaspect=45]{->}(10.75,0.75)%
(10.15,2.2)
\psline[linewidth=0.6pt,linestyle=solid,linecolor=cyan,arrowsize=1.5pt 5]{->}%
(7.6,4)(8.6,4)
\pstextpath[c]{\psline[linewidth=0.6pt,linestyle=solid,linecolor=cyan](8.6,4)%
(10.5,4)}{\makebox(0,0){$\pmb{\gamma^{o,2}_{a_{j}}}$}}
\psline[linewidth=0.6pt,linestyle=solid,linecolor=magenta](12.5,4)(13.3,4)
\pstextpath[c]{\psline[linewidth=0.6pt,linestyle=solid,linecolor=magenta,%
arrowsize=1.5pt 5]{->}(10.5,4)(12.5,4)}{\makebox(0,0){$\pmb{\gamma^{o,4}_{%
a_{j}}}$}}
\psline[linewidth=0.6pt,linestyle=solid,linecolor=blue,arrowsize=1.5pt 5]{->}%
(8.4,6.1)(9.3,5.2)
\pstextpath[c]{\psline[linewidth=0.6pt,linestyle=solid,linecolor=blue]%
(9.3,5.2)(10.5,4)}{\makebox(0,0){$\pmb{\gamma^{o,1}_{a_{j}}}$}}
\psline[linewidth=0.6pt,linestyle=solid,linecolor=green,arrowsize=1.5pt 5]{->}%
(8.4,1.9)(9.3,2.8)
\pstextpath[c]{\psline[linewidth=0.6pt,linestyle=solid,linecolor=green]%
(9.3,2.8)(10.5,4)}{\makebox(0,0){$\pmb{\gamma^{o,3}_{a_{j}}}$}}
\rput(11.5,4.85){\makebox(0,0){$\pmb{\widehat{\Omega}^{o,1}_{a_{j}}}$}}
\rput(9,4.65){\makebox(0,0){$\pmb{\widehat{\Omega}^{o,2}_{a_{j}}}$}}
\rput(9,3.35){\makebox(0,0){$\pmb{\widehat{\Omega}^{o,3}_{a_{j}}}$}}
\rput(11.5,3.15){\makebox(0,0){$\pmb{\widehat{\Omega}^{o,4}_{a_{j}}}$}}
\psdots[dotstyle=*,dotscale=1.5](3.5,4)
\psdots[dotstyle=*,dotscale=1.5](10.5,4)
\rput(3.5,3.7){\makebox(0,0){$\pmb{a_{j}^{o}}$}}
\rput(10.5,3.7){\makebox(0,0){$\pmb{0}$}}
\end{pspicture}
\end{center}
\caption{The conformal mapping $\zeta \! = \! \Phi^{o}_{a_{j}}(z) \! := \!
(\tfrac{3}{2}(n \! + \! \tfrac{1}{2}) \int_{a_{j}^{o}}^{z}(R_{o}(s))^{1/2}h_{
V}^{o}(s))^{2/3}$, $j \! = \! 1,\dotsc,N \! + \! 1$, where $(\Phi^{o}_{a_{j}}
)^{-1}$ denotes the inverse mapping}
\end{figure}
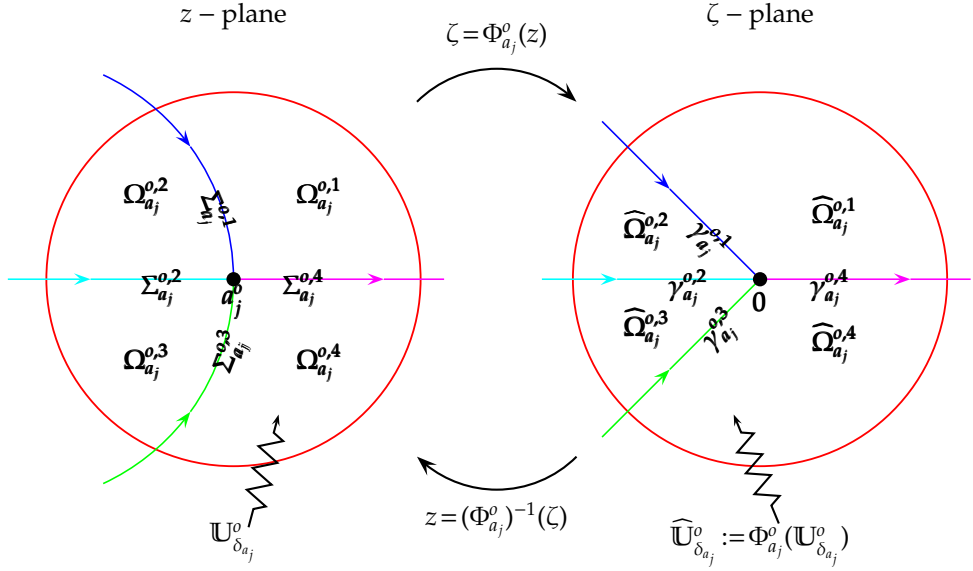
\begin{figure}[bht]
\begin{center}
\vspace{0.55cm}
\begin{pspicture}(0,0)(14,8)
\psset{xunit=1cm,yunit=1cm}
\pscircle[linewidth=0.7pt,linestyle=solid,linecolor=red](3.5,4){2.5}
\pscircle[linewidth=0.7pt,linestyle=solid,linecolor=red](10.5,4){2.5}
\psarcn[linewidth=0.6pt,linestyle=solid,linecolor=blue](6.5,4){3}{136}{120}
\pstextpath[c]{\psarcn[linewidth=0.6pt,linestyle=solid,linecolor=blue,%
arrowsize=1.5pt 5]{->}(6.5,4){3}{180}{136}}{\makebox(0,0){$\pmb{\Sigma^{o,%
1}_{b_{j-1}}}$}}
\psarc[linewidth=0.6pt,linestyle=solid,linecolor=green](6.5,4){3}{228}{240}
\pstextpath[c]{\psarc[linewidth=0.6pt,linestyle=solid,linecolor=green,%
arrowsize=1.5pt 5]{->}(6.5,4){3}{180}{228}}{\makebox(0,0){$\pmb{\Sigma^{o,%
3}_{b_{j-1}}}$}}
\psline[linewidth=0.6pt,linestyle=solid,linecolor=cyan,arrowsize=1.5pt 5]{->}%
(0.5,4)(1.6,4)
\pstextpath[c]{\psline[linewidth=0.6pt,linestyle=solid,linecolor=cyan](1.6,4)%
(3.5,4)}{\makebox(0,0){$\pmb{\Sigma^{o,4}_{b_{j-1}}}$}}
\psline[linewidth=0.6pt,linestyle=solid,linecolor=magenta](5.5,4)(6.3,4)
\pstextpath[c]{\psline[linewidth=0.6pt,linestyle=solid,linecolor=magenta,%
arrowsize=1.5pt 5]{->}(3.5,4)(5.5,4)}{\makebox(0,0){%
$\pmb{\Sigma^{o,2}_{b_{j-1}}}$}}
\rput(3.5,7.5){\makebox(0,0){$z-\text{plane}$}}
\rput(3.5,0.5){\makebox(0,0){$\mathbb{U}^{o}_{\delta_{b_{j-1}}}$}}
\pszigzag[coilwidth=0.3cm,coilarm=0.25cm,coilaspect=45]{->}(3.7,0.7)(3.9,2.2)
\rput(4.8,5.1){\makebox(0,0){$\pmb{\Omega^{o,2}_{b_{j-1}}}$}}
\rput(2.35,5.1){\makebox(0,0){$\pmb{\Omega^{o,1}_{b_{j-1}}}$}}
\rput(2.35,2.9){\makebox(0,0){$\pmb{\Omega^{o,4}_{b_{j-1}}}$}}
\rput(4.8,2.9){\makebox(0,0){$\pmb{\Omega^{o,3}_{b_{j-1}}}$}}
\rput(7,7.2){\makebox(0,0){$\zeta \! = \! \Phi^{o}_{b_{j-1}}(z)$}}
\rput(7,0.8){\makebox(0,0){$z \! = \! (\Phi^{o}_{b_{j-1}})^{-1}(\zeta)$}}
\psarcn[linewidth=0.8pt,linestyle=solid,linecolor=black,arrowsize=1.5pt 5]%
{->}(7,5.3){1.5}{135}{45}
\psarc[linewidth=0.8pt,linestyle=solid,linecolor=black,arrowsize=1.5pt 5]%
{<-}(7,2.7){1.5}{225}{315}
\rput(10.5,7.5){\makebox(0,0){$\zeta -\text{plane}$}}
\rput(10.5,0.5){\makebox(0,0){$\widehat{\mathbb{U}}^{o}_{\delta_{b_{j-1}}} \!
:= \! \Phi^{o}_{b_{j-1}}(\mathbb{U}^{o}_{\delta_{b_{j-1}}})$}}
\pszigzag[coilwidth=0.3cm,coilarm=0.25cm,coilaspect=45]{->}(10.75,0.75)%
(10.55,2.2)
\psline[linewidth=0.6pt,linestyle=solid,linecolor=magenta,arrowsize=1.5pt 5]%
{-<}(7.6,4)(8.7,4)
\pstextpath[c]{\psline[linewidth=0.6pt,linestyle=solid,linecolor=magenta]%
(8.6,4)(10.5,4)}{\makebox(0,0){$\pmb{\gamma^{o,2}_{b_{j-1}}}$}}
\psline[linewidth=0.6pt,linestyle=solid,linecolor=cyan](12.5,4)(13.3,4)
\pstextpath[c]{\psline[linewidth=0.6pt,linestyle=solid,linecolor=cyan,%
arrowsize=1.5pt 5]{-<}(10.5,4)(12.6,4)}{\makebox(0,0){$\pmb{\gamma^{o,4}_{%
b_{j-1}}}$}}
\psline[linewidth=0.6pt,linestyle=solid,linecolor=blue,arrowsize=1.5pt 5]{-<}%
(8.4,6.1)(9.2,5.3)
\pstextpath[c]{\psline[linewidth=0.6pt,linestyle=solid,linecolor=blue]%
(9.1,5.4)(10.5,4)}{\makebox(0,0){$\pmb{\gamma^{o,1}_{b_{j-1}}}$}}
\psline[linewidth=0.6pt,linestyle=solid,linecolor=green,arrowsize=1.5pt 5]{-<}%
(8.4,1.9)(9.2,2.7)
\pstextpath[c]{\psline[linewidth=0.6pt,linestyle=solid,linecolor=green]%
(9.1,2.6)(10.5,4)}{\makebox(0,0){$\pmb{\gamma^{o,3}_{b_{j-1}}}$}}
\rput(11.5,4.85){\makebox(0,0){$\pmb{\widehat{\Omega}^{o,1}_{b_{j-1}}}$}}
\rput(8.9,4.65){\makebox(0,0){$\pmb{\widehat{\Omega}^{o,2}_{b_{j-1}}}$}}
\rput(8.9,3.35){\makebox(0,0){$\pmb{\widehat{\Omega}^{o,3}_{b_{j-1}}}$}}
\rput(11.5,3.15){\makebox(0,0){$\pmb{\widehat{\Omega}^{o,4}_{b_{j-1}}}$}}
\psdots[dotstyle=*,dotscale=1.5](3.5,4)
\psdots[dotstyle=*,dotscale=1.5](10.5,4)
\rput(3.5,3.6){\makebox(0,0){$\pmb{b_{j-1}^{o}}$}}
\rput(10.5,3.7){\makebox(0,0){$\pmb{0}$}}
\end{pspicture}
\end{center}
\caption{The conformal mapping $\zeta \! = \! \Phi^{o}_{b_{j-1}}(z) \! := \!
(-\tfrac{3}{2}(n \! + \! \tfrac{1}{2}) \int_{z}^{b_{j-1}^{o}}(R_{o}(s))^{1/2}
h_{V}^{o}(s))^{2/3}$, $j \! = \! 1,\dotsc,N \! + \! 1$, where $(\Phi^{o}_{b_{
j-1}})^{-1}$ denotes the inverse mapping}
\end{figure}
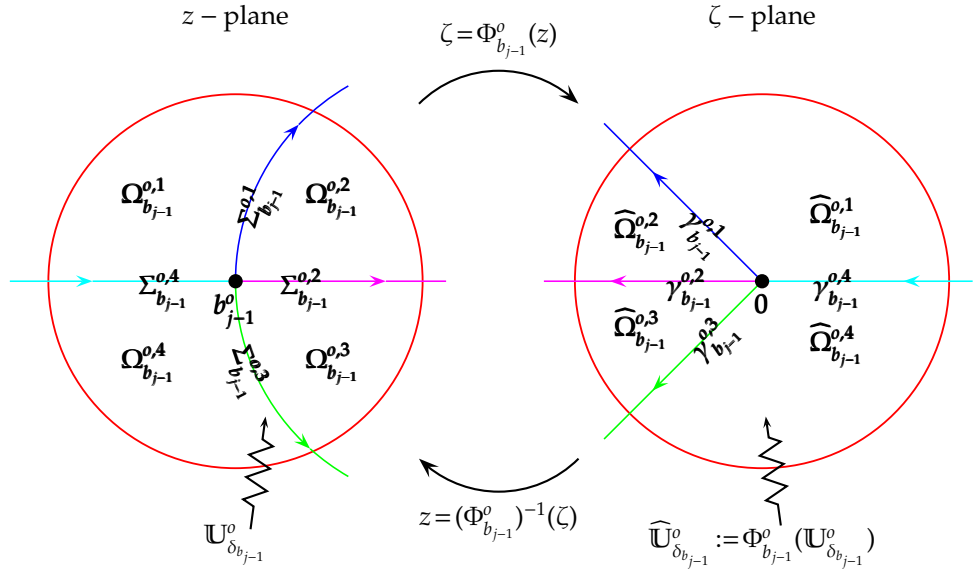
\begin{description}
\item[\pmb{(i)}] $\Phi^{o}_{a_{j}}(a_{j}^{o}) \! = \! 0$ (resp., $\Phi^{o}_{
b_{j-1}}(b_{j-1}^{o}) \! = \! 0)$;
\item[\pmb{(ii)}] $\Phi^{o}_{a_{j}} \colon \mathbb{U}^{o}_{\delta_{a_{j}}} \!
\to \! \widehat{\mathbb{U}}^{o}_{\delta_{a_{j}}} \! := \! \Phi^{o}_{a_{j}}
(\mathbb{U}^{o}_{\delta_{a_{j}}})$ (resp., $\Phi^{o}_{b_{j-1}} \colon \mathbb{
U}^{o}_{\delta_{b_{j-1}}} \! \to \! \widehat{\mathbb{U}}^{o}_{\delta_{b_{j-1}
}} \! := \! \Phi^{o}_{b_{j-1}}(\mathbb{U}^{o}_{\delta_{b_{j-1}}}))$;
\item[\pmb{(iii)}] $\Phi^{o}_{a_{j}}(\mathbb{U}^{o}_{\delta_{a_{j}}} \cap
\Sigma^{o,l}_{a_{j}}) \! = \! \Phi^{o}_{a_{j}}(\mathbb{U}^{o}_{\delta_{a_{j}}
}) \cap \gamma^{o,l}_{a_{j}}$ (resp., $\Phi^{o}_{b_{j-1}}(\mathbb{U}^{o}_{
\delta_{b_{j-1}}} \cap \Sigma^{o,l}_{b_{j-1}}) \! = \! \Phi^{o}_{b_{j-1}}
(\mathbb{U}^{o}_{\delta_{b_{j-1}}}) \cap \gamma^{o,l}_{b_{j-1}})$;
\item[\pmb{(iv)}] $\Phi^{o}_{a_{j}}(\mathbb{U}^{o}_{\delta_{a_{j}}} \cap
\Omega^{o,l}_{a_{j}}) \! = \! \Phi^{o}_{a_{j}}(\mathbb{U}^{o}_{\delta_{a_{j}}
}) \cap \widehat{\Omega}^{o,l}_{a_{j}}$ (resp., $\Phi^{o}_{b_{j-1}}(\mathbb{
U}^{o}_{\delta_{b_{j-1}}} \cap \Omega^{o,l}_{b_{j-1}}) \! = \! \Phi^{o}_{b_{j
-1}}(\mathbb{U}^{o}_{\delta_{b_{j-1}}}) \cap \widehat{\Omega}^{o,l}_{b_{j-1}
})$, with $\widehat{\Omega}^{o,1}_{a_{j}}$ (and $\widehat{\Omega}^{o,1}_{b_{j
-1}})$ $= \! \lbrace \mathstrut \zeta \! \in \! \mathbb{C}; \, \arg (\zeta) \!
\in \! (0,2 \pi/3) \rbrace$, $\widehat{\Omega}^{o,2}_{a_{j}}$ (and $\widehat{
\Omega}^{o,2}_{b_{j-1}})$ $= \! \lbrace \mathstrut \zeta \! \in \! \mathbb{C};
\, \arg (\zeta) \! \in \! (2 \pi/3,\pi) \rbrace$, $\widehat{\Omega}^{o,3}_{a_{
j}}$ (and $\widehat{\Omega}^{o,3}_{b_{j-1}})$ $= \! \lbrace \mathstrut \zeta
\! \in \! \mathbb{C}; \, \arg (\zeta) \! \in \! (-\pi,-2 \pi/3) \rbrace$, and
$\widehat{\Omega}^{o,4}_{a_{j}}$ (and $\widehat{\Omega}^{o,4}_{b_{j-1}})$ $=
\! \lbrace \mathstrut \zeta \! \in \! \mathbb{C}; \, \arg (\zeta) \! \in \!
(-2 \pi/3,0) \rbrace$\footnote{The precise angles between the sectors are 
not absolutely important; one could, for example, replace $2 \pi/3$ by any 
angle strictly between $0$ and $\pi$ \cite{a2,a46,a47,a49,a79}.}.
\end{description}

Introduce, now, the Airy function, $\operatorname{Ai}(\cdot)$, which appears 
in several of the final results of this work: $\operatorname{Ai}(\cdot)$ 
is determined (uniquely) as the solution of the second-order, non-constant 
coefficient, homogeneous ODE (see, for example, Chapter~10 of \cite{a82})
\begin{equation*}
\operatorname{Ai}^{\prime \prime}(z) \! - \! z \operatorname{Ai}(z) \! = \! 
0,
\end{equation*}
with asymptotics (at infinity)
\begin{equation}
\begin{split}
\operatorname{Ai}(z) \underset{\underset{\vert \arg z \vert < \pi}{z \to
\infty}}{\sim}& \, \dfrac{1}{2 \sqrt{\smash[b]{\pi}}}z^{-1/4} \, \me^{-
\widehat{\zeta}(z)} \sum_{k=0}^{\infty}(-1)^{k}s_{k}(\widehat{\zeta}(z))^{-k},
\qquad \widehat{\zeta}(z) \! := \! \dfrac{2}{3}z^{3/2}, \\
\operatorname{Ai}^{\prime}(z) \underset{\underset{\vert \arg z \vert < \pi}{z
\to \infty}}{\sim}& \, -\dfrac{1}{2 \sqrt{\smash[b]{\pi}}}z^{1/4} \, \me^{-
\widehat{\zeta}(z)} \sum_{k=0}^{\infty}(-1)^{k}t_{k}(\widehat{\zeta}(z))^{-k},
\end{split}
\end{equation}
where $s_{0} \! = \! t_{0} \! = \! 1$,
\begin{equation*}
s_{k} \! = \! \dfrac{\Gamma (3k \! + \! 1/2)}{54^{k}k! \Gamma (k \! + \! 1/2)}
\! = \! \dfrac{(2k \! + \! 1)(2k \! + \! 3) \cdots (6k \! - \! 1)}{216^{k}k!},
\qquad \quad t_{k} \! = \! -\left(\dfrac{6k \! + \! 1}{6k \! -
\! 1} \right) \! s_{k}, \quad k \! \in \! \mathbb{N},
\end{equation*}
and $\Gamma (\pmb{\cdot})$ is the gamma (factorial) function.

In order to present the final asymptotic (as $n \! \to \! \infty)$ results,
and for arbitrary $j \! = \! 1,\dotsc,N \! + \! 1$, consider the following
decomposition (see Figure~7), into bounded and unbounded regions, of $\mathbb{
C}$ and the neighbourhoods of the end-points $b_{i-1}^{o}$, $a_{i}^{o}$, $i \!
= \! 1,\dotsc,N \! + \! 1$ (as per the discussion above, $\mathbb{U}^{o}_{
\delta_{b_{k-1}}} \cap \mathbb{U}^{o}_{\delta_{a_{k}}} \! = \! \varnothing$,
$k \! = \! 1,\dotsc,N \! + \! 1)$. Asymptotics (as $n \! \to \! \infty)$ for
$\boldsymbol{\pi}_{2n+1}(z)$, with $z \! \in \! \cup_{j=1}^{4}(\Upsilon_{j}^{
o} \cup (\cup_{k=1}^{N+1}(\Omega_{b_{k-1}}^{o,j} \cup \Omega_{a_{k}}^{o,j}))
)$, are now presented. These asymptotic expansions are obtained via a union of
the DZ non-linear steepest-descent method \cite{a1,a2} and the extension of
Deift-Venakides-Zhou \cite{a3} (see, also, \cite{a45,a46,a47,a48,a49,a50,a51,%
a52,a53,a54,a55,a56,a57,a58,a59,a60,a61,a62,a63,a66,a67,a68,a69}, and the 
pedagogical exposition \cite{a79}).
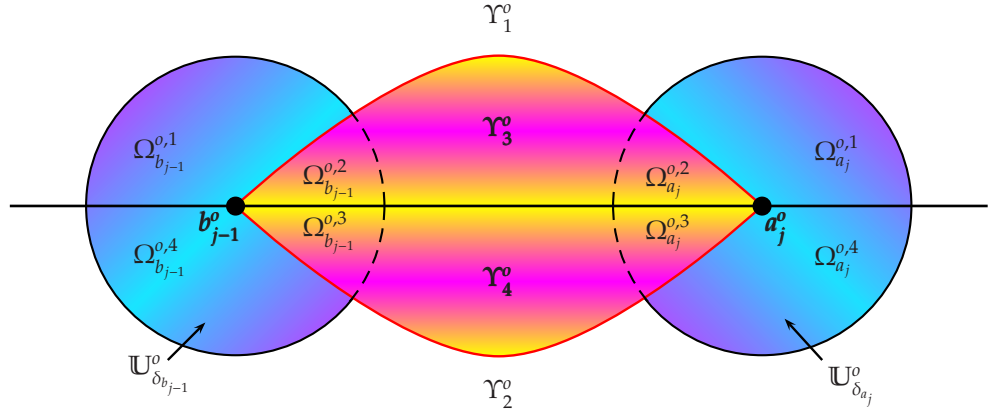
\begin{figure}[tbh]
\begin{center}
\vspace{0.55cm}
\begin{pspicture}(0,0)(14,6)
\psset{xunit=1cm,yunit=1cm}
\pscircle[fillstyle=gradient,gradangle=45,gradlines=600,gradbegin=magenta,%
gradend=cyan,gradmidpoint=0.5](3.5,3){2}
\pscircle[fillstyle=gradient,gradangle=135,gradlines=600,gradbegin=magenta,%
gradend=cyan,gradmidpoint=0.5](10.5,3){2}
\pscurve[linewidth=0.9pt,linestyle=solid,linecolor=red,fillstyle=gradient,%
gradangle=0,gradlines=600,gradbegin=yellow,gradend=magenta,gradmidpoint=0.5]%
(3.5,3)(7,5)(10.5,3)
\pscurve[linewidth=0.9pt,linestyle=solid,linecolor=red,fillstyle=gradient,%
gradangle=0,gradlines=600,gradbegin=yellow,gradend=magenta,gradmidpoint=0.5]%
(3.5,3)(7,1)(10.5,3)
\rput(7,5.5){\makebox(0,0){$\Upsilon^{o}_{1}$}}
\rput(7,0.5){\makebox(0,0){$\Upsilon^{o}_{2}$}}
\rput(7,4){\makebox(0,0){$\pmb{\Upsilon^{o}_{3}}$}}
\rput(7,2){\makebox(0,0){$\pmb{\Upsilon^{o}_{4}}$}}
\rput(2.5,0.75){\makebox(0,0){$\mathbb{U}^{o}_{\delta_{b_{j-1}}}$}}
\rput(11.70,0.60){\makebox(0,0){$\mathbb{U}^{o}_{\delta_{a_{j}}}$}}
\psline[linewidth=0.9pt,linestyle=solid,linecolor=black]{->}(2.6,0.95)%
(3.1,1.45)
\psline[linewidth=0.9pt,linestyle=solid,linecolor=black]{->}(11.3,0.95)%
(10.8,1.50)
\psline[linewidth=1.0pt,linestyle=solid,linecolor=black](0.5,3)(13.5,3)
\rput(2.5,3.7){\makebox(0,0){$\Omega^{o,1}_{b_{j-1}}$}}
\rput(2.5,2.3){\makebox(0,0){$\Omega^{o,4}_{b_{j-1}}$}}
\rput(4.75,3.35){\makebox(0,0){$\Omega^{o,2}_{b_{j-1}}$}}
\rput(4.75,2.65){\makebox(0,0){$\Omega^{o,3}_{b_{j-1}}$}}
\rput(9.25,3.35){\makebox(0,0){$\Omega^{o,2}_{a_{j}}$}}
\rput(9.25,2.65){\makebox(0,0){$\Omega^{o,3}_{a_{j}}$}}
\rput(11.5,3.7){\makebox(0,0){$\Omega^{o,1}_{a_{j}}$}}
\rput(11.5,2.3){\makebox(0,0){$\Omega^{o,4}_{a_{j}}$}}
\psarcn[linewidth=0.8pt,linestyle=dashed,linecolor=black](3.5,3){1.98}{40}{0}
\psarc[linewidth=0.8pt,linestyle=dashed,linecolor=black](3.5,3){1.98}{320}{360}
\psarc[linewidth=0.8pt,linestyle=dashed,linecolor=black](10.5,3){1.98}{140}%
{180}
\psarc[linewidth=0.8pt,linestyle=dashed,linecolor=black](10.5,3){1.98}{180}%
{220}
\psdots[dotstyle=*,dotscale=2](3.5,3)
\psdots[dotstyle=*,dotscale=2](10.5,3)
\rput(3.3,2.7){\makebox(0,0){$\pmb{b_{j-1}^{o}}$}}
\rput(10.7,2.7){\makebox(0,0){$\pmb{a_{j}^{o}}$}}
\end{pspicture}
\end{center}
\caption{Region-by-region decomposition of $\mathbb{C}$ and the neighbourhoods
surrounding the end-points of the support of the `odd' equilibrium measure,
$\lbrace b_{j-1}^{o},a_{j}^{o} \rbrace_{j=1}^{N+1}$}
\end{figure}
\begin{eeee}
In order to eschew a flood of superfluous notation, the simplified `notation'
$\mathcal{O}((n \! + \! 1/2)^{-2})$ is maintained throughout Theorem~2.3.1
(see below), and is to be understood in the following, \emph{normal} sense:
for a compact subset, $\mathfrak{D}$, say, of $\mathbb{C}$, and uniformly with
respect to $z \! \in \! \mathfrak{D}$, $\mathcal{O}((n \! + \! 1/2)^{-2}) \!
:= \! \mathcal{O}(c^{\natural}(z,n)(n \! + \! 1/2)^{-2})$, where $\norm{c^{
\natural}(\boldsymbol{\cdot},n)}_{\mathcal{L}^{p}(\mathfrak{D})} \! =_{n \to
\infty} \! \mathcal{O}(1)$, $p \! \in \! \lbrace 1,2,\infty \rbrace$, and
$\exists \, \, K_{\mathfrak{D}} \! > \! 0$ (and finite) such that, $\forall
\, \, z \! \in \! \mathfrak{D}$, $\vert c^{\natural}(z,n) \vert \! \leqslant_{
n \to \infty} \! K_{\mathfrak{D}}$. \hfill $\blacksquare$
\end{eeee}
\begin{dddd}
Let the external field $\widetilde{V} \colon \mathbb{R} \setminus \{0\} \! \to
\! \mathbb{R}$ satisfy conditions~{\rm (2.3)--(2.5)}. Set
\begin{equation}
\md \mu_{V}^{o}(x) \! := \! \psi_{V}^{o}(x) \, \md x \! = \! \dfrac{1}{2 \pi
\mi}(R_{o}(x))^{1/2}_{+}h_{V}^{o}(x) \pmb{1}_{J_{o}}(x) \, \md x,
\end{equation}
where
\begin{equation}
(R_{o}(z))^{1/2} \! := \! \left(\prod_{k=1}^{N+1}(z \! - \! b_{k-1}^{o})(z \!
- \! a_{k}^{o}) \right)^{1/2},
\end{equation}
with $(R_{o}(x))^{1/2}_{\pm} \! := \! \lim_{\varepsilon \downarrow 0}(R_{o}(x
\! \pm \! \mi \varepsilon))^{1/2}$, $x \! \in \! J_{o} \! := \! \operatorname{
supp}(\mu_{V}^{o}) \! = \! \cup_{j=1}^{N+1}(b_{j-1}^{o},a_{j}^{o})$ $(\subset
\mathbb{R} \setminus \{0\})$, $N \! \in \! \mathbb{N}$ (and finite), $b_{0}^{
o} \! := \! \min \lbrace \operatorname{supp}(\mu_{V}^{o}) \rbrace \! \notin \!
\lbrace -\infty,0 \rbrace$, $a_{N+1}^{o} \! := \! \max \lbrace \operatorname{
supp}(\mu_{V}^{o}) \rbrace \! \notin \! \lbrace 0,+\infty \rbrace$, and
$-\infty \! < \! b_{0}^{o} \! < \! a_{1}^{o} \! < \! b_{1}^{o} \! < \! a_{
2}^{o} \! < \! \cdots \! < \! b_{N}^{o} \! < \! a_{N+1}^{o} \! < \! +\infty$,
and the branch of the square root is chosen so that $z^{-(N+1)}(R_{o}(z))^{
1/2} \! \sim_{\underset{z \in \mathbb{C}_{\pm}}{z \to \infty}} \! \pm 1$,
\begin{equation}
h_{V}^{o}(z) \! := \! \dfrac{1}{2} \! \left(2 \! + \! \dfrac{1}{n} \right)^{-
1} \oint_{C_{\mathrm{R}}^{o}} \dfrac{\left(\frac{2 \mi}{\pi s} \! + \! \frac{
\mi \widetilde{V}^{\prime}(s)}{\pi} \right)}{\sqrt{\smash[b]{R_{o}(s)}} \, (s
\! - \! z)} \, \md s
\end{equation}
(real analytic for $z \! \in \! \mathbb{R} \setminus \{0\})$, $C_{\mathrm{R}
}^{o}$ $(\subset \mathbb{C}^{\ast})$ is the boundary of any open
doubly-connected annular region of the type $\lbrace \mathstrut z^{\prime} \!
\in \! \mathbb{C}; \, 0 \! < \! r^{\flat} \! < \! \vert z^{\prime} \vert \! <
\! R^{\flat} \! < \! +\infty \rbrace$, where the simple outer (resp., inner)
boundary $\lbrace \mathstrut z^{\prime} \! = \! R^{\flat} \me^{\mi \vartheta},
\, 0 \! \leqslant \! \vartheta \! \leqslant \! 2 \pi \rbrace$ (resp.,
$\lbrace \mathstrut z^{\prime} \! = \! r^{\flat} \me^{\mi \vartheta}, \, 0 \!
\leqslant \! \vartheta \! \leqslant \! 2 \pi \rbrace)$ is traversed clockwise
(resp., counter-clockwise), with the numbers $0 \! < \! r^{\flat} \! < \! R^{
\flat} \! < \! +\infty$ chosen such that, for (any) non-real $z$ in the domain
of analyticity of $\widetilde{V}$ (that is, $\mathbb{C}^{\ast})$, $\mathrm{
int}(C_{\mathrm{R}}^{o}) \! \supset \! J_{o} \cup \{z\}$, $\pmb{1}_{J_{o}}
(x)$ denotes the indicator (characteristic) function of the set $J_{o}$, and
$\lbrace b_{j-1}^{o},a_{j}^{o} \rbrace_{j=1}^{N+1}$ satisfy the following
$n$-dependent and (locally) solvable system of $2(N \! + \! 1)$ moment
conditions:
\begin{gather}
\begin{split}
\int_{J_{o}} \dfrac{(2s^{-1} \! + \! \widetilde{V}^{\prime}(s))s^{j}}{(R_{o}
(s))^{1/2}_{+}} \, \md s \! = \! 0, \quad j \! = \! 0,\dotsc,N, \qquad \qquad
\int_{J_{o}} \dfrac{(2s^{-1} \! + \! \widetilde{V}^{\prime}(s))s^{N+1}}{(R_{o}
(s))^{1/2}_{+}} \, \dfrac{\md s}{2 \pi \mi} \! = \! -\left(2 \! + \! \dfrac{
1}{n} \right), \\
\int_{a_{j}^{o}}^{b_{j}^{o}} \! \left(\dfrac{\mi (R_{o}(s))^{1/2}}{2 \pi}
\int_{J_{o}} \dfrac{(2 \xi^{-1} \! + \! \widetilde{V}^{\prime}(\xi))}{(R_{o}
(\xi))^{1/2}_{+}(\xi \! - \! s)} \, \md \xi \right) \! \md s \! = \! \ln \!
\left\vert \dfrac{a_{j}^{o}}{b_{j}^{o}} \right\vert \! + \! \dfrac{1}{2} \!
\left(\widetilde{V}(a_{j}^{o}) \! - \! \widetilde{V}(b_{j}^{o}) \right), \quad
j \! = \! 1,\dotsc,N.
\end{split}
\end{gather}
Suppose, furthermore, that $\widetilde{V} \colon \mathbb{R} \setminus \{0\} \!
\to \! \mathbb{R}$ is regular, namely:
\begin{description}
\item[{\rm (i)}] $h_{V}^{o}(x) \! \not\equiv \! 0$ on $\overline{J_{o}}:= \!
J_{o} \cup \! \left(\cup_{k=1}^{N+1}\{b_{k-1}^{o},a_{k}^{o}\} \right);$
\item[{\rm (ii)}]
\begin{equation}
2 \! \left(2 \! + \! \dfrac{1}{n} \right) \! \int_{J_{o}} \ln (\vert x \! - \!
s \vert) \, \md \mu_{V}^{o}(s) \! - \! 2 \ln \vert x \vert \! - \! \widetilde{
V}(x) \! - \! \ell_{o} \! - \! 2 \! \left(2 \! + \! \dfrac{1}{n} \right) \!
Q_{o} \! = \! 0, \quad x \! \in \! \overline{J_{o}},
\end{equation}
which defines the `odd' variational constant $\ell_{o} \! \in \! \mathbb{R}$
(the same on each---compact---interval $[b_{j-1}^{o},a_{j}^{o}]$, $j \! = \! 
1,\dotsc,N \! + \! 1)$, where
\begin{equation}
Q_{o} \! := \! \int_{J_{o}} \ln (\lvert s \rvert) \, \md \mu_{V}^{o}(s),
\end{equation}
and
\begin{equation*}
2 \! \left(2 \! + \! \dfrac{1}{n} \right) \! \int_{J_{o}} \ln (\vert x \! - \!
s \vert) \, \md \mu_{V}^{o}(s) \! - \! 2 \ln \vert x \vert \! - \! \widetilde{
V}(x) \! - \! \ell_{o} \! - \! 2 \! \left(2 \! + \! \dfrac{1}{n} \right) \!
Q_{o} \! < \! 0, \quad x \! \in \! \mathbb{R} \setminus \overline{J_{o}} \,;
\end{equation*}
\item[{\rm (iii)}]
\begin{equation*}
g^{o}_{+}(z) \! + \! g^{o}_{-}(z) \! - \! \widetilde{V}(z) \! - \! \ell_{o} \!
- \! \left(\mathfrak{Q}^{+}_{\mathscr{A}} \! + \! \mathfrak{Q}^{-}_{\mathscr{
A}} \right) \! < \! 0, \quad z \! \in \! \mathbb{R} \setminus \overline{J_{o}},
\end{equation*}
where
\begin{equation}
g^{o}(z) \! := \! \int_{J_{o}} \ln \! \left((z \! - \! s)^{2+\frac{1}{n}}
(zs)^{-1} \right) \md \mu_{V}^{o}(s), \quad z \! \in \! \mathbb{C} \setminus
(-\infty,\max \{0,a_{N+1}^{o}\}),
\end{equation}
and
\begin{equation*}
\mathfrak{Q}^{\pm}_{\mathscr{A}} \! := \! \left(1 \! + \! \dfrac{1}{n} \right)
\! \int_{J_{o}} \ln (\lvert s \rvert) \psi_{V}^{o}(s) \, \md s \! - \! \mi \pi
\int_{J_{o} \cap \mathbb{R}_{-}} \psi_{V}^{o}(s) \, \md s \! \pm \! \mi \pi
\! \left(2 \! + \! \dfrac{1}{n} \right) \! \int_{J_{o} \cap \mathbb{R}_{+}}
\psi_{V}^{o}(s) \, \md s,
\end{equation*}
with
\begin{align*}
\int_{J_{o} \cap \mathbb{R}_{-}} \psi_{V}^{o}(s) \, \md s =& \,
\begin{cases}
0, &\text{$J_{o} \! \subset \! \mathbb{R}_{+}$,} \\
1, &\text{$J_{o} \! \subset \! \mathbb{R}_{-}$,} \\
\int_{b_{0}^{o}}^{a_{j}^{o}} \psi_{V}^{o}(s) \, \md s, &\text{$(a_{j}^{o},
b_{j}^{o}) \! \ni \! 0, \quad j \! = \! 1,\dotsc,N$,}
\end{cases} \\
\int_{J_{o} \cap \mathbb{R}_{+}} \psi_{V}^{o}(s) \, \md s =& \,
\begin{cases}
0, &\text{$J_{o} \! \subset \! \mathbb{R}_{-}$,} \\
1, &\text{$J_{o} \! \subset \! \mathbb{R}_{+}$,} \\
\int_{b_{j}^{o}}^{a_{N+1}^{o}} \psi_{V}^{o}(s) \, \md s, &\text{$(a_{j}^{o},
b_{j}^{o}) \! \ni \! 0, \quad j \! = \! 1,\dotsc,N$;}
\end{cases}
\end{align*}
\item[{\rm (iv)}]
\begin{equation*}
\mi \! \left(g^{o}_{+}(z) \! - \! g^{o}_{-}(z) \! - \! \mathfrak{Q}^{+}_{
\mathscr{A}} \! + \! \mathfrak{Q}^{-}_{\mathscr{A}} \right)^{\prime} \! > \!
0, \quad z \! \in \! J_{o}.
\end{equation*}
\end{description}

Set
\begin{equation}
\overset{o}{m}^{\raise-1.0ex\hbox{$\scriptstyle \infty$}}(z) \! = \!
\begin{cases}
\overset{o}{\mathfrak{M}}^{\raise-1.0ex\hbox{$\scriptstyle \infty$}}(z),
&\text{$z \! \in \! \mathbb{C}_{+}$,} \\
-\mi \, \overset{o}{\mathfrak{M}}^{\raise-1.0ex\hbox{$\scriptstyle \infty$}}
(z) \sigma_{2}, &\text{$z \in \! \mathbb{C}_{-}$,}
\end{cases}
\end{equation}
where $(\det (\overset{o}{m}^{\raise-1.0ex\hbox{$\scriptstyle \infty$}}(z))
\! = \! 1)$
\begin{gather}
\overset{o}{\mathfrak{M}}^{\raise-1.0ex\hbox{$\scriptstyle \infty$}}(z) \!
= \!
\begin{pmatrix}
\frac{((\gamma^{o}(0))^{-1} \gamma^{o}(z)+\gamma^{o}(0)(\gamma^{o}(z))^{-1})}{
2} \mathfrak{m}^{o}_{11}(z) & -\frac{((\gamma^{o}(0))^{-1} \gamma^{o}(z)-
\gamma^{o}(0)(\gamma^{o}(z))^{-1})}{2 \mi} \mathfrak{m}^{o}_{12}(z) \\
\frac{((\gamma^{o}(0))^{-1} \gamma^{o}(z)-\gamma^{o}(0)(\gamma^{o}(z))^{-1})}{
2 \mi} \mathfrak{m}^{o}_{21}(z) & \frac{((\gamma^{o}(0))^{-1} \gamma^{o}(z)+
\gamma^{o}(0)(\gamma^{o}(z))^{-1})}{2} \mathfrak{m}^{o}_{22}(z)
\end{pmatrix}, \\
\gamma^{o}(z) \! := \! \left(\! \left(\dfrac{z \! - \! b_{0}^{o}}{z \! - \!
a_{N+1}^{o}} \right) \! \prod_{k=1}^{N} \! \left(\dfrac{z \! - \! b_{k}^{o}}{z
\! - \! a_{k}^{o}} \right) \right)^{1/4}, \qquad \qquad \gamma^{o}(0) \! := \!
\left(\prod_{k=1}^{N+1} \dfrac{b_{k-1}^{o}}{a_{k}^{o}} \right)^{1/4} \quad
(> \! 0), \\
\mathfrak{m}^{o}_{11}(z) \! := \! \dfrac{1}{\mathbb{E}} \dfrac{\boldsymbol{
\theta}^{o}(\boldsymbol{u}^{o}_{+}(0) \! + \! \boldsymbol{d}_{o}) \boldsymbol{
\theta}^{o}(\boldsymbol{u}^{o}(z) \! - \! \frac{1}{2 \pi}(n \! + \! \frac{1}{
2}) \boldsymbol{\Omega}^{o} \! + \! \boldsymbol{d}_{o})}{\boldsymbol{\theta}^{
o}(\boldsymbol{u}^{o}_{+}(0) \! - \! \frac{1}{2 \pi}(n \! + \! \frac{1}{2})
\boldsymbol{\Omega}^{o} \! + \! \boldsymbol{d}_{o}) \boldsymbol{\theta}^{o}
(\boldsymbol{u}^{o}(z) \! + \! \boldsymbol{d}_{o})}, \\
\mathfrak{m}^{o}_{12}(z) \! := \! \dfrac{1}{\mathbb{E}} \dfrac{\boldsymbol{
\theta}^{o}(\boldsymbol{u}^{o}_{+}(0) \! + \! \boldsymbol{d}_{o}) \boldsymbol{
\theta}^{o}(-\boldsymbol{u}^{o}(z) \! - \! \frac{1}{2 \pi}(n \! + \! \frac{1}{
2}) \boldsymbol{\Omega}^{o} \! + \! \boldsymbol{d}_{o})}{\boldsymbol{\theta}^{
o}(\boldsymbol{u}^{o}_{+}(0) \! - \! \frac{1}{2 \pi}(n \! + \! \frac{1}{2})
\boldsymbol{\Omega}^{o} \! + \! \boldsymbol{d}_{o}) \boldsymbol{\theta}^{o}
(-\boldsymbol{u}^{o}(z) \! + \! \boldsymbol{d}_{o})}, \\
\mathfrak{m}^{o}_{21}(z) \! := \! \mathbb{E} \dfrac{\boldsymbol{\theta}^{o}
(\boldsymbol{u}^{o}_{+}(0) \! + \! \boldsymbol{d}_{o}) \boldsymbol{\theta}^{o}
(\boldsymbol{u}^{o}(z) \! - \! \frac{1}{2 \pi}(n \! + \! \frac{1}{2})
\boldsymbol{\Omega}^{o} \! - \! \boldsymbol{d}_{o})}{\boldsymbol{\theta}^{o}
(-\boldsymbol{u}^{o}_{+}(0) \! - \! \frac{1}{2 \pi}(n \! + \! \frac{1}{2})
\boldsymbol{\Omega}^{o} \! - \! \boldsymbol{d}_{o}) \boldsymbol{\theta}^{o}
(\boldsymbol{u}^{o}(z) \! - \! \boldsymbol{d}_{o})}, \\
\mathfrak{m}^{o}_{22}(z) \! := \! \mathbb{E} \dfrac{\boldsymbol{\theta}^{o}
(\boldsymbol{u}^{o}_{+}(0) \! + \! \boldsymbol{d}_{o}) \boldsymbol{\theta}^{o}
(-\boldsymbol{u}^{o}(z) \! - \! \frac{1}{2 \pi}(n \! + \! \frac{1}{2})
\boldsymbol{\Omega}^{o} \! - \! \boldsymbol{d}_{o})}{\boldsymbol{\theta}^{o}
(-\boldsymbol{u}^{o}_{+}(0) \! - \! \frac{1}{2 \pi}(n \! + \! \frac{1}{2})
\boldsymbol{\Omega}^{o} \! - \! \boldsymbol{d}_{o}) \boldsymbol{\theta}^{o}
(\boldsymbol{u}^{o}(z) \! + \! \boldsymbol{d}_{o})},
\end{gather}
with
\begin{gather*}
\mathbb{E} \! = \! \exp \! \left(\mi 2 \pi \! \left(n \! + \! \dfrac{1}{2}
\right) \! \int_{J_{o} \cap \mathbb{R}_{+}} \psi_{V}^{o}(s) \, \md s \right),
\\
\boldsymbol{u}^{o}(z) \! = \! \int_{a_{N+1}^{o}}^{z} \boldsymbol{\omega}^{o},
\qquad \qquad \boldsymbol{u}^{o}_{+}(0) \! = \! \int_{a_{N+1}^{o}}^{0^{+}}
\boldsymbol{\omega}^{o},
\end{gather*}
$\boldsymbol{\Omega}^{o} \! = \! (\Omega^{o}_{1},\Omega^{o}_{2},\dotsc,
\Omega^{o}_{N})^{\operatorname{T}}$ $(\in \! \mathbb{R}^{N})$, where
\begin{equation*}
\Omega^{o}_{j} \! := \! 4 \pi \int_{b_{j}^{o}}^{a_{N+1}^{o}} \psi_{V}^{o}(s)
\, \md s, \quad j \! = \! 1,\dotsc,N,
\end{equation*}
and
\begin{equation*}
\boldsymbol{d}_{o} \! \equiv \! \sum_{j=1}^{N} \int_{a_{j}^{o}}^{z_{j}^{o,+}}
\boldsymbol{\omega}^{o} \quad \left(= \! -\sum_{j=1}^{N+1} \int_{a_{j}^{o}}^{
z_{j}^{o,-}} \boldsymbol{\omega}^{o} \right),
\end{equation*}
where a set of $N$ upper-edge and lower-edge finite-length-gap roots/zeros
are
\begin{equation*}
\left\lbrace z_{j}^{o,\pm} \right\rbrace_{j=1}^{N} \! = \! \left\lbrace
\mathstrut z^{\pm} \! \in \! \mathbb{C}_{\pm}; \, ((\gamma^{o}(0))^{-1}
\gamma^{o}(z) \! \mp \! \gamma^{o}(0)(\gamma^{o}(z))^{-1}) \vert_{z=z^{\pm}}
\! = \! 0 \right\rbrace,
\end{equation*}
with $z_{j}^{o,\pm} \! \in \! (a_{j}^{o},b_{j}^{o})^{\pm}$ $(\subset \!
\mathbb{C}_{\pm})$, $j \! = \! 1,\dotsc,N$.

Let $\overset{o}{\operatorname{Y}} \colon \mathbb{C} \setminus \mathbb{R}
\! \to \! \operatorname{SL}_{2}(\mathbb{C})$ be the unique solution of
{\rm \pmb{RHP2}} whose integral representations are given in Lemma
{\rm 2.2.2;} in particular, $z \boldsymbol{\pi}_{2n+1}(z) \! := \! (\overset{
o}{\operatorname{Y}}(z))_{11}$. Then:\\
{\rm \pmb{(1)}} for $z \! \in \! \Upsilon^{o}_{1}$ $(\subset \! \mathbb{C}_{
+})$,
\begin{align}
z \boldsymbol{\pi}_{2n+1}(z) \underset{n \to \infty}{=}& \, \mathbb{E} \exp
\! \left(n(g^{o}(z) \! - \! \mathfrak{Q}^{+}_{\mathscr{A}}) \right) \! \left(
(\overset{o}{m}^{\raise-1.0ex\hbox{$\scriptstyle \infty$}}(z))_{11} \! \left(
1 \! + \! \dfrac{1}{n \! + \! \frac{1}{2}} \! \left(\mathscr{R}^{o}_{0}(z)
\right)_{11} \right. \right. \nonumber \\
+&\left. \left. \, \mathcal{O} \! \left(\dfrac{1}{(n \! + \! \frac{1}{2})^{2}}
\right) \right) \! + \!
(\overset{o}{m}^{\raise-1.0ex\hbox{$\scriptstyle \infty$}}(z))_{21} \! \left(
\dfrac{1}{n \! + \! \frac{1}{2}} \! \left(\mathscr{R}^{o}_{0}(z) \right)_{12}
\! + \! \mathcal{O} \! \left(\dfrac{1}{(n \! + \! \frac{1}{2})^{2}} \right)
\right) \right),
\end{align}
and
\begin{align}
z \int_{\mathbb{R}} \dfrac{(s \boldsymbol{\pi}_{2n+1}(s)) \me^{-n \widetilde{V}
(s)}}{s(s \! - \! z)} \, \dfrac{\md s}{2 \pi \mi} \underset{n \to \infty}{=}&
\, \dfrac{1}{\mathbb{E}} \exp \! \left(-n(g^{o}(z) \! - \! \ell_{o} \! - \!
\mathfrak{Q}^{+}_{\mathscr{A}}) \right) \nonumber \\
\times& \, \left((\overset{o}{m}^{\raise-1.0ex\hbox{$\scriptstyle \infty$}}
(z))_{12} \! \left(1 \! + \! \dfrac{1}{n \! + \! \frac{1}{2}} \! \left(
\mathscr{R}^{o}_{0}(z) \right)_{11} \! + \! \mathcal{O} \! \left(\dfrac{1}{(n
\! + \! \frac{1}{2})^{2}} \right) \right) \right. \nonumber \\
+&\left. \, (\overset{o}{m}^{\raise-1.0ex\hbox{$\scriptstyle \infty$}}(z))_{2
2} \! \left(\dfrac{1}{n \! + \! \frac{1}{2}} \! \left(\mathscr{R}^{o}_{0}(z)
\right)_{12} \! + \! \mathcal{O} \! \left(\dfrac{1}{(n \! + \! \frac{1}{2})^{
2}} \right) \right) \right),
\end{align}
where
\begin{align}
\mathscr{R}^{o}_{0}(z) :=& \, \sum_{j=1}^{N+1} \! \left(\dfrac{(\mathscr{B}^{o}
(a_{j}^{o}) \widehat{\alpha}_{0}^{o}(a_{j}^{o}) \! - \! \mathscr{A}^{o}(a_{j}^{
o})(\widehat{\alpha}_{1}^{o}(a_{j}^{o}) \! + \! (a_{j}^{o})^{-1} \widehat{
\alpha}_{0}^{o}(a_{j}^{o})))}{(\widehat{\alpha}_{0}^{o}(a_{j}^{o}))^{2}a_{j}^{
o}} \right. \nonumber \\
+&\left. \, \dfrac{(\mathscr{B}^{o}(b_{j-1}^{o}) \widehat{\alpha}_{0}^{o}(b_{j
-1}^{o}) \! - \! \mathscr{A}^{o}(b_{j-1}^{o})(\widehat{\alpha}_{1}^{o}(b_{j-
1}^{o}) \! + \! (b_{j-1}^{o})^{-1} \widehat{\alpha}_{0}^{o}(b_{j-1}^{o})))}{
(\widehat{\alpha}_{0}^{o}(b_{j-1}^{o}))^{2}b_{j-1}^{o}} \right. \nonumber \\
+&\left. \, \dfrac{1}{(z \! - \! b_{j-1}^{o})} \! \left(\dfrac{\mathscr{A}^{o}
(b_{j-1}^{o})}{\widehat{\alpha}_{0}^{o}(b_{j-1}^{o})(z \! - \! b_{j-1}^{o})}
\! + \! \dfrac{(\mathscr{B}^{o}(b_{j-1}^{o}) \widehat{\alpha}_{0}^{o}(b_{j-1}^{
o}) \! - \! \mathscr{A}^{o}(b_{j-1}^{o}) \widehat{\alpha}_{1}^{o}(b_{j-1}^{o})
)}{(\widehat{\alpha}_{0}^{o}(b_{j-1}^{o}))^{2}} \right) \right. \nonumber \\
+&\left. \, \dfrac{1}{(z \! - \! a_{j}^{o})} \! \left(\dfrac{\mathscr{A}^{o}
(a_{j}^{o})}{\widehat{\alpha}_{0}^{o}(a_{j}^{o})(z \! - \! a_{j}^{o})} \! + \!
\dfrac{(\mathscr{B}^{o}(a_{j}^{o}) \widehat{\alpha}_{0}^{o}(a_{j}^{o}) \! - \!
\mathscr{A}^{o}(a_{j}^{o}) \widehat{\alpha}_{1}^{o}(a_{j}^{o}))}{(\widehat{
\alpha}_{0}^{o}(a_{j}^{o}))^{2}} \right) \right),
\end{align}
with, for $j \! = \! 1,\dotsc,N \! + \! 1$,
\begin{gather}
\mathscr{A}^{o}(b_{j-1}^{o}) \! = \! -\dfrac{s_{1}(\gamma^{o}(0))^{2} \me^{
\mi (n+\frac{1}{2}) \mho_{j-1}^{o}}}{Q_{0}^{o}(b_{j-1}^{o})} \!
\begin{pmatrix}
\varkappa^{o}_{1}(b_{j-1}^{o}) \varkappa^{o}_{2}(b_{j-1}^{o}) & \mi 
(\varkappa^{o}_{1}(b_{j-1}^{o}))^{2} \\
\mi (\varkappa^{o}_{2}(b_{j-1}^{o}))^{2} & -\varkappa^{o}_{1}(b_{j-1}^{o}) 
\varkappa^{o}_{2}(b_{j-1}^{o})
\end{pmatrix}, \\
\mathscr{A}^{o}(a_{j}^{o}) \! = \! \dfrac{s_{1}Q_{0}^{o}(a_{j}^{o}) \me^{
\mi (n+\frac{1}{2}) \mho_{j}^{o}}}{(\gamma^{o}(0))^{2}} \!
\begin{pmatrix}
-\varkappa^{o}_{1}(a_{j}^{o}) \varkappa^{o}_{2}(a_{j}^{o}) & \mi 
(\varkappa^{o}_{1}(a_{j}^{o}))^{2} \\
\mi (\varkappa^{o}_{2}(a_{j}^{o}))^{2} & \varkappa^{o}_{1}(a_{j}^{o}) 
\varkappa^{o}_{2}(a_{j}^{o})
\end{pmatrix}, \\
\dfrac{\mathscr{B}^{o}(b_{j-1}^{o})}{\me^{\mi (n+\frac{1}{2}) \mho_{j-1}^{
o}}} \! = \!
\begin{pmatrix}
\boxed{\begin{matrix} \varkappa_{1}^{o}(b_{j-1}^{o}) \varkappa_{2}^{o}
(b_{j-1}^{o}) \! \left(-\frac{s_{1}(\gamma^{o}(0))^{2}}{Q_{0}^{o}(b_{j-1}^{
o})} \right. \\
\left. \times \left\{\daleth^{1}_{1}(b_{j-1}^{o}) \! + \! \daleth^{1}_{-1}
(b_{j-1}^{o}) \! - \! Q_{1}^{o}(b_{j-1}^{o}) \right. \right. \\
\left. \left. \times \, (Q_{0}^{o}(b_{j-1}^{o}))^{-1} \right\} \! - \! t_{1}
(\gamma^{o}(0))^{2} \! \left\{Q_{0}^{o}(b_{j-1}^{o}) \right. \right. \\
\left. \left. + \, (Q_{0}^{o}(b_{j-1}^{o}))^{-1} \aleph^{1}_{1}(b_{j-1}^{o}) 
\aleph^{1}_{-1}(b_{j-1}^{o}) \right\} \right. \\
\left. + \, \mi (s_{1} \! + \! t_{1}) \! \left\{\aleph^{1}_{-1}(b_{j-1}^{o}) 
\! - \! \aleph^{1}_{1}(b_{j-1}^{o}) \right\} \right)
\end{matrix}} & 
\boxed{\begin{matrix} (\varkappa_{1}^{o}(b_{j-1}^{o}))^{2} \! \left(-\frac{
\mi s_{1}(\gamma^{o}(0))^{2}}{Q_{0}^{o}(b_{j-1}^{o})} \! \left\{2 \daleth^{
1}_{1}(b_{j-1}^{o}) \right. \right. \\
\left. \left. -\, Q_{1}^{o}(b_{j-1}^{o})(Q_{0}^{o}(b_{j-1}^{o}))^{-1} 
\right\} \! + \! \mi t_{1} \! \left\{Q_{0}^{o}(b_{j-1}^{o}) \right. \right. \\
\left. \left. \times \, (\gamma^{o}(0))^{-2} \! - \! (Q_{0}^{o}(b_{j-1}^{o}
))^{-1}(\gamma^{o}(0))^{2} \right. \right. \\
\left. \left. \times \, (\aleph^{1}_{1}(b_{j-1}^{o}))^{2} \right\} \! + \! 
2(s_{1} \! - \! t_{1}) \aleph^{1}_{1}(b_{j-1}^{o}) \right)
\end{matrix}} \\
\boxed{\begin{matrix} (\varkappa_{2}^{o}(b_{j-1}^{o}))^{2} \! \left(-\frac{
\mi s_{1}(\gamma^{o}(0))^{2}}{Q_{0}^{o}(b_{j-1}^{o})} \! \left\{2 \daleth^{
1}_{-1}(b_{j-1}^{o}) \right. \right. \\
\left. \left. -\, Q_{1}^{o}(b_{j-1}^{o})(Q_{0}^{o}(b_{j-1}^{o}))^{-1} 
\right\} \! + \! \mi t_{1} \! \left\{Q_{0}^{o}(b_{j-1}^{o}) \right. \right. \\
\left. \left. \times \, (\gamma^{o}(0))^{-2} \! - \! (Q_{0}^{o}(b_{j-1}^{o}
))^{-1}(\gamma^{o}(0))^{2} \right. \right. \\
\left. \left. \times \, (\aleph^{1}_{-1}(b_{j-1}^{o}))^{2} \right\} \!- \! 
2(s_{1} \! - \! t_{1}) \aleph^{1}_{-1}(b_{j-1}^{o}) \right)
\end{matrix}} & 
\boxed{\begin{matrix} \varkappa_{1}^{o}(b_{j-1}^{o}) \varkappa_{2}^{o}
(b_{j-1}^{o}) \! \left(\frac{s_{1}(\gamma^{o}(0))^{2}}{Q_{0}^{o}(b_{j-1}^{
o})} \right. \\
\left. \times \left\{\daleth^{1}_{1}(b_{j-1}^{o}) \! + \! \daleth^{1}_{-1}
(b_{j-1}^{o}) \! - \! Q_{1}^{o}(b_{j-1}^{o}) \right. \right. \\
\left. \left. \times \, (Q_{0}^{o}(b_{j-1}^{o}))^{-1} \right\} \! + \! t_{1}
(\gamma^{o}(0))^{2} \! \left\{Q_{0}^{o}(b_{j-1}^{o}) \right. \right. \\
\left. \left. + \, (Q_{0}^{o}(b_{j-1}^{o}))^{-1} \aleph^{1}_{1}(b_{j-1}^{o}) 
\aleph^{1}_{-1}(b_{j-1}^{o}) \right\} \right. \\
\left. + \, \mi (s_{1} \! + \! t_{1}) \! \left\{\aleph^{1}_{1}(b_{j-1}^{o}) 
\! - \! \aleph^{1}_{-1}(b_{j-1}^{o}) \right\} \right)
\end{matrix}}
\end{pmatrix}, \\
\dfrac{\mathscr{B}^{o}(a_{j}^{o})}{\me^{\mi (n+\frac{1}{2}) \mho_{j}^{o}}} \! 
= \!
\begin{pmatrix}
\boxed{\begin{matrix} \varkappa_{1}^{o}(a_{j}^{o}) \varkappa_{2}^{o}(a_{j}^{
o}) \! \left(-\frac{s_{1}}{(\gamma^{o}(0))^{2}} \! \left\{Q_{1}^{o}(a_{j}^{
o}) \right. \right. \\
\left. \left. + \, Q_{0}^{o}(a_{j}^{o}) \left[\daleth^{1}_{1}(a_{j}^{o}) \! 
+ \! \daleth^{1}_{-1}(a_{j}^{o}) \right] \right\} \! - \! t_{1} \right. \\
\left. \times \left\{(\gamma^{o}(0))^{-2}Q_{0}^{o}(a_{j}^{o}) \aleph^{1}_{1}
(a_{j}^{o}) \aleph^{1}_{-1}(a_{j}^{o}) \right. \right. \\
\left. \left. + \, (\gamma^{o}(0))^{2}(Q_{0}^{o}(a_{j}^{o}))^{-1} \right\} 
\right. \\
\left. + \, \mi (s_{1} \! + \! t_{1}) \! \left\{\aleph^{1}_{-1}(a_{j}^{o}) 
\! - \! \aleph^{1}_{1}(a_{j}^{o}) \right\} \right)
\end{matrix}} & 
\boxed{\begin{matrix} (\varkappa_{1}^{o}(a_{j}^{o}))^{2} \! \left(\frac{\mi 
s_{1}}{(\gamma^{o}(0))^{2}} \! \left\{Q_{1}^{o}(a_{j}^{o}) \! + \! 2Q_{0}^{o}
(a_{j}^{o}) \right. \right. \\
\left. \left. \times \, \daleth^{1}_{1}(a_{j}^{o}) \right\} \! + \! \mi t_{1} 
\! \left\{Q_{0}^{o}(a_{j}^{o})(\aleph^{1}_{1}(a_{j}^{o}))^{2} \right. \right. 
\\
\left. \left. \times \, (\gamma^{o}(0))^{-2} \! - \! (Q_{0}^{o}(a_{j}^{o}))^{
-1}(\gamma^{o}(0))^{2} \right\} \right. \\
\left. - \, 2(s_{1} \! - \! t_{1}) \aleph^{1}_{1}(a_{j}^{o}) \right)
\end{matrix}} \\
\boxed{\begin{matrix} (\varkappa_{2}^{o}(a_{j}^{o}))^{2} \! \left(\frac{\mi 
s_{1}}{(\gamma^{o}(0))^{2}} \! \left\{Q_{1}^{o}(a_{j}^{o}) \! + \! 2Q_{0}^{o}
(a_{j}^{o}) \right. \right. \\
\left. \left. \times \, \daleth^{1}_{-1}(a_{j}^{o}) \right\} \! + \! \mi t_{1} 
\! \left\{Q_{0}^{o}(a_{j}^{o})(\aleph^{1}_{-1}(a_{j}^{o}))^{2} \right. 
\right. \\
\left. \left. \times \, (\gamma^{o}(0))^{-2} \! - \! (Q_{0}^{o}(a_{j}^{o}))^{
-1}(\gamma^{o}(0))^{2} \right\} \right. \\
\left. + \, 2(s_{1} \! - \! t_{1}) \aleph^{1}_{-1}(a_{j}^{o}) \right)
\end{matrix}} & 
\boxed{\begin{matrix} \varkappa_{1}^{o}(a_{j}^{o}) \varkappa_{2}^{o}(a_{j}^{
o}) \! \left(\frac{s_{1}}{(\gamma^{o}(0))^{2}} \! \left\{Q_{1}^{o}(a_{j}^{o}) 
\right. \right. \\
\left. \left. + \, Q_{0}^{o}(a_{j}^{o}) \left[\daleth^{1}_{1}(a_{j}^{o}) \! + 
\! \daleth^{1}_{-1}(a_{j}^{o}) \right] \right\} \! + \! t_{1} \right. \\
\left. \times \left\{(\gamma^{o}(0))^{-2}Q_{0}^{o}(a_{j}^{o}) \aleph^{1}_{1}
(a_{j}^{o}) \aleph^{1}_{-1}(a_{j}^{o}) \right. \right. \\
\left. \left. + \, (\gamma^{o}(0))^{2}(Q_{0}^{o}(a_{j}^{o}))^{-1} \right\} 
\right. \\
\left. + \, \mi (s_{1} \! + \! t_{1}) \! \left\{\aleph^{1}_{1}(a_{j}^{o}) 
\! - \! \aleph^{1}_{-1}(a_{j}^{o}) \right\} \right)
\end{matrix}}
\end{pmatrix}, \\
s_{1} \! = \! \dfrac{5}{72}, \qquad \qquad \quad t_{1} \! = \! -\dfrac{7}{72}, 
\qquad \qquad \quad \mho_{i}^{o} \! := \!
\begin{cases}
\Omega_{i}^{o}, &\text{$i \! = \! 1,\dotsc,N$,} \\
0, &\text{$i \! = \! 0,N \! + \! 1$,}
\end{cases} \\
Q_{0}^{o}(b_{0}^{o}) \! = \! -\mi \left(\! (a_{N+1}^{o} \! - \! b_{0}^{o})^{-
1} \prod_{k=1}^{N} \! \left(\dfrac{b_{k}^{o} \! - \! b_{0}^{o}}{a_{k}^{o} \! -
\! b_{0}^{o}} \right) \right)^{1/2}, \\
Q_{1}^{o}(b_{0}^{o}) \! = \! \dfrac{1}{2}Q_{0}^{o}(b_{0}^{o}) \! \left(\sum_{k
=1}^{N} \! \left(\dfrac{1}{b_{0}^{o} \! - \! b_{k}^{o}} \! - \! \dfrac{1}{b_{
0}^{o} \! - \! a_{k}^{o}} \right) \! - \! \dfrac{1}{b_{0}^{o} \! - \! a_{N+
1}^{o}} \right), \\
Q_{0}^{o}(a_{N+1}^{o}) \! = \! \left(\! (a_{N+1}^{o} \! - \! b_{0}^{o}) \prod_{
k=1}^{N} \! \left(\dfrac{a_{N+1}^{o} \! - \! b_{k}^{o}}{a_{N+1}^{o} \! - \!
a_{k}^{o}} \right) \right)^{1/2}, \\
Q_{1}^{o}(a_{N+1}^{o}) \! = \! \dfrac{1}{2}Q_{0}^{o}(a_{N+1}^{o}) \! \left(
\sum_{k=1}^{N} \! \left(\dfrac{1}{a_{N+1}^{o} \! - \! b_{k}^{o}} \! - \!
\dfrac{1}{a_{N+1}^{o} \! - \! a_{k}^{o}} \right) \! + \! \dfrac{1}{a_{N+1}^{
o} \! - \! b_{0}^{o}} \right), \\
Q_{0}^{o}(b_{j}^{o}) \! = \! -\mi \left(\! \dfrac{(b_{j}^{o} \! - \! b_{0}^{o})
}{(a_{N+1}^{o} \! - \! b_{j}^{o})(b_{j}^{o} \! - \! a_{j}^{o})} \prod_{k=1}^{j
-1} \! \left(\dfrac{b_{j}^{o} \! - \! b_{k}^{o}}{b_{j}^{o} \! - \! a_{k}^{o}}
\right) \! \prod_{l=j+1}^{N} \! \left(\! \dfrac{b_{l}^{o} \! - \! b_{j}^{o}}{
a_{l}^{o} \! - \! b_{j}^{o}} \right) \right)^{1/2}, \\
Q_{1}^{o}(b_{j}^{o}) \! = \! \dfrac{1}{2}Q_{0}^{o}(b_{j}^{o}) \! \left(\sum_{
\substack{k=1\\k \not= j}}^{N} \! \left(\! \dfrac{1}{b_{j}^{o} \! - \! b_{k}^{
o}} \! - \! \dfrac{1}{b_{j}^{o} \! - \! a_{k}^{o}} \right) \! + \! \dfrac{1}{
b_{j}^{o} \! - \! b_{0}^{o}} \! - \! \dfrac{1}{b_{j}^{o} \! - \! a_{N+1}^{o}}
\! - \! \dfrac{1}{b_{j}^{o} \! - \! a_{j}^{o}} \right), \\
Q_{0}^{o}(a_{j}^{o}) \! = \! \left(\! \dfrac{(a_{j}^{o} \! - \! b_{0}^{o})(b_{
j}^{o} \! - \! a_{j}^{o})}{(a_{N+1}^{o} \! - \! a_{j}^{o})} \prod_{k=1}^{j-1}
\! \left(\dfrac{a_{j}^{o} \! - \! b_{k}^{o}}{a_{j}^{o} \! - \! a_{k}^{o}}
\right) \! \prod_{l=j+1}^{N} \! \left(\! \dfrac{b_{l}^{o} \! - \! a_{j}^{o}}{
a_{l}^{o} \! - \! a_{j}^{o}} \right) \right)^{1/2}, \\
Q_{1}^{o}(a_{j}^{o}) \! = \! \dfrac{1}{2}Q_{0}^{o}(a_{j}^{o}) \! \left(\sum_{
\substack{k=1\\k \not= j}}^{N} \! \left(\! \dfrac{1}{a_{j}^{o} \! - \! b_{k}^{
o}} \! - \! \dfrac{1}{a_{j}^{o} \! - \! a_{k}^{o}} \right) \! + \! \dfrac{1}{
a_{j}^{o} \! - \! b_{0}^{o}} \! - \! \dfrac{1}{a_{j}^{o} \! - \! a_{N+1}^{o}}
\! + \! \dfrac{1}{a_{j}^{o} \! - \! b_{j}^{o}} \right),
\end{gather}
where $\mi Q_{0}^{o}(b_{j-1}^{o}),Q^{o}_{0}(a_{j}^{o}) \! > \! 0$, $j \! = \! 
1,\dotsc,N \! + \! 1$,
\begin{gather}
\varkappa_{1}^{o}(\xi) \! = \! \dfrac{1}{\mathbb{E}} \dfrac{\bm{\theta}^{o}
(\bm{u}^{o}_{+}(0) \! + \! \bm{d}_{o}) \bm{\theta}^{o}(\bm{u}^{o}_{+}(\xi) 
\! - \! \frac{1}{2 \pi}(n \! + \! \frac{1}{2}) \bm{\Omega}^{o} \! + \! \bm{
d}_{o})}{\bm{\theta}^{o}(\bm{u}^{o}_{+}(0) \! - \! \frac{1}{2 \pi}(n \! + \! 
\frac{1}{2}) \bm{\Omega}^{o} \! + \! \bm{d}_{o}) \bm{\theta}^{o}(\bm{u}^{o}_{
+}(\xi) \! + \! \bm{d}_{o})}, \\
\varkappa_{2}^{o}(\xi) \! = \! \mathbb{E} \dfrac{\bm{\theta}^{o}(-\bm{u}^{o}_{
+}(0) \! - \! \bm{d}_{o}) \bm{\theta}^{o}(\bm{u}^{o}_{+}(\xi) \! - \! \frac{
1}{2 \pi}(n \! + \! \frac{1}{2}) \bm{\Omega}^{o} \! - \! \bm{d}_{o})}{\bm{
\theta}^{o}(-\bm{u}^{o}_{+}(0) \! - \! \frac{1}{2 \pi}(n \! + \! \frac{1}{2}) 
\bm{\Omega}^{o} \! - \! \bm{d}_{o}) \bm{\theta}^{o}(\bm{u}^{o}_{+}(\xi) \! - 
\! \bm{d}_{o})}, \\
\aleph^{\varepsilon_{1}}_{\varepsilon_{2}}(\xi) \! = \! -\dfrac{\mathfrak{
u}^{o}(\varepsilon_{1},\varepsilon_{2},\bm{0};\xi)}{\bm{\theta}^{o}
(\varepsilon_{1} \bm{u}^{o}_{+}(\xi) \! + \! \varepsilon_{2} \bm{d}_{o})} \! 
+ \! \dfrac{\mathfrak{u}^{o}(\varepsilon_{1},\varepsilon_{2},\bm{\Omega}^{o};
\xi)}{\bm{\theta}^{o}(\varepsilon_{1} \bm{u}^{o}_{+}(\xi) \! - \! \frac{1}{2 
\pi}(n \! + \! \frac{1}{2}) \bm{\Omega}^{o} \! + \! \varepsilon_{2} \bm{d}_{
o})}, \quad \varepsilon_{1},\varepsilon_{2} \! = \! \pm 1,
\end{gather}
\begin{align}
\daleth^{\varepsilon_{1}}_{\varepsilon_{2}}(\xi) =& \, -\dfrac{\mathfrak{
v}^{o}(\varepsilon_{1},\varepsilon_{2},\bm{0};\xi)}{\bm{\theta}^{o}
(\varepsilon_{1} \bm{u}^{o}_{+}(\xi) \! + \! \varepsilon_{2} \bm{d}_{o})} \! 
+ \! \dfrac{\mathfrak{v}^{o}(\varepsilon_{1},\varepsilon_{2},\bm{\Omega}^{o};
\xi)}{\bm{\theta}^{o}(\varepsilon_{1} \bm{u}^{o}_{+}(\xi) \! - \! \frac{1}{2 
\pi}(n \! + \! \frac{1}{2}) \bm{\Omega}^{o} \! + \! \varepsilon_{2} \bm{d}_{
o})} \! - \! \left(\dfrac{\mathfrak{u}^{o}(\varepsilon_{1},\varepsilon_{2},
\bm{0};\xi)}{\bm{\theta}^{o}(\varepsilon_{1} \bm{u}^{o}_{+}(\xi) \! + \! 
\varepsilon_{2} \bm{d}_{o})} \right)^{2} \nonumber \\
+& \, \dfrac{\mathfrak{u}^{o}(\varepsilon_{1},\varepsilon_{2},\bm{0};\xi) 
\mathfrak{u}^{o}(\varepsilon_{1},\varepsilon_{2},\bm{\Omega}^{o};\xi)}{\bm{
\theta}^{o}(\varepsilon_{1} \bm{u}^{o}_{+}(\xi) \! + \! \varepsilon_{2} 
\bm{d}_{o}) \bm{\theta}^{o}(\varepsilon_{1} \bm{u}^{o}_{+}(\xi) \! - \! 
\frac{1}{2 \pi}(n \! + \! \frac{1}{2}) \bm{\Omega}^{o} \! + \! \varepsilon_{
2} \bm{d}_{o})},
\end{align}
\begin{gather}
\mathfrak{u}^{o}(\varepsilon_{1},\varepsilon_{2},\bm{\Omega}^{o},\xi) \! := 
\! 2 \pi \Lambda^{\raise-1.0ex\hbox{$\scriptstyle 1$}}_{o}(\varepsilon_{1},
\varepsilon_{2},\bm{\Omega}^{o},\xi), \qquad \mathfrak{v}^{o}(\varepsilon_{
1},\varepsilon_{2},\bm{\Omega}^{o},\xi) \! := \! -2 \pi^{2} 
\Lambda^{\raise-1.0ex\hbox{$\scriptstyle 2$}}_{o}(\varepsilon_{1},
\varepsilon_{2},\bm{\Omega}^{o},\xi), \\
\Lambda^{\raise-1.0ex\hbox{$\scriptstyle j^{\prime}$}}_{o}(\varepsilon_{1},
\varepsilon_{2},\bm{\Omega}^{o},\xi) \! = \! \sum_{m \in \mathbb{Z}^{N}}
(\mathfrak{r}_{o}(\xi))^{j^{\prime}} \me^{2 \pi \mi (m,\varepsilon_{1} 
\bm{u}^{o}_{+}(\xi)-\frac{1}{2 \pi}(n+\frac{1}{2}) \bm{\Omega}^{o}+
\varepsilon_{2} \bm{d}_{o})+ \pi \mi(m,\bm{\tau}^{o}m)}, \quad j^{\prime} \! 
= \! 1,2, \\
\mathfrak{r}_{o}(\xi) \! := \! \dfrac{2(m,\vec{\moo}_{o}(\xi))}{
\leftthreetimes^{\raise+0.3ex\hbox{$\scriptstyle o$}}(\xi)}, \qquad \qquad 
\vec{\moo}_{o}(\xi) \! = \! 
(\rightthreetimes^{\raise-0.9ex\hbox{$\scriptstyle o$}}_{1}(\xi),
\rightthreetimes^{\raise-0.9ex\hbox{$\scriptstyle o$}}_{2}(\xi),\dotsc,
\rightthreetimes^{\raise-0.9ex\hbox{$\scriptstyle o$}}_{N}(\xi)), \\
\rightthreetimes^{\raise-0.9ex\hbox{$\scriptstyle o$}}_{j^{\prime}}(\xi) \! 
:= \! \sum_{k=1}^{N}c_{j^{\prime}k}^{o} \xi^{N-k}, \quad j^{\prime} \! = \! 
1,\dotsc,N, \\
\leftthreetimes^{\raise+0.3ex\hbox{$\scriptstyle o$}}(b_{0}^{o}) \! = \! \mi 
(-1)^{N} \eta_{b_{0}^{o}}, \quad 
\leftthreetimes^{\raise+0.3ex\hbox{$\scriptstyle o$}}(a_{N+1}^{o}) \! = \! 
\eta_{a_{N+1}^{o}}, \quad 
\leftthreetimes^{\raise+0.3ex\hbox{$\scriptstyle o$}}(b_{j}^{o}) \! = \! \mi 
(-1)^{N-j} \eta_{b_{j}^{o}}, \quad 
\leftthreetimes^{\raise+0.3ex\hbox{$\scriptstyle o$}}(a_{j}^{o}) \! = \! 
(-1)^{N-j+1} \eta_{a_{j}^{o}}, \\
\eta_{b_{0}^{o}} \! := \! \left((a_{N+1}^{o} \! - \! b_{0}^{o}) \prod_{k=1}^{
N}(b_{k}^{o} \! - \! b_{0}^{o})(a_{k}^{o} \! - \! b_{0}^{o}) \right)^{1/2}, \\
\eta_{a_{N+1}^{o}} \! := \! \left((a_{N+1}^{o} \! - \! b_{0}^{o}) \prod_{k=
1}^{N}(a_{N+1}^{o} \! - \! b_{k}^{o})(a_{N+1}^{o} \! - \! a_{k}^{o}) 
\right)^{1/2}, \\
\eta_{b_{j}^{o}} \! := \! \left(\! (b_{j}^{o} \! - \! a_{j}^{o})(a_{N+1}^{o} 
\! - \! b_{j}^{o})(b_{j}^{o} \! - \! b_{0}^{o}) \prod_{k=1}^{j-1}(b_{j}^{o} 
\! - \! b_{k}^{o})(b_{j}^{o} \! - \! a_{k}^{o}) \prod_{l=j+1}^{N}(b_{l}^{o} 
\! - \! b_{j}^{o})(a_{l}^{o} \! - \! b_{j}^{o}) \! \right)^{1/2}, \\
\eta_{a_{j}^{o}} \! := \! \left(\! (b_{j}^{o} \! - \! a_{j}^{o})(a_{N+1}^{o} 
\! - \! a_{j}^{o})(a_{j}^{o} \! - \! b_{0}^{o}) \prod_{k=1}^{j-1}(a_{j}^{o} 
\! - \! b_{k}^{o})(a_{j}^{o} \! - \! a_{k}^{o}) \prod_{l=j+1}^{N}(b_{l}^{o} 
\! - \! a_{j}^{o})(a_{l}^{o} \! - \! a_{j}^{o}) \! \right)^{1/2},
\end{gather}
where $c^{o}_{j^{\prime}k^{\prime}}$, $j^{\prime},k^{\prime} \! = \! 1,\dotsc,
N$, are obtained {}from Equations~{\rm (O1)} and~{\rm (O2)}, $\eta_{b_{j-1}^{
o}},\eta_{a_{j}^{o}} \! > \! 0$, $j \! = \! 1, \dotsc,N \! + \! 1$, and
\begin{align}
\widehat{\alpha}^{o}_{0}(b_{0}^{o}) \! =& \, \dfrac{4}{3} \mi (-1)^{N}h_{V}^{o}
(b_{0}^{o}) \eta_{b_{0}^{o}}, \\
\widehat{\alpha}^{o}_{1}(b_{0}^{o}) \! =& \, \mi (-1)^{N} \! \left(\dfrac{2}{5}
h_{V}^{o}(b_{0}^{o}) \eta_{b_{0}^{o}} \! \left(\sum_{l=1}^{N} \! \left(\dfrac{
1}{b_{0}^{o} \! - \! b_{l}^{o}} \! + \! \dfrac{1}{b_{0}^{o} \! - \! a_{l}^{o}}
\right) \! + \! \dfrac{1}{b_{0}^{o} \! - \! a_{N+1}^{o}} \right) \! + \!
\dfrac{4}{5}(h_{V}^{o}(b_{0}^{o}))^{\prime} \eta_{b_{0}^{o}} \right), \\
\widehat{\alpha}^{o}_{0}(a_{N+1}^{o}) \! =& \, \dfrac{4}{3}h_{V}^{o}(a_{N+1}^{
o}) \eta_{a_{N+1}^{o}}, \\
\widehat{\alpha}^{o}_{1}(a_{N+1}^{o}) \! =& \, \dfrac{2}{5}h_{V}^{o}(a_{N+1}^{
o}) \eta_{a_{N+1}^{o}} \! \left(\sum_{l=1}^{N} \! \left(\dfrac{1}{a_{N+1}^{o}
\! - \! b_{l}^{o}} \! + \! \dfrac{1}{a_{N+1}^{o} \! - \! a_{l}^{o}} \right) \!
+ \! \dfrac{1}{a_{N+1}^{o} \! - \! b_{0}^{o}} \right) \! + \! \dfrac{4}{5}
(h_{V}^{o}(a_{N+1}^{o}))^{\prime} \eta_{a_{N+1}^{o}}, \\
\widehat{\alpha}^{o}_{0}(b_{j}^{o}) \! =& \, \dfrac{4}{3} \mi (-1)^{N-j}h_{V}^{
o}(b_{j}^{o}) \eta_{b_{j}^{o}}, \\
\widehat{\alpha}^{o}_{1}(b_{j}^{o}) \! =& \, \mi (-1)^{N-j} \! \left(\dfrac{
2}{5}h_{V}^{o}(b_{j}^{o}) \eta_{b_{j}^{o}} \! \left(\sum_{\substack{k=1\\k
\not= j}}^{N} \! \left(\dfrac{1}{b_{j}^{o} \! - \! b_{k}^{o}} \! + \! \dfrac{
1}{b_{j}^{o} \! - \! a_{k}^{o}} \right) \! + \! \dfrac{1}{b_{j}^{o} \! - \!
a_{j}^{o}} \! + \! \dfrac{1}{b_{j}^{o} \! - \! a_{N+1}^{o}} \! + \! \dfrac{1}{
b_{j}^{o} \! - \! b_{0}^{o}} \right) \! + \! \dfrac{4}{5}(h_{V}^{o}(b_{j}^{o}
))^{\prime} \eta_{b_{j}^{o}} \right), \\
\widehat{\alpha}^{o}_{0}(a_{j}^{o}) \! =& \, \dfrac{4}{3}(-1)^{N-j+1}h_{V}^{o}
(a_{j}^{o}) \eta_{a_{j}^{o}}, \\
\widehat{\alpha}^{o}_{1}(a_{j}^{o}) \! =& \, (-1)^{N-j+1} \! \left(\dfrac{2}{5}
h_{V}^{o}(a_{j}^{o}) \eta_{a_{j}^{o}} \! \left(\sum_{\substack{k=1\\k \not= j}
}^{N} \! \left(\dfrac{1}{a_{j}^{o} \! - \! b_{k}^{o}} \! + \! \dfrac{1}{a_{j}^{
o} \! - \! a_{k}^{o}} \right) \! + \! \dfrac{1}{a_{j}^{o} \! - \! b_{j}^{o}}
\! + \! \dfrac{1}{a_{j}^{o} \! - \! a_{N+1}^{o}} \! + \! \dfrac{1}{a_{j}^{o}
\! - \! b_{0}^{o}} \right) \! + \! \dfrac{4}{5}(h_{V}^{o}(a_{j}^{o}))^{\prime} 
\eta_{a_{j}^{o}} \right);
\end{align}
{\rm \pmb{(2)}} for $z \! \in \! \Upsilon^{o}_{2}$ $(\subset \! \mathbb{C}_{
-})$,
\begin{align}
z \boldsymbol{\pi}_{2n+1}(z) \underset{n \to \infty}{=}& \, \dfrac{1}{\mathbb{
E}} \exp \! \left(n(g^{o}(z) \! - \! \mathfrak{Q}^{-}_{\mathscr{A}}) \right)
\! \left((\overset{o}{m}^{\raise-1.0ex\hbox{$\scriptstyle \infty$}}(z))_{11}
\! \left(1 \! + \! \dfrac{1}{n \! + \! \frac{1}{2}} \! \left(\mathscr{R}^{o}_{
0}(z) \right)_{11} \right. \right. \nonumber \\
+&\left. \left. \, \mathcal{O} \! \left(\dfrac{1}{(n \! + \! \frac{1}{2})^{2}}
\right) \right) \! + \!
(\overset{o}{m}^{\raise-1.0ex\hbox{$\scriptstyle \infty$}}(z))_{21} \! \left(
\dfrac{1}{n \! + \! \frac{1}{2}} \! \left(\mathscr{R}^{o}_{0}(z) \right)_{12}
\! + \! \mathcal{O} \! \left(\dfrac{1}{(n \! + \! \frac{1}{2})^{2}} \right)
\right) \right),
\end{align}
and
\begin{align}
z \int_{\mathbb{R}} \dfrac{(s \boldsymbol{\pi}_{2n+1}(s)) \me^{-n \widetilde{V}
(s)}}{s(s \! - \! z)} \, \dfrac{\md s}{2 \pi \mi} \underset{n \to \infty}{=}&
\, \mathbb{E} \exp \! \left(-n(g^{o}(z) \! - \! \ell_{o} \! - \! \mathfrak{Q}^{
-}_{\mathscr{A}}) \right) \nonumber \\
\times& \, \left(\! (\overset{o}{m}^{\raise-1.0ex\hbox{$\scriptstyle \infty$}}
(z))_{12} \! \left(1 \! + \! \dfrac{1}{n \! + \! \frac{1}{2}} \! \left(
\mathscr{R}^{o}_{0}(z) \right)_{11} \! + \! \mathcal{O} \! \left(\dfrac{1}{(n
\! + \! \frac{1}{2})^{2}} \right) \right) \right. \nonumber \\
+&\left. \, (\overset{o}{m}^{\raise-1.0ex\hbox{$\scriptstyle \infty$}}(z))_{2
2} \! \left(\dfrac{1}{n \! + \! \frac{1}{2}} \! \left(\mathscr{R}^{o}_{0}(z)
\right)_{12} \! + \! \mathcal{O} \! \left(\dfrac{1}{(n \! + \! \frac{1}{2})^{
2}} \right) \right) \right);
\end{align}
{\rm \pmb{(3)}} for $z \! \in \! \Upsilon^{o}_{3}$ $(\subset \! \cup_{j=1}^{N
+1} \left\lbrace \mathstrut z \! \in \! \mathbb{C}^{\ast}; \, \Re (z) \! \in
\! (b_{j-1}^{o},a_{j}^{o}), \, \inf_{q \in (b_{j-1}^{o},a_{j}^{o})} \vert z \!
- \! q \vert \! < \! 2^{-1/2} \min \lbrace \delta^{o}_{b_{j-1}},\delta^{o}_{
a_{j}} \rbrace \right\rbrace \! \subset \! \mathbb{C}_{+})$,
\begin{align}
z \boldsymbol{\pi}_{2n+1}(z) \underset{n \to \infty}{=}& \, \mathbb{E} \exp \!
\left(n(g^{o}(z) \! - \! \mathfrak{Q}^{+}_{\mathscr{A}}) \right) \! \left(\!
\left(\! (\overset{o}{m}^{\raise-1.0ex\hbox{$\scriptstyle \infty$}}(z))_{11}
\! + \! (\overset{o}{m}^{\raise-1.0ex\hbox{$\scriptstyle \infty$}}(z))_{12}
\me^{-4(n+\frac{1}{2}) \pi \mi \int_{z}^{a_{N+1}^{o}} \psi_{V}^{o}(s) \, \md
s} \right) \right. \nonumber \\
\times&\left. \, \left(1 \! + \! \dfrac{1}{n \! + \! \frac{1}{2}} \! \left(
\mathscr{R}^{o}_{0}(z) \right)_{11} \! + \! \mathcal{O} \! \left(\dfrac{1}{(
n \! + \! \frac{1}{2})^{2}} \right) \right) \! + \!
\left(\! (\overset{o}{m}^{\raise-1.0ex\hbox{$\scriptstyle \infty$}}(z))_{22}
\me^{-4(n+\frac{1}{2}) \pi \mi \int_{z}^{a_{N+1}^{o}} \psi_{V}^{o}(s) \, \md
s} \right. \right. \nonumber \\
+&\left. \left. \, (\overset{o}{m}^{\raise-1.0ex\hbox{$\scriptstyle \infty$}}
(z))_{21} \right) \! \left(\dfrac{1}{n \! + \! \frac{1}{2}} \! \left(\mathscr{
R}^{o}_{0}(z) \right)_{12} \! + \! \mathcal{O} \! \left(\dfrac{1}{(n \! + \!
\frac{1}{2})^{2}} \right) \right) \right),
\end{align}
and
\begin{align}
z \int_{\mathbb{R}} \dfrac{(s \boldsymbol{\pi}_{2n+1}(s)) \me^{-n \widetilde{V}
(s)}}{s(s \! - \! z)} \, \dfrac{\md s}{2 \pi \mi} \underset{n \to \infty}{=}&
\dfrac{1}{\mathbb{E}} \exp \! \left(-n(g^{o}(z) \! - \! \ell_{o} \! - \!
\mathfrak{Q}^{+}_{\mathscr{A}}) \right) \! \left(\!
(\overset{o}{m}^{\raise-1.0ex\hbox{$\scriptstyle \infty$}}(z))_{12} \right.
\nonumber \\
\times&\left. \, \left(1 \! + \! \dfrac{1}{n \! + \! \frac{1}{2}} \! \left(
\mathscr{R}^{o}_{0}(z) \right)_{11} \! + \! \mathcal{O} \! \left(\dfrac{1}{(
n \! + \! \frac{1}{2})^{2}} \right) \right) \! + \!
(\overset{o}{m}^{\raise-1.0ex\hbox{$\scriptstyle \infty$}}(z))_{22} \right.
\nonumber \\
\times&\left. \, \left(\dfrac{1}{n \! + \! \frac{1}{2}} \! \left(\mathscr{R}^{
o}_{0}(z) \right)_{12} \! + \! \mathcal{O} \! \left(\dfrac{1}{(n \! + \! \frac{
1}{2})^{2}} \right) \right) \right);
\end{align}
{\rm \pmb{(4)}} for $z \! \in \! \Upsilon^{o}_{4}$ $(\subset \! \cup_{j=1}^{N
+1} \left\lbrace \mathstrut z \! \in \! \mathbb{C}^{\ast}; \, \Re (z) \! \in
\! (b_{j-1}^{o},a_{j}^{o}), \, \inf_{q \in (b_{j-1}^{o},a_{j}^{o})} \vert z \!
- \! q \vert \! < \! 2^{-1/2} \min \lbrace \delta^{o}_{b_{j-1}},\delta^{o}_{
a_{j}} \rbrace \right\rbrace \! \subset \! \mathbb{C}_{-})$,
\begin{align}
z \boldsymbol{\pi}_{2n+1}(z) \underset{n \to \infty}{=}& \, \dfrac{1}{\mathbb{
E}} \exp \! \left(n(g^{o}(z) \! - \! \mathfrak{Q}^{-}_{\mathscr{A}}) \right)
\! \left(\! \left(\!
(\overset{o}{m}^{\raise-1.0ex\hbox{$\scriptstyle \infty$}}(z))_{11} \! - \!
(\overset{o}{m}^{\raise-1.0ex\hbox{$\scriptstyle \infty$}}(z))_{12} \me^{4(n+
\frac{1}{2}) \pi \mi \int_{z}^{a_{N+1}^{o}} \psi_{V}^{o}(s) \, \md s} \right)
\right. \nonumber \\
\times&\left. \, \left(1 \! + \! \dfrac{1}{n \! + \! \frac{1}{2}} \! \left(
\mathscr{R}^{o}_{0}(z) \right)_{11} \! + \! \mathcal{O} \! \left(\dfrac{1}{(n
\! + \! \frac{1}{2})^{2}} \right) \right) \! + \!
\left(-(\overset{o}{m}^{\raise-1.0ex\hbox{$\scriptstyle \infty$}}(z))_{22}
\me^{4(n+\frac{1}{2}) \pi \mi \int_{z}^{a_{N+1}^{o}} \psi_{V}^{o}(s) \, \md s}
\right. \right. \nonumber \\
+&\left. \left. \, (\overset{o}{m}^{\raise-1.0ex\hbox{$\scriptstyle \infty$}}
(z))_{21} \right) \! \left(\dfrac{1}{n \! + \! \frac{1}{2}} \! \left(\mathscr{
R}^{o}_{0}(z) \right)_{12} \! + \! \mathcal{O} \! \left(\dfrac{1}{(n \! + \!
\frac{1}{2})^{2}} \right) \right) \right),
\end{align}
and
\begin{align}
z \int_{\mathbb{R}} \dfrac{(s \boldsymbol{\pi}_{2n+1}(s)) \me^{-n \widetilde{V}
(s)}}{s(s \! - \! z)} \, \dfrac{\md s}{2 \pi \mi} \underset{n \to \infty}{=}&
\, \mathbb{E} \exp \! \left(-n(g^{o}(z) \! - \! \ell_{o} \! - \! \mathfrak{Q}^{
-}_{\mathscr{A}}) \right) \! \left(\!
(\overset{o}{m}^{\raise-1.0ex\hbox{$\scriptstyle \infty$}}(z))_{12} \right.
\nonumber \\
\times&\left. \, \left(1 \! + \! \dfrac{1}{n \! + \! \frac{1}{2}} \! \left(
\mathscr{R}^{o}_{0}(z) \right)_{11} \! + \! \mathcal{O} \! \left(\dfrac{1}{(
n \! + \! \frac{1}{2})^{2}} \right) \right) \! + \!
(\overset{o}{m}^{\raise-1.0ex\hbox{$\scriptstyle \infty$}}(z))_{22} \right.
\nonumber \\
\times&\left. \, \left(\dfrac{1}{n \! + \! \frac{1}{2}} \! \left(\mathscr{R}^{
o}_{0}(z) \right)_{12} \! + \! \mathcal{O} \! \left(\dfrac{1}{(n \! + \! \frac{
1}{2})^{2}} \right) \right) \right);
\end{align}
{\rm \pmb{(5)}} for $z \! \in \! \Omega^{o,1}_{b_{j-1}}$ $(\subset \! \mathbb{
C}_{+} \cap \mathbb{U}^{o}_{\delta_{b_{j-1}}})$, $j \! = \! 1,\dotsc,N \! + \!
1$,
\begin{align}
z \boldsymbol{\pi}_{2n+1}(z) \underset{n \to \infty}{=}& \, \mathbb{E} \exp \!
\left(n(g^{o}(z) \! - \! \mathfrak{Q}^{+}_{\mathscr{A}}) \right) \! \left(\!
(m^{b,1}_{p}(z))_{11} \! \left(1 \! + \! \dfrac{1}{n \! + \! \frac{1}{2}} \!
\left(\mathscr{R}^{o}_{0}(z) \! - \! \widetilde{\mathscr{R}}^{o}_{0}(z)
\right)_{11} \right. \right. \nonumber \\
+&\left. \left. \, \mathcal{O} \! \left(\dfrac{1}{(n \! + \! \frac{1}{2})^{2}}
\right) \right) \! + \! (m^{b,1}_{p}(z))_{21} \! \left(\dfrac{1}{n \! + \!
\frac{1}{2}} \! \left(\mathscr{R}^{o}_{0}(z) \! - \! \widetilde{\mathscr{R}}^{
o}_{0}(z) \right)_{12} \! + \! \mathcal{O} \! \left(\dfrac{1}{(n \! + \! \frac{
1}{2})^{2}} \right) \right) \right),
\end{align}
and
\begin{align}
z \int_{\mathbb{R}} \dfrac{(s \boldsymbol{\pi}_{2n+1}(s)) \me^{-n \widetilde{V}
(s)}}{s(s \! - \! z)} \, \dfrac{\md s}{2 \pi \mi} \underset{n \to \infty}{=}&
\, \dfrac{1}{\mathbb{E}} \exp \! \left(-n(g^{o}(z) \! - \! \ell_{o} \! - \!
\mathfrak{Q}^{+}_{\mathscr{A}}) \right) \! \left((m^{b,1}_{p}(z))_{12} \right.
\nonumber \\
\times&\left. \, \left(1 \! + \! \dfrac{1}{n \! + \! \frac{1}{2}} \! \left(
\mathscr{R}^{o}_{0}(z) \! - \! \widetilde{\mathscr{R}}^{o}_{o}(z) \right)_{11}
\! + \! \mathcal{O} \! \left(\dfrac{1}{(n \! + \! \frac{1}{2})^{2}} \right)
\right) \! + \! (m^{b,1}_{p}(z))_{22} \right. \nonumber \\
\times&\left. \, \left(\dfrac{1}{n \! + \! \frac{1}{2}} \! \left(\mathscr{R}^{
o}_{0}(z) \! - \! \widetilde{\mathscr{R}}^{o}_{0}(z) \right)_{12} \! + \!
\mathcal{O} \! \left(\dfrac{1}{(n \! + \! \frac{1}{2})^{2}} \right) \right)
\right),
\end{align}
where
\begin{align}
(m^{b,1}_{p}(z))_{11} :=& \, -\mi \sqrt{\smash[b]{\pi}} \, \me^{\frac{1}{2}
(n+\frac{1}{2}) \xi^{o}_{b_{j-1}}(z)} \! \left(\mi \! \left(\operatorname{Ai}
(p_{b})(p_{b})^{1/4} \! - \! \operatorname{Ai}^{\prime}(p_{b})(p_{b})^{-1/4}
\right) \! (\overset{o}{m}^{\raise-1.0ex\hbox{$\scriptstyle \infty$}}(z))_{11}
\right. \nonumber \\
-&\left. \left(\operatorname{Ai}(p_{b})(p_{b})^{1/4} \! + \!
\operatorname{Ai}^{\prime}(p_{b})(p_{b})^{-1/4} \right) \!
(\overset{o}{m}^{\raise-1.0ex\hbox{$\scriptstyle \infty$}}(z))_{12} \me^{-\mi
(n+\frac{1}{2}) \mho^{o}_{j-1}} \right), \\
(m^{b,1}_{p}(z))_{12} :=& \, \sqrt{\smash[b]{\pi}} \, \me^{-\frac{\mi \pi}{6}}
\me^{-\frac{1}{2}(n+\frac{1}{2}) \xi^{o}_{b_{j-1}}(z)} \! \left(\mi \! \left(-
\operatorname{Ai}(\omega^{2}p_{b})(p_{b})^{1/4} \! + \! \omega^{2}
\operatorname{Ai}^{\prime}(\omega^{2}p_{b})(p_{b})^{-1/4} \right) \!
(\overset{o}{m}^{\raise-1.0ex\hbox{$\scriptstyle \infty$}}(z))_{11} \right.
\nonumber \\
\times &\left. \me^{\mi (n+\frac{1}{2}) \mho^{o}_{j-1}} \! + \! \left(
\operatorname{Ai}(\omega^{2}p_{b})(p_{b})^{1/4} \! + \! \omega^{2}
\operatorname{Ai}^{\prime}(\omega^{2}p_{b})(p_{b})^{-1/4} \right) \!
(\overset{o}{m}^{\raise-1.0ex\hbox{$\scriptstyle \infty$}}(z))_{12} \right), \\
(m^{b,1}_{p}(z))_{21} :=& \, -\mi \sqrt{\smash[b]{\pi}} \, \me^{\frac{1}{2}
(n+\frac{1}{2}) \xi^{o}_{b_{j-1}}(z)} \! \left(\mi \! \left(\operatorname{Ai}
(p_{b})(p_{b})^{1/4} \! - \! \operatorname{Ai}^{\prime}(p_{b})(p_{b})^{-1/4}
\right) \! (\overset{o}{m}^{\raise-1.0ex\hbox{$\scriptstyle \infty$}}(z))_{21}
\right. \nonumber \\
-&\left. \left(\operatorname{Ai}(p_{b})(p_{b})^{1/4} \! + \!
\operatorname{Ai}^{\prime}(p_{b})(p_{b})^{-1/4} \right) \!
(\overset{o}{m}^{\raise-1.0ex\hbox{$\scriptstyle \infty$}}(z))_{22} \me^{-\mi
(n+\frac{1}{2}) \mho^{o}_{j-1}} \right), \\
(m^{b,1}_{p}(z))_{22} :=& \, \sqrt{\smash[b]{\pi}} \, \me^{-\frac{\mi \pi}{6}}
\me^{-\frac{1}{2}(n+\frac{1}{2}) \xi^{o}_{b_{j-1}}(z)} \! \left(\mi \! \left(-
\operatorname{Ai}(\omega^{2}p_{b})(p_{b})^{1/4} \! + \! \omega^{2}
\operatorname{Ai}^{\prime}(\omega^{2}p_{b})(p_{b})^{-1/4} \right) \!
(\overset{o}{m}^{\raise-1.0ex\hbox{$\scriptstyle \infty$}}(z))_{21} \right.
\nonumber \\
\times &\left. \me^{\mi (n+\frac{1}{2}) \mho^{o}_{j-1}} \! + \! \left(
\operatorname{Ai}(\omega^{2}p_{b})(p_{b})^{1/4} \! + \! \omega^{2}
\operatorname{Ai}^{\prime}(\omega^{2}p_{b})(p_{b})^{-1/4} \right) \!
(\overset{o}{m}^{\raise-1.0ex\hbox{$\scriptstyle \infty$}}(z))_{22} \right),
\end{align}
with $\omega \! = \! \exp (2 \pi \mi/3)$,
\begin{gather}
\widetilde{\mathscr{R}}^{o}_{0}(z) \! := \! \sum_{j=1}^{N+1} \! \left(\mathscr{
R}^{0}_{b_{j-1}^{o}}(z) \pmb{1}_{\mathbb{U}^{o}_{\delta_{b_{j-1}}}}(z) \! + \!
\mathscr{R}^{0}_{a_{j}^{o}}(z) \pmb{1}_{\mathbb{U}^{o}_{\delta_{a_{j}}}}(z)
\right), \\
\xi^{o}_{b_{j-1}}(z) \! = \! -2 \int_{z}^{b^{o}_{j-1}}(R_{o}(s))^{1/2}h_{V}^{o}
(s) \, \md s, \qquad \quad p_{b} \! = \! \left(\dfrac{3}{4} \! \left(n \! + \!
\frac{1}{2} \right) \! \xi^{o}_{b_{j-1}}(z) \right)^{2/3}, \\
\mathscr{R}^{0}_{b_{j-1}^{o}}(z) \! = \! \dfrac{1}{\xi_{b_{j-1}}^{o}(z)}
\overset{o}{m}^{\raise-1.0ex\hbox{$\scriptstyle \infty$}}(z) \!
\begin{pmatrix}
-(s_{1}+t_{1}) & -\mi (s_{1}-t_{1}) \me^{\mi (n+\frac{1}{2}) \mho_{j-1}^{o}} \\
-\mi (s_{1}-t_{1}) \me^{-\mi (n+\frac{1}{2}) \mho_{j-1}^{o}} & (s_{1}+t_{1})
\end{pmatrix} \!
(\overset{o}{m}^{\raise-1.0ex\hbox{$\scriptstyle \infty$}}(z))^{-1}, \\
\mathscr{R}_{a_{j}^{o}}^{0}(z) \! = \! \dfrac{1}{\xi_{a_{j}}^{o}(z)}
\overset{o}{m}^{\raise-1.0ex\hbox{$\scriptstyle \infty$}}(z) \!
\begin{pmatrix}
-(s_{1}+t_{1}) & \mi (s_{1}-t_{1}) \me^{\mi (n+\frac{1}{2}) \mho_{j}^{o}} \\
\mi (s_{1}-t_{1}) \me^{-\mi (n+\frac{1}{2}) \mho_{j}^{o}} & (s_{1}+t_{1})
\end{pmatrix} \!
(\overset{o}{m}^{\raise-1.0ex\hbox{$\scriptstyle \infty$}}(z))^{-1}, \\
\xi_{a_{j}}^{o}(z) \! = \! 2 \int_{a_{j}^{o}}^{z}(R_{o}(s))^{1/2}h_{V}^{o}(s)
\, \md s,
\end{gather}
and $\pmb{1}_{\mathbb{U}^{o}_{\delta_{b_{j-1}}}}(z)$ (resp., $\pmb{1}_{\mathbb{
U}^{o}_{\delta_{a_{j}}}}(z))$ the indicator (characteristic) function of the
set $\mathbb{U}^{o}_{\delta_{b_{j-1}}}$ (resp., $\mathbb{U}^{o}_{\delta_{a_{j}
}});$\\
{\rm \pmb{(6)}} for $z \! \in \! \Omega^{o,2}_{b_{j-1}}$ $(\subset \! \mathbb{
C}_{+} \cap \mathbb{U}^{o}_{\delta_{b_{j-1}}})$, $j \! = \! 1,\dotsc,N \! + \!
1$,
\begin{align}
z \boldsymbol{\pi}_{2n+1}(z) \underset{n \to \infty}{=}& \, \mathbb{E} \exp \!
\left(n(g^{o}(z) \! - \! \mathfrak{Q}^{+}_{\mathscr{A}}) \right) \! \left(\!
\left(\! (m^{b,2}_{p}(z))_{11} \! + \! (m^{b,2}_{p}(z))_{12} \me^{-4(n+\frac{
1}{2}) \pi \mi \int_{z}^{a_{N+1}^{o}} \psi_{V}^{o}(s) \, \md s} \right)
\right. \nonumber \\
\times&\left. \, \left(1 \! + \! \dfrac{1}{n \! + \! \frac{1}{2}} \! \left(
\mathscr{R}^{o}_{0}(z) \! - \! \widetilde{\mathscr{R}}^{o}_{0}(z) \right)_{11}
\! + \! \mathcal{O} \! \left(\dfrac{1}{(n \! + \! \frac{1}{2})^{2}} \right)
\right) \! + \! \left((m^{b,2}_{p}(z))_{21} \! + \! (m^{b,2}_{p}(z))_{22}
\right. \right. \nonumber \\
\times&\left. \left. \, \me^{-4(n+\frac{1}{2}) \pi \mi \int_{z}^{a_{N+1}^{o}}
\psi_{V}^{o}(s) \, \md s} \right) \! \left(\dfrac{1}{n \! + \! \frac{1}{2}}
\! \left(\mathscr{R}^{o}_{0}(z) \! - \! \widetilde{\mathscr{R}}^{o}_{0}(z)
\right)_{12} \! + \! \mathcal{O} \! \left(\dfrac{1}{(n \! + \! \frac{1}{2})^{
2}} \right) \right) \right),
\end{align}
and
\begin{align}
z \int_{\mathbb{R}} \dfrac{(s \boldsymbol{\pi}_{2n+1}(s)) \me^{-n \widetilde{V}
(s)}}{s(s \! - \! z)} \, \dfrac{\md s}{2 \pi \mi} \underset{n \to \infty}{=}&
\, \dfrac{1}{\mathbb{E}} \exp \! \left(-n(g^{o}(z) \! - \! \ell_{o} \! - \!
\mathfrak{Q}^{+}_{\mathscr{A}}) \right) \! \left((m^{b,2}_{p}(z))_{12} \right.
\nonumber \\
\times&\left. \, \left(1 \! + \! \dfrac{1}{n \! + \! \frac{1}{2}} \! \left(
\mathscr{R}^{o}_{0}(z) \! - \! \widetilde{\mathscr{R}}^{o}_{0}(z) \right)_{11}
\! + \! \mathcal{O} \! \left(\dfrac{1}{(n \! + \! \frac{1}{2})^{2}} \right)
\right) \! + \! (m^{b,2}_{p}(z))_{22} \right. \nonumber \\
\times&\left. \, \left(\dfrac{1}{n \! + \! \frac{1}{2}} \! \left(\mathscr{R}^{
o}_{0}(z) \! - \! \widetilde{\mathscr{R}}^{o}_{0}(z) \right)_{12} \! + \!
\mathcal{O} \! \left(\dfrac{1}{(n \! + \! \frac{1}{2})^{2}} \right) \right)
\right),
\end{align}
where
\begin{align}
(m^{b,2}_{p}(z))_{11} :=& \, -\mi \sqrt{\smash[b]{\pi}} \, \me^{\frac{1}{2}
(n+\frac{1}{2}) \xi^{o}_{b_{j-1}}(z)} \! \left(\mi \! \left(-\omega
\operatorname{Ai}(\omega p_{b})(p_{b})^{1/4} \! + \! \omega^{2} \operatorname{
Ai}^{\prime}(\omega p_{b})(p_{b})^{-1/4} \right) \!
(\overset{o}{m}^{\raise-1.0ex\hbox{$\scriptstyle \infty$}}(z))_{11} \right.
\nonumber \\
+&\left. \left(\omega \operatorname{Ai}(\omega p_{b})(p_{b})^{1/4} \! + \!
\omega^{2} \operatorname{Ai}^{\prime}(\omega p_{b})(p_{b})^{-1/4} \right) \!
(\overset{o}{m}^{\raise-1.0ex\hbox{$\scriptstyle \infty$}}(z))_{12} \me^{-\mi
(n+\frac{1}{2}) \mho^{o}_{j-1}} \right), \\
(m^{b,2}_{p}(z))_{12} :=& \, \sqrt{\smash[b]{\pi}} \, \me^{-\frac{\mi \pi}{6}}
\me^{-\frac{1}{2}(n+\frac{1}{2}) \xi^{o}_{b_{j-1}}(z)} \! \left(\mi \! \left(-
\operatorname{Ai}(\omega^{2}p_{b})(p_{b})^{1/4} \! + \! \omega^{2}
\operatorname{Ai}^{\prime}(\omega p_{b})(p_{b})^{-1/4} \right) \!
(\overset{o}{m}^{\raise-1.0ex\hbox{$\scriptstyle \infty$}}(z))_{11} \right.
\nonumber \\
\times &\left. \me^{\mi (n+\frac{1}{2}) \mho^{o}_{j-1}} \! + \! \left(
\operatorname{Ai}(\omega^{2} p_{b})(p_{b})^{1/4} \! + \! \omega^{2}
\operatorname{Ai}^{\prime}(\omega p_{b})(p_{b})^{-1/4} \right) \!
(\overset{o}{m}^{\raise-1.0ex\hbox{$\scriptstyle \infty$}}(z))_{12} \right), \\
(m^{b,2}_{p}(z))_{21} :=& \, -\mi \sqrt{\smash[b]{\pi}} \, \me^{\frac{1}{2}
(n+\frac{1}{2}) \xi^{o}_{b_{j-1}}(z)} \! \left(\mi \! \left(-\omega
\operatorname{Ai}(\omega p_{b})(p_{b})^{1/4} \! + \! \omega^{2} \operatorname{
Ai}^{\prime}(\omega p_{b})(p_{b})^{-1/4} \right) \!
(\overset{o}{m}^{\raise-1.0ex\hbox{$\scriptstyle \infty$}}(z))_{21} \right.
\nonumber \\
+&\left. \left(\omega \operatorname{Ai}(\omega p_{b})(p_{b})^{1/4} \! + \!
\omega^{2} \operatorname{Ai}^{\prime}(\omega p_{b})(p_{b})^{-1/4} \right) \!
(\overset{o}{m}^{\raise-1.0ex\hbox{$\scriptstyle \infty$}}(z))_{22} \me^{-\mi
(n+\frac{1}{2}) \mho^{o}_{j-1}} \right), \\
(m^{b,2}_{p}(z))_{22} :=& \, \sqrt{\smash[b]{\pi}} \, \me^{-\frac{\mi \pi}{6}}
\me^{-\frac{1}{2}(n+\frac{1}{2}) \xi^{o}_{b_{j-1}}(z)} \! \left(\mi \! \left(-
\operatorname{Ai}(\omega^{2}p_{b})(p_{b})^{1/4} \! + \! \omega^{2}
\operatorname{Ai}^{\prime}(\omega p_{b})(p_{b})^{-1/4} \right) \!
(\overset{o}{m}^{\raise-1.0ex\hbox{$\scriptstyle \infty$}}(z))_{21} \right.
\nonumber \\
\times &\left. \me^{\mi (n+\frac{1}{2}) \mho^{o}_{j-1}} \! + \! \left(
\operatorname{Ai}(\omega^{2}p_{b})(p_{b})^{1/4} \! + \! \omega^{2}
\operatorname{Ai}^{\prime}(\omega p_{b})(p_{b})^{-1/4} \right) \!
(\overset{o}{m}^{\raise-1.0ex\hbox{$\scriptstyle \infty$}}(z))_{22} \right);
\end{align}
{\rm \pmb{(7)}} for $z \! \in \! \Omega^{o,3}_{b_{j-1}}$ $(\subset \! \mathbb{
C}_{-} \cap \mathbb{U}^{o}_{\delta_{b_{j-1}}})$, $j \! = \! 1,\dotsc,N \! + \!
1$,
\begin{align}
z \boldsymbol{\pi}_{2n+1}(z) \underset{n \to \infty}{=}& \, \dfrac{1}{\mathbb{
E}} \exp \! \left(n(g^{o}(z) \! - \! \mathfrak{Q}^{-}_{\mathscr{A}}) \right)
\! \left(\! \left(\! (m^{b,3}_{p}(z))_{11} \! - \! (m^{b,3}_{p}(z))_{12} \me^{
4(n+\frac{1}{2}) \pi \mi \int_{z}^{a_{N+1}^{o}} \psi_{V}^{o}(s) \, \md s}
\right) \right. \nonumber \\
\times&\left. \, \left(1 \! + \! \dfrac{1}{n \! + \! \frac{1}{2}} \! \left(
\mathscr{R}^{o}_{0}(z) \! - \! \widetilde{\mathscr{R}}^{o}_{0}(z) \right)_{11}
\! + \! \mathcal{O} \! \left(\dfrac{1}{(n \! + \! \frac{1}{2})^{2}} \right)
\right) \! + \! \left((m^{b,3}_{p}(z))_{21} \! - \! (m^{b,3}_{p}(z))_{22}
\right. \right. \nonumber \\
\times&\left. \left. \, \me^{4(n+\frac{1}{2}) \pi \mi \int_{z}^{a_{N+1}^{o}}
\psi_{V}^{o}(s) \, \md s} \right) \! \left(\dfrac{1}{n \! + \! \frac{1}{2}}
\! \left(\mathscr{R}^{o}_{0}(z) \! - \! \widetilde{\mathscr{R}}^{o}_{0}(z)
\right)_{12} \! + \! \mathcal{O} \! \left(\dfrac{1}{(n \! + \! \frac{1}{2})^{
2}} \right) \right) \right),
\end{align}
and
\begin{align}
z \int_{\mathbb{R}} \dfrac{(s \boldsymbol{\pi}_{2n+1}(s)) \me^{-n \widetilde{V}
(s)}}{s(s \! - \! z)} \, \dfrac{\md s}{2 \pi \mi} \underset{n \to \infty}{=}&
\, \mathbb{E} \exp \! \left(-n(g^{o}(z) \! - \! \ell_{o} \! - \! \mathfrak{Q}^{
-}_{\mathscr{A}}) \right) \! \left((m^{b,3}_{p}(z))_{12} \right. \nonumber \\
\times&\left. \, \left(1 \! + \! \dfrac{1}{n \! + \! \frac{1}{2}} \! \left(
\mathscr{R}^{o}_{0}(z) \! - \! \widetilde{\mathscr{R}}^{o}_{0}(z) \right)_{11}
\! + \! \mathcal{O} \! \left(\dfrac{1}{(n \! + \! \frac{1}{2})^{2}} \right)
\right) \! + \! (m^{b,3}_{p}(z))_{22} \right. \nonumber \\
\times&\left. \, \left(\dfrac{1}{n \! + \! \frac{1}{2}} \! \left(\mathscr{R}^{
o}_{0}(z) \! - \! \widetilde{\mathscr{R}}^{o}_{0}(z) \right)_{12} \! + \!
\mathcal{O} \! \left(\dfrac{1}{(n \! + \! \frac{1}{2})^{2}} \right) \right)
\right),
\end{align}
where
\begin{align}
(m^{b,3}_{p}(z))_{11} :=& \, -\mi \sqrt{\smash[b]{\pi}} \, \me^{\frac{1}{2}
(n+\frac{1}{2}) \xi^{o}_{b_{j-1}}(z)} \! \left(\mi \! \left(-\omega^{2}
\operatorname{Ai}(\omega^{2} p_{b})(p_{b})^{1/4} \! + \! \omega \operatorname{
Ai}^{\prime}(\omega^{2} p_{b})(p_{b})^{-1/4} \right) \!
(\overset{o}{m}^{\raise-1.0ex\hbox{$\scriptstyle \infty$}}(z))_{11} \right.
\nonumber \\
+&\left. \left(\omega^{2} \operatorname{Ai}(\omega^{2}p_{b})(p_{b})^{1/4} \! +
\! \omega \operatorname{Ai}^{\prime}(\omega^{2} p_{b})(p_{b})^{-1/4} \right)
\! (\overset{o}{m}^{\raise-1.0ex\hbox{$\scriptstyle \infty$}}(z))_{12} \me^{
\mi (n+\frac{1}{2}) \mho^{o}_{j-1}} \right), \\
(m^{b,3}_{p}(z))_{12} :=& \, \sqrt{\smash[b]{\pi}} \, \me^{-\frac{\mi \pi}{6}}
\me^{-\frac{1}{2}(n+\frac{1}{2}) \xi^{o}_{b_{j-1}}(z)} \! \left(\mi \! \left(
\omega^{2} \operatorname{Ai}(\omega p_{b})(p_{b})^{1/4} \! - \! \operatorname{
Ai}^{\prime}(\omega p_{b})(p_{b})^{-1/4} \right) \!
(\overset{o}{m}^{\raise-1.0ex\hbox{$\scriptstyle \infty$}}(z))_{11} \right.
\nonumber \\
\times &\left. \me^{-\mi (n+\frac{1}{2}) \mho^{o}_{j-1}} \! - \! \left(\omega^{
2} \operatorname{Ai}(\omega p_{b})(p_{b})^{1/4} \! + \! \operatorname{Ai}^{
\prime}(\omega p_{b})(p_{b})^{-1/4} \right) \!
(\overset{o}{m}^{\raise-1.0ex\hbox{$\scriptstyle \infty$}}(z))_{12} \right), \\
(m^{b,3}_{p}(z))_{21} :=& \, -\mi \sqrt{\smash[b]{\pi}} \, \me^{\frac{1}{2}
(n+\frac{1}{2}) \xi^{o}_{b_{j-1}}(z)} \! \left(\mi \! \left(-\omega^{2}
\operatorname{Ai}(\omega^{2}p_{b})(p_{b})^{1/4} \! + \! \omega \operatorname{
Ai}^{\prime}(\omega^{2}p_{b})(p_{b})^{-1/4} \right) \!
(\overset{o}{m}^{\raise-1.0ex\hbox{$\scriptstyle \infty$}}(z))_{21} \right.
\nonumber \\
+&\left. \left(\omega^{2} \operatorname{Ai}(\omega^{2}p_{b})(p_{b})^{1/4} \! +
\! \omega \operatorname{Ai}^{\prime}(\omega^{2}p_{b})(p_{b})^{-1/4} \right) \!
(\overset{o}{m}^{\raise-1.0ex\hbox{$\scriptstyle \infty$}}(z))_{22} \me^{\mi
(n+\frac{1}{2}) \mho^{o}_{j-1}} \right), \\
(m^{b,3}_{p}(z))_{22} :=& \, \sqrt{\smash[b]{\pi}} \, \me^{-\frac{\mi \pi}{6}}
\me^{-\frac{1}{2}(n+\frac{1}{2}) \xi^{o}_{b_{j-1}}(z)} \! \left(\mi \! \left(
\omega^{2} \operatorname{Ai}(\omega p_{b})(p_{b})^{1/4} \! - \! \operatorname{
Ai}^{\prime}(\omega p_{b})(p_{b})^{-1/4} \right) \!
(\overset{o}{m}^{\raise-1.0ex\hbox{$\scriptstyle \infty$}}(z))_{21} \right.
\nonumber \\
\times &\left. \me^{-\mi (n+\frac{1}{2}) \mho^{o}_{j-1}} \! - \! \left(\omega^{
2} \operatorname{Ai}(\omega p_{b})(p_{b})^{1/4} \! + \! \operatorname{Ai}^{
\prime}(\omega p_{b})(p_{b})^{-1/4} \right) \!
(\overset{o}{m}^{\raise-1.0ex\hbox{$\scriptstyle \infty$}}(z))_{22} \right);
\end{align}
{\rm \pmb{(8)}} for $z \! \in \! \Omega^{o,4}_{b_{j-1}}$ $(\subset \! \mathbb{
C}_{-} \cap \mathbb{U}^{o}_{\delta_{b_{j-1}}})$, $j \! = \! 1,\dotsc,N \! + \!
1$,
\begin{align}
z \boldsymbol{\pi}_{2n+1}(z) \underset{n \to \infty}{=}& \, \dfrac{1}{\mathbb{
E}} \exp \! \left(n(g^{o}(z) \! - \! \mathfrak{Q}^{-}_{\mathscr{A}}) \right)
\! \left(\! (m^{b,4}_{p}(z))_{11} \! \left(1 \! + \! \dfrac{1}{n \! + \! \frac{
1}{2}} \! \left(\mathscr{R}^{o}_{0}(z) \! - \! \widetilde{\mathscr{R}}^{o}_{0}
(z) \right)_{11} \right. \right. \nonumber \\
+&\left. \left. \, \mathcal{O} \! \left(\dfrac{1}{(n \! + \! \frac{1}{2})^{2}}
\right) \right) \! + \! (m^{b,4}_{p}(z))_{21} \! \left(\dfrac{1}{n \! + \!
\frac{1}{2}} \! \left(\mathscr{R}^{o}_{0}(z) \! - \! \widetilde{\mathscr{R}}^{
o}_{0}(z) \right)_{12} \! + \! \mathcal{O} \! \left(\dfrac{1}{(n \! + \!
\frac{1}{2})^{2}} \right) \right) \right),
\end{align}
and
\begin{align}
z \int_{\mathbb{R}} \dfrac{(s \boldsymbol{\pi}_{2n+1}(s)) \me^{-n \widetilde{V}
(s)}}{s(s \! - \! z)} \, \dfrac{\md s}{2 \pi \mi} \underset{n \to \infty}{=}&
\, \mathbb{E} \exp \! \left(-n(g^{o}(z) \! - \! \ell_{o} \! - \! \mathfrak{Q}^{
-}_{\mathscr{A}}) \right) \! \left((m^{b,4}_{p}(z))_{12} \right. \nonumber \\
\times&\left. \, \left(1 \! + \! \dfrac{1}{n \! + \! \frac{1}{2}} \! \left(
\mathscr{R}^{o}_{0}(z) \! - \! \widetilde{\mathscr{R}}^{o}_{0}(z) \right)_{11}
\! + \! \mathcal{O} \! \left(\dfrac{1}{(n \! + \! \frac{1}{2})^{2}} \right)
\right) \! + \! (m^{b,4}_{p}(z))_{22} \right. \nonumber \\
\times&\left. \, \left(\dfrac{1}{n \! + \! \frac{1}{2}} \! \left(\mathscr{R}^{
o}_{0}(z) \! - \! \widetilde{\mathscr{R}}^{o}_{0}(z) \right)_{12} \! + \!
\mathcal{O} \! \left(\dfrac{1}{(n \! + \! \frac{1}{2})^{2}} \right) \right)
\right),
\end{align}
where
\begin{align}
(m^{b,4}_{p}(z))_{11} :=& \, -\mi \sqrt{\smash[b]{\pi}} \, \me^{\frac{1}{2}
(n+\frac{1}{2}) \xi^{o}_{b_{j-1}}(z)} \! \left(\mi \! \left(\operatorname{Ai}
(p_{b})(p_{b})^{1/4} \! - \! \operatorname{Ai}^{\prime}(p_{b})(p_{b})^{-1/4}
\right) \! (\overset{o}{m}^{\raise-1.0ex\hbox{$\scriptstyle \infty$}}(z))_{11}
\right. \nonumber \\
-&\left. \left(\operatorname{Ai}(p_{b})(p_{b})^{1/4} \! + \!
\operatorname{Ai}^{\prime}(p_{b})(p_{b})^{-1/4} \right) \!
(\overset{o}{m}^{\raise-1.0ex\hbox{$\scriptstyle \infty$}}(z))_{12} \me^{\mi
(n+\frac{1}{2}) \mho^{o}_{j-1}} \right), \\
(m^{b,4}_{p}(z))_{12} :=& \, \sqrt{\smash[b]{\pi}} \, \me^{-\frac{\mi \pi}{6}}
\me^{-\frac{1}{2}(n+\frac{1}{2}) \xi^{o}_{b_{j-1}}(z)} \! \left(\mi \! \left(
\omega^{2} \operatorname{Ai}(\omega p_{b})(p_{b})^{1/4} \! - \! \operatorname{
Ai}^{\prime}(\omega p_{b})(p_{b})^{-1/4} \right) \!
(\overset{o}{m}^{\raise-1.0ex\hbox{$\scriptstyle \infty$}}(z))_{11} \right.
\nonumber \\
\times &\left. \me^{-\mi (n+\frac{1}{2}) \mho^{o}_{j-1}} \! - \! \left(\omega^{
2} \operatorname{Ai}(\omega p_{b})(p_{b})^{1/4} \! + \! \operatorname{Ai}^{
\prime}(\omega p_{b})(p_{b})^{-1/4} \right) \!
(\overset{o}{m}^{\raise-1.0ex\hbox{$\scriptstyle \infty$}}(z))_{12} \right), \\
(m^{b,4}_{p}(z))_{21} :=& \, -\mi \sqrt{\smash[b]{\pi}} \, \me^{\frac{1}{2}
(n+\frac{1}{2}) \xi^{o}_{b_{j-1}}(z)} \! \left(\mi \! \left(\operatorname{Ai}
(p_{b})(p_{b})^{1/4} \! - \! \operatorname{Ai}^{\prime}(p_{b})(p_{b})^{-1/4}
\right) \! (\overset{o}{m}^{\raise-1.0ex\hbox{$\scriptstyle \infty$}}(z))_{21}
\right. \nonumber \\
-&\left. \left(\operatorname{Ai}(p_{b})(p_{b})^{1/4} \! + \!
\operatorname{Ai}^{\prime}(p_{b})(p_{b})^{-1/4} \right) \!
(\overset{o}{m}^{\raise-1.0ex\hbox{$\scriptstyle \infty$}}(z))_{22} \me^{\mi
(n+\frac{1}{2}) \mho^{o}_{j-1}} \right), \\
(m^{b,4}_{p}(z))_{22} :=& \, \sqrt{\smash[b]{\pi}} \, \me^{-\frac{\mi \pi}{6}}
\me^{-\frac{1}{2}(n+\frac{1}{2}) \xi^{o}_{b_{j-1}}(z)} \! \left(\mi \!
\left(\omega^{2} \operatorname{Ai}(\omega p_{b})(p_{b})^{1/4} \! - \!
\operatorname{Ai}^{\prime}(\omega p_{b})(p_{b})^{-1/4} \right) \!
(\overset{o}{m}^{\raise-1.0ex\hbox{$\scriptstyle \infty$}}(z))_{21} \right.
\nonumber \\
\times &\left. \me^{-\mi (n+\frac{1}{2}) \mho^{o}_{j-1}} \! - \! \left(\omega^{
2} \operatorname{Ai}(\omega p_{b})(p_{b})^{1/4} \! + \! \operatorname{Ai}^{
\prime}(\omega p_{b})(p_{b})^{-1/4} \right) \!
(\overset{o}{m}^{\raise-1.0ex\hbox{$\scriptstyle \infty$}}(z))_{22}
\right);
\end{align}
{\rm \pmb{(9)}} for $z \! \in \! \Omega^{o,1}_{a_{j}}$ $(\subset \! \mathbb{
C}_{+} \cap \mathbb{U}^{o}_{\delta_{a_{j}}})$, $j \! = \! 1,\dotsc,N \! + \!
1$,
\begin{align}
z \boldsymbol{\pi}_{2n+1}(z) \underset{n \to \infty}{=}& \, \mathbb{E} \exp \!
\left(n(g^{o}(z) \! - \! \mathfrak{Q}^{+}_{\mathscr{A}}) \right) \! \left(\!
(m^{a,1}_{p}(z))_{11} \! \left(1 \! + \! \dfrac{1}{n \! + \! \frac{1}{2}} \!
\left(\mathscr{R}^{o}_{0}(z) \! - \! \widetilde{\mathscr{R}}^{o}_{0}(z)
\right)_{11} \right. \right. \nonumber \\
+&\left. \left. \, \mathcal{O} \! \left(\dfrac{1}{(n \! + \! \frac{1}{2})^{2}}
\right) \right) \! + \! (m^{a,1}_{p}(z))_{21} \! \left(\dfrac{1}{n \! + \!
\frac{1}{2}} \! \left(\mathscr{R}^{o}_{0}(z) \! - \! \widetilde{\mathscr{R}}^{
o}_{0}(z) \right)_{12} \! + \! \mathcal{O} \! \left(\dfrac{1}{(n \! + \!
\frac{1}{2})^{2}} \right) \right) \right),
\end{align}
and
\begin{align}
z \int_{\mathbb{R}} \dfrac{(s \boldsymbol{\pi}_{2n+1}(s)) \me^{-n \widetilde{V}
(s)}}{s(s \! - \! z)} \, \dfrac{\md s}{2 \pi \mi} \underset{n \to \infty}{=}&
\, \dfrac{1}{\mathbb{E}} \exp \! \left(-n(g^{o}(z) \! - \! \ell_{o} \! - \!
\mathfrak{Q}^{+}_{\mathscr{A}}) \right) \! \left((m^{a,1}_{p}(z))_{12} \right.
\nonumber \\
\times&\left. \, \left(1 \! + \! \dfrac{1}{n \! + \! \frac{1}{2}} \! \left(
\mathscr{R}^{o}_{0}(z) \! - \! \widetilde{\mathscr{R}}^{o}_{0}(z) \right)_{11}
\! + \! \mathcal{O} \! \left(\dfrac{1}{(n \! + \! \frac{1}{2})^{2}} \right)
\right) \! + \! (m^{a,1}_{p}(z))_{22} \right. \nonumber \\
\times&\left. \, \left(\dfrac{1}{n \! + \! \frac{1}{2}} \! \left(\mathscr{R}^{
o}_{0}(z) \! - \! \widetilde{\mathscr{R}}^{o}_{0}(z) \right)_{12} \! + \!
\mathcal{O} \! \left(\dfrac{1}{(n \! + \! \frac{1}{2})^{2}} \right) \right)
\right),
\end{align}
where
\begin{align}
(m^{a,1}_{p}(z))_{11} :=& \, -\mi \sqrt{\smash[b]{\pi}} \, \me^{\frac{1}{2}
(n+\frac{1}{2}) \xi^{o}_{a_{j}}(z)} \! \left(\mi \! \left(\operatorname{Ai}
(p_{a})(p_{a})^{1/4} \! - \! \operatorname{Ai}^{\prime}(p_{a})(p_{a})^{-1/4}
\right) \! (\overset{o}{m}^{\raise-1.0ex\hbox{$\scriptstyle \infty$}}(z))_{11}
\right. \nonumber \\
+&\left. \left(\operatorname{Ai}(p_{a})(p_{a})^{1/4} \! + \!
\operatorname{Ai}^{\prime}(p_{a})(p_{a})^{-1/4} \right) \!
(\overset{o}{m}^{\raise-1.0ex\hbox{$\scriptstyle \infty$}}(z))_{12} \me^{-\mi
(n+\frac{1}{2}) \mho^{o}_{j}} \right), \\
(m^{a,1}_{p}(z))_{12} :=& \, \sqrt{\smash[b]{\pi}} \, \me^{-\frac{\mi \pi}{6}}
\me^{-\frac{1}{2}(n+\frac{1}{2}) \xi^{o}_{a_{j}}(z)} \! \left(\mi \! \left(
\operatorname{Ai}(\omega^{2}p_{a})(p_{a})^{1/4} \! - \! \omega^{2}
\operatorname{Ai}^{\prime}(\omega^{2}p_{a})(p_{a})^{-1/4} \right) \!
(\overset{o}{m}^{\raise-1.0ex\hbox{$\scriptstyle \infty$}}(z))_{11} \right.
\nonumber \\
\times &\left. \me^{\mi (n+\frac{1}{2}) \mho^{o}_{j}} \! + \! \left(
\operatorname{Ai}(\omega^{2}p_{a})(p_{a})^{1/4} \! + \! \omega^{2}
\operatorname{Ai}^{\prime}(\omega^{2}p_{a})(p_{a})^{-1/4} \right) \!
(\overset{o}{m}^{\raise-1.0ex\hbox{$\scriptstyle \infty$}}(z))_{12} \right), \\
(m^{a,1}_{p}(z))_{21} :=& \, -\mi \sqrt{\smash[b]{\pi}} \, \me^{\frac{1}{2}
(n+\frac{1}{2}) \xi^{o}_{a_{j}}(z)} \! \left(\mi \! \left(\operatorname{Ai}
(p_{a})(p_{a})^{1/4} \! - \! \operatorname{Ai}^{\prime}(p_{a})(p_{a})^{-1/4}
\right) \! (\overset{o}{m}^{\raise-1.0ex\hbox{$\scriptstyle \infty$}}(z))_{21}
\right. \nonumber \\
+&\left. \left(\operatorname{Ai}(p_{a})(p_{a})^{1/4} \! + \!
\operatorname{Ai}^{\prime}(p_{a})(p_{a})^{-1/4} \right) \!
(\overset{o}{m}^{\raise-1.0ex\hbox{$\scriptstyle \infty$}}(z))_{22} \me^{-\mi
(n+\frac{1}{2}) \mho^{o}_{j}} \right), \\
(m^{a,1}_{p}(z))_{22} :=& \, \sqrt{\smash[b]{\pi}} \, \me^{-\frac{\mi \pi}{6}}
\me^{-\frac{1}{2}(n+\frac{1}{2}) \xi^{o}_{a_{j}}(z)} \! \left(\mi \! \left(
\operatorname{Ai}(\omega^{2}p_{a})(p_{a})^{1/4} \! - \! \omega^{2}
\operatorname{Ai}^{\prime}(\omega^{2}p_{a})(p_{a})^{-1/4} \right) \!
(\overset{o}{m}^{\raise-1.0ex\hbox{$\scriptstyle \infty$}}(z))_{21} \right.
\nonumber \\
\times &\left. \me^{\mi (n+\frac{1}{2}) \mho^{o}_{j}} \! + \! \left(
\operatorname{Ai}(\omega^{2}p_{a})(p_{a})^{1/4} \! + \! \omega^{2}
\operatorname{Ai}^{\prime}(\omega^{2}p_{a})(p_{a})^{-1/4} \right) \!
(\overset{o}{m}^{\raise-1.0ex\hbox{$\scriptstyle \infty$}}(z))_{22} \right),
\end{align}
with
\begin{equation}
p_{a} \! = \! \left(\dfrac{3}{4} \! \left(n \! + \! \dfrac{1}{2} \right) \!
\xi^{o}_{a_{j}}(z) \right)^{2/3};
\end{equation}
{\rm \pmb{(10)}} for $z \! \in \! \Omega^{o,2}_{a_{j}}$ $(\subset \! \mathbb{
C}_{+} \cap \mathbb{U}^{o}_{\delta_{a_{j}}})$, $j \! = \! 1,\dotsc,N \! + \!
1$,
\begin{align}
z \boldsymbol{\pi}_{2n+1}(z) \underset{n \to \infty}{=}& \, \mathbb{E} \exp \!
\left(n(g^{0}(z) \! - \! \mathfrak{Q}^{+}_{\mathscr{A}}) \right) \! \left(\!
\left(\! (m^{a,2}_{p}(z))_{11} \! + \! (m^{a,2}_{p}(z))_{12} \me^{-4(n+\frac{
1}{2}) \pi \mi \int_{z}^{a_{N+1}^{o}} \psi_{V}^{o}(s) \, \md s} \right)
\right. \nonumber \\
\times&\left. \, \left(1 \! + \! \dfrac{1}{n \! + \! \frac{1}{2}} \! \left(
\mathscr{R}^{o}_{0}(z) \! - \! \widetilde{\mathscr{R}}^{o}_{0}(z) \right)_{11}
\! + \! \mathcal{O} \! \left(\dfrac{1}{(n \! + \! \frac{1}{2})^{2}} \right)
\right) \! + \! \left((m^{a,2}_{p}(z))_{21} \! + \! (m^{a,2}_{p}(z))_{22}
\right. \right. \nonumber \\
\times&\left. \left. \, \me^{-4(n+\frac{1}{2}) \pi \mi \int_{z}^{a_{N+1}^{o}}
\psi_{V}^{o}(s) \, \md s} \right) \! \left(\dfrac{1}{n} \! \left(\mathscr{R}^{
o}_{0}(z) \! - \! \widetilde{\mathscr{R}}^{o}_{0}(z) \right)_{12} \! + \!
\mathcal{O} \! \left(\dfrac{1}{(n \! + \! \frac{1}{2})^{2}} \right) \right)
\right),
\end{align}
and
\begin{align}
z \int_{\mathbb{R}} \dfrac{(s \boldsymbol{\pi}_{2n+1}(s)) \me^{-n \widetilde{V}
(s)}}{s(s \! - \! z)} \, \dfrac{\md s}{2 \pi \mi} \underset{n \to \infty}{=}&
\, \dfrac{1}{\mathbb{E}} \exp \! \left(-n(g^{o}(z) \! - \! \ell_{o} \! - \!
\mathfrak{Q}^{+}_{\mathscr{A}}) \right) \! \left((m^{a,2}_{p}(z))_{12} \right.
\nonumber \\
\times&\left. \, \left(1 \! + \! \dfrac{1}{n \! + \! \frac{1}{2}} \! \left(
\mathscr{R}^{o}_{0}(z) \! - \! \widetilde{\mathscr{R}}^{o}_{0}(z) \right)_{11}
\! + \! \mathcal{O} \! \left(\dfrac{1}{(n \! + \! \frac{1}{2})^{2}} \right)
\right) \! + \! (m^{a,2}_{p}(z))_{22} \right. \nonumber \\
\times&\left. \, \left(\dfrac{1}{n \! + \! \frac{1}{2}} \! \left(\mathscr{R}^{
o}_{0}(z) \! - \! \widetilde{\mathscr{R}}^{o}_{0}(z) \right)_{12} \! + \!
\mathcal{O} \! \left(\dfrac{1}{(n \! + \! \frac{1}{2})^{2}} \right) \right)
\right),
\end{align}
where
\begin{align}
(m^{a,2}_{p}(z))_{11} :=& \, -\mi \sqrt{\smash[b]{\pi}} \, \me^{\frac{1}{2}
(n+\frac{1}{2}) \xi^{o}_{a_{j}}(z)} \! \left(\mi \! \left(-\omega
\operatorname{Ai}(\omega p_{a})(p_{a})^{1/4} \! + \! \omega^{2}
\operatorname{Ai}^{\prime}(\omega p_{a})(p_{a})^{-1/4} \right) \!
(\overset{o}{m}^{\raise-1.0ex\hbox{$\scriptstyle \infty$}}(z))_{11} \right.
\nonumber \\
-&\left. \left(\omega \operatorname{Ai}(\omega p_{a})(p_{a})^{1/4} \! + \!
\omega^{2} \operatorname{Ai}^{\prime}(\omega p_{a})(p_{a})^{-1/4} \right) \!
(\overset{o}{m}^{\raise-1.0ex\hbox{$\scriptstyle \infty$}}(z))_{12} \me^{-\mi
(n+\frac{1}{2}) \mho^{o}_{j}} \right), \\
(m^{a,2}_{p}(z))_{12} :=& \, \sqrt{\smash[b]{\pi}} \, \me^{-\frac{\mi \pi}{6}}
\me^{-\frac{1}{2}(n+\frac{1}{2}) \xi^{o}_{a_{j}}(z)} \! \left(\mi \! \left(
\operatorname{Ai}(\omega^{2}p_{a})(p_{a})^{1/4} \! - \! \omega^{2}
\operatorname{Ai}^{\prime}(\omega^{2}p_{a})(p_{a})^{-1/4} \right) \!
(\overset{o}{m}^{\raise-1.0ex\hbox{$\scriptstyle \infty$}}(z))_{11} \right.
\nonumber \\
\times &\left. \me^{\mi (n+\frac{1}{2}) \mho^{o}_{j}} \! + \! \left(
\operatorname{Ai}(\omega^{2}p_{a})(p_{a})^{1/4} \! + \! \omega^{2}
\operatorname{Ai}^{\prime}(\omega^{2}p_{a})(p_{a})^{-1/4} \right) \!
(\overset{o}{m}^{\raise-1.0ex\hbox{$\scriptstyle \infty$}}(z))_{12} \right), \\
(m^{a,2}_{p}(z))_{21} :=& \, -\mi \sqrt{\smash[b]{\pi}} \, \me^{\frac{1}{2}
(n+\frac{1}{2}) \xi^{o}_{a_{j}}(z)} \! \left(\mi \! \left(-\omega
\operatorname{Ai}(\omega p_{a})(p_{a})^{1/4} \! + \! \omega^{2}
\operatorname{Ai}^{\prime}(\omega p_{a})(p_{a})^{-1/4} \right) \!
(\overset{o}{m}^{\raise-1.0ex\hbox{$\scriptstyle \infty$}}(z))_{21} \right.
\nonumber \\
-&\left. \left(\omega \operatorname{Ai}(\omega p_{a})(p_{a})^{1/4} \! + \!
\omega^{2} \operatorname{Ai}^{\prime}(\omega p_{a})(p_{a})^{-1/4} \right) \!
(\overset{o}{m}^{\raise-1.0ex\hbox{$\scriptstyle \infty$}}(z))_{22} \me^{-\mi
(n+\frac{1}{2}) \mho^{o}_{j}} \right), \\
(m^{a,2}_{p}(z))_{22} :=& \, \sqrt{\smash[b]{\pi}} \, \me^{-\frac{\mi \pi}{6}}
\me^{-\frac{1}{2}(n+\frac{1}{2}) \xi^{o}_{a_{j}}(z)} \! \left(\mi \! \left(
\operatorname{Ai}(\omega^{2}p_{a})(p_{a})^{1/4} \! - \! \omega^{2}
\operatorname{Ai}^{\prime}(\omega^{2}p_{a})(p_{a})^{-1/4} \right) \!
(\overset{o}{m}^{\raise-1.0ex\hbox{$\scriptstyle \infty$}}(z))_{21} \right.
\nonumber \\
\times &\left. \me^{\mi (n+\frac{1}{2}) \mho^{o}_{j}} \! + \! \left(
\operatorname{Ai}(\omega^{2}p_{a})(p_{a})^{1/4} \! + \! \omega^{2}
\operatorname{Ai}^{\prime}(\omega^{2}p_{a})(p_{a})^{-1/4} \right) \!
(\overset{o}{m}^{\raise-1.0ex\hbox{$\scriptstyle \infty$}}(z))_{22} \right);
\end{align}
{\rm \pmb{(11)}} for $z \! \in \! \Omega^{o,3}_{a_{j}}$ $(\subset \! \mathbb{
C}_{-} \cap \mathbb{U}^{o}_{\delta_{a_{j}}})$, $j \! = \! 1,\dotsc,N \! + \!
1$,
\begin{align}
z \boldsymbol{\pi}_{2n+1}(z) \underset{n \to \infty}{=}& \, \dfrac{1}{\mathbb{
E}} \exp \! \left(n(g^{o}(z) \! - \! \mathfrak{Q}^{-}_{\mathscr{A}}) \right)
\! \left(\! \left(\! (m^{a,3}_{p}(z))_{11} \! - \! (m^{a,3}_{p}(z))_{12} \me^{
4(n+\frac{1}{2}) \pi \mi \int_{z}^{a_{N+1}^{o}} \psi_{V}^{o}(s) \, \md s}
\right) \right. \nonumber \\
\times&\left. \, \left(1 \! + \! \dfrac{1}{n \! + \! \frac{1}{2}} \! \left(
\mathscr{R}^{o}_{0}(z) \! - \! \widetilde{\mathscr{R}}^{o}_{0}(z) \right)_{11}
\! + \! \mathcal{O} \! \left(\dfrac{1}{(n \! + \! \frac{1}{2})^{2}} \right)
\right) \! + \! \left((m^{a,3}_{p}(z))_{21} \! - \! (m^{a,3}_{p}(z))_{22}
\right. \right. \nonumber \\
\times&\left. \left. \, \me^{4(n+\frac{1}{2}) \pi \mi \int_{z}^{a_{N+1}^{o}}
\psi_{V}^{o}(s) \, \md s} \right) \! \left(\dfrac{1}{n \! + \! \frac{1}{2}}
\! \left(\mathscr{R}^{o}_{0}(z) \! - \! \widetilde{\mathscr{R}}^{o}_{0}(z)
\right)_{12} \! + \! \mathcal{O} \! \left(\dfrac{1}{(n \! + \! \frac{1}{2})^{
2}} \right) \right) \right),
\end{align}
and
\begin{align}
z \int_{\mathbb{R}} \dfrac{(s \boldsymbol{\pi}_{2n+1}(s)) \me^{-n \widetilde{V}
(s)}}{s(s \! - \! z)} \, \dfrac{\md s}{2 \pi \mi} \underset{n \to \infty}{=}&
\, \mathbb{E} \exp \! \left(-n(g^{o}(z) \! - \! \ell_{o} \! - \! \mathfrak{Q}^{
-}_{\mathscr{A}}) \right) \! \left((m^{a,3}_{p}(z))_{12} \right. \nonumber \\
\times&\left. \, \left(1 \! + \! \dfrac{1}{n \! + \! \frac{1}{2}} \! \left(
\mathscr{R}^{o}_{0}(z) \! - \! \widetilde{\mathscr{R}}^{o}_{0}(z) \right)_{11}
\! + \! \mathcal{O} \! \left(\dfrac{1}{(n \! + \! \frac{1}{2})^{2}} \right)
\right) \! + \! (m^{a,3}_{p}(z))_{22} \right. \nonumber \\
\times&\left. \, \left(\dfrac{1}{n \! + \! \frac{1}{2}} \! \left(\mathscr{R}^{
o}_{0}(z) \! - \! \widetilde{\mathscr{R}}^{o}_{0}(z) \right)_{12} \! + \!
\mathcal{O} \! \left(\dfrac{1}{(n \! + \! \frac{1}{2})^{2}} \right) \right)
\right),
\end{align}
where
\begin{align}
(m^{a,3}_{p}(z))_{11} :=& \, -\mi \sqrt{\smash[b]{\pi}} \, \me^{\frac{1}{2}
(n+\frac{1}{2}) \xi^{o}_{a_{j}}(z)} \! \left(\mi \! \left(-\omega^{2}
\operatorname{Ai}(\omega^{2}p_{a})(p_{a})^{1/4} \! + \! \omega
\operatorname{Ai}^{\prime}(\omega^{2}p_{a})(p_{a})^{-1/4} \right) \!
(\overset{o}{m}^{\raise-1.0ex\hbox{$\scriptstyle \infty$}}(z))_{11} \right.
\nonumber \\
-&\left. \left(\omega^{2} \operatorname{Ai}(\omega^{2}p_{a})(p_{a})^{1/4} \! +
\! \omega \operatorname{Ai}^{\prime}(\omega^{2}p_{a})(p_{a})^{-1/4} \right) \!
(\overset{o}{m}^{\raise-1.0ex\hbox{$\scriptstyle \infty$}}(z))_{12} \me^{\mi
(n+\frac{1}{2}) \mho^{o}_{j}} \right), \\
(m^{a,3}_{p}(z))_{12} :=& \, \sqrt{\smash[b]{\pi}} \, \me^{-\frac{\mi \pi}{6}}
\me^{-\frac{1}{2}(n+\frac{1}{2}) \xi^{o}_{a_{j}}(z)} \! \left(\mi \! \left(-
\omega^{2} \operatorname{Ai}(\omega p_{a})(p_{a})^{1/4} \! + \! \operatorname{
Ai}^{\prime}(\omega p_{a})(p_{a})^{-1/4} \right) \!
(\overset{o}{m}^{\raise-1.0ex\hbox{$\scriptstyle \infty$}}(z))_{11} \right.
\nonumber \\
\times &\left. \me^{-\mi (n+\frac{1}{2}) \mho^{o}_{j}} \! - \! \left(\omega^{
2} \operatorname{Ai}(\omega p_{a})(p_{a})^{1/4} \! + \! \operatorname{Ai}^{
\prime}(\omega p_{a})(p_{a})^{-1/4} \right) \!
(\overset{o}{m}^{\raise-1.0ex\hbox{$\scriptstyle \infty$}}(z))_{12} \right), \\
(m^{a,3}_{p}(z))_{21} :=& \, -\mi \sqrt{\smash[b]{\pi}} \, \me^{\frac{1}{2}
(n+\frac{1}{2}) \xi^{o}_{a_{j}}(z)} \! \left(\mi \! \left(-\omega^{2}
\operatorname{Ai}(\omega^{2}p_{a})(p_{a})^{1/4} \! + \! \omega
\operatorname{Ai}^{\prime}(\omega^{2}p_{a})(p_{a})^{-1/4} \right) \!
(\overset{o}{m}^{\raise-1.0ex\hbox{$\scriptstyle \infty$}}(z))_{21} \right.
\nonumber \\
-&\left. \left(\omega^{2} \operatorname{Ai}(\omega^{2}p_{a})(p_{a})^{1/4} \! +
\! \omega \operatorname{Ai}^{\prime}(\omega^{2}p_{a})(p_{a})^{-1/4} \right) \!
(\overset{o}{m}^{\raise-1.0ex\hbox{$\scriptstyle \infty$}}(z))_{22} \me^{\mi
(n+\frac{1}{2}) \mho^{o}_{j}} \right), \\
(m^{a,3}_{p}(z))_{22} :=& \, \sqrt{\smash[b]{\pi}} \, \me^{-\frac{\mi \pi}{6}}
\me^{-\frac{1}{2}(n+\frac{1}{2}) \xi^{o}_{a_{j}}(z)} \! \left(\mi \! \left(-
\omega^{2} \operatorname{Ai}(\omega p_{a})(p_{a})^{1/4} \! + \! \operatorname{
Ai}^{\prime}(\omega p_{a})(p_{a})^{-1/4} \right) \!
(\overset{o}{m}^{\raise-1.0ex\hbox{$\scriptstyle \infty$}}(z))_{21} \right.
\nonumber \\
\times &\left. \me^{-\mi (n+\frac{1}{2}) \mho^{o}_{j}} \! - \! \left(\omega^{
2} \operatorname{Ai}(\omega p_{a})(p_{a})^{1/4} \! + \! \operatorname{Ai}^{
\prime}(\omega p_{a})(p_{a})^{-1/4} \right) \!
(\overset{o}{m}^{\raise-1.0ex\hbox{$\scriptstyle \infty$}}(z))_{22} \right);
\end{align}
and {\rm \pmb{(12)}} for $z \! \in \! \Omega^{o,4}_{a_{j}}$ $(\subset \!
\mathbb{C}_{-} \cap \mathbb{U}^{o}_{\delta_{a_{j}}})$, $j \! = \! 1,\dotsc,N
\! + \! 1$,
\begin{align}
z \boldsymbol{\pi}_{2n+1}(z) \underset{n \to \infty}{=}& \, \dfrac{1}{\mathbb{
E}} \exp \! \left(n(g^{o}(z) \! - \! \mathfrak{Q}^{-}_{\mathscr{A}}) \right)
\! \left(\! (m^{a,4}_{p}(z))_{11} \! \left(1 \! + \! \dfrac{1}{n \! + \! \frac{
1}{2}} \! \left(\mathscr{R}^{o}_{0}(z) \! - \! \widetilde{\mathscr{R}}^{o}_{0}
(z) \right)_{11} \right. \right. \nonumber \\
+&\left. \left. \, \mathcal{O} \! \left(\dfrac{1}{(n \! + \! \frac{1}{2})^{2}}
\right) \right) \! + \! (m^{a,4}_{p}(z))_{21} \! \left(\dfrac{1}{n \! + \!
\frac{1}{2}} \! \left(\mathscr{R}^{o}_{0}(z) \! - \! \widetilde{\mathscr{R}}^{
o}_{0}(z) \right)_{12} \! + \! \mathcal{O} \! \left(\dfrac{1}{(n \! + \!
\frac{1}{2}^{2}} \right) \right) \right),
\end{align}
and
\begin{align}
z \int_{\mathbb{R}} \dfrac{(s \boldsymbol{\pi}_{2n+1}(s)) \me^{-n \widetilde{V}
(s)}}{s(s \! - \! z)} \, \dfrac{\md s}{2 \pi \mi} \underset{n \to \infty}{=}&
\, \mathbb{E} \exp \! \left(-n(g^{o}(z) \! - \! \ell_{o} \! - \! \mathfrak{Q}^{
-}_{\mathscr{A}}) \right) \! \left((m^{a,4}_{p}(z))_{12} \right. \nonumber \\
\times&\left. \, \left(1 \! + \! \dfrac{1}{n \! + \! \frac{1}{2}} \! \left(
\mathscr{R}^{o}_{0}(z) \! - \! \widetilde{\mathscr{R}}^{o}_{0}(z) \right)_{11}
\! + \! \mathcal{O} \! \left(\dfrac{1}{(n \! + \! \frac{1}{2})^{2}} \right)
\right) \! + \! (m^{a,4}_{p}(z))_{22} \right. \nonumber \\
\times&\left. \, \left(\dfrac{1}{n \! + \! \frac{1}{2}} \! \left(\mathscr{R}^{
o}_{0}(z) \! - \! \widetilde{\mathscr{R}}^{o}_{0}(z) \right)_{12} \! + \!
\mathcal{O} \! \left(\dfrac{1}{(n \! + \! \frac{1}{2})^{2}} \right) \right)
\right),
\end{align}
where
\begin{align}
(m^{a,4}_{p}(z))_{11} :=& \, -\mi \sqrt{\smash[b]{\pi}} \, \me^{\frac{1}{2}
(n+\frac{1}{2}) \xi^{o}_{a_{j}}(z)} \! \left(\mi \! \left(\operatorname{Ai}
(p_{a})(p_{a})^{1/4} \! - \! \operatorname{Ai}^{\prime}(p_{a})(p_{a})^{-1/4}
\right) \! (\overset{o}{m}^{\raise-1.0ex\hbox{$\scriptstyle \infty$}}(z))_{11}
\right. \nonumber \\
+&\left. \left(\operatorname{Ai}(p_{a})(p_{a})^{1/4} \! + \!
\operatorname{Ai}^{\prime}(p_{a})(p_{a})^{-1/4} \right) \!
(\overset{o}{m}^{\raise-1.0ex\hbox{$\scriptstyle \infty$}}(z))_{12} \me^{\mi
(n+\frac{1}{2}) \mho^{o}_{j}} \right), \\
(m^{a,4}_{p}(z))_{12} :=& \, \sqrt{\smash[b]{\pi}} \, \me^{-\frac{\mi \pi}{6}}
\me^{-\frac{1}{2}(n+\frac{1}{2}) \xi^{o}_{a_{j}}(z)} \! \left(\mi \! \left(-
\omega^{2} \operatorname{Ai}(\omega p_{a})(p_{a})^{1/4} \! + \! \operatorname{
Ai}^{\prime}(\omega p_{a})(p_{a})^{-1/4} \right) \!
(\overset{o}{m}^{\raise-1.0ex\hbox{$\scriptstyle \infty$}}(z))_{11} \right.
\nonumber \\
\times &\left. \me^{-\mi (n+\frac{1}{2}) \mho^{o}_{j}} \! - \! \left(\omega^{
2} \operatorname{Ai}(\omega p_{a})(p_{a})^{1/4} \! + \! \operatorname{Ai}^{
\prime}(\omega p_{a})(p_{a})^{-1/4} \right) \!
(\overset{o}{m}^{\raise-1.0ex\hbox{$\scriptstyle \infty$}}(z))_{12} \right), \\
(m^{a,4}_{p}(z))_{21} :=& \, -\mi \sqrt{\smash[b]{\pi}} \, \me^{\frac{1}{2}
(n+\frac{1}{2}) \xi^{o}_{a_{j}}(z)} \! \left(\mi \! \left(\operatorname{Ai}
(p_{a})(p_{a})^{1/4} \! - \! \operatorname{Ai}^{\prime}(p_{a})(p_{a})^{-1/4}
\right) \! (\overset{o}{m}^{\raise-1.0ex\hbox{$\scriptstyle \infty$}}(z))_{21}
\right. \nonumber \\
+&\left. \left(\operatorname{Ai}(p_{a})(p_{a})^{1/4} \! + \!
\operatorname{Ai}^{\prime}(p_{a})(p_{a})^{-1/4} \right) \!
(\overset{o}{m}^{\raise-1.0ex\hbox{$\scriptstyle \infty$}}(z))_{22} \me^{\mi
(n+\frac{1}{2}) \mho^{o}_{j}} \right), \\
(m^{a,4}_{p}(z))_{22} :=& \, \sqrt{\smash[b]{\pi}} \, \me^{-\frac{\mi \pi}{6}}
\me^{-\frac{1}{2}(n+\frac{1}{2}) \xi^{o}_{a_{j}}(z)} \! \left(\mi \! \left(-
\omega^{2} \operatorname{Ai}(\omega p_{a})(p_{b})^{1/4} \! + \! \operatorname{
Ai}^{\prime}(\omega p_{a})(p_{a})^{-1/4} \right) \!
(\overset{o}{m}^{\raise-1.0ex\hbox{$\scriptstyle \infty$}}(z))_{21} \right.
\nonumber \\
\times &\left. \me^{-\mi (n+\frac{1}{2}) \mho^{o}_{j}} \! - \! \left(\omega^{
2} \operatorname{Ai}(\omega p_{a})(p_{a})^{1/4} \! + \! \operatorname{Ai}^{
\prime}(\omega p_{a})(p_{a})^{-1/4} \right) \!
(\overset{o}{m}^{\raise-1.0ex\hbox{$\scriptstyle \infty$}}(z))_{22} \right).
\end{align}
\end{dddd}
\begin{eeee}
Using limiting values, if necessary, all of the above (asymptotic) formulae 
for $\boldsymbol{\pi}_{2n+1}(z)$ and $z \int_{\mathbb{R}}(s \boldsymbol{\pi}_{
2n+1}(s)) \me^{-n \widetilde{V}(s)}(s(s \! - \! z))^{-1} \, \tfrac{\md s}{2 
\pi \mi}$ have a natural interpretation on the real and imaginary axes. \hfill 
$\blacksquare$
\end{eeee}
\begin{dddd}
Let all the conditions stated in Theorem~{\rm 2.3.1} be valid, and let 
$\overset{o}{\operatorname{Y}} \colon \mathbb{C} \setminus \mathbb{R} \! 
\to \! \operatorname{SL}_{2}(\mathbb{C})$ be the unique solution of 
{\rm \pmb{RHP2}}. Let $H^{(m)}_{k}$, $(m,k) \! \in \! \mathbb{Z} \times 
\mathbb{N}$, be the Hankel determinants associated with the bi-infinite, 
real-valued, strong moment sequence $\left\lbrace c_{k} \! = \! \int_{\mathbb{
R}}s^{k} \me^{-n \widetilde{V}(s)} \, \md s, \, n \! \in \! \mathbb{N} 
\right\rbrace_{k \in \mathbb{Z}}$ defined in Equations~{\rm (1.1)}, and 
let $\boldsymbol{\pi}_{2n+1}(z)$ be the odd degree monic orthogonal 
$L$-polynomial defined in Lemma~{\rm 2.2.2}, that is, $z \boldsymbol{
\pi}_{2n+1}(z) \linebreak[4]
:= \! (\overset{o}{\mathrm{Y}}(z))_{11}$, with $n \! \to \! \infty$ 
asymptotics (in the entire complex plane) given in Theorem~{\rm 2.3.1}. 
Then,
\begin{equation}
(\xi^{(2n+1)}_{-n-1})^{2} \! = \! \dfrac{1}{\norm{\boldsymbol{\pi}_{2n+1}
(\cdot)}_{\mathscr{L}}^{2}} \! = \dfrac{H^{(-2n)}_{2n+1}}{H^{(-2n-2)}_{2n+1}}
\underset{n \to \infty}{=} \dfrac{\me^{-n \ell_{o}}}{\pi} \Xi^{\natural} \!
\left(1 \! + \! \frac{1}{n \! + \! \frac{1}{2}} \Xi^{\natural}(\mathfrak{Q}^{
\natural})_{12} \! + \! \mathcal{O} \! \left(\dfrac{c^{\natural}(n)}{(n \! +
\! \frac{1}{2})^{2}} \right) \right),
\end{equation}
where
\begin{align}
\Xi^{\natural} :=& \, 2 \mathbb{E}^{2} \! \left(\sum_{k=1}^{N+1} \! \left((b_{k
-1}^{o})^{-1} \! - \! (a_{k}^{o})^{-1} \right) \right)^{-1} \dfrac{\boldsymbol{
\theta}^{o}(\boldsymbol{u}^{o}_{+}(0) \! - \! \frac{1}{2 \pi}(n \! + \! \frac{
1}{2}) \boldsymbol{\Omega}^{o} \! + \! \boldsymbol{d}_{o}) \boldsymbol{
\theta}^{o}(-\boldsymbol{u}^{o}_{+}(0) \! + \! \boldsymbol{d}_{o})}{
\boldsymbol{\theta}^{o}(-\boldsymbol{u}^{o}_{+}(0) \! - \! \frac{1}{2 \pi}(n
\! + \! \frac{1}{2}) \boldsymbol{\Omega}^{o} \! + \! \boldsymbol{d}_{o})
\boldsymbol{\theta}^{o}(\boldsymbol{u}^{o}_{+}(0) \! + \! \boldsymbol{d}_{o}
)}, \\
\mathfrak{Q}^{\natural} :=& \, -2 \mi \sum_{j=1}^{N+1} \! \left(\dfrac{(
\mathscr{A}^{o}(a_{j}^{o})(\widehat{\alpha}^{o}_{1}(a_{j}^{o}) \! + \! 2(a_{
j}^{o})^{-1} \widehat{\alpha}^{o}_{0}(a_{j}^{o})) \! - \! \mathscr{B}^{o}
(a_{j}^{o}) \widehat{\alpha}^{o}_{0}(a_{j}^{o}))}{(a_{j}^{o})^{2}(\widehat{
\alpha}^{o}_{0}(a_{j}^{o}))^{2}} \right. \nonumber \\
+&\left. \, \dfrac{(\mathscr{A}^{o}(b_{j-1}^{o})(\widehat{\alpha}^{o}_{1}(b_{j
-1}^{o}) \! + \! 2(b_{j-1}^{o})^{-1} \widehat{\alpha}^{o}_{0}(b_{j-1}^{o})) \!
- \! \mathscr{B}^{o}(b_{j-1}^{o}) \widehat{\alpha}^{o}_{0}(b_{j-1}^{o}))}{(b_{
j-1}^{o})^{2}(\widehat{\alpha}^{o}_{0}(b_{j-1}^{o}))^{2}} \right),
\end{align}
$(\mathfrak{Q}^{\natural})_{12}$ denotes the $(1 \, 2)$-element of $\mathfrak{
Q}^{\natural}$, $c^{\natural}(n) \! =_{n \to \infty} \! \mathcal{O}(1)$, and 
all relevant parameters are defined in Theorem~{\rm 2.3.1:} asymptotics for 
$\xi^{(2n+1)}_{-n-1}$ are obtained by taking the positive square root of 
both sides of Equation~{\rm (2.118)}. Furthermore, the $n \! \to \! \infty$ 
asymptotic expansion (in the entire complex plane) for the odd degree 
orthonormal $L$-polynomial,
\begin{equation}
\phi_{2n+1}(z) \! = \! \xi^{(2n+1)}_{-n-1} \boldsymbol{\pi}_{2n+1}(z),
\end{equation}
to $\mathcal{O} \! \left((n \! + \! 1/2)^{-2} \right)$, is given by the 
(scalar) multiplication of the $n \! \to \! \infty$ asymptotics of 
$\boldsymbol{\pi}_{2n+1}(z)$ and $\xi^{(2n+1)}_{-n-1}$ stated, respectively, 
in Theorem~{\rm 2.3.1} and Equations~{\rm (2.118)--(2.120)}.
\end{dddd}
\begin{eeee}
Since, {}from general theory (cf. Section~1), and, by construction (cf. 
Equations~(1.3) and~(1.8)), $\xi^{(2n+1)}_{-n-1} \! > \! 0$, it follows, 
incidentally, {}from Theorem~2.3.2, Equations~(2.118)--(2.120) that: (i) 
$\Xi^{\natural} \! > \! 0$; and (ii) $\Im ((\mathfrak{Q}^{\natural})_{12}) 
\! = \! 0$. \hfill $\blacksquare$
\end{eeee}
\section{The Equilibrium Measure, the Variational Problem, and the 
Tr\-a\-n\-s\-f\-o\-r\-m\-e\-d RHP}
In this section, the detailed analysis of the `odd degree' variational 
problem, and the associated `odd' equilibrium measure, is undertaken (see 
Lemmas~3.1--3.3 and Lemma~3.5), including the discussion of the corresponding 
$g$-function, denoted, herein, as $g^{o}$, and \textbf{RHP2}, that is, 
$(\overset{o}{\operatorname{Y}}(z),\mathrm{I} \! + \! \exp (-n \widetilde{V}
(z)) \sigma_{+},\linebreak[4]
\mathbb{R})$, is reformulated as an equivalent\footnote{If there are two RHPs, 
$(\mathcal{Y}_{1}(z),\upsilon_{1}(z),\varGamma_{1})$ and $(\mathcal{Y}_{2}(z),
\upsilon_{2}(z),\varGamma_{2})$, say, with $\varGamma_{2} \subset \varGamma_{
1}$ and $\upsilon_{1}(z) \! \! \upharpoonright_{\varGamma_{1} \setminus 
\varGamma_{2}} \! =_{n \to \infty} \! \mathrm{I} \! + \! o(1)$, then, within 
the BC framework \cite{a74}, and modulo $o(1)$ estimates, their solutions, 
$\mathcal{Y}_{1}$ and $\mathcal{Y}_{2}$, respectively, are (asymptotically) 
equal.}, auxiliary RHP (see Lemma~3.4). The proofs of Lemmas~3.1--3.3 are 
modelled on the calculations of Saff-Totik (\cite{a43}, Chapter~1), Deift 
(\cite{a79}, Chapter~6), and Johansson \cite{a80}.

One begins by establishing the existence of the `odd' equilibrium measure, 
$\mu_{V}^{o}$ $(\in \! \mathcal{M}_{1}(\mathbb{R}))$.
\begin{ccccc}
Let the external field $\widetilde{V} \colon \mathbb{R} \setminus \{0\} \! \to
\! \mathbb{R}$ satisfy conditions~{\rm (2.3)--(2.5)}, and set $w^{o}(z) \! :=
\! \me^{-\widetilde{V}(z)}$. For $\mu^{o} \! \in \! \mathcal{M}_{1}(\mathbb{
R})$, define the weighted energy functional $\mathrm{I}_{V}^{o}[\mu^{o}]
\colon \mathcal{M}_{1}(\mathbb{R}) \! \to \! \mathbb{R}$,
\begin{equation*}
\mathrm{I}_{V}^{o}[\mu^{o}] \! := \! \iint_{\mathbb{R}^{2}} \ln \! \left(\vert
s \! - \! t \vert^{2+\frac{1}{n}} \vert st \vert^{-1}w^{o}(s)w^{o}(t) \right)^{
-1} \md \mu^{o}(s) \, \md \mu^{o}(t), \quad n \! \in \! \mathbb{N},
\end{equation*}
and consider the minimisation problem
\begin{equation*}
E_{V}^{o} \! = \! \inf \! \left\lbrace \mathstrut \mathrm{I}_{V}^{o}[\mu^{o}]; 
\, \mu^{o} \! \in \! \mathcal{M}_{1}(\mathbb{R}) \right\rbrace.
\end{equation*}
Then: {\rm (1)} $E_{V}^{o}$ is finite; {\rm (2)} $\exists \, \, \mu_{V}^{o}
\in \! \mathcal{M}_{1}(\mathbb{R})$ such that $\mathrm{I}_{V}^{o}[\mu_{V}^{o}]
\! = \! E_{V}^{o}$ (the infimum is attained), and $\mu_{V}^{o}$ has finite
weighted logarithmic energy $(-\infty \! < \! \mathrm{I}_{V}^{o}[\mu_{V}^{o}]
\! < \! +\infty);$ and {\rm (3)} $J_{o} \! := \! \operatorname{supp}(\mu_{
V}^{o})$ is compact, $J_{o} \subset \lbrace \mathstrut z; \, w^{o}(z) \! >
\! 0 \rbrace$, and $J_{o}$ has positive logarithmic capacity, that is,
$\operatorname{cap}(J_{o}) \! := \! \exp \! \left(-\inf \lbrace \mathrm{I}_{
V}^{o}[\mu^{o}]; \, \mu^{o} \! \in \! \mathcal{M}_{1}(J_{o}) \rbrace \right)
\! > \! 0$.
\end{ccccc}

\emph{Proof.} Let $\mu^{o} \! \in \! \mathcal{M}_{1}(\mathbb{R})$, and
set\footnote{All the introduced variables in the proof depend on $n$, which
would necessitate the introduction of additional, superfluous notation to
encode this $n$ dependence; however, for simplicity of notation, such
cumbersome $n$-dependencies will not be introduced, and the reader should be
cognizant of this fact: \emph{mutatis mutandis} for the remainder of the
paper.} $w^{o}(z) \! := \! \exp (-\widetilde{V}(z))$, where $\widetilde{V}
\colon \mathbb{R} \setminus \{0\} \! \to \! \mathbb{R}$ satisfies
conditions~(2.3)--(2.5). {}From the definition of $\mathrm{I}_{V}^{o}
[\mu^{o}]$ given in the Lemma, one shows that, for $n \! \in \! \mathbb{N}$,
\begin{align*}
\mathrm{I}_{V}^{o}[\mu^{o}] =& \iint_{\mathbb{R}^{2}} \! \left(\! \left(1 \!
+ \! \dfrac{1}{n} \right) \! \ln (\vert s \! - \! t \vert^{-1}) \! + \! \ln(
\vert s^{-1} \! - \! t^{-1} \vert^{-1}) \right) \md \mu^{o}(s) \, \md \mu^{o}
(t) \! + \! 2 \int_{\mathbb{R}} \widetilde{V}(s) \, \md \mu^{o}(s) \\
=:& \iint_{\mathbb{R}^{2}}K_{V,n}^{o}(s,t) \, \md \mu^{o}(s) \, \md \mu^{o}(t),
\end{align*}
where (the $n$-dependent symmetric kernel)
\begin{equation*}
K_{V,n}^{o}(s,t) \! = \! K_{V,n}^{o}(t,s) \! := \! \left(1 \! + \! \dfrac{1}{
n} \right) \! \ln (\vert s \! - \! t \vert^{-1}) \! + \! \ln (\vert s^{-1} \!
- \! t^{-1} \vert^{-1}) \! + \! \widetilde{V}(s) \! + \! \widetilde{V}(t)
\end{equation*}
(of course, the definition of $\mathrm{I}_{V}^{o}[\mu^{o}]$ only makes sense
provided both integrals exist and are finite). Recall the following
inequalities (see, for example, Chapter~6 of \cite{a79}): $\vert s \! - \! t
\vert \! \leqslant \! (1 \! + \! s^{2})^{1/2}(1 \! + \! t^{2})^{1/2}$ and
$\vert s^{-1} \! - \! t^{-1} \vert \! \leqslant \! (1 \! + \! s^{-2})^{1/2}(1
\! + \! t^{-2})^{1/2}$, $s,t \! \in \! \mathbb{R}$, whence $\ln (\vert s \! -
\! t \vert^{-1}) \! \geqslant \! -\tfrac{1}{2} \ln (1 \! + \! s^{2}) \! - \!
\tfrac{1}{2} \ln (1 \! + \! t^{2})$ and $\ln (\vert s^{-1} \! - \! t^{-1}
\vert^{-1}) \! \geqslant \! -\tfrac{1}{2} \ln (1 \! + \! s^{-2}) \! - \!
\tfrac{1}{2} \ln (1 \! + \! t^{-2})$; thus,
\begin{equation*}
K_{V,n}^{o}(s,t) \! \geqslant \! \dfrac{1}{2} \! \left(2 \widetilde{V}(s) \! -
\! \left(1 \! + \! \dfrac{1}{n} \right) \! \ln (s^{2} \! + \! 1) \! - \! \ln
(s^{-2} \! + \! 1) \right) \! + \! \dfrac{1}{2} \! \left(2 \widetilde{V}(t) \!
- \! \left(1 \! + \! \dfrac{1}{n} \right) \! \ln (t^{2} \! + \! 1) \! - \! \ln
(t^{-2} \! + \! 1) \right).
\end{equation*}
Recalling conditions~(2.3)--(2.5) for the external field $\widetilde{V} \colon
\mathbb{R} \setminus \{0\} \! \to \! \mathbb{R}$, in particular, $\exists \,
\, \delta_{1}$ $(= \! \delta_{1}(n))$ $> \! 0$ (resp., $\exists \, \, \delta_{
2}$ $(= \! \delta_{2}(n))$ $> \! 0)$ such that $\widetilde{V}(x) \! \geqslant
\! (1 \! + \! \delta_{1})(1 \! + \! \tfrac{1}{n}) \ln (x^{2} \! + \! 1)$
(resp., $\widetilde{V}(x) \! \geqslant \! (1 \! + \! \delta_{2}) \ln (x^{-2}
\! + \! 1))$ for sufficiently large $\vert x \vert$ (resp., small $\vert x
\vert)$, it follows that $2 \widetilde{V}(x) \! - \! (1 \! + \! \tfrac{1}{n})
\ln (x^{2} \! + \! 1) \! - \! \ln (x^{-2} \! + \! 1) \! \geqslant \! C_{V}^{o}
\! > \! -\infty$, whence $K_{V,n}^{o}(s,t) \! \geqslant \! C_{V}^{o}$ $(> \!
-\infty)$, which shows that $K_{V,n}^{o}(s,t)$ is bounded {}from below (on
$\mathbb{R}^{2})$; hence,
\begin{equation*}
\mathrm{I}_{V}^{o}[\mu^{o}] \! \geqslant \! \iint_{\mathbb{R}^{2}}C_{V}^{o} \,
\md \mu^{o}(s) \, \md \mu^{o}(t) \! = \! C_{V}^{o} \underbrace{\int_{\mathbb{
R}} \md \mu^{o}(s)}_{= \, 1} \, \, \underbrace{\int_{\mathbb{R}} \md \mu^{o}
(t)}_{= \, 1} \! \geqslant \! C_{V}^{o} \quad (> \! -\infty).
\end{equation*}
It follows {}from the above inequality and the definition of $E_{V}^{o}$ 
stated in the Lemma that, $\forall \, \, \mu^{o} \! \in \! \mathcal{M}_{1}
(\mathbb{R})$, $E_{V}^{o} \! \geqslant \! C_{V}^{o} \! > \! -\infty$, which 
shows that $E_{V}^{o}$ is bounded {}from below. Let $\varepsilon$ be an 
arbitrarily fixed, sufficiently small positive real number, and set $\Sigma_{
o,\varepsilon} \! := \! \lbrace \mathstrut z; \, w^{o}(z) \! \geqslant \! 
\varepsilon \rbrace$; then $\Sigma_{o,\varepsilon}$ is compact, and $\Sigma_{
o,0} \! := \! \cup_{l=1}^{\infty} \Sigma_{o,1/l} \! = \! \cup_{l=1}^{\infty} 
\lbrace \mathstrut z; \, w^{o}(z) \! > \! l^{-1} \rbrace \! = \! \lbrace 
\mathstrut z; \, w^{o}(z) \! > \! 0 \rbrace$. Since, for $\widetilde{V} 
\colon \mathbb{R} \setminus \{0\} \! \to \! \mathbb{R}$ satisfying 
conditions~(2.3)--(2.5), $w^{o}$ is an \emph{admissible weight} \cite{a43}, 
in which case $\Sigma_{o,0}$ has positive logarithmic capacity, that is, 
$\operatorname{cap}(\Sigma_{o,0}) \! = \! \exp (-\inf \lbrace \mathstrut 
\mathrm{I}_{V}^{o}[\mu^{o}]; \, \mu^{o} \! \in \! \mathcal{M}_{1}(\Sigma_{o,
0}) \rbrace) \! > \! 0$, it follows that $\exists \, \, l^{\ast} \! \in \! 
\mathbb{N}$ such that $\operatorname{cap}(\Sigma_{o,1/l^{\ast}}) \! = \! \exp 
(-\inf \lbrace \mathstrut \mathrm{I}_{V}^{o}[\mu^{o}]; \, \mu^{o} \! \in \! 
\mathcal{M}_{1}(\Sigma_{o,1/l^{\ast}}) \rbrace) \! > \! 0$, which, in turn, 
means that there exists a probability measure, $\mu^{o}_{l^{\ast}}$, say, with 
$\operatorname{supp}(\mu^{o}_{l^{\ast}}) \subseteq \Sigma_{o,1/l^{\ast}}$, 
such that $\iint_{\Sigma_{o,1/l^{\ast}}^{2}} \ln (\vert s \! - \! t \vert^{-
(2+\frac{1}{n})} \vert st \vert) \, \md \mu^{o}_{l^{\ast}}(s) \, \md \mu^{
o}_{l^{\ast}}(t) \! < \! +\infty$, where $\Sigma_{o,1/l^{\ast}}^{2} \! = \! 
\Sigma_{o,1/l^{\ast}} \times \Sigma_{o,1/l^{\ast}}$ $(\subseteq \mathbb{R}^{
2})$. For $z \! \in \! \operatorname{supp}(\mu^{o}_{l^{\ast}}) \! \subseteq 
\! \Sigma_{o,1/l^{\ast}}$, it follows that $w^{o}(z) \! \geqslant \! 1/l^{
\ast}$, whence $\iint_{\Sigma_{o,1/l^{\ast}}^{2}} \ln (w^{o}(s)w^{o}(t))^{-1} 
\, \md \mu^{o}_{l^{\ast}}(s) \, \md \mu^{o}_{l^{\ast}}(t) \! \leqslant \! 2 
\ln (l^{\ast}) \! < \! +\infty$ $\Rightarrow$
\begin{equation*}
\mathrm{I}_{V}^{o}[\mu_{l^{\ast}}^{o}] \! = \! \iint_{\Sigma_{o,1/l^{\ast}}^{
2}} \ln \! \left(\vert s \! - \! t \vert^{2+\frac{1}{n}} \vert st \vert^{-1}
w^{o}(s)w^{o}(t) \right)^{-1} \, \md \mu^{o}_{l^{\ast}}(s) \, \md \mu^{o}_{
l^{\ast}}(t) \! < \! +\infty;
\end{equation*}
thus, it follows that $E_{V}^{o} \! := \! \inf \lbrace \mathstrut \mathrm{
I}_{V}^{o}[\mu^{o}]; \, \mu^{o} \! \in \! \mathcal{M}_{1}(\mathbb{R}) \rbrace$
is finite (see, also, below).

Choose a sequence of probability measures $\{\mu^{o}_{m}\}_{m=1}^{\infty}$ in
$\mathcal{M}_{1}(\mathbb{R})$ such that $\mathrm{I}_{V}^{o}[\mu^{o}_{m}] \!
\leqslant \! E_{V}^{o} \! + \! \tfrac{1}{m}$. {}From the analysis above, it
follows that
\begin{align*}
\mathrm{I}_{V}^{o}[\mu_{m}^{o}] =& \iint_{\mathbb{R}^{2}}K_{V,n}^{o}(s,t) \,
\md \mu^{o}_{m}(s) \, \md \mu^{o}_{m}(t) \! \geqslant \! \iint_{\mathbb{R}^{2}
} \! \left(\dfrac{1}{2} \! \left(2 \widetilde{V}(s) \! - \! \left(1 \! + \!
\dfrac{1}{n} \right) \! \ln (s^{2} \! + \! 1) \! - \! \ln (s^{-2} \! + \! 1)
\right) \right. \\
+& \left. \dfrac{1}{2} \! \left(2 \widetilde{V}(t) \! - \! \left(1 \! + \!
\dfrac{1}{n} \right) \! \ln (t^{2} \! + \! 1) \! - \! \ln (t^{-2} \! + \! 1)
\right) \right) \md \mu_{m}^{o}(s) \, \md \mu_{m}^{o}(t).
\end{align*}
Set
\begin{equation*}
\widehat{\psi}_{V}^{o}(z) \! := \! 2 \widetilde{V}(z) \! - \! \left(1 \! + \!
\dfrac{1}{n} \right) \! \ln (z^{2} \! + \! 1) \! - \! \ln (z^{-2} \! + \! 1).
\end{equation*}
Then $\mathrm{I}_{V}^{o}[\mu_{m}^{o}] \! \geqslant \! \int_{\mathbb{R}}
\widehat{\psi}_{V}^{o}(s) \, \md \mu_{m}^{o}(s)$ $\Rightarrow$ $E_{V}^{o} \! +
\! \tfrac{1}{m} \! \geqslant \! \int_{\mathbb{R}} \widehat{\psi}_{V}^{o}(s) \,
\md \mu_{m}^{o}(s)$. Recalling that $\exists \, \, \delta_{1} \! > \! 0$
(resp., $\exists \, \, \delta_{2} \! > \! 0)$ such that $\widetilde{V}(x) \!
\geqslant \! (1 \! + \! \delta_{1})(1 \! + \! \tfrac{1}{n}) \ln (x^{2} \! + \!
1)$ (resp., $\widetilde{V}(x) \! \geqslant \! (1 \! + \! \delta_{2}) \ln (x^{-
2} \! + \! 1))$ for sufficiently large $\vert x \vert$ (resp., small $\vert x
\vert)$, it follows that, for any $b_{o} \! > \! 0$, $\exists \, \, M_{o} \! >
\! 1$ such that $\widehat{\psi}_{V}^{o}(z) \! > \! b_{o} \, \, \forall \, \, z
\! \in \! \lbrace \vert z \vert \! \geqslant \! M_{o} \rbrace \cup \lbrace
\vert z \vert \! \leqslant \! M_{o}^{-1} \rbrace \! =: \! \mathfrak{D}_{o}$,
which implies that
\begin{align*}
E_{V}^{o} \! + \! \dfrac{1}{m} \geqslant& \int_{\mathbb{R}} \widehat{\psi}_{
V}^{o}(s) \, \md \mu_{m}^{o}(s) \! = \! \int_{\mathfrak{D}_{o}} \underbrace{
\widehat{\psi}_{V}^{o}(s)}_{> \, b_{o}} \, \md \mu_{m}^{o}(s) \! + \! \int_{
\mathbb{R} \setminus \mathfrak{D}_{o}} \underbrace{\widehat{\psi}_{V}^{o}(s)
}_{\geqslant \, -\vert C_{V}^{o} \vert} \, \md \mu_{m}^{o}(s) \\
\geqslant& \, b_{o} \int_{\mathfrak{D}_{o}} \md \mu^{o}_{m}(s) \! - \! \vert
C_{V}^{o} \vert \underbrace{\int_{\mathbb{R} \setminus \mathfrak{D}_{o}} \md
\mu_{m}^{o}(s)}_{\in \, [0,1]} \! \geqslant \! b_{o} \int_{\mathfrak{D}_{o}}
\md \mu_{m}^{o}(s) \! - \! \vert C_{V}^{o} \vert;
\end{align*}
thus,
\begin{equation*}
\int_{\mathfrak{D}_{o}} \md \mu_{m}^{o}(s) \! \leqslant \! b_{o}^{-1} \!
\left(E_{V}^{o} \! + \! \vert C_{V}^{o} \vert \! + \! \dfrac{1}{m} \right),
\end{equation*}
whence
\begin{equation*}
\limsup_{m \to \infty} \int_{\mathfrak{D}_{o}} \md \mu_{m}^{o}(s) \! \leqslant
\! \limsup_{m \to \infty} \! \left(b_{o}^{-1} \! \left(E_{V}^{o} \! + \! \vert
C_{V}^{o} \vert \! + \! \dfrac{1}{m} \right) \right).
\end{equation*}
By the Archimedean property, it follows that, $\forall \, \, \epsilon_{o} \! 
> \! 0$, $\exists \, \, N \! \in \! \mathbb{N}$ such that, $\forall \, \, m 
\! > \! N \Rightarrow m^{-1} \! < \! \epsilon_{o}$; thus, choosing $b_{o} \! 
= \! \epsilon^{-1}(E_{V}^{o} \! + \! \vert C_{V}^{o} \vert \! + \! \epsilon_{
o})$, where $\epsilon$ is some arbitrarily fixed, sufficiently small positive 
real number, it follows that $\limsup_{m \to \infty} \int_{\mathfrak{D}_{o}} 
\md \mu_{m}^{o}(s) \! \leqslant \! \epsilon$ $\Rightarrow$ the sequence of 
probability measures $\lbrace \mu_{m}^{o} \rbrace_{m=1}^{\infty}$ in 
$\mathcal{M}_{1}(\mathbb{R})$ is \emph{tight} \cite{a80} (that is, given 
$\epsilon \! > \! 0$, $\exists \, \, M \! > \! 1$ such that $\limsup_{m \to 
\infty} \mu_{m}^{o}(\lbrace \vert s \vert \! \geqslant \! M \rbrace \cup 
\lbrace \vert s \vert \! \leqslant \! M^{-1} \rbrace) \! := \! \limsup_{m \to 
\infty} \int_{\lbrace \vert s \vert \geqslant M \rbrace \cup \lbrace \vert 
s \vert \leqslant M^{-1} \rbrace} \md \mu_{m}^{o}(s) \! \leqslant \! 
\epsilon)$. Since the sequence of probabilty measures $\lbrace \mu_{m}^{o} 
\rbrace_{m=1}^{\infty}$ in $\mathcal{M}_{1}(\mathbb{R})$ is tight, by a 
Helly Selection Theorem, there exists a $(\mathrm{weak}^{\ast}$ convergent) 
subsequence of probability measures $\lbrace \mu^{o}_{m_{k}} \rbrace_{k=1}^{
\infty}$ in $\mathcal{M}_{1}(\mathbb{R})$ converging (weakly) to a probability 
measure $\mu^{o} \! \in \! \mathcal{M}_{1}(\mathbb{R})$, symbolically 
$\mu^{o}_{m_{k}} \! \overset{\ast}{\to} \! \mu^{o}$ as $k \! \to \! 
\infty$\footnote{A sequence of probability measures $\lbrace \mu_{m} 
\rbrace_{m=1}^{\infty}$ in $\mathcal{M}_{1}(D)$ is said to \emph{converge 
weakly} as $m \! \to \! \infty$ to $\mu \! \in \! \mathcal{M}_{1}(D)$, 
symbolically $\mu_{m} \! \overset{\ast}{\to} \! \mu$, if $\mu_{m}(f) \! := \! 
\int_{D}f(s) \, \md \mu_{m}(s) \! \to \! \int_{D}f(s) \, \md \mu (s) \! =: \! 
\mu (f)$ as $m \! \to \! \infty$ $\forall \, \, f \! \in \! \mathbf{\mathrm{
C}_{b}^{0}}(D)$, where $\mathbf{\mathrm{C}_{b}^{0}}(D)$ denotes the set of 
all bounded, continuous functions on $D$ with compact support.}. One now shows 
that, if $\mu_{m}^{o} \! \overset{\ast}{\to} \! \mu^{o}$, $\mu_{m}^{o},\mu^{o} 
\! \in \! \mathcal{M}_{1}(\mathbb{R})$, then $\liminf_{m \to \infty} \mathrm{
I}_{V}^{o}[\mu_{m}^{o}] \! \geqslant \! \mathrm{I}_{V}^{o}[\mu^{o}]$. Since 
$w^{o}$ is continuous, thus upper semi-continuous \cite{a43}, there exists 
a sequence $\lbrace w^{o}_{m} \rbrace_{m=1}^{\infty}$ (resp., $\lbrace 
\widetilde{V}_{m} \rbrace_{m=1}^{\infty})$ of continuous functions on 
$\mathbb{R}$ such that $w^{o}_{m+1} \! \leqslant \! w^{o}_{m}$ (resp., 
$\widetilde{V}_{m+1} \! \geqslant \! \widetilde{V}_{m})$\footnote{Adding a 
suitable constant, if necessary, which does not change $\mu_{m}^{o}$, or the 
regularity of $\widetilde{V} \colon \mathbb{R} \setminus \lbrace 0 \rbrace 
\! \to \! \mathbb{R}$, one may assume that $\widetilde{V} \! \geqslant \! 
0$ and $\widetilde{V}_{m} \! \geqslant \! 0$, $m \! \in \! \mathbb{N}$.}, 
$m \! \in \! \mathbb{N}$, and $w^{o}_{m}(z) \! \searrow \! w^{o}(z)$ (resp., 
$\widetilde{V}_{m}(z) \! \nearrow \! \widetilde{V}(z))$ as $m \! \to \! 
\infty$ for every $z \! \in \! \mathbb{R}$; in particular,
\begin{equation*}
\mathrm{I}_{V}^{o}[\mu_{k}^{o}] \! = \! \iint_{\mathbb{R}^{2}}K_{V,n}^{o}(s,t)
\, \md \mu_{k}^{o}(s) \, \md \mu_{k}^{o}(t) \! \geqslant \! \iint_{\mathbb{
R}^{2}}K_{V_{m},n}^{o}(s,t) \, \md \mu_{k}^{o}(s) \, \md \mu_{k}^{o}(t).
\end{equation*}
For arbitrary $q \! \in \! \mathbb{R}$, $\mathrm{I}_{V}^{o}[\mu_{k}^{o}] \!
\geqslant \! \iint_{\mathbb{R}^{2}}p^{o}(s,t;n) \, \md \mu_{k}^{o}(s) \,
\md \mu_{k}^{o}(t)$, where $p^{o}(s,t;n) \! := \! \min \left\lbrace q,
K_{V_{m},n}^{o} (s,t) \right\rbrace \! = \! p^{o}(t,s;n)$ (bounded and
continuous on $\mathbb{R}^{2})$. Recall that $\lbrace \mu_{m}^{o} \rbrace_{m
=1}^{\infty}$ is tight in $\mathcal{M}_{1}(\mathbb{R})$. For $M_{o} \! > \!
1$, let $h_{M}^{o}(x) \! \in \! \mathbf{\mathrm{C}_{b}^{0}}(\mathbb{R})$ be
such that:
\begin{compactenum}
\item[(i)] $h_{M}^{o}(x) \! = \! 1$, $x \! \in \! [-M_{o},-M_{o}^{-1}] \cup
[M_{o}^{-1},M_{o}] \! =: \! \mathfrak{D}_{M_{o}}$;
\item[(ii)] $h_{M}^{o}(x) \! = \! 0$, $x \! \in \! \mathbb{R} \setminus
\mathfrak{D}_{M_{o}+1}$; and
\item[(iii)] $0 \! \leqslant \! h_{M}^{o}(x) \! \leqslant \! 1$, $x \! \in \!
\mathbb{R}$.
\end{compactenum}
Note the decomposition $\iint_{\mathbb{R}^{2}}p^{o}(t,s;n) \, \md \mu_{k}^{o}
(t) \, \md \mu_{k}^{o}(s) \! = \! I_{a} \! + \! I_{b} \! + \! I_{c}$, where
\begin{align*}
I_{a} :=& \iint_{\mathbb{R}^{2}}p^{o}(t,s;n)(1 \! - \! h_{M}^{o}(s)) \, \md
\mu^{o}_{k}(t) \, \md \mu_{k}^{o}(s), \\
I_{b} :=& \iint_{\mathbb{R}^{2}}p^{o}(t,s;n)h_{M}^{o}(s)(1 \! - \! h_{M}^{o}
(t)) \, \md \mu^{o}_{k}(t) \, \md \mu_{k}^{o}(s), \\
I_{c} :=& \iint_{\mathbb{R}^{2}}p^{o}(t,s;n)h_{M}^{o}(t)h_{M}^{o}(s) \, \md
\mu^{o}_{k}(t) \, \md \mu_{k}^{o}(s).
\end{align*}
One shows that, for $n \! \in \! \mathbb{N}$,
\begin{align*}
\vert I_{a} \vert \leqslant& \iint_{\mathbb{R}^{2}} \vert p^{o}(t,s;n) \vert
(1 \! - \! h_{M}^{o}(s)) \, \md \mu_{k}^{o}(t) \, \md \mu_{k}^{o}(s) \\
\leqslant& \, \sup_{(t,s) \in \mathbb{R}^{2}} \vert p^{o}(t,s;n) \vert
\underbrace{\int_{\mathbb{R}} \md \mu_{k}^{o}(t)}_{= \, 1} \! \left(\int_{
\mathfrak{D}_{M_{o}}} \underbrace{(1 \! - \! h_{M}^{o}(s))}_{= \, 0} \, \md
\mu_{k}^{o}(s) \! + \! \int_{\mathbb{R} \setminus \mathfrak{D}_{M_{o}+1}}(1
\! - \! \underbrace{h_{M}^{o}(s)}_{= \, 0}) \, \md \mu_{k}^{o}(s) \right),
\end{align*}
whence
\begin{equation*}
\limsup_{k \to \infty} \vert I_{a} \vert \! \leqslant \! \sup_{(t,s) \in
\mathbb{R}^{2}} \vert p^{o}(t,s;n) \vert \underbrace{\limsup_{k \to \infty}
\int_{\mathbb{R} \setminus \mathfrak{D}_{M_{o}+1}} \, \md \mu_{k}^{o}(s)}_{
\leqslant \, \epsilon} \! \leqslant \! \epsilon \sup_{(t,s) \in \mathbb{R}^{
2}} \vert p^{o}(t,s;n) \vert;
\end{equation*}
similarly,
\begin{equation*}
\limsup_{k \to \infty} \vert I_{b} \vert \! \leqslant \! \epsilon \sup_{(t,s)
\in \mathbb{R}^{2}} \vert p^{o}(t,s;n) \vert.
\end{equation*}
Since, for $n \! \in \! \mathbb{N}$, $p^{o}(t,s;n)$ is continuous and bounded
on $\mathbb{R}^{2}$, there exists, by a generalisation of the
Stone-Weierstrass Theorem (for the single-variable case), a polynomial in two
variables (with $n$-dependent coefficients), $p(t,s;n)$, say, with $p(t,s;n)
\! = \! \sum_{i \geqslant i_{o}} \sum_{j \geqslant j_{o}} \gamma_{ij}(n)t^{i}
s^{j}$, such that $\vert p^{o}(t,s;n) \! - \! p(t,s;n) \vert \! \leqslant \!
\epsilon (n) \! := \! \epsilon$; thus,
\begin{equation*}
\vert h_{M}^{o}(t)h_{M}^{o}(s)p^{o}(t,s;n) \! - \! h_{M}^{o}(t)h_{M}^{o}(s)p
(t,s;n) \vert \! \leqslant \! \epsilon, \quad t,s \! \in \! \mathbb{R}.
\end{equation*}
Rewrite $I_{c}$ as
\begin{align*}
I_{c} &= \iint_{\mathbb{R}^{2}}h_{M}^{o}(s)h_{M}^{o}(t)p(t,s;n) \, \md \mu_{
k}^{o}(t) \, \md \mu_{k}^{o}(s) \! + \! \iint_{\mathbb{R}^{2}}h_{M}^{o}(s)h_{
M}^{o}(t)(p^{o}(t,s;n) \! - \! p(t,s;n)) \, \md \mu_{k}^{o}(t) \, \md \mu_{
k}^{o}(s) \\
&=: I_{c}^{\alpha} \! + \! I_{c}^{\beta}.
\end{align*}
One now shows that
\begin{align*}
\vert I_{c}^{\beta} \vert \leqslant& \iint_{\mathbb{R}^{2}}h_{M}^{o}(s)h_{
M}^{o}(t) \underbrace{\vert p^{o}(t,s;n) \! - \! p(t,s;n) \vert}_{\leqslant \,
\epsilon} \, \md \mu_{k}^{o}(t) \, \md \mu_{k}^{o}(s) \! \leqslant \! \epsilon
\int_{\mathbb{R}}h_{M}^{o}(s) \, \md \mu_{k}^{o}(s) \int_{\mathbb{R}}h_{M}^{o}
(t) \, \md \mu_{k}^{o}(t) \\
\leqslant& \, \epsilon \! \left(\int_{\mathfrak{D}_{M_{o}}} \underbrace{h_{
M}^{o}(s)}_{= \, 1} \, \md \mu_{k}^{o}(s) \! + \! \int_{\mathbb{R} \setminus
\mathfrak{D}_{M_{o}+1}} \underbrace{h_{M}^{o}(s)}_{= \, 0} \, \md \mu_{k}^{o}
(s) \right)^{2} \! \leqslant \! \epsilon \! \left(\int_{\mathfrak{D}_{M_{o}}}
\md \mu_{k}^{o}(s) \right)^{2} \\
\leqslant& \, \epsilon \left(\int_{\mathbb{R}} \md \mu_{k}^{o}(s) \right)^{2}
\! \leqslant \! \epsilon,
\end{align*}
and
\begin{align*}
I_{c}^{\alpha}=& \iint_{\mathbb{R}^{2}}h_{M}^{o}(s)h_{M}^{o}(t) \sum_{i
\geqslant i_{o}} \sum_{j \geqslant j_{o}} \gamma_{ij}(n)t^{i}s^{j} \, \md
\mu_{k}^{o}(t) \, \md \mu_{k}^{o}(s) \\
=& \, \sum_{i \geqslant i_{o}} \sum_{j \geqslant j_{o}} \gamma_{ij}(n) \!
\left(\int_{\mathbb{R}}h_{M}^{o}(t)t^{i} \, \md \mu_{k}^{o}(t) \right) \!
\left(\int_{\mathbb{R}}h_{M}^{o}(s)s^{j} \, \md \mu_{k}^{o}(s) \right) \\
\to& \, \sum_{i \geqslant i_{o}} \sum_{j \geqslant j_{o}} \gamma_{ij}(n) \!
\left(\int_{\mathbb{R}}h_{M}^{o}(t)t^{i} \, \md \mu^{o}(t) \right) \! \left(
\int_{\mathbb{R}}h_{M}^{o}(s)s^{j} \, \md \mu^{o}(s) \right) \qquad
(\text{since} \, \, \mu_{k}^{o} \! \overset{\ast}{\to} \! \mu^{o} \, \,
\text{as} \, \, k \! \to \! \infty) \\
=& \, \iint_{\mathbb{R}^{2}} \! \left(\sum_{i \geqslant i_{o}} \sum_{j
\geqslant j_{o}} \gamma_{ij}(n)t^{i}s^{j} \right) \! h_{M}^{o}(t)h_{M}^{o}(s)
\, \md \mu^{o}(t) \, \md \mu^{o}(s),
\end{align*}
whence, recalling that $p(t,s;n) \! = \! \sum_{i \geqslant i_{o}} \sum_{j
\geqslant j_{o}} \gamma_{ij}(n)t^{i}s^{j}$, it follows that
\begin{equation*}
I_{c}^{\alpha} \! = \! \iint_{\mathbb{R}^{2}}p(t,s;n)h_{M}^{o}(t)h_{M}^{o}(s)
\, \md \mu^{o}(t) \, \md \mu^{o}(s).
\end{equation*}
Furthermore, for $n \! \in \! \mathbb{N}$,
\begin{align*}
I_{c}^{\alpha} \leqslant& \iint_{\mathbb{R}^{2}}p^{o}(t,s;n)h_{M}^{o}(t)h_{M}^{
o}(s) \, \md \mu^{o}(t) \, \md \mu^{o}(s) \! + \! \epsilon \underbrace{\int_{
\mathbb{R}}h_{M}^{o}(t) \, \md \mu^{o}(t)}_{\leqslant \, 1} \, \, \underbrace{
\int_{\mathbb{R}}h_{M}^{o}(s) \, \md \mu^{o}(s)}_{\leqslant \, 1} \Rightarrow
\\
I_{c}^{\alpha} \leqslant& \, \iint_{\mathbb{R}^{2}}p^{o}(t,s;n) \vert 1 \! +
\! (h_{M}^{o}(t) \! - \! 1) \vert \vert 1 \! + \! (h_{M}^{o}(s) \! - \! 1)
\vert \, \md \mu^{o}(t) \, \md \mu^{o}(s) \! + \! \epsilon \\
\leqslant& \, \iint_{\mathbb{R}^{2}}p^{o}(t,s;n) \, \md \mu^{o}(t) \, \md
\mu^{o}(s) \! + \! \iint_{\mathbb{R}^{2}}p^{o}(t,s;n) \vert h_{M}^{o}(s) \! -
\! 1 \vert \, \md \mu^{o}(t) \, \md \mu^{o}(s) \! + \! \epsilon \\
+& \, \iint_{\mathbb{R}^{2}}p^{o}(t,s;n) \vert h_{M}^{o}(t) \! - \! 1 \vert \,
\md \mu^{o}(t) \, \md \mu^{o}(s) \! + \! \iint_{\mathbb{R}^{2}}p^{o}(t,s;n)
\vert h_{M}^{o}(t) \! - \! 1 \vert \vert h_{M}^{o}(s) \! - \! 1 \vert \, \md
\mu^{o}(t) \, \md \mu^{o}(s) \\
\leqslant& \, \iint_{\mathbb{R}^{2}}p^{o}(t,s;n) \, \md \mu^{o}(t) \, \md
\mu^{o}(s) \! + \! 2 \sup_{(t,s) \in \mathbb{R}^{2}} \vert p^{o}(t,s;n) \vert
\underbrace{\int_{(\mathbb{R} \setminus \mathfrak{D}_{M_{o}}) \cup \mathfrak{
D}_{M_{o}}} \vert h_{M}^{o}(s) \! - \! 1 \vert \, \md \mu^{o}(s)}_{\leqslant
\, \epsilon} \, \, \underbrace{\int_{\mathbb{R}} \md \mu^{o}(t)}_{= \, 1} \\
+& \, \sup_{(t,s) \in \mathbb{R}^{2}} \vert p^{o}(t,s;n) \vert \left(
\underbrace{\int_{(\mathbb{R} \setminus \mathfrak{D}_{M_{o}}) \cup \mathfrak{
D}_{M_{o}}} \vert h_{M}^{o}(t) \! - \! 1 \vert \, \md \mu^{o}(t)}_{\leqslant
\, \epsilon} \right)^{2} \! + \! \epsilon \\
\leqslant& \, \iint_{\mathbb{R}^{2}}p^{o}(t,s;n) \, \md \mu^{o}(t) \, \md \mu^{
o}(s) \! + \! \epsilon \! \left(1 \! + \! 2 \sup_{(t,s) \in \mathbb{R}^{2}}
\vert p^{o}(t,s;n) \vert \right) \! + \! \mathcal{O}(\epsilon^{2}),
\end{align*}
whereupon, neglecting the $\mathcal{O}(\epsilon^{2})$ term, and setting
$\varkappa^{\natural}_{n} \! := \! 1 \! + \! 2 \sup_{(t,s) \in \mathbb{R}^{2}}
\vert p^{o}(t,s;n) \vert$, one obtains
\begin{equation*}
I_{c}^{\alpha} \! \leqslant \! \iint_{\mathbb{R}^{2}}p^{o}(t,s;n) \, \md
\mu^{o}(t) \, \md \mu^{o}(s) \! + \! \varkappa^{\natural}_{n} \epsilon.
\end{equation*}
Hence, assembling the above-derived bounds for $I_{a}$, $I_{b}$, $I_{c}^{
\beta}$, and $I_{c}^{\alpha}$, one arrives at, for $n \! \in \! \mathbb{N}$,
upon setting $\epsilon^{\natural}_{n} \! := \! 2 \varkappa^{\natural}_{n}
\epsilon$,
\begin{equation*}
\iint_{\mathbb{R}^{2}}p^{o}(t,s;n) \, \md \mu_{k}^{o}(t) \, \md \mu^{o}_{k}(s)
\! - \! \iint_{\mathbb{R}^{2}}p^{o}(t,s;n) \, \md \mu^{o}(t) \, \md \mu^{o}(s)
\! \leqslant \! \epsilon^{\natural}_{n};
\end{equation*}
thus,
\begin{equation*}
\iint_{\mathbb{R}^{2}}p^{o}(t,s;n) \, \md \mu_{k}^{o}(t) \, \md \mu^{o}_{k}(s)
\! \to \! \iint_{\mathbb{R}^{2}}p^{o}(t,s;n) \, \md \mu^{o}(t) \, \md \mu^{o}
(s) \quad \text{as} \, \, k \! \to \! \infty.
\end{equation*}
Recalling that $p^{o}(t,s;n) \! := \! \min \left\lbrace q,K_{V_{m},n}^{o}(t,s)
\right\rbrace$, $(q,m) \! \in \! \mathbb{R} \times \mathbb{N}$, it follows
{}from the above analysis that, for $n \! \in \! \mathbb{N}$,
\begin{equation*}
\liminf_{k \to \infty} \mathrm{I}_{V}^{o}[\mu_{k}^{o}] \! \geqslant \! \iint_{
\mathbb{R}^{2}} \min \left\lbrace q,K_{V_{m},n}^{o}(t,s) \right\rbrace \md
\mu^{o}(t) \, \md \mu^{o}(s):
\end{equation*}
letting $q \! \uparrow \! \infty$ and $m \! \to \! \infty$, and using the
Monotone Convergence Theorem, one arrives at, upon noting that $\min
\left\lbrace q,K_{V_{m},n}^{o}(t,s) \right\rbrace \! \to \! K_{V,n}^{o}(t,s)$,
\begin{equation*}
\liminf_{k \to \infty} \mathrm{I}_{V}^{o}[\mu_{k}^{o}] \! \geqslant \! \iint_{
\mathbb{R}^{2}}K_{V,n}^{o}(t,s) \, \md \mu^{o}(t) \, \md \mu^{o}(s) \! = \!
\mathrm{I}_{V}^{o}[\mu^{o}], \quad \mu_{k}^{o},\mu^{o} \! \in \! \mathcal{
M}_{1}(\mathbb{R}).
\end{equation*}
Since, {}from the analysis above, it was shown that there exists a weakly 
$(\mathrm{weak}^{\ast})$ convergent subsequence (of probability measures) 
$\lbrace \mu_{m_{k}}^{o} \rbrace_{k=1}^{\infty}$ $(\subset \! \mathcal{M}_{1}
(\mathbb{R}))$ of $\lbrace \mu_{m}^{o} \rbrace_{m=1}^{\infty}$ $(\subset \! 
\mathcal{M}_{1}(\mathbb{R}))$ with a weak limit $\mu^{o} \! \in \! \mathcal{
M}_{1}(\mathbb{R})$, namely, $\mu_{m_{k}}^{o} \! \to \! \mu^{o}$ as $k \! \to 
\! \infty$, upon recalling that $\mathrm{I}_{V}^{o}[\mu_{m}^{o}] \! \leqslant 
\! E_{V}^{o} \! + \! \tfrac{1}{m}$, $m \! \in \! \mathbb{N}$, it follows 
that, in the limit as $m \! \to \! \infty$, $\mathrm{I}_{V}^{o}[\mu^{o}] \! 
\leqslant \! E_{V}^{o} \! := \! \inf \lbrace \mathstrut \mathrm{I}_{V}^{o}
[\mu^{o}]; \, \mu^{o} \! \in \! \mathcal{M}_{1}(\mathbb{R}) \rbrace$; {}from 
the latter two inequalities, it follows, thus, that $\exists \, \, \mu^{o} 
\! := \! \mu_{V}^{o} \! \in \! \mathcal{M}_{1}(\mathbb{R})$, the `odd' 
equilibrium measure, such that $\mathrm{I}_{V}^{o}[\mu_{V}^{o}] \! = \! \inf 
\lbrace \mathstrut \mathrm{I}_{V}^{o}[\mu^{o}]; \, \mu^{o} \! \in \! \mathcal{
M}_{1}(\mathbb{R}) \rbrace$, that is, the infimum is attained (the uniqueness 
of $\mu_{V}^{o} \! \in \! \mathcal{M}_{1}(\mathbb{R})$ is proven in Lemma~3.3 
below).

The compactness of $\operatorname{supp}(\mu_{V}^{o}) \! =: \! J_{o}$ is now 
established: actually, the following proof is true for any $\mu \! \in \! 
\mathcal{M}_{1}(\mathbb{R})$ achieving the above minimum; in particular, for 
$\mu \! = \! \mu_{V}^{o}$. Without loss of generality, therefore, let $\mu_{w} 
\! \in \! \mathcal{M}_{1}(\mathbb{R})$ be such that $\mathrm{I}_{V}^{o}[\mu_{
w}] \! = \! E_{V}^{o}$, and let $D$ be any proper subset of $\mathbb{R}$ for 
which $\mu_{w}(D) \! := \! \int_{D} \md \mu_{w}(s) \! > \! 0$. As in 
\cite{a80}, set
\begin{equation*}
\mu_{w}^{\varepsilon}(z) \! := \! \left(1 \! + \! \varepsilon \mu_{w}(D)
\right)^{-1}(\mu_{w}(z) \! + \! \varepsilon (\mu_{w} \! \! \upharpoonright_{
D})(z)),  \quad \varepsilon \! \in \! (-1,1),
\end{equation*}
where $\mu_{w} \! \! \upharpoonright_{D}$ denotes the restriction of $\mu_{w}$
to $D$ (note, also, that $\mu^{\varepsilon}_{w} \! > \! 0$ and bounded, and
$\int_{\mathbb{R}} \md \mu_{w}^{\epsilon}(s) \! = \! 1)$. Using the fact that
$K_{V,n}^{o}(s,t) \! = \! K_{V,n}^{o}(t,s)$, one shows that, for $n \! \in \!
\mathbb{N}$,
\begin{align*}
\mathrm{I}_{V}^{o}[\mu_{w}^{\varepsilon}] =& \iint_{\mathbb{R}^{2}}K_{V,n}^{o}
(s,t) \, \md \mu^{\varepsilon}_{w}(s) \, \md \mu^{\varepsilon}_{w}(t) \\
=& \, (1 \! + \! \varepsilon \mu_{w}(D))^{-2} \iint_{\mathbb{R}^{2}}K_{V,n}^{o}
(s,t)(\md \mu_{w}(s) \! + \! \varepsilon \md (\mu_{w} \! \! \upharpoonright_{
D})(s))(\md \mu_{w}(t) \! + \! \varepsilon \md (\mu_{w} \! \! \upharpoonright_{
D})(t)) \\
=& \, (1 \! + \! \varepsilon \mu_{w}(D))^{-2} \! \left(\mathrm{I}_{V}^{o}[
\mu_{w}] \! + \! 2 \varepsilon \iint_{\mathbb{R}^{2}}K_{V,n}^{o}(s,t) \, \md
\mu_{w}(s) \, \md (\mu_{w} \! \! \upharpoonright_{D})(t) \right. \\
+& \left. \, \varepsilon^{2} \iint_{\mathbb{R}^{2}}K_{V,n}^{o}(s,t) \, \md
(\mu_{w} \! \! \upharpoonright_{D})(t) \, \md (\mu_{w} \! \! \upharpoonright_{
D})(s) \right).
\end{align*}
(Note that all the above integrals are finite due to the argument at the 
beginning of the proof.) By the minimal property of $\mu_{w} \! \in \! 
\mathcal{M}_{1}(\mathbb{R})$, it follows that $\partial_{\varepsilon} \mathrm{
I}_{V}^{o}[\mu_{w}^{\varepsilon}] \! = \! 0$, which implies that, for $n \!
\in \! \mathbb{N}$,
\begin{equation*}
\iint_{\mathbb{R}^{2}}(K_{V,n}^{o}(s,t) \! - \! \mathrm{I}_{V}^{o}[\mu_{w}])
\, \md \mu_{w}(s) \, \md (\mu_{w} \! \! \upharpoonright_{D})(t) \! = \! 0;
\end{equation*}
but, recalling that, for $\widehat{\psi}_{V}^{o}(z) \! := \! 2 \widetilde{V}
(z) \! - \! (1 \! + \! \tfrac{1}{n}) \ln (z^{2} \! + \! 1) \! - \! \ln (z^{-2}
\! + \! 1)$, $K_{V,n}^{o}(t,s) \! \geqslant \! \tfrac{1}{2} \psi_{V}^{o}(s) \!
+ \! \tfrac{1}{2} \psi_{V}^{o}(t)$, it follows from the above that,
\begin{gather*}
\iint_{\mathbb{R}^{2}} \mathrm{I}_{V}^{o}[\mu_{w}] \, \md \mu_{w}(s) \, \md
(\mu_{w} \! \! \upharpoonright_{D})(t) \! \geqslant \! \iint_{\mathbb{R}^{2}}
\! \left(\tfrac{1}{2} \widehat{\psi}_{V}^{o}(s) \! + \! \tfrac{1}{2} \widehat{
\psi}_{V}^{o}(t) \right) \! \md \mu_{w}(s) \, \md (\mu_{w} \! \!
\upharpoonright_{D})(t) \Rightarrow \\
0 \! \geqslant \! \iint_{\mathbb{R}^{2}} \! \left(\tfrac{1}{2} \widehat{
\psi}_{V}^{o}(s) \! + \! \tfrac{1}{2} \widehat{\psi}_{V}^{o}(t) \! - \!
\mathrm{I}_{V}^{o}[\mu_{w}] \right) \! \md \mu_{w}(s) \, \md (\mu_{w} \! \!
\upharpoonright_{D})(t),
\end{gather*}
whence
\begin{equation*}
\int_{\mathbb{R}} \! \left(\widehat{\psi}_{V}^{o}(t) \! + \! \left(\int_{
\mathbb{R}} \widehat{\psi}_{V}^{o}(s) \, \md \mu_{w}(s) \right) \! - \! 2
\mathrm{I}_{V}^{o}[\mu_{w}] \right) \! \md (\mu_{w} \! \! \upharpoonright_{D}
)(t) \! \leqslant \! 0.
\end{equation*}
Recalling that
\begin{equation*}
\widehat{\psi}_{V}^{o}(x) \! := \! 2 \widetilde{V}(x) \! - \! \left(1 \! + \!
\dfrac{1}{n} \right) \! \ln (x^{2} \! + \! 1) \! - \! \ln (x^{-2} \! + \! 1)
\! = \!
\begin{cases}
+\infty, &\text{$\vert x \vert \! \to \! \infty$,} \\
+\infty, &\text{$\vert x \vert \! \to \! 0$,}
\end{cases}
\end{equation*}
it follows that, $\exists \, \, T_{m} \! := \! T_{m}(n) \! > \! 1$ such that
\begin{equation*}
\widehat{\psi}_{V}^{o}(t) \! + \! \int_{\mathbb{R}} \widehat{\psi}_{V}^{o}(s)
\, \md \mu_{w}(s) \! - \! 2 \mathrm{I}_{V}^{o}[\mu_{w}] \! \geqslant \! 1
\quad \text{for} \quad t \! \in \! \left((-T_{m},-T_{m}^{-1}) \cup (T_{m}^{-
1},T_{m}) \right)^{c}
\end{equation*}
(note, also, that $+\infty \! > \! \mathrm{I}_{V}^{o}[\mu_{w}] \! = \! \iint_{
\mathbb{R}^{2}}K_{V,n}^{o}(t,s) \, \md \mu_{w}(t) \, \md \mu_{w}(s) \! = \! 
\int_{\mathbb{R}} \widehat{\psi}_{V}^{o}(\xi) \, \md \mu_{w}(\xi) \! =$ a 
finite real number). Hence, if $D$ $(\subset \mathbb{R})$ is such that $D 
\subset (\lbrace \vert x \vert \! \geqslant \! T_{m} \rbrace \cup \lbrace 
\vert x \vert \! \leqslant \! T_{m}^{-1} \rbrace)$, $T_{m} \! > \! 1$, it 
follows {}from the above calculations that, for $n \! \in \! \mathbb{N}$,
\begin{equation*}
0 \! \geqslant \! \int_{\mathbb{R}} \! \left(\widehat{\psi}_{V}^{o}(t) \! + \!
\left(\int_{\mathbb{R}} \widehat{\psi}_{V}^{o}(s) \, \md \mu_{w}(s) \right)
\! - \! 2 \mathrm{I}_{V}^{o}[\mu_{w}] \right) \! \md (\mu_{w} \! \!
\upharpoonright_{D})(t) \! \geqslant \! 1,
\end{equation*}
which is a contradiction; hence, $\operatorname{supp}(\mu_{w}) \subseteq [-
T_{m},-T_{m}^{-1}] \cup [T_{m}^{-1},T_{m}]$, $T_{m} \! > \! 1$; in particular, 
$J_{o} \! := \! \operatorname{supp}(\mu_{V}^{o}) \subseteq [-T_{m},-T_{m}^{-
1}] \cup [T_{m}^{-1},T_{m}]$, $T_{m} \! > \! 1$, which establishes the 
compactness of the support of the `odd' equilibrium measure $\mu_{V}^{o} \! 
\in \! \mathcal{M}_{1}(\mathbb{R})$. Furthermore, it is worth noting that, 
since $J_{o} \! := \! \operatorname{supp}(\mu_{V}^{o}) \! =$ compact $(\subset 
\overline{\mathbb{R}} \setminus \{0,\pm \infty\})$, and $\widetilde{V} \colon 
\mathbb{R} \setminus \{0\} \! \to \! \mathbb{R}$ is real analytic on $J_{o}$, 
for $n \! \in \! \mathbb{N}$,
\begin{align*}
+\infty \! > \! E_{V}^{o} \, (=& \, \mathrm{I}_{V}^{o}[\mu_{V}^{o}]) \,
\geqslant \, \iint_{\mathbb{R}^{2}} \ln \! \left(\vert s \! - \! t \vert^{2+
\frac{1}{n}} \vert st \vert^{-1}w^{o}(s)w^{o}(t) \right)^{-1} \md \mu_{V}^{o}
(s) \, \md \mu_{V}^{o}(t) \\
=& \, \iint_{J_{o}^{2}} \ln \! \left(\vert s \! - \! t \vert^{2+\frac{1}{n}}
\vert st \vert^{-1}w^{o}(s)w^{o}(t) \right)^{-1} \md \mu_{V}^{o}(s) \, \md
\mu_{V}^{o}(t) \! > \! -\infty;
\end{align*}
moreover, a straightforward consequence of the fact just established is that
$J_{o}$ has positive logarithmic capacity, that is, $\operatorname{cap}
(J_{o}) \! = \! \exp (-E_{V}^{o}) \! > \! 0$. \hfill $\qed$
\begin{eeeee}
It is important to note from the latter part of the proof of Lemma~3.1 that 
$J_{o} \! \not\supseteq \! \lbrace 0,\pm \infty \rbrace$. This can also be 
seen as follows. For $\varepsilon$ some arbitrarily fixed, sufficiently small 
positive real number and $\Sigma_{\varepsilon} \! := \! \lbrace \mathstrut 
z; \, w^{o}(z) \! \geqslant \! \varepsilon \rbrace$, if $(s,t) \! \notin \! 
\Sigma_{\varepsilon} \times \Sigma_{\varepsilon}$, then, for $n \! \in \! 
\mathbb{N}$, $\ln (\vert s \! - \! t \vert^{2+\frac{1}{n}} \vert st \vert^{-1}
w^{o}(s)w^{o}(t))^{-1} \! =: \! K_{V,n}^{o}(s,t)$ $(= \! K_{V,n}^{o}(t,s))$ 
$> \! E_{V}^{o} \! + \! 1$, which is a contradiction, since it was established 
above that the minimum is attained $\Leftrightarrow \! (s,t) \! \in \! 
\Sigma_{\varepsilon} \times \Sigma_{\varepsilon}$. Towards this end, it is 
enough to show that (see, for example, \cite{a43}), if $\lbrace (s_{m},t_{m}) 
\rbrace_{m=1}^{\infty}$ is a sequence with $\lim \min_{m \to \infty} \lbrace 
w^{o}(s_{m}),w^{o}(t_{m}) \rbrace \! = \! 0$, then, for $n \! \in \! \mathbb{
N}$, $\lim_{m \to \infty} \ln (\vert s_{m} \! - \! t_{m} \vert^{2+\frac{1}{n}
} \vert s_{m}t_{m} \vert^{-1}w^{o}(s_{m})w^{o}(t_{m}))^{-1} \! = \! \lim_{m 
\to \infty}K_{V,n}^{o}(s_{m},t_{m}) \! = \! +\infty$. Without loss of 
generality, one can assume that $s_{m} \! \to \! s$ and $t_{m} \! \to \! t$ 
as $m \! \to \! \infty$, where $s$, $t$, or both may be infinite; thus, there 
are several cases to consider:
\begin{compactenum}
\item[(i)] if $s$ and $t$ are finite, then, from $\lim \min_{m \to \infty}
\lbrace w^{o}(s_{m}),w^{o}(t_{m}) \rbrace \! = \! \min \lbrace w^{o}(s),w^{o}
(t) \rbrace \! = \! 0$, it is clear that $\lim_{m \to \infty}K_{V,n}^{o}(s_{
m},t_{m}) \! = \! +\infty$;
\item[(ii)] if $\vert s \vert \! = \! \infty$ (resp., $\vert t \vert \! = \!
\infty)$ but $t \! = \! \text{finite}$ (resp., $s \! = \! \text{finite})$,
then, due to the fact that $\widetilde{V} \colon \mathbb{R} \setminus \{0\} \!
\to \! \mathbb{R}$ satisfies (for $n \! \in \! \mathbb{N})$ the conditions
\begin{equation*}
2 \widetilde{V}(x) \! - \! \left(1 \! + \! \dfrac{1}{n} \right) \! \ln (x^{2}
\! + \! 1) \! - \! \ln (x^{-2} \! + \! 1) \! = \!
\begin{cases}
+\infty, &\text{$\vert x \vert \! \to \! \infty$,} \\
+\infty, &\text{$\vert x \vert \! \to \! 0$,}
\end{cases}
\end{equation*}
it follows that $\lim_{m \to \infty}K_{V,n}^{o}(s_{m},t_{m}) \! = \! +\infty$;
\item[(iii)] if $\vert s \vert \! = \! 0$ (resp., $\vert t \vert \! = \! 0)$
but $t \! = \! \text{finite}$ (resp., $s \! = \! \text{finite})$, then, as
a result of the above conditions for $\widetilde{V}$, it follows that $\lim_{
m \to \infty}K_{V,n}^{o}(s_{m},t_{m}) \! = \! +\infty$;
\item[(iv)] if $\vert s \vert \! = \! \infty$ and $\vert t \vert \! = \!
\infty$, then, again due to the above conditions for $\widetilde{V}$, it
follows that $\lim_{m \to \infty}K_{V,n}^{o} \linebreak[4]
(s_{m},t_{m}) \! = \! +\infty$; and
\item[(v)] if $\vert s \vert \! = \! 0$ and $\vert t \vert \! = \! 0$, then,
again, as above, it follows that $\lim_{m \to \infty}K_{V,n}^{o}(s_{m},t_{m})
\! = \! +\infty$.
\end{compactenum}
Hence, for $n \! \in \! \mathbb{N}$, $K_{V,n}^{o}(s,t) \! > \! E_{V}^{o} \! 
+ \! 1$ if $(s,t) \! \notin \! \Sigma_{\varepsilon} \times \Sigma_{
\varepsilon}$, that is, if $s$, $t$, or both $\in \! \lbrace 0,\pm \infty 
\rbrace$ (which can not be the case, as the infimum $E_{V}^{o}$ is attained 
$\Leftrightarrow \! (s,t) \! \in \! \Sigma_{\varepsilon} \times \Sigma_{
\varepsilon}$, whence $\operatorname{supp}(\mu_{V}^{o}) \! =: \! J_{o} \! 
\not\supseteq \! \{0,\pm \infty\})$. \hfill $\blacksquare$
\end{eeeee}

In order to demonstrate the uniqueness of the `odd' equilibrium measure, 
$\mu_{V}^{o}$ $(\in \! \mathcal{M}_{1}(\mathbb{R}))$, the following lemma is 
requisite.
\begin{ccccc}
Let $\mu \! := \! \mu_{1} \! - \! \mu_{2}$, where $\mu_{1},\mu_{2}$ are 
non-negative, finite-moment $(\int_{\operatorname{supp}(\mu_{j})}s^{m} \, 
\md \mu_{j}(s) \! < \! \infty$, $m \! \in \! \mathbb{Z}$, $j \! = \! 1,2)$ 
measures on $\mathbb{R}$ supported on distinct sets $(\operatorname{supp}
(\mu_{1}) \cap \operatorname{supp}(\mu_{2}) \! = \! \varnothing)$, be the 
(unique) Jordan decomposition of the finite-moment signed measure on $\mathbb{
R}$ with mean zero, that is, $\int_{\operatorname{supp}(\mu)} \md \mu (s) \! 
= \! 0$, and with $\operatorname{supp}(\mu) \! = \! \text{compact}$. Suppose 
that $-\infty \! < \! \iint_{\mathbb{R}^{2}} \ln (\vert s \! - \! t \vert^{-
(2+\frac{1}{n})} \vert st \vert) \, \md \mu_{j}(s) \, \md \mu_{j}(t) \! < \! 
+\infty$, $n \! \in \! \mathbb{N}$, $j \! =\! 1,2$. Then, for $n \! \in \! 
\mathbb{N}$,
\begin{equation*}
\iint_{\mathbb{R}^{2}} \ln \! \left(\dfrac{\vert st \vert}{\vert s \! - \! t
\vert^{2+\frac{1}{n}}} \right) \! \md \mu (s) \, \md \mu (t) \! = \! \iint_{
\mathbb{R}^{2}} \ln \! \left(\dfrac{\vert s \! - \! t \vert^{2+\frac{1}{n}}}{
\vert st \vert}w^{o}(s)w^{o}(t) \right)^{-1} \md \mu (s) \, \md \mu (t) \!
\geqslant \! 0,
\end{equation*}
where equality holds if, and only if, $\mu \! = \! 0$.
\end{ccccc}

\emph{Proof.} Recall the following identity \cite{a79}\footnote{One could also
carry out the proof via the following identity: for $s \! \in \! \mathbb{R}$,
$n \! \in \! \mathbb{N}$, and any $\epsilon \! > \! 0$, $\ln \! \left(s^{2+
1/n} \! + \! \epsilon^{2} \right) \! = \! \ln (\epsilon^{2}) \! + \! 2 \,
\Im \! \left(\int_{0}^{+\infty}(\mi u)^{-1}(\me^{\mi us^{1+1/2n}}-1) \me^{-
\epsilon u} \, \md u \right)$.} (see pg.~147, Equation~(6.44)): for $\xi \!
\in \! \mathbb{R}$ and any $\varepsilon \! > \! 0$,
\begin{equation*}
\ln (\xi^{2} \! + \! \varepsilon^{2}) \! = \! \ln (\varepsilon^{2}) \! + \! 2
\, \Im \! \left(\int_{0}^{+\infty} \! \left(\dfrac{\me^{\mi \xi v} \! - \!
1}{\mi v} \right) \! \me^{-\varepsilon v} \, \md v \right);
\end{equation*}
thus, it follows that, for $n \! \in \! \mathbb{N}$,
\begin{align*}
\left(1 \! + \! \dfrac{1}{2n} \right) \! \iint_{\mathbb{R}^{2}} \ln ((s \! 
- \! t)^{2}+ &\varepsilon^{2}) \, \md \mu (s) \, \md \mu (t) \! = \! \left(1 
\! + \! \dfrac{1}{2n} \right) \! \iint_{\mathbb{R}^{2}} \ln (\varepsilon^{2}) 
\, \md \mu (s) \, \md \mu (t) \\
+& \, 2 \! \left(1 \! + \! \dfrac{1}{2n} \right) \! \iint_{\mathbb{R}^{2}} \! 
\left(\Im \! \left(\int_{0}^{+\infty} \! \left(\dfrac{\me^{\mi (s-t)v} \! - 
\! 1}{\mi v} \right) \! \me^{-\varepsilon v} \, \md v \right) \! \right) \! 
\md \mu (s) \, \md \mu (t), \\
\iint_{\mathbb{R}^{2}} \ln (s^{2} \! + \! \varepsilon^{2}) \, &\md \mu (s) \, 
\md \mu (t) = \iint_{\mathbb{R}^{2}} \ln (\varepsilon^{2}) \, \md \mu (s) \, 
\md \mu (t) \\
+& \, \iint_{\mathbb{R}^{2}} \! \left(2 \, \Im \! \left(\int_{0}^{+\infty} \!
\left(\dfrac{\me^{\mi sv} \! - \! 1}{\mi v} \right) \! \me^{-\varepsilon v} \,
\md v \right) \! \right) \! \md \mu (s) \, \md \mu (t), \\
\iint_{\mathbb{R}^{2}} \ln (t^{2} \! + \! \varepsilon^{2}) \, &\md \mu (s) \,
\md \mu (t) = \iint_{\mathbb{R}^{2}} \ln (\varepsilon^{2}) \, \md \mu (s)
\, \md \mu (t) \\
+& \, \iint_{\mathbb{R}^{2}} \! \left(2 \, \Im \! \left(\int_{0}^{+\infty} \!
\left(\dfrac{\me^{\mi tv} \! - \! 1}{\mi v} \right) \! \me^{-\varepsilon v} \,
\md v \right) \! \right) \! \md \mu (s) \, \md \mu (t);
\end{align*}
but, since $\iint_{\mathbb{R}^{2}} \md \mu (s) \, \md \mu (t) \! = \! \left(
\int_{\mathbb{R}} \md \mu (s) \right)^{2} \! = \! 0$, one obtains, after some 
rearrangement,
\begin{align*}
\left(1 \! + \! \dfrac{1}{2n} \right) \! \iint_{\mathbb{R}^{2}} \ln ((s \! -
\! t)^{2} \! + \! \varepsilon^{2}) \, \md \mu (s) \, \md \mu (t) \! =& \, 2
\! \left(1 \! + \! \dfrac{1}{2n} \right) \! \Im \! \left(\int_{0}^{+\infty}
\me^{-\varepsilon v} \! \left(\iint_{\mathbb{R}^{2}} \! \left(\dfrac{\me^{\mi
(s-t)v} \! - \! 1}{\mi v} \right) \right. \right. \\
\times&\left. \left. \md \mu (s) \, \md \mu (t) \right) \! \md v \right), \\
\iint_{\mathbb{R}^{2}} \ln (s^{2} \! + \! \varepsilon^{2}) \, \md \mu (s) \,
\md \mu (t) \! =& \, 2 \, \Im \left(\int_{0}^{+\infty} \me^{-\varepsilon v} \!
\left(\iint_{\mathbb{R}^{2}} \! \left(\dfrac{\me^{\mi sv} \! - \! 1}{\mi v}
\right) \! \md \mu (s) \, \md \mu (t) \right) \! \md v \right), \\
\iint_{\mathbb{R}^{2}} \ln (t^{2} \! + \! \varepsilon^{2}) \, \md \mu (s) \,
\md \mu (t) \! =& \, 2 \, \Im \left(\int_{0}^{+\infty} \me^{-\varepsilon v} \!
\left(\iint_{\mathbb{R}^{2}} \! \left(\dfrac{\me^{\mi tv} \! - \! 1}{\mi v}
\right) \! \md \mu (s) \, \md \mu (t) \right) \! \md v \right).
\end{align*}
Noting that
\begin{align*}
\iint_{\mathbb{R}^{2}} \! \left(\dfrac{\me^{\mi (s-t)v} \! - \! 1}{\mi v}
\right) \! \md \mu (s) \, \md \mu (t) =& \, \dfrac{1}{\mi v} \iint_{\mathbb{
R}^{2}} \me^{\mi (s-t)v} \, \md \mu (s) \, \md \mu (t) \! - \! \dfrac{1}{\mi
v} \underbrace{\iint_{\mathbb{R}^{2}} \md \mu (s) \, \md \mu (t)}_{= \, 0} \\
=& \, \dfrac{1}{\mi v} \int_{\mathbb{R}} \me^{\mi sv} \, \md \mu (s) \int_{
\mathbb{R}} \me^{-\mi tv} \, \md \mu (t),
\end{align*}
and setting $\widehat{\mu}(z) \! := \! \int_{\mathbb{R}} \me^{\mi \xi z} \,
\md \mu (\xi)$, one gets that
\begin{equation*}
\iint_{\mathbb{R}^{2}} \! \left(\dfrac{\me^{\mi (s-t)v} \! - \! 1}{\mi v}
\right) \! \md \mu (s) \, \md \mu (t) \! = \! \dfrac{1}{\mi v} \vert \widehat{
\mu}(v) \vert^{2}:
\end{equation*}
also,
\begin{equation*}
\iint_{\mathbb{R}^{2}} \! \left(\dfrac{\me^{\mi sv} \! - \! 1}{\mi v} \right)
\! \md \mu (s) \, \md \mu (t) \! = \! \dfrac{1}{\mi v} \int_{\mathbb{R}} \me^{
\mi sv} \, \md \mu (s) \underbrace{\int_{\mathbb{R}} \md \mu (t)}_{= \, 0}-
\dfrac{1}{\mi v} \underbrace{\int_{\mathbb{R}} \md \mu (s)}_{= \, 0} \, \,
\underbrace{\int_{\mathbb{R}} \md \mu (t)}_{= \, 0} \! = \! 0;
\end{equation*}
similarly,
\begin{equation*}
\iint_{\mathbb{R}^{2}} \! \left(\dfrac{\me^{\mi tv} \! - \! 1}{\mi v} \right)
\! \md \mu (s) \, \md \mu (t) \! = \! 0.
\end{equation*}
Hence,
\begin{gather*}
\left(1 \! + \! \dfrac{1}{2n} \right) \! \iint_{\mathbb{R}^{2}} \ln ((s \! -
\! t)^{2} \! + \! \varepsilon^{2}) \, \md \mu (s) \, \md \mu (t) \! = \! 2 \!
\left(1 \! + \! \dfrac{1}{2n} \right) \! \Im \! \left(\int_{0}^{+\infty}
\dfrac{\vert \widehat{\mu}(v) \vert^{2}}{\mi v} \me^{-\varepsilon v} \, \md v
\right), \\
\iint_{\mathbb{R}^{2}} \ln (s^{2} \! + \! \varepsilon^{2}) \, \md \mu (s) \,
\md \mu (t) \! = \! \iint_{\mathbb{R}^{2}} \ln (t^{2} \! + \! \varepsilon^{2})
\, \md \mu (s) \, \md \mu (t) \! = \! 0.
\end{gather*}
Noting that $\widehat{\mu}(0) \! = \! \int_{\mathbb{R}} \md \mu (\xi) \! =
\! 0$, a Taylor expansion about $v \! = \! 0$ shows that $\widehat{\mu}(v) \!
=_{v \to 0} \! \widehat{\mu}^{\prime}(0)v \! + \! \mathcal{O}(v^{2})$, where
$\widehat{\mu}^{\prime}(0) \! := \! \partial_{v} \widehat{\mu}(v) \vert_{v=
0}$; thus, $v^{-1} \vert \widehat{\mu}(v) \vert^{2} \! =_{v \to 0} \! \vert
\widehat{\mu}^{\prime}(0) \vert^{2}v \! + \! \mathcal{O}(v^{2})$, which means
that there is no singularity in the integrand as $v \! \to \! 0$ (in fact,
$v^{-1} \vert \widehat{\mu}(v) \vert^{2}$ is real analytic in a neighbourhood
of the origin), whence
\begin{equation*}
\left(1 \! + \! \dfrac{1}{2n} \right) \! \iint_{\mathbb{R}^{2}} \ln ((s \! -
\! t)^{2} \! + \! \varepsilon^{2}) \, \md \mu (s) \, \md \mu (t) \! = \! -2
\! \left(1 \! + \! \dfrac{1}{2n} \right) \! \int_{0}^{+\infty}v^{-1} \vert
\widehat{\mu}(v) \vert^{2} \me^{-\varepsilon v} \, \md v.
\end{equation*}
Recalling that $\iint_{\mathbb{R}^{2}} \ln (\ast^{2} \! + \! \varepsilon^{2})
\, \md \mu (s) \, \md \mu (t) \! = \! 0$, $\ast \! \in \! \{s,t\}$, and
adding, it follows that
\begin{equation*}
\iint_{\mathbb{R}^{2}} \ln \! \left(\dfrac{(s^{2} \! + \! \varepsilon^{2})^{
1/2}(t^{2} \! + \! \varepsilon^{2})^{1/2}}{((s \! - \! t)^{2} \! + \!
\varepsilon^{2})^{1+\frac{1}{2n}}} \right) \! \md \mu (s) \, \md \mu (t) \!
= \! 2 \! \left(1 \! + \! \dfrac{1}{2n} \right) \! \int_{0}^{+\infty}v^{-1}
\vert \widehat{\mu}(v) \vert^{2} \me^{-\varepsilon v} \, \md v.
\end{equation*}
Now, using the fact that $\ln ((s \! - \! t)^{2} \! + \! \varepsilon^{2})^{-
1}$ (resp., $\ln (s^{2} \! + \! \varepsilon^{2})^{1/2}$ and $\ln (t^{2} \!
+ \! \varepsilon^{2})^{1/2})$ is (resp., are) bounded below (resp., above)
uniformly with respect to $\varepsilon$ and that the measures have compact
support, letting $\varepsilon \! \downarrow \! 0$ and using the Monotone
Convergence Theorem, one arrives at
\begin{align*}
\left(1 \! + \! \dfrac{1}{2n} \right) \! \iint_{\mathbb{R}^{2}} \ln \! \left(
\dfrac{(s^{2} \! + \! \varepsilon^{2})^{1/2}(t^{2} \! + \! \varepsilon^{2})^{
1/2}}{((s \! - \! t)^{2} \! + \! \varepsilon^{2})^{1+\frac{1}{2n}}} \right) \!
\md \mu (s) \, \md \mu (t) \underset{\varepsilon \downarrow 0}{=}& \, \iint_{
\mathbb{R}^{2}} \ln \! \left(\dfrac{\vert st \vert}{\vert s \! - \! t \vert^{
2+\frac{1}{n}}} \right) \! \md \mu (s) \, \md \mu (t) \\
=& \, 2 \! \left(1 \! + \! \dfrac{1}{2n} \right) \! \int_{0}^{+\infty}v^{-1}
\vert \widehat{\mu}(v) \vert^{2} \, \md v \! \geqslant \! 0,
\end{align*}
where, trivially, equality holds if, and only if, $\mu \! = \! 0$.
Furthermore, noting that, since $\int_{\mathbb{R}} \md \mu (\xi) \! = \! 0$,
$\iint_{\mathbb{R}^{2}} \ln (w^{e}(\ast))^{-1} \, \md \mu (s) \, \md \mu (t)
\! = \! 0$, $\ast \! \in \! \{s,t\}$, letting $\varepsilon \! \downarrow \!
0$ and using monotone convergence, one also arrives at
\begin{align*}
\iint_{\mathbb{R}^{2}} \ln \! \left(\dfrac{(s^{2} \! + \! \varepsilon^{2})^{
1/2}(t^{2} \! + \! \varepsilon^{2})^{1/2}}{((s \! - \! t)^{2} \! + \!
\varepsilon^{2})^{1+\frac{1}{2n}}w^{o}(s)w^{o}(t)} \right) \! \md \mu (s) \,
\md \mu (t) \underset{\varepsilon \downarrow 0}{=}& \, \iint_{\mathbb{R}^{2}}
\ln \! \left(\dfrac{\vert st \vert}{\vert s \! - \! t \vert^{2+\frac{1}{n}}
w^{o}(s)w^{o}(t)} \right) \md \mu (s) \, \md \mu (t) \\
=& \, 2 \! \left(1 \! + \! \dfrac{1}{2n} \right) \! \int_{0}^{+\infty}v^{-1}
\vert \widehat{\mu}(v) \vert^{2} \, \md v \! \geqslant \! 0,
\end{align*}
where, again, and trivially, equality holds if, and only if, $\mu \! = \! 0$.
\hfill $\square$

The uniqueness of $\mu_{V}^{o}$ $(\in \! \mathcal{M}_{1}(\mathbb{R}))$ will
now be established.
\begin{ccccc}
Let the external field $\widetilde{V} \colon \mathbb{R} \setminus \{0\} \! \to
\! \mathbb{R}$ satisfy conditions~{\rm (2.3)--(2.5)}. Set $w^{o}(z) \! := \!
\exp (-\widetilde{V}(z))$, and define, for $n \! \in \! \mathbb{N}$,
\begin{equation*}
\mathrm{I}_{V}^{o}[\mu^{o}] \colon \mathcal{M}_{1}(\mathbb{R}) \! \to \!
\mathbb{R}, \, \, \mu^{o} \! \mapsto \! \iint_{\mathbb{R}^{2}} \ln \! \left(
\vert s \! - \! t \vert^{2+\frac{1}{n}} \vert st \vert^{-1}w^{o}(s)w^{o}(t)
\right)^{-1} \md \mu^{o}(s) \, \md \mu^{o}(t),
\end{equation*}
and consider the minimisation problem $E_{V}^{o} \! = \! \inf \left\lbrace 
\mathstrut \mathrm{I}_{V}^{o}[\mu^{o}]; \, \mu^{o} \! \in \! \mathcal{M}_{1}
(\mathbb{R}) \right\rbrace$. Then, $\exists ! \, \, \mu_{V}^{o} \in \! 
\mathcal{M}_{1}(\mathbb{R})$ such that $\mathrm{I}_{V}^{o}[\mu_{V}^{o}] \! = 
\! E_{V}^{o}$.
\end{ccccc}

\emph{Proof.} It was shown in Lemma~3.1 that $\exists \, \, \mu_{V}^{o} \! 
\in \! \mathcal{M}_{1}(\mathbb{R})$, the `odd' equilibrium measure, such that 
$\mathrm{I}_{V}^{o}[\mu^{o}] \! = \! E_{V}^{o}$; therefore, it remains to 
establish the uniqueness of the `odd' equilibrium measure. Let $\widetilde{
\mu}_{V}^{o} \! \in \! \mathcal{M}_{1}(\mathbb{R})$ be a second probability 
measure for which $\mathrm{I}_{V}^{o}[\widetilde{\mu}_{V}^{o}] \! = \! E_{
V}^{o} \! = \! \mathrm{I}_{V}^{o}[\mu_{V}^{o}]$: the argument in Lemma~3.1 
shows that $\widetilde{J}_{o} \! := \! \operatorname{supp}(\widetilde{\mu}_{
V}^{o}) \! = \! \text{compact} \subset \! \overline{\mathbb{R}} \setminus 
\lbrace 0,\pm \infty \rbrace$, and that $\mathrm{I}_{V}^{o}[\widetilde{\mu}_{
V}^{o}] \! < \! +\infty$. Define the finite-moment, signed measure $\mu^{
\sharp} \! := \! \widetilde{\mu}_{V}^{o} \! - \! \mu_{V}^{o}$, where 
$\widetilde{\mu}_{V}^{o},\mu_{V}^{o} \! \in \! \mathcal{M}_{1}(\mathbb{R})$, 
and $\widetilde{J}_{o} \cap J_{o} \! = \! \varnothing$, with (cf. Lemma~3.1), 
$J_{o} \! = \! \operatorname{supp}(\mu_{V}^{o}) \! = \! \text{compact} \! 
\subset \overline{\mathbb{R}} \setminus \lbrace 0,\pm \infty \rbrace$; thus, 
{}from Lemma~3.2 (with $\mu \! \to \! \mu^{\sharp})$, namely,
\begin{equation*}
\iint_{\mathbb{R}^{2}} \ln \! \left(\vert s \! - \! t \vert^{-(2+\frac{1}{n})
} \vert st \vert \right) \md \mu^{\sharp}(s) \, \md \mu^{\sharp}(t) \! = \!
\iint_{\mathbb{R}^{2}} \ln \! \left(\vert s \! - \! t \vert^{2+\frac{1}{n}}
\vert st \vert^{-1}w^{o}(s)w^{o}(t) \right)^{-1} \md \mu^{\sharp}(s) \, \md
\mu^{\sharp}(t) \! \geqslant \! 0,
\end{equation*}
it follows that
\begin{align*}
\iint_{\mathbb{R}^{2}} \ln \! \left(\dfrac{\vert st \vert}{\vert s \! - \! t
\vert^{2+\frac{1}{n}}} \right) \! \left(\md \widetilde{\mu}_{V}^{o}(s) \, \md
\widetilde{\mu}_{V}^{o}(t) \! + \! \md \mu_{V}^{o}(s) \, \md \mu_{V}^{o}(t)
\right) \geqslant& \, \iint_{\mathbb{R}^{2}} \ln \! \left(\dfrac{\vert st
\vert}{\vert s \! - \! t \vert^{2+\frac{1}{n}}} \right) \! \left(\md
\widetilde{\mu}_{V}^{o}(s) \, \md \mu_{V}^{o}(t) \right. \\
+& \left. \, \md \mu_{V}^{o}(s) \md \widetilde{\mu}_{V}^{o}(t) \right),
\end{align*}
or, via a straightforward symmetry argument,
\begin{align*}
\iint_{\mathbb{R}^{2}} \ln \! \left(\dfrac{\vert st \vert}{\vert s \! - \! t
\vert^{2+\frac{1}{n}}} \right) \! \left(\md \widetilde{\mu}_{V}^{o}(s) \, \md
\widetilde{\mu}_{V}^{o}(t) \! + \! \md \mu_{V}^{o}(s) \, \md \mu_{V}^{o}(t)
\right) \geqslant& \, 2 \iint_{\mathbb{R}^{2}} \ln \! \left(\dfrac{\vert st
\vert}{\vert s \! - \! t \vert^{2+\frac{1}{n}}} \right) \! \md \widetilde{
\mu}_{V}^{o}(s) \, \md \mu_{V}^{o}(t) \\
=& \, 2 \iint_{\mathbb{R}^{2}} \ln \! \left(\dfrac{\vert st \vert}{\vert s \!
- \! t \vert^{2+\frac{1}{n}}} \right) \! \md \mu_{V}^{o}(s) \, \md \widetilde{
\mu}_{V}^{o}(t).
\end{align*}
The above shows that (since both $\mathrm{I}_{V}^{o}[\mu_{V}^{o}]$ and
$\mathrm{I}_{V}^{o}[\widetilde{\mu}_{V}^{o}] \! < \! +\infty)$ $\ln (\vert st
\vert \vert s \! - \! t \vert^{-(2+\frac{1}{n})})$ is integrable with respect
to both $\md \widetilde{\mu}_{V}^{o}(s) \, \md \mu_{V}^{o}(t)$ and $\md \mu_{
V}^{o}(s) \, \md \widetilde{\mu}_{V}^{o}(t)$. {}From an argument on pg.~149 of
\cite{a79}, it follows that $\ln (\vert st \vert \vert s \! - \! t \vert^{-(2
+\frac{1}{n})})$ is integrable with respect to (the measure) $\md \mu_{t}^{o}
(s) \, \md \mu_{t}^{o}(t^{\prime})$, where $\mu_{t}^{o}(z) \! := \! \mu_{V}^{
o}(z) \! + \! t(\widetilde{\mu}_{V}^{o}(z) \! - \! \mu_{V}^{o}(z))$, $(z,t)
\! \in \! \mathbb{R} \times [0,1]$. Set
\begin{equation*}
\mathscr{F}_{\mu}(t) \! := \! \iint_{\mathbb{R}^{2}} \ln \! \left(\vert st^{
\prime} \vert \vert s \! - \! t^{\prime} \vert^{-(2+\frac{1}{n})}(w^{o}(s)
w^{o}(t^{\prime}))^{-1} \right) \md \mu_{t}^{o}(s) \, \md \mu_{t}^{o}(t^{
\prime})
\end{equation*}
$(= \! \mathrm{I}_{V}^{o}[\mu_{t}^{o}])$. Noting that
\begin{align*}
\md \mu_{t}^{o}(s) \, \md \mu_{t}^{o}(t^{\prime}) =& \, \md \mu_{V}^{o}(s) \,
\md \mu_{V}^{o}(t^{\prime}) \! + \! t \md \mu_{V}^{o}(s)(\md \widetilde{\mu}_{
V}^{o}(t^{\prime}) \! - \! \md \mu_{V}^{o}(t^{\prime})) \! + \! t \md \mu_{
V}^{o}(t^{\prime})(\md \widetilde{\mu}_{V}^{o}(s) \! - \! \md \mu_{V}^{o}(s))
\\
+& \, t^{2}(\md \widetilde{\mu}_{V}^{o}(s) \! - \! \md \mu_{V}^{o}(s))(\md
\widetilde{\mu}_{V}^{o}(t^{\prime}) \! - \! \md \mu_{V}^{o}(t^{\prime})),
\end{align*}
it follows that
\begin{align*}
\mathscr{F}_{\mu}(t) =& \, \mathrm{I}_{V}^{o}[\mu_{V}^{o}] \! + \! 2t \iint_{
\mathbb{R}^{2}} \ln \! \left(\dfrac{\vert st^{\prime} \vert}{\vert s \! - \!
t^{\prime} \vert^{2+\frac{1}{n}}}(w^{o}(s)w^{o}(t^{\prime}))^{-1} \right) \!
\md \mu_{V}^{o}(s) \, (\md \widetilde{\mu}_{V}^{o}(t^{\prime}) \! - \! \md
\mu_{V}^{o}(t^{\prime})) \\
+& \, t^{2} \iint_{\mathbb{R}^{2}} \ln \! \left(\dfrac{\vert st^{\prime}
\vert}{\vert s \! - \! t^{\prime} \vert^{2+\frac{1}{n}}}(w^{o}(s)w^{o}(t^{
\prime}))^{-1} \right) \! (\md \widetilde{\mu}_{V}^{o}(s) \! - \! \md \mu_{
V}^{o}(s)) \, (\md \widetilde{\mu}_{V}^{o}(t^{\prime}) \! - \! \md \mu_{V}^{o}
(t^{\prime})).
\end{align*}
Since $\mu^{\sharp} \! \in \! \mathcal{M}_{1}(\mathbb{R})$ is finite-moment
signed measure with mean zero, that is, $\int_{\mathbb{R}} \md \mu^{\sharp}
(\xi) \! = \! \int_{\mathbb{R}} \md (\widetilde{\mu}_{V}^{o} \! - \! \mu_{
V}^{o})(\xi) \! = \! 0$, and compact support, it follows {}from the analysis
above and the result of Lemma~3.2 that $\mathscr{F}_{\mu}(t)$ is
convex\footnote{If $f$ is twice differentiable on $(a,b)$, then $f^{\prime
\prime}(x) \! \geqslant \! 0$ on $(a,b)$ is both a necessary and sufficient
condition that $f$ be convex on $(a,b)$.}; thus, for $t \! \in \! [0,1]$,
\begin{align*}
\mathrm{I}_{V}^{o}[\mu_{V}^{o}] \leqslant& \, \mathscr{F}_{\mu}(t) \! = \!
\mathrm{I}_{V}^{o}[\mu_{t}^{o}] \! = \! \mathscr{F}_{\mu}(t \! + \! (1 \! - \!
t)0) \! \leqslant \! t \mathscr{F}_{\mu}(1) \! + \! (1 \! - \! t) \mathscr{
F}_{\mu}(0) \\
=& \, t \mathrm{I}_{V}^{o}[\widetilde{\mu}_{V}^{o}] \! + \! (1 \! - \! t)
\mathrm{I}_{V}^{o}[\mu_{V}^{o}] \! = \! t \mathrm{I}_{V}^{o}[\mu_{V}^{o}] \!
+ \! (1 \! - \! t) \mathrm{I}_{V}^{o}[\mu_{V}^{o}] \Rightarrow \\
\mathrm{I}_{V}^{o}[\mu_{V}^{o}] \leqslant& \, \mathrm{I}_{V}^{o}[\mu_{t}^{o}]
\! \leqslant \! \mathrm{I}_{V}^{o}[\mu_{V}^{o}],
\end{align*}
whence $\mathrm{I}_{V}^{o}[\mu_{t}^{o}] \! = \! \mathrm{I}_{V}^{o}[\mu_{V}^{
o}] \! := \! E_{V}^{o}$ $(= \! \text{const.})$. Since $\mathrm{I}_{V}^{o}
[\mu_{t}^{o}] \! = \! \mathscr{F}_{\mu}(t) \! = \! E_{V}^{o}$, it follows, in
particular, that $\mathscr{F}_{\mu}^{\prime \prime}(0) \! = \! 0 \Rightarrow$
\begin{align*}
0 =& \, \iint_{\mathbb{R}^{2}} \ln \! \left(\dfrac{\vert st^{\prime} \vert}{
\vert s \! - \! t^{\prime} \vert^{2+\frac{1}{n}}}(w^{o}(s)w^{o}(t^{\prime}))^{
-1} \right) \! (\md \widetilde{\mu}_{V}^{o}(s) \! - \! \md \mu_{V}^{o}(s)) \,
(\md \widetilde{\mu}_{V}^{o}(t^{\prime}) \! - \! \md \mu_{V}^{o}(t^{\prime}))
\\
=& \, \iint_{\mathbb{R}^{2}} \ln \! \left(\dfrac{\vert st^{\prime} \vert}{
\vert s \! - \! t^{\prime} \vert^{2+\frac{1}{n}}} \right) \! (\md \widetilde{
\mu}_{V}^{o}(s) \! - \! \md \mu_{V}^{o}(s)) \, (\md \widetilde{\mu}_{V}^{o}
(t^{\prime}) \! - \! \md \mu_{V}^{o}(t^{\prime})) \\
+& \, 2 \int_{\mathbb{R}} \widetilde{V}(t^{\prime}) \, \md (\widetilde{\mu}_{
V}^{o} \! - \! \mu_{V}^{o})(t^{\prime}) \, \, \underbrace{\int_{\mathbb{R}}
\md (\widetilde{\mu}_{V}^{o} \! - \! \mu_{V}^{o})(s)}_{= \, 0} \Rightarrow \\
0 =& \, \iint_{\mathbb{R}^{2}} \ln \! \left(\dfrac{\vert st^{\prime} \vert}{
\vert s \! - \! t^{\prime} \vert^{2+\frac{1}{n}}} \right) \! \md (\widetilde{
\mu}_{V}^{o} \! - \! \mu_{V}^{o})(s) \, \md (\widetilde{\mu}_{V}^{o} \! - \!
\mu_{V}^{o})(t^{\prime});
\end{align*}
but, in Lemma~3.2, it was shown that,
\begin{equation*}
\iint_{\mathbb{R}^{2}} \ln \! \left(\dfrac{\vert st^{\prime} \vert}{\vert s
\! - \! t^{\prime} \vert^{2+\frac{1}{n}}} \right) \! \md (\widetilde{\mu}_{
V}^{o} \! - \! \mu_{V}^{o})(s) \, \md (\widetilde{\mu}_{V}^{o} \! - \! \mu_{
V}^{o})(t^{\prime}) \! = \! \left(2 \! + \! \dfrac{1}{n} \right) \! \int_{
0}^{+\infty} \xi^{-1} \vert (\widehat{\widetilde{\mu}_{V}^{o}} \! - \!
\widehat{\mu_{V}^{o}})(\xi) \vert^{2} \, \md \xi \quad (\geqslant \! 0),
\end{equation*}
whence $\int_{0}^{+\infty} \xi^{-1} \vert (\widehat{\widetilde{\mu}_{V}^{o}}
\! - \! \widehat{\mu_{V}^{o}})(\xi) \vert^{2} \, \md \xi \! = \! 0$
$\Rightarrow$ $\widehat{\widetilde{\mu}_{V}^{o}}(\xi) \! = \! \widehat{\mu_{
V}^{o}}(\xi)$, $\xi \! \geqslant \! 0$. Noting that
\begin{equation*}
\widehat{\widetilde{\mu}_{V}^{o}}(-\xi) \! = \! \int_{\mathbb{R}} \me^{\mi
s(-\xi)} \, \md \widetilde{\mu}_{V}^{o}(s) \! = \! \overline{\widehat{
\widetilde{\mu}_{V}^{o}}(\xi)} \qquad \text{and} \qquad \widehat{\mu_{V}^{o}}
(-\xi) \! = \! \int_{\mathbb{R}} \me^{\mi s(-\xi)} \, \md \mu_{V}^{o}(s) \! =
\! \overline{\widehat{\mu_{V}^{o}}(\xi)},
\end{equation*}
it follows {}from $\widehat{\widetilde{\mu}_{V}^{o}}(\xi) \! = \! \widehat{
\mu_{V}^{o}}(\xi)$, $\xi \! \geqslant \! 0$, via a complex-conjugation
argument, that $\widehat{\widetilde{\mu}_{V}^{o}}(-\xi) \! = \! \widehat{\mu_{
V}^{o}}(-\xi)$, $\xi \! \geqslant \! 0$; hence, $\widehat{\widetilde{\mu}_{
V}^{o}}(\xi) \! = \! \widehat{\mu_{V}^{o}}(\xi)$, $\xi \! \in \! \mathbb{R}$.
The latter relation shows that $\int_{\mathbb{R}} \me^{\mi s \xi} \, \md
(\widetilde{\mu}_{V}^{o} \! - \! \mu_{V}^{o})(s) \! = \! 0$ $\Rightarrow$
$\widetilde{\mu}_{V}^{o} \! = \! \mu_{V}^{o}$; thus the uniqueness of the
`odd' equilibrium measure. \hfill $\qed$

Before proceeding to Lemma~3.4, the following observations, which are
interesting, non-trivial and important results in their own right, should be
noted. Let $\widetilde{V} \colon \mathbb{R} \setminus \{0\} \! \to \! \mathbb{
R}$ satisfy conditions~(2.3)--(2.5). For each $m \! \in \! \mathbb{Z}_{0}^{+}$
and any $(2m+1)$-tuple $(x_{1},x_{2},\dotsc,x_{2m+1})$ of distinct, finite
and non-zero real numbers, let, for $n \! \in \! \mathbb{N}$,
\begin{equation*}
\mathfrak{d}_{o,m}^{\widetilde{V}}(n) \! := \! \dfrac{1}{2m(2m \! + \! 1)} \,
\inf_{\{x_{1},x_{2},\dotsc,x_{2m+1}\} \subset \mathbb{R} \setminus \{0\}}
\left(\sum_{\substack{j,k=1\\j \neq k}}^{2m+1} \ln \! \left(\left\vert x_{j}
\! - \! x_{k} \right\vert^{1+\frac{1}{n}} \! \left\vert x_{k}^{-1} \! - \!
x_{j}^{-1} \right\vert \right)^{-1} \! + \! 4m \sum_{i=1}^{2m+1} \widetilde{V}
(x_{i}) \right).
\end{equation*}
For each $m \! \in \! \mathbb{Z}_{0}^{+}$, a set $\left\lbrace x_{1}^{\flat},
x_{2}^{\flat},\dotsc,x_{2m+1}^{\flat} \right\rbrace$ which realizes the above
infimum, that is, for $n \! \in \! \mathbb{N}$,
\begin{equation*}
\mathfrak{d}_{o,m}^{\widetilde{V}}(n) \! = \! \dfrac{1}{2m(2m \! + \! 1)} \!
\left(\sum_{\substack{j,k=1\\j \neq k}}^{2m+1} \ln \! \left(\left\vert x_{j}^{
\flat} \! - \! x_{k}^{\flat} \right\vert^{1+\frac{1}{n}} \! \left\vert (x_{
k}^{\flat})^{-1} \! - \! (x_{j}^{\flat})^{-1} \right\vert \right)^{-1} \! + \!
4m \sum_{i=1}^{2m+1} \widetilde{V}(x_{i}^{\flat}) \right),
\end{equation*}
will be called (with slight abuse of nomenclature) a \emph{generalised
weighted $(2m \! + \! 1)$-Fekete set}, and the points $x_{1}^{\flat},x_{2}^{
\flat},\dotsc,x_{2m+1}^{\flat}$ will be called \emph{generalised weighted
Fekete points}. For $\left\lbrace x_{1}^{\flat},x_{2}^{\flat},\dotsc,x_{2m+
1}^{\flat} \right\rbrace$ a generalised weighted $(2m \! + \! 1)$-Fekete set,
denote by
\begin{equation*}
\mu_{\mathbf{x}^{\flat}}^{o} \! := \! \dfrac{1}{2m \! + \! 1} \sum_{j=1}^{2m
+1} \delta_{x_{j}^{\flat}},
\end{equation*}
where $\delta_{x_{j}^{\flat}}$, $j \! = \! 1,\dotsc,2m \! + \! 1$, is the 
Dirac delta measure (atomic mass) concentrated at $x_{j}^{\flat}$, the 
\emph{normalised counting measure}, that is, $\int_{\mathbb{R}} \md \mu_{
\mathbf{x}^{\flat}}^{o}(s) \! = \! 1$. Then, mimicking the calculations in 
Chapter~6 of \cite{a79} and the techniques used to prove Theorem~1.34 in 
\cite{a44} (see, in particular, Section~2 of \cite{a44}), one proves that, 
for $n \! \in \! \mathbb{N}$ (the details are left to the interested reader):
\begin{compactenum}
\item[\textbullet] $\lim_{m \to \infty} \mathfrak{d}_{o,m}^{\widetilde{V}}(n)$
exists, more precisely,
\begin{equation*}
\lim_{m \to \infty} \mathfrak{d}_{o,m}^{\widetilde{V}}(n) \! = \! E_{V}^{o} \! 
= \! \inf \! \left\lbrace \mathstrut \mathrm{I}_{V}^{o}[\mu^{o}]; \, \mu^{o} 
\! \in \! \mathcal{M}_{1}(\mathbb{R}) \right\rbrace,
\end{equation*}
where (the functional) $\mathrm{I}_{V}^{o}[\mu^{o}] \colon \mathcal{M}_{1}
(\mathbb{R}) \! \to \! \mathbb{R}$ is defined in Lemma~3.1, and $\lim_{m 
\to \infty} \exp (-\mathfrak{d}_{o,m}^{\widetilde{V}}(n) \! ) \linebreak[4]
= \! \exp (-E_{V}^{o})$ is positive and finite;
\item[\textbullet] $\mu_{\mathbf{x}^{\flat}}^{o}$ converges weakly (in the 
weak-$\ast$ topology of measures) to the `odd' equilibrium measure $\mu_{V}^{
o}$, that is, $\mu_{\mathbf{x}^{\flat}}^{o} \overset{\ast}{\to} \mu_{V}^{o}$ 
as $m \! \to \! \infty$.
\end{compactenum}

\textbf{RHP2}, that is, $(\overset{o}{\operatorname{Y}}(z),\mathrm{I} \! + \!
\exp (-n \widetilde{V}(z)) \sigma_{+},\mathbb{R})$, is now reformulated as an
equivalent, auxiliary RHP normalised at zero.
\begin{notrem}
For completeness, the integrand appearing in the definition of $g^{o}(z)$ (see
Lemma~3.4 below) is defined as follows: $\ln ((z \! - \! s)^{2+\frac{1}{n}}
(zs)^{-1}) \! := \! (2 \! + \! \tfrac{1}{n}) \ln (z \! - \! s) \! - \! \ln z
\! - \! \ln s$, where, for $s \! < \! 0$, $\ln s \! := \! \ln \lvert s \rvert
\! + \! \mi \pi$. \hfill $\blacksquare$
\end{notrem}
\begin{ccccc}
Let the external field $\widetilde{V} \colon \mathbb{R} \setminus \{0\} \! \to
\! \mathbb{R}$ satisfy conditions~{\rm (2.3)--(2.5)}. For the associated `odd'
equilibrium measure, $\mu_{V}^{o} \! \in \! \mathcal{M}_{1}(\mathbb{R})$, set
$J_{o} \! := \! \operatorname{supp}(\mu_{V}^{o})$, where $J_{o}$ $(=$ compact)
$\subset \! \overline{\mathbb{R}} \setminus \lbrace 0,\pm \infty \rbrace$,
and let $\overset{o}{\mathrm{Y}} \colon \mathbb{C} \setminus \mathbb{R} \! \to
\! \mathrm{SL}_{2}(\mathbb{C})$ be the (unique) solution of {\rm \pmb{RHP2}}.
Let
\begin{equation*}
\overset{o}{\mathscr{M}}(z) \! := \! \me^{-\frac{n \ell_{o}}{2} \mathrm{ad}
(\sigma_{3})} \overset{o}{\mathrm{Y}}(z) \me^{-n(g^{o}(z)-\mathfrak{Q}_{
\mathscr{A}}) \sigma_{3}},
\end{equation*}
where $g^{o}(z)$, the `odd' $g$-function, is defined by, for $n \! \in \! 
\mathbb{N}$,
\begin{equation*}
g^{o}(z) \! := \! \int_{J_{o}} \ln \! \left((z \! - \! s)^{2+\frac{1}{n}}
(zs)^{-1} \right) \md \mu_{V}^{o}(s), \quad z \! \in \! \mathbb{C} \setminus
(-\infty,\max \lbrace 0,\max \lbrace \operatorname{supp}(\mu_{V}^{o}) \rbrace
\rbrace),
\end{equation*}
$\ell_{o}$ $(\in \! \mathbb{R})$, the `odd' variational constant, is given in
{\rm Lemma~3.6} below, and
\begin{equation*}
\mathfrak{Q}_{\mathscr{A}} \! = \!
\begin{cases}
\mathfrak{Q}_{\mathscr{A}}^{+} \! := \! (1 \! + \! \frac{1}{n}) \int_{J_{o}}
\ln (\lvert s \rvert) \, \md \mu_{V}^{o}(s) \! - \! \mi \pi \int_{J_{o} \cap
\mathbb{R}_{-}} \md \mu_{V}^{o}(s) \! + \! \mi \pi (2 \! + \! \frac{1}{n})
\int_{J_{o} \cap \mathbb{R}_{+}} \md \mu_{V}^{o}(s), &\text{$z \! \in \!
\mathbb{C}_{+}$,} \\
\mathfrak{Q}_{\mathscr{A}}^{-} \! := \! (1 \! + \! \frac{1}{n}) \int_{J_{o}}
\ln (\lvert s \rvert) \, \md \mu_{V}^{o}(s) \! - \! \mi \pi \int_{J_{o} \cap
\mathbb{R}_{-}} \md \mu_{V}^{o}(s) \! - \! \mi \pi (2 \! + \! \frac{1}{n})
\int_{J_{o} \cap \mathbb{R}_{+}} \md \mu_{V}^{o}(s), &\text{$z \! \in \!
\mathbb{C}_{-}$,}
\end{cases}
\end{equation*}
with (see Lemma~{\rm 3.5}, item~{\rm (1)}, below)
\begin{equation*}
\int_{J_{o} \cap \mathbb{R}_{-}} \md \mu_{V}^{o}(s) \! = \!
\begin{cases}
0, &\text{$J_{o} \! \subset \! \mathbb{R}_{+}$,} \\
1, &\text{$J_{o} \! \subset \! \mathbb{R}_{-}$,} \\
\int_{b_{0}^{o}}^{a_{j}^{o}} \md \mu_{V}^{o}(s), &\text{$(a_{j}^{o},b_{j}^{o})
\! \ni \! 0, \quad j \! = \! 1,\dotsc,N$,}
\end{cases}
\end{equation*}
and
\begin{equation*}
\int_{J_{o} \cap \mathbb{R}_{+}} \md \mu_{V}^{o}(s) \! = \!
\begin{cases}
0, &\text{$J_{o} \! \subset \! \mathbb{R}_{-}$,} \\
1, &\text{$J_{o} \! \subset \! \mathbb{R}_{+}$,} \\
\int_{b_{j}^{o}}^{a_{N+1}^{o}} \md \mu_{V}^{o}(s), &\text{$(a_{j}^{o},b_{j}^{o}
) \! \ni \! 0, \quad j \! = \! 1,\dotsc,N$.}
\end{cases}
\end{equation*}
Then $\overset{o}{\mathscr{M}} \colon \mathbb{C} \setminus \mathbb{R} \! \to 
\! \mathrm{SL}_{2}(\mathbb{C})$ solves the following (normalised at zero) 
{\rm RHP:} {\rm (i)} $\overset{o}{\mathscr{M}}(z)$ is holomorphic for $z \! 
\in \! \mathbb{C} \setminus \mathbb{R};$ {\rm (ii)} the boundary values 
$\overset{o}{\mathscr{M}}_{\pm}(z) := \! \lim_{\underset{\pm \Im (z^{\prime})
>0}{z^{\prime} \to z}} \overset{o}{\mathscr{M}}(z^{\prime})$ satisfy the jump 
condition
\begin{equation*}
\overset{o}{\mathscr{M}}_{+}(z) \! = \! \overset{o}{\mathscr{M}}_{-}(z) \!
\begin{pmatrix}
\me^{-n(g^{o}_{+}(z)-g^{o}_{-}(z)-\mathfrak{Q}^{+}_{\mathscr{A}}+\mathfrak{Q}^{
-}_{\mathscr{A}})} & \me^{n(g^{o}_{+}(z)+g^{o}_{-}(z)-\widetilde{V}(z)- \ell_{
o}-\mathfrak{Q}^{+}_{\mathscr{A}}-\mathfrak{Q}^{-}_{\mathscr{A}})} \\
0 & \me^{n(g^{o}_{+}(z)-g^{o}_{-}(z)-\mathfrak{Q}^{+}_{\mathscr{A}}+\mathfrak{
Q}^{-}_{\mathscr{A}})}
\end{pmatrix}, \quad z \! \in \! \mathbb{R},
\end{equation*}
with $g^{o}_{\pm}(z) \! := \! \lim_{\varepsilon \downarrow 0}g^{o}(z \! \pm
\! \mi \varepsilon);$ {\rm (iii)} $\overset{o}{\mathscr{M}}(z) \! =_{\underset{
z \in \mathbb{C} \setminus \mathbb{R}}{z \to 0}} \! \mathrm{I} \! + \!
\mathcal{O}(z);$ and {\rm (iv)} $\overset{o}{\mathscr{M}}(z) \! =_{\underset{z
\in \mathbb{C} \setminus \mathbb{R}}{z \to \infty}} \! \mathcal{O}(1)$.
\end{ccccc}

\emph{Proof.} For (arbitrary) $z_{1},z_{2} \! \in \! \mathbb{C}_{\pm}$, note
that, {}from the definition of $g^{o}(z)$ stated in the Lemma, $g^{o}(z_{2})
\! - \! g^{o}(z_{1}) \! = \! \mi \pi \int_{z_{1}}^{z_{2}} \mathscr{F}^{o}(s)
\, \md s$, where
\begin{equation*}
\mathscr{F}^{o} \colon \mathbb{C} \setminus (\operatorname{supp}(\mu_{V}^{o})
\cup \{0\}) \! \to \! \mathbb{C}, \, \, z \! \mapsto \! -\dfrac{1}{\mi \pi} \!
\left(\dfrac{1}{z} \! + \! \left(2 \! + \! \dfrac{1}{n} \right) \! \int_{J_{
o}} \dfrac{\md \mu_{V}^{o}(s)}{s \! - \! z} \right),
\end{equation*}
with $\mathscr{F}^{o}(z) \! =_{z \to 0} \! -\tfrac{1}{\pi \mi z} \! + \!
\mathcal{O}(1)$ (since $\mu_{V}^{o} \! \in \! \mathcal{M}_{1}(\mathbb{R})$; in
particular, $\int_{\mathbb{R}}s^{m} \, \md \mu_{V}^{o}(s) \! < \! \infty$, $m
\! \in \! \mathbb{Z})$; thus, $\vert g^{o}(z_{2}) \! - \! g^{o}(z_{1}) \vert
\! \leqslant \! \pi \sup_{z \in \mathbb{C}_{\pm}} \vert \mathscr{F}^{o}(z)
\vert \vert z_{2} \! - \! z_{1} \vert$, that is, $g^{o}(z)$ is uniformly
Lipschitz continuous in $\mathbb{C}_{\pm}$. Thus, {}from the definition of
$g^{o}(z)$ stated in the Lemma:
\begin{compactenum}
\item[(1)] for $s \! \in \! J_{o}$, $z \! \in \! \mathbb{C} \setminus (-
\infty,\max \lbrace 0,\max \lbrace \operatorname{supp}(\mu_{V}^{o}) \rbrace
\rbrace)$, with $\vert s/z \vert \! \ll \! 1$, and $\mu_{V}^{o} \! \in \!
\mathcal{M}_{1}(\mathbb{R})$, in particular, $\int_{\mathbb{R}} \md \mu_{V}^{
o}(s)$ $(= \! \int_{J_{o}} \md \mu_{V}^{o}(s))$ $= \! 1$ and $\int_{\mathbb{
R}} s^{m} \, \md \mu_{V}^{o}(s)$ $(= \! \int_{J_{o}} s^{m} \, \md \mu_{V}^{o}
(s))$ $< \! \infty$, $m \! \in \! \mathbb{N}$, it follows {}from the
expansions $\tfrac{1}{s-z} \! = \! -\sum_{k=0}^{l} \tfrac{s^{k}}{z^{k+1}} \!
+ \! \tfrac{s^{l+1}}{z^{l+1}(s-z)}$, $l \! \in \! \mathbb{Z}_{0}^{+}$, and
$\ln (z \! - \! s) \! =_{\vert z \vert \to \infty} \! \ln (z) \! - \! \sum_{k
=1}^{\infty} \tfrac{1}{k}(\tfrac{s}{z})^{k}$, that
\begin{equation*}
g^{o}(z) \underset{\mathbb{C} \setminus \mathbb{R} \ni z \to \infty}{=} \left(
1 \! + \! \dfrac{1}{n} \right) \! \ln (z) \! - \! \int_{J_{o}} \ln (\lvert s 
\rvert) \, \md \mu_{V}^{o}(s) \! - \! \mi \pi \int_{J_{o} \cap \mathbb{R}_{-}} 
\md \mu_{V}^{o}(s) \! + \! \mathcal{O}(z^{-1}),
\end{equation*}
where $\int_{J_{o} \cap \mathbb{R}_{-}} \md \mu_{V}^{o}(s)$ is given in the 
Lemma;
\item[(2)] for $s \! \in \! J_{o}$, $z \! \in \! \mathbb{C} \setminus (-
\infty,\max \lbrace 0,\max \lbrace \operatorname{supp}(\mu_{V}^{o}) \rbrace 
\rbrace)$, with $\vert z/s \vert \! \ll \! 1$, and $\mu_{V}^{o} \! \in \! 
\mathcal{M}_{1}(\mathbb{R})$, in particular, $\int_{\mathbb{R}}s^{-m} \, 
\md \mu_{V}^{o}(s)$ $(= \! \int_{J_{o}}s^{-m} \, \md \mu_{V}^{o}(s))$ $< \! 
\infty$, $m \! \in \! \mathbb{N}$, it follows {}from the expansions $\tfrac{
1}{z-s} \! = \! -\sum_{k=0}^{l} \tfrac{z^{k}}{s^{k+1}} \! + \! \tfrac{z^{l+1}
}{s^{l+1}(z-s)}$, $l \! \in \! \mathbb{Z}_{0}^{+}$, and $\ln (s \! - \! z) \! 
=_{\vert z \vert \to 0} \! \ln (s) \! - \! \sum_{k=1}^{\infty} \tfrac{1}{k}
(\tfrac{z}{s})^{k}$, that
\begin{align*}
g^{o}(z) \underset{\mathbb{C}_{\pm} \ni z \to 0}{=}& \, -\ln (z) \! + \! 
\left(1 \! + \! \dfrac{1}{n} \right) \int_{J_{o}} \ln (\lvert s \rvert) \, 
\md \mu_{V}^{o}(s) \! - \! \mi \pi \int_{J_{o} \cap \mathbb{R}_{-}} \md 
\mu_{V}^{o}(s) \\
\pm& \, \mi \pi \! \left(2 \! + \! \dfrac{1}{n} \right) \! \int_{J_{o} \cap
\mathbb{R}_{+}} \md \mu_{V}^{o}(s) \! + \! \mathcal{O}(z),
\end{align*}
where $\int_{J_{o} \cap \mathbb{R}_{+}} \md \mu_{V}^{o}(s)$ is given in the
Lemma.
\end{compactenum}
Items~(i)--(iv) now follow {}from the definitions of $\overset{o}{\mathscr{M}}
(z)$ (in terms of $\overset{o}{\mathrm{Y}}(z))$ and $g^{o}(z)$ stated in the
Lemma, and the above two asymptotic expansions. \hfill $\qed$
\begin{ccccc}
Let the external field $\widetilde{V} \colon \mathbb{R} \setminus \{0\} \! \to
\! \mathbb{R}$ satisfy conditions~{\rm (2.3)--(2.5)}. For $\mu_{V}^{o} \! \in
\! \mathcal{M}_{1}(\mathbb{R})$, the associated `odd' equilibrium measure, set
$J_{o} \! := \! \operatorname{supp}(\mu_{V}^{o})$, where $J_{o}$ $(=$ compact)
$\subset \! \overline{\mathbb{R}} \setminus \lbrace 0,\pm \infty \rbrace$. 
Then: {\rm (1)} $J_{o} \! = \! \cup_{j=1}^{N+1}(b_{j-1}^{o},a_{j}^{o})$, 
with $N \! \in \! \mathbb{N}$ and finite, $b_{0}^{o} \! := \! \min \lbrace 
\operatorname{supp}(\mu_{V}^{o}) \rbrace \! \notin \lbrace -\infty,0 
\rbrace$, $a_{N+1}^{o} \! := \! \max \lbrace \operatorname{supp}(\mu_{V}^{o}) 
\rbrace \! \notin \! \lbrace 0,+\infty \rbrace$, and $-\infty \! < \! b_{0}^{
o} \! < \! a_{1}^{o} \! < \! b_{1}^{o} \! < \! a_{2}^{o} \! < \! \cdots \! < 
\! b_{N}^{o} \! < \! a_{N+1}^{o} \! < \! +\infty$, and $\lbrace b_{j-1}^{o},
a_{j}^{o} \rbrace_{j=1}^{N+1}$ satisfy the $n$-dependent and (locally) 
solvable system of $2(N \! + \! 1)$ moment conditions
\begin{gather*}
\int_{J_{o}} \dfrac{(\frac{2 \mi}{\pi s} \! + \! \frac{\mi \widetilde{V}^{
\prime}(s)}{\pi})s^{j}}{(R_{o}(s))^{1/2}_{+}} \, \dfrac{\md s}{2 \pi \mi} \! =
\! 0, \quad j \! = \! 0,\dotsc,N, \qquad \qquad \int_{J_{o}} \dfrac{(\frac{2
\mi}{\pi s} \! + \! \frac{\mi \widetilde{V}^{\prime}(s)}{\pi})s^{N+1}}{(R_{o}
(s))^{1/2}_{+}} \, \dfrac{\md s}{2 \pi \mi} \! = \! \dfrac{1}{\mi \pi} \!
\left(2 \! + \! \dfrac{1}{n} \right), \\
\int_{a_{j}^{o}}^{b_{j}^{o}} \! \left(\mi (R_{o}(s))^{1/2} \int_{J_{o}}
\dfrac{(\frac{\mi}{\pi \xi} \! + \! \frac{\mi \widetilde{V}^{\prime}(\xi)}{2
\pi})}{(R_{o}(\xi))^{1/2}_{+}(\xi \! - \! s)} \, \dfrac{\md \xi}{2 \pi \mi}
\right) \! \md s \! = \! \dfrac{1}{2 \pi} \ln \! \left\vert \dfrac{a_{j}^{o}}{
b_{j}^{o}} \right\vert \! + \! \dfrac{1}{4 \pi} \! \left(\widetilde{V}(a_{j}^{
o}) \! - \! \widetilde{V}(b_{j}^{o}) \right), \quad j \! = \! 1,\dotsc,N,
\end{gather*}
where $(R_{o}(z))^{1/2}$ is defined in Theorem~{\rm 2.3.1},
Equation~{\rm (2.8)}, with $(R_{o}(z))^{1/2}_{\pm} \! := \! \lim_{\varepsilon
\downarrow 0}(R_{o}(z \! \pm \! \mi \varepsilon))^{1/2}$, and the branch of
the square root chosen so that $z^{-(N+1)}(R_{o}(z))^{1/2} \! \sim_{\underset{
z \in \mathbb{C}_{\pm}}{z \to \infty}} \! \pm 1;$ and {\rm (2)} the density of
the `odd' equilibrium measure, which is absolutely continuous with respect to
Lebesgue measure, is given by
\begin{equation*}
\md \mu_{V}^{o}(x) \! := \! \psi_{V}^{o}(x) \, \md x \! = \! \dfrac{1}{2 \pi
\mi}(R_{o}(x))^{1/2}_{+}h_{V}^{o}(x) \pmb{1}_{J_{o}}(x) \, \md x,
\end{equation*}
where
\begin{equation*}
h_{V}^{o}(z) \! := \! \dfrac{1}{2} \! \left(2 \! + \! \dfrac{1}{n} \right)^{-
1} \oint_{C_{\mathrm{R}}^{o}} \dfrac{(\frac{2 \mi}{\pi s} \! + \! \frac{\mi
\widetilde{V}^{\prime}(s)}{\pi})}{(R_{o}(s))^{1/2}(s \! - \! z)} \, \md s
\end{equation*}
(real analytic for $z \! \in \! \mathbb{R} \setminus \{0\})$, with $C_{
\mathrm{R}}^{o}$ $(\subset \mathbb{C}^{\ast})$ the boundary of any open
doubly-connected annular region of the type $\lbrace \mathstrut z^{\prime} \!
\in \! \mathbb{C}; \, 0 \! < \! r \! < \! \vert z^{\prime} \vert \! < \! R
\! < \! +\infty \rbrace$, where the simple outer (resp., inner) boundary
$\lbrace \mathstrut z^{\prime} \! = \! R \me^{\mi \vartheta}, \, 0 \!
\leqslant \! \vartheta \! \leqslant \! 2 \pi \rbrace$ (resp., $\lbrace
\mathstrut z^{\prime} \! = \! r \me^{\mi \vartheta}, \, 0 \! \leqslant \!
\vartheta \! \leqslant \! 2 \pi \rbrace)$ is traversed clockwise (resp.,
counter-clockwise), with the numbers $0 \! < \! r \! < \! R \! < \! +\infty$
chosen such that, for (any) non-real $z$ in the domain of analyticity of
$\widetilde{V}$ (that is, $\mathbb{C}^{\ast})$, $\mathrm{int}(C_{\mathrm{R}
}^{o}) \! \supset \! J_{o} \cup \{z\}$, $\pmb{1}_{J_{o}}(x)$ is the indicator
(characteristic) function of the set $J_{o}$, and $\psi_{V}^{o}(x) \!
\geqslant \! 0$ (resp., $\psi_{V}^{o}(x) \! > \! 0)$ $\forall \, \, x \! \in
\! \overline{J_{o}}:= \! \cup_{j=1}^{N+1}[b_{j-1}^{o},a_{j}^{o}]$ (resp.,
$\forall \, \, x \! \in \! J_{o})$.
\end{ccccc}

\emph{Proof.} One begins by showing that the support of the `odd' equilibrium 
measure, $\operatorname{supp}(\mu_{V}^{o}) \! =: \! J_{o}$, consists of the 
union of a finite number of disjoint and bounded (real) intervals. Recall 
{}from Lemma~3.1 that $J_{o} \! = \text{compact} \subset \overline{\mathbb{R}} 
\setminus \lbrace 0,\pm \infty \rbrace$, and that $\widetilde{V}$ is real 
analytic on $\mathbb{R} \setminus \{0\}$, thus real analytic on $J_{o}$, with 
an analytic continuation to the following (open) neighbourhood of $J_{o}$, 
$\mathbb{U} \! := \! \lbrace \mathstrut z \! \in \! \mathbb{C}; \, \inf_{q 
\in J_{o}} \vert z \! - \! q \vert \! < \! r \! \in \! (0,1) \rbrace \setminus 
\{0\}$. In analogy with Equation~(2.1) of \cite{a44}, for each $m \! \in \! 
\mathbb{Z}_{0}^{+}$ and any $2m \! + \! 1$-tuple $(x_{1},x_{2},\dotsc,x_{2m
+1})$ of distinct, finite and non-zero real numbers, let, for $n \! \in \! 
\mathbb{N}$,
\begin{align*}
\mathrm{d}_{\widetilde{V},m}^{o}(n) :=& \, \left(\sup_{\{x_{1},x_{2},\dotsc,
x_{2m+1}\} \subset \mathbb{R} \setminus \{0\}} \, \prod_{\substack{j,k=1\\j<
k}}^{2m+1} \left\vert x_{j} \! - \! x_{k} \right\vert^{2+\frac{2}{n}} \!
\left\vert x_{k}^{-1} \! - \! x_{j}^{-1} \right\vert^{2} \, \me^{-2
\widetilde{V}(x_{j})} \me^{-2 \widetilde{V}(x_{k})} \right)^{\frac{1}{2m
(2m+1)}} \\
=& \, \left(\sup_{\{x_{1},x_{2},\dotsc,x_{2m+1}\} \subset \mathbb{R} \setminus
\{0\}} \, \prod_{\substack{j,k=1\\j<k}}^{2m+1} \left\vert x_{j} \! - \! x_{k}
\right\vert^{2+\frac{2}{n}} \! \left\vert x_{k}^{-1} \! - \! x_{j}^{-1}
\right\vert^{2} \, \me^{-4m \sum_{i=1}^{2m+1} \widetilde{V}(x_{i})} \right)^{
\frac{1}{2m(2m+1)}},
\end{align*}
where $\prod_{\substack{j,k=1\\j<k}}^{2m+1}(\star) \! = \! \prod_{j=1}^{2m}
\prod_{k=j+1}^{2m+1}(\star)$. Denote by $\left\lbrace x^{\ast}_{1},x^{\ast}_{
2},\dotsc,x^{\ast}_{2m+1} \right\rbrace$, with $ x^{\ast}_{i} \! < \!
x^{\ast}_{j} \, \, \forall \, \, i \! < \! j \! \in \! \lbrace 1,\dotsc,2m \!
+ \! 1\rbrace$, the associated generalised weighted $(2m \! + \! 1)$-Fekete
set (see the discussion preceding Lemma~3.4), that is, for $n \! \in \!
\mathbb{N}$,
\begin{equation*}
\mathrm{d}_{\widetilde{V},m}^{o}(n) \! = \! \left(\prod_{\substack{j,k=1\\j<
k}}^{2m+1} \left\vert x_{j}^{\ast} \! - \! x_{k}^{\ast} \right\vert^{2+\frac{
2}{n}} \! \left\vert (x_{k}^{\ast})^{-1} \! - \! (x_{j}^{\ast})^{-1}
\right\vert^{2} \, \me^{-4m \sum_{i=1}^{2m+1} \widetilde{V}(x_{i}^{\ast})}
\right)^{\frac{1}{2m(2m+1)}}.
\end{equation*}
Proceeding, now, as in the proof of Theorem~1.34, Equation~(1.35), of 
\cite{a44}, in particular, mimicking the calculations on pp.~408--413 of 
\cite{a44} (for the proofs of Lemmae~2.3 and~2.15 therein), namely, using 
those techniques to show that, in the present case, the nearest-neighbour 
distances $\lbrace x_{j+1}^{\ast} \! - \! x_{j}^{\ast} \rbrace_{j=1}^{2m}$ 
are not `too small' as $m \! \to \! \infty$, and the calculations on 
pp.~413--415 of \cite{a44} (for the proof of Lemma~2.26 therein), one shows 
that, for the regular case considered herein (cf. Subsection~2.2), the `odd' 
equilibrium measure, $\mu_{V}^{o}$ $(\in \! \mathcal{M}_{1}(\mathbb{R}))$, 
is absolutely continuous with respect to Lebesgue measure, that is, the 
density of the `odd' equilibrium measure has the representation $\md \mu_{
V}^{o}(x) \! := \! \psi_{V}^{o}(x) \, \md x$, $x \! \in \! \operatorname{supp}
(\mu_{V}^{o})$, where $\psi_{V}^{o}(x) \! \geqslant \! 0$ on $\overline{J_{
o}}$, with $\psi_{V}^{o}(\pmb{\cdot})$ determined (explicitly) 
below\footnote{The analysis of \cite{a44} is, in some sense, more complicated 
than the one of the present paper, because, unlike the `real-line' case 
considered herein, that is, $\operatorname{supp}(\mu_{V}^{o}) \! =: \! J_{o} 
\subset \overline{\mathbb{R}} \setminus \lbrace 0,\pm \infty \rbrace$, the 
end-point effects at $\pm 1$ in \cite{a44} require special consideration (see, 
also, Section~4 of \cite{a44}).}.

Set
\begin{equation}
\mathscr{H}^{o}(z) \! := \! (\mathscr{F}^{o}(z))^{2} \! - \! \int_{J_{o}}
\dfrac{(4 \mi (2 \! + \! \frac{1}{n})^{2} \psi_{V}^{o}(\xi)(\mathcal{H} \psi_{
V}^{o})(\xi) \! - \! \frac{4 \mi}{\pi \xi}(2 \! + \! \frac{1}{n}) \psi_{V}^{o}
(\xi))}{(\xi \! - \! z)} \, \dfrac{\md \xi}{2 \pi \mi}, \quad z \! \in \!
\mathbb{C} \setminus (J_{o} \cup \{0\}),
\end{equation}
where, {}from the proof of Lemma~3.4,
\begin{equation}
\mathscr{F}^{o}(z) \! = \! -\dfrac{1}{\mi \pi} \! \left(\dfrac{1}{z} \! + \!
\left(2 \! + \! \dfrac{1}{n} \right) \! \int_{J_{o}} \dfrac{\md \mu_{V}^{o}
(s)}{s-z} \right),
\end{equation}
with $\int_{J_{o}} \tfrac{\md \mu_{V}^{o}(s)}{s-z}$ the Stieltjes transform
of the `odd' equilibrium measure, and
\begin{equation*}
\mathcal{H} \colon \mathcal{L}^{2}_{\mathrm{M}_{2}(\mathbb{C})} \! \to \!
\mathcal{L}^{2}_{\mathrm{M}_{2}(\mathbb{C})}, \, \, f \! \mapsto \! (\mathcal{
H}f)(z) \! := \! \vip_{\raise-1.95ex\hbox{$\scriptstyle{} \mathbb{R}$}}
\dfrac{f(s)}{z \! - \! s} \, \dfrac{\md s}{\pi}
\end{equation*}
denotes the Hilbert transform, with $\pvi_{}$ denoting the principle value
integral. Via the distributional identities $\tfrac{1}{x-(x_{0} \pm \mi 0)} \!
= \! \tfrac{1}{x-x_{0}} \! \pm \! \mi \pi \delta (x \! - \! x_{0})$, with
$\delta (\pmb{\cdot})$ the Dirac delta function, and $\int_{\xi_{1}}^{\xi_{2}}
f(\xi) \delta (\xi \! - \! x) \md \xi \! = \!
\begin{cases}
f(x), &\text{$x \! \in \! (\xi_{1},\xi_{2})$,} \\
0, &\text{$x \! \in \! \mathbb{R} \setminus (\xi_{1},\xi_{2}),$}
\end{cases}$ it follows that
\begin{equation*}
\mathscr{H}^{o}_{\pm}(z) \! = \!
\begin{cases}
(\mathscr{F}^{o}_{\pm}(z))^{2} \! - \! \pvi_{\raise-0.95ex\hbox{$\scriptstyle{}
J_{o}$}} \tfrac{(4 \mi (2+\frac{1}{n})^{2} \psi_{V}^{o}(\xi)(\mathcal{H} \psi_{
V}^{o})(\xi)-\frac{4 \mi}{\pi \xi}(2+\frac{1}{n}) \psi_{V}^{o}(\xi))}{(\xi -
z)} \, \tfrac{\md \xi}{2 \pi \mi} \! \mp \! \mathfrak{Y}(z), &\text{$z \! \in
\! J_{o}$,} \\
(\mathscr{F}^{o}_{\pm}(z))^{2} \! - \! \int_{J_{o}} \tfrac{(4 \mi (2+\frac{1}{
n})^{2} \psi_{V}^{o}(\xi)(\mathcal{H} \psi_{V}^{o})(\xi)-\frac{4 \mi}{\pi \xi}
(2+\frac{1}{n}) \psi_{V}^{o}(\xi))}{(\xi -z)} \, \tfrac{\md \xi}{2 \pi \mi},
&\text{$z \! \notin \! J_{o}$,}
\end{cases}
\end{equation*}
where
\begin{equation*}
\mathfrak{Y}(z) \! := \! \dfrac{1}{2} \! \left(4 \mi \! \left(2 \! + \! \dfrac{
1}{n} \right)^{2} \psi_{V}^{o}(z)(\mathcal{H} \psi_{V}^{o})(z) \! - \! \dfrac{
4 \mi}{\pi z} \! \left(2 \! + \! \dfrac{1}{n} \right) \! \psi_{V}^{o}(z)
\right),
\end{equation*}
and $\star^{o}_{\pm}(z) \! := \! \lim_{\varepsilon \downarrow 0} \star^{o}(z
\! \pm \! \mi 0)$, $\star \! \in \! \{\mathcal{H},\mathscr{F}\}$. Recall the
definition of $g^{o}(z)$ given in Lemma~3.4: for $n \! \in \! \mathbb{N}$,
\begin{equation*}
g^{o}(z) \! := \! \int_{J_{o}} \ln \! \left(\dfrac{(z \! - \! s)^{2+\frac{1}{
n}}}{zs} \right) \! \md \mu_{V}^{o}(s) \! = \! \int_{J_{o}} \ln \! \left(
\dfrac{(z \! - \! s)^{2+\frac{1}{n}}}{zs} \right) \! \psi_{V}^{o}(s) \, \md s,
\quad z \! \in \! \mathbb{C} \setminus (-\infty,\max \{0,\max \{J_{o}\}\});
\end{equation*}
noting the above distributional identities and the fact that $\int_{J_{o}}
\psi_{V}^{o}(s) \, \md s \! = \! 1$, one shows that
\begin{equation*}
(g^{o}_{\pm}(z))^{\prime} \! := \! \lim_{\varepsilon \downarrow 0}(g^{o})^{
\prime}(z \! \pm \! \mi \varepsilon) \! = \!
\begin{cases}
-\tfrac{1}{z} \! - \! (2 \! + \! \tfrac{1}{n})
\pvi_{\raise-1.05ex\hbox{$\scriptstyle{}J_{o}$}} \tfrac{\psi_{V}^{o}(s)}{s-z}
\, \md s \! \mp \! (2 \! + \! \tfrac{1}{n}) \pi \mi \psi_{V}^{o}(z),
&\text{$z \! \in \! J_{o}$,} \\
-\tfrac{1}{z} \! - \! (2 \! + \! \tfrac{1}{n}) \int_{J_{o}} \tfrac{\psi_{V}^{o}
(s)}{s-z} \, \md s, &\text{$z \! \notin \! J_{o}$,}
\end{cases}
\end{equation*}
whence one concludes that
\begin{align*}
(g^{o}_{+} \! + \! g^{o}_{-})^{\prime}(z) =& \, -\dfrac{2}{z} \! - \! 2 \!
\left(2 \! + \! \dfrac{1}{n} \right) \!
\vip_{\raise-1.95ex\hbox{$\scriptstyle{}J_{o}$}} \dfrac{\psi_{V}^{o}(s)}{s \!
- \! z} \, \md s \! = \! -\dfrac{2}{z} \! + \! 2 \! \left(2 \! + \! \dfrac{1}{
n} \right) \! \pi (\mathcal{H} \psi_{V}^{o})(z), \quad z \! \in \! J_{o}, \\
(g^{o}_{+} \! - \! g^{o}_{-})^{\prime}(z) =& \,
\begin{cases}
-2 \! \left(2 \! + \! \dfrac{1}{n} \right) \! \pi \mi \psi_{V}^{o}(z),
&\text{$z \! \in \! J_{o}$,} \\
0, &\text{$z \! \notin \! J_{o}$.}
\end{cases}
\end{align*}
Demanding that (see Lemma~3.6 below) $(g^{o}_{+} \! + \! g^{o}_{-})^{\prime}
(z) \! = \! \widetilde{V}^{\prime}(z)$, $z \! \in \! J_{o}$, one shows {}from
the above that, for $J_{o} \! \ni \! z$, $((g^{o}(z))^{\prime} \! + \! \tfrac{
1}{z})_{+} \! + \! ((g^{o}(z))^{\prime} \! + \! \tfrac{1}{z})_{-} \! = \! 2(2
\! + \! \tfrac{1}{n}) \pi (\mathcal{H} \psi_{V}^{o})(z) \! = \! \tfrac{2}{z}
\! + \! \widetilde{V}^{\prime}(z)$ $\Rightarrow$
\begin{equation}
(\mathcal{H} \psi_{V}^{o})(z) \! = \! \dfrac{1}{2(2 \! + \! \frac{1}{n}) \pi}
\! \left(\dfrac{2}{z} \! + \! \widetilde{V}^{\prime}(z) \right), \quad z \!
\in \! J_{o}.
\end{equation}
{}From Equation~(3.2) and the distributional identities above, one shows that
\begin{equation}
\mathscr{F}^{o}_{\pm}(z) \! := \! \lim_{\varepsilon \downarrow 0} \mathscr{
F}^{o}(z \! \pm \! \mi \varepsilon) \! = \!
\begin{cases}
-\tfrac{1}{\pi \mi z} \! - \! \mi (2 \! + \! \tfrac{1}{n})(\mathcal{H} \psi_{
V}^{o})(z) \! \mp \! (2 \! + \! \tfrac{1}{n}) \psi_{V}^{o}(z), &\text{$z \!
\in \! J_{o}$,} \\
-\tfrac{1}{\pi \mi} \! \left(\tfrac{1}{z} \! + \! (2 \! + \! \tfrac{1}{n})
\int_{J_{o}} \tfrac{\psi_{V}^{o}(s)}{s-z} \, \md s \right), &\text{$z \!
\notin \! J_{o}$;}
\end{cases}
\end{equation}
thus, for $z \! \in \! \mathbb{R} \setminus (J_{o} \cup \{0\})$, $\mathscr{F}^{
o}_{+}(z) \! = \! \mathscr{F}^{o}_{-}(z) \! = \! -\tfrac{1}{\pi \mi}(\tfrac{1}{
z} \! + \! (2 \! + \! \tfrac{1}{n}) \int_{J_{o}} \tfrac{\psi_{V}^{o}(s)}{s-z}
\, \md s)$. Hence, for $z \! \notin \! J_{o} \cup \{0\}$, one deduces that
$\mathscr{H}^{o}_{+}(z) \! = \! \mathscr{H}^{o}_{-}(z)$. For $z \! \in \! J_{
o}$, one notes that
\begin{equation*}
\mathscr{H}^{o}_{+}(z) \! - \! \mathscr{H}^{o}_{-}(z) \! = \! (\mathscr{F}^{
o}_{+}(z))^{2} \! - \! (\mathscr{F}^{o}_{-}(z))^{2} \! - \! 4 \mi \! \left(2
\! + \! \dfrac{1}{n} \right)^{2} \! \psi_{V}^{o}(z)(\mathcal{H} \psi_{V}^{o})
(z) \! + \! \dfrac{4 \mi}{\pi z} \! \left(2 \! + \! \dfrac{1}{n} \right) \!
\psi_{V}^{o}(z),
\end{equation*}
and
\begin{align*}
(\mathscr{F}^{o}_{\pm}(z))^{2} &= -\dfrac{1}{\pi^{2}z^{2}} \! + \! \dfrac{2}{
\pi z} \! \left(2 \! + \! \dfrac{1}{n} \right) \! (\mathcal{H} \psi_{V}^{o})
(z) \! \mp \! \dfrac{2 \mi}{\pi z} \! \left(2 \! + \! \dfrac{1}{n} \right) \!
\psi_{V}^{o}(z) \! - \! \left(2 \! + \! \dfrac{1}{n} \right)^{2} \! ((
\mathcal{H} \psi_{V}^{o})(z))^{2} \\
&\pm \, 2 \mi \! \left(2 \! + \! \dfrac{1}{n} \right)^{2} \! \psi_{V}^{o}(z)
(\mathcal{H} \psi_{V}^{o})(z) \! + \! \left(2 \! + \! \dfrac{1}{n} \right)^{2}
\! (\psi_{V}^{o}(z))^{2},
\end{align*}
whence $(\mathscr{F}^{o}_{+}(z))^{2} \! - \! (\mathscr{F}^{o}_{-}(z))^{2} \!
= \! -\tfrac{4 \mi}{\pi z}(2 \! + \! \tfrac{1}{n}) \psi_{V}^{o}(z) \! + \! 4
\mi (2 \! + \! \tfrac{1}{n})^{2} \psi_{V}^{o}(z)(\mathcal{H} \psi_{V}^{o})(z)$
$\Rightarrow$ $\mathscr{H}^{o}_{+}(z) \! - \! \mathscr{H}^{o}_{-}(z) \! = \!
0$; thus, for $z \! \in \! J_{o}$, $\mathscr{H}^{o}_{+}(z) \! = \! \mathscr{
H}^{o}_{-}(z)$. The above argument shows, therefore, that $\mathscr{H}^{o}(z)$
is analytic across $\mathbb{R} \setminus \{0\}$; in fact, $\mathscr{H}^{o}(z)$
is entire for $z \! \in \! \mathbb{C}^{\ast}$. Recalling that $\mu_{V}^{o} \!
\in \! \mathcal{M}_{1}(\mathbb{R})$, in particular, $\int_{J_{o}}s^{-m} \,
\md \mu_{V}^{o}(s) \! = \! \int_{J_{o}}s^{-m} \psi_{V}^{o}(s) \, \md s \! <
\! \infty$, $m \! \in \! \mathbb{N}$, one shows that, for $\vert z/s \vert \!
\ll \! 1$, with $s \! \in \! J_{o}$ and $z \! \notin \! J_{o}$, via the
expansion $\tfrac{1}{z-s} \! = \! -\sum_{k=0}^{l} \tfrac{z^{k}}{s^{k+1}} \!
+ \! \tfrac{z^{l+1}}{s^{l+1}(z-s)}$, $l \! \in \! \mathbb{Z}_{0}^{+}$,
\begin{equation*}
(\mathscr{F}^{o}(z))^{2} \! \underset{z \to 0}{=} \! -\dfrac{1}{\pi^{2}z^{2}}
\! - \! \dfrac{1}{z} \! \left(\dfrac{2}{\pi^{2}} \! \left(2 \! + \! \dfrac{1}{
n} \right) \! \int_{J_{o}}s^{-1} \, \md \mu_{V}^{o}(s) \right) \! + \!
\mathcal{O}(1),
\end{equation*}
whence, upon recalling the definition of $\mathscr{H}^{o}(z)$, in particular,
for $\vert z/\xi \vert \! \ll \! 1$, with $\xi \! \in \! J_{o}$ and $z \!
\notin \! J_{o}$, via the expansion $\tfrac{1}{z-\xi} \! = \! -\sum_{k=0}^{l}
\tfrac{z^{k}}{\xi^{k+1}} \! + \! \tfrac{z^{l+1}}{\xi^{l+1}(z-\xi)}$, $l \! \in
\! \mathbb{Z}_{0}^{+}$,
\begin{equation*}
\int_{J_{o}} \dfrac{(4 \mi (2 \! + \! \frac{1}{n})^{2} \psi_{V}^{o}(\xi)
(\mathcal{H} \psi_{V}^{o})(\xi) \! - \! \frac{4 \mi}{\pi \xi}(2 \! + \! \frac{
1}{n}) \psi_{V}^{o}(\xi))}{(\xi \! - \! z)} \, \dfrac{\md \xi}{2 \pi \mi}
\underset{z \to 0}{=} \mathcal{O}(1),
\end{equation*}
it follows that
\begin{equation*}
\mathscr{H}^{o}(z) \underset{z \to 0}{=} -\dfrac{1}{\pi^{2}z^{2}} \! - \!
\dfrac{1}{z} \! \left(\dfrac{2}{\pi^{2}} \! \left(2 \! + \! \dfrac{1}{n}
\right) \! \int_{J_{o}}s^{-1} \, \md \mu_{V}^{o}(s) \right) \! + \! \mathcal{O}
(1),
\end{equation*}
which shows that $\mathscr{H}^{o}(z)$ has a pole of order $2$ at $z \! = \!
0$, with $\operatorname{Res}(\mathscr{H}^{o}(z);0) \! = \! -\tfrac{2}{\pi^{2}}
(2 \! + \! \tfrac{1}{n}) \int_{J_{o}}s^{-1} \, \md \mu_{V}^{o}(s)$. One learns
{}from the above analysis that $z^{2} \mathscr{H}^{o}(z)$ is entire: look, in
particular, at the behaviour of $z^{2} \mathscr{H}^{o}(z)$ as $\vert z \vert
\! \to \! \infty$. Recalling Equations~(3.1) and~(3.2), one shows that, for
$\mu_{V}^{o} \! \in \! \mathcal{M}_{1}(\mathbb{R})$, in particular, $\int_{J_{
o}} \md \mu_{V}^{o}(s) \! = \! 1$ and $\int_{J_{o}}s^{m} \, \md \mu_{V}^{o}(s)
\! < \! \infty$, $m \! \in \! \mathbb{N}$, for $\vert s/z \vert \! \ll \! 1$,
with $s \! \in \! J_{o}$ and $z \! \notin \! J_{o}$, via the expansion $\tfrac{
1}{s-z} \! = \! -\sum_{k=0}^{l} \tfrac{s^{k}}{z^{k+1}} \! + \! \tfrac{s^{l+1}}{
z^{l+1}(s-z)}$, $l \! \in \! \mathbb{Z}_{0}^{+}$,
\begin{align*}
z^{2} \mathscr{H}^{o}(z) \, +& \, \dfrac{1}{\pi^{2}} \! \left(1 \! + \! \tfrac{
1}{n} \right)^{2} \! - \! \int_{J_{o}}s \! \left(4 \mi \! \left(2 \! + \!
\tfrac{1}{n} \right)^{2} \! \psi_{V}^{o}(s)(\mathcal{H} \psi_{V}^{o})(s) \! -
\! \tfrac{4 \mi}{\pi s} \! \left(2 \! + \! \tfrac{1}{n} \right) \! \psi_{V}^{o}
(s) \right) \! \dfrac{\md s}{2 \pi \mi} \\
-& \, z \int_{J_{o}} \! \left(4 \mi \! \left(2 \! + \! \tfrac{1}{n} \right)^{
2} \! \psi_{V}^{o}(s)(\mathcal{H} \psi_{V}^{o})(s) \! - \! \tfrac{4 \mi}{\pi
s} \! \left(2 \! + \! \tfrac{1}{n} \right) \! \psi_{V}^{o}(s) \right) \!
\dfrac{\md s}{2 \pi \mi} \underset{\vert z \vert \to \infty}{=} \mathcal{O}
(z^{-1});
\end{align*}
thus, due to the entirety of $\mathscr{H}^{o}(z)$, it follows, by a 
generalisation of Liouville's Theorem, that
\begin{align*}
z^{2} \mathscr{H}^{o}(z) &+ \! \dfrac{1}{\pi^{2}} \! \left(1 \! + \! \tfrac{
1}{n} \right)^{2} \! - \! \int_{J_{o}} \! s \! \left(4 \mi \! \left(2 \! + \!
\tfrac{1}{n} \right)^{2} \! \psi_{V}^{o}(s)(\mathcal{H} \psi_{V}^{o})(s) \! -
\! \tfrac{4 \mi}{\pi s} \! \left(2 \! + \! \tfrac{1}{n} \right) \! \psi_{V}^{o}
(s) \right) \! \dfrac{\md s}{2 \pi \mi} \\
&- \, z \int_{J_{o}} \! \left(4 \mi \! \left(2 \! + \! \tfrac{1}{n} \right)^{
2} \! \psi_{V}^{o}(s)(\mathcal{H} \psi_{V}^{o})(s) \! - \! \tfrac{4 \mi}{\pi
s} \! \left(2 \! + \! \tfrac{1}{n} \right) \! \psi_{V}^{o}(s) \right) \!
\dfrac{\md s}{2 \pi \mi} \! = \! 0.
\end{align*}
Substituting Equation~(3.1) into the above formula, one notes that
\begin{gather*}
(\mathscr{F}^{o}(z))^{2} \! - \! \int_{J_{o}} \dfrac{(4 \mi (2 \! + \! \frac{
1}{n})^{2} \psi_{V}^{o}(\xi)(\mathcal{H} \psi_{V}^{o})(\xi) \! - \! \frac{4
\mi}{\pi \xi}(2 \! + \! \frac{1}{n}) \psi_{V}^{o}(\xi))}{(\xi \! - \! z)} \,
\dfrac{\md \xi}{2 \pi \mi} \! - \! \dfrac{1}{z} \int_{J_{o}} \! \left(4 \mi \!
\left(2 \! + \! \tfrac{1}{n} \right)^{2} \! \psi_{V}^{o}(\xi)(\mathcal{H}
\psi_{V}^{o})(\xi) \right. \\
\left. - \, \tfrac{4 \mi}{\pi \xi} \! \left(2 \! + \! \tfrac{1}{n} \right) \!
\psi_{V}^{o}(\xi) \right) \! \tfrac{\md \xi}{2 \pi \mi} \! + \! \tfrac{(1+
\frac{1}{n})^{2}}{\pi^{2}z^{2}} \! - \! \tfrac{1}{z^{2}} \int_{J_{o}} \xi \!
\left(4 \mi \! \left(2 \! + \! \tfrac{1}{n} \right)^{2} \! \psi_{V}^{o}(\xi)
(\mathcal{H} \psi_{V}^{o})(\xi) \! - \! \tfrac{4 \mi}{\pi \xi} \! \left(2 \! +
\! \tfrac{1}{n} \right) \! \psi_{V}^{o}(\xi) \right) \! \tfrac{\md \xi}{2 \pi
\mi} \! = \! 0.
\end{gather*}
Via Equation~(3.3), it follows that $4 \mi (2 \! + \! \tfrac{1}{n})^{2} \psi_{
V}^{o}(\xi)(\mathcal{H} \psi_{V}^{o})(\xi) \! - \! \tfrac{4 \mi}{\pi \xi}(2
\! + \! \tfrac{1}{n}) \psi_{V}^{o}(\xi) \! = \! \tfrac{2 \mi}{\pi}(2 \! + \!
\tfrac{1}{n}) \psi_{V}^{o}(\xi) \widetilde{V}^{\prime}(\xi)$; substituting the
latter expression into the above equation, and re-arranging, one obtains,
\begin{equation}
(\mathscr{F}^{o}(z))^{2} \! - \! \tfrac{(2+\frac{1}{n})}{\pi^{2}} \int_{J_{o}}
\tfrac{\widetilde{V}^{\prime}(\xi) \psi_{V}^{o}(\xi)}{\xi -z} \, \md \xi \! +
\! \tfrac{(1+\frac{1}{n})^{2}}{\pi^{2}z^{2}} \! - \! \tfrac{(2+\frac{1}{n})}{
\pi^{2}z^{2}} \int_{J_{o}} \xi \widetilde{V}^{\prime}(\xi) \psi_{V}^{o}(\xi)
\, \md \xi \! - \! \tfrac{(2+\frac{1}{n})}{\pi^{2}z} \int_{J_{o}} \widetilde{
V}^{\prime}(\xi) \psi_{V}^{o}(\xi) \, \md \xi \! = \! 0.
\end{equation}
But
\begin{align*}
\dfrac{(2 \! + \! \frac{1}{n})}{\pi^{2}} \int_{J_{o}} \dfrac{\widetilde{V}^{
\prime}(\xi) \psi_{V}^{o}(\xi)}{\xi \! - \! z} \, \md \xi =& \, \dfrac{(2 \! +
\! \frac{1}{n})}{\pi^{2}} \int_{J_{o}} \dfrac{(\widetilde{V}^{\prime}(\xi) \!
- \! \widetilde{V}^{\prime}(z)) \psi_{V}^{o}(\xi)}{\xi \! - \! z} \, \md \xi
\! + \! \dfrac{(2 \! + \! \frac{1}{n})}{\pi^{2}} \int_{J_{o}} \dfrac{
\widetilde{V}^{\prime}(z) \psi_{V}^{o}(\xi)}{\xi \! - \! z} \, \md \xi \\
=& \, \dfrac{(2 \! + \! \frac{1}{n})}{\pi^{2}} \int_{J_{o}} \dfrac{(\widetilde{
V}^{\prime}(\xi) \! - \! \widetilde{V}^{\prime}(z)) \psi_{V}^{o}(\xi)}{\xi \!
- \! z} \, \md \xi \! + \! \dfrac{\widetilde{V}^{\prime}(z)}{\pi^{2}} \!
\underbrace{\left(\! \left(2 \! + \! \tfrac{1}{n} \right) \! \int_{J_{o}}
\dfrac{\psi_{V}^{o}(\xi)}{\xi \! - \! z} \, \md \xi \right)}_{= \, -\mi \pi
\mathscr{F}^{o}(z)-z^{-1}} \\
=& \, \dfrac{(2 \! + \! \frac{1}{n})}{\pi^{2}} \int_{J_{o}} \dfrac{(\widetilde{
V}^{\prime}(\xi) \! - \! \widetilde{V}^{\prime}(z)) \psi_{V}^{o}(\xi)}{\xi \!
- \! z} \, \md \xi \! - \! \dfrac{\mi \widetilde{V}^{\prime}(z) \mathscr{F}^{o}
(z)}{\pi} \! - \! \dfrac{\widetilde{V}^{\prime}(z)}{\pi^{2}z}:
\end{align*}
substituting the above into Equation~(3.5), one arrives at, upon completing
the square and re-arran\-g\-i\-n\-g terms,
\begin{equation}
\left(\mathscr{F}^{o}(z) \! + \! \dfrac{\mi \widetilde{V}^{\prime}(z)}{2 \pi}
\right)^{2} \! + \! \dfrac{\mathfrak{q}_{V}^{o}(z)}{\pi^{2}} \! = \! 0,
\end{equation}
where
\begin{align*}
\mathfrak{q}_{V}^{o}(z) :=& \, \left(\dfrac{\widetilde{V}^{\prime}(z)}{2}
\right)^{2} \! + \! \dfrac{\widetilde{V}^{\prime}(z)}{z} \! - \! \left(2 \! +
\! \dfrac{1}{n} \right) \! \int_{J_{o}} \dfrac{(\widetilde{V}^{\prime}(\xi) \!
- \! \widetilde{V}^{\prime}(z)) \psi_{V}^{o}(\xi)}{\xi \! - \! z} \, \md \xi \\
+& \, \dfrac{1}{z^{2}} \! \left(\! \left(1 \! + \! \dfrac{1}{n} \right)^{2} \!
- \! \left(2 \! + \! \dfrac{1}{n} \right) \! \int_{J_{o}}(\xi \! + \! z)
\widetilde{V}^{\prime}(\xi) \psi_{V}^{o}(\xi) \, \md \xi \right).
\end{align*}
(Equation~(3.6) above generalises Equation~(3.5) for $q^{(0)}(x)$ in 
\cite{a46} for the case when $\widetilde{V} \colon \mathbb{R} \setminus 
\{0\} \! \to \! \mathbb{R}$ is real analytic; moreover, it is analogous to 
Equation~(1.37) of \cite{a44}.) Note that, since $\widetilde{V} \colon \mathbb{
R} \setminus \{0\} \! \to \! \mathbb{R}$ satisfies conditions (2.3)--(2.5), 
it follows {}from $\alpha^{l} \! - \! \beta^{l} \! = \! (\alpha \! - \! \beta)
(\alpha^{l-1} \! + \! \alpha^{l-2} \beta \! + \! \cdots \! + \! \alpha \beta^{
l-2} \! + \! \beta^{l-1})$, $l \! \in \! \mathbb{N}$, that $\mathfrak{q}_{V}^{
o}(z)$ is real analytic on $J_{o}$ (and real analytic on $\mathbb{R} \setminus
\{0\})$. For $x \! \in \! J_{o}$, set $z \! := \! x \! + \! \mi \varepsilon$,
and consider the $\varepsilon \! \downarrow \! 0$ limit of Equation~(3.6):
$\lim_{\varepsilon \downarrow 0}(\mathscr{F}^{o}(x \! + \! \mi \varepsilon) \!
+ \! \tfrac{\mi \widetilde{V}^{\prime}(x+\mi \varepsilon)}{2 \pi})^{2} \! = \!
(\mathscr{F}^{o}_{+}(x) \! + \! \tfrac{\mi \widetilde{V}^{\prime}(x)}{2 \pi})^{
2}$ (as $\widetilde{V}$ is real analytic on $J_{o})$; recalling that $\mathscr{
F}^{o}_{+}(x) \! = \! -\tfrac{1}{\pi \mi x} \! - \! \mi (2 \! + \! \tfrac{1}{n}
)(\mathcal{H} \psi_{V}^{o})(x) \! - \! (2 \! + \! \tfrac{1}{n}) \psi_{V}^{o}
(x)$, via Equation~(3.3), it follows that $\mathscr{F}^{o}_{+}(x) \! = \!
-\tfrac{\mi \widetilde{V}^{\prime}(x)}{2 \pi} \! - \! (2 \! + \! \tfrac{1}{n})
\psi_{V}^{o}(x)$ $\Rightarrow$ $(\mathscr{F}^{o}_{+}(x) \! + \! \tfrac{\mi
\widetilde{V}^{\prime}(x)}{2 \pi})^{2} \! = \! ((2 \! + \! \tfrac{1}{n}) \psi_{
V}^{o}(x))^{2}$, whence $(\psi_{V}^{o}(x))^{2} \! = \! -\mathfrak{q}_{V}^{o}
(x)/((2 \! + \! \tfrac{1}{n}) \pi)^{2}$, $x \! \in \! J_{o}$, whereupon, using
the fact that (see above) $\psi_{V}^{o}(x) \! \geqslant \! 0 \, \, \forall \,
\, x \! \in \! \overline{J_{o}}$, it follows that $\mathfrak{q}_{V}^{o}(x) \!
\leqslant \! 0$, $x \! \in \! J_{o}$; moreover, as a by-product, decomposing
$\mathfrak{q}_{V}^{o}(x)$, for $x \! \in \! J_{o}$, into positive and negative
parts, that is, $\mathfrak{q}_{V}^{o}(x) \! = \! (\mathfrak{q}_{V}^{o}(x))^{+}
\! - \! (\mathfrak{q}_{V}^{o}(x))^{-}$, $x \! \in \! J_{o}$, where
$(\mathfrak{q}_{V}^{o}(x))^{\pm} \! := \! \max \left\lbrace \pm \mathfrak{q}_{
V}^{o}(x),0 \right\rbrace$ $(\geqslant \! 0)$, one learns {}from the above
analysis that, for $x \! \in \! J_{o}$, $(\mathfrak{q}_{V}^{o}(x))^{+} \!
\equiv \! 0$ and $\psi_{V}^{o}(x) \! = \! \tfrac{1}{(2+\frac{1}{n}) \pi}
((\mathfrak{q}_{V}^{o}(x))^{-})^{1/2}$; and, since $\int_{J_{o}} \! \psi_{V}^{
o}(s) \, \md s \! = \! 1$, it follows that $\tfrac{1}{(2+\frac{1}{n}) \pi}
\int_{J_{o}}((\mathfrak{q}_{V}^{o}(s))^{-})^{1/2} \, \md s \! = \! 1$, which
gives rise to the interesting fact that the function $(\mathfrak{q}_{V}^{o}
(x))^{-} \! \not\equiv \! 0$ on $J_{o}$. (Even though $(\mathfrak{q}_{V}^{o}
(x))^{-}$ depends on $\md \mu_{V}^{o}(x) \! = \! \psi_{V}^{o}(x) \, \md x$,
and thus $\psi_{V}^{o}(x) \! = \! \tfrac{1}{(2+\frac{1}{n}) \pi}((\mathfrak{
q}_{V}^{o}(x))^{-})^{1/2}$ is an implicit representation for $\psi_{V}^{o}$,
it is still a useful relation which can be used to obtain additional, valuable
information about $\psi_{V}^{o}$.) For $x \! \notin \! J_{o}$, set $z \! := \!
x \! + \! \mi \varepsilon$, and (again) study the $\varepsilon \! \downarrow
\! 0$ limit of Equation~(3.6): in this case, $\lim_{\varepsilon \downarrow 0}
(\mathscr{F}^{o}(x \! + \! \mi \varepsilon) \! + \! \tfrac{\mi \widetilde{V}^{
\prime}(x+\mi \varepsilon)}{2 \pi})^{2} \! = \! (\mathscr{F}^{o}_{+}(x) \! +
\! \tfrac{\mi \widetilde{V}^{\prime}(x)}{2 \pi})^{2} \! = \! (\mathscr{F}^{o}
(x) \! + \! \tfrac{\mi \widetilde{V}^{\prime}(x)}{2 \pi})^{2}$; recalling
that, for $x \! \notin \! J_{o}$, $\mathscr{F}^{o}(x) \! = \! -\tfrac{1}{\pi
\mi}(\tfrac{1}{x} \! + \! (2 \! + \! \tfrac{1}{n}) \int_{J_{o}} \tfrac{\psi_{
V}^{o}(s)}{s-x} \, \md s) \! = \! \tfrac{\mi}{\pi x} \! - \! \mi (2 \! + \!
\tfrac{1}{n})(\mathcal{H} \psi_{V}^{o})(x)$, substituting the latter
expression into Equation~(3.6), one arrives at $(\tfrac{1}{\pi x} \! - \! (2
\! + \! \tfrac{1}{n})(\mathcal{H} \psi_{V}^{o})(x) \! + \! \tfrac{\widetilde{
V}^{\prime}(x)}{2 \pi})^{2} \! = \! \mathfrak{q}_{V}^{o}(x)/\pi^{2}$, $x
\! \notin \! J_{o}$ (since $\widetilde{V}^{\prime}$ is real analytic on
$(\mathbb{R} \setminus \{0\}) \setminus J_{o}$, it follows that $\mathfrak{
q}_{V}^{o}(x)$, too, is real analytic on $(\mathbb{R} \setminus \{0\})
\setminus J_{o}$, in which case, this latter relation merely states that, for
$x \! = \! 0$, $+\infty \! = \! +\infty)$, whence $\mathfrak{q}_{V}^{o}(x) 
\! \geqslant \! 0 \, \, \forall \, \, x \! \notin \! J_{o}$.

Now, recalling that, on a compact subset of $\mathbb{R}$, an analytic function
changes sign an at most countable number of times, it follows {}from the above
argument, the fact that $\widetilde{V} \colon \mathbb{R} \setminus \{0\}
\! \to \! \mathbb{R}$ satisfying conditions (2.3)--(2.5) is regular (cf.
Subsection~2.2), in particular, $\widetilde{V}$ is real analytic in the (open)
neighbourhood $\mathbb{U} \! := \! \lbrace \mathstrut z \! \in \! \mathbb{C};
\, \inf_{q \in J_{o}} \vert z \! - \! q \vert \! < \! r \! \in \! (0,1)
\rbrace \setminus \{0\}$, $\mu_{V}^{o}$ has compact support, and mimicking a
part of the calculations subsumed in the proof of Theorem~1.38 in \cite{a44},
that $J_{o} \! := \! \operatorname{supp}(\mu_{V}^{o}) \! = \! \lbrace
\mathstrut x \! \in \! \mathbb{R}; \, q_{V}^{o}(x) \! \leqslant \! 0 \rbrace$
consists of the disjoint union of a finite number of bounded (real) intervals,
with representation $J_{o} \! := \! \cup_{j=1}^{N+1}J_{j}^{o}$, where $J_{j}^{
o} \! := \! [b_{j-1}^{o},a_{j}^{o}]$, with $N \! \in \! \mathbb{N}$ and
finite, $b_{0}^{o} \! := \! \min \lbrace J_{o} \rbrace \! \notin \! \lbrace
-\infty,0 \rbrace$, $a_{N+1}^{o} \! := \! \max \lbrace J_{o} \rbrace \!
\notin \! \lbrace 0,+\infty \rbrace$, and $-\infty \! < \! b_{0}^{o} \! < \!
a_{1}^{o} \! < \! b_{1}^{o} \! < \! a_{2}^{o} \! < \! \cdots \! < \! b_{N}^{o}
\! < \! a_{N+1}^{o} \! < \! +\infty$. (One notes that $\widetilde{V}$ is real
analytic in, say, the open neighbourhood $\widetilde{\mathbb{U}} \! := \!
\cup_{j=1}^{N+1} \widetilde{\mathbb{U}}_{j}$, where $\widetilde{\mathbb{U}}_{
j} \! := \! \lbrace \mathstrut z \! \in \! \mathbb{C}^{\ast}; \, \inf_{q \in
J_{j}^{o}} \vert z \! - \! q \vert \! < \! r_{j} \! \in \! (0,1) \rbrace$,
with $\widetilde{\mathbb{U}}_{i} \cap \widetilde{\mathbb{U}}_{j} \! = \!
\varnothing$, $i \! \not= \! j \! = \! 1,\dotsc,N \! + \! 1$.) Furthermore, as
a by-product of the above representation for $J_{o}$, it follows that, since
$J_{i}^{o} \cap J_{j}^{o} \! = \! \varnothing$, $i \! \not= \! j \! = \! 1,
\dotsc,N \! + \! 1$, $\text{meas}(J_{o}) \! = \! \sum_{j=1}^{N+1} \vert b_{j
-1}^{o} \! - \! a_{j}^{o} \vert \! < \! +\infty$.

It remains, still, to determine the $2(N \! + \! 1)$ conditions satisfied by
the end-points of the support of the `odd' equilibrium measure, $\lbrace b_{j
-1}^{o},a_{j}^{o} \rbrace_{j=1}^{N+1}$. Towards this end, one proceeds as
follows. {}From the formula for $\mathscr{F}^{o}(z)$ given in Equation~(3.2):
\begin{compactenum}
\item[(i)] for $\mu_{V}^{o} \! \in \! \mathcal{M}_{1}(\mathbb{R})$, in
particular, $\int_{\mathbb{R}} \md \mu_{V}^{o}(s) \! = \! 1$ and $\int_{
\mathbb{R}}s^{m} \, \md \mu_{V}^{o}(s) \! < \! \infty$, $m \! \in \! \mathbb{
N}$, $s \! \in \! J_{o}$ and $z \! \notin \! J_{o}$, with $\vert s/z \vert \!
\ll \! 1$ (e.g., $\vert z \vert \! \gg \! \max_{j=1,\dotsc,N+1} \lbrace \vert
b_{j-1}^{o} \! - \! a_{j}^{o} \vert \rbrace)$, via the expansion $\tfrac{1}{s
-z} \! = \! -\sum_{k=0}^{l} \tfrac{s^{k}}{z^{k+1}} \! + \! \tfrac{s^{l+1}}{z^{
l+1}(s-z)}$, $l \! \in \! \mathbb{Z}_{0}^{+}$, one gets that $\mathscr{F}^{o}
(z) \! =_{z \to \infty} \! \tfrac{1}{\pi \mi}(1 \! + \! \tfrac{1}{n}) \tfrac{
1}{z} \! + \! \mathcal{O}(z^{-2})$;
\item[(ii)] for $\mu_{V}^{o} \! \in \! \mathcal{M}_{1}(\mathbb{R})$, in 
particular, $\int_{\mathbb{R}}s^{-m} \, \md \mu_{V}^{o}(s) \! < \! \infty$,
$m \! \in \! \mathbb{N}$, $s \! \in \! J_{o}$ and $z \! \notin \! J_{o}$,
with $\vert z/s \vert \! \ll \! 1$ (e.g., $\vert z \vert \! \ll \! \min_{j=1,
\dotsc,N+1} \lbrace \vert b_{j-1}^{o} \! - \! a_{j}^{o} \vert \rbrace)$, via
the expansion $\tfrac{1}{z-s} \! = \! -\sum_{k=0}^{l} \tfrac{z^{k}}{s^{k+1}}
\! + \! \tfrac{z^{l+1}}{s^{l+1}(z-s)}$, $l \! \in \! \mathbb{Z}_{0}^{+}$, one
gets that $\mathscr{F}^{o}(z) \! =_{z \to 0} \! -\tfrac{1}{\pi \mi z} \! + \!
\mathcal{O}(1)$.
\end{compactenum}
Recalling, also, the formulae for $\mathscr{F}^{o}_{\pm}(z)$ given in
Equation~(3.4), one deduces that $\mathscr{F}^{o}_{+}(z) \! + \! \mathscr{F}^{
o}_{-}(z) \! = \! -\mi \widetilde{V}^{\prime}(z)/\pi$, $z \! \in \! J_{o}$,
and $\mathscr{F}^{o}_{+}(z) \! - \! \mathscr{F}^{o}_{-}(z) \! = \! 0$, $z \!
\notin \! J_{o}$; thus, one learns that $\mathscr{F}^{o} \colon \mathbb{C}
\setminus (J_{o} \cup \{0\}) \! \to \! \mathbb{C}$ solves the following
(scalar and homogeneous) RHP:
\begin{compactenum}
\item[(1)] $\mathscr{F}^{o}(z)$ is holomorphic (resp., meromorphic) for $z \!
\in \! \mathbb{C} \setminus (J_{o} \cup \{0\})$ (resp., $z \! \in \! \mathbb{
C} \setminus J_{o})$;
\item[(2)] $\mathscr{F}^{o}_{\pm}(z) \! := \! \lim_{\varepsilon \downarrow 0}
\mathscr{F}^{o}(z \! \pm \! \mi \varepsilon)$ satisfy the boundary condition
$\mathscr{F}^{o}_{+}(z) \! + \! \mathscr{F}^{o}_{-}(z) \! = \! -\mi
\widetilde{V}^{\prime}(z)/\pi$, $z \! \in \! J_{o}$, with $\mathscr{F}^{o}_{
+}(z) \! = \! \mathscr{F}^{o}_{-}(z) \! := \! \mathscr{F}^{o}(z)$ for $z \!
\notin \! J_{o}$;
\item[(3)] $\mathscr{F}^{o}(z) \! =_{\underset{z \in \mathbb{C} \setminus
\mathbb{R}}{z \to \infty}} \! \tfrac{1}{\pi \mi}(1 \! + \! \tfrac{1}{n})
\tfrac{1}{z} \! + \! \mathcal{O}(z^{-2})$; and
\item[(4)] $\operatorname{Res}(\mathscr{F}^{o}(z);0) \! = \! -\tfrac{1}{\pi
\mi}$.
\end{compactenum}
The solution of this RHP is (see, for example, \cite{a83})
\begin{equation*}
\mathscr{F}^{o}(z) \! = \! -\dfrac{1}{\pi \mi z} \! + \! (R_{o}(z))^{1/2}
\int_{J_{o}} \dfrac{(\frac{2}{\mi \pi s} \! + \! \frac{\widetilde{V}^{\prime}
(s)}{\mi \pi})}{(R_{o}(s))^{1/2}_{+}(s \! - \! z)} \, \dfrac{\md s}{2 \pi
\mi}, \quad z \! \in \! \mathbb{C} \setminus (J_{o} \cup \{0\}),
\end{equation*}
where $(R_{o}(z))^{1/2}$ is defined in the Lemma, with $(R_{o}(z))^{1/2}_{\pm}
\! := \! \lim_{\varepsilon \downarrow 0}(R_{o}(z \! \pm \! \mi \varepsilon))^{
1/2}$, and the branch of the square root is chosen so that $z^{-(N+1)}(R_{o}
(z))^{1/2} \! \sim_{\underset{z \in \mathbb{C}_{\pm}}{z \to \infty}} \! \pm
1$. (Note that $(R_{o}(z))^{1/2}$ is pure imaginary on $J_{o}$.) It follows
{}from the above integral representation for $\mathscr{F}^{o}(z)$ that, for
$s \! \in \! J_{o}$ and $z \! \notin \! J_{o}$, with $\vert s/z \vert \! \ll
\! 1$ (e.g., $\vert z \vert \! \gg \! \max_{j=1,\dotsc,N+1} \lbrace \vert b_{
j-1}^{o} \! - \! a_{j}^{o} \vert \rbrace)$, via the expansion $\tfrac{1}{s-z}
\! = \! -\sum_{k=0}^{l} \tfrac{s^{k}}{z^{k+1}} \! + \! \tfrac{s^{l+1}}{z^{l+1}
(s-z)}$, $l \! \in \! \mathbb{Z}_{0}^{+}$,
\begin{equation*}
\mathscr{F}^{o}(z) \underset{z \to \infty}{=} -\dfrac{1}{\mi \pi z} \! + \!
\dfrac{(z^{N+1} \! + \! \dotsb)}{z} \int_{J_{o}} \! \dfrac{(\frac{2 \mi}{\pi
s} \! + \! \frac{\mi \widetilde{V}^{\prime}(s)}{\pi})}{(R_{o}(s))^{1/2}_{+}}
\! \left(1 \! + \! \dfrac{s}{z} \! + \! \cdots \! + \! \dfrac{s^{N}}{z^{N}}
\! + \! \dfrac{s^{N+1}}{z^{N+1}} \! + \! \dotsb \right) \! \dfrac{\md s}{2
\pi \mi}:
\end{equation*}
now, recalling {}from above that $\mathscr{F}^{o}(z) \! =_{\underset{z \in
\mathbb{C}_{+}}{z \to \infty}} \! \tfrac{1}{\pi \mi}(1 \! + \! \tfrac{1}{n})
\tfrac{1}{z} \! + \! \mathcal{O}(z^{-2})$, it follows that, upon removing the
secular (growing) terms,
\begin{equation*}
\int_{J_{o}} \! \left(\dfrac{2 \mi}{\pi s} \! + \! \dfrac{\mi \widetilde{V}^{
\prime}(s)}{\pi} \right) \! \dfrac{s^{j}}{(R_{o}(s))^{1/2}_{+}} \, \dfrac{\md
s}{2 \pi \mi} \! = \! 0, \quad j \! = \! 0,\dotsc,N
\end{equation*}
(which gives $N \! + \! 1$ (real) moment conditions), and, upon equating $z^{
-1}$ terms,
\begin{equation*}
\int_{J_{o}} \! \left(\dfrac{2 \mi}{\pi s} \! + \! \dfrac{\mi \widetilde{V}^{
\prime}(s)}{\pi} \right) \! \dfrac{s^{N+1}}{(R_{o}(s))^{1/2}_{+}} \, \dfrac{
\md s}{2 \pi \mi} \! = \! \dfrac{1}{\mi \pi} \! \left(2 \! + \! \dfrac{1}{n}
\right);
\end{equation*}
it remains, thus, to determine an additional $2(N \! + \! 1) \! - \! (N
\! + \! 1) \! - \! 1 \! = \! N$ moment conditions. {}From the integral
representation for $\mathscr{F}^{o}(z)$, a residue calculation shows that
\begin{equation}
\mathscr{F}^{o}(z) \! = \! -\dfrac{\mi \widetilde{V}^{\prime}(z)}{2 \pi} \!
- \! \dfrac{(R_{o}(z))^{1/2}}{2} \oint_{C_{\mathrm{R}}^{o}} \! \dfrac{(\frac{
2 \mi}{\pi s} \! + \! \frac{\mi \widetilde{V}^{\prime}(s)}{\pi})}{(R_{o}(s))^{
1/2}(s \! - \! z)} \, \dfrac{\md s}{2 \pi \mi},
\end{equation}
where $C_{\mathrm{R}}^{o}$ $(\subset \mathbb{C}^{\ast})$ denotes the boundary
of any open doubly-connected annular region of the type $\lbrace \mathstrut
z^{\prime} \! \in \! \mathbb{C}; \, 0 \! < \! r \! < \! \vert z^{\prime} \vert
\! < \! R \! < \! +\infty \rbrace$, where the simple outer (resp., inner)
boundary $\lbrace \mathstrut z^{\prime} \! = \! R \me^{\mi \vartheta}, \, 0
\! \leqslant \! \vartheta \! \leqslant \! 2 \pi \rbrace$ (resp., $\lbrace
\mathstrut z^{\prime} \! = \! r \me^{\mi \vartheta}, \, 0 \! \leqslant \!
\vartheta \! \leqslant \! 2 \pi \rbrace)$ is traversed clockwise (resp.,
counter-clockwise), with the numbers $0 \! < \! r \! < \! R \! < \! +\infty$
chosen such that, for (any) non-real $z$ in the domain of analyticity of
$\widetilde{V}$ (that is, $\mathbb{C}^{\ast})$, $\mathrm{int}(C_{\mathrm{R}}^{
o}) \! \supset \! J_{o} \cup \{z\}$. Recall {}from Equation~(3.4) that, for
$z \! \in \! \mathbb{R} \setminus \overline{J_{o}}$ $(\supset \cup_{j=1}^{N}
(a_{j}^{o},b_{j}^{o}))$, $\mathscr{F}^{o}_{+}(z) \! = \! \mathscr{F}^{o}_{-}
(z) \! = \! -\tfrac{1}{\pi \mi}(\tfrac{1}{z} \! + \! (2 \! + \! \tfrac{1}{n})
\int_{J_{o}} \tfrac{\psi_{V}^{o}(s)}{s-z} \, \md s)$, whence $\mathscr{F}^{o}
(z) \! + \! \tfrac{1}{\pi \mi z} \! = \! -\mi (2 \! + \! \tfrac{1}{n})
(\mathcal{H} \psi_{V}^{o})(z)$; thus, using Equation~(3.7), one arrives at
\begin{equation*}
\left(2 \! + \! \dfrac{1}{n} \right) \! (\mathcal{H} \psi_{V}^{o})(z) \! =
\! \dfrac{\widetilde{V}^{\prime}(z)}{2 \pi} \! + \! \dfrac{1}{\pi z} \! + \!
\dfrac{\mi (R_{o}(z))^{1/2}}{2} \oint_{C_{\mathrm{R}}^{o}} \! \dfrac{(\frac{
2}{\pi \mi \xi} \! + \! \frac{\widetilde{V}^{\prime}(\xi)}{\pi \mi})}{(R_{o}
(\xi))^{1/2}(\xi \! - \! z)} \, \dfrac{\md \xi}{2 \pi \mi}, \quad z \! \in \!
\cup_{j=1}^{N}(a_{j}^{o},b_{j}^{o}).
\end{equation*}
A contour integration argument shows that
\begin{equation*}
\int_{a_{j}^{o}}^{b_{j}^{o}} \! \left((\mathcal{H} \psi_{V}^{o})(s) \! - \!
\dfrac{1}{2(2+\frac{1}{n}) \pi} \! \left(\dfrac{2}{s} \! + \! \widetilde{V}^{
\prime}(s) \right) \right) \md s \! = \! 0, \quad j \! = \! 1,\dotsc,N,
\end{equation*}
whence, using the above expression for $(\mathcal{H} \psi_{V}^{o})(z)$, $z \!
\in \! \cup_{j=1}^{N}(a_{j}^{o},b_{j}^{o})$, it follows that
\begin{equation}
\int_{a_{j}^{o}}^{b_{j}^{o}} \! \left(\dfrac{\mi (R_{o}(s))^{1/2}}{2} \oint_{
C_{\mathrm{R}}^{o}} \! \dfrac{(\frac{2}{\pi \mi \xi} \! + \! \frac{\widetilde{
V}^{\prime}(\xi)}{\pi \mi})}{(R_{o}(\xi))^{1/2}(\xi \! - \! s)} \, \dfrac{\md
\xi}{2 \pi \mi} \right) \! \md s \! = \! 0, \quad j \! = \! 1,\dotsc,N:
\end{equation}
now, `collapsing' the contour $C_{\mathrm{R}}^{o}$ down to $\mathbb{R}
\setminus \{0\}$ and using the Residue Theorem, one shows that
\begin{equation*}
\dfrac{\mi (R_{o}(z))^{1/2}}{2} \oint_{C_{\mathrm{R}}^{o}} \! \dfrac{(\frac{
2}{\pi \mi \xi} \! + \! \frac{\widetilde{V}^{\prime}(\xi)}{\pi \mi})}{(R_{o}
(\xi))^{1/2}(\xi \! - \! z)} \, \dfrac{\md \xi}{2 \pi \mi} \! = \! -\dfrac{1}{
\pi z} \! - \! \dfrac{\widetilde{V}^{\prime}(z)}{2 \pi} \! + \! \mi (R_{o}(z)
)^{1/2} \int_{J_{o}} \! \dfrac{(\frac{2}{\pi \mi \xi} \! + \! \frac{\widetilde{
V}^{\prime}(\xi)}{\pi \mi})}{(R_{o}(\xi))^{1/2}_{+}(\xi \! - \! z)} \, \dfrac{
\md \xi}{2 \pi \mi};
\end{equation*}
substituting the latter relation into Equation~(3.8), one arrives at, after
straightforward integration and using the Fundamental Theorem of Calculus,
for $j \! = \! 1,\dotsc,N$,
\begin{equation*}
\int_{a_{j}^{o}}^{b_{j}^{o}} \! \left(\mi (R_{o}(s))^{1/2} \int_{J_{o}} \!
\left(\dfrac{\mi}{\pi \xi} \! + \! \dfrac{\mi \widetilde{V}^{\prime}(\xi)}{2
\pi} \right) \! \dfrac{1}{(R_{o}(\xi))^{1/2}_{+}(\xi \! - \! s)} \, \dfrac{\md
\xi}{2 \pi \mi} \right) \md s \! = \! \dfrac{1}{2 \pi} \ln \! \left\vert
\dfrac{a_{j}^{o}}{b_{j}^{o}} \right\vert \! + \! \dfrac{1}{4 \pi} \!
\left(\widetilde{V}(a_{j}^{o}) \! - \! \widetilde{V}(b_{j}^{o}) \right),
\end{equation*}
which give the remaining $N$ moment conditions determining the end-points of
the support of the `odd' equilibrium measure, $\lbrace b_{j-1}^{o},a_{j}^{o}
\rbrace_{j=1}^{N+1}$. Since $J_{o} \not\supseteq \lbrace 0,\pm \infty
\rbrace$ and $\widetilde{V}$ is real analytic on $J_{o}$,
\begin{equation*}
(R_{o}(s))^{1/2} \! =_{s \downarrow b_{j-1}^{o}} \! \mathcal{O} \! \left((s \!
- \! b_{j-1}^{o})^{1/2} \right) \qquad \text{and} \qquad (R_{o}(s))^{1/2} \!
=_{s \uparrow a_{j}^{o}} \! \mathcal{O} \! \left((a_{j}^{o} \! - \! s)^{1/2}
\right), \quad j \! = \! 1,\dotsc,N \! + \! 1,
\end{equation*}
which shows that all the integrals above constituting the $n$-dependent system
of $2(N \! + \! 1)$ moment conditions for the end-points of the support of
$\mu_{V}^{o}$ have removable singularities at $b_{j-1}^{o},a_{j}^{o}$, $j \!
= \! 1,\dotsc,N \! + \! 1$.

Recall {}from Equation~(3.4) that, for $z \! \in \! J_{o}$, $\mathscr{F}^{o}_{
\pm}(z) \! = \! -\tfrac{1}{\pi \mi z} \! - \! \mi (2 \! + \! \tfrac{1}{n})
(\mathcal{H} \psi_{V}^{o})(z) \! \mp \! (2 \! + \! \tfrac{1}{n}) \psi_{V}^{o}
(z)$: using the fact that, {}from Equation~(3.3), for $z \! \in \! J_{o}$,
$(\mathcal{H} \psi_{V}^{o})(z) \! = \! \tfrac{1}{2(2+\frac{1}{n}) \pi}
(\tfrac{2}{z} \! + \! \widetilde{V}^{\prime}(z))$, it follows that
\begin{equation*}
\mathscr{F}^{o}_{\pm}(z) \! = \! \dfrac{\widetilde{V}^{\prime}(z)}{2 \pi \mi}
\! \mp \! \left(2 \! + \! \tfrac{1}{n} \right) \! \psi_{V}^{o}(z), \quad z \!
\in \! J_{o}.
\end{equation*}
{}From Equation~(3.7), it follows that
\begin{equation*}
\mathscr{F}^{o}_{\pm}(z) \! = \! \dfrac{\widetilde{V}^{\prime}(z)}{2 \pi
\mi} \! + \! \dfrac{(R_{o}(z))^{1/2}_{\pm}}{2} \oint_{C_{\mathrm{R}}^{o}} \!
\dfrac{(\frac{2}{\pi \mi s} \! + \! \frac{\widetilde{V}^{\prime}(s)}{\mi \pi}
)}{(R_{o}(s))^{1/2}(s \! - \! z)} \, \dfrac{\md s}{2 \pi \mi};
\end{equation*}
thus, equating the above two expressions for $\mathscr{F}^{o}_{\pm}(z)$, one
arrives at $\psi_{V}^{o}(x) \! = \! \tfrac{1}{2 \pi \mi}(R_{o}(x))^{1/2}_{+}
h_{V}^{o}(x) \pmb{1}_{J_{o}}(x)$, where $h_{V}^{o}(z)$ is defined in the
Lemma, and $\pmb{1}_{J_{o}}(x)$ is the characteristic function of the set
$J_{o}$, which gives rise to the formula for the density of the `odd'
equilibrium measure, $\md \mu_{V}^{o}(x) \! = \! \psi_{V}^{o}(x) \, \md x$
(the integral representation of $h_{V}^{o}(z)$ shows that it is analytic in
some open subset of $\mathbb{C}^{\ast}$ containing $J_{o})$. Now, recalling
that $\widetilde{V} \colon \mathbb{R} \setminus \{0\} \! \to \! \mathbb{R}$
satisfying conditions~(2.3)--(2.5) is regular, and that, for $s \! \in \!
J_{o}$ (resp., $s \! \in \! \overline{J_{o}})$, $\psi_{V}^{o}(s) \! > \! 0$
(resp., $\psi_{V}^{o}(s) \! \geqslant \! 0)$ and $(R_{o}(s))^{1/2}_{+} \! =
\! \mi (\vert R_{o}(s) \vert)^{1/2} \! \in \! \mi \mathbb{R}_{\pm}$ (resp.,
$(R_{o}(s))^{1/2}_{+} \! = \! \mi (\vert R_{o}(s) \vert)^{1/2} \! \in \! \mi
\mathbb{R})$, it follows {}from the formula $\psi_{V}^{o}(s) \! = \! \tfrac{
1}{2 \pi \mi}(R_{o}(s))^{1/2}_{+}h_{V}^{o}(s) \pmb{1}_{J_{o}}(s)$ and the
regularity assumption, namely, $h_{V}^{o}(z) \! \not\equiv \! 0$ for $z \!
\in \! \overline{J_{o}}$, that $(\vert R_{o}(s) \vert)^{1/2}h_{V}^{o}(s) \! >
\! 0$, $s \! \in \! J_{o}$ (resp., $(\vert R_{o}(s) \vert)^{1/2}h_{V}^{o}(s)
\! \geqslant \! 0$, $s \! \in \! \overline{J_{o}})$.

Finally, it will be shown that, if $\overline{J_{o}} \! := \! \cup_{j=1}^{N
+1}[b_{j-1}^{o},a_{j}^{o}]$, the end-points of the support of the `odd' 
equilibrium measure, which satisfy the $n$-dependent system of $2(N \! + \! 
1)$ moment conditions stated in the Lemma, are (real) analytic functions of 
$z_{o}$, thus proving the (local) solvability of the $n$-dependent $2(N \! 
+ \! 1)$ moment conditions. Towards this end, one follows closely the idea 
of the proof of Theorem~1.3~(iii) in \cite{a81} (see, also, Section~8 of 
\cite{a44}, and \cite{a84}). Recall {}from Subsection~2.2 that $\widetilde{V}
(z) \! := \! z_{o}V(z)$, where $z_{o} \colon \mathbb{N} \times \mathbb{N} \! 
\to \! \mathbb{R}_{+}$, $(n,\mathscr{N}) \! \mapsto \! z_{o} \! := \! 
\mathscr{N}/n$, and, in the double-scaling limit as $\mathscr{N},n \! \to \! 
\infty$, $z_{o} \! = \! 1 \! + \! o(1)$. Furthermore, {}from the analysis 
above, it was shown that the end-points of the support of the `odd' 
equilibrium measure were the simple zeros/roots of the function $\mathfrak{
q}_{V}^{o}(z)$, that is (up to re-arrangement), $\lbrace b_{0}^{o},b_{1}^{o},
\dotsc,b_{N}^{o},a_{1}^{o},a_{2}^{o},\dotsc,a_{N+1}^{o} \rbrace \! = \! 
\lbrace \mathstrut x \! \in \! \mathbb{R}; \, \mathfrak{q}_{V}^{o}(x) \! = 
\! 0 \rbrace$ (these are the only roots for the regular case studied in this 
work). The function $\mathfrak{q}_{V}^{o}(x) \! \in \! \mathbb{R}(x)$ (the 
algebra of rational functions in $x$ with coefficients in $\mathbb{R})$ is 
real rational on $\mathbb{R}$ and real analytic on $\mathbb{R} \setminus 
\{0\}$, it has analytic continuation to $\lbrace \mathstrut z \! \in \! 
\mathbb{C}; \, \inf_{p \in \mathbb{R}} \vert z \! - \! p \vert \! < \! r \! 
\in \! (0,1) \rbrace \setminus \{0\}$ (independent of $z_{o})$, and depends 
continuously on $z_{o}$; thus, its simple zeros/roots, that is, $b_{k-1}^{o} 
\! = \! b_{k-1}^{o}(z_{o})$ and $a_{k}^{o} \! = \! a_{k}^{o}(z_{o})$, $k \! 
= \! 1,\dotsc,N \! + \! 1$, are continuous functions of $z_{o}$.

Write the large-$z$ (e.g., $\vert z \vert \! \gg \! \max_{j=1,\dotsc,N+1} 
\lbrace \vert b_{j-1}^{o} \! - \! a_{j}^{o} \vert \rbrace)$ asymptotic 
expansion for $\mathscr{F}^{o}(z)$ given above as follows:
\begin{equation*}
\mathscr{F}^{o}(z) \underset{z \to \infty}{=} -\dfrac{1}{\mi \pi z} \! - \!
\dfrac{(R_{o}(z))^{1/2}}{2 \pi \mi z} \sum_{j=0}^{\infty} \mathcal{T}_{j}^{o}
z^{-j},
\end{equation*}
where
\begin{equation*}
\mathcal{T}_{j}^{o} \! := \! \int_{J_{o}} \! \left(\dfrac{2}{\mi \pi s} \! +
\! \dfrac{\widetilde{V}^{\prime}(s)}{\mi \pi} \right) \! \dfrac{s^{j}}{(R_{o}
(s))^{1/2}_{+}} \, \md s, \quad j \! \in \! \mathbb{Z}_{0}^{+}.
\end{equation*}
Set, for $n \! \in \! \mathbb{N}$,
\begin{equation*}
\mathcal{N}_{j}^{o} \! := \! \int_{a_{j}^{o}}^{b_{j}^{o}} \! \left((\mathcal{
H} \psi_{V}^{o})(s) \! - \! \dfrac{1}{2(2 \! + \! \frac{1}{n}) \pi} \! \left(
\dfrac{2}{s} \! + \! \widetilde{V}^{\prime}(s) \right) \right) \md s, \quad j
\! = \! 1,\dotsc,N.
\end{equation*}
The $(n$-dependent) $2(N \! + \! 1)$ moment conditions are, thus, equivalent 
to $\mathcal{T}_{j}^{o} \! = \! 0$, $j \! = \! 0,\dotsc,N$, $\mathcal{T}_{N+
1}^{o} \! = \! -2(2 \! + \! \tfrac{1}{n})$, and $\mathcal{N}_{j}^{o} \! = 
\! 0$, $j \! = \! 1,\dotsc,N$. It will first be shown that, for regular 
$\widetilde{V} \colon \mathbb{R} \setminus \{0\} \! \to \! \mathbb{R}$ 
satisfying conditions~(2.3)--(2.5), the Jacobian of the transformation 
$\lbrace b_{0}^{o}(z_{o}),\dotsc,b_{N}^{o}(z_{o}),a_{1}^{o}(z_{o}),\dotsc,
a_{N+1}^{o}(z_{o}) \rbrace \! \mapsto \! \lbrace \mathcal{T}_{0}^{o},\dotsc,
\mathcal{T}_{N+1}^{o},\mathcal{N}_{1}^{o},\dotsc,\mathcal{N}_{N}^{o} \rbrace$, 
that is, $\operatorname{Jac}(\mathcal{T}_{0}^{o},\dotsc,\mathcal{T}_{N+1}^{
o},\mathcal{N}_{1}^{o},\dotsc,\mathcal{N}_{N}^{o}) \! := \! \tfrac{\partial 
(\mathcal{T}_{0}^{o},\dotsc,\mathcal{T}_{N+1}^{o},\mathcal{N}_{1}^{o},\dotsc,
\mathcal{N}_{N}^{o})}{\partial (b_{0}^{o},\dotsc,b_{N}^{o},a_{1}^{o},\dotsc,
a_{N+1}^{o})}$, is non-zero w\-h\-e\-n\-e\-v\-e\-r $b_{j-1}^{o} \! = \! 
b_{j-1}^{o}(z_{o})$ and $a_{j}^{o} \! = \! a_{j}^{o}(z_{o})$, $j \! = \! 1,
\dotsc,N \! + \! 1$, are chosen so that $\overline{J_{o}} \! = \! \cup_{j=
1}^{N+1}[b_{j-1}^{o},a_{j}^{o}]$. Using the equation $(\mathcal{H} \psi_{V}^{
o})(z) \! = \! (2 \! + \! \tfrac{1}{n})^{-1}(\mi \mathscr{F}^{o}(z) \! + \! 
\tfrac{1}{\pi z})$ (cf. Equation~(3.2)), one follows the analysis on 
pp.~778--779 of \cite{a81} (see, also, Section~3, Lemma~3.5, of \cite{a38}) 
to show that, for $k \! = \! 1,\dotsc,N \! + \! 1$:
\begin{gather}
\dfrac{\partial \mathcal{T}_{j}^{o}}{\partial b_{k-1}^{o}} \! = \! b_{k-1}^{o}
\dfrac{\partial \mathcal{T}_{j-1}^{o}}{\partial b_{k-1}^{o}} \! + \! \dfrac{
1}{2} \mathcal{T}_{j-1}^{o}, \quad j \! \in \! \mathbb{N}, \tag{T1} \\
\dfrac{\partial \mathcal{T}_{j}^{o}}{\partial a_{k}^{o}} \! = \! a_{k}^{o}
\dfrac{\partial \mathcal{T}_{j-1}^{o}}{\partial a_{k}^{o}} \! + \! \dfrac{1}{
2} \mathcal{T}_{j-1}^{o}, \quad j \! \in \! \mathbb{N}, \tag{T2} \\
\dfrac{\partial \mathscr{F}^{o}(z)}{\partial b_{k-1}^{o}} \! = \! -\dfrac{1}{2
\pi \mi} \! \left(\dfrac{\partial \mathcal{T}_{0}^{o}}{\partial b_{k-1}^{o}}
\right) \! \dfrac{(R_{o}(z))^{1/2}}{z \! - \! b_{k-1}^{o}}, \quad z \! \in \!
\mathbb{C} \setminus (\overline{J_{o}} \cup \{0\}), \tag{F1} \\
\dfrac{\partial \mathscr{F}^{o}(z)}{\partial a_{k}^{o}} \! = \! -\dfrac{1}{2
\pi \mi} \! \left(\dfrac{\partial \mathcal{T}_{0}^{o}}{\partial a_{k}^{o}}
\right) \! \dfrac{(R_{o}(z))^{1/2}}{z \! - \! a_{k}^{o}}, \quad z \! \in \!
\mathbb{C} \setminus (\overline{J_{o}} \cup \{0\}), \tag{F2} \\
\dfrac{\partial \mathcal{N}_{j}^{o}}{\partial b_{k-1}^{o}} \! = \! -\dfrac{1}{
2(2 \! + \! \frac{1}{n}) \pi} \! \left(\dfrac{\partial \mathcal{T}_{0}^{o}}{
\partial b_{k-1}^{o}} \right) \! \int_{a_{j}^{o}}^{b_{j}^{o}} \dfrac{(R_{o}
(s))^{1/2}}{s \! - \! b_{k-1}^{o}} \, \md s, \quad j \! = \! 1,\dotsc,N,
\tag{N1} \\
\dfrac{\partial \mathcal{N}_{j}^{o}}{\partial a_{k}^{o}} \! = \! -\dfrac{1}{2
(2 \! + \! \frac{1}{n}) \pi} \! \left(\dfrac{\partial \mathcal{T}_{0}^{o}}{
\partial a_{k}^{o}} \right) \! \int_{a_{j}^{o}}^{b_{j}^{o}} \dfrac{(R_{o}(s)
)^{1/2}}{s \! - \! a_{k}^{o}} \, \md s, \quad j \! = \! 1,\dotsc,N; \tag{N2}
\end{gather}
moreover, if one evaluates Equations~(T1) and~(T2) on the solution of the
$n$-dependent system of $2(N \! + \! 1)$ moment conditions, that is, $\mathcal{
T}_{j}^{o} \! = \! 0$, $j \! = \! 0,\dotsc,N$, $\mathcal{T}_{N+1}^{o} \! = \!
-2(2 \! + \! \tfrac{1}{n})$, and $\mathcal{N}_{i}^{o} \! = \! 0$, $i \! = \!
1,\dotsc,N$, one arrives at
\begin{equation}
\dfrac{\partial \mathcal{T}_{j}^{o}}{\partial b_{k-1}^{o}} \! = \! (b_{k-1}^{
o})^{j} \dfrac{\partial \mathcal{T}_{0}^{o}}{\partial b_{k-1}^{o}}, \qquad
\qquad \dfrac{\partial \mathcal{T}_{j}^{o}}{\partial a_{k}^{o}} \! = \!
(a_{k}^{o})^{j} \dfrac{\partial \mathcal{T}_{0}^{o}}{\partial a_{k}^{o}},
\quad j \! = \! 0,\dotsc,N \! + \! 1. \tag{S1}
\end{equation}
Via Equations~(N1), (N2), and~(S1), one now computes the Jacobian of the
transformation $\lbrace b_{0}^{o}(z_{o}),\dotsc,\linebreak[4]
b_{N}^{o}(z_{o}),a_{1}^{o}(z_{o}),\dotsc,a_{N+1}^{o}(z_{o}) \rbrace \!
\mapsto \! \lbrace \mathcal{T}_{0}^{o},\dotsc,\mathcal{T}_{N+1}^{o},\mathcal{
N}_{1}^{o},\dotsc,\mathcal{N}_{N}^{o} \rbrace$ on the solution of the
$n$-dependent system of $2(N \! + \! 1)$ moment conditions:
\begin{align*}
&\operatorname{Jac}(\mathcal{T}_{0}^{o},\dotsc,\mathcal{T}_{N+1}^{o},\mathcal{
N}_{1}^{o},\dotsc,\mathcal{N}_{N}^{o}) \! := \! \dfrac{\partial (\mathcal{T}_{
0}^{o},\dotsc,\mathcal{T}_{N+1}^{o},\mathcal{N}_{1}^{o},\dotsc,\mathcal{N}_{
N}^{o})}{\partial (b_{0}^{o},\dotsc,b_{N}^{o},a_{1}^{o},\dotsc,a_{N+1}^{o})} \\
&=
\begin{vmatrix}
\frac{\partial \mathcal{T}_{0}^{o}}{\partial b_{0}^{o}} & \frac{\partial
\mathcal{T}_{0}^{o}}{\partial b_{1}^{o}} & \cdots & \frac{\partial \mathcal{
T}_{0}^{o}}{\partial b_{N}^{o}} & \frac{\partial \mathcal{T}_{0}^{o}}{\partial
a_{1}^{o}} & \frac{\partial \mathcal{T}_{0}^{o}}{\partial a_{2}^{o}} & \cdots
& \frac{\partial \mathcal{T}_{0}^{o}}{\partial a_{N+1}^{o}} \\
\frac{\partial \mathcal{T}_{1}^{o}}{\partial b_{0}^{o}} & \frac{\partial
\mathcal{T}_{1}^{o}}{\partial b_{1}^{o}} & \cdots & \frac{\partial \mathcal{
T}_{1}^{o}}{\partial b_{N}^{o}} & \frac{\partial \mathcal{T}_{1}^{o}}{\partial
a_{1}^{o}} & \frac{\partial \mathcal{T}_{1}^{o}}{\partial a_{2}^{o}} & \cdots
& \frac{\partial \mathcal{T}_{1}^{o}}{\partial a_{N+1}^{o}} \\
\vdots & \vdots & \ddots & \vdots & \vdots & \vdots & \ddots & \vdots \\
\frac{\partial \mathcal{T}_{N+1}^{o}}{\partial b_{0}^{o}} & \frac{\partial
\mathcal{T}_{N+1}^{o}}{\partial b_{1}^{o}} & \cdots & \frac{\partial \mathcal{
T}_{N+1}^{o}}{\partial b_{N}^{o}} & \frac{\partial \mathcal{T}_{N+1}^{o}}{
\partial a_{1}^{o}} & \frac{\partial \mathcal{T}_{N+1}^{o}}{\partial a_{2}^{
o}} & \cdots & \frac{\partial \mathcal{T}_{N+1}^{o}}{\partial a_{N+1}^{o}} \\
\frac{\partial \mathcal{N}_{1}^{o}}{\partial b_{0}^{o}} & \frac{\partial
\mathcal{N}_{1}^{o}}{\partial b_{1}^{o}} & \cdots & \frac{\partial \mathcal{
N}_{1}^{o}}{\partial b_{N}^{o}} & \frac{\partial \mathcal{N}_{1}^{o}}{\partial
a_{1}^{o}} & \frac{\partial \mathcal{N}_{1}^{o}}{\partial a_{2}^{o}} & \cdots
& \frac{\partial \mathcal{N}_{1}^{o}}{\partial a_{N+1}^{o}} \\
\frac{\partial \mathcal{N}_{2}^{o}}{\partial b_{0}^{o}} & \frac{\partial
\mathcal{N}_{2}^{o}}{\partial b_{1}^{o}} & \cdots & \frac{\partial \mathcal{
N}_{2}^{o}}{\partial b_{N}^{o}} & \frac{\partial \mathcal{N}_{2}^{o}}{\partial
a_{1}^{o}} & \frac{\partial \mathcal{N}_{2}^{o}}{\partial a_{2}^{o}} & \cdots
& \frac{\partial \mathcal{N}_{2}^{o}}{\partial a_{N+1}^{o}} \\
\vdots & \vdots & \ddots & \vdots & \vdots & \vdots & \ddots & \vdots \\
\frac{\partial \mathcal{N}_{N}^{o}}{\partial b_{0}^{o}} & \frac{\partial
\mathcal{N}_{N}^{o}}{\partial b_{1}^{o}} & \cdots & \frac{\partial \mathcal{
N}_{N}^{o}}{\partial b_{N}^{o}} & \frac{\partial \mathcal{N}_{N}^{o}}{\partial
a_{1}^{o}} & \frac{\partial \mathcal{N}_{N}^{o}}{\partial a_{2}^{o}} & \cdots
& \frac{\partial \mathcal{N}_{N}^{o}}{\partial a_{N+1}^{o}}
\end{vmatrix} \\
&= \dfrac{(-1)^{N}}{(2(2 \! + \! \frac{1}{n}) \pi)^{N}} \! \left(\prod_{k=
1}^{N+1} \dfrac{\partial \mathcal{T}_{0}^{o}}{\partial b_{k-1}^{o}} \dfrac{
\partial \mathcal{T}_{0}^{o}}{\partial a_{k}^{o}} \right) \! \left(\prod_{j=
1}^{N} \int_{a_{j}^{o}}^{b_{j}^{o}}(R_{o}(s_{j}))^{1/2} \, \md s_{j} \right)
\! \Delta_{d}^{o},
\end{align*}
where
\begin{equation*}
\Delta_{d}^{o} \! := \!
\begin{vmatrix}
1 & 1 & \cdots & 1 & 1 & 1 & \cdots & 1 \\
b_{0}^{o} & b_{1}^{o} & \cdots & b_{N}^{o} & a_{1}^{o} & a_{2}^{o} & \cdots &
a_{N+1}^{o} \\
\vdots & \vdots & \ddots & \vdots & \vdots & \vdots & \ddots & \vdots \\
(b_{0}^{o})^{N+1} & (b_{1}^{o})^{N+1} & \cdots & (b_{N}^{o})^{N+1} & (a_{1}^{
o})^{N+1} & (a_{2}^{o})^{N+1} & \cdots & (a_{N+1}^{o})^{N+1} \\
\frac{1}{s_{1}-b_{0}^{o}} & \frac{1}{s_{1}-b_{1}^{o}} & \cdots & \frac{1}{s_{1}
-b_{N}^{o}} & \frac{1}{s_{1}-a_{1}^{o}} & \frac{1}{s_{1}-a_{2}^{o}} & \cdots &
\frac{1}{s_{1}-a_{N+1}^{o}} \\
\frac{1}{s_{2}-b_{0}^{o}} & \frac{1}{s_{2}-b_{1}^{o}} & \cdots & \frac{1}{s_{2}
-b_{N}^{o}} & \frac{1}{s_{2}-a_{1}^{o}} & \frac{1}{s_{2}-a_{2}^{o}} & \cdots &
\frac{1}{s_{2}-a_{N+1}^{o}} \\
\vdots & \vdots & \ddots & \vdots & \vdots & \vdots & \ddots & \vdots \\
\frac{1}{s_{N}-b_{0}^{o}} & \frac{1}{s_{N}-b_{1}^{o}} & \cdots & \frac{1}{s_{N}
-b_{N}^{o}} & \frac{1}{s_{N}-a_{1}^{o}} & \frac{1}{s_{N}-a_{2}^{o}} & \cdots &
\frac{1}{s_{N}-a_{N+1}^{o}}
\end{vmatrix}.
\end{equation*}
The above determinant, that is, $\Delta_{d}^{o}$, has been calculated on 
pg.~780 of \cite{a81} (see, also, Section~5.3, Equations~(5.148) and~(5.149) 
of \cite{a50}), namely,
\begin{equation*}
\Delta_{d}^{o} \! = \! \dfrac{\left(\prod_{j=1}^{N+1} \prod_{k=1}^{N+1}(b_{k-
1}^{o} \! - \! a_{j}^{o}) \right) \! \left(\prod_{\substack{j,k=1\\j<k}}^{N+
1}(a_{k}^{o} \! - \! a_{j}^{o})(b_{k-1}^{o} \! - \! b_{j-1}^{o}) \right) \!
\left(\prod_{\substack{j,k=1\\j<k}}^{N}(s_{k} \! - \! s_{j}) \right)}{(-1)^{
N} \prod_{j=1}^{N} \prod_{k=1}^{N+1}(s_{j} \! - \! a_{k}^{o})(s_{j} \! - \!
b_{k-1}^{o})};
\end{equation*}
but, for $-\infty \! < \! b_{0}^{o} \! < \! a_{1}^{o} \! < \! s_{1} \! < \!
b_{1}^{o} \! < \! a_{2}^{o} \! < \! s_{2} \! < \! b_{2}^{o} \! < \! \dotsb \!
< \! b_{N-1}^{o} \! < \! a_{N}^{o} \! < \! s_{N} \! < \! b_{N}^{o} \! < \!
a_{N+1}^{o} \! < \! +\infty$, $\Delta_{d}^{o} \! \not= \! 0$ (which means that
it is of a fixed sign), and $\int_{a_{j}^{o}}^{b_{j}^{o}}(R_{o}(s_{j}))^{1/2}
\, \md s_{j} \! > \! 0$, $j \! = \! 1,\dotsc,N$, whence
\begin{equation*}
\left(\prod_{j=1}^{N} \int_{a_{j}^{o}}^{b_{j}^{o}}(R_{o}(s_{j}))^{1/2} \, \md
s_{j} \right) \! \Delta_{d}^{o} \! \not= \! 0.
\end{equation*}
It remains to show that $\partial \mathcal{T}_{0}^{o}/\partial b_{k-1}^{o}$
and $\partial \mathcal{T}_{0}^{o}/\partial a_{k}^{o}$, $k \! = \! 1,\dotsc,N
\! + \! 1$, too, are non-zero; for this purpose, one exploits the fact that
$\mathcal{T}_{0}^{o} \! = \! (\mi \pi)^{-1} \int_{J_{o}}(2s^{-1} \! + \!
\widetilde{V}^{\prime}(s))(R_{o}(s))^{-1/2}_{+} \, \md s$ is independent of
$z$. It follows {}from Equation~(3.7), the integral representation for $h_{
V}^{o}(z)$ given in the Lemma, and Equations~(F1) and~(F2) that
\begin{gather*}
\dfrac{(z \! - \! b_{k-1}^{o})}{\sqrt{\smash[b]{R_{o}(z)}}} \dfrac{\partial
\mathscr{F}^{o}(z)}{\partial b_{k-1}^{o}} \! = \! -\dfrac{1}{2 \pi \mi} \!
\left(2 \! + \! \dfrac{1}{n} \right) \! \left((z \! - \! b_{k-1}^{o}) \dfrac{
\partial h_{V}^{o}(z)}{\partial b_{k-1}^{o}} \! - \! \dfrac{1}{2}h_{V}^{o}(z)
\right), \quad k \! = \! 1,\dotsc,N \! + \! 1, \\
\dfrac{(z \! - \! a_{k}^{o})}{\sqrt{\smash[b]{R_{o}(z)}}} \dfrac{\partial
\mathscr{F}^{o}(z)}{\partial a_{k}^{o}} \! = \! -\dfrac{1}{2 \pi \mi} \!
\left(2 \! + \! \dfrac{1}{n} \right) \! \left((z \! - \! a_{k}^{o}) \dfrac{
\partial h_{V}^{o}(z)}{\partial a_{k}^{o}} \! - \! \dfrac{1}{2}h_{V}^{o}(z)
\right), \quad k \! = \! 1,\dotsc,N \! + \! 1:
\end{gather*}
using, now, the $z$-independence of $\mathcal{T}_{0}^{o}$, and the fact that,
for the case of regular $\widetilde{V} \colon \mathbb{R} \setminus \{0\} \!
\to \! \mathbb{R}$ satisfying conditions~(2.3)--(2.5), $h_{V}^{o}(b_{j-1}^{
o}),h_{V}^{o}(a_{j}^{o}) \! \not= \! 0$, $j \! = \! 1,\dotsc,N \! + \! 1$,
one shows that
\begin{gather*}
\left. \dfrac{(z \! - \! b_{k-1}^{o})}{\sqrt{\smash[b]{R_{o}(z)}}} \dfrac{
\partial \mathscr{F}^{o}(z)}{\partial b_{k-1}^{o}} \right\vert_{z=b_{k-1}^{o}}
\! = \! \dfrac{1}{4 \pi \mi} \! \left(2 \! + \! \dfrac{1}{n} \right) \!
h_{V}^{o}(b_{k-1}^{o}) \! \not= \! 0, \quad k \! = \! 1,\dotsc,N \! + \! 1, \\
\left. \dfrac{(z \! - \! a_{k}^{o})}{\sqrt{\smash[b]{R_{o}(z)}}} \dfrac{
\partial \mathscr{F}^{o}(z)}{\partial a_{k}^{o}} \right\vert_{z=a_{k}^{o}} \!
= \! \dfrac{1}{4 \pi \mi} \! \left(2 \! + \! \dfrac{1}{n} \right) \! h_{V}^{o}
(a_{k}^{o}) \! \not= \! 0, \quad k \! = \! 1,\dotsc,N \! + \! 1;
\end{gather*}
thus, via Equations~(F1) and~(F2), one arrives at
\begin{equation*}
\dfrac{\partial \mathcal{T}_{0}^{o}}{\partial b_{k-1}^{o}} \! = \! -\dfrac{1}{
2} \! \left(2 \! + \! \dfrac{1}{n} \right) \! h_{V}^{o}(b_{k-1}^{o}) \! \not=
\! 0 \qquad \text{and} \qquad \dfrac{\partial \mathcal{T}_{0}^{o}}{\partial
a_{k}^{o}} \! = \! -\dfrac{1}{2} \! \left(2 \! + \! \dfrac{1}{n} \right) \!
h_{V}^{o}(a_{k}^{o}) \! \not= \! 0, \quad k \! = \! 1,\dotsc,N \! + \! 1,
\end{equation*}
whence
\begin{equation*}
\prod_{k=1}^{N+1} \dfrac{\partial \mathcal{T}_{0}^{o}}{\partial b_{k-1}^{o}}
\dfrac{\partial \mathcal{T}_{0}^{o}}{\partial a_{k}^{o}} \! = \! \left(\dfrac{
1}{2} \! \left(2 \! + \! \dfrac{1}{n} \right) \right)^{2(N+1)} \, \prod_{k=
1}^{N+1}h_{V}^{o}(b_{k-1}^{o})h_{V}^{o}(a_{k}^{o}) \! \not= \! 0.
\end{equation*}
Hence, $\operatorname{Jac}(\mathcal{T}_{0}^{o},\dotsc,\mathcal{T}_{N+1}^{o},
\mathcal{N}_{1}^{o},\dotsc,\mathcal{N}_{N}^{o}) \! \not= \! 0$.

It remains, still, to show that $\mathcal{T}_{j}^{o}$, $j \! = \! 0,\dotsc,
N \! + \! 1$, and $\mathcal{N}_{i}^{o}$, $i \! = \! 1,\dotsc,N$, are (real)
analytic functions of $\lbrace b_{j-1}^{o},a_{j}^{o} \rbrace_{j=1}^{N+1}$.
{}From the definition of $\mathcal{T}_{j}^{o}$, $j \! \in \! \mathbb{Z}_{0}^{
+}$, above, using the fact that they are independent of $z$, thus giving rise
to zero residue contributions, a straightforward residue calculus calculation
shows that, equivalently,
\begin{equation*}
\mathcal{T}_{j}^{o} \! = \! \dfrac{1}{2} \oint_{C_{\mathrm{R}}^{o}} \! \left(
\dfrac{2}{\mi \pi s} \! + \! \dfrac{\widetilde{V}^{\prime}(s)}{\mi \pi}
\right) \! \dfrac{s^{j}}{(R_{o}(s))^{1/2}} \, \md s, \quad j \! \in \!
\mathbb{Z}_{0}^{+},
\end{equation*}
where (the closed contour) $C_{\mathrm{R}}^{o}$ has been defined above: the
only factor depending on $\lbrace b_{k-1}^{o},a_{k}^{o} \rbrace_{k=1}^{N+1}$
is $\sqrt{\smash[b]{R_{o}(z)}}$. As $\sqrt{\smash[b]{R_{o}(z)}}$ is analytic
$\forall \, \, z \! \in \! \mathbb{C} \setminus \cup_{j=1}^{N+1}[b_{j-1}^{o},
a_{j}^{o}]$, and since $C_{\mathrm{R}}^{o} \subset \mathbb{C} \setminus \cup_{
j=1}^{N+1}[b_{j-1}^{o},a_{j}^{o}]$, with $\operatorname{int}(C_{\mathrm{R}}^{
o}) \supset \overline{J_{o}} \cup \{z\}$, it follows that, in particular,
$\sqrt{\smash[b]{R_{o}(z)}} \! \! \upharpoonright_{C_{\mathrm{R}}^{o}}$ is
an analytic function of $\lbrace b_{j-1}^{o},a_{j}^{o} \rbrace_{j=1}^{N+1}$,
which implies, via the above (equivalent) contour integral representation of
$\mathcal{T}_{j}^{o}$, $j \! \in \! \mathbb{Z}_{0}^{+}$, that $\mathcal{T}_{
k}^{o}$, $k \! = \! 0,\dotsc,N \! + \! 1$, are (real) analytic functions
of $\lbrace b_{j-1}^{o},a_{j}^{o} \rbrace_{j=1}^{N+1}$. Recalling that
$(\mathcal{H} \psi_{V}^{o})(z) \! = \! (2(2 \! + \! \tfrac{1}{n}) \pi)^{-1}
(2z^{-1} \! + \! \widetilde{V}^{\prime}(z)) \! - \! \tfrac{1}{2 \pi}(R_{o}
(z))^{1/2}h_{V}^{o}(z)$, it follows {}from the definition of $\mathcal{N}_{
j}^{o}$, $j \! = \! 1,\dotsc,N$, that
\begin{equation*}
\mathcal{N}_{j}^{o} \! = \! -\dfrac{1}{2 \pi} \int_{a_{j}^{o}}^{b_{j}^{o}}
(R_{o}(s))^{1/2}h_{V}^{o}(s) \, \md s, \quad j \! = \! 1,\dotsc,N:
\end{equation*}
making the linear change of variables $u_{j} \colon \mathbb{C} \! \to \!
\mathbb{C}$, $s \! \mapsto \! u_{j}(s) \! := \! (b_{j}^{o} \! - \! a_{j}^{o}
)^{-1}(s \! - \! a_{j}^{o})$, $j \! = \! 1,\dotsc,N$, which take each of the
(compact) intervals $[a_{j}^{o},b_{j}^{o}]$, $j \! = \! 1,\dotsc,N$, onto
$[0,1]$, and setting
\begin{equation*}
\sqrt{\smash[b]{\widehat{R}_{o}(z)}}:= \! \left(\prod_{k_{1}=1}^{j}(z \! - \!
b_{k_{1}-1}^{o}) \prod_{k_{2}=1}^{j-1}(z \! - \! a_{k_{2}}^{o}) \prod_{k_{3}=j
+1}^{N+1}(a_{k_{3}}^{o} \! - \! z) \prod_{k_{4}=j+2}^{N+1}(b_{k_{4}-1}^{o} \!
- \! z) \right)^{1/2},
\end{equation*}
one arrives at
\begin{equation*}
\mathcal{N}_{j}^{o} \! = \! -\dfrac{1}{2 \pi} \! \left(b_{j}^{o} \! - \! a_{
j}^{o} \right)^{2} \! \int_{0}^{1} \! \left(u_{j}(1 \! - \! u_{j}) \right)^{
1/2} \left(\widehat{R}_{o}((b_{j}^{o} \! - \! a_{j}^{o})u_{j} \! + \! a_{j}^{
o}) \right)^{1/2}h_{V}^{o}((b_{j}^{o} \! - \! a_{j}^{o})u_{j} \! + \! a_{j}^{
o}) \, \md u_{j}, \quad j \! = \! 1,\dotsc,N.
\end{equation*}
Recalling that $h_{V}^{o}(z)$ is analytic on $\mathbb{R} \setminus \{0\}$, in
particular, $h_{V}^{o}(b_{j-1}^{o}),h_{V}^{o}(a_{j}^{o}) \! \not= \! 0$, $j \!
= \! 1,\dotsc,N \! + \! 1$, and that it is an analytic function of $\lbrace
b_{k-1}^{o}(z_{o}),a_{k}^{o}(z_{o}) \rbrace_{k=1}^{N+1}$ (since $-\infty \!
< \! b_{0}^{o} \! < \! a_{1}^{o} \! < \! b_{1}^{o} \! < \! a_{2}^{o} \! < \!
\dotsb \! < \! b_{N}^{o} \! < \! a_{N+1}^{o} \! < \! +\infty)$, and noting
{}from the definition of $\sqrt{\smash[b]{\widehat{R}_{o}(z)}}$ above that,
it, too, is an analytic function of $(b_{j-1}^{o} \! - \! a_{j}^{o})u_{j} \!
+ \! a_{j}^{o}$, $(j,u_{j}) \! \in \! \{1,\dotsc,N\} \times [0,1]$, and thus
an analytic function of $\lbrace b_{j-1}^{o}(z_{o}),a_{j}^{o}(z_{o}) \rbrace_{
j=1}^{N+1}$, it follows that $\mathcal{N}_{j}^{o}$, $j \! = \! 1,\dotsc,N$,
are (real) analytic functions of $\lbrace b_{j-1}^{o}(z_{o}),a_{j}^{o}(z_{o})
\rbrace_{j=1}^{N+1}$.

Thus, as the Jacobian of the transformation $\lbrace b_{0}^{o}(z_{o}),\dotsc,
b_{N}^{o}(z_{o}),a_{1}^{o}(z_{o}),\dotsc,a_{N+1}^{o}(z_{o}) \rbrace \! \mapsto
\! \lbrace \mathcal{T}_{0}^{o},\dotsc,\linebreak[4]
\mathcal{T}_{N+1}^{o},\mathcal{N}_{1}^{o},\dotsc,\mathcal{N}_{N}^{o} \rbrace$
is non-zero whenever $\lbrace b_{j-1}^{o}(z_{o}),a_{j}^{o}(z_{o}) \rbrace_{j=
1}^{N+1}$, the end-points of the support of the `odd' equilibrium measure, are
chosen so that, for regular $\widetilde{V} \colon \mathbb{R} \setminus \{0\}
\! \to \! \mathbb{R}$ satisfying conditions~(2.3)--(2.5), $\overline{J_{o}}=
\! \cup_{j=1}^{N+1}[b_{j-1}^{o},a_{j}^{o}]$, and $\mathcal{T}_{j}^{o}$, $j \!
= \! 0,\dotsc,N \! + \! 1$, and $\mathcal{N}_{k}^{o}$, $k \! = \! 1,\dotsc,N$,
are (real) analytic functions of $\lbrace b_{j-1}^{o}(z_{o}),a_{j}^{o}(z_{o})
\rbrace_{j=1}^{N+1}$, it follows, via the Implicit Function Theorem, that
$b_{j-1}^{o}(z_{o}),a_{j}^{o}(z_{o})$, $j \! = \! 1,\dotsc,N \! + \! 1$, are
real analytic functions of $z_{o}$. \hfill $\qed$
\begin{eeeee}
It turns out that, for $\widetilde{V} \colon \mathbb{R} \setminus \{0\} \! \to
\! \mathbb{R}$ (satisfying conditions~(2.3)--(2.5)) of the form
\begin{equation*}
\widetilde{V}(z) \! = \! \sum^{2m_{2}}_{k=-2m_{1}} \widetilde{\varrho}_{k}
z^{k},
\end{equation*}
with $\widetilde{\varrho}_{k} \! \in \! \mathbb{R}$, $k \! = \! -2m_{1},
\dotsc,2m_{2}$, $m_{1,2} \! \in \! \mathbb{N}$, and (since $\widetilde{V}
(\pm \infty),\widetilde{V}(0) \! > \! 0)$ $\widetilde{\varrho}_{-2m_{1}},
\widetilde{\varrho}_{2m_{2}} \! > \! 0$, the integral for $h_{V}^{o}(z)$, that
is, $h_{V}^{o}(z) \! = \! \tfrac{1}{2}(2 \! + \! \tfrac{1}{n})^{-1} \oint_{
C_{\mathrm{R}}^{o}}(R_{o}(s))^{-1/2}(\frac{2 \mi}{\pi s} \! + \! \frac{\mi
\widetilde{V}^{\prime}(s)}{\pi})(s \! - \! z)^{-1} \, \md s$, can be evaluated
explicitly. Let $C_{\mathrm{R}}^{o} \! = \! \widetilde{\Gamma}_{\infty}^{o}
\cup \widetilde{\Gamma}_{0}^{o}$, where $\widetilde{\Gamma}_{\infty}^{o} \!
:= \! \lbrace \mathstrut z^{\prime} \! = \! R \me^{\mi \vartheta}, \, R \! >
\! 1/\varepsilon, \, \vartheta \! \in \! [0,2 \pi] \rbrace$ (oriented
clockwise), and $\widetilde{\Gamma}_{0}^{o} \! := \! \lbrace \mathstrut
z^{\prime} \! = \! r \me^{\mi \vartheta}, \, 0 \! < \! r \! < \! \varepsilon,
\, \vartheta \! \in \! [0,2 \pi] \rbrace$ (oriented counter-clockwise), with
$\varepsilon$ some arbitrarily fixed, sufficiently small positive real number
chosen such that: (i) $\partial \lbrace \mathstrut z^{\prime} \! \in \!
\mathbb{C}; \, \vert z^{\prime} \vert \! = \! \varepsilon \rbrace \cap
\partial \lbrace \mathstrut z^{\prime} \! \in \! \mathbb{C}; \, \vert
z^{\prime} \vert \! = \! 1/\varepsilon \rbrace \! = \! \varnothing$; (ii)
$\{\mathstrut z^{\prime} \! \in \! \mathbb{C}; \, \vert z^{\prime} \vert \! <
\! \varepsilon\} \cap (J_{o} \cup \{z\}) \! = \! \varnothing$; (iii) $\lbrace
\mathstrut z^{\prime} \! \in \! \mathbb{C}; \, \vert z^{\prime} \vert \! > \!
1/\varepsilon \rbrace \cap (J_{o} \cup \{z\}) \! = \! \varnothing$; and (iv)
$\lbrace \mathstrut z^{\prime} \! \in \! \mathbb{C}; \, \varepsilon \! < \!
\vert z^{\prime} \vert \! < \! 1/\varepsilon \rbrace \! \supset \! J_{o} \cup
\{z\}$. A tedious, but otherwise straightforward, residue calculus calculation
shows that
\begin{align*}
h_{V}^{o}(z) \! =& \, \dfrac{z^{2m_{2}-N-2}}{(2 \! + \! \frac{1}{n})} \, \,
\sum_{j=0}^{2m_{2}-N-2} \underset{\substack{0 \leqslant \vert k \vert + \vert
l \vert \leqslant 2m_{2}-j-N-2\\k_{i} \geqslant 0, \, \, l_{i} \geqslant 0,
\, \, i \in \{0,\dotsc,N\}}}{\sideset{}{'}{\sum}_{k_{0},\dotsc,k_{N}} \, \,
\sideset{}{'}{\sum}_{l_{0},\dotsc,l_{N}}}(2m_{2} \! - \! j) \widetilde{
\varrho}_{2m_{2}-j} \! \left(\prod_{p=0}^{N} \prod_{j_{p}=0}^{k_{p}-1} \!
\left(\dfrac{1}{2} \! + \! j_{p} \right) \! \right) \\
\times& \left(\prod_{q=0}^{N} \prod_{\widetilde{m}_{q}=0}^{l_{q}-1} \! \left(
\dfrac{1}{2} \! + \! \widetilde{m}_{q} \right) \! \right) \! \dfrac{\left(
\prod_{p^{\prime}=0}^{N}(b_{p^{\prime}}^{o})^{k_{p^{\prime}}} \right) \! \left(
\prod_{q^{\prime}=0}^{N}(a_{q^{\prime}+1}^{o})^{l_{q^{\prime}}} \right)}{\left(
\prod_{l^{\prime}=0}^{N}k_{l^{\prime}}! \right) \! \left(\prod_{\widetilde{m}^{
\prime}=0}^{N}l_{\widetilde{m}^{\prime}}! \right)} \, z^{-(j+\vert k \vert +
\vert l \vert)} \\
+& \, \dfrac{(-1)^{\mathcal{N}_{+}}(\prod_{k=1}^{N+1} \vert b_{k-1}^{o}a_{k}^{
o} \vert)^{-1/2}}{(2 \! + \! \frac{1}{n})z^{2m_{1}+1}} \, \, \sum_{j=-2m_{1}+
1}^{0} \underset{\substack{0 \leqslant \vert k \vert + \vert l \vert \leqslant
2m_{1}+j\\k_{i} \geqslant 0, \, \, l_{i} \geqslant 0, \, \, i \in \{0,\dotsc,
N\}}}{\sideset{}{''}{\sum}_{k_{0},\dotsc,k_{N}} \, \, \sideset{}{''}{\sum}_{
l_{0},\dotsc,l_{N}}}(-2m_{1} \! - \! j) \widetilde{\varrho}_{-2m_{1}-j} \\
\times& \left(\prod_{p=0}^{N} \prod_{j_{p}=0}^{k_{p}-1} \! \left(\dfrac{1}{2}
\! + \! j_{p} \right) \! \right) \! \left(\prod_{q=0}^{N} \prod_{\widetilde{
m}_{q}=0}^{l_{q}-1} \! \left(\dfrac{1}{2} \! + \! \widetilde{m}_{q} \right) \!
\right) \! \dfrac{\left(\prod_{p^{\prime}=0}^{N}(b_{p^{\prime}}^{o})^{k_{p^{
\prime}}} \right)^{-1} \! \left(\prod_{q^{\prime}=0}^{N}(a_{q^{\prime}+1}^{
o})^{l_{q^{\prime}}} \right)^{-1}}{\left(\prod_{l^{\prime}=0}^{N}k_{l^{\prime}
}! \right) \! \left(\prod_{\widetilde{m}^{\prime}=0}^{N}l_{\widetilde{m}^{
\prime}}! \right)} \\
\times& \, z^{\vert k \vert +\vert l \vert -j} \! + \! \dfrac{2(-1)^{\mathcal{
N}_{+}}(\prod_{k=1}^{N+1} \vert b_{k-1}^{o}a_{k}^{o} \vert)^{-1/2}}{(2 \! +
\! \frac{1}{n})z},
\end{align*}
where $\mathcal{N}_{+} \! \in \! \lbrace 0,\dotsc,N \! + \! 1 \rbrace$ is the
number of bands to the right of $z \! = \! 0$, $\vert k \vert \! := \! k_{0}
\! + \! k_{1} \! + \! \dotsb \! + \! k_{N}$ $(\geqslant \! 0)$, $\vert l \vert
\! := \! l_{0} \! + \! l_{1} \! + \! \dotsb \! + \! l_{N}$ $(\geqslant \! 0)$,
and the primes (resp., double primes) on the summations mean that all possible
sums over $\{k_{l}\}_{l=0}^{N}$ and $\{l_{k}\}_{k=0}^{N}$ must be taken for
which $0 \! \leqslant \! k_{0} \! + \! \cdots \! + \! k_{N} \! + \! l_{0} \! +
\! \cdots \! + \! l_{N} \! \leqslant \! 2m_{2} \! - \! j \! - \! N \! - \! 2$,
$j \! = \! 0,\dotsc,2m_{2} \! - \! N \! - \! 2$, $k_{i} \! \geqslant \! 0$,
$l_{i} \! \geqslant \! 0$, $i \! = \! 0,\dotsc,N$ (resp., $0 \! \leqslant \!
k_{0} \! + \! \cdots \! + \! k_{N} \! + \! l_{0} \! + \! \cdots \! + \! l_{N}
\! \leqslant \! 2m_{1} \! + \! j$, $j \! = \! -2m_{1} \! + \! 1,\dotsc,0$,
$k_{i} \! \geqslant \! 0$, $l_{i} \! \geqslant \! 0$, $i \! = \! 0,\dotsc,N)$.
It is important to note that all of the above sums are finite sums: any sums
for which the upper limit is less than the lower limit are defined to be zero,
and any products in which the upper limit is less than the lower limit are
defined to be one; for example, $\sum_{j=0}^{-1}(\ast) \! := \! 0$ and
$\prod_{j=0}^{-1}(\ast) \! := \! 1$.

It is interesting to note that one may derive explicit formulae for the
various moments of the `odd' equilibrium measure, that is, $\int_{\mathbb{R}}
s^{\pm m} \, \md \mu_{V}^{o}(s) \! = \! \int_{J_{o}}s^{\pm m} \psi_{V}^{o}(s)
\, \md s$, $m \! \in \! \mathbb{N}$, in terms of the external field and the
function $(R_{o}(z))^{1/2}$; without loss of generality, and for demonstrative
purposes only, consider, say, the following moments: $\int_{J_{o}}s^{\pm j}
\, \md \mu_{V}^{o}(s)$, $j \! = \! 1,2,3$ (the calculations below
straightforwardly generalise to $\int_{J_{o}}s^{\pm (k+3)} \, \md \mu_{V}^{o}
(s)$, $k \! \in \! \mathbb{N})$. Recall the following formulae for $\mathscr{
F}^{o}(z)$ given in Lemma~3.5:
\begin{gather*}
\mathscr{F}^{o}(z) \! = \! -\dfrac{1}{\pi \mi z} \! - \! \dfrac{1}{\pi \mi} \!
\left(2 \! + \! \dfrac{1}{n} \right) \! \int_{J_{o}} \! \dfrac{\md \mu_{V}^{o}
(s)}{s \! - \! z}, \quad z \! \in \! \mathbb{C} \setminus (J_{o} \cup \{0\}),
\\
\mathscr{F}^{o}(z) \! = \! -\dfrac{1}{\pi \mi z} \! - \! (R_{o}(z))^{1/2}
\int_{J_{o}} \! \dfrac{(\frac{2 \mi}{\pi s} \! + \! \frac{\mi \widetilde{V}^{
\prime}(s)}{\pi})}{(R_{o}(s))^{1/2}_{+}(s \! - \! z)} \, \dfrac{\md s}{2 \pi
\mi}, \quad z \! \in \! \mathbb{C} \setminus (J_{o} \cup \{0\}).
\end{gather*}
One derives the following asymptotic expansions: (1) for $\mu_{V}^{o} \! \in
\! \mathcal{M}_{1}(\mathbb{R})$, in particular, $\int_{J_{o}}s^{-m} \, \md
\mu_{V}^{o}(s) \! < \! \infty$, $m \! \in \! \mathbb{N}$, $s \! \in \! J_{o}$
and $z \! \notin \! J_{o}$, with $\vert z/s \vert \! \ll \! 1$ (e.g., $\vert z
\vert \! \ll \! \min_{j=1,\dotsc,N+1} \lbrace \vert b_{j-1}^{o} \! - \! a_{
j}^{o} \vert \rbrace)$, via the expansions $\tfrac{1}{z-s} \! = \! -\sum_{k=
0}^{l} \tfrac{z^{k}}{s^{k+1}} \! + \! \tfrac{z^{l+1}}{s^{l+1}(z-s)}$, $l \!
\in \! \mathbb{Z}_{0}^{+}$, and $\ln (1 \! - \! \ast) \! = \! -\sum_{k=1}^{
\infty} \tfrac{\ast^{k}}{k}$, $\vert \ast \vert \! \ll \! 1$,
\begin{align*}
\mathscr{F}^{o}(z) \underset{z \to 0}{=}& \, -\dfrac{1}{\pi \mi z} \! - \!
\dfrac{1}{\pi \mi} \! \left(2 \! + \! \dfrac{1}{n} \right) \! \int_{J_{o}}
s^{-1} \, \md \mu_{V}^{o}(s) \! + \! z \! \left(-\dfrac{1}{\pi \mi} \! \left(
2 \! + \! \dfrac{1}{n} \right) \! \int_{J_{o}}s^{-2} \, \md \mu_{V}^{o}(s)
\right) \\
+& \, z^{2} \! \left(-\dfrac{1}{\pi \mi} \! \left(2 \! + \! \dfrac{1}{n}
\right) \! \int_{J_{o}}s^{-3} \, \md \mu_{V}^{o}(s) \right) \! + \!
\mathcal{O} \! \left(z^{3} \right),
\end{align*}
and
\begin{equation*}
\mathscr{F}^{o}(z) \underset{z \to 0}{=} \! -\dfrac{1}{\pi \mi z} \! + \!
\gamma_{V}^{o} \! \left(\check{Q}_{0}^{o} \! + \! z(\check{Q}_{1}^{o} \! - \!
\check{P}_{0}^{o} \check{Q}_{0}^{o}) \! + \! z^{2}(\check{Q}_{2}^{o} \! - \!
\check{P}_{0}^{o} \check{Q}_{1}^{o} \! + \! \check{P}_{1}^{o} \check{Q}_{0}^{
o}) \! + \! \mathcal{O}(z^{3}) \right),
\end{equation*}
where
\begin{gather*}
\gamma_{V}^{o} \! := \! (-1)^{\mathcal{N}_{+}} \! \left(\prod_{j=1}^{N+1}
\left\vert b_{j-1}^{o}a_{j}^{o} \right\vert \right)^{1/2}, \qquad \qquad
\check{P}_{0}^{o} \! := \! \dfrac{1}{2} \sum_{j=1}^{N+1} \! \left(\dfrac{1}{
b_{j-1}^{o}} \! + \! \dfrac{1}{a_{j}^{o}} \right), \\
\check{P}_{1}^{o} \! := \! \dfrac{1}{2}(\check{P}_{0}^{o})^{2} \! - \! \dfrac{
1}{4} \sum_{j=1}^{N+1} \! \left(\dfrac{1}{(b_{j-1}^{o})^{2}} \! + \! \dfrac{
1}{(a_{j}^{o})^{2}} \right), \qquad \check{Q}_{j}^{o} \! := \! -\int_{J_{o}}
\! \dfrac{(\frac{2 \mi}{\pi s} \! + \! \frac{\mi \widetilde{V}^{\prime}(s)}{
\pi})}{(R_{o}(s))^{1/2}_{+}s^{j+1}} \, \dfrac{\md s}{2 \pi \mi}, \quad j \! =
\! 0,1,2;
\end{gather*}
and (2) for $\mu_{V}^{o} \! \in \! \mathcal{M}_{1}(\mathbb{R})$, in
particular, $\int_{J_{o}} \md \mu_{V}^{o}(s) \! = \! 1$ and $\int_{J_{o}}s^{m}
\, \md \mu_{V}^{o}(s) \! < \! \infty$, $m \! \in \! \mathbb{N}$, $s \! \in \!
J_{o}$ and $z \! \notin \! J_{o}$, with $\vert s/z \vert \! \ll \! 1$ (e.g.,
$\vert z \vert \! \gg \! \max_{j=1,\dotsc,N+1} \lbrace \vert b_{j-1}^{o} \!
- \! a_{j}^{o} \vert \rbrace)$, via the expansions $\tfrac{1}{s-z} \! = \! -
\sum_{k=0}^{l} \tfrac{s^{k}}{z^{k+1}} \! + \! \tfrac{s^{l+1}}{z^{l+1}(s-z)}$,
$l \! \in \! \mathbb{Z}_{0}^{+}$, and $\ln (1 \! - \! \ast) \! = \! -\sum_{k=
1}^{\infty} \tfrac{\ast^{k}}{k}$, $\vert \ast \vert \! \ll \! 1$,
\begin{align*}
\mathscr{F}^{o}(z) \underset{z \to \infty}{=}& \, \dfrac{1}{\pi \mi} \! \left(
1 \! + \! \dfrac{1}{n} \right) \! \dfrac{1}{z} \! + \! \dfrac{1}{z^{2}} \!
\left(\dfrac{1}{\pi \mi} \! \left(2 \! + \! \dfrac{1}{n} \right) \! \int_{J_{
o}}s \, \md \mu_{V}^{o}(s) \right) \! + \! \dfrac{1}{z^{3}} \! \left(\dfrac{
1}{\pi \mi} \! \left(2 \! + \! \dfrac{1}{n} \right) \! \int_{J_{o}}s^{2} \,
\md \mu_{V}^{o}(s) \right) \\
+& \, \dfrac{1}{z^{4}} \! \left(\dfrac{1}{\pi \mi} \! \left(2 \! + \! \dfrac{
1}{n} \right) \! \int_{J_{o}}s^{3} \, \md \mu_{V}^{o}(s) \right) \! + \!
\mathcal{O} \! \left(\dfrac{1}{z^{5}} \right),
\end{align*}
and
\begin{equation*}
\mathscr{F}^{o}(z) \underset{z \to \infty}{=} -\dfrac{1}{\pi \mi z} \! + \!
z^{N} \! \left(1 \! + \! \dfrac{\check{\alpha}}{z} \! + \! \dfrac{\widetilde{
P}_{0}^{o}}{z^{2}} \! + \! \dfrac{\widetilde{P}_{1}^{o}}{z^{3}} \! + \! \dotsb
\right) \! \int_{J_{o}} \! \dfrac{(\frac{2 \mi}{\pi s} \! + \! \frac{\mi
\widetilde{V}^{\prime}(s)}{\pi})}{(R_{o}(s))^{1/2}_{+}} \! \left(1 \! + \!
\cdots \! + \! \dfrac{s^{N}}{z^{N}} \! + \! \dfrac{s^{N+1}}{z^{N+1}} \! + \!
\dotsb \right) \! \dfrac{\md s}{2 \pi \mi},
\end{equation*}
where
\begin{gather*}
\check{\alpha} \! := \! -\dfrac{1}{2} \sum_{j=1}^{N+1} \! \left(b_{j-1}^{
o} \! + \! a_{j}^{o} \right), \qquad \qquad \widetilde{P}_{0}^{o} \! := \!
\dfrac{1}{2}(\check{\alpha})^{2} \! - \! \dfrac{1}{4} \sum_{j=1}^{N+1} \!
\left((b_{j-1}^{o})^{2} \! + \! (a_{j}^{o})^{2} \right), \\
\widetilde{P}_{1}^{o} \! := \! -\dfrac{1}{3!} \sum_{j=1}^{N+1} \! \left((b_{j
-1}^{o})^{3} \! + \! (a_{j}^{o})^{3} \right) \! + \! \dfrac{(\check{\alpha})^{
3}}{3!} \! - \! \dfrac{\check{\alpha}}{4} \sum_{j=1}^{N+1} \! \left((b_{j-1}^{
o})^{2} \! + \! (a_{j}^{o})^{2} \right).
\end{gather*}
Recalling the following $(n$-dependent) $N \! + \! 2$ moment conditions stated
in Lemma~3.5,
\begin{equation*}
\int_{J_{o}} \dfrac{(\frac{2 \mi}{\pi s} \! + \! \frac{\mi \widetilde{V}^{
\prime}(s)}{\pi})s^{j}}{(R_{o}(s))^{1/2}_{+}} \, \md s \! = \! 0, \quad j \!
= \! 0,\dotsc,N, \qquad \text{and} \qquad \int_{J_{o}} \dfrac{(\frac{2 \mi}{
\pi s} \! + \! \frac{\mi \widetilde{V}^{\prime}(s)}{\pi})s^{N+1}}{(R_{o}(s)
)^{1/2}_{+}} \, \md s \! = \! 2 \! \left(2 \! + \! \dfrac{1}{n} \right),
\end{equation*}
and equating the respective pairs of asymptotic expansions above (as $z \! \to
\! 0$ and $z \! \to \! \infty)$ for $\mathscr{F}^{o}(z)$, one arrives at the
following expressions for the first three `positive' and `negative' moments
of the `odd' equilibrium measure:
\begin{align*}
\int_{J_{o}}s \, \md \mu_{V}^{o}(s)=& \, \dfrac{1}{2(2 \! + \! \frac{1}{n})}
\! \left(\int_{J_{o}} \dfrac{(\frac{2 \mi}{\pi s} \! + \! \frac{\mi \widetilde{
V}^{\prime}(s)}{\pi})s^{N+2}}{(R_{o}(s))^{1/2}_{+}} \, \md s \! - \! \left(2
\! + \! \dfrac{1}{n} \right) \! \sum_{j=1}^{N+1}(b_{j-1}^{o} \! + \!
a_{j}^{o}) \right) \! , \\
\int_{J_{o}}s^{2} \, \md \mu_{V}^{o}(s)=& \, \dfrac{1}{2(2 \! + \! \frac{1}{
n})} \! \left(\int_{J_{o}} \dfrac{(\frac{2 \mi}{\pi s} \! + \! \frac{\mi
\widetilde{V}^{\prime}(s)}{\pi})s^{N+3}}{(R_{o}(s))^{1/2}_{+}} \, \md s \!
- \! \dfrac{1}{2} \! \left(\sum_{j=1}^{N+1}(b_{j-1}^{o} \! + \! a_{j}^{o})
\right) \! \int_{J_{o}} \dfrac{(\frac{2 \mi}{\pi s} \! + \! \frac{\mi
\widetilde{V}^{\prime}(s)}{\pi})s^{N+2}}{(R_{o}(s))^{1/2}_{+}} \, \md s
\right. \\
+&\left. \, \dfrac{1}{2} \! \left(2 \! + \! \dfrac{1}{n} \right) \! \left(
\dfrac{1}{2} \! \left(\sum_{j=1}^{N+1}(b_{j-1}^{o} \! + \! a_{j}^{o}) \right)^{
2} \! - \! \sum_{j=1}^{N+1}((b_{j-1}^{o})^{2} \! + \! (a_{j}^{o})^{2}) \right)
\right), \\
\int_{J_{o}}s^{3} \, \md \mu_{V}^{o}(s)=& \, \dfrac{1}{2(2 \! + \! \frac{1}{n}
)} \! \left(\int_{J_{o}} \dfrac{(\frac{2 \mi}{\pi s} \! + \! \frac{\mi
\widetilde{V}^{\prime}(s)}{\pi})s^{N+4}}{(R_{o}(s))^{1/2}_{+}} \, \md s \!
- \! \dfrac{1}{2} \! \left(\sum_{j=1}^{N+1}(b_{j-1}^{o} \! + \! a_{j}^{o})
\right) \! \int_{J_{o}} \dfrac{(\frac{2 \mi}{\pi s} \! + \! \frac{\mi
\widetilde{V}^{\prime}(s)}{\pi})s^{N+3}}{(R_{o}(s))^{1/2}_{+}} \, \md s
\right. \\
+&\left. \, \dfrac{1}{4} \! \left(\dfrac{1}{2} \! \left(\sum_{j=1}^{N+1}(b_{j
-1}^{o} \! + \! a_{j}^{o}) \right)^{2} \! - \! \sum_{j=1}^{N+1}((b_{j-1}^{o})^{
2} \! + \! (a_{j}^{o})^{2}) \right) \! \int_{J_{o}} \dfrac{(\frac{2 \mi}{\pi
s} \! + \! \frac{\mi \widetilde{V}^{\prime}(s)}{\pi})s^{N+2}}{(R_{o}(s))^{
1/2}_{+}} \, \md s \right. \\
-&\left. \, \dfrac{1}{4} \! \left(2 \! + \! \dfrac{1}{n} \right) \! \left(
\dfrac{1}{3!} \! \left(\sum_{j=1}^{N+1}(b_{j-1}^{o} \! + \! a_{j}^{o})
\right)^{3} \! + \! \dfrac{4}{3} \sum_{j=1}^{N+1}((b_{j-1}^{o})^{3} \! + \!
(a_{j}^{o})^{3}) \! - \! \sum_{j=1}^{N+1}(b_{j-1}^{o} \! + \! a_{j}^{o})
\right. \right. \\
\times& \left. \left. \sum_{k=1}^{N+1}((b_{k-1}^{o})^{2} \! + \! (a_{k}^{o})^{
2}) \right) \right), \\
\int_{J_{o}}s^{-1} \, \md \mu_{V}^{o}(s)=& \, \dfrac{(-1)^{\mathcal{N}_{+}}
\left(\prod_{j=1}^{N+1} \vert b_{j-1}^{o}a_{j}^{o} \vert \right)^{1/2}}{2(2 \!
+ \! \frac{1}{n})} \int_{J_{o}} \dfrac{(\frac{2 \mi}{\pi s} \! + \! \frac{\mi
\widetilde{V}^{\prime}(s)}{\pi})}{s(R_{o}(s))^{1/2}_{+}} \, \md s, \\
\int_{J_{o}}s^{-2} \, \md \mu_{V}^{o}(s)=& \, \dfrac{(-1)^{\mathcal{N}_{+}}
\left(\prod_{j=1}^{N+1} \vert b_{j-1}^{o}a_{j}^{o} \vert \right)^{1/2}}{2(2
\! + \! \frac{1}{n})} \left(\int_{J_{o}} \dfrac{(\frac{2 \mi}{\pi s} \! + \!
\frac{\mi \widetilde{V}^{\prime}(s)}{\pi})}{s^{2}(R_{o}(s))^{1/2}_{+}} \, \md
s \! - \! \dfrac{1}{2} \! \left(\sum_{j=1}^{N+1} \! \left(\dfrac{1}{b_{j-1}^{
o}} \! + \! \dfrac{1}{a_{j}^{o}} \right) \! \right) \right. \\
\times& \left. \int_{J_{o}} \dfrac{(\frac{2 \mi}{\pi s} \! + \! \frac{\mi
\widetilde{V}^{\prime}(s)}{\pi})}{s(R_{o}(s))^{1/2}_{+}} \, \md s \right), \\
\int_{J_{o}}s^{-3} \, \md \mu_{V}^{o}(s)=& \, \dfrac{(-1)^{\mathcal{N}_{+}}
\left(\prod_{j=1}^{N+1} \vert b_{j-1}^{o}a_{j}^{o} \vert \right)^{1/2}}{2(2
\! + \! \frac{1}{n})} \left(\int_{J_{o}} \dfrac{(\frac{2 \mi}{\pi s} \! + \!
\frac{\mi \widetilde{V}^{\prime}(s)}{\pi})}{s^{3}(R_{o}(s))^{1/2}_{+}} \, \md
s \! - \! \dfrac{1}{2} \! \left(\sum_{j=1}^{N+1} \! \left(\dfrac{1}{b_{j-1}^{
o}} \! + \! \dfrac{1}{a_{j}^{o}} \right) \! \right) \right. \\
\times& \left. \int_{J_{o}} \dfrac{(\frac{2 \mi}{\pi s} \! + \! \frac{\mi
\widetilde{V}^{\prime}(s)}{\pi})}{s^{2}(R_{o}(s))^{1/2}_{+}} \, \md s \! + \!
\left(\dfrac{1}{8} \! \left(\sum_{j=1}^{N+1} \! \left(\dfrac{1}{b_{j-1}^{o}}
\! + \! \dfrac{1}{a_{j}^{o}} \right) \right)^{2} \! - \! \dfrac{1}{4} \sum_{
j=1}^{N+1} \! \left(\dfrac{1}{(b_{j-1}^{o})^{2}} \! + \! \dfrac{1}{(a_{j}^{o}
)^{2}} \right) \right) \right. \\
\times& \left. \int_{J_{o}} \dfrac{(\frac{2 \mi}{\pi s} \! + \! \frac{\mi
\widetilde{V}^{\prime}(s)}{\pi})}{s(R_{o}(s))^{1/2}_{+}} \, \md s \right).
\end{align*}
It is important to note that all of the above integrals are real valued
(since, for $s \! \in \! \overline{J_{o}}$, $(R_{o}(s))^{1/2}_{+} \! = \! \mi
(\vert R_{o}(s) \vert)^{1/2} \! \in \! \mi \mathbb{R})$ and bounded (since,
for $j \! = \! 1,\dotsc,N \! + \! 1$, $(R_{o}(s))^{1/2} \! =_{s \downarrow
b_{j-1}^{o}} \! \mathcal{O}((s \! - \! b_{j-1}^{o})^{1/2})$ and $(R_{o}(s))^{
1/2} \! =_{s \uparrow a_{j}^{o}} \! \mathcal{O}((a_{j}^{o} \! - \! s)^{1/2})$,
that is, there are removable singularities at the end-points of the support of
the `odd' equilibrium measure). \hfill $\blacksquare$
\end{eeeee}
\begin{ccccc}
Let the external field $\widetilde{V} \colon \mathbb{R} \setminus \{0\} \!
\to \! \mathbb{R}$ satisfy conditions~{\rm (2.3)--(2.5)}. Let the `odd'
equilibrium measure, $\mu_{V}^{o}$, and its support, $\operatorname{supp}
(\mu_{V}^{o}) \! =: \! J_{o}$ $(\subset \overline{\mathbb{R}} \setminus
\lbrace 0,\pm \infty \rbrace)$, be as described in Lemma~{\rm 3.5}, and let
there exist $\ell_{o}$ $(\in \! \mathbb{R})$, the `odd' variational constant,
such that
\begin{equation}
\begin{aligned}
2 \! \left(2 \! + \! \dfrac{1}{n} \right) \! \int_{J_{o}} \ln (\vert x \! -
\! s \vert) \psi_{V}^{o}(s) \, \md s \! - \! 2 \ln \vert x \vert \! - \!
\widetilde{V}(x) \! - \! \ell_{o} \! - \! 2 \! \left(2 \! + \! \dfrac{1}{n}
\right) \! Q_{o} =& \, 0, \quad x \! \in \! \overline{J_{o}}, \\
2 \! \left(2 \! + \! \dfrac{1}{n} \right) \! \int_{J_{o}} \ln (\vert x \! -
\! s \vert) \psi_{V}^{o}(s) \, \md s \! - \! 2 \ln \vert x \vert \! - \!
\widetilde{V}(x) \! - \! \ell_{o} \! - \! 2 \! \left(2 \! + \! \dfrac{1}{n}
\right) \! Q_{o} \leqslant& \, 0, \quad x \! \in \! \mathbb{R} \setminus
\overline{J_{o}},
\end{aligned}
\end{equation}
where
\begin{equation*}
Q_{o} \! := \! \int_{J_{o}} \ln (\lvert s \rvert) \psi_{V}^{o}(s) \, \md s,
\end{equation*}
and, for $\widetilde{V}$ regular, the inequality in the second of
Equations~{\rm (3.9)} is strict. Then:
\begin{compactenum}
\item[{\rm (1)}] $g^{o}_{+}(z) \! + \! g^{o}_{-}(z) \! - \! \widetilde{V}(z)
\! - \! \ell_{o} \! - \! (\mathfrak{Q}^{+}_{\mathscr{A}} \! + \! \mathfrak{
Q}^{-}_{\mathscr{A}}) \! = \! 0$, $z \! \in \! \overline{J_{o}}$, where $g^{
o}_{\pm}(z) \! := \! \lim_{\varepsilon \downarrow 0}g^{o}(z \! \pm \! \mi
\varepsilon);$
\item[{\rm (2)}] $g^{o}_{+}(z) \! + \! g^{o}_{-}(z) \! - \! \widetilde{V}(z)
\! - \! \ell_{o} \! - \! (\mathfrak{Q}^{+}_{\mathscr{A}} \! + \! \mathfrak{
Q}^{-}_{\mathscr{A}}) \! \leqslant \! 0$, $z \! \in \! \mathbb{R} \setminus
\overline{J_{o}}$, where equality holds for at most a finite number of points,
and, for $\widetilde{V}$ regular, the inequality is strict;
\item[{\rm (3)}] $g^{o}_{+}(z) \! - \! g^{o}_{-}(z) \! - \! \mathfrak{Q}^{+}_{
\mathscr{A}} \! + \! \mathfrak{Q}^{-}_{\mathscr{A}} \! \in \! \mi f_{g^{o}}^{
\mathbb{R}}$, $z \! \in \! \mathbb{R}$, where $f_{g^{o}}^{\mathbb{R}} \colon
\mathbb{R} \! \to \! \mathbb{R}$ is some bounded function, and, in particular,
$g^{o}_{+}(z) \! - \! g^{o}_{-}(z) \! - \! \mathfrak{Q}^{+}_{\mathscr{A}} \!
+ \! \mathfrak{Q}^{-}_{\mathscr{A}} \! = \! \mi \operatorname{const.}$, $z \!
\in \! \mathbb{R} \setminus \overline{J_{o}}$, where $\operatorname{const.} \!
\in \! \mathbb{R};$
\item[{\rm (4)}] $\mi (g^{o}_{+}(z) \! - \! g^{o}_{-}(z) \! - \! \mathfrak{Q}^{
+}_{\mathscr{A}} \! + \! \mathfrak{Q}^{-}_{\mathscr{A}})^{\prime} \! \geqslant
\! 0$, $z \! \in \! J_{o}$, and where, for $\widetilde{V}$ regular, equality
holds for at most a finite number of points.
\end{compactenum}
\end{ccccc}

\emph{Proof.} Set (cf. Lemma~3.5) $J_{o} \! := \! \cup_{j=1}^{N+1}J_{j}^{o}$,
where $J_{j}^{o} \! = \! (b_{j-1}^{o},a_{j}^{o}) \! =$ the $j$th `band', 
with $N \! \in \! \mathbb{N}$ and finite, $b_{0}^{o} \! := \! \min \lbrace 
\operatorname{supp}(\mu_{V}^{o}) \rbrace \! \notin \! \lbrace -\infty,0 
\rbrace$, $a_{N+1}^{o} \! := \! \max \lbrace \operatorname{supp}(\mu_{V}^{o}) 
\rbrace \! \notin \! \lbrace 0,+\infty \rbrace$, and $-\infty \! < \! b_{0}^{
o} \! < \! a_{1}^{o} \! < \! b_{1}^{o} \! < \! a_{2}^{o} \! < \! \cdots \! < 
\! b_{N}^{o} \! < \! a_{N+1}^{o} \! < \! +\infty$, and $\lbrace b_{j-1}^{o},
a_{j}^{o} \rbrace_{j=1}^{N+1}$ satisfy the $n$-dependent and (locally) 
solvable system of $2(N \! + \! 1)$ moment conditions given in Lemma~3.5. 
Consider the following cases: \textbf{(1)} $z \! \in \! \overline{J_{j}^{o}}
:= \! [b_{j-1}^{o},a_{j}^{o}]$, $j \! = \! 1,\dotsc,N \! + \! 1$; \textbf{(2)} 
$z \! \in \! (a_{j}^{o},b_{j}^{o}) \! =$ the $j$th `gap', $j \! = \! 1,\dotsc,
N$; \textbf{(3)} $z\! \in \! (a_{N+1}^{o},+\infty)$; and \textbf{(4)} $z \! 
\in \! (-\infty,b_{0}^{o})$.

$\mathbf{(1)}$ Recall the definition of $g^{o}(z)$ given in Lemma~3.4, namely,
$g^{o}(z) \! := \! \int_{J_{o}} \ln ((z \! - \! s)^{2+\frac{1}{n}}(zs)^{-1})
\psi_{V}^{o} \linebreak[4]
(s) \, \md s$, $z \! \in \! \mathbb{C} \setminus (-\infty,\max \lbrace 0,a_{N
+1}^{o} \rbrace)$, where the representation (cf. Lemma~3.5) $\md \mu_{V}^{o}
(s) \! = \! \psi_{V}^{o}(s) \, \md s$, $s \! \in \! J_{o}$, was substituted
into the latter. For $z \! \in \! \overline{J_{j}^{o}}$, $j \! = \! 1,\dotsc,
N \! + \! 1$, one shows that
\begin{align*}
g^{o}_{\pm}(z)=& \left(2 \! + \! \dfrac{1}{n} \right) \! \int_{J_{o}} \ln
(\vert z \! - \! s \vert) \psi_{V}^{o}(s) \, \md s \! \pm \! \mi \pi \! \left(
2 \! + \! \dfrac{1}{n} \right) \! \int_{z}^{a_{N+1}^{o}} \psi_{V}^{o}(s) \,
\md s \! - \! \int_{J_{o}} \ln (\lvert s \rvert) \psi_{V}^{o}(s) \, \md s \\
-& \, \mi \pi \int_{J_{o} \cap \mathbb{R}_{-}} \psi_{V}^{o}(s) \, \md s \! -
\!
\begin{cases}
\ln \vert z \vert, &\text{$z \! > \! 0,$} \\
\ln \vert z \vert \! \pm \! \mi \pi, &\text{$z \! < \! 0$,}
\end{cases}
\end{align*}
where $g^{o}_{\pm}(z) \! := \! \lim_{\varepsilon \downarrow 0}g^{o}(z \! \pm
\! \mi \varepsilon)$, whence
\begin{align*}
g^{o}_{+}(z) \! - \! g^{o}_{-}(z) \! - \! \mathfrak{Q}^{+}_{\mathscr{A}} \!
+ \! \mathfrak{Q}^{-}_{\mathscr{A}}=& \, 2 \! \left(2 \! + \! \dfrac{1}{n}
\right) \! \pi \mi \int_{z}^{a_{N+1}^{o}} \psi_{V}^{o}(s) \, \md s \! - \! 2
\! \left(2 \! + \! \dfrac{1}{n} \right) \! \pi \mi \int_{J_{o} \cap \mathbb{
R}_{+}} \psi_{V}^{o}(s) \, \md s \\
+& \,
\begin{cases}
0, &\text{$z \! > \! 0$,} \\
-2 \pi \mi, &\text{$z \! < \! 0$,}
\end{cases}
\end{align*}
which shows that $g^{o}_{+}(z) \! - \! g^{o}_{-}(z) \! - \! \mathfrak{Q}^{+}_{
\mathscr{A}} \! + \! \mathfrak{Q}^{-}_{\mathscr{A}} \! \in \! \mi \mathbb{R}$,
and $\Re (g^{o}_{+}(z) \! - \! g^{o}_{-}(z) \! - \! \mathfrak{Q}^{+}_{\mathscr{
A}} \! + \! \mathfrak{Q}^{-}_{\mathscr{A}}) \! = \! 0$; moreover, using the
Fundamental Theorem of Calculus, one shows that $(g^{o}_{+}(z) \! - \! g^{o}_{
-}(z) \! - \! \mathfrak{Q}^{+}_{\mathscr{A}} \! + \! \mathfrak{Q}^{-}_{
\mathscr{A}})^{\prime} \! = \! -2(2 \! + \! \tfrac{1}{n}) \pi \mi \psi_{V}^{o}
(z)$, whence $\mi (g^{o}_{+}(z) \! - \! g^{o}_{-}(z) \! - \! \mathfrak{Q}^{+}_{
\mathscr{A}} \! + \! \mathfrak{Q}^{-}_{\mathscr{A}})^{\prime} \! = \! 2(2 \! +
\! \tfrac{1}{n}) \pi \psi_{V}^{o}(z) \! \geqslant \! 0$, since $\psi_{V}^{o}
(z) \! \geqslant \! 0 \, \, \forall \, \, z \! \in \! \overline{J_{o}}$
$(\supset \overline{J_{j}^{o}}$, $j \! = \! 1,\dotsc,N \! + \! 1)$.
Furthermore, using the first of Equations~(3.9), one shows that
\begin{align*}
g^{o}_{+}(z) \! + \! g^{o}_{-}(z) \! - \! \widetilde{V}(z) \! - \! \ell_{o} \!
- \! (\mathfrak{Q}^{+}_{\mathscr{A}} \! + \! \mathfrak{Q}^{-}_{\mathscr{A}})
=& \, 2 \! \left(2 \! + \! \dfrac{1}{n} \right) \! \int_{J_{o}} \ln (\vert z
\! - \! s \vert) \psi_{V}^{o}(s) \, \md s \! - \! 2 \ln \vert z \vert \! - \!
\widetilde{V}(z) \\
-& \, \ell_{o} \! - \! 2 \! \left(2 \! + \! \dfrac{1}{n} \right) \! Q_{o} \!
= \! 0,
\end{align*}
which gives the formula for the ($n$-dependent) `odd' variational constant 
$\ell_{o}$ $(\in \! \mathbb{R})$, which is the same \cite{a81,a85} (see, also, 
Section~7 of \cite{a44}) for each compact interval $\overline{J_{j}^{o}}$, $j
\! = \! 1,\dotsc,N \! + \! 1$; in particular,
\begin{align*}
\ell_{o}=& \, \dfrac{1}{\pi} \! \left(2 \! + \! \dfrac{1}{n} \right) \! \sum_{
j=1}^{N+1} \int_{b_{j-1}^{o}}^{a_{j}^{o}} \ln \! \left(\left\vert \! \left(
\tfrac{1}{2}(b_{N}^{o} \! + \! a_{N+1}^{o})-s \right) \! s^{-1} \right\vert
\right) \! \left(\vert R_{o}(s) \vert \right)^{1/2} \! h_{V}^{o}(s) \, \md
s \! - \! 2 \ln \! \left\vert \tfrac{1}{2}(b_{N}^{o} \! + \! a_{N+1}^{o})
\right\vert \\
-& \, \widetilde{V} \! \left(\tfrac{1}{2}(b_{N}^{o} \! + \! a_{N+1}^{o})
\right),
\end{align*}
where $(\vert R_{o}(s) \vert)^{1/2}h_{V}^{o}(s) \! \geqslant \! 0$, $j \! = \!
1,\dotsc,N \! + \! 1$, and where there are no singularities in the integrand,
since, for (any) $r \! > \! 0$, $\lim_{\vert x \vert \to 0} \vert x \vert^{r}
\ln \vert x \vert \! = \! 0$.

$\mathbf{(2)}$ For $z \! \in \! (a_{j}^{o},b_{j}^{o})$, $j \! = \! 1,\dotsc,
N$, one shows that
\begin{align*}
g^{o}_{\pm}(z) =& \left(2 \! + \! \dfrac{1}{n} \right) \! \int_{J_{o}} \ln
(\vert z \! - \! s \vert) \psi_{V}^{o}(s) \, \md s \! \pm \! \left(2 \! + \!
\dfrac{1}{n} \right) \! \pi \mi \sum_{k=j+1}^{N+1} \int_{b_{k-1}^{o}}^{a_{k}^{
o}} \psi_{V}^{o}(s) \, \md s \! - \! \int_{J_{o}} \ln (\lvert s \rvert) \psi_{
V}^{o}(s) \, \md s \\
-& \, \mi \pi \int_{J_{o} \cap \mathbb{R}_{-}} \psi_{V}^{o}(s) \, \md s \! -
\!
\begin{cases}
\ln \vert z \vert, &\text{$z \! > \! 0,$} \\
\ln \vert z \vert \! \pm \! \mi \pi, &\text{$z \! < \! 0$,}
\end{cases}
\end{align*}
whence
\begin{align*}
g^{o}_{+}(z) \! - \! g^{o}_{-}(z) \! - \! \mathfrak{Q}^{+}_{\mathscr{A}} \! +
\! \mathfrak{Q}^{-}_{\mathscr{A}}=& \, 2 \! \left(2 \! + \! \dfrac{1}{n}
\right) \! \pi \mi \int_{b_{j}^{o}}^{a_{N+1}^{o}} \psi_{V}^{o}(s) \, \md s \!
- \! 2 \! \left(2 \! + \! \dfrac{1}{n} \right) \! \pi \mi \int_{J_{o} \cap
\mathbb{R}_{+}} \psi_{V}^{o}(s) \, \md s \\
+& \,
\begin{cases}
0, &\text{$z \! > \! 0$,} \\
-2 \pi \mi, &\text{$z \! < \! 0$,}
\end{cases}
\end{align*}
which shows that $g^{o}_{+}(z) \! - \! g^{o}_{-}(z) \! - \! \mathfrak{Q}^{+}_{
\mathscr{A}} \! + \! \mathfrak{Q}^{-}_{\mathscr{A}} \! = \! \mi \operatorname{
const.}$, with $\operatorname{const.} \! \in \! \mathbb{R}$, and $\Re (g^{o}_{
+}(z) \! - \! g^{o}_{-}(z) \! - \! \mathfrak{Q}^{+}_{\mathscr{A}} \! + \!
\mathfrak{Q}^{-}_{\mathscr{A}}) \! = \! 0$; moreover, $\mi (g^{o}_{+}(z) \! -
\! g^{o}_{-}(z) \! - \! \mathfrak{Q}^{+}_{\mathscr{A}} \! + \! \mathfrak{Q}^{
-}_{\mathscr{A}})^{\prime} \! = \! 0$. One notes {}from the above formulae
for $g^{o}_{\pm}(z)$ that
\begin{align*}
g^{o}_{+}(z) \! + \! g^{o}_{-}(z) \! - \! \widetilde{V}(z) \! - \! \ell_{o} \!
- \! (\mathfrak{Q}^{+}_{\mathscr{A}} \! + \! \mathfrak{Q}^{-}_{\mathscr{A}})=&
\, 2 \! \left(2 \! + \! \dfrac{1}{n} \right) \! \int_{J_{o}} \ln (\vert z \!
- \! s \vert) \psi_{V}^{o}(s) \, \md s \! - \! 2 \ln \vert z \vert \! - \!
\widetilde{V}(z) \\
-& \, \ell_{o} \! - \! 2 \! \left(2 \! + \! \dfrac{1}{n} \right) \! Q_{o}.
\end{align*}
Recalling that (cf. Lemma~3.5) $\mathcal{H} \colon \mathcal{L}^{2}_{\mathrm{
M}_{2}(\mathbb{C})} \! \to \! \mathcal{L}^{2}_{\mathrm{M}_{2}(\mathbb{C})}$,
$f \! \mapsto \! (\mathcal{H}f)(z) \! := \! \pvi \, \tfrac{f(s)}{z-s} \,
\tfrac{\md s}{\pi}$, where $\pvi$ denotes the principal value integral, one
shows that, for $z \! \in \! (a_{j}^{o},b_{j}^{o})$, $j \! = \! 1,\dotsc,N$,
\begin{equation*}
2 \! \left(2 \! + \! \dfrac{1}{n} \right) \! \int_{J_{o}} \ln (\vert z \! - \!
s \vert) \psi_{V}^{o}(s) \, \md s \! = \! 2 \! \left(2 \! + \! \dfrac{1}{n}
\right) \! \pi \int_{a_{j}^{o}}^{z}(\mathcal{H} \psi_{V}^{o})(s) \, \md s \! +
\! 2 \! \left(2 \! + \! \dfrac{1}{n} \right) \! \int_{J_{o}} \ln (\vert a_{
j}^{o} \! - \! s \vert) \psi_{V}^{o}(s) \, \md s;
\end{equation*}
thus,
\begin{align*}
g^{o}_{+}(z) \! + \! g^{o}_{-}(z) \! - \! \widetilde{V}(z) \! - \! \ell_{o} \!
- \! (\mathfrak{Q}^{+}_{\mathscr{A}} \! + \! \mathfrak{Q}^{-}_{\mathscr{A}})
=& \, 2(2 \! + \! \tfrac{1}{n}) \pi \int_{a_{j}^{o}}^{z}(\mathcal{H} \psi_{V}^{
o})(s) \, \md s \! + \! 2(2 \! + \! \tfrac{1}{n}) \int_{J_{o}} \ln (\vert a_{
j}^{o} \! - \! s \vert) \psi_{V}^{o}(s) \, \md s \\
-& \, 2 \ln \vert z \vert \! - \! \widetilde{V}(z) \! - \! \ell_{o} \! - \! 2
(2 \! + \! \tfrac{1}{n})Q_{o} \\
=& \, 2(2 \! + \! \tfrac{1}{n}) \int_{J_{o}} \ln (\vert a_{j}^{o} \! - \! s
\vert) \psi_{V}^{o}(s) \, \md s \! + \! 2(2 \! + \! \tfrac{1}{n}) \pi \int_{
a_{j}^{o}}^{z}(\mathcal{H} \psi_{V}^{o})(s) \, \md s \\
-& \, \int_{a_{j}^{o}}^{z} \widetilde{V}^{\prime}(s) \, \md s \! - \! 2 \int_{
a_{j}^{o}}^{z} \dfrac{1}{s} \, \md s \! - \! 2 \ln \vert a_{j}^{o} \vert \! -
\! \widetilde{V}(a_{j}^{o}) \! - \! \ell_{o} \! - \! 2(2 \! + \! \tfrac{1}{n})
Q_{o} \\
=& \, \int_{a_{j}^{o}}^{z} \! \left(2 \! \left(2 \! + \! \dfrac{1}{n} \right)
\! \pi (\mathcal{H} \psi_{V}^{o})(s) \! - \! \widetilde{V}^{\prime}(s) \! - \!
\dfrac{2}{s} \right) \! \md s,
\end{align*}
since
\begin{equation*}
2 \! \left(2 \! + \! \dfrac{1}{n} \right) \! \int_{J_{o}} \ln (\vert a_{j}^{o}
\! - \! s \vert) \psi_{V}^{o}(s) \, \md s \! - \! 2 \ln \vert a_{j}^{o} \vert
\! - \! \widetilde{V}(a_{j}^{o}) \! - \! \ell_{o} \! - \! 2 \! \left(2 \! + \!
\dfrac{1}{n} \right) \! Q_{o} \! = \! 0,
\end{equation*}
whence, for $j \! = \! 1,\dotsc,N \! + \! 1$,
\begin{equation*}
g^{o}_{+}(z) \! + \! g^{o}_{-}(z) \! - \! \widetilde{V}(z) \! - \! \ell_{o} \!
- \! (\mathfrak{Q}^{+}_{\mathscr{A}} \! + \! \mathfrak{Q}^{-}_{\mathscr{A}})
\! = \! \int_{a_{j}^{o}}^{z} \! \left(2 \! \left(2 \! + \! \dfrac{1}{n}
\right) \! \pi (\mathcal{H} \psi_{V}^{o})(s) \! - \! \widetilde{V}^{\prime}(s)
\! - \! \dfrac{2}{s} \right) \! \md s, \quad z \! \in \! (a_{j}^{o},b_{j}^{o}).
\end{equation*}
It was shown in the proof of Lemma~3.5 that $2(2 \! + \! \tfrac{1}{n}) \pi
(\mathcal{H} \psi_{V}^{o})(s) \! = \! \widetilde{V}^{\prime}(s) \! + \!
\tfrac{2}{s} \! - \! (2 \! + \! \tfrac{1}{n})(R_{o}(s))^{1/2}h_{V}^{o}(s)$,
$s \! \in \! (a_{j}^{o},b_{j}^{o})$, $j \! = \! 1,\dotsc,N$, whence
\begin{equation*}
g^{o}_{+}(z) \! + \! g^{o}_{-}(z) \! - \! \widetilde{V}(z) \! - \! \ell_{o} \!
- \! (\mathfrak{Q}^{+}_{\mathscr{A}} \! + \! \mathfrak{Q}^{-}_{\mathscr{A}})
\! = \! - \! \left(2 \! + \! \dfrac{1}{n} \right) \! \int_{a_{j}^{o}}^{z}(R_{o}
(s))^{1/2}h_{V}^{o}(s) \, \md s \! < \! 0, \quad z \! \in \! \cup_{j=1}^{N}
(a_{j}^{o},b_{j}^{o}):
\end{equation*}
since $h_{V}^{o}(z)$ is real analytic on $\mathbb{R} \setminus \{0\}$ and
$(R_{o}(s))^{1/2}h_{V}^{o}(s) \! > \! 0 \, \, \forall \, \, s \! \in \!
\cup_{j=1}^{N}(a_{j}^{o},b_{j}^{o})$, it follows that one has equality only
at points $z \! \in \! \cup_{j=1}^{N}(a_{j}^{o},b_{j}^{o})$ for which $h_{
V}^{o}(z) \! = \! 0$, which are finitely denumerable. (Note that, for $z \!
\in \! \cup_{j=1}^{N}(a_{j}^{o},b_{j}^{o})$, $(R_{o}(s))^{1/2}_{+} \! = \!
(R_{o}(s))^{1/2}_{-} \! = \! (R_{o}(s))^{1/2}$.)

$\mathbf{(3)}$ For $z \! \in \! (a_{N+1}^{o},+\infty)$, one shows that
\begin{align*}
g^{o}_{\pm}(z)=& \left(2 \! + \! \dfrac{1}{n} \right) \! \int_{J_{o}} \ln
(\vert z \! - \! s \vert) \psi_{V}^{o}(s) \, \md s \! - \! \int_{J_{o}} \ln
(\lvert s \rvert) \psi_{V}^{o}(s) \, \md s \! - \! \mi \pi \int_{J_{o} \cap
\mathbb{R}_{-}} \psi_{V}^{o}(s) \, \md s \\
-& \,
\begin{cases}
\ln \vert z \vert, &\text{$z \! > \! 0,$} \\
\ln \vert z \vert \! \pm \! \mi \pi, &\text{$z \! < \! 0$,}
\end{cases}
\end{align*}
whence
\begin{equation*}
g^{o}_{+}(z) \! - \! g^{o}_{-}(z) \! - \! \mathfrak{Q}^{+}_{\mathscr{A}} \! +
\! \mathfrak{Q}^{-}_{\mathscr{A}} \! = \! -2 \! \left(2 \! + \! \dfrac{1}{n}
\right) \! \pi \mi \int_{J_{o} \cap \mathbb{R}_{+}} \psi_{V}^{o}(s) \, \md s
\! + \!
\begin{cases}
0, &\text{$z \! > \! 0$,} \\
-2 \pi \mi, &\text{$z \! < \! 0$,}
\end{cases}
\end{equation*}
which shows that $g^{o}_{+}(z) \! - \! g^{o}_{-}(z) \! - \! \mathfrak{Q}^{+}_{
\mathscr{A}} \! + \! \mathfrak{Q}^{-}_{\mathscr{A}}$ is pure imaginary, and
$\mi (g^{o}_{+}(z) \! - \! g^{o}_{-}(z) \! - \! \mathfrak{Q}^{+}_{\mathscr{A}}
\! + \! \mathfrak{Q}^{-}_{\mathscr{A}})^{\prime} \! = \! 0$. Also, one shows
that
\begin{align*}
g^{o}_{+}(z) \! + \! g^{o}_{-}(z) \! - \! \widetilde{V}(z) \! - \! \ell_{o} \!
- \! (\mathfrak{Q}^{+}_{\mathscr{A}} \! + \! \mathfrak{Q}^{-}_{\mathscr{A}})=&
\, 2 \! \left(2 \! + \! \dfrac{1}{n} \right) \! \int_{J_{o}} \ln (\vert z \!
- \! s \vert) \psi_{V}^{o}(s) \, \md s \! - \! 2 \ln \vert z \vert \! - \!
\widetilde{V}(z) \\
-& \, \ell_{o} \! - \! 2 \! \left(2 \! + \! \dfrac{1}{n} \right) \! Q_{o};
\end{align*}
and, following the analysis of case~$\mathbf{(2)}$ above, one shows that, for
$z \! \in \! (a_{N+1}^{o},+\infty)$,
\begin{align*}
2 \! \left(2 \! + \! \dfrac{1}{n} \right) \! \int_{J_{o}} \ln (\vert z \! -
\! s \vert) \psi_{V}^{o}(s) \, \md s \! - \! 2 \ln \vert z \vert \! - \!
\widetilde{V}(z) \! - \! \ell_{o} \! - \! 2 \! \left(2 \! + \! \dfrac{1}{n}
\right) \! Q_{o}=& \, \int_{a_{N+1}^{o}}^{z} \! \left(2 \! \left(2 \! + \!
\dfrac{1}{n} \right) \! \pi (\mathcal{H} \psi_{V}^{o})(s) \right. \\
-&\left. \, \widetilde{V}^{\prime}(s) \! - \! \dfrac{2}{s} \right) \md s,
\end{align*}
thus, via the relation (cf. case~$\mathbf{(2)}$ above) $2(2 \! + \! \tfrac{1}{
n}) \pi (\mathcal{H} \psi_{V}^{o})(s) \! = \! \widetilde{V}^{\prime}(s) \! +
\! \tfrac{2}{s} \! - \! (2 \! + \! \tfrac{1}{n})(R_{o}(s))^{1/2}h_{V}^{o}(s)$,
$s \! \in \! (a_{N+1}^{o},+\infty)$, one arrives at
\begin{equation*}
g^{o}_{+}(z) \! + \! g^{o}_{-}(z) \! - \! \widetilde{V}(z) \! - \! \ell_{o} \!
- \! (\mathfrak{Q}^{+}_{\mathscr{A}} \! + \! \mathfrak{Q}^{-}_{\mathscr{A}})
\! = \! - \! \left(2 \! + \! \dfrac{1}{n} \right) \! \int_{a_{N+1}^{o}}^{z}
(R_{o}(s))^{1/2}h_{V}^{o}(s) \, \md s \! < \! 0, \quad z \! \in \! (a_{N+1}^{
o},+\infty).
\end{equation*}
If: (1) $z \! \to \! +\infty$ (e.g., $\vert z \vert \! \gg \! \max_{j=1,
\dotsc,N+1} \lbrace \vert b_{j-1}^{o},a_{j}^{o} \vert \rbrace)$, $s \! \in \!
J_{o}$, and $\vert s/z \vert \! \ll \! 1$, {}from $\mu_{V}^{o} \! \in \!
\mathcal{M}_{1}(\mathbb{R})$, in particular, $\int_{J_{o}} \md \mu_{V}^{o}(s)
\! = \! 1$ and $\int_{J_{o}}s^{m} \, \md \mu_{V}^{o}(s) \! < \! \infty$, $m \!
\in \! \mathbb{N}$, the formula for $g^{o}_{+}(z) \! + \! g^{o}_{-}(z) \! - \!
\widetilde{V}(z) \! - \! \ell_{o} \! - \! (\mathfrak{Q}^{+}_{\mathscr{A}} \! +
\! \mathfrak{Q}^{-}_{\mathscr{A}})$ above, and the expansions $\tfrac{1}{s-z}
\! = \! -\sum_{k=0}^{l} \tfrac{s^{k}}{z^{k+1}} \! + \! \tfrac{s^{l+1}}{z^{l+1}
(s-z)}$, $l \! \in \! \mathbb{Z}_{0}^{+}$, and $\ln (z \! - \! s) \! =_{\vert
z \vert \to \infty} \! \ln (z) \! - \! \sum_{k=1}^{\infty} \tfrac{1}{k}
(\tfrac{s}{z})^{k}$, one shows that
\begin{equation*}
g^{o}_{+}(z) \! + \! g^{o}_{-}(z) \! - \! \widetilde{V}(z) \! - \! \ell_{o} \!
- \! (\mathfrak{Q}^{+}_{\mathscr{A}} \! + \! \mathfrak{Q}^{-}_{\mathscr{A}})
\underset{z \to +\infty}{=} \left(1 \! + \! \dfrac{1}{n} \right) \! \ln (z^{
2} \! + \! 1) \! - \! \widetilde{V}(z) \! + \! \mathcal{O}(1),
\end{equation*}
which, upon recalling that (cf. condition~(2.4)) $\lim_{\vert x \vert \to
\infty}(\widetilde{V}(x)/\ln (x^{2} \! + \! 1)) \! = \! +\infty$, shows that
$g^{o}_{+}(z) \! + \! g^{o}_{-}(z) \! - \! \widetilde{V}(z) \! - \! \ell_{o}
\! - \! (\mathfrak{Q}^{+}_{\mathscr{A}} \! + \! \mathfrak{Q}^{-}_{\mathscr{A}
}) \! < \! 0$; and (2) $\vert z \vert \! \to \! 0$ (e.g., $\vert z \vert \!
\ll \! \min_{j=1,\dotsc,N+1} \lbrace \vert b_{j-1}^{o},a_{j}^{o} \vert
\rbrace)$, $s \! \in \! J_{o}$, and $\vert z/s \vert \! \ll \! 1$, {}from
$\mu_{V}^{o} \! \in \! \mathcal{M}_{1}(\mathbb{R})$, in particular, $\int_{
J_{o}}s^{-m} \, \md \mu_{V}^{o}(s) \! < \! \infty$, $m \! \in \! \mathbb{N}$,
the above formula for $g^{o}_{+}(z) \! + \! g^{o}_{-}(z) \! - \! \widetilde{
V}(z) \! - \! \ell_{o} \! - \! (\mathfrak{Q}^{+}_{\mathscr{A}} \! + \!
\mathfrak{Q}^{-}_{\mathscr{A}})$, and the expansions $\tfrac{1}{z-s} \! = \!
-\sum_{k=0}^{l} \tfrac{z^{k}}{s^{k+1}} \! + \! \tfrac{z^{l+1}}{s^{l+1}(z-s)}$,
$l \! \in \! \mathbb{Z}_{0}^{+}$, and $\ln (s \! - \! z) \! =_{\vert z \vert
\to 0} \! \ln (s) \! - \! \sum_{k=1}^{\infty} \tfrac{1}{k}(\tfrac{z}{s})^{k}$,
one shows that
\begin{equation*}
g^{o}_{+}(z) \! + \! g^{o}_{-}(z) \! - \! \widetilde{V}(z) \! - \! \ell_{o} \!
- \! (\mathfrak{Q}^{+}_{\mathscr{A}} \! + \! \mathfrak{Q}^{-}_{\mathscr{A}})
\underset{\vert z \vert \to 0}{=} \ln (z^{-2} \! + \! 1) \! - \! \widetilde{V}
(z) \! + \! \mathcal{O}(1),
\end{equation*}
whereupon, recalling that (cf. condition~(2.5)) $\lim_{\vert x \vert \to 0}
(\widetilde{V}(x)/\ln (x^{-2} \! + \! 1)) \! = \! +\infty$, it follows that
$g^{o}_{+}(z) \! + \! g^{o}_{-}(z) \! - \! \widetilde{V}(z) \! - \! \ell_{o}
\! - \! (\mathfrak{Q}^{+}_{\mathscr{A}} \! + \! \mathfrak{Q}^{-}_{\mathscr{A}}
) \! < \! 0$.

$\mathbf{(4)}$ For $z \! \in \! (-\infty,b_{0}^{o})$, one shows that
\begin{align*}
g^{o}_{\pm}(z)=& \left(2 \! + \! \dfrac{1}{n} \right) \! \int_{J_{o}} \ln
(\vert z \! - \! s \vert) \psi_{V}^{o}(s) \, \md s \! \pm \! \left(2 \! + \!
\dfrac{1}{n} \right) \! \pi \mi \! - \! \int_{J_{o}} \ln (\lvert s \rvert)
\psi_{V}^{o}(s) \, \md s \\
-& \, \mi \pi \int_{J_{o} \cap \mathbb{R}_{-}} \psi_{V}^{o}(s) \, \md s \! -
\!
\begin{cases}
\ln \vert z \vert, &\text{$z \! > \! 0,$} \\
\ln \vert z \vert \! \pm \! \mi \pi, &\text{$z \! < \! 0$,}
\end{cases}
\end{align*}
whence
\begin{equation*}
g^{o}_{+}(z) \! - \! g^{o}_{-}(z) \! - \! \mathfrak{Q}^{+}_{\mathscr{A}} \! +
\! \mathfrak{Q}^{-}_{\mathscr{A}} \! = \! 2 \! \left(2 \! + \! \dfrac{1}{n}
\right) \! \pi \mi \! - \! 2 \! \left(2 \! + \! \dfrac{1}{n} \right) \! \pi
\mi \int_{J_{o} \cap \mathbb{R}_{+}} \psi_{V}^{o}(s) \, \md s \! + \!
\begin{cases}
0, &\text{$z \! > \! 0$,} \\
-2 \pi \mi, &\text{$z \! < \! 0$,}
\end{cases}
\end{equation*}
which shows that $g^{o}_{+}(z) \! - \! g^{o}_{-}(z) \! - \! \mathfrak{Q}^{+}_{
\mathscr{A}} \! + \! \mathfrak{Q}^{-}_{\mathscr{A}}$ is pure imaginary, and
$\mi (g^{o}_{+}(z) \! - \! g^{o}_{-}(z) \! - \! \mathfrak{Q}^{+}_{\mathscr{A}}
\! + \! \mathfrak{Q}^{-}_{\mathscr{A}})^{\prime} \! = \! 0$. Also,
\begin{align*}
g^{o}_{+}(z) \! + \! g^{o}_{-}(z) \! - \! \widetilde{V}(z) \! - \! \ell_{o} \!
- \! (\mathfrak{Q}^{+}_{\mathscr{A}} \! + \! \mathfrak{Q}^{-}_{\mathscr{A}})=&
\, 2 \! \left(2 \! + \! \dfrac{1}{n} \right) \! \int_{J_{o}} \ln (\vert z \!
- \! s \vert) \psi_{V}^{o}(s) \, \md s \! - \! 2 \ln \vert z \vert \! - \!
\widetilde{V}(z) \\
-& \, \ell_{o} \! - \! 2 \! \left(2 \! + \! \dfrac{1}{n} \right) \! Q_{o}:
\end{align*}
proceeding as in the asymptotic analysis for case~$\mathbf{(3)}$ above, one
arrives at
\begin{equation*}
g^{o}_{+}(z) \! + \! g^{o}_{-}(z) \! - \! \widetilde{V}(z) \! - \! \ell_{o} \!
- \! (\mathfrak{Q}^{+}_{\mathscr{A}} \! + \! \mathfrak{Q}^{-}_{\mathscr{A}})
\underset{z \to -\infty}{=} \left(1 \! + \! \dfrac{1}{n} \right) \! \ln (z^{
2} \! + \! 1) \! - \! \widetilde{V}(z) \! + \! \mathcal{O}(1),
\end{equation*}
and
\begin{equation*}
g^{o}_{+}(z) \! + \! g^{o}_{-}(z) \! - \! \widetilde{V}(z) \! - \! \ell_{o} \!
- \! (\mathfrak{Q}^{+}_{\mathscr{A}} \! + \! \mathfrak{Q}^{-}_{\mathscr{A}})
\underset{\vert z \vert \to 0}{=} \ln (z^{-2} \! + \! 1) \! - \! \widetilde{
V}(z) \! + \! \mathcal{O}(1),
\end{equation*}
whence, via conditions~(2.4) and~(2.5), $g^{o}_{+}(z) \! + \! g^{o}_{-}(z) \!
- \! \widetilde{V}(z) \! - \! \ell_{o} \! - \! (\mathfrak{Q}^{+}_{\mathscr{A}}
\! + \! \mathfrak{Q}^{-}_{\mathscr{A}}) \! < \! 0$, $z \! \in \! (-\infty,
b_{0}^{o})$. \hfill $\qed$
\section{The Model RHP and Parametrices}
In this section, the (normalised at zero) auxiliary RHP for $\overset{o}{
\mathscr{M}} \colon \mathbb{C} \setminus \mathbb{R} \! \to \! \operatorname{S
L}_{2}(\mathbb{C})$ formulated in Lemma 3.4 is augmented, by means of a
sequence of contour deformations and transformations \emph{\`{a} la}
Deift-Venakides-Zhou \cite{a1,a2,a3}, into simpler, `model' (matrix) RHPs
which, as $n \! \to \! \infty$, are solved explicitly (in closed form) in
terms of Riemann theta functions (associated with the underlying genus-$N$
hyperelliptic Riemann surface) and Airy functions, and which give rise to the
leading $(\mathcal{O}(1))$ terms of asymptotics for $\boldsymbol{\pi}_{2n+1}
(z)$, $\xi^{(2n+1)}_{-n-1}$ and $\phi_{2n+1}(z)$ stated, respectively, in
Theorems~2.3.1 and~2.3.2, and the asymptotic (as $n \! \to \! \infty)$
analysis of the parametrices, which are `approximate' solutions of the RHP for
$\overset{o}{\mathscr{M}} \colon \mathbb{C} \setminus \mathbb{R} \! \to \!
\operatorname{SL}_{2}(\mathbb{C})$ in neighbourhoods of the end-points of the
support of the `odd' equilibrium measure, and which give rise to the
$\mathcal{O}(n \! + \! 1/2)$ (and $\mathcal{O}((n \! + \! 1/2)^{2}))$
corrections for $\boldsymbol{\pi}_{2n+1}(z)$, $\xi^{(2n+1)}_{-n-1}$ and
$\phi_{2n+1}(z)$ stated, respectively, in Theorems~2.3.1 and~2.3.2, is
undertaken.
\begin{ccccc}
Let the external field $\widetilde{V} \colon \mathbb{R} \setminus \{0\} \!
\to \! \mathbb{R}$ satisfy conditions~{\rm (2.3)--(2.5);} furthermore, let
$\widetilde{V}$ be regular. Let the `odd' equilibrium measure, $\mu_{V}^{o}$,
and its support, $\operatorname{supp}(\mu_{V}^{o}) \! =: \! J_{o} \! = \!
\cup_{j=1}^{N+1}J_{j}^{o} \! := \! \cup_{j=1}^{N+1}(b_{j-1}^{o},a_{j}^{o})$,
be as described in Lemma~{\rm 3.5}, and, along with $\ell_{o}$ $(\in \!
\mathbb{R})$, the `odd' variational constant, satisfy the variational
conditions given in Lemma~{\rm 3.6}, Equations~{\rm (3.9);} moreover, let
conditions~{\rm (1)}--{\rm (4)} stated in Lemma~{\rm 3.6} be valid. Then the
{\rm RHP} for $\overset{o}{\mathscr{M}} \colon \mathbb{C} \setminus \mathbb{R}
\! \to \! \mathrm{SL}_{2}(\mathbb{C})$ formulated in Lemma~{\rm 3.4} can be
equivalently reformulated as follows: {\rm (1)} $\overset{o}{\mathscr{M}}(z)$
is holomorphic for $z \! \in \! \mathbb{C} \setminus \mathbb{R};$ {\rm (2)}
$\overset{o}{\mathscr{M}}_{\pm}(z) \! := \! \lim_{\underset{\pm \Im (z^{
\prime})>0}{z^{\prime} \to z}} \overset{o}{\mathscr{M}}(z^{\prime})$ satisfy 
the boundary condition
\begin{equation*}
\overset{o}{\mathscr{M}}_{+}(z) \! = \! \overset{o}{\mathscr{M}}_{-}(z)
\overset{o}{\upsilon}(z), \quad z \! \in \! \mathbb{R},
\end{equation*}
where, for $i \! = \! 1,\dotsc,N \! + \! 1$ and $j \! = \! 1,\dotsc,N$,
\begin{equation*}
\overset{o}{\upsilon}(z) \! = \!
\begin{cases}
\begin{pmatrix}
\me^{-4(n+\frac{1}{2}) \pi \mi \int_{z}^{a_{N+1}^{o}} \psi_{V}^{o}(s) \, \md
s} \, \me^{\mi \mathfrak{q}_{o}} & 1 \\
0 & \me^{4(n+\frac{1}{2}) \pi \mi \int_{z}^{a_{N+1}^{o}} \psi_{V}^{o}(s) \,
\md s} \, \me^{-\mi \mathfrak{q}_{o}}
\end{pmatrix}, &\text{$z \! \in \! (b_{i-1}^{o},a_{i}^{o})$,} \\
\begin{pmatrix}
\me^{-4(n+\frac{1}{2}) \pi \mi \int_{b_{j}^{o}}^{a_{N+1}^{o}} \psi_{V}^{o}(s)
\, \md s} \, \me^{\mi \mathfrak{q}_{o}} & \me^{n(g^{o}_{+}(z)+g^{o}_{-}(z)-
\widetilde{V}(z)-\ell_{o}-\mathfrak{Q}^{+}_{\mathscr{A}}-\mathfrak{Q}^{-}_{
\mathscr{A}})} \\
0 & \me^{4(n+\frac{1}{2}) \pi \mi \int_{b_{j}^{o}}^{a_{N+1}^{o}} \psi_{V}^{o}
(s) \, \md s} \, \me^{-\mi \mathfrak{q}_{o}}
\end{pmatrix}, &\text{$z \! \in \! (a_{j}^{o},b_{j}^{o})$,} \\
\begin{pmatrix}
\me^{\mi \mathfrak{q}_{o}} & \me^{n(g^{o}_{+}(z)+g^{o}_{-}(z)-\widetilde{V}(z)
-\ell_{o}-\mathfrak{Q}^{+}_{\mathscr{A}}-\mathfrak{Q}^{-}_{\mathscr{A}})} \\
0 & \me^{-\mi \mathfrak{q}_{o}}
\end{pmatrix}, &\text{$z \! \in \! \mathfrak{I}$,}
\end{cases}
\end{equation*}
with
\begin{equation*}
\mathfrak{q}_{o} \! := \! 4 \pi \! \left(n \! + \! \dfrac{1}{2} \right) \!
\int_{J_{o} \cap \mathbb{R}_{+}} \psi_{V}^{o}(s) \, \md s,
\end{equation*}
$\mathfrak{I} \! := \! (-\infty,b_{0}^{o}) \! \cup \! (a_{N+1}^{o},+\infty)$,
$g^{o}(z)$ and $\mathfrak{Q}^{\pm}_{\mathscr{A}}$ defined in Lemma~{\rm 3.4},
\begin{equation*}
\pm \Re \! \left(\mi \int_{z}^{a_{N+1}^{o}} \psi_{V}^{o}(s) \, \md s \right)
\! > \! 0, \quad z \! \in \! \mathbb{C}_{\pm} \cap (\cup_{j=1}^{N+1} \mathbb{
U}_{j}^{o}),
\end{equation*}
where $\mathbb{U}_{j}^{o} \! := \! \lbrace \mathstrut z \! \in \! \mathbb{C}^{
\ast}; \, \Re (z) \! \in \! (b_{j-1}^{o},a_{j}^{o}), \, \inf_{q \in J_{j}^{o}}
\vert z \! - \! q \vert \! < \! r_{j} \! \in \! (0,1) \rbrace$, $j \! = \! 1,
\dotsc,N \! + \! 1$, with $\mathbb{U}_{i}^{o} \cap \mathbb{U}_{j}^{o} \! = \!
\varnothing$, $i \! \not= \! j \! = \! 1,\dotsc,N \! + \! 1$, and $g^{o}_{+}
(z) \! + \! g^{o}_{-}(z) \! - \! \widetilde{V}(z) \! - \! \ell_{o} \! - \!
\mathfrak{Q}^{+}_{\mathscr{A}} \! - \! \mathfrak{Q}^{-}_{\mathscr{A}} \! < \!
0$, $z \! \in \! \mathfrak{I} \cup (\cup_{j=1}^{N+1}(a_{j}^{o},b_{j}^{o}));$
{\rm (3)} $\overset{o}{\mathscr{M}}(z) \! =_{\underset{z \in \mathbb{C}
\setminus \mathbb{R}}{z \to 0}} \! \mathrm{I} \! + \! \mathcal{O}(z);$ and
{\rm (4)} $\overset{o}{\mathscr{M}}(z) \! =_{\underset{z \in \mathbb{C}
\setminus \mathbb{R}}{z \to \infty}} \! \mathcal{O}(1)$.
\end{ccccc}

\emph{Proof.} Item~(1) stated in the Lemma is simply a re-statement of
item~(1) of Lemma~3.4. Write $\mathbb{R} \! = \! (-\infty,b_{0}^{o}) \cup
(a_{N+1}^{o},+\infty) \cup (\cup_{j=1}^{N+1}J_{j}^{o}) \cup (\cup_{k=1}^{N}
(a_{k}^{o},b_{k}^{o})) \cup (\cup_{l=1}^{N+1} \lbrace b_{l-1}^{o},a_{l}^{o}
\rbrace)$, where $J_{j}^{o} \! := \! (b_{j-1}^{o},a_{j}^{o})$, $j \! = \!
1,\dotsc,N \! + \! 1$. Recall {}from the proof of Lemma~3.6 that, for
$\widetilde{V}$, $\mu_{V}^{o}$, and $\ell_{o}$ described therein (and in the
Lemma): (1) $g^{o}_{+}(z) \! - \! g^{o}_{-}(z) \! - \! \mathfrak{Q}^{+}_{
\mathscr{A}} \! + \! \mathfrak{Q}^{-}_{\mathscr{A}} \! =$
\begin{equation*}
\begin{cases}
2(2 \! + \! \tfrac{1}{n}) \pi \mi \int_{z}^{a_{N+1}^{o}} \psi_{V}^{o}(s) \,
\md s \! - \! 2(2 \! + \! \tfrac{1}{n}) \pi \mi \int_{J_{o} \cap \mathbb{R}_{
+}} \psi_{V}^{o}(s) \, \md s \! + \!
\begin{cases}
0, &\text{$z \! \in \! \mathbb{R}_{+} \cap \overline{J_{j}^{o}}, \quad j \! =
\! 1,\dotsc,N \! + \! 1$,} \\
-2 \pi \mi, &\text{$z \! \in \! \mathbb{R}_{-} \cap \overline{J_{j}^{o}},
\quad j \! = \! 1,\dotsc,N \! + \! 1$,}
\end{cases} \\
2(2 \! + \! \tfrac{1}{n}) \pi \mi \int_{b_{j}^{o}}^{a_{N+1}^{o}} \psi_{V}^{o}
(s) \, \md s \! - \! 2(2 \! + \! \tfrac{1}{n}) \pi \mi \int_{J_{o} \cap
\mathbb{R}_{+}} \psi_{V}^{o}(s) \, \md s \! + \!
\begin{cases}
0, &\text{$z \! \in \! \mathbb{R}_{+} \cap (a_{j}^{o},b_{j}^{o}), \quad j \!
= \! 1,\dotsc,N$,} \\
-2 \pi \mi, &\text{$z \! \in \! \mathbb{R}_{-} \cap (a_{j}^{o},b_{j}^{o}),
\quad j \! = \! 1,\dotsc,N$,}
\end{cases} \\
-2(2 \! + \! \tfrac{1}{n}) \pi \mi \int_{J_{o} \cap \mathbb{R}_{+}} \psi_{V}^{
o}(s) \, \md s \! + \!
\begin{cases}
0, &\text{$z \! \in \! \mathbb{R}_{+} \cap (a_{N+1}^{o},+\infty)$,} \\
-2 \pi \mi, &\text{$z \! \in \! \mathbb{R}_{-} \cap (a_{N+1}^{o},+\infty)$,}
\end{cases} \\
2(2 \! + \! \tfrac{1}{n}) \pi \mi \! - \! 2(2 \! + \! \tfrac{1}{n}) \pi \mi
\int_{J_{o} \cap \mathbb{R}_{+}} \psi_{V}^{o}(s) \, \md s \! + \!
\begin{cases}
0, &\text{$z \! \in \! \mathbb{R}_{+} \cap (-\infty,b_{0}^{o})$,} \\
-2 \pi \mi, &\text{$z \! \in \! \mathbb{R}_{-} \cap (-\infty,b_{0}^{o})$,}
\end{cases}
\end{cases}
\end{equation*}
where $\overline{J_{j}^{o}}$ $(:= \! J_{j}^{o} \cup \partial J_{j}^{o}) \! =
\! [b_{j-1}^{o},a_{j}^{o}]$, $j \! = \! 1,\dotsc,N \! + \! 1$; and (2) for $j
\! = \! 1,\dotsc,N \! + \! 1$,
\begin{equation*}
g^{o}_{+}(z) \! + \! g^{o}_{-}(z) \! - \! \widetilde{V}(z) \! - \! \ell_{o} \!
- \! \mathfrak{Q}^{+}_{\mathscr{A}} \! - \! \mathfrak{Q}^{-}_{\mathscr{A}} \!
= \!
\begin{cases}
0, &\text{$z \! \in \! \cup_{j=1}^{N+1} \overline{J_{j}^{o}}$,} \\
-(2 \! + \! \tfrac{1}{n}) \int_{a_{j}^{o}}^{z}(R_{o}(s))^{1/2}h_{V}^{o}(s) \,
\md s \! < \! 0,
&\text{$z \! \in \! (a_{j}^{o},b_{j}^{o})$,} \\
-(2 \! + \! \tfrac{1}{n}) \int_{a_{N+1}^{o}}^{z}(R_{o}(s))^{1/2}h_{V}^{o}(s)
\, \md s \! < \! 0, &\text{$z \! \in \! (a_{N+1}^{o},+\infty)$,} \\
(2 \! + \! \tfrac{1}{n}) \int_{z}^{b_{0}^{o}}(R_{o}(s))^{1/2}h_{V}^{o}(s) \,
\md s \! < \! 0, &\text{$z \! \in \! (-\infty,b_{0}^{o})$.}
\end{cases}
\end{equation*}
Recall, also, the formula for the `jump matrix' given in Lemma~3.4, namely,
\begin{equation*}
\begin{pmatrix}
\me^{-n(g^{o}_{+}(z)-g^{o}_{-}(z)-\mathfrak{Q}^{+}_{\mathscr{A}}+\mathfrak{Q}^{
-}_{\mathscr{A}})} & \me^{n(g^{o}_{+}(z)+g^{o}_{-}(z)-\widetilde{V}(z)-\ell_{o}
-\mathfrak{Q}^{+}_{\mathscr{A}}-\mathfrak{Q}^{-}_{\mathscr{A}})} \\
0 & \me^{n(g^{o}_{+}(z)-g^{o}_{-}(z)-\mathfrak{Q}^{+}_{\mathscr{A}}+\mathfrak{
Q}^{-}_{\mathscr{A}})}
\end{pmatrix}.
\end{equation*}
Partitioning $\mathbb{R}$ as given above, one obtains the formula for
$\overset{o}{\upsilon}(z)$ stated in the Lemma, thus item~(2); moreover,
items~(3) and~(4) are re-statements of the respective items of Lemma~3.4. It
remains, therefore, to show that $\Re (\mi \int_{z}^{a_{N+1}^{o}} \psi_{V}^{o}
(s) \, \md s)$ satisfies the inequalities stated in the Lemma. Recall {}from
the proof of Lemma~3.4 that $g^{o}(z)$ is uniformly Lipschitz continuous in
$\mathbb{C}_{\pm}$; moreover, via the Cauchy-Riemann conditions, item~(4) of
Lemma~3.6, that is, $\mi (g^{o}_{+}(z) \! - \! g^{o}_{-}(z) \! - \! \mathfrak{
Q}^{+}_{\mathscr{A}} \! + \! \mathfrak{Q}^{-}_{\mathscr{A}})^{\prime} \!
\geqslant \! 0$, $z \! \in \! J_{o}$, implies that the quantity $g^{o}_{+}(z)
\! - \! g^{o}_{-}(z) \! - \! \mathfrak{Q}^{+}_{\mathscr{A}} \! + \! \mathfrak{
Q}^{-}_{\mathscr{A}}$ has an analytic continuation, $\mathscr{G}^{o}(z)$, say,
to an open neighbourhood, $\mathbb{U}_{V}^{o}$, say, of $J_{o} \! = \!
\cup_{j=1}^{N+1}(b_{j-1}^{o},a_{j}^{o})$, where $\mathbb{U}_{V}^{o} \! :=
\! \cup_{j=1}^{N+1} \mathbb{U}_{j}^{o}$, with $\mathbb{U}_{j}^{o} \! := \!
\lbrace \mathstrut z \! \in \! \mathbb{C}^{\ast}; \, \Re (z) \! \in \! (b_{j
-1}^{o},a_{j}^{o}), \, \inf_{q \in J_{j}^{o}} \vert z \! - \! q \vert \! <
\! r_{j} \! \in \! (0,1) \rbrace$, $j \! = \! 1,\dotsc,N \! + \! 1$, and
$\mathbb{U}_{i}^{o} \cap \mathbb{U}_{j}^{o} \! = \! \varnothing$, $i \! \not=
\! j \! = \! 1,\dotsc,N \! + \! 1$, with the property that $\pm \Re (\mathscr{
G}^{o}(z)) \! > \! 0$, $z \! \in \! \mathbb{C}_{\pm} \cap \mathbb{U}_{V}^{
o}$. \hfill $\qed$
\begin{eeeee}
Recalling that the external field $\widetilde{V} \colon \mathbb{R} \setminus
\{0\} \! \to \! \mathbb{R}$ is regular, that is, $h_{V}^{o}(z) \! \not\equiv
\! 0 \, \, \forall \, \, z \! \in \! \overline{J_{j}^{o}} \! := \! \cup_{j=
1}^{N+1}[b_{j-1}^{o},a_{j}^{o}]$, the second inequality in Equations~(3.9) is
strict, namely, $2(2 \! + \! \tfrac{1}{n}) \int_{J_{o}} \ln (\vert x \! - \! s
\vert) \psi_{V}^{o}(s) \, \md s \! - \! 2 \ln \vert x \vert \! - \! \widetilde{
V}(z) \! - \! \ell_{o} \! - \! 2(2 \! + \! \tfrac{1}{n})Q_{o} \! < \! 0$, $x
\! \in \! \mathbb{R} \setminus \overline{J_{o}}$, and ({}from the proof of
Lemma~4.1) that $g^{o}_{+}(z) \! + \! g^{o}_{-}(z) \! - \! \widetilde{V}(z)
\! - \! \ell_{o} \! - \! 2(2 \! + \! \tfrac{1}{n})Q_{o} \! < \! 0$, $z \! \in
\! (-\infty,b_{0}^{o}) \cup (a_{N+1}^{o},+\infty) \cup (\cup_{j=1}^{N}(a_{j}^{
o},b_{j}^{o}))$, it follows that
\begin{equation*}
\overset{o}{\upsilon}(z) \underset{n \to \infty}{=}
\begin{cases}
\me^{-(4(n+\frac{1}{2}) \pi \mi \int_{b_{j}^{o}}^{a_{N+1}^{o}} \psi_{V}^{o}
(s) \, \md s) \sigma_{3}} \, \me^{\mi \mathfrak{q}_{o} \sigma_{3}} \! \left(
\mathrm{I} \! + \! o(1) \sigma_{+} \right),
&\text{$z \! \in \! (a_{j}^{o},b_{j}^{o}), \quad j \! = \! 1,\dotsc,N$,} \\
\me^{\mi \mathfrak{q}_{o} \sigma_{3}} \! \left(\mathrm{I} \! + \! o(1)
\sigma_{+} \right), &\text{$z \! \in \! (-\infty,b_{0}^{o}) \cup (a_{N+1}^{o},
+\infty)$,}
\end{cases}
\end{equation*}
where $o(1)$ denotes terms that are exponentially small. \hfill $\blacksquare$
\end{eeeee}
\begin{bbbbb}
Let the external field $\widetilde{V} \colon \mathbb{R} \setminus \{0\} \!
\to \! \mathbb{R}$ satisfy conditions~{\rm (2.3)--(2.5);} furthermore, let
$\widetilde{V}$ be regular. Let the `odd' equilibrium measure, $\mu_{V}^{o}$,
and its support, $\operatorname{supp}(\mu_{V}^{o}) \! =: \! J_{o} \! = \!
\cup_{j=1}^{N+1}J_{j}^{o} \! := \! \cup_{j=1}^{N+1}(b_{j-1}^{o},a_{j}^{o})$,
be as described in Lemma~{\rm 3.5}, and, along with $\ell_{o}$ $(\in \!
\mathbb{R})$, the `odd' variational constant, satisfy the variational
conditions given in Lemma~{\rm 3.6}, Equations~{\rm (3.9);} moreover, let
conditions~{\rm (1)}--{\rm (4)} stated in Lemma~{\rm 3.6} be valid. Let
$\overset{o}{\mathscr{M}}(z) \colon \mathbb{C} \setminus \mathbb{R} \! \to
\! \operatorname{SL}_{2}(\mathbb{C})$ solve the {\rm RHP} formulated in
Lemma~{\rm 4.1}. Set
\begin{equation*}
\overset{o}{\mathscr{M}}^{\raise-1.0ex\hbox{$\scriptstyle \flat$}}(z) \! = \!
\begin{cases}
\overset{o}{\mathscr{M}}(z) \mathbb{E}^{-\sigma_{3}}, &\text{$z \! \in \!
\mathbb{C}_{+}$,} \\
\overset{o}{\mathscr{M}}(z) \mathbb{E}^{\sigma_{3}}, &\text{$z \! \in \!
\mathbb{C}_{-}$,}
\end{cases}
\end{equation*}
where
\begin{equation*}
\mathbb{E} \! := \! \exp \! \left(\mi \mathfrak{q}_{o}/2 \right) \! = \!
\exp \! \left(\mi 2 \pi \! \left(n \! + \! \dfrac{1}{2} \right) \! \int_{J_{
o} \cap \mathbb{R}_{+}} \psi_{V}^{o}(s) \, \md s \right).
\end{equation*}
Then $\overset{o}{\mathscr{M}}^{\raise-1.0ex\hbox{$\scriptstyle \flat$}}(z)
\colon \mathbb{C} \setminus \mathbb{R} \! \to \! \operatorname{SL}_{2}
(\mathbb{C})$ solves the following {\rm RHP:} {\rm (1)}
$\overset{o}{\mathscr{M}}^{\raise-1.0ex\hbox{$\scriptstyle \flat$}}(z)$ is
holomorphic for $z \! \in \! \mathbb{C} \setminus \mathbb{R};$
{\rm (2)}
$\overset{o}{\mathscr{M}}^{\raise-1.0ex\hbox{$\scriptstyle \flat$}}_{\pm}(z)
\! := \! \lim_{\underset{\pm \Im (z^{\prime})>0}{z^{\prime} \to z}} 
\overset{o}{\mathscr{M}}^{\raise-1.0ex\hbox{$\scriptstyle \flat$}}
(z^{\prime})$ satisfy the boundary condition
\begin{equation*}
\overset{o}{\mathscr{M}}^{\raise-1.0ex\hbox{$\scriptstyle \flat$}}_{+}(z) \! =
\! \overset{o}{\mathscr{M}}^{\raise-1.0ex\hbox{$\scriptstyle \flat$}}_{-}(z)
\mathcal{V}_{
\overset{o}{\mathscr{M}}^{\raise-1.0ex\hbox{$\scriptstyle \flat$}}}(z), \quad
z \! \in \! \mathbb{R},
\end{equation*}
where, for $i \! = \! 1,\dotsc,N \! + \! 1$ and $j \! = \! 1,\dotsc,N$,
\begin{equation*}
\mathcal{V}_{
\overset{o}{\mathscr{M}}^{\raise-1.0ex\hbox{$\scriptstyle \flat$}}}(z) \! = \!
\begin{cases}
\begin{pmatrix}
\me^{-4(n+\frac{1}{2}) \pi \mi \int_{z}^{a_{N+1}^{o}} \psi_{V}^{o}(s) \, \md
s} & 1 \\
0 & \me^{4(n+\frac{1}{2}) \pi \mi \int_{z}^{a_{N+1}^{o}} \psi_{V}^{o}(s) \,
\md s}
\end{pmatrix}, &\text{$z \! \in \! (b_{i-1}^{o},a_{i}^{o})$,} \\
\begin{pmatrix}
\me^{-4(n+\frac{1}{2}) \pi \mi \int_{b_{j}^{o}}^{a_{N+1}^{o}} \psi_{V}^{o}(s)
\, \md s} & \me^{n(g^{o}_{+}(z)+g^{o}_{-}(z)-\widetilde{V}(z)-\ell_{o}-
\mathfrak{Q}^{+}_{\mathscr{A}}-\mathfrak{Q}^{-}_{\mathscr{A}})} \\
0 & \me^{4(n+\frac{1}{2}) \pi \mi \int_{b_{j}^{o}}^{a_{N+1}^{o}} \psi_{V}^{o}
(s) \, \md s}
\end{pmatrix}, &\text{$z \! \in \! (a_{j}^{o},b_{j}^{o})$,} \\
\mathrm{I} \! + \! \me^{n(g^{o}_{+}(z)+g^{o}_{-}(z)-\widetilde{V}(z)-\ell_{o}-
\mathfrak{Q}^{+}_{\mathscr{A}}-\mathfrak{Q}^{-}_{\mathscr{A}})} \sigma_{+},
&\text{$z \! \in \! \mathfrak{I}$,}
\end{cases}
\end{equation*}
with $\mathfrak{I} \! := \! (-\infty,b_{0}^{o}) \! \cup \! (a_{N+1}^{o},+
\infty)$, $g^{o}(z)$ and $\mathfrak{Q}^{\pm}_{\mathscr{A}}$ defined in
Lemma~{\rm 3.4},
\begin{equation*}
\pm \Re \! \left(\mi \int_{z}^{a_{N+1}^{o}} \psi_{V}^{o}(s) \, \md s \right)
\! > \! 0, \quad z \! \in \! \mathbb{C}_{\pm} \cap (\cup_{j=1}^{N+1} \mathbb{
U}_{j}^{o}),
\end{equation*}
where $\mathbb{U}_{j}^{o} \! := \! \lbrace \mathstrut z \! \in \! \mathbb{C}^{
\ast}; \, \Re (z) \! \in \! (b_{j-1}^{o},a_{j}^{o}), \, \inf_{q \in J_{j}^{o}}
\vert z \! - \! q \vert \! < \! r_{j} \! \in \! (0,1) \rbrace$, $j \! = \! 1,
\dotsc,N \! + \! 1$, with $\mathbb{U}_{i}^{o} \cap \mathbb{U}_{j}^{o} \! = \!
\varnothing$, $i \! \not= \! j \! = \! 1,\dotsc,N \! + \! 1$, and $g^{o}_{+}
(z) \! + \! g^{o}_{-}(z) \! - \! \widetilde{V}(z) \! - \! \ell_{o} \! - \!
\mathfrak{Q}^{+}_{\mathscr{A}} \! - \! \mathfrak{Q}^{-}_{\mathscr{A}} \! < \!
0$, $z \! \in \! \mathfrak{I} \cup (\cup_{j=1}^{N+1}(a_{j}^{o},b_{j}^{o}));$
{\rm (3)} $\overset{o}{\mathscr{M}}^{\raise-1.0ex\hbox{$\scriptstyle \flat$}}
(z) \! =_{\underset{z \in \mathbb{C}_{+}}{z \to 0}} \! (\mathrm{I} \! + \!
\mathcal{O}(z)) \mathbb{E}^{-\sigma_{3}}$ and
$\overset{o}{\mathscr{M}}^{\raise-1.0ex\hbox{$\scriptstyle \flat$}}(z) \! =_{
\underset{z \in \mathbb{C}_{-}}{z \to 0}} \! (\mathrm{I} \! + \! \mathcal{O}
(z)) \mathbb{E}^{\sigma_{3}};$ and {\rm (4)}
$\overset{o}{\mathscr{M}}^{\raise-1.0ex\hbox{$\scriptstyle \flat$}}(z) \! =_{
\underset{z \in \mathbb{C} \setminus \mathbb{R}}{z \to \infty}} \! \mathcal{O}
(1)$.
\end{bbbbb}

\emph{Proof.} Follows {}from the definition of
$\overset{o}{\mathscr{M}}^{\raise-1.0ex\hbox{$\scriptstyle \flat$}}(z)$ in
terms of $\overset{o}{\mathscr{M}}(z)$ given in the Proposition and the RHP
for $\overset{o}{\mathscr{M}}(z)$ formulated in Lemma~4.1. \hfill $\qed$
\begin{ccccc}
Let the external field $\widetilde{V} \colon \mathbb{R} \setminus \{0\} \!
\to \! \mathbb{R}$ satisfy conditions~{\rm (2.3)--(2.5);} furthermore, let
$\widetilde{V}$ be regular. Let the `odd' equilibrium measure, $\mu_{V}^{o}$,
and its support, $\operatorname{supp}(\mu_{V}^{o}) \! =: \! J_{o} \! = \!
\cup_{j=1}^{N+1}J_{j}^{o} \! := \! \cup_{j=1}^{N+1}(b_{j-1}^{o},a_{j}^{o})$,
be as described in Lemma~{\rm 3.5}, and, along with $\ell_{o}$ $(\in \!
\mathbb{R})$, the `odd' variational constant, satisfy the variational
conditions given in Lemma~{\rm 3.6}, Equations~{\rm (3.9);} moreover, let
conditions~{\rm (1)}--{\rm (4)} stated in Lemma~{\rm 3.6} be valid. Let
$\overset{o}{\mathscr{M}}^{\raise-1.0ex\hbox{$\scriptstyle \flat$}}(z) \colon
\mathbb{C} \setminus \mathbb{R} \! \to \! \mathrm{SL}_{2}(\mathbb{C})$ solve
the {\rm RHP} formulated in Proposition~{\rm 4.1}, and let the deformed (and
oriented) contour $\Sigma^{\sharp}_{o} \! := \! \mathbb{R} \cup (\cup_{j=
1}^{N+1}(J_{j}^{o,\smallfrown} \cup J_{j}^{o,\smallsmile}))$ be as in
Figure~{\rm 8} below; furthermore, $\cup_{j=1}^{N+1}(\Omega_{j}^{o,
\smallfrown} \cup \Omega_{j}^{o,\smallsmile} \cup J_{j}^{o,\smallfrown} 
\cup J_{j}^{o,\smallsmile}) \subset \cup_{j=1}^{N+1} \mathbb{U}_{j}^{o}$
(Figure~{\rm 8)}, where $\mathbb{U}_{j}^{o}$, $j \! = \! 1,\dotsc,N \! + \! 
1$, is defined in Lemma~{\rm 4.1}. Set
\begin{equation*}
\overset{o}{\mathscr{M}}^{\raise-1.0ex\hbox{$\scriptstyle \sharp$}}(z) \! 
:= \!
\begin{cases}
\overset{o}{\mathscr{M}}^{\raise-1.0ex\hbox{$\scriptstyle \flat$}}(z),
&\text{$z \! \in \! \mathbb{C} \setminus (\Sigma_{o}^{\sharp} \cup (\cup_{j=
1}^{N+1}(\Omega_{j}^{o,\smallfrown} \cup \Omega_{j}^{o,\smallsmile})))$,} \\
\overset{o}{\mathscr{M}}^{\raise-1.0ex\hbox{$\scriptstyle \flat$}}(z) \! \left(
\mathrm{I} \! - \! \me^{-4(n+\frac{1}{2}) \pi \mi \int_{z}^{a_{N+1}^{o}} \psi_{
V}^{o}(s) \, \md s} \sigma_{-} \right), &\text{$z \! \in \! \mathbb{C}_{+}
\cap (\cup_{j=1}^{N+1} \Omega_{j}^{o,\smallfrown})$,} \\
\overset{o}{\mathscr{M}}^{\raise-1.0ex\hbox{$\scriptstyle \flat$}}(z) \! \left(
\mathrm{I} \! + \! \me^{4(n+\frac{1}{2}) \pi \mi \int_{z}^{a_{N+1}^{o}} \psi_{
V}^{o}(s) \, \md s} \sigma_{-} \right), &\text{$z \! \in \! \mathbb{C}_{-}
\cap (\cup_{j=1}^{N+1} \Omega_{j}^{o,\smallsmile})$.}
\end{cases}
\end{equation*}
Then $\overset{o}{\mathscr{M}}^{\raise-1.0ex\hbox{$\scriptstyle \sharp$}}
\colon \mathbb{C} \setminus \Sigma_{o}^{\sharp} \! \to \! \mathrm{SL}_{2}(
\mathbb{C})$ solves the following, equivalent {\rm RHP:} {\rm (1)} $\overset{
o}{\mathscr{M}}^{\raise-1.0ex\hbox{$\scriptstyle \sharp$}}(z)$ is holomorphic
for $z \! \in \! \mathbb{C} \setminus \Sigma_{o}^{\sharp};$ {\rm (2)}
$\overset{o}{\mathscr{M}}^{\raise-1.0ex\hbox{$\scriptstyle \sharp$}}_{\pm}
(z) \! := \! \lim_{\underset{z^{\prime} \, \in \, \pm \, \mathrm{side} \,
\mathrm{of} \, \Sigma_{o}^{\sharp}}{z^{\prime} \to z}} 
\overset{o}{\mathscr{M}}^{\raise-1.0ex\hbox{$\scriptstyle \sharp$}}
(z^{\prime})$ satisfy the boundary condition
\begin{equation*}
\overset{o}{\mathscr{M}}^{\raise-1.0ex\hbox{$\scriptstyle \sharp$}}_{+}(z) \!
= \! \overset{o}{\mathscr{M}}^{\raise-1.0ex\hbox{$\scriptstyle \sharp$}}_{-}
(z) \overset{o}{\upsilon}^{\raise-1.0ex\hbox{$\scriptstyle \sharp$}}(z), \quad
z \! \in \! \textstyle \Sigma_{o}^{\sharp},
\end{equation*}
where, for $i \! = \! 1,\dotsc,N \! + \! 1$ and $j \! = \! 1,\dotsc,N$,
\begin{equation*}
\overset{o}{\upsilon}^{\raise-1.0ex\hbox{$\scriptstyle \sharp$}}(z) \! = \!
\begin{cases}
\mi \sigma_{2}, &\text{$z \! \in \! J_{i}^{o}$,} \\
\mathrm{I} \! + \! \me^{-4(n+\frac{1}{2}) \pi \mi \int_{z}^{a_{N+1}^{o}} \psi_{
V}^{o}(s) \, \md s} \sigma_{-}, &\text{$z \! \in \! J_{i}^{o,\smallfrown}$,} \\
\mathrm{I} \! + \! \me^{4(n+\frac{1}{2}) \pi \mi \int_{z}^{a_{N+1}^{o}} \psi_{
V}^{o}(s) \, \md s} \sigma_{-}, &\text{$z \! \in \! J_{i}^{o,\smallsmile}$,} \\
\begin{pmatrix}
\me^{-4(n+\frac{1}{2}) \pi \mi \int_{b_{j}^{o}}^{a_{N+1}^{o}} \psi_{V}^{o}
(s) \, \md s} & \me^{n(g^{o}_{+}(z)+g^{o}_{-}(z)-\widetilde{V}(z)-\ell_{o}-
\mathfrak{Q}^{+}_{\mathscr{A}}-\mathfrak{Q}^{-}_{\mathscr{A}})} \\
0 & \me^{4(n+\frac{1}{2}) \pi \mi \int_{b_{j}^{o}}^{a_{N+1}^{o}} \psi_{V}^{o}
(s) \, \md s}
\end{pmatrix}, &\text{$z \! \in \! (a_{j}^{o},b_{j}^{o})$,} \\
\mathrm{I} \! + \! \me^{n(g^{o}_{+}(z)+g^{o}_{-}(z)-\widetilde{V}(z)-\ell_{o}-
\mathfrak{Q}^{+}_{\mathscr{A}}-\mathfrak{Q}^{-}_{\mathscr{A}})} \sigma_{+},
&\text{$z \! \in \! \mathfrak{I}$,}
\end{cases}
\end{equation*}
with $\mathfrak{I} \! := \! (-\infty,b_{0}^{o}) \cup (a_{N+1}^{o},+\infty)$,
and $\Re (\mi \int_{z}^{a_{N+1}^{o}} \psi_{V}^{o}(s) \, \md s) \! > \! 0$
(resp., $\Re (\mi \int_{z}^{a_{N+1}^{o}} \psi_{V}^{o}(s) \, \md s) \! < \!
0)$, $z \! \in \! \Omega_{j}^{o,\smallfrown}$ (resp., $z \! \in \! \Omega_{
j}^{o,\smallsmile});$ {\rm (3)}
\begin{align*}
\overset{o}{\mathscr{M}}^{\raise-1.0ex\hbox{$\scriptstyle \sharp$}}(z)
\underset{\underset{z \in \mathbb{C}_{+} \setminus \cup_{j=1}^{N+1}(J_{j}^{o,
\smallfrown} \cup \Omega_{j}^{o,\smallfrown})}{z \to 0}}{=}& \, \left(\mathrm{
I} \! + \! \mathcal{O}(z) \right) \! \mathbb{E}^{-\sigma_{3}}, \\
\overset{o}{\mathscr{M}}^{\raise-1.0ex\hbox{$\scriptstyle \sharp$}}(z)
\underset{\underset{z \in \mathbb{C}_{-} \setminus \cup_{j=1}^{N+1}(J_{j}^{o,
\smallsmile} \cup \Omega_{j}^{o,\smallsmile})}{z \to 0}}{=}& \, \left(\mathrm{
I} \! + \! \mathcal{O}(z) \right) \! \mathbb{E}^{\sigma_{3}};
\end{align*}
and {\rm (4)}
\begin{equation*}
\overset{o}{\mathscr{M}}^{\raise-1.0ex\hbox{$\scriptstyle \sharp$}}(z) \!
\underset{\underset{z \in \mathbb{C} \setminus (\Sigma_{o}^{\sharp} \cup
(\cup_{j=1}^{N+1}(\Omega_{j}^{o,\smallfrown} \cup \Omega_{j}^{o,\smallsmile}))
)}{z \to \infty}}{=} \!
\mathcal{O}(1).
\end{equation*}
\end{ccccc}
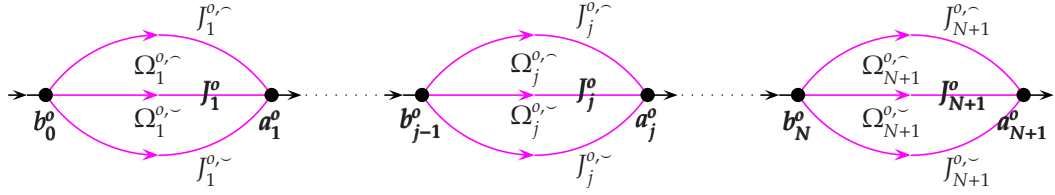
\begin{figure}[tbh]
\begin{center}
\vspace{-0.40cm}
\begin{pspicture}(0,0)(14,5)
\psset{xunit=1cm,yunit=1cm}
\psarcn[linewidth=0.6pt,linestyle=solid,linecolor=magenta,arrowsize=1.5pt 5]%
{->}(2,1.5){1.8}{146}{90}
\psarcn[linewidth=0.6pt,linestyle=solid,linecolor=magenta](2,1.5){1.8}{90}{34}
\psarc[linewidth=0.6pt,linestyle=solid,linecolor=magenta,arrowsize=1.5pt 5]%
{->}(2,3.5){1.8}{214}{270}
\psarc[linewidth=0.6pt,linestyle=solid,linecolor=magenta](2,3.5){1.8}{270}{326}
\psline[linewidth=0.6pt,linestyle=solid,linecolor=black,arrowsize=1.5pt 4]%
{->}(0,2.5)(0.25,2.5)
\psline[linewidth=0.6pt,linestyle=solid,linecolor=black](0.25,2.5)(0.5,2.5)
\psline[linewidth=0.6pt,linestyle=solid,linecolor=magenta,arrowsize=1.5pt 5]%
{->}(0.5,2.5)(2,2.5)
\psline[linewidth=0.6pt,linestyle=solid,linecolor=magenta](2,2.5)(3.5,2.5)
\psline[linewidth=0.6pt,linestyle=solid,linecolor=black,arrowsize=1.5pt 4]%
{->}(3.5,2.5)(3.9,2.5)
\psline[linewidth=0.6pt,linestyle=solid,linecolor=black,arrowsize=1.5pt 4]%
{->}(5,2.5)(5.25,2.5)
\psline[linewidth=0.6pt,linestyle=solid,linecolor=black](5.25,2.5)(5.5,2.5)
\psline[linewidth=0.6pt,linestyle=solid,linecolor=magenta,arrowsize=1.5pt 5]%
{->}(5.5,2.5)(7,2.5)
\psline[linewidth=0.6pt,linestyle=solid,linecolor=magenta](7,2.5)(8.5,2.5)
\psline[linewidth=0.6pt,linestyle=solid,linecolor=black,arrowsize=1.5pt 4]%
{->}(8.5,2.5)(8.9,2.5)
\psarcn[linewidth=0.6pt,linestyle=solid,linecolor=magenta,arrowsize=1.5pt 5]%
{->}(7,1.5){1.8}{146}{90}
\psarcn[linewidth=0.6pt,linestyle=solid,linecolor=magenta](7,1.5){1.8}{90}{34}
\psarc[linewidth=0.6pt,linestyle=solid,linecolor=magenta,arrowsize=1.5pt 5]%
{->}(7,3.5){1.8}{214}{270}
\psarc[linewidth=0.6pt,linestyle=solid,linecolor=magenta](7,3.5){1.8}{270}{326}
\psline[linewidth=0.6pt,linestyle=solid,linecolor=black,arrowsize=1.5pt 4]%
{->}(10,2.5)(10.25,2.5)
\psline[linewidth=0.6pt,linestyle=solid,linecolor=black](10.25,2.5)(10.5,2.5)
\psline[linewidth=0.6pt,linestyle=solid,linecolor=magenta,arrowsize=1.5pt 5]%
{->}(10.5,2.5)(12,2.5)
\psline[linewidth=0.6pt,linestyle=solid,linecolor=magenta](12,2.5)(13.5,2.5)
\psline[linewidth=0.6pt,linestyle=solid,linecolor=black,arrowsize=1.5pt 4]%
{->}(13.5,2.5)(13.9,2.5)
\psarcn[linewidth=0.6pt,linestyle=solid,linecolor=magenta,arrowsize=1.5pt 5]%
{->}(12,1.5){1.8}{146}{90}
\psarcn[linewidth=0.6pt,linestyle=solid,linecolor=magenta](12,1.5){1.8}{90}{34}
\psarc[linewidth=0.6pt,linestyle=solid,linecolor=magenta,arrowsize=1.5pt 5]%
{->}(12,3.5){1.8}{214}{270}
\psarc[linewidth=0.6pt,linestyle=solid,linecolor=magenta](12,3.5){1.8}{270}%
{326}
\psline[linewidth=0.7pt,linestyle=dotted,linecolor=darkgray](3.95,2.5)%
(4.9,2.5)
\psline[linewidth=0.7pt,linestyle=dotted,linecolor=darkgray](8.95,2.5)%
(9.9,2.5)
\psdots[dotstyle=*,dotscale=1.5](0.5,2.5)
\psdots[dotstyle=*,dotscale=1.5](3.5,2.5)
\psdots[dotstyle=*,dotscale=1.5](5.5,2.5)
\psdots[dotstyle=*,dotscale=1.5](8.5,2.5)
\psdots[dotstyle=*,dotscale=1.5](10.5,2.5)
\psdots[dotstyle=*,dotscale=1.5](13.5,2.5)
\rput(0.5,2.1){\makebox(0,0){$\pmb{b_{0}^{o}}$}}
\rput(3.5,2.1){\makebox(0,0){$\pmb{a_{1}^{o}}$}}
\rput(5.5,2.1){\makebox(0,0){$\pmb{b_{j-1}^{o}}$}}
\rput(8.5,2.1){\makebox(0,0){$\pmb{a_{j}^{o}}$}}
\rput(10.5,2.1){\makebox(0,0){$\pmb{b_{N}^{o}}$}}
\rput(13.5,2.1){\makebox(0,0){$\pmb{a_{N+1}^{o}}$}}
\rput(2.7,2.5){\makebox(0,0){$\pmb{J_{1}^{o}}$}}
\rput(7.7,2.5){\makebox(0,0){$\pmb{J_{j}^{o}}$}}
\rput(12.7,2.5){\makebox(0,0){$\pmb{J_{N+1}^{o}}$}}
\rput(2.75,3.5){\makebox(0,0){$J_{1}^{o,\smallfrown}$}}
\rput(2.75,1.5){\makebox(0,0){$J_{1}^{o,\smallsmile}$}}
\rput(2,2.85){\makebox(0,0){$\Omega_{1}^{o,\smallfrown}$}}
\rput(2,2.15){\makebox(0,0){$\Omega_{1}^{o,\smallsmile}$}}
\rput(7.8,3.5){\makebox(0,0){$J_{j}^{o,\smallfrown}$}}
\rput(7.8,1.5){\makebox(0,0){$J_{j}^{o,\smallsmile}$}}
\rput(7,2.85){\makebox(0,0){$\Omega_{j}^{o,\smallfrown}$}}
\rput(7,2.15){\makebox(0,0){$\Omega_{j}^{o,\smallsmile}$}}
\rput(12.75,3.5){\makebox(0,0){$J_{N+1}^{o,\smallfrown}$}}
\rput(12.75,1.5){\makebox(0,0){$J_{N+1}^{o,\smallsmile}$}}
\rput(11.75,2.85){\makebox(0,0){$\Omega_{N+1}^{o,\smallfrown}$}}
\rput(11.75,2.15){\makebox(0,0){$\Omega_{N+1}^{o,\smallsmile}$}}
\end{pspicture}
\end{center}
\vspace{-1.00cm}
\caption{Oriented/deformed contour $\Sigma_{o}^{\sharp} \! := \! \mathbb{R}
\cup (\cup_{j=1}^{N+1}(J_{j}^{o,\smallfrown} \cup J_{j}^{o,\smallsmile}))$}
\end{figure}

\emph{Proof.} Items~(1), (3), and~(4) in the formulation of the RHP for
$\overset{o}{\mathscr{M}}^{\raise-1.0ex\hbox{$\scriptstyle \sharp$}} \colon
\mathbb{C} \setminus \Sigma_{o}^{\sharp} \! \to \! \mathrm{SL}_{2}(\mathbb{
C})$ follow {}from the definition of
$\overset{o}{\mathscr{M}}^{\raise-1.0ex\hbox{$\scriptstyle \sharp$}}(z)$ (in
terms of $\overset{o}{\mathscr{M}}^{\raise-1.0ex\hbox{$\scriptstyle \flat$}}
(z))$ given in the Lemma and the respective items~(1), (3), and~(4) for the
RHP for $\overset{o}{\mathscr{M}}^{\raise-1.0ex\hbox{$\scriptstyle \flat$}}
\colon \mathbb{C} \setminus \mathbb{R} \! \to \! \mathrm{SL}_{2}(\mathbb{C})$
stated in Proposition~4.1; it remains, therefore, to verify item~(2), that is,
the formula for
$\overset{o}{\upsilon}^{\raise-1.0ex\hbox{$\scriptstyle \sharp$}}(z)$. Recall
{}from item~(2) of Proposition~4.1 that, for $z \! \in \! (b_{j-1}^{o},
a_{j}^{o})$, $j \! = \! 1,\dotsc,N \! + \! 1$,
$\overset{o}{\mathscr{M}}^{\raise-1.0ex\hbox{$\scriptstyle \flat$}}_{+}(z) \!
= \! \overset{o}{\mathscr{M}}^{\raise-1.0ex\hbox{$\scriptstyle \flat$}}_{-}(z)
\mathcal{V}_{\overset{o}{\mathscr{M}}^{\raise-1.0ex\hbox{$\scriptstyle
\flat$}}}(z)$, where $\mathcal{V}_{
\overset{o}{\mathscr{M}}^{\raise-1.0ex\hbox{$\scriptstyle \flat$}}}(z) \! = \!
\left(
\begin{smallmatrix}
\me^{-4(n+\frac{1}{2}) \pi \mi \int_{z}^{a_{N+1}^{o}} \psi_{V}^{o}(s) \, \md
s} & 1 \\
0 & \me^{4(n+\frac{1}{2}) \pi \mi \int_{z}^{a_{N+1}^{o}} \psi_{V}^{o}(s) \,
\md s}
\end{smallmatrix}
\right)$: noting the matrix factorisation
\begin{align*}
\begin{pmatrix}
\me^{-4(n+\frac{1}{2}) \pi \mi \int_{z}^{a_{N+1}^{o}} \psi_{V}^{o}(s) \, \md
s} & 1 \\
0 & \me^{4(n+\frac{1}{2}) \pi \mi \int_{z}^{a_{N+1}^{o}} \psi_{V}^{o}(s) \,
\md s}
\end{pmatrix} =& \,
\begin{pmatrix}
1 & 0 \\
\me^{4(n+\frac{1}{2}) \pi \mi \int_{z}^{a_{N+1}^{o}} \psi_{V}^{o}(s) \, \md s}
& 1
\end{pmatrix} \!
\begin{pmatrix}
0 & 1 \\
-1 & 0
\end{pmatrix} \\
\times& \,
\begin{pmatrix}
1 & 0 \\
\me^{-4(n+\frac{1}{2}) \pi \mi \int_{z}^{a_{N+1}^{o}} \psi_{V}^{o}(s) \, \md
s} & 1
\end{pmatrix},
\end{align*}
it follows that, $z \! \in \! (b_{j-1}^{o},a_{j}^{o})$, $j \! = \! 1,\dotsc,N
\! + \! 1$,
\begin{equation*}
\overset{o}{\mathscr{M}}^{\raise-1.0ex\hbox{$\scriptstyle \flat$}}_{+}(z) \!
\begin{pmatrix}
1 & 0 \\
-\me^{-4(n+\frac{1}{2}) \pi \mi \int_{z}^{a_{N+1}^{o}} \psi_{V}^{o}(s) \, \md
s} & 1
\end{pmatrix} \! =
\overset{o}{\mathscr{M}}^{\raise-1.0ex\hbox{$\scriptstyle \flat$}}_{-}(z) \!
\begin{pmatrix}
1 & 0 \\
\me^{4(n+\frac{1}{2}) \pi \mi \int_{z}^{a_{N+1}^{o}} \psi_{V}^{o}(s) \, \md
s} & 1
\end{pmatrix} \mi \sigma_{2}.
\end{equation*}
It was shown in Lemma~4.1 that $\pm \Re (\mi \int_{z}^{a_{N+1}^{o}} \psi_{V}^{
o}(s) \, \md s) \! > \! 0$ for $z \! \in \! \mathbb{C}_{\pm} \cap \mathbb{U}_{
j}^{o}$, where $\mathbb{U}_{j}^{o} \! := \! \lbrace \mathstrut z \! \in \!
\mathbb{C}^{\ast}; \, \Re (z) \! \in \! (b_{j-1}^{o},a_{j}^{o}), \, \inf_{q
\in J_{j}^{o}} \vert z \! - \! q \vert \! < \! r_{j} \! \in \! (0,1) \rbrace$,
$j \! = \! 1,\dotsc,N \! + \! 1$, with $\mathbb{U}_{i}^{o} \cap \mathbb{U}_{
j}^{o} \! = \! \varnothing$, $i \! \not= \! j \! = \! 1,\dotsc,N \! + \! 1$,
and $J_{j}^{o} \! := \! (b_{j-1}^{o},a_{j}^{o})$, $j \! = \! 1,\dotsc,N \! +
\! 1$. (One notes that the terms $\pm 4(n \! + \! \tfrac{1}{2}) \pi \mi \int_{
z}^{a_{N+1}^{o}} \psi_{V}^{o}(s) \, \md s$, which are pure imaginary for $z
\! \in \! \mathbb{R}$, and corresponding to which $\exp (\pm 4(n \! + \!
\tfrac{1}{2}) \pi \mi \int_{z}^{a_{N+1}^{o}} \psi_{V}^{o}(s) \, \md s)$ are
undulatory, are continued analytically to $\mathbb{C}_{\pm} \cap (\cup_{j=1}^{
N+1} \mathbb{U}_{j}^{o})$, respectively, corresponding to which $\exp (\pm 4
(n \! + \! \tfrac{1}{2}) \pi \mi \int_{z}^{a_{N+1}^{o}} \psi_{V}^{o}(s) \,
\md s)$ are exponentially decreasing as $n \! \to \! \infty)$. As per the DZ
non-linear steepest-descent method \cite{a1,a2} (see, also, the extension
\cite{a3}), one now `deforms' the original (and oriented) contour $\mathbb{R}$
to the deformed, or extended, (and oriented) contour/skeleton $\Sigma_{o}^{
\sharp} \! := \! \mathbb{R} \cup (\cup_{j=1}^{N+1}(J_{j}^{o,\smallfrown} \cup
J_{j}^{o,\smallsmile}))$ (Figure~8) in such a way that the upper (resp.,
lower) `lips' of the `lenses' $J_{j}^{o,\smallfrown}$ (resp., $J_{j}^{o,
\smallsmile})$, $j \! = \! 1,\dotsc,N \! + \! 1$, which are the boundaries of
$\Omega_{j}^{o,\smallfrown}$ (resp., $\Omega_{j}^{o,\smallsmile})$, $j \! =
\! 1,\dotsc,N \! + \! 1$, respectively, lie within the domain of analytic
continuation of $g^{o}_{+}(z) \! - \! g^{o}_{-}(z) \! - \! \mathfrak{Q}^{
+}_{\mathscr{A}} \! + \! \mathfrak{Q}^{-}_{\mathscr{A}}$ (cf. the proof of
Lemma~4.1), that is, $\cup_{j=1}^{N+1}(\Omega_{j}^{o,\smallfrown} \cup
\Omega_{j}^{o,\smallsmile} \cup J_{j}^{o,\smallfrown} \cup J_{j}^{o,
\smallsmile}) \subset \cup_{j=1}^{N+1} \mathbb{U}_{j}^{o}$; in particular,
each (oriented) interval $J_{j}^{o} \! = \! (b_{j-1}^{o},a_{j}^{o})$, $j \! =
\! 1,\dotsc,N \! + \! 1$, in the original (and oriented) contour $\mathbb{R}$
is `split' (or branched) into three, and the new (and oriented) contour
$\Sigma_{o}^{\sharp}$ is the old contour $(\mathbb{R})$ together with the
(oriented) boundary of $N \! + \! 1$ lens-shaped regions, one region
surrounding each (bounded and oriented) interval $J_{j}^{o}$. Now,
recalling the definition of
$\overset{o}{\mathscr{M}}^{\raise-1.0ex\hbox{$\scriptstyle \sharp$}}(z)$ (in
terms of $\overset{o}{\mathscr{M}}^{\raise-1.0ex\hbox{$\scriptstyle \flat$}}
(z))$ stated in the Lemma, and the expressions for (the jump matrix)
$\mathcal{V}_{\overset{o}{\mathscr{M}}^{\raise-1.0ex\hbox{$\scriptstyle
\flat$}}}(z)$ given in Proposition~4.1, one arrives at the formula for
$\overset{o}{\upsilon}^{\raise-1.0ex\hbox{$\scriptstyle \sharp$}}(z)$ given
in item~(2) of the Lemma. \hfill $\qed$
\begin{eeeee}
The jump condition stated in item~(2) of Lemma~4.2, that is,
$\overset{o}{\mathscr{M}}^{\raise-1.0ex\hbox{$\scriptstyle \sharp$}}_{+}(z) \!
= \! \overset{o}{\mathscr{M}}^{\raise-1.0ex\hbox{$\scriptstyle \sharp$}}_{-}
(z) \overset{o}{\upsilon}^{\raise-1.0ex\hbox{$\scriptstyle \sharp$}}(z)$,
$z \! \in \! \Sigma_{o}^{\sharp}$, with
$\overset{o}{\upsilon}^{\raise-1.0ex\hbox{$\scriptstyle \sharp$}}(z)$ given
therein, should, of course, be understood as follows: the $\mathrm{SL}_{2}(
\mathbb{C})$-valued functions
$\overset{o}{\mathscr{M}}^{\raise-1.0ex\hbox{$\scriptstyle \sharp$}} \! \!
\upharpoonright_{\mathbb{C}_{\pm} \setminus \Sigma_{o}^{\sharp}}$ have a
continuous extension to $\Sigma_{o}^{\sharp}$ with boundary values
$\overset{o}{\mathscr{M}}^{\raise-1.0ex\hbox{$\scriptstyle \sharp$}}_{\pm}
(z) \! := \! \lim_{\underset{z^{\prime} \, \in \, \pm \, \mathrm{side} \,
\mathrm{of} \, \Sigma_{o}^{\sharp}}{z^{\prime} \to z \in \Sigma_{o}^{\sharp}}}
\! \overset{o}{\mathscr{M}}^{\raise-1.0ex\hbox{$\scriptstyle \sharp$}}
(z^{\prime})$ satisfying the above jump relation
($\overset{o}{\mathscr{M}}^{\raise-1.0ex\hbox{$\scriptstyle \sharp$}}(z)$ is
continuous in each component of $\mathbb{C} \setminus \Sigma_{o}^{\sharp}$ up
to the boundary with boundary values
$\overset{o}{\mathscr{M}}^{\raise-1.0ex\hbox{$\scriptstyle \sharp$}}_{\pm}(z)$
satisfying the above jump relation on $\Sigma_{o}^{\sharp})$. \hfill
$\blacksquare$
\end{eeeee}

Recalling {}from Proposition~4.1 that, for $z \! \in \! (-\infty,b_{0}^{o})
\cup (a_{N+1}^{o},+\infty) \cup (\cup_{j=1}^{N}(a_{j}^{o},b_{j}^{o}))$, $g^{
o}_{+}(z) \! + \! g^{o}_{-}(z) \! - \! \widetilde{V}(z) \! - \! \ell_{o} \! -
\! \mathfrak{Q}^{+}_{\mathscr{A}} \! - \! \mathfrak{Q}^{-}_{\mathscr{A}} \!
< \! 0$, and, {}from Lemma~4.2, $\Re (\mi \int_{z}^{a_{N+1}^{o}} \psi_{V}^{o}
(s) \, \md s) \! > \! 0$ for $z \! \in \! J_{j}^{o,\smallfrown}$ (resp., $\Re
(\mi \int_{z}^{a_{N+1}^{o}} \psi_{V}^{o}(s) \, \md s) \! < \! 0$ for $z \! \in
\! J_{j}^{o,\smallsmile})$, $j \! = \! 1,\dotsc,N \! + \! 1$, one arrives at
the following large-$n$ asymptotic behaviour for the jump matrix
$\overset{o}{\upsilon}^{\raise-1.0ex\hbox{$\scriptstyle \sharp$}}(z)$: for 
$i \! = \! 1,\dotsc,N \! + \! 1$ and $j \! = \! 1,\dotsc,N$,
\begin{equation*}
\overset{o}{\upsilon}^{\raise-1.0ex\hbox{$\scriptstyle \sharp$}}(z) \!
\underset{n \to \infty}{=} \!
\begin{cases}
\mi \sigma_{2}, &\text{$z \! \in \! J_{i}^{o}$,} \\
\mathrm{I} \! + \! \mathcal{O} \! \left(\me^{-(n+\frac{1}{2})c \vert z \vert}
\right) \! \sigma_{-}, &\text{$z \! \in \! J_{i}^{o,\smallfrown} \cup J_{i}^{
o,\smallsmile}$,} \\
\me^{-(4(n+\frac{1}{2}) \pi \mi \int_{b_{j}^{o}}^{a_{N+1}^{o}} \psi_{V}^{o}(s)
\, \md s) \sigma_{3}} \! \left(\mathrm{I} \! + \! \mathcal{O}(\me^{-(n+\frac{
1}{2})c \vert z-a_{j}^{o} \vert}) \sigma_{+} \right), &\text{$z \! \in \! (a_{
j}^{o},b_{j}^{o}) \setminus \widehat{\mathbb{U}}_{\delta_{0}^{o}}(0)$,} \\
\me^{-(4(n+\frac{1}{2}) \pi \mi \int_{b_{j}^{o}}^{a_{N+1}^{o}} \psi_{V}^{o}(s)
\, \md s) \sigma_{3}} \! \left(\mathrm{I} \! + \! \mathcal{O}(\me^{-(n+\frac{
1}{2})c \vert z \vert^{-1}}) \sigma_{+} \right), &\text{$z \! \in \! (a_{j}^{
o},b_{j}^{o}) \cap \widehat{\mathbb{U}}_{\delta_{0}^{o}}(0)$,} \\
\mathrm{I} \! + \! \mathcal{O} \! \left(\me^{-(n+\frac{1}{2})c \vert z \vert}
\right) \! \sigma_{+}, &\text{$z \! \in \! \mathfrak{I} \setminus \widehat{
\mathbb{U}}_{\delta_{0}^{o}}(0)$,} \\
\mathrm{I} \! + \! \mathcal{O} \! \left(\me^{-(n+\frac{1}{2})c \vert z \vert^{
-1}} \right) \! \sigma_{+}, &\text{$z \! \in \! \mathfrak{I} \cap \widehat{
\mathbb{U}}_{\delta_{0}^{o}}(0)$,}
\end{cases}
\end{equation*}
where $c$ (some generic number) $> \! 0$, $\widehat{\mathbb{U}}_{\delta_{0}^{
o}}(0) \! := \! \lbrace \mathstrut z \! \in \! \mathbb{C}; \, \vert z \vert \!
< \! \delta_{0}^{o} \rbrace$, with $\delta_{0}^{o}$ some arbitrarily fixed,
sufficiently small positive real number, $\mathfrak{I} \! := \! (-\infty,b_{
0}^{o}) \cup (a_{N+1}^{o},+\infty)$, and where the respective convergences
are normal, that is, uniform in (respective) compact subsets (see Section~5
below).

Recall {}from Lemma~2.56 of \cite{a1} that, for an oriented skeleton in
$\mathbb{C}$ on which the jump matrix of an RHP is defined, one may always
choose to add or delete a portion of the skeleton on which the jump matrix
equals $\mathrm{I}$ without altering the RHP in the operator sense; hence,
neglecting those jumps on $\Sigma_{o}^{\sharp}$ tending exponentially quickly
(as $n \! \to \! \infty)$ to $\mathrm{I}$, and removing the corresponding
oriented skeletons from $\Sigma_{o}^{\sharp}$, it becomes more or less
transparent how to construct a parametrix, that is, an approximate solution,
of the RHP for
$\overset{o}{\mathscr{M}}^{\raise-1.0ex\hbox{$\scriptstyle \sharp$}}
\colon \mathbb{C} \setminus \Sigma_{o}^{\sharp} \! \to \! \mathrm{SL}_{2}(
\mathbb{C})$ stated in Lemma~4.2, namely, the large-$n$ solution of the RHP
for $\overset{o}{\mathscr{M}}^{\raise-1.0ex\hbox{$\scriptstyle \sharp$}}(z)$
formulated in Lemma~4.2 should be `close to' the solution of the following
limiting, or model, RHP (for
$\overset{o}{m}^{\raise-1.0ex\hbox{$\scriptstyle \infty$}}(z))$.
\begin{ccccc}
Let the external field $\widetilde{V} \colon \mathbb{R} \setminus \{0\} \!
\to \! \mathbb{R}$ satisfy conditions~{\rm (2.3)--(2.5);} furthermore, let
$\widetilde{V}$ be regular. Let the `odd' equilibrium measure, $\mu_{V}^{o}$,
and its support, $\operatorname{supp}(\mu_{V}^{o}) \! =: \! J_{o} \! = \!
\cup_{j=1}^{N+1}J_{j}^{o} \! := \! \cup_{j=1}^{N+1}(b_{j-1}^{o},a_{j}^{o})$,
be as described in Lemma~{\rm 3.5}, and, along with $\ell_{o}$ $(\in \!
\mathbb{R})$, the `odd' variational constant, satisfy the variational
conditions given in Lemma~{\rm 3.6}, Equations~{\rm (3.9);} moreover, let
conditions~{\rm (1)}--{\rm (4)} stated in Lemma~{\rm 3.6} be valid. Then
$\overset{o}{m}^{\raise-1.0ex\hbox{$\scriptstyle \infty$}} \colon \mathbb{C}
\setminus J_{o}^{\infty} \to \! \mathrm{SL}_{2}(\mathbb{C})$, where $J_{o}^{
\infty} \! := \! J_{o} \cup (\cup_{j=1}^{N}(a_{j}^{o},b_{j}^{o}))$, solves
the following (model) {\rm RHP:} {\rm (1)}
$\overset{o}{m}^{\raise-1.0ex\hbox{$\scriptstyle \infty$}}(z)$ is holomorphic
for $z \! \in \! \mathbb{C} \setminus J_{o}^{\infty};$ {\rm (2)}
$\overset{o}{m}^{\raise-1.0ex\hbox{$\scriptstyle \infty$}}_{\pm}(z) \! := \!
\lim_{\underset{z^{\prime} \, \in \, \pm \, \mathrm{side} \, \mathrm{of} \,
J_{o}^{\infty}}{z^{\prime} \to z}} 
\overset{o}{m}^{\raise-1.0ex\hbox{$\scriptstyle \infty$}}(z^{\prime})$
satisfy the boundary condition
\begin{equation*}
\overset{o}{m}^{\raise-1.0ex\hbox{$\scriptstyle \infty$}}_{+}(z) \! = \!
\overset{o}{m}^{\raise-1.0ex\hbox{$\scriptstyle \infty$}}_{-}(z)
\overset{o}{\upsilon}^{\raise-1.0ex\hbox{$\scriptstyle \infty$}}(z), \quad
z \! \in \! J_{o}^{\infty},
\end{equation*}
where
\begin{equation*}
\overset{o}{\upsilon}^{\raise-1.0ex\hbox{$\scriptstyle \infty$}}(z) \! = \!
\begin{cases}
\mi \sigma_{2}, &\text{$z\! \in \! (b_{i-1}^{o},a_{i}^{o}), \quad i \! = \! 1,
\dotsc,N \! + \! 1$,} \\
\me^{-(4(n+\frac{1}{2}) \pi \mi \int_{b_{j}^{o}}^{a_{N+1}^{o}} \psi_{V}^{o}(s)
\, \md s) \sigma_{3}}, &\text{$z \! \in \! (a_{j}^{o},b_{j}^{o}), \quad j \! =
\! 1,\dotsc,N;$}
\end{cases}
\end{equation*}
{\rm (3)} $\overset{o}{m}^{\raise-1.0ex\hbox{$\scriptstyle \infty$}}(z) \! =_{
\underset{z \in \mathbb{C}_{+} \setminus J_{o}^{\infty}}{z \to 0}} \! \left(
\mathrm{I} \! + \! \mathcal{O}(z) \right) \! \mathbb{E}^{-\sigma_{3}}$ and
$\overset{o}{m}^{\raise-1.0ex\hbox{$\scriptstyle \infty$}}(z) \! =_{\underset{
z \in \mathbb{C}_{-} \setminus J_{o}^{\infty}}{z \to 0}} \! \left(\mathrm{I}
\! + \! \mathcal{O}(z) \right) \! \mathbb{E}^{\sigma_{3}};$ and {\rm (4)}
$\overset{o}{m}^{\raise-1.0ex\hbox{$\scriptstyle \infty$}}(z) \! =_{\underset{
z \in \mathbb{C} \setminus J_{o}^{\infty}}{z \to \infty}} \! \mathcal{O}(1)$.
\end{ccccc}

The model RHP for $\overset{o}{m}^{\raise-1.0ex\hbox{$\scriptstyle \infty$}}
\colon \mathbb{C} \setminus J_{o}^{\infty} \! \to \! \operatorname{SL}_{2}
(\mathbb{C})$ formulated in Lemma~4.3 is (explicitly) solvable in terms of
Riemann theta functions (see, for example, Section~3 of \cite{a45}; see, 
also, Section~4.2 of \cite{a46}): the solution is succinctly presented below.
\begin{ccccc}
Let $\gamma^{o} \colon \mathbb{C} \setminus ((-\infty,b_{0}^{o}) \cup (a_{N+
1}^{o},+\infty) \cup (\cup_{j=1}^{N}(a_{j}^{o},b_{j}^{o}))) \! \to \! \mathbb{
C}$ be defined by
\begin{equation*}
\gamma^{o}(z) \! := \!
\begin{cases}
\left(\left(\dfrac{z \! - \! b_{0}^{o}}{z \! - \! a_{N+1}^{o}} \right) \!
\mathlarger{\prod_{k=1}^{N}} \! \left(\dfrac{z \! - \! b_{k}^{o}}{z \! - \!
a_{k}^{o}} \right) \right)^{1/4}, &\text{$z \! \in \! \mathbb{C}_{+}$,} \\
-\mi \left(\left(\dfrac{z \! - \! b_{0}^{o}}{z \! - \! a_{N+1}^{o}} \right) \!
\mathlarger{\prod_{k=1}^{N}} \! \left(\dfrac{z \! - \! b_{k}^{o}}{z \! - \!
a_{k}^{o}} \right) \right)^{1/4}, &\text{$z \! \in \! \mathbb{C}_{-}$,}
\end{cases}
\end{equation*}
and set
\begin{equation*}
\gamma^{o}(0) \! := \! \left(\prod_{k=1}^{N+1}b_{k-1}^{o}(a_{k}^{o})^{-1}
\right)^{1/4} \quad (> \! 0).
\end{equation*}
If $(\gamma^{o}(0))^{4} \!= \! 1$, then, on the lower edge of each
finite-length gap, that is, $(a_{j}^{o},b_{j}^{o})^{-}$, $j \! = \! 1,\dotsc,
N$, $(\gamma^{o}(0))^{-1} \gamma^{o} \linebreak[4]
(z) \! + \! \gamma^{o}(0)(\gamma^{o}(z))^{-1} \! = \! 0$ has exactly one
root/zero, and, on the upper edge of each finite-length gap, that is, $(a_{
j}^{o},b_{j}^{o})^{+}$, $j \! = \! 1,\dotsc,N$, $(\gamma^{o}(0))^{-1} \gamma^{
o}(z) \! - \! \gamma^{o}(0)(\gamma^{o}(z))^{-1} \! = \! 0$ has exactly one
root/zero; otherwise, if $(\gamma^{o}(0))^{4} \! \not= \! 1$, there is an
additional root/zero in the exterior/unbounded gap $(-\infty,b_{0}^{o}) \cup
(a_{N+1}^{o},+\infty)$. For both cases, label a set of $N$ of the lower-edge
and upper-edge finite-length-gap roots/zeros as
\begin{equation*}
\left\lbrace z_{j}^{o,\pm} \! \in \! (a_{j}^{o},b_{j}^{o})^{\pm} \subset
\mathbb{C}_{\pm}, \, j \! = \! 1,\dotsc,N; \, ((\gamma^{o}(0))^{-1} \gamma^{o}
(z) \! \mp \! \gamma^{o}(0)(\gamma^{o}(z))^{-1}) \vert_{z=z_{j}^{o,\pm}} \! =
\! 0 \right\rbrace
\end{equation*}
(in the plane, $z_{j}^{o,+} \! = \! z_{j}^{o,-} \! := \! z_{j}^{o} \! \in \!
(a_{j}^{o},b_{j}^{o})$, $j \! = \! 1,\dotsc,N)$. Furthermore, $\gamma^{o}(z)$
solves the following (scalar) {\rm RHP:}
\begin{compactenum}
\item[{\rm (1)}] $\gamma^{o}(z)$ is holomorphic for $z \! \in \! \mathbb{C}
\setminus ((-\infty,b_{0}^{o}) \cup (a_{N+1}^{o},+\infty) \cup (\cup_{j=1}^{N}
(a_{j}^{o},b_{j}^{o})));$
\item[{\rm (2)}] $\gamma^{o}_{+}(z) \! = \! \gamma^{o}_{-}(z) \mi$, $z \! \in
\! (-\infty,b_{0}^{o}) \cup (a_{N+1}^{o},+\infty) \cup (\cup_{j=1}^{N}(a_{j}^{
o},b_{j}^{o}));$
\item[{\rm (3)}] $\gamma^{o}(z) \! =_{\underset{z \in \mathbb{C}_{\pm}}{z \to
0}} \! (-\mi)^{(1 \mp 1)/2} \gamma^{o}(0)(1 \! + \! \mathcal{O}(z));$ and
\item[{\rm (4)}] $\gamma^{o}(z) \! =_{\underset{z \in \mathbb{C}_{\pm}}{z \to
\infty}} \! \mathcal{O}(1)$.
\end{compactenum}
\end{ccccc}

\emph{Proof.} Define $\gamma^{o}(z)$ as in the Lemma: then one notes that
$(\gamma^{o}(0))^{-1} \gamma^{o}(z) \! \mp \! \gamma^{o}(0)(\gamma^{o}(z))^{-
1} \! = \! 0 \! \Leftrightarrow \! (\gamma^{o}(z))^{2} \! \mp \! (\gamma^{o}
(0))^{2} \! = \! 0 \! \Rightarrow \! (\gamma^{o}(z))^{4} \! - \! (\gamma^{o}
(0))^{4} \! = \! 0 \! \Leftrightarrow \! \mathscr{Q}^{o}(z)$ $(\in \! \mathbb{
R}[z])$ $:= \! (z \! - \! b_{0}^{o}) \prod_{k=1}^{N}(z \! - \! b_{k}^{o}) \!
- \! (\gamma^{o}(0))^{4}(z \! - \! a_{N+1}^{o}) \prod_{k=1}^{N}(z \! - \! a_{
k}^{o}) \! = \! 0$, whence, via a straightforward calculation, and using the
fact that $(\gamma^{o}(0))^{4} \! = \! \prod_{k=1}^{N+1}b_{k-1}^{o}(a_{k}^{o}
)^{-1} \! > \! 0$, one shows that $\mathscr{Q}^{o}(a_{j}^{o}) \! = \! (-1)^{N-
j+1} \widehat{\mathscr{Q}}^{o}_{a_{j}^{o}}$, $j \! = \! 1,\dotsc,N$, where
$\widehat{\mathscr{Q}}^{o}_{a_{j}^{o}} \! := \! (b_{j}^{o} \! - \! a_{j}^{o})
(a_{j}^{o} \! - \! b_{0}^{o}) \prod_{k=1}^{j-1}(a_{j}^{o} \! - \! b_{k}^{o})
\prod_{l=j+1}^{N}(b_{l}^{o} \! - \! a_{j}^{o}) \! > \! 0$, and $\mathscr{Q}^{o}
(b_{j}^{o}) \! = \! -(-1)^{N-j+1} \widehat{\mathscr{Q}}^{o}_{b_{j}^{o}}$, $j
\! = \! 1,\dotsc,N$, where $\widehat{\mathscr{Q}}^{o}_{b_{j}^{o}} \! := \!
(\gamma^{o}(0))^{4}(b_{j}^{o} \! - \! a_{j}^{o})(a_{N+1}^{o} \! - \! b_{j}^{o}
) \prod_{k=1}^{j-1}(b_{j}^{o} \! - \! a_{k}^{o}) \prod_{l=j+1}^{N}(a_{l}^{o}
\! - \! b_{j}^{o}) \! > \! 0$; thus, $\mathscr{Q}^{o}(a_{j}^{o}) \mathscr{Q}^{
o}(b_{j}^{o}) \! < \! 0$, $j \! = \! 1,\dotsc,N$, which shows that: (i) for
$(\gamma^{o}(0))^{4} \! \not= \! 1$, since $\deg (\mathscr{Q}^{o}(z)) \! = \!
N \! + \! 1$, there are $N \! + \! 1$ (simple) roots/zeros of $\mathscr{Q}^{o}
(z)$, one in each (open) finite-length gap $(a_{j}^{o},b_{j}^{o})$, $j \! = \!
1,\dotsc,N$, and one in the (open) unbounded/exterior gap $(-\infty,b_{0}^{
o}) \cup (a_{N+1}^{o},+\infty)$; and (ii) for $(\gamma^{o}(0))^{4} \! = \!
1$, since $\deg (\mathscr{Q}^{o}(z)) \! = \! N$, there are $N$ (simple)
roots/zeros of $\mathscr{Q}^{o}(z)$, one in each (open) finite-length gap
$(a_{j}^{o},b_{j}^{o})$, $j \! = \! 1,\dotsc,N$. For both cases, label a set
of $N$ of the roots/zeros of $\mathscr{Q}^{o}(z)$ as $\lbrace z_{j}^{o}
\rbrace_{j=1}^{N}$. A straightforward analysis of the branch cuts shows that,
for $z \! \in \! \cup_{j=1}^{N}(a_{j}^{o},b_{j}^{o})^{\pm}$, $\pm (\gamma^{o}
(z))^{2} \! > \! 0$, whence $\lbrace z_{j}^{o,\pm} \rbrace_{j=1}^{N} \! =
\! \lbrace \mathstrut z^{\pm} \! \in \! (a_{j}^{o},b_{j}^{o})^{\pm} \subset
\mathbb{C}_{\pm}, \, j \! = \! 1,\dotsc,N; \, ((\gamma^{o}(0))^{-1} \gamma^{o}
(z) \! \mp \! \gamma^{o}(0)(\gamma^{o}(z))^{-1}) \vert_{z=z^{\pm}} \! = \! 0
\rbrace$. Setting $\widetilde{J}^{o} \! := \! (-\infty,b_{0}^{o}) \cup (a_{N+
1}^{o},+\infty) \cup (\cup_{j=1}^{N}(a_{j}^{o},b_{j}^{o}))$, one shows, upon
performing a straightforward analysis of the branch cuts, that $\gamma^{o}
(z)$ solves the RHP $(\gamma^{o}(z),\mi,\widetilde{J}^{o})$ formulated in
the Lemma. \hfill $\qed$

All of the notation/nomenclature used in Lemma~4.5 below has been defined at
the end of Subsection~2.1; the reader, therefore, is advised to peruse the
relevant notations(s), etc., before proceeding to Lemma~4.5. Let $\mathcal{Y}_{
o}$ denote the Riemann surface of $y^{2} \! = \! R_{o}(z) \! = \! \prod_{k=1}^{
N+1}(z \! - \! b_{k-1}^{o})(z \! - \! a_{k}^{o})$, where the single-valued
branch of the square root is chosen so that $z^{-(N+1)}(R_{o}(z))^{1/2} \!
\sim_{\underset{z \in \mathbb{C}_{\pm}}{z \to \infty}} \! \pm 1$. Let
$\mathscr{P} \! := \! (y,z)$ denote a point on the Riemann surface $\mathcal{
Y}_{o}$ $(:= \! \lbrace \mathstrut (y,z); \, y^{2} \! = \! R_{o}(z) \rbrace)$.
The notation $0^{\pm}$ (used in Lemma~4.5 below) means: $\mathscr{P} \! \to \!
0^{\pm} \! \Leftrightarrow \! z \! \to \! 0, y \! \sim \! \pm (-1)^{\mathcal{
N}_{+}}(\prod_{k=1}^{N+1} \vert b_{k-1}^{o}a_{k}^{o} \vert)^{1/2}$, where
$\mathcal{N}_{+} \! \in \! \lbrace 0,\dotsc,N \! + \! 1 \rbrace$ is the number
of bands to the right of $z \! = \! 0$.
\begin{ccccc}
Let $\overset{o}{m}^{\raise-1.0ex\hbox{$\scriptstyle \infty$}} \colon \mathbb{
C} \setminus J_{o}^{\infty} \to \! \mathrm{SL}_{2}(\mathbb{C})$ solve the
{\rm RHP} formulated in Lemma~{\rm 4.3}. Then,
\begin{equation*}
\overset{o}{m}^{\raise-1.0ex\hbox{$\scriptstyle \infty$}}(z) \! = \!
\begin{cases}
\overset{o}{\mathfrak{M}}^{\raise-1.0ex\hbox{$\scriptstyle \infty$}}(z),
&\text{$z \! \in \! \mathbb{C}_{+}$,} \\
-\mi \, \overset{o}{\mathfrak{M}}^{\raise-1.0ex\hbox{$\scriptstyle \infty$}}
(z) \sigma_{2}, &\text{$z \in \! \mathbb{C}_{-}$,}
\end{cases}
\end{equation*}
where
\begin{equation*}
\overset{o}{\mathfrak{M}}^{\raise-1.0ex\hbox{$\scriptstyle \infty$}}(z) \! :=
\! \mathbb{E}^{-\sigma_{3}} \!
\begin{pmatrix}
\frac{\boldsymbol{\theta}^{o}(\boldsymbol{u}^{o}_{+}(0)+\boldsymbol{d}_{o})}{
\boldsymbol{\theta}^{o}(\boldsymbol{u}^{o}_{+}(0)-\frac{1}{2 \pi}(n+\frac{1}{
2}) \boldsymbol{\Omega}^{o}+\boldsymbol{d}_{o})} & 0 \\
0 & \frac{\boldsymbol{\theta}^{o}(\boldsymbol{u}^{o}_{+}(0)+\boldsymbol{d}_{
o})}{\boldsymbol{\theta}^{o}(-\boldsymbol{u}^{o}_{+}(0)-\frac{1}{2 \pi}(n+
\frac{1}{2}) \boldsymbol{\Omega}^{o}-\boldsymbol{d}_{o})}\end{pmatrix} \!
\overset{o}{\boldsymbol{\Theta}}^{\raise-1.0ex\hbox{$\scriptstyle \infty$}}_{
\natural}(z),
\end{equation*}
and
\begin{gather*}
\overset{o}{\boldsymbol{\Theta}}^{\raise-1.0ex\hbox{$\scriptstyle \infty$}}_{
\natural}(z) \! = \!
\begin{pmatrix}
\frac{((\gamma^{o}(0))^{-1} \gamma^{o}(z)+\gamma^{o}(0)(\gamma^{o}(z))^{-1})}{
2} \overset{o}{\boldsymbol{\Theta}}^{\raise-1.0ex\hbox{$\scriptstyle \infty$}
}_{11}(z) & -\frac{((\gamma^{o}(0))^{-1} \gamma^{o}(z)-\gamma^{o}(0)(\gamma^{
o}(z))^{-1})}{2 \mi}
\overset{o}{\boldsymbol{\Theta}}^{\raise-1.0ex\hbox{$\scriptstyle \infty$}}_{1
2}(z) \\
\frac{((\gamma^{o}(0))^{-1} \gamma^{o}(z)-\gamma^{o}(0)(\gamma^{o}(z))^{-1})}{
2 \mi}
\overset{o}{\boldsymbol{\Theta}}^{\raise-1.0ex\hbox{$\scriptstyle \infty$}}_{
21}(z) & \frac{((\gamma^{o}(0))^{-1} \gamma^{o}(z)+\gamma^{o}(0)(\gamma^{o}
(z))^{-1})}{2}
\overset{o}{\boldsymbol{\Theta}}^{\raise-1.0ex\hbox{$\scriptstyle \infty$}}_{
22}(z)
\end{pmatrix}, \\
\overset{o}{\boldsymbol{\Theta}}^{\raise-1.0ex\hbox{$\scriptstyle \infty$}}_{
11}(z) \! := \! \dfrac{\boldsymbol{\theta}^{o}(\boldsymbol{u}^{o}(z) \! - \!
\frac{1}{2 \pi}(n \! + \! \frac{1}{2}) \boldsymbol{\Omega}^{o} \! + \!
\boldsymbol{d}_{o})}{\boldsymbol{\theta}^{o}(\boldsymbol{u}^{o}(z) \! + \!
\boldsymbol{d}_{o})}, \qquad \quad
\overset{o}{\boldsymbol{\Theta}}^{\raise-1.0ex\hbox{$\scriptstyle \infty$}}_{
12}(z) \! := \! \dfrac{\boldsymbol{\theta}^{o}(-\boldsymbol{u}^{o}(z) \! - \!
\frac{1}{2 \pi}(n \! + \! \frac{1}{2}) \boldsymbol{\Omega}^{o} \! + \!
\boldsymbol{d}_{o})}{\boldsymbol{\theta}^{o}(-\boldsymbol{u}^{o}(z) \! + \!
\boldsymbol{d}_{o})}, \\
\overset{o}{\boldsymbol{\Theta}}^{\raise-1.0ex\hbox{$\scriptstyle \infty$}}_{
21}(z) \! := \! \dfrac{\boldsymbol{\theta}^{o}(\boldsymbol{u}^{o}(z) \! - \!
\frac{1}{2 \pi}(n \! + \! \frac{1}{2}) \boldsymbol{\Omega}^{o} \! - \!
\boldsymbol{d}_{o})}{\boldsymbol{\theta}^{o}(\boldsymbol{u}^{o}(z) \! - \!
\boldsymbol{d}_{o})}, \qquad \quad
\overset{o}{\boldsymbol{\Theta}}^{\raise-1.0ex\hbox{$\scriptstyle \infty$}}_{
22}(z) \! := \! \dfrac{\boldsymbol{\theta}^{o}(-\boldsymbol{u}^{o}(z) \! - \!
\frac{1}{2 \pi}(n \! + \! \frac{1}{2}) \boldsymbol{\Omega}^{o} \! - \!
\boldsymbol{d}_{o})}{\boldsymbol{\theta}^{o}(\boldsymbol{u}^{o}(z) \! + \!
\boldsymbol{d}_{o})},
\end{gather*}
with $\gamma^{o}(z)$ characterised completely in Lemma~{\rm 4.4}, $\boldsymbol{
\Omega}^{o} \! := \! (\Omega_{1}^{o},\Omega_{2}^{o},\dotsc,\Omega_{N}^{o})^{
\mathrm{T}}$ $(\in \! \mathbb{R}^{N})$, where $\Omega_{j}^{o} \! = \! 4 \pi
\int_{b_{j}^{o}}^{a_{N+1}^{o}} \linebreak[4]
\psi_{V}^{o}(s) \, \md s$, $j \! = \! 1,\dotsc,N$, and ${}^{\mathrm{T}}$
denotes transposition, $\boldsymbol{d}_{o} \! \equiv \! -\sum_{j=1}^{N} \int_{
a_{j}^{o}}^{z_{j}^{o,-}} \! \boldsymbol{\omega}^{o}$ $(= \! \sum_{j=1}^{N}
\int_{a_{j}^{o}}^{z_{j}^{o,+}} \! \boldsymbol{\omega}^{o})$, $\lbrace z_{j}^{
o,\pm} \rbrace_{j=1}^{N}$ are characterised completely in Lemma~{\rm 4.4},
$\boldsymbol{\omega}^{o}$ is the associated normalised basis of holomorphic
one-forms of $\mathcal{Y}_{o}$, $\boldsymbol{u}^{o}(z) \! := \! \int_{a_{N+1
}^{o}}^{z} \boldsymbol{\omega}^{o}$ $(\in \! \operatorname{Jac}(\mathcal{Y}_{
o}))$, and $\boldsymbol{u}^{o}_{+}(0) \! := \! \int_{a_{N+1}^{o}}^{0^{+}}
\boldsymbol{\omega}^{o};$ furthermore, the solution is unique.
\end{ccccc}

\emph{Proof.} Let $\overset{o}{m}^{\raise-1.0ex\hbox{$\scriptstyle \infty$}}
\colon \mathbb{C} \setminus J_{o}^{\infty} \to \! \mathrm{SL}_{2}(\mathbb{C}
)$ solve the RHP formulated in Lemma~4.3, and define
$\overset{o}{m}^{\raise-1.0ex\hbox{$\scriptstyle \infty$}}(z)$, in terms of
$\overset{o}{\mathfrak{M}}^{\raise-1.0ex\hbox{$\scriptstyle \infty$}}(z)$, as
in the Lemma. A straightforward calculation shows that
$\overset{o}{\mathfrak{M}}^{\raise-1.0ex\hbox{$\scriptstyle \infty$}} \colon
\mathbb{C} \setminus \mathbb{R} \! \to \! \mathrm{SL}_{2}(\mathbb{C})$ solves
the following `twisted' RHP: (i)
$\overset{o}{\mathfrak{M}}^{\raise-1.0ex\hbox{$\scriptstyle \infty$}}(z)$ is
holomorphic for $z \! \in \! \mathbb{C} \setminus \widetilde{J}^{o}$, where
$\widetilde{J}^{o} \! := \! (-\infty,b_{0}^{o}) \cup (a_{N+1}^{o},+\infty)
\cup (\cup_{j=1}^{N}(a_{j}^{o},b_{j}^{o}))$; (ii)
$\overset{o}{\mathfrak{M}}^{\raise-1.0ex\hbox{$\scriptstyle \infty$}}_{\pm}
(z) \! := \! \lim_{\underset{z^{\prime} \, \in \, \pm \, \mathrm{side} \,
\mathrm{of} \, \widetilde{J}^{o}}{z^{\prime} \to z}} \! 
\overset{o}{\mathfrak{M}}^{\raise-1.0ex\hbox{$\scriptstyle \infty$}}
(z^{\prime})$ satisfy the boundary condition 
$\overset{o}{\mathfrak{M}}^{\raise-1.0ex\hbox{$\scriptstyle \infty$}}_{+}
\linebreak[4]
(z) \! = \!
\overset{o}{\mathfrak{M}}^{\raise-1.0ex\hbox{$\scriptstyle \infty$}}_{-}(z)
\overset{o}{\mathscr{V}}^{\raise-1.0ex\hbox{$\scriptstyle \infty$}}(z)$, $z \!
\in \! \widetilde{J}^{o}$, where
\begin{equation}
\overset{o}{\mathscr{V}}^{\raise-1.0ex\hbox{$\scriptstyle \infty$}}(z) \! :=
\!
\begin{cases}
\mathrm{I}, &\text{$z \! \in \! J_{o}$,} \\
-\mi \sigma_{2}, &\text{$z \! \in \! (-\infty,b_{0}^{o}) \cup (a_{N+1}^{o},
+\infty)$,} \\
-\mi \sigma_{2} \me^{-\mi (n+\frac{1}{2}) \Omega_{j}^{o} \sigma_{3}},
&\text{$z \! \in \! (a_{j}^{o},b_{j}^{o}), \quad j \! = \! 1,\dotsc,N$,}
\end{cases}
\end{equation}
with $\Omega_{j}^{o} \! = \! 4 \pi \int_{b_{j}^{o}}^{a_{N+1}^{o}} \psi_{V}^{
o}(s) \, \md s$, $j \! = \! 1,\dotsc,N$; (iii)
$\overset{o}{\mathfrak{M}}^{\raise-1.0ex\hbox{$\scriptstyle \infty$}}(z) \!
=_{\underset{z \in \mathbb{C}_{+}}{z \to 0}} \! \mathbb{E}^{-\sigma_{3}} \! +
\! \mathcal{O}(z)$ and
$\overset{o}{\mathfrak{M}}^{\raise-1.0ex\hbox{$\scriptstyle \infty$}}
(z) \! =_{\underset{z \in \mathbb{C}_{-}}{z \to 0}} \! \mi \mathbb{E}^{\sigma_{
3}} \sigma_{2} \! + \! \mathcal{O}(z)$; and (iv)
$\overset{o}{\mathfrak{M}}^{\raise-1.0ex\hbox{$\scriptstyle \infty$}}(z) \!
=_{\underset{z \in \mathbb{C} \setminus \widetilde{J}^{o}}{z \to \infty}} \!
\mathcal{O}(1)$. The solution of this latter (twisted) RHP for
$\overset{o}{\mathfrak{M}}^{\raise-1.0ex\hbox{$\scriptstyle \infty$}}(z)$ is
constructed out of the solution of two, simpler RHPs: $(\mathscr{N}^{o}(z),
-\mi \sigma_{2},\widetilde{J}^{o})$ and
$(\overset{o}{\mathfrak{m}}^{\raise-1.0ex\hbox{$\scriptstyle \infty$}}(z),
\overset{o}{\mathscr{U}}^{\raise-1.0ex\hbox{$\scriptstyle \infty$}}(z),
\widetilde{J}^{o})$, where
$\overset{o}{\mathscr{U}}^{\raise-1.0ex\hbox{$\scriptstyle \infty$}}(z)$
equals $\exp (\mi (n \! + \! 1/2) \Omega_{j}^{o} \sigma_{3}) \sigma_{1}$ for
$z \! \in \! (a_{j}^{o},b_{j}^{o})$, $j \! = \! 1,\dotsc,N$, and equals
$\mathrm{I}$ for $z \! \in \! (-\infty,b_{0}^{o}) \cup (a_{N+1}^{o},+
\infty)$. The RHP $(\mathscr{N}^{o}(z),-\mi \sigma_{2},\widetilde{J}^{o})$ 
is solved explicitly by diagonalising the jump matrix, and thus reduced to 
two scalar RHPs \cite{a2} (see, also, \cite{a45,a47,a79}): the solution is
\begin{equation*}
\mathscr{N}^{o}(z) \! = \!
\begin{pmatrix}
\frac{1}{2}((\gamma^{o}(0))^{-1} \gamma^{o}(z)+\gamma^{o}(0)(\gamma^{o}(z))^{-
1}) & -\frac{1}{2 \mi}((\gamma^{o}(0))^{-1} \gamma^{o}(z)-\gamma^{o}(0)
(\gamma^{o}(z))^{-1}) \\
\frac{1}{2 \mi}((\gamma^{o}(0))^{-1} \gamma^{o}(z)-\gamma^{o}(0)(\gamma^{o}
(z))^{-1}) & \frac{1}{2}((\gamma^{o}(0))^{-1} \gamma^{o}(z)+\gamma^{o}(0)
(\gamma^{o}(z))^{-1})
\end{pmatrix},
\end{equation*}
where $\gamma^{o} \colon \mathbb{C} \setminus \widetilde{J}^{o} \! \to \!
\mathbb{C}$ is characterised completely in Lemma~4.4; furthermore, $\mathscr{
N}^{o}(z)$ is piecewise (sectionally) holomorphic for $z \! \in \! \mathbb{C}
\setminus \widetilde{J}^{o}$, and $\mathscr{N}^{o}(z) \! =_{\underset{z \in
\mathbb{C}_{+}}{z \to 0}} \! \mathrm{I} \! + \! \mathcal{O}(z)$ and $\mathscr{
N}^{o}(z) \! =_{\underset{z \in \mathbb{C}_{-}}{z \to 0}} \! \mi \sigma_{2} \!
+ \! \mathcal{O}(z)$\footnote{Note that, strictly speaking, $\mathscr{N}^{o}
(z)$, as given above, does not solve the RHP $(\mathscr{N}^{o}(z),-\mi
\sigma_{2},\widetilde{J}^{o})$ in the sense defined heretofore, as $\mathscr{
N}^{o} \! \! \upharpoonright_{\mathbb{C}_{\pm}}$ can not be extended
continuously to $\overline{\mathbb{C}}_{\pm}$; however, $\mathscr{N}^{o}
(\pmb{\cdot} \! \pm \! \mi \varepsilon)$ converge in $\mathcal{L}^{2}_{
\mathrm{M}_{2}(\mathbb{C}),\mathrm{loc}}(\mathbb{R})$ as $\varepsilon \!
\downarrow \! 0$ to $\mathrm{SL}_{2}(\mathbb{C})$-valued functions $\mathscr{
N}^{o}(z)$ in $\mathcal{L}^{2}_{\mathrm{M}_{2}(\mathbb{C})}(\widetilde{J}^{o}
)$ that satisfy $\mathscr{N}^{o}_{+}(z) \! = \! \mathscr{N}^{o}_{-}(z)(-\mi
\sigma_{2})$ $\mathrm{a.e.}$ on $\widetilde{J}^{o}$: one then shows that
$\mathscr{N}^{o}(z)$ is the unique solution of the RHP $(\mathscr{N}^{o}(z),
-\mi \sigma_{2},\widetilde{J}^{o})$, where the latter boundary/jump condition
is interpreted in the $\mathcal{L}^{2}_{\mathrm{M}_{2}(\mathbb{C}),\mathrm{loc}
}$ sense.}.

Consider, now, the functions $\boldsymbol{\theta}^{o}(\boldsymbol{u}^{o}(z) \!
\pm \! \boldsymbol{d}_{o})$, where $\boldsymbol{u}^{o}(z) \colon z \! \to \!
\operatorname{Jac}(\mathcal{Y}_{o})$, $z \! \mapsto \! \boldsymbol{u}^{o}(z)
\! := \! \int_{a_{N+1}^{o}}^{z} \boldsymbol{\omega}^{o}$, with $\boldsymbol{
\omega}^{o}$ the associated normalised basis of holomorphic one-forms of
$\mathcal{Y}_{o}$, $\boldsymbol{d}_{o} \! \equiv \! -\sum_{j=1}^{N} \int_{a_{
j}^{o}}^{z_{j}^{o,-}} \boldsymbol{\omega}^{o} \! = \! \sum_{j=1}^{N} \int_{a_{
j}^{o}}^{z_{j}^{o,+}} \boldsymbol{\omega}^{o}$, where $\equiv$ denotes
equivalence modulo the period lattice, and $\lbrace z_{j}^{o,\pm} \rbrace_{j=
1}^{N}$ are characterised completely in Lemma~4.4. {}From the general theory
of theta functions on Riemann surfaces (see, for example, \cite{a77,a78}),
$\boldsymbol{\theta}^{o}(\boldsymbol{u}^{o}(z) \! + \! \boldsymbol{d}_{o})$,
for $z \! \in \! \mathcal{Y}_{o} \! := \! \lbrace \mathstrut (y,z); \, y^{2}
\! = \! \prod_{k=1}^{N+1}(z \! - \! b_{k-1}^{o})(z \! - \! a_{k}^{o})
\rbrace$, is either identically zero on $\mathcal{Y}_{o}$ or has precisely 
$N$ (simple) zeros (the generic case). In this case, since the divisors 
$\prod_{j=1}^{N}z_{j}^{o,-}$ and $\prod_{j=1}^{N}z_{j}^{o,+}$ are non-special, 
one uses Lemma~3.27 of \cite{a45} (see, also, Lemma~4.2 of \cite{a46}) and 
the representation \cite{a78} $\boldsymbol{K}_{o} \! = \! \sum_{j=1}^{N} 
\int_{a_{j}^{o}}^{a_{N+1}^{o}} \boldsymbol{\omega}^{o}$, for the `odd' 
vector of Riemann constants, with $2 \boldsymbol{K}_{o} \! = \! 0$ and $s 
\boldsymbol{K}_{o} \! \not= \! 0$, $0 \! < \! s \! < \! 2$, to arrive at
\begin{align*}
\boldsymbol{\theta}^{o}(\boldsymbol{u}^{o}(z) \! + \! \boldsymbol{d}_{o}) =&
\, \boldsymbol{\theta}^{o} \! \left(\boldsymbol{u}^{o}(z) \! - \! \sum_{j=1}^{
N} \int_{a_{j}^{o}}^{z_{j}^{o,-}} \boldsymbol{\omega}^{o} \right) \! = \!
\boldsymbol{\theta}^{o} \! \left(\int_{a_{N+1}^{o}}^{z} \boldsymbol{\omega}^{
o} \! - \! \boldsymbol{K}_{o} \! - \! \sum_{j=1}^{N} \int_{a_{N+1}^{o}}^{z_{
j}^{o,-}} \boldsymbol{\omega}^{o} \right) \! = \! 0 \\
\Leftrightarrow& \, z \! \in \! \left\{z_{1}^{o,-},z_{2}^{o,-},\dotsc,
z_{N}^{o,-} \right\}, \\
\boldsymbol{\theta}^{o}(\boldsymbol{u}^{o}(z) \! - \! \boldsymbol{d}_{o}) =&
\, \boldsymbol{\theta}^{o} \! \left(\boldsymbol{u}^{o}(z) \! - \! \sum_{j=1}^{
N} \int_{a_{j}^{o}}^{z_{j}^{o,+}} \boldsymbol{\omega}^{o} \right) \! = \!
\boldsymbol{\theta}^{o} \! \left(\int_{a_{N+1}^{o}}^{z} \boldsymbol{\omega}^{
o} \! - \! \boldsymbol{K}_{o} \! - \! \sum_{j=1}^{N} \int_{a_{N+1}^{o}}^{z_{
j}^{o,+}} \boldsymbol{\omega}^{o} \right) \! = \! 0 \\
\Leftrightarrow& \, z \! \in \! \left\{z_{1}^{o,+},z_{2}^{o,+},\dotsc,
z_{N}^{o,+} \right\}.
\end{align*}
Following Lemma~3.21 of \cite{a45}, set
\begin{equation*}
\overset{o}{\mathfrak{m}}^{\raise-1.0ex\hbox{$\scriptstyle \infty$}}(z) \! :=
\!
\begin{pmatrix}
\frac{\boldsymbol{\theta}^{o}(\boldsymbol{u}^{o}(z)-\frac{1}{2 \pi}(n+\frac{
1}{2}) \boldsymbol{\Omega}^{o}+\boldsymbol{d}_{o})}{\boldsymbol{\theta}^{o}
(\boldsymbol{u}^{o}(z)+\boldsymbol{d}_{o})} &
\frac{\boldsymbol{\theta}^{o}(-\boldsymbol{u}^{o}(z)-\frac{1}{2 \pi}(n+\frac{
1}{2}) \boldsymbol{\Omega}^{o}+\boldsymbol{d}_{o})}{\boldsymbol{\theta}^{o}(-
\boldsymbol{u}^{o}(z)+\boldsymbol{d}_{o})} \\
\frac{\boldsymbol{\theta}^{o}(\boldsymbol{u}^{o}(z)-\frac{1}{2 \pi}(n+\frac{
1}{2}) \boldsymbol{\Omega}^{o}-\boldsymbol{d}_{o})}{\boldsymbol{\theta}^{o}
(\boldsymbol{u}^{o}(z)-\boldsymbol{d}_{o})} &
\frac{\boldsymbol{\theta}^{o}(-\boldsymbol{u}^{o}(z)-\frac{1}{2 \pi}(n+\frac{
1}{2}) \boldsymbol{\Omega}^{o}-\boldsymbol{d}_{o})}{\boldsymbol{\theta}^{o}
(\boldsymbol{u}^{o}(z)+\boldsymbol{d}_{o})}
\end{pmatrix},
\end{equation*}
where $\boldsymbol{\Omega}^{o} \! := \! (\Omega_{1}^{o},\Omega_{2}^{o},\dotsc,
\Omega_{N}^{o})^{\mathrm{T}}$ $(\in \! \mathbb{R}^{N})$, with $\Omega_{j}^{
o}$, $j \! = \! 1,\dotsc,N$, given above, and ${}^{\mathrm{T}}$ denoting
transposition. Using Lemma~3.18 of \cite{a45} (or, equivalently, 
Equations~(4.65) and~(4.66) of \cite{a46}), that is, for $z \! \in \! (a_{j}^{
o},b_{j}^{o})$, $j \! = \! 1,\dotsc,N$, $\boldsymbol{u}^{o}_{+}(z) \! + \!
\boldsymbol{u}^{o}_{-}(z) \! \equiv \! -\tau^{o}_{j}$ $(:= \! -\tau^{o}e_{j}
)$, $j \! = \! 1,\dotsc,N$, with $\tau^{o} \! := \! (\tau^{o})_{i,j=1,\dotsc,
N} \! := \! (\oint_{\boldsymbol{\beta}_{j}^{o}} \omega^{o}_{i})_{i,j=1,\dotsc,
N}$ (the associated matrix of Riemann periods), and, for $z \! \in \! (-\infty,
b_{0}^{o}) \cup (a_{N+1}^{o},+\infty)$, $\boldsymbol{u}^{o}_{+}(z) \! + \!
\boldsymbol{u}^{o}_{-}(z) \! \equiv \! 0$, where $\boldsymbol{u}^{o}_{\pm}(z)
\! := \! \int_{a_{N+1}^{o}}^{z^{\pm}} \boldsymbol{\omega}^{o}$, with $z^{\pm}
\! \in \! (a_{j}^{o},b_{j}^{o})^{\pm}$, $j \! = \! 1,\dotsc,N$, and the
evenness and (quasi-) periodicity properties of $\boldsymbol{\theta}^{o}$, one
shows that, for $z \! \in \! (a_{j}^{o},b_{j}^{o})$, $j \! = \! 1,\dotsc,N$,
\begin{gather*}
\dfrac{\boldsymbol{\theta}^{o}(\boldsymbol{u}^{o}_{+}(z) \! - \! \frac{1}{2
\pi}(n \! + \! \frac{1}{2}) \boldsymbol{\Omega}^{o} \! + \! \boldsymbol{d}_{o}
)}{\boldsymbol{\theta}^{o}(\boldsymbol{u}^{o}_{+}(z) \! + \! \boldsymbol{d}_{
o})} \! = \! \me^{-\mi (n+\frac{1}{2}) \Omega_{j}^{o}} \dfrac{\boldsymbol{
\theta}^{o}(-\boldsymbol{u}^{o}_{-}(z) \! - \! \frac{1}{2 \pi}(n \! + \! \frac{
1}{2}) \boldsymbol{\Omega}^{o} \! + \! \boldsymbol{d}_{o})}{\boldsymbol{
\theta}^{o}(-\boldsymbol{u}^{o}_{-}(z) \! + \! \boldsymbol{d}_{o})}, \\
\dfrac{\boldsymbol{\theta}^{o}(\boldsymbol{u}^{o}_{+}(z) \! - \! \frac{1}{2
\pi}(n \! + \! \frac{1}{2}) \boldsymbol{\Omega}^{o} \! - \! \boldsymbol{d}_{o}
)}{\boldsymbol{\theta}^{o}(\boldsymbol{u}^{o}_{+}(z) \! - \! \boldsymbol{d}_{
o})} \! = \! \me^{-\mi (n+\frac{1}{2}) \Omega_{j}^{o}} \dfrac{\boldsymbol{
\theta}^{o}(-\boldsymbol{u}^{o}_{-}(z) \! - \! \frac{1}{2 \pi}(n \! + \! \frac{
1}{2}) \boldsymbol{\Omega}^{o} \! - \! \boldsymbol{d}_{o})}{\boldsymbol{
\theta}^{o}(\boldsymbol{u}^{o}_{-}(z) \! + \! \boldsymbol{d}_{o})}, \\
\dfrac{\boldsymbol{\theta}^{o}(-\boldsymbol{u}^{o}_{+}(z) \! - \! \frac{1}{2
\pi}(n \! + \! \frac{1}{2}) \boldsymbol{\Omega}^{o} \! + \! \boldsymbol{d}_{o}
)}{\boldsymbol{\theta}^{o}(-\boldsymbol{u}^{o}_{+}(z) \! + \! \boldsymbol{d}_{
o})} \! = \! \me^{\mi (n+\frac{1}{2}) \Omega_{j}^{o}} \dfrac{\boldsymbol{
\theta}^{o}(\boldsymbol{u}^{o}_{-}(z) \! - \! \frac{1}{2 \pi}(n \! + \! \frac{
1}{2}) \boldsymbol{\Omega}^{o} \! + \! \boldsymbol{d}_{o})}{\boldsymbol{
\theta}^{o}(\boldsymbol{u}^{o}_{-}(z) \! + \! \boldsymbol{d}_{o})}, \\
\dfrac{\boldsymbol{\theta}^{o}(-\boldsymbol{u}^{o}_{+}(z) \! - \! \frac{1}{2
\pi}(n \! + \! \frac{1}{2}) \boldsymbol{\Omega}^{o} \! - \! \boldsymbol{d}_{o}
)}{\boldsymbol{\theta}^{o}(\boldsymbol{u}^{o}_{+}(z) \! + \! \boldsymbol{d}_{
o})} \! = \! \me^{\mi (n+\frac{1}{2}) \Omega_{j}^{o}} \dfrac{\boldsymbol{
\theta}^{o}(\boldsymbol{u}^{o}_{-}(z) \! - \! \frac{1}{2 \pi}(n \! + \! \frac{
1}{2}) \boldsymbol{\Omega}^{o} \! - \! \boldsymbol{d}_{o})}{\boldsymbol{
\theta}^{o}(\boldsymbol{u}^{o}_{-}(z) \! - \! \boldsymbol{d}_{o})},
\end{gather*}
and, for $z \! \in \! (-\infty,b_{0}^{o}) \cup (a_{N+1}^{o},+\infty)$, one 
obtains the same relations as above but with the changes $\exp (\mp \mi (n 
\! + \! 1/2) \Omega_{j}^{o}) \! \to \! 1$. Set, as in Proposition~3.31 of 
\cite{a45},
\begin{equation*}
\overset{o}{\mathcal{Q}}^{\raise-1.0ex\hbox{$\scriptstyle \infty$}}(z) \! :=
\!
\begin{pmatrix}
(\mathscr{N}^{o}(z))_{11}
(\overset{o}{\mathfrak{m}}^{\raise-1.0ex\hbox{$\scriptstyle \infty$}}(z))_{11}
& (\mathscr{N}^{o}(z))_{12}
(\overset{o}{\mathfrak{m}}^{\raise-1.0ex\hbox{$\scriptstyle \infty$}}(z))_{12}
\\(\mathscr{N}^{o}(z))_{21}
(\overset{o}{\mathfrak{m}}^{\raise-1.0ex\hbox{$\scriptstyle \infty$}}(z))_{21}
& (\mathscr{N}^{o}(z))_{22}
(\overset{o}{\mathfrak{m}}^{\raise-1.0ex\hbox{$\scriptstyle \infty$}}(z))_{22}
\end{pmatrix},
\end{equation*}
where $(\ast)_{ij}$, $i,j \! = \! 1,2$, denotes the $(i \, j)$-element of
$(\ast)$. Recalling that $\mathscr{N}^{o} \colon \mathbb{C} \setminus
\widetilde{J}^{o} \! \to \! \mathrm{SL}_{2}(\mathbb{C})$ solves the RHP
$(\mathscr{N}^{o}(z),-\mi \sigma_{2},\widetilde{J}^{o})$, using the above
theta-functional relations and the small-$z$ asymptotic expansion of
$\boldsymbol{u}^{o}(z)$ (see Section~5, the proof of Proposition~5.3), one
shows that
$\overset{o}{\mathcal{Q}}^{\raise-1.0ex\hbox{$\scriptstyle \infty$}}(z)$
solves the following RHP: (i)
$\overset{o}{\mathcal{Q}}^{\raise-1.0ex\hbox{$\scriptstyle \infty$}}(z)$ is
holomorphic for $z \! \in \! \mathbb{C} \setminus \widetilde{J}^{o}$; (ii)
$\overset{o}{\mathcal{Q}}^{\raise-1.0ex\hbox{$\scriptstyle \infty$}}_{\pm}
(z)\! := \! \lim_{\underset{z^{\prime} \, \in \, \pm \, \mathrm{side} \, 
\mathrm{of} \, \widetilde{J}^{o}}{z^{\prime} \to z}} 
\overset{o}{\mathcal{Q}}^{\raise-1.0ex\hbox{$\scriptstyle \infty$}}
(z^{\prime})$ satisfy the boundary condition
$\overset{o}{\mathcal{Q}}^{\raise-1.0ex\hbox{$\scriptstyle \infty$}}_{+}(z) \!
= \! \overset{o}{\mathcal{Q}}^{\raise-1.0ex\hbox{$\scriptstyle \infty$}}_{-}
(z) \overset{o}{\mathscr{V}}^{\raise-1.0ex\hbox{$\scriptstyle \infty$}}(z)$,
$z \! \in \! \widetilde{J}^{o}$, where
$\overset{o}{\mathscr{V}}^{\raise-1.0ex\hbox{$\scriptstyle \infty$}}(z)$ is
defined in Equation~(4.1); (iii)
\begin{align*}
\overset{o}{\mathcal{Q}}^{\raise-1.0ex\hbox{$\scriptstyle \infty$}}(z)
\underset{\underset{z \in \mathbb{C}_{+}}{z \to 0}}{=}&
\begin{pmatrix}
\frac{\boldsymbol{\theta}^{o}(\boldsymbol{u}^{o}_{+}(0)-\frac{1}{2 \pi}(n+
\frac{1}{2}) \boldsymbol{\Omega}^{o}+\boldsymbol{d}_{o})}{\boldsymbol{\theta}^{
o}(\boldsymbol{u}^{o}_{+}(0)+\boldsymbol{d}_{o})} & 0 \\
0 & \frac{\boldsymbol{\theta}^{o}(-\boldsymbol{u}^{o}_{+}(0)-\frac{1}{2 \pi}(n
+\frac{1}{2}) \boldsymbol{\Omega}^{o}-\boldsymbol{d}_{o})}{\boldsymbol{
\theta}^{o}(\boldsymbol{u}^{o}_{+}(0)+\boldsymbol{d}_{o})}
\end{pmatrix} \! + \! \mathcal{O}(z), \\
\overset{o}{\mathcal{Q}}^{\raise-1.0ex\hbox{$\scriptstyle \infty$}}(z)
\underset{\underset{z \in \mathbb{C}_{-}}{z \to 0}}{=}&
\begin{pmatrix}
0 & \frac{\boldsymbol{\theta}^{o}(-\boldsymbol{u}^{o}_{-}(0)-\frac{1}{2 \pi}(n
+\frac{1}{2}) \boldsymbol{\Omega}^{o}+\boldsymbol{d}_{o})}{\boldsymbol{
\theta}^{o}(-\boldsymbol{u}^{o}_{-}(0)+\boldsymbol{d}_{o})} \\
-\frac{\boldsymbol{\theta}^{o}(\boldsymbol{u}^{o}_{-}(0)-\frac{1}{2 \pi}(n+
\frac{1}{2}) \boldsymbol{\Omega}^{o}-\boldsymbol{d}_{o})}{\boldsymbol{\theta}^{
o}(\boldsymbol{u}^{o}_{-}(0)-\boldsymbol{d}_{o})} & 0
\end{pmatrix} \! + \! \mathcal{O}(z),
\end{align*}
where $\boldsymbol{u}^{o}_{\pm}(0) \! := \! \int_{a_{N+1}^{o}}^{0^{\pm}}
\boldsymbol{\omega}^{o}$; and (iv)
$\overset{o}{\mathcal{Q}}^{\raise-1.0ex\hbox{$\scriptstyle \infty$}}(z) \!
=_{\underset{z \in \mathbb{C} \setminus \widetilde{J}^{o}}{z \to \infty}} \!
\mathcal{O}(1)$. Now, using the fact that, for $z \! \in \! (a_{j}^{o},b_{j}^{
o})$, $j \! = \! 1,\dotsc,N$, $\boldsymbol{u}^{o}_{+}(0) \! + \! \boldsymbol{
u}^{o}_{-}(0) \! = \! \int_{a_{N+1}^{o}}^{0^{+}} \boldsymbol{\omega}^{o} \! +
\! \int_{a_{N+1}^{o}}^{0^{-}} \boldsymbol{\omega}^{o} \! \equiv \! -\tau_{j}^{
o}$, $j \! = \! 1,\dotsc,N$, and, for $z \! \in \! (-\infty,b_{0}^{o}) \cup
(a_{N+1}^{o},+\infty)$, $\boldsymbol{u}^{o}_{+}(0) \! + \! \boldsymbol{u}^{o}_{
-}(0) \! = \! \int_{a_{N+1}^{o}}^{0^{+}} \boldsymbol{\omega}^{o} \! + \! \int_{
a_{N+1}^{o}}^{0^{-}} \boldsymbol{\omega}^{o} \! \equiv \! 0$, upon multiplying
$\overset{o}{\mathcal{Q}}^{\raise-1.0ex\hbox{$\scriptstyle \infty$}}(z)$ on
the left by
\begin{equation*}
\operatorname{diag} \! \left(\! \dfrac{\boldsymbol{\theta}^{o}(\boldsymbol{u}^{
o}_{+}(0) \! + \! \boldsymbol{d}_{o}) \mathbb{E}^{-1}}{\boldsymbol{\theta}^{o}
(\boldsymbol{u}^{o}_{+}(0) \! - \! \frac{1}{2 \pi}(n \! + \! \frac{1}{2})
\boldsymbol{\Omega}^{o} \! + \! \boldsymbol{d}_{o})},\dfrac{\boldsymbol{
\theta}^{o}(\boldsymbol{u}^{o}_{+}(0) \! + \! \boldsymbol{d}_{o}) \mathbb{E}}{
\boldsymbol{\theta}^{o}(-\boldsymbol{u}^{o}_{+}(0) \! - \! \frac{1}{2 \pi}(n
\! + \! \frac{1}{2}) \boldsymbol{\Omega}^{o} \! - \! \boldsymbol{d}_{o})}
\right) \! =: \!
\overset{o}{\mathfrak{c}}^{\raise-1.0ex\hbox{$\scriptstyle \infty$}},
\end{equation*}
that is, $\overset{o}{\mathcal{Q}}^{\raise-1.0ex\hbox{$\scriptstyle \infty$}}
(z) \! \to \! \overset{o}{\mathfrak{c}}^{\raise-1.0ex\hbox{$\scriptstyle
\infty$}} \overset{o}{\mathcal{Q}}^{\raise-1.0ex\hbox{$\scriptstyle \infty$}}
(z) \! =: \! \mathcal{M}^{\infty}_{o}(z)$, one shows that $\mathcal{M}^{
\infty}_{o} \colon \mathbb{C} \setminus \widetilde{J}^{o} \! \to \! 
\mathrm{SL}_{2}(\mathbb{C})$ solves the RHP $(\mathcal{M}^{\infty}_{o}(z),
\overset{o}{\mathscr{V}}^{\raise-1.0ex\hbox{$\scriptstyle \infty$}}(z),
\linebreak[4]
\widetilde{J}^{o})$. Using, finally, the formula 
$\overset{o}{m}^{\raise-1.0ex\hbox{$\scriptstyle \infty$}}(z) \! = \!
\begin{cases}
\overset{o}{\mathfrak{M}}^{\raise-1.0ex\hbox{$\scriptstyle \infty$}}(z),
&\text{$z \! \in \! \mathbb{C}_{+}$,} \\
-\mi \, \overset{o}{\mathfrak{M}}^{\raise-1.0ex\hbox{$\scriptstyle \infty$}}
(z) \sigma_{2}, &\text{$z \! \in \! \mathbb{C}_{-}$,}
\end{cases}$ one shows that
$\overset{o}{m}^{\raise-1.0ex\hbox{$\scriptstyle \infty$}} \colon \mathbb{C}
\setminus J_{o}^{\infty} \! \to \! \operatorname{SL}_{2}(\mathbb{C})$ solves
the model RHP formulated in Lemma~4.3. One notes {}from the formula for
$\overset{o}{\mathfrak{M}}^{\raise-1.0ex\hbox{$\scriptstyle \infty$}}(z)$
stated in the Lemma that it is well defined for $\mathbb{C} \setminus \mathbb{
R}$; in particular, it is single valued and analytic (see below) for $z \!
\in \! \mathbb{C} \setminus \widetilde{J}^{o}$ (independently of the path in
$\mathbb{C} \setminus \widetilde{J}^{o}$ chosen to evaluate $\boldsymbol{u}^{
o}(z) \! = \! \int_{a_{N+1}^{o}}^{z} \boldsymbol{\omega}^{o})$. Furthermore
(cf. Lemma~4.4 and the analysis above), since $\{\mathstrut z^{\prime} \! \in
\! \mathbb{C}; \, \boldsymbol{\theta}^{o}(\boldsymbol{u}^{o}(z^{\prime}) \!
\pm \! \boldsymbol{d}_{o}) \! = \! 0 \rbrace \! = \! \lbrace z_{j}^{o,\mp}
\rbrace_{j=1}^{N} \! = \! \lbrace \mathstrut z^{\prime} \! \in \! \mathbb{C};
\, ((\gamma^{o}(0))^{-1} \gamma^{o}(z) \! \pm \! \gamma^{o}(0)(\gamma^{o}(z)
)^{-1}) \vert_{z=z^{\prime}} \! = \! 0 \rbrace$, one notes that the (simple)
poles of
$(\overset{o}{\mathfrak{m}}^{\raise-1.0ex\hbox{$\scriptstyle \infty$}}(z))_{1
1}$ and $(\overset{o}{\mathfrak{m}}^{\raise-1.0ex\hbox{$\scriptstyle \infty$}}
(z))_{22}$ (resp.,
$(\overset{o}{\mathfrak{m}}^{\raise-1.0ex\hbox{$\scriptstyle \infty$}}(z))_{1
2}$ and
$(\overset{o}{\mathfrak{m}}^{\raise-1.0ex\hbox{$\scriptstyle \infty$}}(z))_{
21})$, that is, $\lbrace \mathstrut z^{\prime} \! \in \! \mathbb{C}; \,
\boldsymbol{\theta}^{o}(\boldsymbol{u}^{o}(z^{\prime}) \! + \! \boldsymbol{
d}_{o}) \! = \! 0 \rbrace$ (resp., $\lbrace \mathstrut z^{\prime} \! \in \!
\mathbb{C}; \, \boldsymbol{\theta}^{o}(\boldsymbol{u}^{o}(z^{\prime}) \! - \!
\boldsymbol{d}_{o}) \! = \! 0 \rbrace)$, are exactly cancelled by the (simple)
zeros of $(\gamma^{o}(0))^{-1} \gamma^{o}(z) \! + \! \gamma^{o}(0)(\gamma^{o}
(z))^{-1}$ (resp., $(\gamma^{o}(0))^{-1} \gamma^{o}(z) \! - \! \gamma^{o}(0)
(\gamma^{o}(z))^{-1})$; thus,
$\overset{o}{\mathfrak{M}}^{\raise-1.0ex\hbox{$\scriptstyle \infty$}}(z)$ has
only $\tfrac{1}{4}$-root singularities at the end-points of the support of the
`odd' equilibrium measure, $\lbrace b_{j-1}^{o},a_{j}^{o} \rbrace_{j=1}^{N+
1}$. (This shows that
$\overset{o}{\mathfrak{M}}^{\raise-1.0ex\hbox{$\scriptstyle \infty$}}(z)$
obtains its boundary values,
$\overset{o}{\mathfrak{M}}^{\raise-1.0ex\hbox{$\scriptstyle \infty$}}_{\pm}(z)
\! := \! \lim_{\varepsilon \downarrow 0}
\overset{o}{\mathfrak{M}}^{\raise-1.0ex\hbox{$\scriptstyle \infty$}}(z \! \pm
\! \mi \varepsilon)$, in the $\mathcal{L}^{2}_{\mathrm{M}_{2}(\mathbb{C})}
(\mathbb{R})$ sense.) {}From the definition of
$\overset{o}{m}^{\raise-1.0ex\hbox{$\scriptstyle \infty$}}(z)$ in terms of
$\overset{o}{\mathfrak{M}}^{\raise-1.0ex\hbox{$\scriptstyle \infty$}}(z)$
given in the Lemma, the explicit formula for
$\overset{o}{\mathfrak{M}}^{\raise-1.0ex\hbox{$\scriptstyle \infty$}}(z)$, and
recalling that $\overset{o}{m}^{\raise-1.0ex\hbox{$\scriptstyle \infty$}}(z)$
solves the model RHP formulated in Lemma~4.3, one learns that, as $\det
(\overset{o}{\upsilon}^{\raise-1.0ex\hbox{$\scriptstyle \infty$}}(z)) \! = \!
1$, $\det (\overset{o}{m}^{\raise-1.0ex\hbox{$\scriptstyle \infty$}}_{+}(z))
\! = \! \det (\overset{o}{m}^{\raise-1.0ex\hbox{$\scriptstyle \infty$}}_{-}
(z))$, that is,
$\det (\overset{o}{m}^{\raise-1.0ex\hbox{$\scriptstyle \infty$}}(z))$ has no
`jumps', whence
$\det (\overset{o}{m}^{\raise-1.0ex\hbox{$\scriptstyle \infty$}}(z))$ has, at
worst, (isolated) $\tfrac{1}{2}$-root singularities at $\lbrace b_{j-1}^{o},
a_{j}^{o} \rbrace_{j=1}^{N+1}$, which are removable, which implies that $\det
(\overset{o}{m}^{\raise-1.0ex\hbox{$\scriptstyle \infty$}}(z))$ is entire
and bounded; hence, via a generalisation of Liouville's Theorem, and the
asymptotic relation $\det
(\overset{o}{m}^{\raise-1.0ex\hbox{$\scriptstyle \infty$}}(z)) \! =_{
\underset{z \in \mathbb{C} \setminus \mathbb{R}}{z \to 0}} \! 1 \! + \!
\mathcal{O}(z)$, one arrives at
$\det (\overset{o}{m}^{\raise-1.0ex\hbox{$\scriptstyle \infty$}}(z)) \! = \!
1 \! \Rightarrow \! \overset{o}{m}^{\raise-1.0ex\hbox{$\scriptstyle \infty$}}
\! \in \! \operatorname{SL}_{2}(\mathbb{C})$. Also, {}from the definition of
$\overset{o}{m}^{\raise-1.0ex\hbox{$\scriptstyle \infty$}}(z)$ in terms of
$\overset{o}{\mathfrak{M}}^{\raise-1.0ex\hbox{$\scriptstyle \infty$}}(z)$ 
and the explicit formula for
$\overset{o}{\mathfrak{M}}^{\raise-1.0ex\hbox{$\scriptstyle \infty$}}(z)$, it
follows that both $\overset{o}{m}^{\raise-1.0ex\hbox{$\scriptstyle \infty$}}
(z)$ and $(\overset{o}{m}^{\raise-1.0ex\hbox{$\scriptstyle \infty$}}(z))^{-1}$
are uniformly bounded as functions of $n$ (as $n \! \to \! \infty)$ for $z$ in
compact subsets away {}from $\lbrace b_{j-1}^{o},a_{j}^{o} \rbrace_{j=1}^{N+
1}$.

Let $\mathscr{S}^{\infty}_{o} \colon \mathbb{C} \setminus \mathbb{R} \! 
\to \! \operatorname{SL}_{2}(\mathbb{C})$ be another solution of the RHP
$(\overset{o}{\mathfrak{M}}^{\raise-1.0ex\hbox{$\scriptstyle \infty$}}(z),
\overset{o}{\mathscr{V}}^{\raise-1.0ex\hbox{$\scriptstyle \infty$}}(z),
\mathbb{R})$ formulated at the beginning of the proof. Set $\Delta^{o}(z) \!
:= \! \mathscr{S}^{\infty}_{o}(z)
(\overset{o}{\mathfrak{M}}^{\raise-1.0ex\hbox{$\scriptstyle \infty$}}(z))^{
-1}$. Then $\Delta^{o}_{+}(z) \! = \! (\mathscr{S}^{\infty}_{o}(z))_{+}
(\overset{o}{\mathfrak{M}}^{\raise-1.0ex\hbox{$\scriptstyle \infty$}}_{+}
(z))^{-1} \! = \! (\mathscr{S}^{\infty}_{o}(z))_{-} \linebreak[4]
\cdot \overset{o}{\mathscr{V}}^{\raise-1.0ex\hbox{$\scriptstyle \infty$}}(z)
(\overset{o}{\mathfrak{M}}^{\raise-1.0ex\hbox{$\scriptstyle \infty$}}_{-}(z)
\overset{o}{\mathscr{V}}^{\raise-1.0ex\hbox{$\scriptstyle \infty$}}(z))^{-1}
\! = \! (\mathscr{S}^{\infty}_{o}(z))_{-}
(\overset{o}{\mathfrak{M}}^{\raise-1.0ex\hbox{$\scriptstyle \infty$}}_{-}(z)
)^{-1} \! = \! \Delta^{o}_{-}(z)$, hence $\Delta^{o}(z)$ is analytic across
$\mathbb{R}$; moreover, as $\det
(\overset{o}{\mathfrak{M}}^{\raise-1.0ex\hbox{$\scriptstyle \infty$}}(z)) \!
= \! 1$, it follows that $\Delta^{o}(z)$ has, at worst, $\mathcal{L}^{1}_{
\mathrm{M}_{2}(\mathbb{C})}(\ast)$-singularities at $b_{j-1}^{o},a_{j}^{o}$,
$j \! = \! 1,\dotsc,N \! + \! 1$, which, as per the discussion above, are
removable; thence, noting that $\Delta^{o}(z) \! \to \! \mathrm{I}$ as $z \!
\to \! 0$ $(z \! \in \! \mathbb{C} \setminus \mathbb{R})$, one concludes that
$\Delta^{o}(z) \! = \! \mathrm{I}$, whence $\mathscr{S}^{\infty}_{o}(z) \! =
\! \overset{o}{\mathfrak{M}}^{\raise-1.0ex\hbox{$\scriptstyle \infty$}}(z)$.
\hfill $\qed$

In order to prove that there is a solution of the (full) RHP 
$(\overset{o}{\mathscr{M}}^{\raise-1.0ex\hbox{$\scriptstyle \sharp$}}(z),
\overset{o}{\upsilon}^{\raise-1.0ex\hbox{$\scriptstyle \sharp$}}(z),\Sigma_{
o}^{\sharp})$, formulated in Lemma~4.2, close to the parametrix, one needs 
to know that the parametrix is \emph{uniformly} bounded: more precisely, by 
(certain) general theorems (see, for example, \cite{a86}), one needs to know 
that $\overset{o}{\upsilon}^{\raise-1.0ex\hbox{$\scriptstyle \sharp$}}(z) \! 
\to \! \overset{o}{\upsilon}^{\raise-1.0ex\hbox{$\scriptstyle \infty$}}(z)$ 
as $n \! \to \! \infty$ uniformly for $z \! \in \! \Sigma_{o}^{\sharp}$ in 
the $\mathcal{L}^{2}_{\mathrm{M}_{2}(\mathbb{C})}(\Sigma_{o}^{\sharp}) \cap 
\mathcal{L}^{\infty}_{\mathrm{M}_{2}(\mathbb{C})}(\Sigma_{o}^{\sharp})$ sense, 
that is, uniformly,
\begin{equation*}
\lim_{n \to \infty}
\norm{\overset{o}{\upsilon}^{\raise-1.0ex\hbox{$\scriptstyle \sharp$}}(\cdot)
\! - \! \overset{o}{\upsilon}^{\raise-1.0ex\hbox{$\scriptstyle \infty$}}
(\cdot)}_{\mathcal{L}^{2}_{\mathrm{M}_{2}(\mathbb{C})}(\Sigma_{o}^{\sharp})
\cap \mathcal{L}^{\infty}_{\mathrm{M}_{2}(\mathbb{C})}(\Sigma_{o}^{\sharp})}
\! := \! \lim_{n \to \infty} \sum_{p \in \{2,\infty\}} \norm{\overset{o}{
\upsilon}^{\raise-1.0ex\hbox{$\scriptstyle \sharp$}}(\cdot) \! - \! \overset{
o}{\upsilon}^{\raise-1.0ex\hbox{$\scriptstyle \infty$}}(\cdot)}_{\mathcal{L}^{
p}_{\mathrm{M}_{2}(\mathbb{C})}(\Sigma_{o}^{\sharp})} \! = \! 0;
\end{equation*}
however, notwithstanding the fact that $\widetilde{V} \colon \mathbb{R}
\setminus \{0\} \! \to \! \mathbb{R}$ is regular $(h_{V}^{o}(b_{j-1}^{o}),h_{
V}^{o}(a_{j}^{o}) \! \not= \! 0$, $j \! = \! 1,\dotsc,N \! + \! 1)$, since the
strict inequalities $g^{o}_{+}(z) \! + \! g^{o}_{-}(z) \! - \! \widetilde{V}
(z) \! - \! \ell_{o} \! - \! \mathfrak{Q}^{+}_{\mathscr{A}} \! - \! \mathfrak{
Q}^{-}_{\mathscr{A}} \! < \! 0$, $z \! \in \! (-\infty,b_{0}^{o}) \cup (a_{N+
1}^{o},+\infty) \cup (\cup_{j=1}^{N}(a_{j}^{o},b_{j}^{o}))$, and $\pm \Re
(\mi \int_{z}^{a_{N+1}^{o}} \psi_{V}^{o}(s) \, \md s) \! > \! 0$, $z \! \in
\! \mathbb{C}_{\pm} \cap (\cup_{j=1}^{N+1} \mathbb{U}_{j}^{o})$, fail at the
end-points of the support of the `odd' equilibrium measure, this implies that
$\overset{o}{\upsilon}^{\raise-1.0ex\hbox{$\scriptstyle \sharp$}}(z) \! \to \!
\overset{o}{\upsilon}^{\raise-1.0ex\hbox{$\scriptstyle \infty$}}(z)$ as $n \!
\to \! \infty$ \emph{pointwise}, but not uniformly, for $z \! \in \! \Sigma_{
o}^{\sharp}$, whence, one can not conclude that $\overset{o}{\mathscr{M}
}^{\raise-1.0ex\hbox{$\scriptstyle \sharp$}}(z) \! \to \! \overset{o}{m}^{
\raise-1.0ex\hbox{$\scriptstyle \infty$}}(z)$ as $n \! \to \! \infty$
uniformly for $z \! \in \! \Sigma_{o}^{\sharp}$. The resolution of this lack
of uniformity at the end-points of the support of the `odd' equilibrium
measure constitutes, therefore, the essential analytical obstacle remaining
for the analysis of the RHP
$(\overset{o}{\mathscr{M}}^{\raise-1.0ex\hbox{$\scriptstyle \sharp$}}(z),
\overset{o}{\upsilon}^{\raise-1.0ex\hbox{$\scriptstyle \sharp$}}(z),\Sigma_{
o}^{\sharp})$, and a substantial part of the following analysis is devoted to
overcoming this problem.

The key necessary to remedy (and control) the above-mentioned analytical 
difficulty is to construct parametrices for the solution of the RHP 
$(\overset{o}{\mathscr{M}}^{\raise-1.0ex\hbox{$\scriptstyle \sharp$}}(z),
\overset{o}{\upsilon}^{\raise-1.0ex\hbox{$\scriptstyle \sharp$}}(z),\Sigma_{
o}^{\sharp})$ in `small' neighbourhoods (open discs) about $\lbrace b_{j
-1}^{o},a_{j}^{o} \rbrace_{j=1}^{N+1}$ (where the convergence of $\overset{
o}{\upsilon}^{\raise-1.0ex\hbox{$\scriptstyle \sharp$}}(z)$ to 
$\overset{o}{\upsilon}^{\raise-1.0ex\hbox{$\scriptstyle \infty$}}(z)$ as $n 
\! \to \! \infty$ is not uniform) in such a way that, on the boundary of 
these neighbourhoods, the parametrices `match' with the solution of the model 
RHP, $\overset{o}{m}^{\raise-1.0ex\hbox{$\scriptstyle \infty$}}(z)$, up to $o
(1)$ (in fact, $\mathcal{O}((n \! + \! 1/2)^{-1}))$ as $n \! \to \! \infty$; 
furthermore, in the generic framework considered in this work, namely, 
$\widetilde{V} \colon \mathbb{R} \setminus \{0\} \! \to \! \mathbb{R}$ is 
regular, in which case the (density of the) `odd' equilibrium measure behaves 
as a square root at the end-points of $\operatorname{supp}(\mu_{V}^{o})$, that 
is, $\psi_{V}^{o}(s) \! =_{s \downarrow b_{j-1}^{o}} \! \mathcal{O}((s \! - 
\! b_{j-1}^{o})^{1/2})$ and $\psi_{V}^{o}(s) \! =_{s \uparrow a_{j}^{o}} \! 
\mathcal{O}((a_{j}^{o} \! - \! s)^{1/2})$, $j \! = \! 1,\dotsc,N \! + \! 1$, 
it is well known \cite{a3,a47,a79} that the parametrices can be expressed 
in terms of Airy functions. (The general method used to construct such 
parametrices is via a Vanishing Lemma \cite{a87}.) More precisely, one 
surrounds the end-points of the support of the `odd' equilibrium measure, 
$\lbrace b_{j-1}^{o},a_{j}^{o} \rbrace_{j=1}^{N+1}$, by `small', mutually 
disjoint (open) discs,
\begin{equation*}
\mathbb{D}_{\epsilon}(b_{j-1}^{o}) \! := \! \left\{\mathstrut z \! \in \!
\mathbb{C}; \, \vert z \! - \! b_{j-1}^{o} \vert \! < \! \epsilon_{j}^{b}
\right\} \qquad \text{and} \qquad \mathbb{D}_{\epsilon}(a_{j}^{o}) \! := \!
\left\{\mathstrut z \! \in \! \mathbb{C}; \, \vert z \! - \! a_{j}^{o} \vert
\! < \! \epsilon_{j}^{a} \right\}, \quad j \! = \! 1,\dotsc,N \! + \! 1,
\end{equation*}
where $\epsilon_{j}^{b},\epsilon_{j}^{a}$ are arbitrarily fixed, sufficiently
small positive real numbers chosen so that $\mathbb{D}_{\epsilon}(b_{i-1}^{
o}) \cap \mathbb{D}_{\epsilon}(a_{j}^{o}) \! = \! \varnothing$, $i,j \! = 
\! 1,\dotsc,N \! + \! 1$, and defines $S_{p}^{o}(z)$, the parametrix for
$\overset{o}{\mathscr{M}}^{\raise-1.0ex\hbox{$\scriptstyle \sharp$}}(z)$, 
by $\overset{o}{m}^{\raise-1.0ex\hbox{$\scriptstyle \infty$}}(z)$ for $z \! 
\in \! \mathbb{C} \setminus (\cup_{j=1}^{N+1}(\mathbb{D}_{\epsilon}(b_{j-1}^{
o}) \cup \mathbb{D}_{\epsilon}(a_{j}^{o})))$, and by $m_{p}^{o}(z)$ for $z \! 
\in \! \cup_{j=1}^{N+1}(\mathbb{D}_{\epsilon}(b_{j-1}^{o}) \cup \mathbb{D}_{
\epsilon}(a_{j}^{o}))$, and solves the local RHP for $m_{p}^{o}(z)$ on 
$\cup_{j=1}^{N+1}(\mathbb{D}_{\epsilon}(b_{j-1}^{o}) \cup \mathbb{D}_{\epsilon}
(a_{j}^{o}))$ in such a way (`optimal', in the nomenclature of \cite{a47}) 
that $m_{p}^{o}(z) \! \approx_{n \to \infty}
\overset{o}{m}^{\raise-1.0ex\hbox{$\scriptstyle \infty$}}(z)$ (to $\mathcal{O}
((n \! + \! 1/2)^{-1}))$ for $z \! \in \! \cup_{j=1}^{N+1}(\partial \mathbb{
D}_{\epsilon}(b_{j-1}^{o}) \cup \partial \mathbb{D}_{\epsilon}(a_{j}^{o}))$,
whence $\mathcal{R}^{o}(z) \! := \!
\overset{o}{\mathscr{M}}^{\raise-1.0ex\hbox{$\scriptstyle \sharp$}}(z)(S_{p}^{
o}(z))^{-1} \colon \mathbb{C} \setminus \widetilde{\Sigma}_{o}^{\sharp} \!
\to \! \operatorname{SL}_{2}(\mathbb{C})$, where $\widetilde{\Sigma}_{o}^{
\sharp} \! := \! \Sigma_{o}^{\sharp} \cup (\cup_{j=1}^{N+1}(\partial \mathbb{
D}_{\epsilon}(b_{j-1}^{o}) \cup \partial \mathbb{D}_{\epsilon}(a_{j}^{o})))$,
solves the RHP $(\mathcal{R}^{o}(z),\upsilon^{o}_{\mathcal{R}}(z),\widetilde{
\Sigma}_{o}^{\sharp})$ with $\norm{\upsilon^{o}_{\mathcal{R}}(\cdot) \! - \!
\mathrm{I}}_{\cap_{p \in \{2,\infty\}} \mathcal{L}^{p}_{\mathrm{M}_{2}
(\mathbb{C})}(\widetilde{\Sigma}_{o}^{\sharp})} \! =_{n \to \infty} \!
\mathcal{O}((n \! + \! 1/2)^{-1})$ uniformly; in particular, the error term,
which is $\mathcal{O}((n \! + \! 1/2)^{-1})$ as $n \! \to \! \infty$, is
uniform in $\cap_{p \in \{1,2,\infty\}} \mathcal{L}^{p}_{\mathrm{M}_{2}
(\mathbb{C})}(\widetilde{\Sigma}_{o}^{\sharp})$. By general Riemann-Hilbert
techniques (see, for example, \cite{a86}), $\mathcal{R}^{o}(z)$ (and thus
$\overset{o}{\mathscr{M}}^{\raise-1.0ex\hbox{$\scriptstyle \sharp$}}(z)$
via the relation
$\overset{o}{\mathscr{M}}^{\raise-1.0ex\hbox{$\scriptstyle \sharp$}}(z) \! =
\! \mathcal{R}^{o}(z)S_{p}^{o}(z))$ can be computed to any order of $(n \! +
\! 1/2)^{-1}$ (as $n \! \to \! \infty)$ via a Neumann series expansion (of
the corresponding resolvent kernel). In fact, at the very core of the
above-mentioned discussion, and the analysis that follows, is the following
Corollary (see, for example, \cite{a79}, Corollary~7.108):
\begin{fffff}[Deift {\rm \cite{a79}}]
For an oriented contour $\varSigma \subset \mathbb{C}$, let $m^{\infty} \colon
\mathbb{C} \setminus \varSigma \! \to \! \operatorname{SL}_{2}(\mathbb{C})$
and $m^{(n)} \colon \mathbb{C} \setminus \varSigma \! \to \! \operatorname{S
L}_{2}(\mathbb{C})$, $n \! \in \! \mathbb{N}$, respectively, solve the
following, equivalent {\rm RHPs}, $(m^{\infty}(z),\upsilon^{\infty}(z),
\varSigma)$ and $(m^{(n)}(z),\upsilon^{(n)} \linebreak[4]
(z),\varSigma)$, where
\begin{equation*}
\upsilon^{\infty} \colon \varSigma \! \to \! \operatorname{GL}_{2}(\mathbb{C}
), \, \, z \! \mapsto \! \left(\mathrm{I} \! - \! w^{\infty}_{-}(z) \right)^{
-1} \! \left(\mathrm{I} \! + \! w^{\infty}_{+}(z) \right)
\end{equation*}
and
\begin{equation*}
\upsilon^{(n)} \colon \varSigma \! \to \! \operatorname{GL}_{2}(\mathbb{C}),
\, \, z \! \mapsto \! \left(\mathrm{I} \! - \! w^{(n)}_{-}(z) \right)^{-1} \!
\left(\mathrm{I} \! + \! w^{(n)}_{+}(z) \right),
\end{equation*}
and suppose that $(\id \! - \! C^{\infty}_{w^{\infty}})^{-1}$ exists, where
\begin{equation*}
\mathcal{L}^{2}_{\mathrm{M}_{2}(\mathbb{C})}(\varSigma) \! \ni \! f \! \mapsto
\! C^{\infty}_{w^{\infty}}f \! := \! C^{\infty}_{+}(fw^{\infty}_{-}) \! + \!
C^{\infty}_{-}(fw^{\infty}_{+}),
\end{equation*}
with
\begin{equation*}
C^{\infty}_{\pm} \colon \mathcal{L}^{2}_{\mathrm{M}_{2}(\mathbb{C})}
(\varSigma) \! \to \! \mathcal{L}^{2}_{\mathrm{M}_{2}(\mathbb{C})}(\varSigma
), \, \, f \! \mapsto \! (C^{\infty}_{\pm}f)(z) := \lim_{\underset{z^{\prime}
\, \in \, \pm \, \mathrm{side} \, \mathrm{of} \, \varSigma}{z^{\prime} \to z}
} \int_{\varSigma} \dfrac{f(s)}{s \! - \! z^{\prime}} \, \dfrac{\md s}{2 \pi
\mi},
\end{equation*}
and $\norm{w^{(n)}_{l}(\cdot) \! - \! w^{\infty}_{l}(\cdot)}_{\cap_{p \in \{2,
\infty\}}\mathcal{L}^{p}_{\mathrm{M}_{2}(\mathbb{C})}(\varSigma)} \! \to \! 0$
as $n \! \to \! \infty$, $l \! = \! \pm 1$. Then, $\exists \, \, N^{\ast} \!
\in \! \mathbb{N}$ such that, $\forall \, \, n \! > \! N^{\ast}$, $m^{\infty}
(z)$ and $m^{(n)}(z)$ exist, and $\norm{m^{(n)}_{l}(\cdot) \! - \!
m^{\infty}_{l}(\cdot)}_{\mathcal{L}^{2}_{\mathrm{M}_{2}(\mathbb{C})}
(\varSigma)} \! \to \! 0$ as $n \! \to \! \infty$, $l \! = \! \pm 1$.
\end{fffff}

A detailed exposition, including further motivations, for the construction 
of parametrices of the above-mentioned type can be found in 
\cite{a3,a45,a46,a47,a49,a79}; rather than regurgitating, verbatim, much of 
the analysis that can be found in the latter references, the point of view 
taken here is that one follows the scheme presented therein to obtain the 
results stated below, that is, the parametrix for the RHP 
$(\overset{o}{\mathscr{M}}^{\raise-1.0ex\hbox{$\scriptstyle \sharp$}}(z),
\overset{o}{\upsilon}^{\raise-1.0ex\hbox{$\scriptstyle \sharp$}}(z),\Sigma_{
o}^{\sharp})$ formulated in Lemma~4.2. In the case of the right-most 
end-points of the support of the `odd' equilibrium measure, $\lbrace 
a_{j}^{o} \rbrace_{j=1}^{N+1}$, a terse sketch of a proof is presented for 
the reader's convenience, and the remaining (left-most) end-points, namely, 
$b_{0}^{o},b_{1}^{o},\dotsc,b_{N}^{o}$, are analysed analogously.

The parametrix for the RHP
$(\overset{o}{\mathscr{M}}^{\raise-1.0ex\hbox{$\scriptstyle \sharp$}}(z),
\overset{o}{\upsilon}^{\raise-1.0ex\hbox{$\scriptstyle \sharp$}}(z),
\Sigma_{o}^{\sharp})$ is now presented. By a parametrix of the RHP
$(\overset{o}{\mathscr{M}}^{\raise-1.0ex\hbox{$\scriptstyle \sharp$}}(z),
\overset{o}{\upsilon}^{\raise-1.0ex\hbox{$\scriptstyle \sharp$}}(z),\Sigma_{
o}^{\sharp})$, in the neighbourhoods of the end-points of the support of the
`odd' equilibrium measure, $\lbrace b_{j-1}^{o},a_{j}^{o} \rbrace_{j=1}^{N+
1}$, is meant the solution of the RHPs formulated in the following two Lemmae
(Lemmae~4.6 and~4.7). Define the `small', mutually disjoint (open) discs about
the end-points of the support of the `odd' equilibrium measure as follows:
$\mathbb{U}_{\delta_{b_{j-1}}}^{o} \! := \! \lbrace \mathstrut z \! \in \!
\mathbb{C}; \, \vert z \! - \! b_{j-1}^{o} \vert \! < \! \delta_{b_{j-1}}^{o}
\! \in \! (0,1) \rbrace$ and $\mathbb{U}_{\delta_{a_{j}}}^{o} \! := \! \lbrace
\mathstrut z \! \in \! \mathbb{C}; \, \vert z \! - \! a_{j}^{o} \vert \! < \!
\delta_{a_{j}}^{o} \! \in \! (0,1) \rbrace$, $j \! = \! 1,\dotsc,N \! + \! 1$,
where $\delta_{b_{j-1}}^{o}$ and $\delta_{a_{j}}^{o}$ are sufficiently small,
positive real numbers chosen (amongst other things: see Lemmae~4.6 and~4.7
below) so that $\mathbb{U}_{\delta_{b_{i-1}}}^{o} \cap \mathbb{U}_{\delta_{
a_{j}}}^{o} \! = \! \varnothing$, $i,j \! = \! 1,\dotsc,N \! + \! 1$ (the
corresponding regions $\Omega_{b_{j-1}}^{o,l}$ and $\Omega_{a_{j}}^{o,l}$, and
arcs $\Sigma_{b_{j-1}}^{o,l}$ and $\Sigma_{a_{j}}^{o,l}$, $j \! = \! 1,\dotsc,
N \! + \! 1$, $l \! = \! 1,2,3,4$, respectively, are defined more precisely
below; see, also, Figures~5 and~6).
\begin{eeeee}
In order to simplify the results of Lemmae~4.6 and~4.7 (see below), it is
convenient to introduce the following notation: (i)
\begin{gather*}
\Psi^{o}_{1}(z) \! := \!
\begin{pmatrix}
\operatorname{Ai}(z) & \operatorname{Ai}(\omega^{2}z) \\
\operatorname{Ai}^{\prime}(z) & \omega^{2} \operatorname{Ai}^{\prime}
(\omega^{2}z)
\end{pmatrix} \! \me^{-\frac{\mi \pi}{6} \sigma_{3}}, \qquad \quad \Psi^{o}_{2}
(z) \! := \!
\begin{pmatrix}
\operatorname{Ai}(z) & \operatorname{Ai}(\omega^{2}z) \\
\operatorname{Ai}^{\prime}(z) & \omega^{2} \operatorname{Ai}^{\prime}
(\omega^{2}z)
\end{pmatrix} \! \me^{-\frac{\mi \pi}{6} \sigma_{3}}(\mathrm{I} \! - \!
\sigma_{-}), \\
\Psi^{o}_{3}(z) \! := \!
\begin{pmatrix}
\operatorname{Ai}(z) & -\omega^{2} \operatorname{Ai}(\omega z) \\
\operatorname{Ai}^{\prime}(z) & -\operatorname{Ai}^{\prime}(\omega z)
\end{pmatrix} \! \me^{-\frac{\mi \pi}{6} \sigma_{3}}(\mathrm{I} \! + \!
\sigma_{-}), \qquad \quad \Psi^{o}_{4}(z) \! := \!
\begin{pmatrix}
\operatorname{Ai}(z) & -\omega^{2} \operatorname{Ai}(\omega z) \\
\operatorname{Ai}^{\prime}(z) & -\operatorname{Ai}^{\prime}(\omega z)
\end{pmatrix} \! \me^{-\frac{\mi \pi}{6} \sigma_{3}},
\end{gather*}
where $\operatorname{Ai}(\pmb{\cdot})$ is the Airy function (cf.
Subsection~2.3), and $\omega \! = \! \exp (2 \pi \mi/3)$; and (ii)
\begin{equation*}
\mho^{o}_{j} \! := \!
\begin{cases}
\Omega^{o}_{j}, &\text{$j \! = \! 1,\dotsc,N$,} \\
0, &\text{$j \! = \! 0,N \! + \! 1$,}
\end{cases}
\end{equation*}
where $\Omega^{o}_{j} \! = \! 4 \pi \int_{b_{j}^{o}}^{a_{N+1}^{o}} \psi_{V}^{o}
(s) \, \md s$. \hfill $\blacksquare$
\end{eeeee}
\begin{ccccc}
Let $\overset{o}{\mathscr{M}}^{\raise-1.0ex\hbox{$\scriptstyle \sharp$}}
\colon \mathbb{C} \setminus \Sigma_{o}^{\sharp} \! \to \! \operatorname{SL}_{
2}(\mathbb{C})$ solve the {\rm RHP}
$(\overset{o}{\mathscr{M}}^{\raise-1.0ex\hbox{$\scriptstyle \sharp$}}(z),
\overset{o}{\upsilon}^{\raise-1.0ex\hbox{$\scriptstyle \sharp$}}(z),\Sigma_{
o}^{\sharp})$ formulated in Lemma~{\rm 4.2}, and set
\begin{equation*}
\mathbb{U}_{\delta_{b_{j-1}}}^{o} \! := \! \left\{\mathstrut z \! \in \!
\mathbb{C}; \, \vert z \! - \! b_{j-1}^{o} \vert \! < \! \delta_{b_{j-1}}^{o}
\! \in \! (0,1) \right\}, \quad j \! = \! 1,\dotsc,N \! + \! 1.
\end{equation*}
Let
\begin{equation*}
\Phi_{b_{j-1}}^{o}(z) \! := \! \left(\dfrac{3}{4} \! \left(n \! + \! \dfrac{
1}{2} \right) \! \xi_{b_{j-1}}^{o}(z) \right)^{2/3}, \quad j \! = \! 1,\dotsc,
N \! + \! 1,
\end{equation*}
with
\begin{equation*}
\xi_{b_{j-1}}^{o}(z) \! = \! -2 \int_{z}^{b_{j-1}^{o}}(R_{o}(s))^{1/2}h_{V}^{
o}(s) \, \md s,
\end{equation*}
where, for $z \! \in \! \mathbb{U}_{\delta_{b_{j-1}}}^{o} \setminus (-\infty,
b_{j-1}^{o})$, $\xi_{b_{j-1}}^{o}(z) \! = \! \mathfrak{b}(z \! - \! b_{j-1}^{
o})^{3/2}G_{b_{j-1}}^{o}(z)$, with $\mathfrak{b} \! := \! \pm 1$ for $z \!
\in \! \mathbb{C}_{\pm}$, and $G_{b_{j-1}}^{o}(z)$ analytic, in particular,
\begin{equation*}
G_{b_{j-1}}^{o}(z) \underset{z \to b_{j-1}^{o}}{=} \dfrac{4}{3}f(b_{j-1}^{o})
\! + \! \dfrac{4}{5}f^{\prime}(b_{j-1}^{o})(z \! - \! b_{j-1}^{o}) \! + \!
\dfrac{2}{7}f^{\prime \prime}(b_{j-1}^{o})(z \! - \! b_{j-1}^{o})^{2} \! + \!
\mathcal{O} \! \left((z \! - \! b_{j-1}^{o})^{3} \right),
\end{equation*}
where
\begin{align*}
f(b_{0}^{o})=& \, \mi (-1)^{N}h_{V}^{o}(b_{0}^{o}) \eta_{b_{0}^{o}}, \\
f^{\prime}(b_{0}^{o})=& \, \mi (-1)^{N} \! \left(\dfrac{1}{2}h_{V}^{o}(b_{0}^{
o}) \eta_{b_{0}^{o}} \! \left(\sum_{l=1}^{N} \! \left(\dfrac{1}{b_{0}^{o} \! -
\! b_{l}^{o}} \! + \! \dfrac{1}{b_{0}^{o} \! - \! a_{l}^{o}} \right) \! + \!
\dfrac{1}{b_{0}^{o} \! - \! a_{N+1}^{o}} \right) \! + \! (h_{V}^{o}(b_{0}^{o})
)^{\prime} \eta_{b_{0}^{o}} \right), \\
f^{\prime \prime}(b_{0}^{o})=& \, \mi (-1)^{N} \! \left(\dfrac{h_{V}^{o}(b_{
0}^{o})(h_{V}^{o}(b_{0}^{o}))^{\prime \prime} \! - \! ((h_{V}^{o}(b_{0}^{o}))^{
\prime})^{2}}{h_{V}^{o}(b_{0}^{o})} \eta_{b_{0}^{o}} \! - \! \dfrac{1}{2}h_{
V}^{o}(b_{0}^{o}) \eta_{b_{0}^{o}} \right. \\
\times&\left. \, \left(\sum_{l=1}^{N} \! \left(\dfrac{1}{(b_{0}^{o} \! - \!
b_{l}^{o})^{2}} \! + \! \dfrac{1}{(b_{0}^{o} \! - \! a_{l}^{o})^{2}} \right)
\! + \! \dfrac{1}{(b_{0}^{o} \! - \! a_{N+1}^{o})^{2}} \right) \right. \\
+&\left. \, \left(\dfrac{1}{2} \! \left(\sum_{k=1}^{N} \! \left(\dfrac{1}{
b_{0}^{o} \! - \! b_{k}^{o}} \! + \! \dfrac{1}{b_{0}^{o} \! - \! a_{k}^{o}}
\right) \! + \! \dfrac{1}{b_{0}^{o} \! - \! a_{N+1}^{o}} \right) \! + \!
\dfrac{(h_{V}^{o}(b_{0}^{o}))^{\prime}}{h_{V}^{o}(b_{0}^{o})} \right) \right.
\\
\times&\left. \, \left(\dfrac{1}{2}h_{V}^{o}(b_{0}^{o}) \eta_{b_{0}^{o}} \!
\left(\sum_{l=1}^{N} \! \left(\dfrac{1}{b_{0}^{o} \! - \! b_{l}^{o}} \! + \!
\dfrac{1}{b_{0}^{o} \! - \! a_{l}^{o}} \right) \! + \! \dfrac{1}{b_{0}^{o} \!
- \! a_{N+1}^{o}} \right) \! + \! (h_{V}^{o}(b_{0}^{o}))^{\prime} \eta_{b_{0}^{
o}} \right) \right),
\end{align*}
with
\begin{equation*}
\eta_{b_{0}^{o}} \! := \! \left((a_{N+1}^{o} \! - \! b_{0}^{o}) \prod_{k=1}^{N}
(b_{k}^{o} \! - \! b_{0}^{o})(a_{k}^{o} \! - \! b_{0}^{o}) \right)^{1/2} \quad
(> \! 0),
\end{equation*}
and, for $j \! = \! 1,\dotsc,N$,
\begin{align*}
f(b_{j}^{o})=& \, \mi (-1)^{N-j}h_{V}^{o}(b_{j}^{o}) \eta_{b_{j}^{o}}, \\
f^{\prime}(b_{j}^{o})=& \, \mi (-1)^{N-j} \! \left(\dfrac{1}{2}h_{V}^{o}(b_{
j}^{o}) \eta_{b_{j}^{o}} \! \left(\sum_{\substack{k=1\\k \not= j}}^{N} \!
\left(\dfrac{1}{b_{j}^{o} \! - \! b_{k}^{o}} \! + \! \dfrac{1}{b_{j}^{o} \! -
\! a_{k}^{o}} \right) \! + \! \dfrac{1}{b_{j}^{o} \! - \! a_{j}^{o}} \! + \!
\dfrac{1}{b_{j}^{o} \! - \! a_{N+1}^{o}} \! + \! \dfrac{1}{b_{j}^{o} \! - \!
b_{0}^{o}} \right) \right. \\
+&\left. \, (h_{V}^{o}(b_{j}^{o}))^{\prime} \eta_{b_{j}^{o}} \right), \\
f^{\prime \prime}(b_{j}^{o})=& \, \mi (-1)^{N-j} \! \left(\dfrac{h_{V}^{o}(b_{
j}^{o})(h_{V}^{o}(b_{j}^{o}))^{\prime \prime} \! - \! ((h_{V}^{o}(b_{j}^{o}))^{
\prime})^{2}}{h_{V}^{o}(b_{j}^{o})} \eta_{b_{j}^{o}} \! - \! \dfrac{1}{2}h_{
V}^{o}(b_{j}^{o}) \eta_{b_{j}^{o}} \! \left(\sum_{\substack{k=1\\k \not= j}}^{
N} \! \left(\dfrac{1}{(b_{j}^{o} \! - \! b_{k}^{o})^{2}} \! + \! \dfrac{1}{(b_{
j}^{o} \! - \! a_{k}^{o})^{2}} \right) \right. \right. \\
+&\left. \left. \, \dfrac{1}{(b_{j}^{o} \! - \! a_{j}^{o})^{2}} \! + \! \dfrac{
1}{(b_{j}^{o} \! - \! a_{N+1}^{o})^{2}} \! + \! \dfrac{1}{(b_{j}^{o} \! - \!
b_{0}^{o})^{2}} \right) \! + \! \left(\dfrac{(h_{V}^{o}(b_{j}^{o}))^{\prime}}{
h_{V}^{o}(b_{j}^{o})} \! + \! \dfrac{1}{2} \! \left(\sum_{\substack{k=1\\k
\not= j}}^{N} \! \left(\dfrac{1}{b_{j}^{o} \! - \! b_{k}^{o}} \! + \! \dfrac{
1}{b_{j}^{o} \! - \! a_{k}^{o}} \right) \right. \right. \right. \\
+&\left. \left. \left. \, \dfrac{1}{b_{j}^{o} \! - \! a_{j}^{o}} \! + \!
\dfrac{1}{b_{j}^{o} \! - \! a_{N+1}^{o}} \! + \! \dfrac{1}{b_{j}^{o} \! - \!
b_{0}^{o}} \right) \! \right) \! \left(\dfrac{1}{2}h_{V}^{o}(b_{j}^{o}) \eta_{
b_{j}^{o}} \! \left(\sum_{\substack{k=1\\k \not= j}}^{N} \! \left(\dfrac{1}{
b_{j}^{o} \! - \! b_{k}^{o}} \! + \! \dfrac{1}{b_{j}^{o} \! - \! a_{k}^{o}}
\right) \right. \right. \right. \\
+&\left. \left. \left. \, \dfrac{1}{b_{j}^{o} \! - \! a_{j}^{o}} \! + \!
\dfrac{1}{b_{j}^{o} \! - \! a_{N+1}^{o}} \! + \! \dfrac{1}{b_{j}^{o} \! - \!
b_{0}^{o}} \right) \! + \!  (h_{V}^{o}(b_{j}^{o}))^{\prime} \eta_{b_{j}^{o}}
\right) \right),
\end{align*}
with
\begin{equation*}
\eta_{b_{j}^{o}} \! := \! \left((b_{j}^{o} \! - \! a_{j}^{o})(a_{N+1}^{o} \!
- \! b_{j}^{o})(b_{j}^{o} \! - \! b_{0}^{o}) \prod_{k=1}^{j-1}(b_{j}^{o} \! -
\! b_{k}^{o})(b_{j}^{o} \! - \! a_{k}^{o}) \prod_{l=j+1}^{N}(b_{l}^{o} \! -
\! b_{j}^{o})(a_{l}^{o} \! - \! b_{j}^{o}) \right)^{1/2} \quad (> \! 0),
\end{equation*}
and $((0,1) \! \ni)$ $\delta_{b_{j-1}}^{o}$, $j \! = \! 1,\dotsc,N \! + \!
1$, are chosen sufficiently small so that $\Phi_{b_{j-1}}^{o}(z)$, which are
bi-holomorphic, conformal, and non-orientation preserving, map $\mathbb{U}_{
\delta_{b_{j-1}}}^{o}$ (and, thus, the oriented contours $\Sigma_{b_{j-1}}^{o}
\! := \! \cup_{l=1}^{4} \Sigma_{b_{j-1}}^{o,l}$, $j \! = \! 1,\dotsc,N \! + \!
1:$ Figure~{\rm 6)} injectively onto open $(n$-dependent) neighbourhoods
$\widehat{\mathbb{U}}_{\delta_{b_{j-1}}}^{o}$, $j \! = \! 1,\dotsc,N \! + \!
1$, of $0$ such that $\Phi_{b_{j-1}}^{o}(b_{j-1}^{o}) \! = \! 0$, $\Phi_{b_{j-
1}}^{o} \colon \mathbb{U}_{\delta_{b_{j-1}}}^{o} \! \to \! \widehat{\mathbb{
U}}_{\delta_{b_{j-1}}}^{o} \! := \! \Phi_{b_{j-1}}^{o}(\mathbb{U}_{\delta_{
b_{j-1}}}^{o})$, $\Phi_{b_{j-1}}^{o}(\mathbb{U}_{\delta_{b_{j-1}}}^{o} \cap
\Sigma_{b_{j-1}}^{o,l}) \! = \! \Phi_{b_{j-1}}^{o}(\mathbb{U}_{\delta_{b_{j-
1}}}^{o}) \cap \gamma_{b_{j-1}}^{o,l}$, and $\Phi_{b_{j-1}}^{o}(\mathbb{U}_{
\delta_{b_{j-1}}}^{o} \cap \Omega_{b_{j-1}}^{o,l}) \! = \! \Phi_{b_{j-1}}^{o}
(\mathbb{U}_{\delta_{b_{j-1}}}^{o}) \cap \widehat{\Omega}_{b_{j-1}}^{o,l}$,
$l \! = \! 1,2,3,4$, with $\widehat{\Omega}_{b_{j-1}}^{o,1} \! = \! \lbrace
\mathstrut \zeta \! \in \! \mathbb{C}; \, \arg (\zeta) \! \in \! (0,2 \pi/3)
\rbrace$, $\widehat{\Omega}_{b_{j-1}}^{o,2} \! = \! \lbrace \mathstrut \zeta
\! \in \! \mathbb{C}; \, \arg (\zeta) \! \in \! (2 \pi/3,\pi) \rbrace$,
$\widehat{\Omega}_{b_{j-1}}^{o,3} \! = \! \lbrace \mathstrut \zeta \! \in \!
\mathbb{C}; \, \arg (\zeta) \! \in \! (-\pi,-2 \pi/3) \rbrace$, and $\widehat{
\Omega}_{b_{j-1}}^{o,4} \! = \! \lbrace \mathstrut \zeta \! \in \! \mathbb{
C}; \, \arg (\zeta) \! \in \! (-2 \pi/3,0) \rbrace$.

The parametrix for the {\rm RHP}
$(\overset{o}{\mathscr{M}}^{\raise-1.0ex\hbox{$\scriptstyle \sharp$}}(z),
\overset{o}{\upsilon}^{\raise-1.0ex\hbox{$\scriptstyle \sharp$}}(z),\Sigma_{
o}^{\sharp})$, for $z \! \in \! \mathbb{U}_{\delta_{b_{j-1}}}^{o}$, $j \! = \!
1,\dotsc,N \! + \! 1$, is the solution of the following {\rm RHPs} for
$\mathcal{X}^{o} \colon \mathbb{U}_{\delta_{b_{j-1}}}^{o} \setminus \Sigma_{
b_{j-1}}^{o} \! \to \! \operatorname{SL}_{2}(\mathbb{C})$, $j \! = \! 1,\dotsc,
N \! + \! 1$, where $\Sigma_{b_{j-1}}^{o} \! := \! (\Phi_{b_{j-1}}^{o})^{-1}
(\gamma_{b_{j-1}}^{o})$, with $(\Phi_{b_{j-1}}^{o})^{-1}$ denoting the inverse
mapping, and $\gamma_{b_{j-1}}^{o} \! := \! \cup_{l=1}^{4} \gamma_{b_{j-1}}^{
o,l}:$ {\rm (i)} $\mathcal{X}^{o}(z)$ is holomorphic for $z \! \in \! \mathbb{
U}_{\delta_{b_{j-1}}}^{o} \setminus \Sigma_{b_{j-1}}^{o}$, $j \! = \! 1,
\dotsc,N \! + \! 1;$ {\rm (ii)} $\mathcal{X}^{o}_{\pm}(z) \! := \! \lim_{
\underset{z^{\prime} \, \in \, \pm \, \mathrm{side} \, \mathrm{of} \,
\Sigma_{b_{j-1}}^{o}}{z^{\prime} \to z}} \mathcal{X}^{o}(z^{\prime})$, $j \!
= \! 1,\dotsc,N \! + \! 1$, satisfy the boundary
condition
\begin{equation*}
\mathcal{X}^{o}_{+}(z) \! = \! \mathcal{X}^{o}_{-}(z)
\overset{o}{\upsilon}^{\raise-1.0ex\hbox{$\scriptstyle \sharp$}}(z), \quad z
\! \in \! \mathbb{U}_{\delta_{b_{j-1}}}^{o} \cap \Sigma_{b_{j-1}}^{o}, \quad j
\! = \! 1,\dotsc,N \! + \! 1,
\end{equation*}
where $\overset{o}{\upsilon}^{\raise-1.0ex\hbox{$\scriptstyle \sharp$}}(z)$ is
given in Lemma~{\rm 4.2;} and {\rm (iii)} uniformly for $z \! \in \! \partial
\mathbb{U}_{\delta_{b_{j-1}}}^{o} \! := \! \left\lbrace \mathstrut z \! \in \!
\mathbb{C}; \, \vert z \! - \! b_{j-1}^{o} \vert \! = \! \delta_{b_{j-1}}^{o}
\right\rbrace$, $j \! = \! 1,\dotsc,N \! + \! 1$,
\begin{equation*}
\overset{o}{m}^{\raise-1.0ex\hbox{$\scriptstyle \infty$}}(z)(\mathcal{X}^{o}
(z))^{-1} \! \underset{\underset{z \in \partial \mathbb{U}_{\delta_{b_{j-1}}}^{
o}}{n \to \infty}}{=} \! \mathrm{I} \! + \! \mathcal{O}((n \!+ \! 1/2)^{-1}),
\quad j \! = \! 1,\dotsc,N \! + \! 1.
\end{equation*}
The solutions of the {\rm RHPs} $(\mathcal{X}^{o}(z),
\overset{o}{\upsilon}^{\raise-1.0ex\hbox{$\scriptstyle \sharp$}}(z),
\mathbb{U}_{\delta_{b_{j-1}}}^{o} \cap \Sigma_{b_{j-1}}^{o})$, $j \! = \! 1,
\dotsc,N \! + \! 1$, are:\\
{\rm \pmb{(1)}} for $z \! \in \! \Omega_{b_{j-1}}^{o,1} \! := \! \mathbb{U}_{
\delta_{b_{j-1}}}^{o} \cap (\Phi_{b_{j-1}}^{o})^{-1}(\widehat{\Omega}_{b_{j-
1}}^{o,1})$, $j \! = \! 1,\dotsc,N \! + \! 1$,
\begin{equation*}
\mathcal{X}^{o}(z) \! = \! \sqrt{\smash[b]{\pi}} \, \me^{-\frac{\mi \pi}{3}}
\overset{o}{m}^{\raise-1.0ex\hbox{$\scriptstyle \infty$}}(z) \sigma_{3} \me^{
\frac{\mi}{2}(n+\frac{1}{2}) \mho_{j-1}^{o} \operatorname{ad}(\sigma_{3})} \!
\begin{pmatrix}
\mi & -\mi \\
1 & 1
\end{pmatrix} \! (\Phi_{b_{j-1}}^{o}(z))^{\frac{1}{4} \sigma_{3}} \Psi^{o}_{1}
(\Phi_{b_{j-1}}^{o}(z)) \me^{\frac{1}{2}(n+\frac{1}{2}) \xi_{b_{j-1}}^{o}(z)
\sigma_{3}} \sigma_{3},
\end{equation*}
where $\overset{o}{m}^{\raise-1.0ex\hbox{$\scriptstyle \infty$}}(z)$ is given
in Lemma~{\rm 4.5}, and $\Psi^{o}_{1}(z)$ and $\mho^{o}_{k}$ are defined in
Remark~{\rm 4.4;}\\
{\rm \pmb{(2)}} for $z \! \in \! \Omega_{b_{j-1}}^{o,2} \! := \! \mathbb{U}_{
\delta_{b_{j-1}}}^{o} \cap (\Phi_{b_{j-1}}^{o})^{-1}(\widehat{\Omega}_{b_{j-
1}}^{o,2})$, $j \! = \! 1,\dotsc,N \! + \! 1$,
\begin{equation*}
\mathcal{X}^{o}(z) \! = \! \sqrt{\smash[b]{\pi}} \, \me^{-\frac{\mi \pi}{3}}
\overset{o}{m}^{\raise-1.0ex\hbox{$\scriptstyle \infty$}}(z) \sigma_{3} \me^{
\frac{\mi}{2}(n+\frac{1}{2}) \mho_{j-1}^{o} \operatorname{ad}(\sigma_{3})} \!
\begin{pmatrix}
\mi & -\mi \\
1 & 1
\end{pmatrix} \! (\Phi_{b_{j-1}}^{o}(z))^{\frac{1}{4} \sigma_{3}} \Psi^{o}_{2}
(\Phi_{b_{j-1}}^{o}(z)) \me^{\frac{1}{2}(n+\frac{1}{2}) \xi_{b_{j-1}}^{o}(z)
\sigma_{3}} \sigma_{3},
\end{equation*}
where $\Psi^{o}_{2}(z)$ is defined in Remark~{\rm 4.4;}\\
{\rm \pmb{(3)}} for $z \! \in \! \Omega_{b_{j-1}}^{o,3} \! := \! \mathbb{U}_{
\delta_{b_{j-1}}}^{o} \cap (\Phi_{b_{j-1}}^{o})^{-1}(\widehat{\Omega}_{b_{j-
1}}^{o,3})$, $j \! = \! 1,\dotsc,N \! + \! 1$,
\begin{equation*}
\mathcal{X}^{o}(z) \! = \! \sqrt{\smash[b]{\pi}} \, \me^{-\frac{\mi \pi}{3}}
\overset{o}{m}^{\raise-1.0ex\hbox{$\scriptstyle \infty$}}(z) \sigma_{3} \me^{-
\frac{\mi}{2}(n+\frac{1}{2}) \mho^{o}_{j-1} \operatorname{ad}(\sigma_{3})} \!
\begin{pmatrix}
\mi & -\mi \\
1 & 1
\end{pmatrix} \! (\Phi_{b_{j-1}}^{o}(z))^{\frac{1}{4} \sigma_{3}} \Psi^{o}_{3}
(\Phi_{b_{j-1}}^{o}(z)) \me^{\frac{1}{2}(n+\frac{1}{2}) \xi_{b_{j-1}}^{o}(z)
\sigma_{3}} \sigma_{3},
\end{equation*}
where $\Psi^{o}_{3}(z)$ is defined in Remark~{\rm 4.4;}\\
{\rm \pmb{(4)}} for $z \! \in \! \Omega_{b_{j-1}}^{o,4} \! := \! \mathbb{U}_{
\delta_{b_{j-1}}}^{o} \cap (\Phi_{b_{j-1}}^{o})^{-1}(\widehat{\Omega}_{b_{j-
1}}^{o,4})$, $j \! = \! 1,\dotsc,N \! + \! 1$,
\begin{equation*}
\mathcal{X}^{o}(z) \! = \! \sqrt{\smash[b]{\pi}} \, \me^{-\frac{\mi \pi}{3}}
\overset{o}{m}^{\raise-1.0ex\hbox{$\scriptstyle \infty$}}(z) \sigma_{3} \me^{-
\frac{\mi}{2}(n+\frac{1}{2}) \mho_{j-1}^{o} \operatorname{ad}(\sigma_{3})} \!
\begin{pmatrix}
\mi & -\mi \\
1 & 1
\end{pmatrix} \! (\Phi_{b_{j-1}}^{o}(z))^{\frac{1}{4} \sigma_{3}} \Psi^{o}_{4}
(\Phi_{b_{j-1}}^{o}(z)) \me^{\frac{1}{2}(n+\frac{1}{2}) \xi_{b_{j-1}}^{o}(z)
\sigma_{3}} \sigma_{3},
\end{equation*}
where $\Psi^{o}_{4}(z)$ is defined in Remark~{\rm 4.4}.
\end{ccccc}
\begin{ccccc}
Let $\overset{o}{\mathscr{M}}^{\raise-1.0ex\hbox{$\scriptstyle \sharp$}}
\colon \mathbb{C} \setminus \Sigma_{o}^{\sharp} \! \to \! \operatorname{SL}_{
2}(\mathbb{C})$ solve the {\rm RHP}
$(\overset{o}{\mathscr{M}}^{\raise-1.0ex\hbox{$\scriptstyle \sharp$}}(z),
\overset{o}{\upsilon}^{\raise-1.0ex\hbox{$\scriptstyle \sharp$}}(z),\Sigma_{
o}^{\sharp})$ formulated in Lemma~{\rm 4.2}, and set
\begin{equation*}
\mathbb{U}_{\delta_{a_{j}}}^{o} \! := \! \left\{\mathstrut z \! \in \!
\mathbb{C}; \, \vert z \! - \! a_{j}^{o} \vert \! < \! \delta_{a_{j}}^{o} \!
\in \! (0,1) \right\}, \quad j \! = \! 1,\dotsc,N \! + \! 1.
\end{equation*}
Let
\begin{equation*}
\Phi_{a_{j}}^{o}(z) \! := \! \left(\dfrac{3}{4} \! \left(n \! + \! \dfrac{1}{
2} \right) \! \xi_{a_{j}}^{o}(z) \right)^{2/3}, \quad j \! = \! 1,\dotsc,N \!
+ \! 1,
\end{equation*}
with
\begin{equation*}
\xi_{a_{j}}^{o}(z) \! = \! 2 \int_{a_{j}^{o}}^{z}(R_{o}(s))^{1/2}h_{V}^{o}(s)
\, \md s,
\end{equation*}
where, for $z \! \in \! \mathbb{U}_{\delta_{a_{j}}}^{o} \setminus (-\infty,
a_{j}^{o})$, $\xi_{a_{j}}^{o}(z) \! = \! (z \! - \! a_{j}^{o})^{3/2}G_{a_{j}
}^{o}(z)$, $j \! = \! 1,\dotsc,N \! + \! 1$, with $G_{a_{j}}^{o}(z)$ analytic,
in particular,
\begin{equation*}
G_{a_{j}}^{o}(z) \underset{z \to a_{j}^{o}}{=} \dfrac{4}{3}f(a_{j}^{o}) \! +
\! \dfrac{4}{5}f^{\prime}(a_{j}^{o})(z \! - \! a_{j}^{o}) \! + \! \dfrac{2}{7}
f^{\prime \prime}(a_{j}^{o})(z \! - \! a_{j}^{o})^{2} \! + \! \mathcal{O} \!
\left((z \! - \! a_{j}^{o})^{3} \right),
\end{equation*}
where
\begin{align*}
f(a_{N+1}^{o})=& \, h_{V}^{o}(a_{N+1}^{o}) \eta_{a_{N+1}^{o}}, \\
f^{\prime}(a_{N+1}^{o})=& \, \dfrac{1}{2}h_{V}^{o}(a_{N+1}^{o}) \eta_{a_{N+
1}^{o}} \! \left(\sum_{l=1}^{N} \! \left(\dfrac{1}{a_{N+1}^{o} \! - \! b_{l}^{
o}} \! + \! \dfrac{1}{a_{N+1}^{o} \! - \! a_{l}^{o}} \right) \! + \! \dfrac{
1}{a_{N+1}^{o} \! - \! b_{0}^{o}} \right) \\
+& \, (h_{V}^{o}(a_{N+1}^{o}))^{\prime} \eta_{a_{N+1}^{o}}, \\
f^{\prime \prime}(a_{N+1}^{o})=& \, \dfrac{h_{V}^{o}(a_{N+1}^{o})(h_{V}^{o}
(a_{N+1}^{o}))^{\prime \prime} \! - \! ((h_{V}^{o}(a_{N+1}^{o}))^{\prime})^{
2}}{h_{V}^{o}(a_{N+1}^{o})} \eta_{a_{N+1}^{o}} \! - \! \dfrac{1}{2}h_{V}^{o}
(a_{N+1}^{o}) \eta_{a_{N+1}^{o}} \\
\times& \, \left(\sum_{l=1}^{N} \! \left(\dfrac{1}{(a_{N+1}^{o} \! - \! b_{
l}^{o})^{2}} \! + \! \dfrac{1}{(a_{N+1}^{o} \! - \! a_{l}^{o})^{2}} \right)
\! + \! \dfrac{1}{(a_{N+1}^{o} \! - \! b_{0}^{o})^{2}} \right) \\
+& \, \left(\dfrac{1}{2} \! \left(\sum_{k=1}^{N} \! \left(\dfrac{1}{a_{N+1}^{
o} \! - \! b_{k}^{o}} \! + \! \dfrac{1}{a_{N+1}^{o} \! - \! a_{k}^{o}} \right)
\! + \! \dfrac{1}{a_{N+1}^{o} \! - \! b_{0}^{o}} \right) \! + \! \dfrac{(
h_{V}^{o}(a_{N+1}^{o}))^{\prime}}{h_{V}^{o}(a_{N+1}^{o})} \right) \\
\times& \, \left(\dfrac{1}{2}h_{V}^{o}(a_{N+1}^{o}) \eta_{a_{N+1}^{o}} \!
\left(\sum_{l=1}^{N} \! \left(\dfrac{1}{a_{N+1}^{o} \! - \! a_{l}^{o}} \! + \!
\dfrac{1}{a_{N+1}^{o} \! - \! b_{l}^{o}} \right) \! + \! \dfrac{1}{a_{N+1}^{o}
\! - \! b_{0}^{o}} \right) \right. \\
+& \left. \, (h_{V}^{o}(a_{N+1}^{o}))^{\prime} \eta_{a_{N+1}^{o}} \right),
\end{align*}
with
\begin{equation*}
\eta_{a_{N+1}^{o}} \! := \! \left((a_{N+1}^{o} \! - \! b_{0}^{o}) \prod_{k=1
}^{N}(a_{N+1}^{o} \! - \! b_{k}^{o})(a_{N+1}^{o} \! - \! a_{k}^{o}) \right)^{
1/2} \quad (> \! 0),
\end{equation*}
and, for $j \! = \! 1,\dotsc,N$,
\begin{align*}
f(a_{j}^{o})=& \, (-1)^{N-j+1}h_{V}^{o}(a_{j}^{o}) \eta_{a_{j}^{o}}, \\
f^{\prime}(a_{j}^{o})=& \, (-1)^{N-j+1} \! \left(\dfrac{1}{2}h_{V}^{o}(a_{j}^{
o}) \eta_{a_{j}^{o}} \! \left(\sum_{\substack{k=1\\k \not= j}}^{N} \! \left(
\dfrac{1}{a_{j}^{o} \! - \! b_{k}^{o}} \! + \! \dfrac{1}{a_{j}^{o} \! - \!
a_{k}^{o}} \right) \! + \! \dfrac{1}{a_{j}^{o} \! - \! b_{j}^{o}} \! + \!
\dfrac{1}{a_{j}^{o} \! - \! a_{N+1}^{o}} \! + \! \dfrac{1}{a_{j}^{o} \! - \!
b_{0}^{o}} \right) \right. \\
+&\left. \, (h_{V}^{o}(a_{j}^{o}))^{\prime} \eta_{a_{j}^{o}} \right), \\
f^{\prime \prime}(a_{j}^{o})=& \, (-1)^{N-j+1} \! \left(\dfrac{h_{V}^{o}(a_{
j}^{o})(h_{V}^{o}(a_{j}^{o}))^{\prime \prime} \! - \! ((h_{V}^{o}(a_{j}^{o}))^{
\prime})^{2}}{h_{V}^{o}(a_{j}^{o})} \eta_{a_{j}^{o}} \! - \! \dfrac{1}{2}h_{
V}^{o}(a_{j}^{o}) \eta_{a_{j}^{o}} \! \left(\sum_{\substack{k=1\\k \not= j}}^{
N} \! \left(\dfrac{1}{(a_{j}^{o} \! - \! b_{k}^{o})^{2}} \! + \! \dfrac{1}{(a_{
j}^{o} \! - \! a_{k}^{o})^{2}} \right) \right. \right. \\
+&\left. \left. \, \dfrac{1}{(a_{j}^{o} \! - \! b_{j}^{o})^{2}} \! + \! \dfrac{
1}{(a_{j}^{o} \! - \! a_{N+1}^{o})^{2}} \! + \! \dfrac{1}{(a_{j}^{o} \! - \!
b_{0}^{o})^{2}} \right) \! + \! \left(\dfrac{(h_{V}^{o}(a_{j}^{o}))^{\prime}}{
h_{V}^{o}(a_{j}^{o})} \! + \! \dfrac{1}{2} \! \left(\sum_{\substack{k=1\\k
\not= j}}^{N} \! \left(\dfrac{1}{a_{j}^{o} \! - \! b_{k}^{o}} \! + \! \dfrac{
1}{a_{j}^{o} \! - \! a_{k}^{o}} \right) \right. \right. \right. \\
+&\left. \left. \left. \, \dfrac{1}{a_{j}^{o} \! - \! b_{j}^{o}} \! + \!
\dfrac{1}{a_{j}^{o} \! - \! a_{N+1}^{o}} \! + \! \dfrac{1}{a_{j}^{o} \! - \!
b_{0}^{o}} \right) \! \right) \! \left(\dfrac{1}{2}h_{V}^{o}(a_{j}^{o}) \eta_{
a_{j}^{o}} \! \left(\sum_{\substack{k=1\\k \not= j}}^{N} \! \left(\dfrac{1}{
a_{j}^{o} \! - \! b_{k}^{o}} \! + \! \dfrac{1}{a_{j}^{o} \! - \! a_{k}^{o}}
\right) \right. \right. \right. \\
+&\left. \left. \left. \, \dfrac{1}{a_{j}^{o} \! - \! b_{j}^{o}} \! + \!
\dfrac{1}{a_{j}^{o} \! - \! a_{N+1}^{o}} \! + \! \dfrac{1}{a_{j}^{o} \! - \!
b_{0}^{o}} \right) \! + \!  (h_{V}^{o}(a_{j}^{o}))^{\prime} \eta_{a_{j}^{o}}
\right) \right),
\end{align*}
with
\begin{equation*}
\eta_{a_{j}^{o}} \! := \! \left((b_{j}^{o} \! - \! a_{j}^{o})(a_{N+1}^{o} \!
- \! a_{j}^{o})(a_{j}^{o} \! - \! b_{0}^{o}) \prod_{k=1}^{j-1}(a_{j}^{o} \! -
\! b_{k}^{o})(a_{j}^{o} \! - \! a_{k}^{o}) \prod_{l=j+1}^{N}(b_{l}^{o} \! -
\! a_{j}^{o})(a_{l}^{o} \! - \! a_{j}^{o}) \right)^{1/2} \quad (> \! 0),
\end{equation*}
and $((0,1) \! \ni)$ $\delta_{a_{j}}^{o}$, $j \! = \! 1,\dotsc,N \! + \!
1$, are chosen sufficiently small so that $\Phi_{a_{j}}^{o}(z)$, which are
bi-holomorph\-i\-c, conformal, and orientation preserving, map $\mathbb{U}_{
\delta_{a_{j}}}^{o}$ (and, thus, the oriented contours $\Sigma_{a_{j}}^{o} \!
:= \! \cup_{l=1}^{4} \Sigma_{a_{j}}^{o,l}$, $j \! = \! 1,\dotsc,N \! + \! 1:$
Figure~{\rm 5)} injectively onto open $(n$-dependent) neighbourhoods $\widehat{
\mathbb{U}}_{\delta_{a_{j}}}^{o}$, $j \! = \! 1,\dotsc,N \! + \! 1$, of $0$
such that $\Phi_{a_{j}}^{o}(a_{j}^{o}) \! = \! 0$, $\Phi_{a_{j}}^{o} \colon
\mathbb{U}_{\delta_{a_{j}}}^{o} \! \to \! \widehat{\mathbb{U}}_{\delta_{a_{
j}}}^{o} \! := \! \Phi_{a_{j}}^{o}(\mathbb{U}_{\delta_{a_{j}}}^{o})$, $\Phi_{
a_{j}}^{o}(\mathbb{U}_{\delta_{a_{j}}}^{o} \cap \Sigma_{a_{j}}^{o,l}) \! = \!
\Phi_{a_{j}}^{o}(\mathbb{U}_{\delta_{a_{j}}}^{o}) \cap \gamma_{a_{j}}^{o,l}$,
and $\Phi_{a_{j}}^{o}(\mathbb{U}_{\delta_{a_{j}}}^{o} \cap \Omega_{a_{j}}^{o,
l}) \! = \! \Phi_{a_{j}}^{o}(\mathbb{U}_{\delta_{a_{j}}}^{o}) \cap \widehat{
\Omega}_{a_{j}}^{o,l}$, $l \! = \! 1,2,3,4$, with $\widehat{\Omega}_{a_{j}}^{
o,1} \! = \! \lbrace \mathstrut \zeta \! \in \! \mathbb{C}; \, \arg (\zeta) \!
\in \! (0,2 \pi/3) \rbrace$, $\widehat{\Omega}_{a_{j}}^{o,2} \! = \! \lbrace
\mathstrut \zeta \! \in \! \mathbb{C}; \, \arg (\zeta) \! \in \! (2 \pi/3,\pi)
\rbrace$, $\widehat{\Omega}_{a_{j}}^{o,3} \! = \! \lbrace \mathstrut \zeta \!
\in \! \mathbb{C}; \, \arg (\zeta) \! \in \! (-\pi,-2 \pi/3) \rbrace$, and
$\widehat{\Omega}_{a_{j}}^{o,4} \! = \! \lbrace \mathstrut \zeta \! \in \!
\mathbb{C}; \, \arg (\zeta) \! \in \! (-2 \pi/3,0) \rbrace$.

The parametrix for the {\rm RHP}
$(\overset{o}{\mathscr{M}}^{\raise-1.0ex\hbox{$\scriptstyle \sharp$}}(z),
\overset{o}{\upsilon}^{\raise-1.0ex\hbox{$\scriptstyle \sharp$}}(z),\Sigma_{
o}^{\sharp})$, for $z \! \in \! \mathbb{U}_{\delta_{a_{j}}}^{o}$, $j \! = \!
1,\dotsc,N \! + \! 1$, is the solution of the following {\rm RHPs} for
$\mathcal{X}^{o} \colon \mathbb{U}_{\delta_{a_{j}}}^{o} \setminus \Sigma_{a_{
j}}^{o} \! \to \! \operatorname{SL}_{2}(\mathbb{C})$, $j \! = \! 1,\dotsc,N \!
+ \! 1$, where $\Sigma_{a_{j}}^{o} \! := \! (\Phi_{a_{j}}^{o})^{-1}(\gamma_{
a_{j}}^{o})$, with $(\Phi_{a_{j}}^{o})^{-1}$ denoting the inverse mapping, and
$\gamma_{a_{j}}^{o} \! := \! \cup_{l=1}^{4} \gamma_{a_{j}}^{o,l}:$ {\rm (i)}
$\mathcal{X}^{o}(z)$ is holomorphic for $z \! \in \! \mathbb{U}_{\delta_{a_{
j}}}^{o} \setminus \Sigma_{a_{j}}^{o}$, $j \! = \! 1,\dotsc,N \! + \! 1;$
{\rm (ii)} $\mathcal{X}^{o}_{\pm}(z) \! := \! \lim_{\underset{z^{\prime} \,
\in \, \pm \, \mathrm{side} \, \mathrm{of} \, \Sigma_{a_{j}}^{o}}{z^{\prime}
\to z}} \mathcal{X}^{o}(z^{\prime})$, $j \! = \! 1,\dotsc,N \! + \! 1$,
satisfy the boundary condition
\begin{equation*}
\mathcal{X}^{o}_{+}(z) \! = \! \mathcal{X}^{o}_{-}(z)
\overset{o}{\upsilon}^{\raise-1.0ex\hbox{$\scriptstyle \sharp$}}(z), \quad z
\! \in \! \mathbb{U}_{\delta_{a_{j}}}^{o} \cap \Sigma_{a_{j}}^{o}, \quad j \!
= \! 1,\dotsc,N \! + \! 1,
\end{equation*}
where $\overset{o}{\upsilon}^{\raise-1.0ex\hbox{$\scriptstyle \sharp$}}(z)$ is
given in Lemma~{\rm 4.2;} and {\rm (iii)} uniformly for $z \! \in \! \partial
\mathbb{U}_{\delta_{a_{j}}}^{o} \! := \! \left\lbrace \mathstrut z \! \in \!
\mathbb{C}; \, \vert z \! - \! a_{j}^{o} \vert \! = \! \delta_{a_{j}}^{o}
\right\rbrace$, $j \! = \! 1,\dotsc,N \! + \! 1$,
\begin{equation*}
\overset{o}{m}^{\raise-1.0ex\hbox{$\scriptstyle \infty$}}(z)(\mathcal{X}^{o}
(z))^{-1} \! \underset{\underset{z \in \partial \mathbb{U}_{\delta_{a_{j}}}^{
o}}{n \to \infty}}{=} \! \mathrm{I} \! + \! \mathcal{O}((n \! + \! 1/2)^{-1}),
\quad j \! = \! 1,\dotsc,N \! + \! 1.
\end{equation*}
The solutions of the {\rm RHPs} $(\mathcal{X}^{o}(z),
\overset{o}{\upsilon}^{\raise-1.0ex\hbox{$\scriptstyle \sharp$}}(z),
\mathbb{U}_{\delta_{a_{j}}}^{o} \cap \Sigma_{a_{j}}^{o})$, $j \! = \! 1,\dotsc,
N \! + \! 1$, are:\\
{\rm \pmb{(1)}} for $z \! \in \! \Omega_{a_{j}}^{o,1} \! := \! \mathbb{U}_{
\delta_{a_{j}}}^{o} \cap (\Phi_{a_{j}}^{o})^{-1}(\widehat{\Omega}_{a_{j}}^{o,
1})$, $j \! = \! 1,\dotsc,N \! + \! 1$,
\begin{equation*}
\mathcal{X}^{o}(z) \! = \! \sqrt{\smash[b]{\pi}} \, \me^{-\frac{\mi \pi}{3}}
\overset{o}{m}^{\raise-1.0ex\hbox{$\scriptstyle \infty$}}(z) \me^{\frac{\mi}{2}
(n+\frac{1}{2}) \mho_{j}^{o} \operatorname{ad}(\sigma_{3})} \!
\begin{pmatrix}
\mi & -\mi \\
1 & 1
\end{pmatrix} \! (\Phi_{a_{j}}^{o}(z))^{\frac{1}{4} \sigma_{3}} \Psi^{o}_{1}
(\Phi_{a_{j}}^{o}(z)) \me^{\frac{1}{2}(n+\frac{1}{2}) \xi_{a_{j}}^{o}(z)
\sigma_{3}},
\end{equation*}
where $\overset{o}{m}^{\raise-1.0ex\hbox{$\scriptstyle \infty$}}(z)$ is given
in Lemma~{\rm 4.5}, and $\Psi^{o}_{1}(z)$ and $\mho_{k}^{o}$ are defined in
Remark~{\rm 4.4;}\\
{\rm \pmb{(2)}} for $z \! \in \! \Omega_{a_{j}}^{o,2} \! := \! \mathbb{U}_{
\delta_{a_{j}}}^{o} \cap (\Phi_{a_{j}}^{o})^{-1}(\widehat{\Omega}_{a_{j}}^{o,
2})$, $j \! = \! 1,\dotsc,N \! + \! 1$,
\begin{equation*}
\mathcal{X}^{o}(z) \! = \! \sqrt{\smash[b]{\pi}} \, \me^{-\frac{\mi \pi}{3}}
\overset{o}{m}^{\raise-1.0ex\hbox{$\scriptstyle \infty$}}(z) \me^{\frac{\mi}{2}
(n+\frac{1}{2}) \mho_{j}^{o} \operatorname{ad}(\sigma_{3})} \!
\begin{pmatrix}
\mi & -\mi \\
1 & 1
\end{pmatrix} \! (\Phi_{a_{j}}^{o}(z))^{\frac{1}{4} \sigma_{3}} \Psi^{o}_{2}
(\Phi_{a_{j}}^{o}(z)) \me^{\frac{1}{2}(n+\frac{1}{2}) \xi_{a_{j}}^{o}(z)
\sigma_{3}},
\end{equation*}
where $\Psi^{o}_{2}(z)$ is defined in Remark~{\rm 4.4;}\\
{\rm \pmb{(3)}} for $z \! \in \! \Omega_{a_{j}}^{o,3} \! := \! \mathbb{U}_{
\delta_{a_{j}}}^{o} \cap (\Phi_{a_{j}}^{o})^{-1}(\widehat{\Omega}_{a_{j}}^{o,
3})$, $j \! = \! 1,\dotsc,N \! + \! 1$,
\begin{equation*}
\mathcal{X}^{o}(z) \! = \! \sqrt{\smash[b]{\pi}} \, \me^{-\frac{\mi \pi}{3}}
\overset{o}{m}^{\raise-1.0ex\hbox{$\scriptstyle \infty$}}(z) \me^{-\frac{\mi}{
2}(n+\frac{1}{2}) \mho^{o}_{j} \operatorname{ad}(\sigma_{3})} \!
\begin{pmatrix}
\mi & -\mi \\
1 & 1
\end{pmatrix} \! (\Phi_{a_{j}}^{o}(z))^{\frac{1}{4} \sigma_{3}} \Psi^{o}_{3}
(\Phi_{a_{j}}^{o}(z)) \me^{\frac{1}{2}(n+\frac{1}{2}) \xi_{a_{j}}^{o}(z)
\sigma_{3}},
\end{equation*}
where $\Psi^{o}_{3}(z)$ is defined in Remark~{\rm 4.4;}\\
{\rm \pmb{(4)}} for $z \! \in \! \Omega_{a_{j}}^{o,4} \! := \! \mathbb{U}_{
\delta_{a_{j}}}^{o} \cap (\Phi_{a_{j}}^{o})^{-1}(\widehat{\Omega}_{a_{j}}^{
o,4})$, $j \! = \! 1,\dotsc,N \! + \! 1$,
\begin{equation*}
\mathcal{X}^{o}(z) \! = \! \sqrt{\smash[b]{\pi}} \, \me^{-\frac{\mi \pi}{3}}
\overset{o}{m}^{\raise-1.0ex\hbox{$\scriptstyle \infty$}}(z) \me^{-\frac{\mi}{
2}(n+\frac{1}{2}) \mho_{j}^{o} \operatorname{ad}(\sigma_{3})} \!
\begin{pmatrix}
\mi & -\mi \\
1 & 1
\end{pmatrix} \! (\Phi_{a_{j}}^{o}(z))^{\frac{1}{4} \sigma_{3}} \Psi^{o}_{4}
(\Phi_{a_{j}}^{o}(z)) \me^{\frac{1}{2}(n+\frac{1}{2}) \xi_{a_{j}}^{o}(z)
\sigma_{3}},
\end{equation*}
where $\Psi^{o}_{4}(z)$ is defined in Remark~{\rm 4.4}.
\end{ccccc}
\begin{eeeee}
Perusing Lemmae~4.6 and~4.7, one notes that the normalisation condition at
zero, which is needed in order to guarantee the existence of solutions to the
corresponding (parametrix) RHPs, is absent. The normalisation conditions at
zero are replaced by the (uniform) matching conditions
$\overset{o}{m}^{\raise-1.0ex\hbox{$\scriptstyle \infty$}}(z)(\mathcal{X}^{o}
(z))^{-1} \! =_{\underset{z \in \partial \mathbb{U}_{\delta_{\ast_{j}}}^{o}}{n
\to \infty}} \! \mathrm{I} \! + \! \mathcal{O}((n \! + \! 1/2)^{-1})$, where
$\ast_{j} \! \in \! \{b_{j-1},a_{j}\}$, $j \! = \! 1,\dotsc,N \! + \! 1$, with
$\partial \mathbb{U}_{\delta_{\ast_{j}}}^{o}$ defined in Lemmae~4.6 and~4.7.
\hfill $\blacksquare$
\end{eeeee}

\emph{Sketch of proof of Lemma~{\rm 4.7}.} Let
$(\overset{o}{\mathscr{M}}^{\raise-1.0ex\hbox{$\scriptstyle \sharp$}}(z),
\overset{o}{\upsilon}^{\raise-1.0ex\hbox{$\scriptstyle \sharp$}}(z),\Sigma_{
o}^{\sharp})$ be the RHP formulated in Lemma~4.2, and recall the definitions
stated therein. For each $a_{j}^{o} \! \in \! \operatorname{supp}(\mu_{V}^{
o})$, $j \! = \! 1,\dotsc,N \! + \! 1$, define $\mathbb{U}_{\delta_{a_{j}}}^{
o}$, $j \! = \! 1,\dotsc,N \! + \! 1$, as in the Lemma, that is, surround each
right-most end-point $a_{j}^{o}$ by open discs of radius $\delta_{a_{j}}^{o}
\! \in \! (0,1)$ centred at $a_{j}^{o}$. Recalling the formula for
$\overset{o}{\upsilon}^{\raise-1.0ex\hbox{$\scriptstyle \sharp$}}(z)$ given in
Lemma~4.2, one shows, via the proof of Lemma~4.1, that:
\begin{compactenum}
\item[(1)] $4 \pi \mi \int_{z}^{a_{N+1}^{o}} \psi_{V}^{o}(s) \, \md s \! = \!
4 \pi \mi (\int_{z}^{a_{j}^{o}} \! + \! \int_{a_{j}^{o}}^{b_{j}^{o}} \! + \!
\int_{b_{j}^{o}}^{a_{N+1}^{o}}) \psi_{V}^{o}(s) \, \md s$, whence, recalling
the expression for the density of the `odd' equilibrium measure given in
Lemma~3.5, that is, $\md \mu_{V}^{o}(x) \! := \! \psi_{V}^{o}(x) \, \md x \! =
\! \tfrac{1}{2 \pi \mi}(R_{o}(x))^{1/2}_{+}h_{V}^{o}(x) \pmb{1}_{J_{o}}(x) \,
\md x$, one arrives at, upon considering the analytic continuation of $4 \pi
\mi \linebreak[4]
\cdot \int_{z}^{a_{N+1}^{o}} \psi_{V}^{o}(s) \, \md s$ to $\mathbb{C}
\setminus \mathbb{R}$ (cf. proof of Lemma~4.1), in particular, to the oriented
(open) skeletons $\mathbb{U}_{\delta_{a_{j}}}^{o} \cap (J_{j}^{o,\smallfrown}
\cup J_{j}^{o,\smallsmile})$, $j \! = \! 1,\dotsc,N \! + \! 1$, $4 \pi \mi
\int_{z}^{a_{N+1}^{o}} \psi_{V}^{o}(s) \, \md s \! = \! -\xi_{a_{j}}^{o}(z) \!
+ \! \mi \mho_{j}^{o}$, $j \! = \! 1,\dotsc,N \! + \! 1$, where $\xi_{a_{j}}^{
o}(z) \! = \! 2 \int_{a_{j}^{o}}^{z}(R_{o}(s))^{1/2}h_{V}^{o}(s) \, \md s$,
and $\mho_{j}^{o}$ are defined in Remark~4.4;
\item[(2)] $g^{o}_{+}(z) \! + \! g^{o}_{-}(z) \! - \! \widetilde{V}(z) \! -
\! \ell_{o} \! - \! \mathfrak{Q}^{+}_{\mathscr{A}} \! - \! \mathfrak{Q}^{-}_{
\mathscr{A}} \! = \! -(2 \! + \! \tfrac{1}{n}) \int_{a_{j}^{o}}^{z}(R_{o}(s)
)^{1/2}h_{V}^{o}(s) \, \md s \! < \! 0$, $z \! \in \! (a_{N+1}^{o},+\infty)
\cup (\cup_{j=1}^{N}(a_{j}^{o},b_{j}^{o}))$.
\end{compactenum}
Via the latter formulae, which appear in the $(i \, j)$-elements, $i,j \! =
\! 1,2$, of the jump matrix
$\overset{o}{\upsilon}^{\raise-1.0ex\hbox{$\scriptstyle \sharp$}}(z)$,
denoting $\overset{o}{\mathscr{M}}^{\raise-1.0ex\hbox{$\scriptstyle \sharp$}}
(z)$ by $\mathcal{X}^{o}(z)$ for $z \! \in \! \mathbb{U}_{\delta_{a_{j}}}^{
o}$, $j \! = \! 1,\dotsc,N \! + \! 1$, and defining
\begin{equation*}
\mathscr{P}_{a_{j}}^{o}(z) \! := \!
\begin{cases}
\mathcal{X}^{o}(z) \me^{-\frac{1}{2}(n+\frac{1}{2}) \xi_{a_{j}}^{o}(z) \sigma_{
3}} \, \me^{\frac{\mi}{2}(n+\frac{1}{2}) \mho_{j}^{o} \sigma_{3}}, &\text{$z
\! \in \! \mathbb{C}_{+} \cap \mathbb{U}_{\delta_{a_{j}}}^{o}, \quad j \! = \!
1,\dotsc,N \! + \! 1$,} \\
\mathcal{X}^{o}(z) \me^{-\frac{1}{2}(n+\frac{1}{2}) \xi_{a_{j}}^{o}(z) \sigma_{
3}} \, \me^{-\frac{\mi}{2}(n+\frac{1}{2}) \mho_{j}^{o} \sigma_{3}}, &\text{$z
\! \in \! \mathbb{C}_{-} \cap \mathbb{U}_{\delta_{a_{j}}}^{o}, \quad j \! =
\! 1,\dotsc,N \! + \! 1$,}
\end{cases}
\end{equation*}
one notes that $\mathscr{P}_{a_{j}}^{o} \colon \mathbb{U}_{\delta_{a_{j}}}^{o}
\setminus J_{a_{j}}^{o} \! \to \! \operatorname{GL}_{2}(\mathbb{C})$, where
$J_{a_{j}}^{o} \! := \! J_{j}^{o,\smallfrown} \cup J_{j}^{o,\smallsmile} \cup
(a_{j}^{o} \! - \! \delta_{a_{j}}^{o},a_{j}^{o} \! + \! \delta_{a_{j}}^{o})$,
$j \! = \! 1,\dotsc,N \! + \! 1$, solves the RHP $(\mathscr{P}_{a_{j}}^{o}(z),
\upsilon^{o}_{\mathscr{P}_{a_{j}}}(z),J_{a_{j}}^{o})$, with constant jump
matrices $\upsilon_{\mathscr{P}_{a_{j}}}^{o}(z)$, $j \! = \! 1,\dotsc,N \! +
\! 1$, defined by
\begin{equation*}
\upsilon_{\mathscr{P}_{a_{j}}}^{o}(z) \! := \!
\begin{cases}
\mathrm{I} \! + \! \sigma_{-}, &\text{$z \! \in \! \mathbb{U}_{\delta_{a_{j}}
}^{o} \cap (J_{j}^{o,\smallfrown} \cup J_{j}^{o,\smallsmile}) \! = \! \Sigma^{
o,1}_{a_{j}} \cup \Sigma^{o,3}_{a_{j}}$,} \\
\mathrm{I} \! + \! \sigma_{+}, &\text{$z \! \in \! \mathbb{U}_{\delta_{a_{j}}
}^{o} \cap (a_{j}^{o},a_{j}^{o} \! + \! \delta_{a_{j}}^{o}) \! = \! \Sigma^{o,
4}_{a_{j}}$,} \\
\mi \sigma_{2}, &\text{$z \! \in \! \mathbb{U}_{\delta_{a_{j}}}^{o} \cap (a_{
j}^{o} \! - \! \delta_{a_{j}}^{o},a_{j}^{o}) \! = \! \Sigma^{o,2}_{a_{j}}$,}
\end{cases}
\end{equation*}
subject, still, to the asymptotic matching conditions
$\overset{o}{m}^{\raise-1.0ex\hbox{$\scriptstyle \infty$}}(z)(\mathcal{X}^{o}
(z))^{-1} \! =_{n \to \infty} \! \mathrm{I} \! + \! \mathcal{O}((n \! + \!
1/2)^{-1})$, uniformly for $z \! \in \! \partial \mathbb{U}_{\delta_{a_{j}}}^{
o}$, $j \! = \! 1,\dotsc,N \! + \! 1$.

Set, as in the Lemma, $\Phi_{a_{j}}^{o}(z) \! := \! (\tfrac{3}{4}(n \! + \!
\tfrac{1}{2}) \xi_{a_{j}}^{o}(z))^{2/3}$, $j \! = \! 1,\dotsc,N \! + \! 1$,
with $\xi_{a_{j}}^{o}(z)$ defined above: a careful analysis of the branch cuts
shows that, for $z \! \in \! \mathbb{U}_{\delta_{a_{j}}}^{o}$, $j \! = \! 1,
\dotsc,N \! + \! 1$, $\Phi_{a_{j}}^{o}(z)$ and $\xi_{a_{j}}^{o}(z)$ satisfy
the properties stated in the Lemma; in particular, for $\Phi_{a_{j}}^{o}
\colon \mathbb{U}_{\delta_{a_{j}}}^{o} \! \to \! \mathbb{C}$, $j \! = \! 1,
\dotsc,N \! + \! 1$, $\Phi_{a_{j}}^{o}(z) \! = \! (z \! - \! a_{j}^{o})^{3/2}
G_{a_{j}}^{o}(z)$, with $G_{a_{j}}^{o}(z)$ holomorphic for $z \! \in \!
\mathbb{U}_{\delta_{a_{j}}}^{o}$ and characterised in the Lemma, $\Phi_{a_{
j}}^{o}(a_{j}^{o}) \! = \! 0$, $(\Phi_{a_{j}}^{o}(z))^{\prime} \! \not= \!
0$, $z \! \in \! \mathbb{U}_{\delta_{a_{j}}}^{o}$, and where $(\Phi_{a_{j}}^{
o}(a_{j}^{o}))^{\prime} \! = \! ((n \! + \! \tfrac{1}{2})f(a_{j}^{o}))^{2/3}
\! > \! 0$, with $f(a_{j}^{o})$ given in the Lemma. One now chooses $\delta_{
a_{j}}^{o}$ $(\in \! (0,1))$, $j \! = \! 1,\dotsc,N \! + \! 1$, and the
oriented---open---skeletons (`near' $a_{j}^{o})$ $J_{a_{j}}^{o}$, $j \! = \!
1,\dotsc,N \! + \! 1$, in such a way that their image under the
bi-holomorphic, conformal and orientation-preserving mappings $\Phi_{a_{j}}^{
o}(z)$ are the union of the straight-line segments $\gamma_{a_{j}}^{o,l}$, $l
\! = \! 1,2,3,4$, $j \! = \! 1,\dotsc,N \! + \! 1$. Set $\zeta \! := \! \Phi_{
a_{j}}^{o}(z)$, $j \! = \! 1,\dotsc,N \! + \! 1$, and consider $\mathcal{X}^{
o}(\Phi_{a_{j}}^{o}(z)) \! := \! \Psi^{o}(\zeta)$. Recalling the properties of
$\Phi_{a_{j}}^{o}(z)$, a straightforward calculation shows that $\Psi^{o}
\colon \Phi_{a_{j}}^{o}(\mathbb{U}_{\delta_{a_{j}}}^{o}) \setminus \cup_{l=
1}^{4} \gamma_{a_{j}}^{o,l} \! \to \! \operatorname{GL}_{2}(\mathbb{C})$, $j
\! = \! 1,\dotsc,N \! + \! 1$, solves the RHPs $(\Psi^{o}(\zeta),\upsilon_{
\Psi^{o}}^{o}(\zeta),\cup_{l=1}^{4}\gamma_{a_{j}}^{o,l})$, $j \! = \! 1,
\dotsc,N \! + \! 1$, with constant jump matrices $\upsilon_{\Psi^{o}}^{o}
(\zeta)$, $j \! = \! 1,\dotsc,N \! + \! 1$, defined by
\begin{equation*}
\upsilon_{\Psi^{o}}^{o}(\zeta) \! := \!
\begin{cases}
\mathrm{I} \! + \! \sigma_{-}, &\text{$\zeta \! \in \! \gamma_{a_{j}}^{o,1}
\cup \gamma_{a_{j}}^{o,3}$,} \\
\mathrm{I} \! + \! \sigma_{+}, &\text{$\zeta \! \in \! \gamma_{a_{j}}^{o,4}$,}
\\
\mi \sigma_{2}, &\text{$\zeta \! \in \! \gamma_{a_{j}}^{o,2}$.}
\end{cases}
\end{equation*}
The solution of the latter (yet-to-be normalised) RHPs is well known; in fact,
their solution is expressed in terms of the Airy function, and is given by
(see, for example, \cite{a3,a46,a47,a49,a79})
\begin{equation*}
\Psi^{o}(\zeta) \! = \!
\begin{cases}
\Psi^{o}_{1}(\zeta), &\text{$\zeta \! \in \! \widehat{\Omega}_{a_{j}}^{o,1},
\quad j \! = \! 1,\dotsc,N \! + \! 1$,} \\
\Psi^{o}_{2}(\zeta), &\text{$\zeta \! \in \! \widehat{\Omega}_{a_{j}}^{o,2},
\quad j \! = \! 1,\dotsc,N \! + \! 1$,} \\
\Psi^{o}_{3}(\zeta), &\text{$\zeta \! \in \! \widehat{\Omega}_{a_{j}}^{o,3},
\quad j \! = \! 1,\dotsc,N \! + \! 1$,} \\
\Psi^{o}_{4}(\zeta), &\text{$\zeta \! \in \! \widehat{\Omega}_{a_{j}}^{o,4},
\quad j \! = \! 1,\dotsc,N \! + \! 1$,}
\end{cases}
\end{equation*}
where $\Psi^{o}_{k}(z)$, $k \! = \! 1,2,3,4$, are defined in Remark~4.4.
Recalling that $\Phi_{a_{j}}^{o}(z)$, $j \! = \! 1,\dotsc,N \! + \! 1$, are
bi-holomorphic, and orientation-preserving conformal mappings, with $\Phi_{a_{
j}}^{o}(a_{j}^{o}) \! = \! 0$ and $\Phi_{a_{j}}^{o}$ $(\colon \mathbb{U}_{
\delta_{a_{j}}}^{o} \! \to \! \Phi_{a_{j}}^{o}(\mathbb{U}_{\delta_{a_{j}}}^{
o}) \! =: \! \widehat{\mathbb{U}}_{\delta_{a_{j}}}^{o})$ $\colon \mathbb{U}_{
\delta_{a_{j}}}^{o} \cap J_{a_{j}}^{o} \! \to \! \Phi_{a_{j}}^{o}(\mathbb{U}_{
\delta_{a_{j}}}^{o} \cap J_{a_{j}}^{o}) \! = \! \widehat{\mathbb{U}}_{\delta_{
a_{j}}}^{o} \cap (\cup_{l=1}^{4} \gamma_{a_{j}}^{o,l})$, $j \! = \! 1,\dotsc,
N \! + \! 1$, one notes that, for any analytic maps $E_{a_{j}}^{o} \colon
\mathbb{U}_{\delta_{a_{j}}}^{o} \! \to \! \operatorname{GL}_{2}(\mathbb{C})$,
$j \! = \! 1,\dotsc,N \! + \! 1$, $\mathbb{U}_{\delta_{a_{j}}}^{o} \setminus
J_{a_{j}}^{o} \! \ni \! \zeta \! \mapsto \! E_{a_{j}}^{o}(\zeta) \Psi^{o}
(\zeta)$ also solves the latter RHPs $(\Psi^{o}(\zeta),\upsilon^{o}_{\Psi^{o}
}(\zeta),\cup_{l=1}^{4} \gamma_{a_{j}}^{o,l})$, $j \! = \! 1,\dotsc,N \! +
\! 1$: one uses this `degree of freedom' of `multiplying on the left' by a
non-degenerate, analytic, matrix-valued function in order to satisfy the
remaining asymptotic (as $n \! \to \! \infty)$ matching condition for the
parametrix, namely,
$\overset{o}{m}^{\raise-1.0ex\hbox{$\scriptstyle \infty$}}(z)(\mathcal{X}^{o}
(z))^{-1} \! =_{\underset{z \in \partial \mathbb{U}_{\delta_{a_{j}}}^{
o}}{n \to \infty}} \! \mathrm{I} \! + \! \mathcal{O}((n \! + \! 1/2)^{-1})$,
uniformly for $z \! \in \! \partial \mathbb{U}_{\delta_{a_{j}}}^{o}$, $j \! =
\! 1,\dotsc,N \! + \! 1$.

Consider, say, and without loss of generality, the regions $\Omega_{a_{j}}^{o,
1} \! := \! (\Phi_{a_{j}}^{o})^{-1}(\widehat{\Omega}_{a_{j}}^{o,1})$, $j \! =
\! 1,\dotsc,N \! + \! 1$ (Figure~5). Re-tracing the above transformations, one
shows that, for $z \! \in \! \Omega_{a_{j}}^{o,1}$ $(\subset \! \mathbb{C}_{
+})$, $j \! = \! 1,\dotsc,N \! + \! 1$, $\mathcal{X}^{o}(z) \! = \! E_{a_{j}
}^{o}(z) \Psi^{o}((\tfrac{3}{4}(n \! + \! \tfrac{1}{2}) \xi_{a_{j}}^{o}(z))^{
2/3}) \exp (\tfrac{1}{2}(n \! + \! \tfrac{1}{2})(\xi_{a_{j}}^{o}(z) \! - \!
\mi \mho_{j}^{o}) \sigma_{3})$, whence, using the expression above for $\Psi^{
o}(\zeta)$, $\zeta \! \in \! \mathbb{C}_{+} \cap \widehat{\Omega}_{a_{j}}^{o,
1}$, $j \! = \! 1,\dotsc,N \! + \! 1$, and the asymptotic expansions for
$\operatorname{Ai}(\pmb{\cdot})$ and $\operatorname{Ai}^{\prime}(\pmb{\cdot})$
(as $n \! \to \! \infty)$ given in Equations~(2.6), one arrives at
\begin{align*}
\mathcal{X}^{o}(z) \underset{\underset{z \in \partial \Omega_{a_{j}}^{o,1}
\cap \partial \mathbb{U}_{\delta_{a_{j}}}^{o}}{n \to \infty}}{=}& \, \dfrac{
1}{\sqrt{\smash[b]{2 \pi}}}E_{a_{j}}^{o}(z) \! \left(\! \left(\dfrac{3}{4}
\! \left(n \! + \! \dfrac{1}{2} \right) \! \xi_{a_{j}}^{o}(z) \right)^{2/3}
\right)^{-\frac{1}{4} \sigma_{3}} \!
\begin{pmatrix}
\me^{-\frac{\mi \pi}{6}} & \me^{\frac{\mi \pi}{3}} \\
-\me^{-\frac{\mi \pi}{6}} & -\me^{\frac{4 \pi \mi}{3}}
\end{pmatrix} \! \me^{-\frac{\mi}{2}(n+\frac{1}{2}) \mho_{j}^{o} \sigma_{3}}
\\
\times& \, \left(\mathrm{I} \! + \! \mathcal{O} \! \left((n \! + \! 1/2)^{-1}
\right) \right):
\end{align*}
demanding that, for $z \! \in \! \partial \Omega_{a_{j}}^{o,1} \cap \partial
\mathbb{U}_{\delta_{a_{j}}}^{o}$, $j \! = \! 1,\dotsc,N \! + \! 1$,
$\overset{o}{m}^{\raise-1.0ex\hbox{$\scriptstyle \infty$}}(z)(\mathcal{X}^{o}
(z))^{-1} \! =_{n \to \infty} \! \mathrm{I} \! + \! \mathcal{O}((n \! + \!
1/2)^{-1})$, one gets that
\begin{equation*}
E_{a_{j}}^{o}(z) \! = \! \dfrac{1}{\sqrt{\smash[b]{2 \mi}}}
\overset{o}{m}^{\raise-1.0ex\hbox{$\scriptstyle \infty$}}(z) \me^{\frac{\mi}{2}
(n+\frac{1}{2}) \mho_{j}^{o} \sigma_{3}} \!
\begin{pmatrix}
\mi & -\mi \\
1 & 1
\end{pmatrix} \! \left(\! \left(\dfrac{3}{4} \! \left(n \! + \! \dfrac{1}{2}
\right) \! \xi_{a_{j}}^{o}(z) \right)^{2/3} \right)^{\frac{1}{4} \sigma_{3}},
\quad j \! = \! 1,\dotsc,N \! + \! 1
\end{equation*}
(note that $\det (E_{a_{j}}^{o}(z)) \! = \! 1)$. One mimicks the above
paradigm for the remaining boundary skeletons $\partial \Omega_{a_{j}}^{o,l}
\cap \partial \mathbb{U}_{\delta_{a_{j}}}^{o}$, $l \! = \! 2,3,4$, $j \! = \!
1,\dotsc,N \! + \! 1$, and shows that the exact same formula for $E_{a_{j}}^{
o}(z)$ given above is obtained; thus, for $E_{a_{j}}^{o}(z)$, $j \! = \! 1,
\dotsc,N \! + \! 1$, as given above, one concludes that, uniformly for $z \!
\in \! \partial \mathbb{U}_{\delta_{a_{j}}}^{o}$, $j \! = \! 1,\dotsc,N \! +
\! 1$, $\overset{o}{m}^{\raise-1.0ex\hbox{$\scriptstyle \infty$}}(z)(\mathcal{
X}^{o}(z))^{-1} \! =_{\underset{z \in \partial \mathbb{U}_{\delta_{a_{j}}}^{o}
}{n \to \infty}} \! \mathrm{I} \! + \! \mathcal{O}((n \! + \! 1/2)^{-1})$.
There remains, however, the question of unimodularity, since
\begin{equation*}
\det (\mathcal{X}^{o}(z)) \! = \! \left\vert
\begin{smallmatrix}
\operatorname{Ai}(\Phi_{a_{j}}^{o}(z)) & \operatorname{Ai}(\omega^{2} \Phi_{
a_{j}}^{o}(z)) \\
\operatorname{Ai}^{\prime}(\Phi_{a_{j}}^{o}(z)) & \omega^{2} \operatorname{
Ai}^{\prime}(\omega^{2} \Phi_{a_{j}}^{o}(z))
\end{smallmatrix}
\right\vert \qquad \text{or} \qquad \left\vert
\begin{smallmatrix}
\operatorname{Ai}(\Phi_{a_{j}}^{o}(z)) & -\omega^{2} \operatorname{Ai}(\omega
\Phi_{a_{j}}^{o}(z)) \\
\operatorname{Ai}^{\prime}(\Phi_{a_{j}}^{o}(z)) & -\operatorname{Ai}^{\prime}
(\omega \Phi_{a_{j}}^{o}(z))
\end{smallmatrix}
\right\vert:
\end{equation*}
multiplying $\mathcal{X}^{o}(z)$ on the left by a constant, $\widetilde{c}$,
say, using the Wronskian relations (see Chapter~10 of \cite{a82}) 
$\operatorname{W}(\operatorname{Ai}(\lambda),\operatorname{Ai}(\omega^{2}
\lambda)) \! = \! (2 \pi)^{-1} \exp (\mi \pi/6)$ and $\operatorname{W}
(\operatorname{Ai}(\lambda),\operatorname{Ai}(\omega \lambda)) \! = \!
-(2 \pi)^{-1} \exp (-\mi \pi/6)$, and the linear dependence relation for Airy
functions, $\operatorname{Ai}(\lambda) \! + \! \omega \operatorname{Ai}(\omega
\lambda) \! + \! \omega^{2} \operatorname{Ai}(\omega^{2} \lambda) \! = \! 0$,
one shows that, upon imposing the condition $\det (\mathcal{X}^{o}(z)) \! =
\! 1$, $\widetilde{c} \! = \! (2 \pi)^{1/2} \exp (-\mi \pi/12)$. \hfill $\qed$

The above analyses lead one to the following lemma.
\begin{ccccc}
Let $\overset{o}{\mathscr{M}}^{\raise-1.0ex\hbox{$\scriptstyle \sharp$}}
\colon \mathbb{C} \setminus \Sigma_{o}^{\sharp} \! \to \! \operatorname{SL}_{
2}(\mathbb{C})$ solve the {\rm RHP}
$(\overset{o}{\mathscr{M}}^{\raise-1.0ex\hbox{$\scriptstyle \sharp$}}(z),
\overset{o}{\upsilon}^{\raise-1.0ex\hbox{$\scriptstyle \sharp$}}(z),\Sigma_{
o}^{\sharp})$ formulated in Lemma~{\rm 4.2}. Define
\begin{equation*}
\mathscr{S}_{p}^{o}(z) \! := \!
\begin{cases}
\overset{o}{m}^{\raise-1.0ex\hbox{$\scriptstyle \infty$}}(z), &\text{$z \! \in
\! \mathbb{C} \setminus \cup_{j=1}^{N+1}(\mathbb{U}_{\delta_{b_{j-1}}}^{o}
\cup \mathbb{U}_{\delta_{a_{j}}}^{o})$,} \\
\mathcal{X}^{o}(z), &\text{$z \! \in \! \cup_{j=1}^{N+1}(\mathbb{U}_{\delta_{
b_{j-1}}}^{o} \cup \mathbb{U}_{\delta_{a_{j}}}^{o})$,}
\end{cases}
\end{equation*}
where $\overset{o}{m}^{\raise-1.0ex\hbox{$\scriptstyle \infty$}} \colon
\mathbb{C} \setminus J_{o}^{\infty} \! \to \! \operatorname{SL}_{2}(\mathbb{
C})$ is characterised completely in Lemma~{\rm 4.5}, and: {\rm (1)} for $z \!
\in \! \mathbb{U}_{\delta_{b_{j-1}}}^{o}$, $j \! = \! 1,\dotsc,N \! + \! 1$,
$\mathcal{X}^{o} \colon \mathbb{U}_{\delta_{b_{j-1}}}^{o} \setminus \Sigma_{
b_{j-1}}^{o} \! \to \! \operatorname{SL}_{2}(\mathbb{C})$ solve the {\rm RHPs}
$(\mathcal{X}^{o}(z),
\overset{o}{\upsilon}^{\raise-1.0ex\hbox{$\scriptstyle \sharp$}}(z),\Sigma_{
b_{j-1}}^{o})$, $j \! = \! 1,\dotsc,N \! + \! 1$, formulated in
Lemma~{\rm 4.6;} and {\rm (2)} for $z \! \in \! \mathbb{U}_{\delta_{a_{j}}}^{
o}$, $j \! = \! 1,\dotsc,N \! + \! 1$, $\mathcal{X}^{o} \colon \mathbb{U}_{
\delta_{a_{j}}}^{o} \setminus \Sigma_{a_{j}}^{o} \! \to \! \operatorname{SL}_{
2}(\mathbb{C})$ solve the {\rm RHPs} $(\mathcal{X}^{o}(z),
\overset{o}{\upsilon}^{\raise-1.0ex\hbox{$\scriptstyle \sharp$}}(z),\Sigma_{
a_{j}}^{o})$, $j \! = \! 1,\dotsc,N \! + \! 1$, formulated in Lemma~{\rm 4.7}.
Set
\begin{equation*}
\mathscr{R}^{o}(z) \! := \!
\overset{o}{\mathscr{M}}^{\raise-1.0ex\hbox{$\scriptstyle \sharp$}}(z) \!
\left(\mathscr{S}_{p}^{o}(z) \right)^{-1},
\end{equation*}
and define the augmented contour $\Sigma_{p}^{o} \! := \! \Sigma_{o}^{\sharp}
\cup (\cup_{j=1}^{N+1}(\partial \mathbb{U}_{\delta_{b_{j-1}}}^{o} \cup \mathbb{
U}_{\delta_{a_{j}}}^{o}))$, with the orientation given in Figure~{\rm 9}. Then
$\mathscr{R}^{o} \colon \mathbb{C} \setminus \Sigma_{p}^{o} \! \to \!
\operatorname{SL}_{2}(\mathbb{C})$ solves the following {\rm RHP:} {\rm (i)}
$\mathscr{R}^{o}(z)$ is holomorphic for $z \! \in \! \mathbb{C} \setminus
\Sigma_{p}^{o};$ {\rm (ii)} $\mathscr{R}^{o}_{\pm}(z) \! := \! \lim_{
\underset{z^{\prime} \! \in \, \pm \, \mathrm{side} \, \mathrm{of} \,
\Sigma_{p}^{o}}{z^{\prime} \to z}} \mathscr{R}^{o}(z^{\prime})$ satisfy
the boundary condition
\begin{equation*}
\mathscr{R}^{o}_{+}(z) \! = \! \mathscr{R}^{o}_{-}(z) \upsilon_{\mathscr{R}}^{
o}(z), \quad z \! \in \! \Sigma_{p}^{o},
\end{equation*}
where
\begin{equation*}
\upsilon^{o}_{\mathscr{R}}(z) \! := \!
\begin{cases}
\upsilon_{\mathscr{R}}^{o,1}(z), &\text{$z \! \in \! (-\infty,b_{0}^{o} \! -
\! \delta_{b_{0}}^{o}) \! \cup \! (a_{N+1}^{o} \! + \! \delta_{a_{N+1}}^{o},
+\infty) \! =: \! \Sigma_{p}^{o,1}$,} \\
\upsilon_{\mathscr{R}}^{o,2}(z), &\text{$z \! \in \! (a_{j}^{o} \! + \!
\delta_{a_{j}}^{o},b_{j}^{o} \! - \! \delta_{b_{j}}^{o}) \! =: \! \Sigma_{p,
j}^{o,2} \! \subset \! \cup_{l=1}^{N} \Sigma_{p,l}^{o,2} \! =: \! \Sigma_{p}^{
o,2}$,} \\
\upsilon_{\mathscr{R}}^{o,3}(z), &\text{$z \! \in \! \cup_{j=1}^{N+1}(J_{j}^{
o,\smallfrown} \setminus (\mathbb{C}_{+} \cap (\mathbb{U}_{\delta_{b_{j-1}}}^{
o} \cup \mathbb{U}_{\delta_{a_{j}}}^{o}))) \! =: \! \Sigma_{p}^{o,3}$,} \\
\upsilon_{\mathscr{R}}^{o,4}(z), &\text{$z \! \in \! \cup_{j=1}^{N+1}(J_{j}^{
o,\smallsmile} \setminus (\mathbb{C}_{-} \cap (\mathbb{U}_{\delta_{b_{j-1}}}^{
o} \cup \mathbb{U}_{\delta_{a_{j}}}^{o}))) \! =: \! \Sigma_{p}^{o,4}$,} \\
\upsilon_{\mathscr{R}}^{o,5}(z), &\text{$z \! \in \! \cup_{j=1}^{N+1}(\partial
\mathbb{U}_{\delta_{b_{j-1}}}^{o} \cup \mathbb{U}_{\delta_{a_{j}}}^{o}) \! =:
\! \Sigma_{p}^{o,5}$,} \\
\mathrm{I}, &\text{$z \! \in \! \Sigma_{p}^{o} \setminus \cup_{l=1}^{5}
\Sigma_{p}^{o,l}$,}
\end{cases}
\end{equation*}
with
\begin{align*}
\upsilon_{\mathscr{R}}^{o,1}(z) \! =& \, \mathrm{I} \! + \! \me^{n(g^{o}_{+}
(z)+g^{o}_{-}(z)-\widetilde{V}(z)-\ell_{o}-\mathfrak{Q}^{+}_{\mathscr{A}}-
\mathfrak{Q}^{-}_{\mathscr{A}})} \, 
\overset{o}{m}^{\raise-1.0ex\hbox{$\scriptstyle \infty$}}(z) \sigma_{+}
(\overset{o}{m}^{\raise-1.0ex\hbox{$\scriptstyle \infty$}}(z))^{-1}, \\
\upsilon_{\mathscr{R}}^{o,2}(z) \! =& \, \mathrm{I} \! + \! \me^{-\mi (n+
\frac{1}{2}) \Omega_{j}^{o}+n(g^{o}_{+}(z)+g^{o}_{-}(z)-\widetilde{V}(z)-
\ell_{o}-\mathfrak{Q}^{+}_{\mathscr{A}}-\mathfrak{Q}^{-}_{\mathscr{A}})} \, 
\overset{o}{m}^{\raise-1.0ex\hbox{$\scriptstyle \infty$}}_{-}(z) \sigma_{+}
(\overset{o}{m}^{\raise-1.0ex\hbox{$\scriptstyle \infty$}}_{-}(z))^{-1}, \\
\upsilon_{\mathscr{R}}^{o,3}(z) \! =& \, \mathrm{I} \! + \! \me^{-4(n+
\frac{1}{2}) \pi \mi \int_{z}^{a_{N+1}^{o}} \psi_{V}^{o}(s) \, \md s} \, 
\overset{o}{m}^{\raise-1.0ex\hbox{$\scriptstyle \infty$}}(z) \sigma_{-}
(\overset{o}{m}^{\raise-1.0ex\hbox{$\scriptstyle \infty$}}(z))^{-1}, \\
\upsilon_{\mathscr{R}}^{o,4}(z) \! =& \, \mathrm{I} \! + \! \me^{4(n+
\frac{1}{2}) \pi \mi \int_{z}^{a_{N+1}^{o}} \psi_{V}^{o}(s) \, \md s} \, 
\overset{o}{m}^{\raise-1.0ex\hbox{$\scriptstyle \infty$}}(z) \sigma_{-}
(\overset{o}{m}^{\raise-1.0ex\hbox{$\scriptstyle \infty$}}(z))^{-1}, \\
\upsilon_{\mathscr{R}}^{o,5}(z) \! =& \, \mathcal{X}^{o}(z)
(\overset{o}{m}^{\raise-1.0ex\hbox{$\scriptstyle \infty$}}(z))^{-1}:
\end{align*}
{\rm (iii)} $\mathscr{R}^{o}(z) \! =_{\underset{z \in \mathbb{C} \setminus
\Sigma_{p}^{o}}{z \to 0}} \! \mathrm{I} \! + \! \mathcal{O}(z);$ and
{\rm (iv)} $\mathscr{R}^{o}(z) \! =_{\underset{z \in \mathbb{C} \setminus
\Sigma_{p}^{o}}{z \to \infty}} \! \mathcal{O}(1)$.
\end{ccccc}
\begin{figure}[tbh]
\begin{center}
\begin{pspicture}(0,0)(15,5)
\psset{xunit=1cm,yunit=1cm}
\psline[linewidth=0.6pt,linestyle=solid,linecolor=black](0,2.5)(1,2.5)
\psline[linewidth=0.6pt,linestyle=solid,linecolor=magenta,arrowsize=1.5pt 5]%
{->}(1,2.5)(2.5,2.5)
\psline[linewidth=0.6pt,linestyle=solid,linecolor=magenta](2.4,2.5)(4,2.5)
\psline[linewidth=0.6pt,linestyle=solid,linecolor=black](4,2.5)(4.7,2.5)
\psarcn[linewidth=0.6pt,linestyle=solid,linecolor=magenta,arrowsize=1.5pt 5]%
{->}(2.5,1.5){1.8}{146}{90}
\psarcn[linewidth=0.6pt,linestyle=solid,linecolor=magenta](2.5,1.5){1.8}{90}%
{34}
\psarc[linewidth=0.6pt,linestyle=solid,linecolor=magenta,arrowsize=1.5pt 5]%
{->}(2.5,3.5){1.8}{214}{270}
\psarc[linewidth=0.6pt,linestyle=solid,linecolor=magenta](2.5,3.5){1.8}{270}%
{326}
\psarcn[linewidth=0.6pt,linestyle=solid,linecolor=cyan,arrowsize=1.5pt 5]%
{->}(4,2.5){0.6}{180}{45}
\psarcn[linewidth=0.6pt,linestyle=solid,linecolor=cyan](4,2.5){0.6}{45}{0}
\psarcn[linewidth=0.6pt,linestyle=solid,linecolor=cyan](4,2.5){0.6}{360}{180}
\psarcn[linewidth=0.6pt,linestyle=solid,linecolor=cyan,arrowsize=1.5pt 5]%
{->}(1,2.5){0.6}{180}{135}
\psarcn[linewidth=0.6pt,linestyle=solid,linecolor=cyan](1,2.5){0.6}{135}{0}
\psarcn[linewidth=0.6pt,linestyle=solid,linecolor=cyan](1,2.5){0.6}{360}{180}
\rput(2.5,2.9){\makebox(0,0){$\Omega^{o,\smallfrown}_{1}$}}
\rput(2.5,2.1){\makebox(0,0){$\Omega^{o,\smallsmile}_{1}$}}
\rput(3,3.6){\makebox(0,0){$J^{o,\smallfrown}_{1}$}}
\rput(3,1.4){\makebox(0,0){$J^{o,\smallsmile}_{1}$}}
\rput(1,3.4){\makebox(0,0){$\partial \mathbb{U}^{o}_{\delta_{b_{0}}}$}}
\rput(4.1,3.4){\makebox(0,0){$\partial \mathbb{U}^{o}_{\delta_{a_{1}}}$}}
\psline[linewidth=0.6pt,linestyle=solid,linecolor=black](5.3,2.5)(6,2.5)
\psline[linewidth=0.6pt,linestyle=solid,linecolor=magenta,arrowsize=1.5pt 5]%
{->}(6,2.5)(7.5,2.5)
\psline[linewidth=0.6pt,linestyle=solid,linecolor=magenta](7.4,2.5)(9,2.5)
\psline[linewidth=0.6pt,linestyle=solid,linecolor=black](9,2.5)(9.7,2.5)
\psarcn[linewidth=0.6pt,linestyle=solid,linecolor=magenta,arrowsize=1.5pt 5]%
{->}(7.5,1.5){1.8}{146}{90}
\psarcn[linewidth=0.6pt,linestyle=solid,linecolor=magenta](7.5,1.5){1.8}{90}%
{34}
\psarc[linewidth=0.6pt,linestyle=solid,linecolor=magenta,arrowsize=1.5pt 5]%
{->}(7.5,3.5){1.8}{214}{270}
\psarc[linewidth=0.6pt,linestyle=solid,linecolor=magenta](7.5,3.5){1.8}{270}%
{326}
\psarcn[linewidth=0.6pt,linestyle=solid,linecolor=cyan,arrowsize=1.5pt 5]%
{->}(6,2.5){0.6}{180}{135}
\psarcn[linewidth=0.6pt,linestyle=solid,linecolor=cyan](6,2.5){0.6}{135}{0}
\psarcn[linewidth=0.6pt,linestyle=solid,linecolor=cyan](6,2.5){0.6}{360}{180}
\psarcn[linewidth=0.6pt,linestyle=solid,linecolor=cyan,arrowsize=1.5pt 5]%
{->}(9,2.5){0.6}{180}{45}
\psarcn[linewidth=0.6pt,linestyle=solid,linecolor=cyan](9,2.5){0.6}{45}{0}
\psarcn[linewidth=0.6pt,linestyle=solid,linecolor=cyan](9,2.5){0.6}{360}{180}
\rput(7.5,2.9){\makebox(0,0){$\Omega^{o,\smallfrown}_{j}$}}
\rput(7.5,2.1){\makebox(0,0){$\Omega^{o,\smallsmile}_{j}$}}
\rput(8,3.6){\makebox(0,0){$J^{o,\smallfrown}_{j}$}}
\rput(8,1.4){\makebox(0,0){$J^{o,\smallsmile}_{j}$}}
\rput(6,3.4){\makebox(0,0){$\partial \mathbb{U}^{o}_{\delta_{b_{j-1}}}$}}
\rput(9.1,3.4){\makebox(0,0){$\partial \mathbb{U}^{o}_{\delta_{a_{j}}}$}}
\psline[linewidth=0.6pt,linestyle=solid,linecolor=black](10.3,2.5)(11,2.5)
\psline[linewidth=0.6pt,linestyle=solid,linecolor=magenta,arrowsize=1.5pt 5]%
{->}(11,2.5)(12.5,2.5)
\psline[linewidth=0.6pt,linestyle=solid,linecolor=magenta](12.4,2.5)(14,2.5)
\psline[linewidth=0.6pt,linestyle=solid,linecolor=black](14,2.5)(15,2.5)
\psarcn[linewidth=0.6pt,linestyle=solid,linecolor=magenta,arrowsize=1.5pt 5]%
{->}(12.5,1.5){1.8}{146}{90}
\psarcn[linewidth=0.6pt,linestyle=solid,linecolor=magenta](12.5,1.5){1.8}{90}%
{34}
\psarc[linewidth=0.6pt,linestyle=solid,linecolor=magenta,arrowsize=1.5pt 5]%
{->}(12.5,3.5){1.8}{214}{270}
\psarc[linewidth=0.6pt,linestyle=solid,linecolor=magenta](12.5,3.5){1.8}{270}%
{326}
\psarcn[linewidth=0.6pt,linestyle=solid,linecolor=cyan,arrowsize=1.5pt 5]%
{->}(11,2.5){0.6}{180}{135}
\psarcn[linewidth=0.6pt,linestyle=solid,linecolor=cyan](11,2.5){0.6}{135}{0}
\psarcn[linewidth=0.6pt,linestyle=solid,linecolor=cyan](11,2.5){0.6}{360}{180}
\psarcn[linewidth=0.6pt,linestyle=solid,linecolor=cyan,arrowsize=1.5pt 5]%
{->}(14,2.5){0.6}{180}{45}
\psarcn[linewidth=0.6pt,linestyle=solid,linecolor=cyan](14,2.5){0.6}{45}{0}
\psarcn[linewidth=0.6pt,linestyle=solid,linecolor=cyan](14,2.5){0.6}{360}{180}
\rput(12.4,2.95){\makebox(0,0){$\Omega^{o,\smallfrown}_{N+1}$}}
\rput(12.4,2.05){\makebox(0,0){$\Omega^{o,\smallsmile}_{N+1}$}}
\rput(13,3.6){\makebox(0,0){$J^{o,\smallfrown}_{N+1}$}}
\rput(13,1.4){\makebox(0,0){$J^{o,\smallsmile}_{N+1}$}}
\rput(11,3.4){\makebox(0,0){$\partial \mathbb{U}^{o}_{\delta_{b_{N}}}$}}
\rput(14.35,3.4){\makebox(0,0){$\partial \mathbb{U}^{o}_{\delta_{a_{N+1}}}$}}
\psline[linewidth=0.9pt,linestyle=dotted,linecolor=darkgray](4.8,2.5)(5.2,2.5)
\psline[linewidth=0.9pt,linestyle=dotted,linecolor=darkgray](9.8,2.5)(10.2,2.5)
\rput(1,2.2){\makebox(0,0){$\scriptstyle \pmb{b_{0}^{o}}$}}
\rput(4,2.2){\makebox(0,0){$\scriptstyle \pmb{a_{1}^{o}}$}}
\rput(6,2.2){\makebox(0,0){$\scriptstyle \pmb{b_{j-1}^{o}}$}}
\rput(9,2.2){\makebox(0,0){$\scriptstyle \pmb{a_{j}^{o}}$}}
\rput(11,2.2){\makebox(0,0){$\scriptstyle \pmb{b_{N}^{o}}$}}
\rput(14,2.2){\makebox(0,0){$\scriptstyle \pmb{a_{N+1}^{o}}$}}
\rput(3.15,2.5){\makebox(0,0){$\scriptstyle \pmb{J_{1}^{o}}$}}
\rput(8.15,2.5){\makebox(0,0){$\scriptstyle \pmb{J_{j}^{o}}$}}
\rput(13.05,2.5){\makebox(0,0){$\scriptstyle \pmb{J_{N+1}^{o}}$}}
\psdots[dotstyle=*,dotscale=1.5](1,2.5)
\psdots[dotstyle=*,dotscale=1.5](4,2.5)
\psdots[dotstyle=*,dotscale=1.5](6,2.5)
\psdots[dotstyle=*,dotscale=1.5](9,2.5)
\psdots[dotstyle=*,dotscale=1.5](11,2.5)
\psdots[dotstyle=*,dotscale=1.5](14,2.5)
\end{pspicture}
\end{center}
\vspace{-1.00cm}
\caption{The augmented contour $\Sigma_{p}^{o} \! := \! \Sigma_{o}^{\sharp}
\cup (\cup_{j=1}^{N+1}(\partial \mathbb{U}_{\delta_{b_{j-1}}}^{o} \cup
\partial \mathbb{U}_{\delta_{a_{j}}}^{o}))$}
\end{figure}
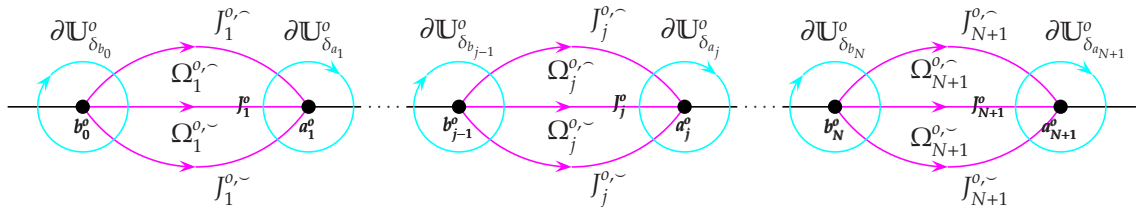

\emph{Proof.} Define the oriented, augmented skeleton $\Sigma_{p}^{o}$ as
in the Lemma: the RHP $(\mathscr{R}^{o}(z),\upsilon^{o}_{\mathscr{R}}(z),
\Sigma_{p}^{o})$ follows {}from the RHPs
$(\overset{o}{\mathscr{M}}^{\raise-1.0ex\hbox{$\scriptstyle \sharp$}}(z),
\overset{o}{\upsilon}^{\raise-1.0ex\hbox{$\scriptstyle \sharp$}}(z),\Sigma_{
o}^{\sharp})$ and $(\overset{o}{m}^{\raise-1.0ex\hbox{$\scriptstyle \infty$}}
(z),\overset{o}{\upsilon}^{\raise-1.0ex\hbox{$\scriptstyle \infty$}}(z),J_{
o}^{\infty})$ formulated in Lemmae~4.2 and~4.3, respectively, upon using the
definitions of $\mathscr{S}_{p}^{o}(z)$ and $\mathscr{R}^{o}(z)$ given in the
Lemma. \hfill $\qed$
\section{Asymptotic (as $n \! \to \! \infty)$ Solution of the RHP for
$\stackrel{o}{\mathrm{Y}} \! (z)$}
In this section, via the Beals-Coifman (BC) construction \cite{a74}, the 
(normalised at zero) RHP $(\mathscr{R}^{o}(z),\upsilon^{o}_{\mathscr{R}}
\linebreak[4]
(z),\Sigma_{p}^{o})$ formulated in Lemma~4.8 is solved asymptotically (as $n 
\! \to \! \infty)$; in particular, it is shown that, uniformly for $z \! \in 
\! \Sigma_{p}^{o}$,
\begin{equation*}
\norm{\upsilon^{o}_{\mathscr{R}}(\cdot) \! - \! \mathrm{I}}_{\cap_{p \in \{
1,2,\infty\}} \mathcal{L}^{p}_{\mathrm{M}_{2}(\mathbb{C})}(\Sigma_{p}^{o})}
\underset{n \to \infty}{=} \mathrm{I} \! + \! \mathcal{O} \! \left(f(n) \!
\left(n \! + \! 1/2 \right)^{-1} \right),
\end{equation*}
where $f(n) \! =_{n \to \infty} \! \mathcal{O}(1)$, and, subsequently, the
original \textbf{RHP2}, that is, $(\overset{o}{\mathrm{Y}}(z),\mathrm{I} \! +
\! \me^{-n \widetilde{V}(z)} \sigma_{+},\mathbb{R})$, is solved asymptotically
by re-tracing the finite sequence of RHP transformations $\mathscr{R}^{o}(z)$
(Lemmae~5.3 and~4.8) $\to \!
\overset{o}{\mathscr{M}}^{\raise-1.0ex\hbox{$\scriptstyle \sharp$}}(z)$
(Lemma~4.2) $\to \!
\overset{o}{\mathscr{M}}^{\raise-1.0ex\hbox{$\scriptstyle \flat$}}(z)$
(Proposition~4.1) $\to \! \overset{o}{\mathscr{M}}(z)$ (Lemma~3.4) $\to \!
\overset{o}{\mathrm{Y}}(z)$. The (unique) solution for $\overset{o}{\mathrm{
Y}}(z)$ then leads to the final asymptotic results for $z \pi_{2n+1}(z)$ (in
the entire complex plane), $\xi_{-n-1}^{(2n+1)}$ and $\phi_{2n+1}(z)$ (in the
entire complex plane) stated, respectively, in Theorems~2.3.1 and~2.3.2.
\begin{bbbbb}
Let $\mathscr{R}^{o} \colon \mathbb{C} \setminus \Sigma_{p}^{o} \! \to \!
\operatorname{SL}_{2}(\mathbb{C})$ solve the {\rm RHP} $(\mathscr{R}^{o}(z),
\upsilon_{\mathscr{R}}^{o}(z),\Sigma_{p}^{o})$ formulated in Lemma {\rm 4.8}.
Then:
\begin{compactenum}
\item[{\rm (1)}] for $z \! \in \! (-\infty,b_{0}^{o} \! - \! \delta_{b_{0}}^{
o}) \cup (a_{N+1}^{o} \! + \! \delta_{a_{N+1}}^{o},+\infty) \! =: \! \Sigma_{
p}^{o,1}$,
\begin{equation*}
\upsilon_{\mathscr{R}}^{o}(z) \underset{n \to \infty}{=}
\begin{cases}
\mathrm{I} \! + \! \mathcal{O}(f_{\infty}(n) \me^{-(n+\frac{1}{2})c_{\infty}
\vert z \vert}), &\text{$z \! \in \! \Sigma_{p}^{o,1} \setminus \mathbb{U}_{
0}^{o}$,} \\
\mathrm{I} \! + \! \mathcal{O}(f_{0}(n) \me^{-(n+\frac{1}{2})c_{0} \vert z
\vert^{-1}}), &\text{$z \! \in \! \Sigma_{p}^{o} \cap \mathbb{U}_{0}^{o}$,}
\end{cases}
\end{equation*}
where $c_{0},c_{\infty} \! > \! 0$, $(f_{\infty}(n))_{ij} \! =_{n \to \infty}
\! \mathcal{O}(1)$, $(f_{0}(n))_{ij} \! =_{n \to \infty} \! \mathcal{O}(1)$,
$i,j \! = \! 1,2$, and $\mathbb{U}_{0}^{o} \! := \! \{\mathstrut z \! \in \!
\mathbb{C}; \, \vert z \vert \! < \! \epsilon\}$, with $\epsilon$ some
arbitrarily fixed, sufficiently small positive real number;
\item[{\rm (2)}] for $z \! \in \! (a_{j}^{o} \! + \! \delta_{a_{j}}^{o},b_{
j}^{o} \! - \! \delta_{b_{j}}^{o}) \! =: \! \Sigma_{p,j}^{o,2} \subset \cup_{
l=1}^{N} \Sigma_{p,l}^{o,2} \! =: \! \Sigma_{p}^{o,2}$, $j \! = \! 1,\dotsc,
N$,
\begin{equation*}
\upsilon_{\mathscr{R}}^{o}(z) \underset{n \to \infty}{=}
\begin{cases}
\mathrm{I} \! + \! \mathcal{O}(f_{j}(n) \me^{-(n+\frac{1}{2})c_{j}(z \! - \!
a_{j}^{o})}), &\text{$z \! \in \! \Sigma^{o,2}_{p,j} \setminus \mathbb{U}^{
o}_{0}$,} \\
\mathrm{I} \! + \! \mathcal{O}(\widetilde{f}_{j}(n) \me^{-(n+\frac{1}{2})
\widetilde{c}_{j} \vert z \vert^{-1}}), &\text{$z \! \in \! \Sigma^{o,2}_{p,
j} \cap \mathbb{U}^{o}_{0}$,}
\end{cases}
\end{equation*}
where $c_{j},\widetilde{c}_{j} \! > \! 0$, $(f_{j}(n))_{kl} \! =_{n \to
\infty} \! \mathcal{O}(1)$, and $(\widetilde{f}_{j}(n))_{kl} \! =_{n \to
\infty} \! \mathcal{O}(1)$, $k,l \! = \! 1,2;$
\item[{\rm (3)}] for $z \! \in \! \cup_{j=1}^{N+1}(J_{j}^{o,\smallfrown}
\setminus (\mathbb{C}_{+} \cap (\mathbb{U}_{\delta_{b_{j-1}}}^{o} \cup
\mathbb{U}_{\delta_{a_{j}}}^{o}))) \! =: \! \Sigma_{p}^{o,3}$,
\begin{equation*}
\upsilon_{\mathscr{R}}^{o}(z) \underset{n \to \infty}{=} \mathrm{I} \! + \!
\mathcal{O}(\overset{\smallfrown}{f}(n) \me^{-(n+\frac{1}{2}) \overset{
\smallfrown}{c} \vert z \vert}),
\end{equation*}
where $\overset{\smallfrown}{c} \! > \! 0$ and $(\overset{\smallfrown}{f}(n))_{
ij} \! =_{n \to \infty} \! \mathcal{O}(1)$, $i,j \! = \! 1,2;$
\item[{\rm (4)}] for $z \! \in \! \cup_{j=1}^{N+1}(J_{j}^{o,\smallsmile}
\setminus (\mathbb{C}_{-} \cap (\mathbb{U}_{\delta_{b_{j-1}}}^{o} \cup
\mathbb{U}_{\delta_{a_{j}}}^{o}))) \! =: \! \Sigma_{p}^{o,4}$,
\begin{equation*}
\upsilon_{\mathscr{R}}^{o}(z) \! \underset{n \to \infty}{=} \! \mathrm{I} \!
+ \! \mathcal{O}(\overset{\smallsmile}{f}(n) \me^{-(n+\frac{1}{2}) \overset{
\smallsmile}{c} \vert z \vert}),
\end{equation*}
where $\overset{\smallsmile}{c} \! > \! 0$ and $(\overset{\smallsmile}{f}(n)
)_{ij} \! =_{n \to \infty} \! \mathcal{O}(1)$, $i,j \! = \! 1,2;$ and
\item[{\rm (5)}]
for $z \! \in \! \cup_{j=1}^{N+1}(\partial \mathbb{U}_{\delta_{b_{j-1}}}^{o}
\cup \partial \mathbb{U}_{\delta_{a_{j}}}^{o}) \! =: \! \Sigma_{p}^{o,5}$,
with $j \! = \! 1,\dotsc,N \! + \! 1$,
\begin{align*}
\upsilon_{\mathscr{R}}^{o}(z) \underset{\underset{z \in \mathbb{C}_{\pm} \cap
\partial \mathbb{U}_{\delta_{b_{j-1}}}^{o}}{n \to \infty}}{=}& \, \mathrm{I}
\! + \! \dfrac{1}{(n \! + \! \frac{1}{2}) \xi_{b_{j-1}}^{o}(z)}
\overset{o}{\mathfrak{M}}^{\raise-1.0ex\hbox{$\scriptstyle \infty$}}(z) \!
\begin{pmatrix}
\mp (s_{1}+t_{1}) & \mp \mi (s_{1}-t_{1}) \me^{\mi (n+\frac{1}{2}) \mho_{j-1}^{
o}} \\
\mp \mi (s_{1}-t_{1}) \me^{-\mi (n+\frac{1}{2}) \mho_{j-1}^{o}} & \pm (s_{1}+
t_{1})
\end{pmatrix} \\
\times& \,
(\overset{o}{\mathfrak{M}}^{\raise-1.0ex\hbox{$\scriptstyle \infty$}}(z))^{-1}
\! + \! \mathcal{O} \! \left(\dfrac{1}{((n \! + \! \frac{1}{2}) \xi_{b_{j-1}
}^{o}(z))^{2}}
\overset{o}{\mathfrak{M}}^{\raise-1.0ex\hbox{$\scriptstyle \infty$}}
(z)f_{b_{j-1}}^{o}(n)
(\overset{o}{\mathfrak{M}}^{\raise-1.0ex\hbox{$\scriptstyle \infty$}}(z))^{-1}
\right),
\end{align*}
where $\overset{o}{\mathfrak{M}}^{\raise-1.0ex\hbox{$\scriptstyle \infty$}}
(z)$ is characterised completely in Lemma {\rm 4.5}, $s_{1} \! = \! 5/72$,
$t_{1} \! = \! -7/72$, for $j \! = \! 1,\dotsc,N \! + \! 1$, $\xi_{b_{j-1}}^{
o}(z) \! = \! -2 \int_{z}^{b_{j-1}^{o}}(R_{o}(s))^{1/2}h_{V}^{o}(s) \, \md s
\! = \! (z \! - \! b_{j-1}^{o})^{3/2}G_{b_{j-1}}^{o}(z)$, with $G_{b_{j-1}}^{
o}(z)$ described completely in Lemma {\rm 4.6}, $\mho_{j-1}^{o}$ is defined
in Remark~{\rm 4.4}, and $(f_{b_{j-1}}^{o}(n))_{kl} \! =_{n \to \infty} \!
\mathcal{O}(1)$, $k,l \! = \! 1,2$, and
\begin{align*}
\upsilon_{\mathscr{R}}^{o}(z) \underset{\underset{z \in \mathbb{C}_{\pm} \cap
\partial \mathbb{U}_{\delta_{a_{j}}}^{o}}{n \to \infty}}{=}& \, \mathrm{I} \!
+ \! \dfrac{1}{(n \! + \! \frac{1}{2}) \xi_{a_{j}}^{o}(z)}
\overset{o}{\mathfrak{M}}^{\raise-1.0ex\hbox{$\scriptstyle \infty$}}(z) \!
\begin{pmatrix}
\mp (s_{1}+t_{1}) & \pm \mi (s_{1}-t_{1}) \me^{\mi (n+\frac{1}{2}) \mho_{j}^{o}
} \\
\pm \mi (s_{1}-t_{1}) \me^{-\mi (n+\frac{1}{2}) \mho_{j}^{o}} & \pm (s_{1}+t_{
1})
\end{pmatrix} \\
\times& \,
(\overset{o}{\mathfrak{M}}^{\raise-1.0ex\hbox{$\scriptstyle \infty$}}(z))^{-1}
\! + \! \mathcal{O} \! \left(\dfrac{1}{((n \! + \! \frac{1}{2}) \xi_{a_{j}}^{o}
(z))^{2}} \overset{o}{\mathfrak{M}}^{\raise-1.0ex\hbox{$\scriptstyle \infty$}}
(z)f_{a_{j}}^{o}(n)
(\overset{o}{\mathfrak{M}}^{\raise-1.0ex\hbox{$\scriptstyle \infty$}}(z))^{-1}
\right), \quad j \! = \! 1,\dotsc,N \! + \! 1,
\end{align*}
where, for $j \! = \! 1,\dotsc,N \! + \! 1$, $\xi_{a_{j}}^{o}(z) \! = \! 2
\int_{a_{j}^{o}}^{z}(R_{o}(s))^{1/2}h_{V}^{o}(s) \, \md s \! = \! (z \! -
\! a_{j}^{o})^{3/2}G_{a_{j}}^{o}(z)$, with $G_{a_{j}}^{o}(z)$ described
completely in Lemma {\rm 4.7}, and $(f_{a_{j}}^{o}(n))_{kl} \! =_{n \to
\infty} \! \mathcal{O}(1)$, $k,l \! = \! 1,2$.
\end{compactenum}
\end{bbbbb}

\emph{Proof.} Recall the definition of $\upsilon_{\mathscr{R}}^{o}(z)$ given
in Lemma~4.8. For $z \! \in \! \Sigma_{p}^{o,1} \! := \! (-\infty,b_{0}^{o} \!
- \! \delta_{b_{0}}^{o}) \cup (a_{N+1}^{o} \! + \! \delta_{a_{N+1}}^{o},+
\infty)$, recall {}from Lemma~4.8 that
\begin{equation*}
\upsilon_{\mathscr{R}}^{o}(z) \! := \! \upsilon_{\mathscr{R}}^{o,1}(z) \! = \!
\mathrm{I} \! + \! \exp \! \left(n(g^{o}_{+}(z) \! + \! g^{o}_{-}(z) \! - \!
\widetilde{V}(z) \! - \! \ell_{o} \! - \! \mathfrak{Q}^{+}_{\mathscr{A}} \! -
\! \mathfrak{Q}^{-}_{\mathscr{A}}) \right) \!
\overset{o}{m}^{\raise-1.0ex\hbox{$\scriptstyle \infty$}}(z) \sigma_{+}
(\overset{o}{m}^{\raise-1.0ex\hbox{$\scriptstyle \infty$}}(z))^{-1},
\end{equation*}
and, {}from the proof of Lemma~4.1, $g^{o}_{+}(z) \! + \! g^{o}_{-}(z) \! - \!
\widetilde{V}(z) \! - \! \ell_{o} \! - \! \mathfrak{Q}^{+}_{\mathscr{A}} \! -
\! \mathfrak{Q}^{-}_{\mathscr{A}}$ equals $-(2 \! + \! \tfrac{1}{n}) \int_{
a_{N+1}^{o}}^{z}(R_{o}(s))^{1/2} \linebreak[4]
\cdot h_{V}^{o}(s) \, \md s$ $(< \! 0)$ for $z \! \in \! (a_{N+1}^{o} \! + \!
\delta_{a_{N+1}}^{o},+\infty)$ and equals $(2 \! + \! \tfrac{1}{n}) \int_{z}^{
b_{0}^{o}}(R_{o}(s))^{1/2}h_{V}^{o}(s) \, \md s$ $(< \! 0)$ for $z \! \in \!
(-\infty,b_{0}^{o} \! - \! \delta_{b_{0}}^{o})$; hence, recalling that
$\widetilde{V} \colon \mathbb{R} \setminus \{0\} \! \to \! \mathbb{R}$, which
is regular, satisfies conditions~(2.3)--(2.5), using the asymptotic expansions
(as $\vert z \vert \! \to \! \infty$ and $\vert z \vert \! \to \! 0)$ for $g^{
o}_{+}(z) \! + \! g^{o}_{-}(z) \! - \! \widetilde{V}(z) \! - \! \ell_{o} \! -
\! \mathfrak{Q}^{+}_{\mathscr{A}} \! - \! \mathfrak{Q}^{-}_{\mathscr{A}}$
given in the proof of Lemma~3.6, that is, $g^{o}_{+}(z) \! + \! g^{o}_{-}(z)
\! - \! \widetilde{V}(z) \! - \! \ell_{o} \! - \! \mathfrak{Q}^{+}_{\mathscr{
A}} \! - \! \mathfrak{Q}^{-}_{\mathscr{A}} \! =_{\vert z \vert \to \infty} \!
(1 \! + \! \tfrac{1}{n}) \ln (z^{2} \! + \! 1) \! - \! \widetilde{V}(z) \! +
\! \mathcal{O}(1)$ and $g^{o}_{+}(z) \! + \! g^{o}_{-}(z) \! - \! \widetilde{
V}(z) \! - \! \ell_{o} \! - \! \mathfrak{Q}^{+}_{\mathscr{A}} \! - \!
\mathfrak{Q}^{-}_{\mathscr{A}} \! =_{\vert z \vert \to 0} \! \ln (z^{-2} \! +
\! 1) \! - \! \widetilde{V}(z) \! + \! \mathcal{O}(1)$, upon recalling the
expression for $\overset{o}{m}^{\raise-1.0ex\hbox{$\scriptstyle \infty$}}(z)$
given in Lemma~4.5 and noting that the respective factors $(\gamma^{o}(0))^{-
1} \gamma^{o}(z) \! \pm \! \gamma^{o}(0)(\gamma^{o}(z))^{-1}$ and $\boldsymbol{
\theta}^{o}(\pm \boldsymbol{u}^{o}(z) \! - \! \tfrac{1}{2 \pi}(n \! + \!
\tfrac{1}{2}) \boldsymbol{\Omega}^{o} \! \pm \! \boldsymbol{d}_{o})$ are
uniformly bounded (with respect to $z)$ in compact subsets outside the open
intervals surrounding the end-points of the suppport of the `odd' equilibrium
measure, defining $\mathbb{U}_{0}^{o}$ as in the Proposition, one arrives at
the asymptotic (as $n \! \to \! \infty)$ estimates for $\upsilon_{\mathscr{R}
}^{o}(z)$ on $\Sigma_{p}^{o,1} \setminus \mathbb{U}_{0}^{o} \! \ni \! z$ and
$\Sigma_{p}^{o,1} \cap \mathbb{U}_{0}^{o} \! \ni \! z$ stated in item~(1) of
the Proposition. (It should be noted that the $n$-dependence of the
$\operatorname{GL}_{2}(\mathbb{C})$-valued factors $f_{\infty}(n)$ and $f_{0}
(n)$ are inherited {}from the bounded $(\mathcal{O}(1))$ $n$-dependence of
the respective Riemann theta functions, whose corresponding series converge
absolutely and uniformly due to the fact that the associated Riemann matrix
of $\boldsymbol{\beta}^{o}$-periods, $\tau^{o}$, is pure imaginary and $-\mi
\tau^{o}$ is positive definite.)

For $z \! \in \! \Sigma_{p,j}^{o,2} \! := \! (a_{j}^{o} \! + \! \delta_{a_{j}
}^{o},b_{j}^{o} \! - \! \delta_{b_{j}}^{o})$, $j \! = \! 1,\dotsc,N$, recall
{}from Lemma~4.8 that
\begin{equation*}
\upsilon_{\mathscr{R}}^{o}(z) \! := \! \upsilon_{\mathscr{R}}^{o,2}(z) \! = \!
\mathrm{I} \! + \! \me^{-\mi (n+\frac{1}{2}) \Omega_{j}^{o}} \exp \! \left(n
(g^{o}_{+}(z) \! + \! g^{o}_{-}(z) \! - \! \widetilde{V}(z) \! - \! \ell_{o}
\! - \! \mathfrak{Q}^{+}_{\mathscr{A}} \! - \! \mathfrak{Q}^{-}_{\mathscr{A}}
) \right) \!
\overset{o}{m}^{\raise-1.0ex\hbox{$\scriptstyle \infty$}}_{-}(z) \sigma_{+}
(\overset{o}{m}^{\raise-1.0ex\hbox{$\scriptstyle \infty$}}_{-}(z))^{-1},
\end{equation*}
and, {}from the proof of Lemma~4.1, $g^{o}_{+}(z) \! + \! g^{o}_{-}(z) \! - \!
\widetilde{V}(z) \! - \! \ell_{o} \! - \! \mathfrak{Q}^{+}_{\mathscr{A}} \! -
\! \mathfrak{Q}^{-}_{\mathscr{A}} \! = \! -(2 \! + \! \tfrac{1}{n}) \int_{a_{
j}^{o}}^{z}(R_{o}(s))^{1/2}h_{V}^{o}(s) \linebreak[4]
\md s$ $(< \! 0)$. Recalling, also, that $(R_{o}(z))^{1/2} \! := \! (\prod_{k
=1}^{N+1}(z \! - \! b_{k-1}^{o})(z \! - \! a_{k}^{o}))^{1/2}$ is continuous
(and bounded) on the compact intervals $[a_{j}^{o},b_{j}^{o}] \supset \Sigma_{
p,j}^{o,2} \! \ni \! z$, $j \! = \! 1,\dotsc,N$, vanishes at the end-points
$\lbrace a_{j}^{o} \rbrace_{j=1}^{N}$ (resp., $\lbrace b_{j}^{o} \rbrace_{j=
1}^{N})$ like $(R_{o}(z))^{1/2} \! =_{z \downarrow a_{j}^{o}} \! \mathcal{O}
((z \! - \! a_{j}^{o})^{1/2})$ (resp., $(R_{o}(z))^{1/2} \! =_{z \uparrow b_{
j}^{o}} \! \mathcal{O}((z \! - \! b_{j}^{o})^{1/2}))$, and is differentiable
on the open intervals $\Sigma_{p,j}^{o,2} \! \ni \! z$, and $h_{V}^{o}(z) \!
= \! \tfrac{1}{2}(2 \! + \! \tfrac{1}{n})^{-1} \oint_{C_{\mathrm{R}}^{o}}
(\tfrac{2 \mi}{\pi s} \! + \! \tfrac{\mi \widetilde{V}^{\prime}(s)}{\pi})
(R_{o}(s))^{-1/2}(s \! - \! z)^{-1} \, \md s$ is analytic, it follows that,
for $z \! \in \! \Sigma_{p,j}^{o,2}$,
\begin{equation*}
\inf_{z \in \Sigma_{p,j}^{o,2}}(R_{o}(z))^{1/2} \! =: \! \widehat{m}_{j} \!
\leqslant \! (R_{o}(z))^{1/2} \! \leqslant \! \widehat{M}_{j} \! := \! \sup_{
z \in \Sigma_{p,j}^{o,2}}(R_{o}(z))^{1/2}, \quad j \! = \! 1,\dotsc,N;
\end{equation*}
thus, recalling the expression for
$\overset{o}{m}^{\raise-1.0ex\hbox{$\scriptstyle \infty$}}(z)$ given in
Lemma~4.5 and noting that the respective factors $(\gamma^{o}(0))^{-1}
\gamma^{o}(z) \! \pm \! \gamma^{o}(0)(\gamma^{o}(z))^{-1}$ and $\boldsymbol{
\theta}^{o}(\pm \boldsymbol{u}^{o}(z) \! - \! \tfrac{1}{2 \pi}(n \! + \!
\tfrac{1}{2}) \boldsymbol{\Omega}^{o} \! \pm \! \boldsymbol{d}_{o})$ are
uniformly bounded (with respect to $z)$ in compact subsets outside the open
intervals surrounding the end-points of the suppport of the `odd' equilibrium
measure, and defining $\mathbb{U}^{o}_{0}$ as in the Proposition, after a
straightforward integration argument, one arrives at the asymptotic (as $n \!
\to \! \infty)$ estimates for $\upsilon_{\mathscr{R}}^{o}(z)$ on $\Sigma_{p,
j}^{o,2} \setminus \mathbb{U}^{o}_{0} \! \ni \! z$ and $\Sigma^{o,2}_{p,j}
\cap \mathbb{U}^{o}_{0} \! \ni \! z$, $j \! = \! 1,\dotsc,N$, stated in
item~(2) of the Proposition (the $n$-dependence of the $\operatorname{GL}_{2}
(\mathbb{C})$-valued factors $f_{j}(n),\widetilde{f}_{j}(n)$, $j \! = \! 1,
\dotsc,N$, is inherited {}from the bounded $(\mathcal{O}(1))$ $n$-dependence
of the respective Riemann theta functions).

For $z \! \in \! \Sigma_{p}^{o,3} \! := \! \cup_{j=1}^{N+1}(J_{j}^{o,
\smallfrown} \setminus (\mathbb{C}_{+} \cap (\mathbb{U}_{\delta_{b_{j-1}}}^{o}
\cup \mathbb{U}_{\delta_{a_{j}}}^{o})))$, recall {}from Lemma~4.1 that $\Re
(\mi \int_{z}^{a_{N+1}^{o}} \psi_{V}^{o}(s) \linebreak[4]
\md s) \! > \! 0$ for $z \! \in \! \mathbb{C}_{+} \cap (\cup_{j=1}^{N+1}
\mathbb{U}_{j}^{o}) \supset \Sigma_{p}^{o,3}$, where $\mathbb{U}_{j}^{o} \!
:= \! \lbrace \mathstrut z \! \in \! \mathbb{C}^{\ast}; \, \Re (z) \! \in \!
(b_{j-1}^{o},a_{j}^{o}), \, \inf_{q \in (b_{j-1}^{o},a_{j}^{o})} \vert z \! -
\! q \vert \! < \! r_{j} \! \in \! (0,1) \rbrace$, $j \! = \! 1,\dotsc,N \! +
\! 1$, with $\mathbb{U}_{i}^{o} \cap \mathbb{U}_{j}^{o} \! = \! \varnothing$,
$i \! \not= \! j \! = \! 1,\dotsc,N \! + \! 1$, and, {}from the proof of
Lemma~4.8,
\begin{equation*}
\upsilon_{\mathscr{R}}^{o}(z) \! := \! \upsilon_{\mathscr{R}}^{o,3}(z) \! = \!
\mathrm{I} \! + \! \exp \! \left(-4 \! \left(n \! + \! \dfrac{1}{2} \right)
\! \pi \mi \int_{z}^{a_{N+1}^{o}} \psi_{V}^{o}(s) \, \md s \right) \!
\overset{o}{m}^{\raise-1.0ex\hbox{$\scriptstyle \infty$}}(z) \sigma_{-}
(\overset{o}{m}^{\raise-1.0ex\hbox{$\scriptstyle \infty$}}(z))^{-1}:
\end{equation*}
using the expression for $\overset{o}{m}^{\raise-1.0ex\hbox{$\scriptstyle
\infty$}}(z)$ given in Lemma~4.5 and noting that the respective factors
$(\gamma^{o}(0))^{-1} \linebreak[4]
\cdot \gamma^{o}(z) \! \pm \! \gamma^{o}(0)(\gamma^{o}(z))^{-1}$ and
$\boldsymbol{\theta}^{o}(\pm \boldsymbol{u}^{o}(z) \! - \! \tfrac{1}{2 \pi}
(n \! + \! \tfrac{1}{2}) \boldsymbol{\Omega}^{o} \! \pm \! \boldsymbol{d}_{
o})$ are uniformly bounded (with respect to $z)$ in compact subsets outside
the open intervals surrounding the end-points of the suppport of the `odd'
equilibrium measure, an arc-length-parametrisation argument, complemented by
an application of the Maximum Length $(\mathrm{ML})$ Theorem, leads one
directly to the asymptotic (as $n \! \to \! \infty)$ estimate for $\upsilon_{
\mathscr{R}}^{o}(z)$ on $\Sigma_{p}^{o,3} \! \ni \! z$ stated in item~(3) of
the Proposition (the $n$-dependence of the $\operatorname{GL}_{2}(\mathbb{
C})$-valued factor $\overset{\smallfrown}{f}(n)$ is inherited {}from the
bounded $(\mathcal{O}(1))$ $n$-dependence of the respective Riemann theta
functions). The above argument applies, \emph{mutatis mutandis}, for the
asymptotic estimate of $\upsilon_{\mathscr{R}}^{o}(z)$ on $\Sigma_{p}^{o,4} \!
:= \! \cup_{j=1}^{N+1}(J_{j}^{o,\smallsmile} \setminus (\mathbb{C}_{-} \cap
(\mathbb{U}_{\delta_{b_{j-1}}}^{o} \cup \mathbb{U}_{\delta_{a_{j}}}^{o}))) \!
\ni \! z$ stated in item~(4) of the Proposition.

Since the estimates in item~(5) of the Proposition are similar, consider,
say, and without loss of generality, the asymptotic (as $n \! \to \! \infty)$
estimate for $\upsilon_{\mathscr{R}}^{o}(z)$ on $\partial \mathbb{U}_{\delta_{
a_{j}}}^{o} \! \ni \! z$, $j \! = \! 1,\dotsc,N \! + \! 1$: this argument
applies, \emph{mutatis mutandis}, for the large-$n$ asymptotics of $\upsilon^{
o}_{\mathscr{R}}(z)$ on $\cup_{j=1}^{N+1} \partial \mathbb{U}_{\delta_{b_{j-1}
}}^{o} \! \ni \! z$. For $z \! \in \! \partial \mathbb{U}_{\delta_{a_{j}}}^{
o}$, $j \! = \! 1,\dotsc,N \! + \! 1$, recall {}from the proof of Lemma~4.8
that $\upsilon^{o}_{\mathscr{R}}(z) \! := \! \upsilon_{\mathscr{R}}^{o,5}(z)
\! = \! \mathcal{X}^{o}(z)
(\overset{o}{m}^{\raise-1.0ex\hbox{$\scriptstyle \infty$}}(z))^{-1}$: using
the expression for the parametrix, $\mathcal{X}^{o}(z)$, given in Lemma~4.7,
and the large-argument asymptotics for the Airy function and its derivative
given in Equations~(2.6), one shows that, for $z \! \in \! \mathbb{C}_{+} 
\cap \partial \mathbb{U}_{\delta_{a_{j}}}^{o}$, $j \! = \! 1,\dotsc,N \! + 
\! 1$,
\begin{align*}
\upsilon_{\mathscr{R}}^{o}(z) \underset{n \to \infty}{=}& \, \mathrm{I} \! 
+ \! \dfrac{\me^{-\mi \pi/3}}{(n \! + \! \frac{1}{2}) \xi_{a_{j}}^{o}(z)} 
\overset{o}{m}^{\raise-1.0ex\hbox{$\scriptstyle \infty$}}(z) \! 
\begin{pmatrix}
\mi \me^{\frac{\mi}{2}(n+\frac{1}{2}) \mho_{j}^{o}} & -\mi \me^{\frac{\mi}{2}
(n+\frac{1}{2}) \mho_{j}^{o}} \\
\me^{-\frac{\mi}{2}(n+\frac{1}{2}) \mho_{j}^{o}} & \me^{-\frac{\mi}{2}(n+
\frac{1}{2}) \mho_{j}^{o}}
\end{pmatrix} \! 
\begin{pmatrix}
-s_{1} \me^{-\frac{\mi \pi}{6}} \me^{-\frac{\mi}{2}(n+\frac{1}{2}) \mho_{j}^{
o}} & s_{1} \me^{\frac{\mi \pi}{3}} \me^{\frac{\mi}{2}(n+\frac{1}{2}) \mho_{
j}^{o}} \\
t_{1} \me^{-\frac{\mi \pi}{6}} \me^{-\frac{\mi}{2}(n+\frac{1}{2}) \mho_{j}^{
o}} & -t_{1} \me^{\frac{4 \pi \mi}{3}} \me^{\frac{\mi}{2}(n+\frac{1}{2}) 
\mho_{j}^{o}}
\end{pmatrix} \\
\times& \, (\overset{o}{m}^{\raise-1.0ex\hbox{$\scriptstyle \infty$}}(z))^{-1} 
\! + \! \mathcal{O} \! \left(\dfrac{1}{((n \! + \! \frac{1}{2}) \xi_{a_{j}}^{
o}(z))^{2}} \overset{o}{m}^{\raise-1.0ex\hbox{$\scriptstyle \infty$}}(z) \! 
\begin{pmatrix}
\ast & \ast \\
\ast & \ast
\end{pmatrix} \! (\overset{o}{m}^{\raise-1.0ex\hbox{$\scriptstyle \infty$}}
(z))^{-1} \right),
\end{align*}
where $\xi_{a_{j}}^{o}(z)$ and $\mho_{j}^{o}$, $j \! = \! 1,\dotsc,N \! + 
\! 1$, and $s_{1}$ and $t_{1}$ are defined in the Proposition, 
$\overset{o}{m}^{\raise-1.0ex\hbox{$\scriptstyle \infty$}}(z)$ is given 
in Lemma~4.5, and
$\left(
\begin{smallmatrix}
\ast & \ast \\
\ast & \ast
\end{smallmatrix}
\right) \! \in \! \operatorname{M}_{2}(\mathbb{C})$, and, for $z \! \in \! 
\mathbb{C}_{-} \cap \partial \mathbb{U}_{\delta_{a_{j}}}^{o}$, $j \! = \! 
1,\dotsc,N \! + \! 1$,
\begin{align*}
\upsilon_{\mathscr{R}}^{o}(z) \underset{n \to \infty}{=}& \, \mathrm{I} \! 
+ \! \dfrac{\me^{-\mi \pi/3}}{(n \! + \! \frac{1}{2}) \xi_{a_{j}}^{o}(z)}
\overset{o}{m}^{\raise-1.0ex\hbox{$\scriptstyle \infty$}}(z) \! 
\begin{pmatrix}
\mi \me^{-\frac{\mi}{2}(n+\frac{1}{2}) \mho_{j}^{o}} & -\mi \me^{-\frac{\mi}{2}
(n+\frac{1}{2}) \mho_{j}^{o}} \\
\me^{\frac{\mi}{2}(n+\frac{1}{2}) \mho_{j}^{o}} & \me^{\frac{\mi}{2}(n+\frac{
1}{2}) \mho_{j}^{o}}
\end{pmatrix} \! 
\begin{pmatrix}
-s_{1} \me^{-\frac{\mi \pi}{6}} \me^{\frac{\mi}{2}(n+\frac{1}{2}) \mho_{j}^{
o}} & s_{1} \me^{\frac{\mi \pi}{3}} \me^{-\frac{\mi}{2}(n+\frac{1}{2}) \mho_{
j}^{o}} \\
t_{1} \me^{-\frac{\mi \pi}{6}} \me^{\frac{\mi}{2}(n+\frac{1}{2}) \mho_{j}^{
o}} & -t_{1} \me^{\frac{4 \pi \mi}{3}} \me^{-\frac{\mi}{2}(n+\frac{1}{2})
\mho_{j}^{o}}
\end{pmatrix} \\
\times& \, (\overset{o}{m}^{\raise-1.0ex\hbox{$\scriptstyle \infty$}}(z))^{-1} 
\! + \! \mathcal{O} \! \left(\dfrac{1}{((n \! + \! \frac{1}{2}) \xi_{a_{j}}^{
o}(z))^{2}} \overset{o}{m}^{\raise-1.0ex\hbox{$\scriptstyle \infty$}}(z) \! 
\begin{pmatrix}
\ast & \ast \\
\ast & \ast
\end{pmatrix} \! (\overset{o}{m}^{\raise-1.0ex\hbox{$\scriptstyle \infty$}}
(z))^{-1} \right).
\end{align*}
Upon recalling the formula for
$\overset{o}{m}^{\raise-1.0ex\hbox{$\scriptstyle \infty$}}(z)$ in terms of
$\overset{o}{\mathfrak{M}}^{\raise-1.0ex\hbox{$\scriptstyle \infty$}}(z)$
given in Lemma~4.5, and noting that the respective factors $(\gamma^{o}(0))^{-
1} \gamma^{o}(z) \! \pm \! \gamma^{o}(0)(\gamma^{o}(z))^{-1}$ and $\boldsymbol{
\theta}^{o}(\pm \boldsymbol{u}^{o}(z) \! - \! \tfrac{1}{2 \pi}(n \! + \!
\tfrac{1}{2}) \boldsymbol{\Omega}^{o} \! \pm \! \boldsymbol{d}_{o})$ are
uniformly bounded (with respect to $z)$ in compact subsets outside the open
intervals surrounding the end-points of the suppport of the `odd' equilibrium
measure, after a straightforward matrix-multiplication argument, one arrives
at the asymptotic (as $n \! \to \! \infty)$ estimates for $\upsilon^{o}_{
\mathscr{R}}(z)$ on $\partial \mathbb{U}_{\delta_{a_{j}}}^{o} \! \ni \! z$,
$j \! = \! 1,\dotsc,N \! + \! 1$, stated in item~(5) of the Proposition (the
$n$-dependence of the $\operatorname{GL}_{2}(\mathbb{C})$-valued factors
$f_{a_{j}}^{o}(n)$, $j \! = \! 1,\dotsc,N \! + \! 1$, is inherited {}from the
bounded $(\mathcal{O}(1))$ $n$-dependence of the respective Riemann theta
functions). \hfill $\qed$
\begin{aaaaa}
For an oriented contour $D \! \subset \! \mathbb{C}$, let $\mathscr{N}_{q}
(D)$ denote the set of all bounded linear operators {}from $\mathcal{L}^{q}_{
\mathrm{M}_{2}(\mathbb{C})}(D)$ into $\mathcal{L}^{q}_{\mathrm{M}_{2}(\mathbb{
C})}(D)$, $q \! \in \! \lbrace 1,2,\infty \rbrace$.
\end{aaaaa}

Since the analysis that follows relies substantially on the BC \cite{a74} 
construction for the solution of a matrix (and suitably normalised) RHP on 
an oriented and unbounded contour, it is convenient to present, with some 
requisite preamble, a succinct and self-contained synopsis of it at this 
juncture. One agrees to call a contour $\Gamma^{\sharp}$ \emph{oriented} if:
\begin{compactenum}
\item[(1)] $\mathbb{C} \setminus \Gamma^{\sharp}$ has finitely many open 
connected components;
\item[(2)] $\mathbb{C} \setminus \Gamma^{\sharp}$ is the disjoint union of 
two, possibly disconnected, open regions, denoted by $\boldsymbol{\mho}^{+}$ 
and $\boldsymbol{\mho}^{-}$;
\item[(3)] $\Gamma^{\sharp}$ may be viewed as either the positively oriented
boundary for $\boldsymbol{\mho}^{+}$ or the negatively oriented boundary for
$\boldsymbol{\mho}^{-}$ ($\mathbb{C} \setminus \Gamma^{\sharp}$ is coloured
by two colours, $\pm)$.
\end{compactenum}
Let $\Gamma^{\sharp}$, as a closed set, be the union of finitely many 
oriented, simple, piecewise-smooth arcs. Denote the set of all 
self-intersections of $\Gamma^{\sharp}$ by $\widehat{\Gamma}^{\sharp}$ 
(with $\mathrm{card}(\widehat{\Gamma}^{\sharp}) \! < \! \infty$ assumed 
throughout). Set $\widetilde{\Gamma}^{\sharp} \! := \! \Gamma^{\sharp} 
\setminus \widehat{\Gamma}^{\sharp}$. The BC \cite{a74} construction for 
the solution of a (matrix) RHP, in the absence of a discrete spectrum and 
spectral singularities \cite{a88} (see, also, \cite{a75,a76,a89,a90,a91}), 
on an oriented contour $\Gamma^{\sharp}$ consists of finding function 
$\mathcal{Y}(z) \colon \mathbb{C} \setminus \Gamma^{\sharp} \! \to \! 
\mathrm{M}_{2}(\mathbb{C})$ such that:
\begin{compactenum}
\item[(1)] $\mathcal{Y}(z)$ is holomorphic for $z \! \in \! \mathbb{C} 
\setminus \Gamma^{\sharp}$, $\mathcal{Y}(z) \! \! \upharpoonright_{\mathbb{C}
\setminus \Gamma^{\sharp}}$ has a continuous extension ({}from `above' and
`below') to $\widetilde{\Gamma}^{\sharp}$, and $\lim_{\genfrac{}{}{0pt}{2}{
z^{\prime} \to z}{z^{\prime} \, \in \, \pm \, \mathrm{side} \, \mathrm{of} \,
\widetilde{\Gamma}^{\sharp}}} \int_{\widetilde{\Gamma}^{\sharp}} \vert
\mathcal{Y}(z^{\prime}) \! - \! \mathcal{Y}_{\pm}(z) \vert^{2} \, \vert \md z
\vert \! = \! 0$;
\item[(2)] $\mathcal{Y}_{\pm}(z) \! := \! \lim_{\underset{z^{\prime} \, \in \,
\pm \, \mathrm{side} \, \mathrm{of} \, \widetilde{\Gamma}^{\sharp}}{z^{\prime}
\to z}} \mathcal{Y}(z^{\prime})$ satisfy $\mathcal{Y}_{+}(z) \! = \! \mathcal{
Y}_{-}(z) \upsilon (z)$, $z \! \in \! \widetilde{\Gamma}^{\sharp}$, for some
(smooth) `jump' matrix $\upsilon \colon \widetilde{\Gamma}^{\sharp} \! \to \!
\mathrm{GL}_{2}(\mathbb{C})$; and
\item[(3)] for arbitrarily fixed $\lambda_{o} \! \in \! \mathbb{C}$, and
uniformly with respect to $z$, $\mathcal{Y}(z) \! =_{\underset{z \in \mathbb{
C} \setminus \Gamma^{\sharp}}{z \to \lambda_{o}}} \! \mathrm{I} \! + \! o(1)$,
where $o(1) \! = \! \mathcal{O}(z \! - \! \lambda_{o})$ if $\lambda_{o}$ is
finite, and $o(1) \! = \! \mathcal{O}(z^{-1})$ if $\lambda_{o}$ is the point
at infinity).
\end{compactenum}
(Condition~(3) is referred to as the \emph{normalisation condition}, and is
necessary in order to prove uniqueness of the associated RHP: one says that
the RHP is `normalised at $\lambda_{o}$'.) Let $\upsilon (z) \! := \!
(\mathrm{I} \! - \! w_{-}(z))^{-1}(\mathrm{I} \! + \! w_{+}(z))$, $z \! \in
\! \widetilde{\Gamma}^{\sharp}$, be a (bounded algebraic) factorisation for
$\upsilon (z)$, where $w_{\pm}(z)$ are some upper/lower, or lower/upper,
triangular matrices (depending on the orientation of $\Gamma^{\sharp})$, and
$w_{\pm}(z) \! \in \! \cap_{p \in \{2,\infty\}} \mathcal{L}^{p}_{\mathrm{M}_{
2}(\mathbb{C})} \linebreak[4]
(\widetilde{\Gamma}^{\sharp})$ (if $\widetilde{\Gamma}^{\sharp}$ is unbounded,
one requires that $w_{\pm}(z) \! =_{\genfrac{}{}{0pt}{2}{z \to \infty}{z \in
\widetilde{\Gamma}^{\sharp}}} \! \pmb{0})$. Define $w(z) \! := \! w_{+}(z)
\! + \! w_{-}(z)$, and introduce the (normalised at $\lambda_{o})$ Cauchy
operators
\begin{equation*}
\mathcal{L}^{2}_{\mathrm{M}_{2}(\mathbb{C})}(\Gamma^{\sharp}) \! \ni \! f \!
\mapsto \! (C^{\lambda_{o}}_{\pm}f)(z) :=
\lim_{\genfrac{}{}{0pt}{2}{z^{\prime} \to z}{z^{\prime} \, \in \, \pm \,
\mathrm{side} \, \mathrm{of} \, \Gamma^{\sharp}}} \int_{\Gamma^{\sharp}}
\dfrac{(z^{\prime} \! - \! \lambda_{o})f(\zeta)}{(\zeta \! - \! \lambda_{o})
(\zeta \! - \! z^{\prime})} \, \dfrac{\md \zeta}{2 \pi \mi},
\end{equation*}
where $\tfrac{(z-\lambda_{o})}{(\zeta -\lambda_{o})(\zeta -z)} \, \tfrac{\md
\zeta}{2 \pi \mi}$ is the Cauchy kernel normalised at $\lambda_{o}$ (which
reduces to the `standard' Cauchy kernel, that is, $\tfrac{1}{\zeta -z} \,
\tfrac{\md \zeta}{2 \pi \mi}$, in the limit $\lambda_{o} \! \to \! \infty)$,
with $C^{\lambda_{o}}_{\pm} \colon \mathcal{L}^{2}_{\mathrm{M}_{2}(\mathbb{
C})}(\Gamma^{\sharp}) \! \to \! \mathcal{L}^{2}_{\mathrm{M}_{2}(\mathbb{C})}
(\Gamma^{\sharp})$ bounded in operator norm\footnote{$\vert \vert C^{\lambda_{
o}}_{\pm} \vert \vert_{\mathscr{N}_{2}(\Gamma^{\sharp})} \! < \! \infty$.},
and $\vert \vert (C^{\lambda_{o}}_{\pm}f)(\cdot) \vert \vert_{\mathcal{L}^{
2}_{\mathrm{M}_{2}(\mathbb{C})}(\Gamma^{\sharp})} \! \leqslant \!
\mathrm{const.} \vert \vert f(\cdot) \vert \vert_{\mathcal{L}^{2}_{\mathrm{
M}_{2}(\mathbb{C})}(\Gamma^{\sharp})}$. Introduce the BC operator $C^{
\lambda_{o}}_{w}$:
\begin{equation*}
\mathcal{L}^{2}_{\mathrm{M}_{2}(\mathbb{C})}(\Gamma^{\sharp}) \! \ni \! f \!
\mapsto \! C^{\lambda_{o}}_{w}f \! := \! C^{\lambda_{o}}_{+}(fw_{-}) \! + \!
C^{\lambda_{o}}_{-}(fw_{+}),
\end{equation*}
which, for $w_{\pm} \! \in \! \mathcal{L}^{\infty}_{\mathrm{M}_{2}(\mathbb{C})}
(\Gamma^{\sharp})$, is bounded {}from $\mathcal{L}^{2}_{\mathrm{M}_{2}(\mathbb{
C})}(\Gamma^{\sharp}) \! \to \! \mathcal{L}^{2}_{\mathrm{M}_{2}(\mathbb{C})}
(\Gamma^{\sharp})$, that is, $\vert \vert C^{\lambda_{o}}_{w} \vert \vert_{
\mathscr{N}_{2}(\Gamma^{\sharp})} \! < \! \infty$; moreover, since $\mathbb{C}
\setminus \Gamma^{\sharp}$ can be coloured by the two colours $\pm$, $C^{
\lambda_{o}}_{\pm}$ are complementary projections \cite{a2,a75,a89,a90}, that
is, $(C^{\lambda_{o}}_{+})^{2} \! = \! C^{\lambda_{o}}_{+}$, $(C^{\lambda_{
o}}_{-})^{2} \! = \! -C^{\lambda_{o}}_{-}$, $C^{\lambda_{o}}_{+}C^{\lambda_{
o}}_{-} \! = \! C^{\lambda_{o}}_{-}C^{\lambda_{o}}_{+} \! = \! \underline{
\pmb{0}}$ (the null operator), and $C^{\lambda_{o}}_{+} \! - \! C^{\lambda_{
o}}_{-} \! = \! \id$ (the identity operator). (In the case that $C^{\lambda_{
o}}_{+}$ and $-C^{\lambda_{o}}_{-}$ are complementary, the contour $\Gamma^{
\sharp}$ can always be oriented in such a way that the $\pm$ regions lie on
the $\pm$ sides of the contour, respectively.) The solution of the above
(normalised at $\lambda_{o})$ RHP is given by the following integral
representation.
\begin{ccccc}[Beals and Coifman {\rm \cite{a74}}]
Set
\begin{equation*}
\mu_{\lambda_{o}}(z) \! = \! \mathcal{Y}_{+}(z) \! \left(\mathrm{I} \! + \!
w_{+}(z) \right)^{-1} \! = \! \mathcal{Y}_{-}(z) \! \left(\mathrm{I} \! - \!
w_{-}(z) \right)^{-1}, \quad z \! \in \! \Gamma^{\sharp}.
\end{equation*}
If $\mu_{\lambda_{o}} \! \in \! \mathrm{I} \! + \! \mathcal{L}^{2}_{\mathrm{
M}_{2}(\mathbb{C})}(\Gamma^{\sharp})$ solves the linear singular integral
equation
\begin{equation*}
(\id \! - \! C^{\lambda_{o}}_{w})(\mu_{\lambda_{o}}(z) \! - \! \mathrm{I}) \!
= \! C_{w}^{\lambda_{o}} \mathrm{I} \! = \! C^{\lambda_{o}}_{+}(w_{-}(z)) \!
+ \! C^{\lambda_{o}}_{-}(w_{+}(z)), \quad z \! \in \! \Gamma^{\sharp},
\end{equation*}
where $\id$ is the identity operator on $\mathcal{L}^{2}_{\mathrm{M}_{2}
(\mathbb{C})}(\Gamma^{\sharp})$, then the solution of the {\rm RHP}
$(\mathcal{Y}(z),\upsilon (z),\Gamma^{\sharp})$ is given by
\begin{equation*}
\mathcal{Y}(z) \! = \! \mathrm{I} \! + \! \int\nolimits_{\Gamma^{\sharp}}
\dfrac{(z \! - \! \lambda_{o}) \mu_{\lambda_{o}}(\zeta)w(\zeta)}{(\zeta \! -
\! \lambda_{o})(\zeta \! - \! z)} \, \dfrac{\md \zeta}{2 \pi \mi}, \quad z
\! \in \! \mathbb{C} \setminus \Gamma^{\sharp},
\end{equation*}
where $\mu_{\lambda_{o}}(z) \! := \! ((\id \! - \! C^{\lambda_{o}}_{w})^{-1}
\mathrm{I})(z)$\footnote{The linear singular integral equation for $\mu_{
\lambda_{o}}(\pmb{\cdot})$ stated in this Lemma~5.1 is well defined in
$\mathcal{L}^{2}_{\mathrm{M}_{2}(\mathbb{C})}(\Gamma^{\sharp})$ provided that
$w_{\pm}(\pmb{\cdot}) \! \in \! \mathcal{L}^{2}_{\mathrm{M}_{2}(\mathbb{C})}
(\Gamma^{\sharp}) \cap \mathcal{L}^{\infty}_{\mathrm{M}_{2}(\mathbb{C})}
(\Gamma^{\sharp})$; furthermore, it is assumed that the associated RHP
$(\mathcal{Y}(z),\upsilon (z),\Gamma^{\sharp})$ is solvable, that is,
$\mathrm{dim} \mathrm{ker}(\id \! - \! C^{\lambda_{o}}_{w}) \! = \!
\mathrm{dim} \left\lbrace \mathstrut \phi \! \in \! \mathcal{L}^{2}_{\mathrm{
M}_{2}(\mathbb{C})}(\Gamma^{\sharp}); \, (\id \! - \! C^{\lambda_{o}}_{w})
\phi \! = \! \underline{\pmb{0}} \right\rbrace \! = \! \mathrm{dim} \,
\varnothing \! = \! 0$ $(\Rightarrow  (\id \! - \! C^{\lambda_{o}}_{w})^{-1}
\! \! \upharpoonright_{\mathcal{L}^{2}_{\mathrm{M}_{2}(\mathbb{C})}(\Gamma^{
\sharp})}$ exists).}.
\end{ccccc}

Recall that $\mathscr{R}^{o} \colon \mathbb{C} \setminus \Sigma_{p}^{o} \! \to
\! \operatorname{SL}_{2}(\mathbb{C})$, which solves the RHP $(\mathscr{R}^{o}
(z),\upsilon_{\mathscr{R}}^{o}(z),\Sigma_{p}^{o})$ formulated in Lemma 4.8,
is normalised at zero, that is, $\mathscr{R}^{o}(0) \! = \! \mathrm{I}$.
Removing {}from the specification of the RHP $(\mathscr{R}^{o}(z),\upsilon_{
\mathscr{R}}^{o}(z),\linebreak[4]
\Sigma_{p}^{o})$ the oriented skeletons on which the jump matrix, $\upsilon_{
\mathscr{R}}^{o}(z)$, is equal to $\mathrm{I}$, in particular (cf. Lemma 4.8),
the oriented skeleton $\Sigma_{p}^{o} \setminus \cup_{l=1}^{5} \Sigma_{p}^{o,
l}$, and setting $\Sigma_{p}^{o} \setminus (\Sigma_{p}^{o} \setminus \cup_{l=
1}^{5} \Sigma_{p}^{o,l}) \! =: \! \widetilde{\Sigma}_{p}^{o}$ (see Figure~10),
\begin{figure}[tbh]
\begin{center}
\begin{pspicture}(0,0)(15,5)
\psset{xunit=1cm,yunit=1cm}
\psdots[dotstyle=*,dotscale=1.5](1.3,2.5)
\psdots[dotstyle=*,dotscale=1.5](3.7,2.5)
\psdots[dotstyle=*,dotscale=1.5](6.3,2.5)
\psdots[dotstyle=*,dotscale=1.5](8.7,2.5)
\psdots[dotstyle=*,dotscale=1.5](11.3,2.5)
\psdots[dotstyle=*,dotscale=1.5](13.7,2.5)
\psline[linewidth=0.6pt,linestyle=solid,linecolor=black](0.4,2.5)(0.7,2.5)
\psline[linewidth=0.6pt,linestyle=solid,linecolor=black,arrowsize=1.5pt 3]%
{->}(0,2.5)(0.4,2.5)
\psline[linewidth=0.6pt,linestyle=solid,linecolor=black,arrowsize=1.5pt 3]%
{->}(4.3,2.5)(4.55,2.5)
\psline[linewidth=0.6pt,linestyle=solid,linecolor=black](4.55,2.5)(4.75,2.5)
\psarcn[linewidth=0.6pt,linestyle=solid,linecolor=magenta,arrowsize=1.5pt 5]%
{->}(2.5,1.5){1.8}{121}{89}
\psarcn[linewidth=0.6pt,linestyle=solid,linecolor=magenta](2.5,1.5){1.8}{90}%
{59}
\psarc[linewidth=0.6pt,linestyle=solid,linecolor=magenta,arrowsize=1.5pt 5]%
{->}(2.5,3.5){1.8}{239}{271}
\psarc[linewidth=0.6pt,linestyle=solid,linecolor=magenta](2.5,3.5){1.8}{270}%
{301}
\psarcn[linewidth=0.6pt,linestyle=solid,linecolor=cyan,arrowsize=1.5pt 5]%
{->}(1.3,2.5){0.6}{180}{135}
\psarcn[linewidth=0.6pt,linestyle=solid,linecolor=cyan](1.3,2.5){0.6}{135}{0}
\psarcn[linewidth=0.6pt,linestyle=solid,linecolor=cyan](1.3,2.5){0.6}{360}{180}
\psarcn[linewidth=0.6pt,linestyle=solid,linecolor=cyan,arrowsize=1.5pt 5]%
{->}(3.7,2.5){0.6}{180}{45}
\psarcn[linewidth=0.6pt,linestyle=solid,linecolor=cyan](3.7,2.5){0.6}{45}{0}
\psarcn[linewidth=0.6pt,linestyle=solid,linecolor=cyan](3.7,2.5){0.6}{360}{180}
\psline[linewidth=0.6pt,linestyle=solid,linecolor=black](5.5,2.5)(5.7,2.5)
\psline[linewidth=0.6pt,linestyle=solid,linecolor=black,arrowsize=1.5pt 3]%
{->}(5.25,2.5)(5.5,2.5)
\psline[linewidth=0.6pt,linestyle=solid,linecolor=black,arrowsize=1.5pt 3]%
{->}(9.3,2.5)(9.55,2.5)
\psline[linewidth=0.6pt,linestyle=solid,linecolor=black](9.55,2.5)(9.75,2.5)
\psarcn[linewidth=0.6pt,linestyle=solid,linecolor=magenta,arrowsize=1.5pt 5]%
{->}(7.5,1.5){1.8}{121}{89}
\psarcn[linewidth=0.6pt,linestyle=solid,linecolor=magenta](7.5,1.5){1.8}{90}%
{59}
\psarc[linewidth=0.6pt,linestyle=solid,linecolor=magenta,arrowsize=1.5pt 5]%
{->}(7.5,3.5){1.8}{239}{271}
\psarc[linewidth=0.6pt,linestyle=solid,linecolor=magenta](7.5,3.5){1.8}{270}%
{301}
\psarcn[linewidth=0.6pt,linestyle=solid,linecolor=cyan,arrowsize=1.5pt 5]%
{->}(6.3,2.5){0.6}{180}{135}
\psarcn[linewidth=0.6pt,linestyle=solid,linecolor=cyan](6.3,2.5){0.6}{135}{0}
\psarcn[linewidth=0.6pt,linestyle=solid,linecolor=cyan](6.3,2.5){0.6}{360}{180}
\psarcn[linewidth=0.6pt,linestyle=solid,linecolor=cyan,arrowsize=1.5pt 5]%
{->}(8.7,2.5){0.6}{180}{45}
\psarcn[linewidth=0.6pt,linestyle=solid,linecolor=cyan](8.7,2.5){0.6}{45}{0}
\psarcn[linewidth=0.6pt,linestyle=solid,linecolor=cyan](8.7,2.5){0.6}{360}{180}
\psline[linewidth=0.6pt,linestyle=solid,linecolor=black](10.5,2.5)(10.7,2.5)
\psline[linewidth=0.6pt,linestyle=solid,linecolor=black,arrowsize=1.5pt 3]%
{->}(10.25,2.5)(10.5,2.5)
\psline[linewidth=0.6pt,linestyle=solid,linecolor=black,arrowsize=1.5pt 3]%
{->}(14.3,2.5)(14.65,2.5)
\psline[linewidth=0.6pt,linestyle=solid,linecolor=black](14.65,2.5)(14.95,2.5)
\psarcn[linewidth=0.6pt,linestyle=solid,linecolor=magenta,arrowsize=1.5pt 5]%
{->}(12.5,1.5){1.8}{121}{89}
\psarcn[linewidth=0.6pt,linestyle=solid,linecolor=magenta](12.5,1.5){1.8}{90}%
{59}
\psarc[linewidth=0.6pt,linestyle=solid,linecolor=magenta,arrowsize=1.5pt 5]%
{->}(12.5,3.5){1.8}{239}{271}
\psarc[linewidth=0.6pt,linestyle=solid,linecolor=magenta](12.5,3.5){1.8}{270}%
{301}
\psarcn[linewidth=0.6pt,linestyle=solid,linecolor=cyan,arrowsize=1.5pt 5]%
{->}(11.3,2.5){0.6}{180}{135}
\psarcn[linewidth=0.6pt,linestyle=solid,linecolor=cyan](11.3,2.5){0.6}{135}{0}
\psarcn[linewidth=0.6pt,linestyle=solid,linecolor=cyan](11.3,2.5){0.6}{360}%
{180}
\psarcn[linewidth=0.6pt,linestyle=solid,linecolor=cyan,arrowsize=1.5pt 5]%
{->}(13.7,2.5){0.6}{180}{45}
\psarcn[linewidth=0.6pt,linestyle=solid,linecolor=cyan](13.7,2.5){0.6}{45}{0}
\psarcn[linewidth=0.6pt,linestyle=solid,linecolor=cyan](13.7,2.5){0.6}{360}%
{180}
\psline[linewidth=0.9pt,linestyle=dotted,linecolor=darkgray](4.8,2.5)(5.2,2.5)
\psline[linewidth=0.9pt,linestyle=dotted,linecolor=darkgray](9.8,2.5)(10.2,2.5)
\rput(1.3,2.2){\makebox(0,0){$\scriptstyle \pmb{b_{0}^{o}}$}}
\rput(3.7,2.2){\makebox(0,0){$\scriptstyle \pmb{a_{1}^{o}}$}}
\rput(6.3,2.2){\makebox(0,0){$\scriptstyle \pmb{b_{j-1}^{o}}$}}
\rput(8.7,2.2){\makebox(0,0){$\scriptstyle \pmb{a_{j}^{o}}$}}
\rput(11.3,2.2){\makebox(0,0){$\scriptstyle \pmb{b_{N}^{o}}$}}
\rput(13.7,2.2){\makebox(0,0){$\scriptstyle \pmb{a_{N+1}^{o}}$}}
\end{pspicture}
\end{center}
\vspace{-1.05cm}
\caption{Oriented skeleton $\widetilde{\Sigma}_{p}^{o} \! := \! \Sigma_{p}^{o}
\setminus (\Sigma_{p}^{o} \setminus \cup_{l=1}^{5} \Sigma_{p}^{o,l})$}
\end{figure}
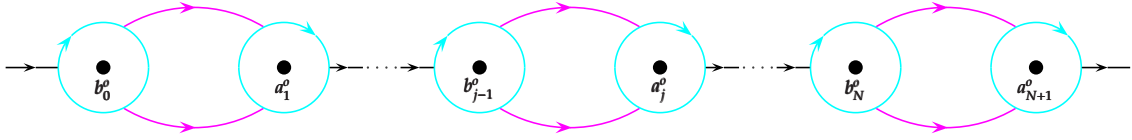
one arrives at the equivalent RHP $(\mathscr{R}^{o}(z),\upsilon_{\mathscr{R}
}^{o}(z),\widetilde{\Sigma}_{p}^{o})$ for $\mathscr{R}^{o} \colon \mathbb{C}
\setminus \widetilde{\Sigma}_{p}^{o} \! \to \! \operatorname{SL}_{2}(\mathbb{
C})$ (the normalisation at zero, of course, remains unchanged). Via the BC
\cite{a74} construction discussed above, write, for $\upsilon_{\mathscr{R}}^{
o} \colon \widetilde{\Sigma}_{p}^{o} \! \to \! \operatorname{SL}_{2}(\mathbb{
C})$, the (bounded algebraic) factorisation
\begin{equation*}
\upsilon_{\mathscr{R}}^{o}(z) \! := \! \left(\mathrm{I} \! - \! w_{-}^{
\Sigma_{\mathscr{R}}^{o}}(z) \right)^{-1} \! \left(\mathrm{I} \! + \! w_{+}^{
\Sigma_{\mathscr{R}}^{o}}(z) \right), \quad z \! \in \! \widetilde{\Sigma}_{
p}^{o}:
\end{equation*}
taking the (so-called) trivial factorisation \cite{a76} (see pp.~293 and~294, 
\emph{Proof of Theorem}~3.14 and \emph{Proposition}~1.9; see, also, 
\cite{a90,a91}) $w^{\Sigma_{\mathscr{R}}^{o}}_{-}(z) \! \equiv \! \mathbf{0}$, 
whence $\upsilon_{\mathscr{R}}^{o}(z) \! = \! \mathrm{I} \! + \! w^{\Sigma_{
\mathscr{R}}^{o}}_{+}(z)$, $z \! \in \! \widetilde{\Sigma}_{p}^{o}$, it 
follows {}from Lemma~5.1 that, upon normalising the Cauchy (integral) 
operator(s) at zero (take the limit $\lambda_{o} \! \to \! 0$ in Lemma~5.1), 
the $(\operatorname{SL}_{2}(\mathbb{C})$-valued) integral representation 
for the---unique---solution of the equivalent RHP $(\mathscr{R}^{o}(z),
\linebreak[4]
\upsilon_{\mathscr{R}}^{o}(z),\widetilde{\Sigma}_{p}^{o})$ is
\begin{equation}
\mathscr{R}^{o}(z) \! = \! \mathrm{I} \! + \! \int_{\widetilde{\Sigma}_{p}^{
o}} \dfrac{z \mu^{\Sigma_{\mathscr{R}}^{o}}(s) w^{\Sigma_{\mathscr{R}}^{o}}_{
+}(s)}{s(s \! - \! z)} \, \dfrac{\md s}{2 \pi \mi}, \quad z \! \in \! \mathbb{
C} \setminus \widetilde{\Sigma}_{p}^{o},
\end{equation}
where $\mu^{\Sigma_{\mathscr{R}}^{o}}(\pmb{\cdot}) \! \in \! \mathrm{I} \! +
\! \mathcal{L}^{2}_{\mathrm{M}_{2}(\mathbb{C})}(\widetilde{\Sigma}_{p}^{o})$
solves the (linear) singular integral equation
\begin{equation*}
(\id \! - \! C^{0}_{w^{\Sigma_{\mathscr{R}}^{o}}}) \mu^{\Sigma_{\mathscr{R}}^{
o}}(z) \! = \! \mathrm{I}, \quad z \! \in \! \widetilde{\Sigma}_{p}^{o},
\end{equation*}
with
\begin{equation*}
\mathcal{L}^{2}_{\mathrm{M}_{2}(\mathbb{C})}(\widetilde{\Sigma}_{p}^{o}) \!
\ni \! f \! \mapsto \! C^{0}_{w^{\Sigma_{\mathscr{R}}^{o}}}f \! := \! C_{-}^{
0}(fw^{\Sigma_{\mathscr{R}}^{o}}_{+}),
\end{equation*}
and
\begin{equation*}
\mathcal{L}^{2}_{\mathrm{M}_{2}(\mathbb{C})}(\widetilde{\Sigma}_{p}^{o}) \!
\ni \! f \! \mapsto \! (C_{\pm}^{0}f)(z) := \lim_{\underset{z^{\prime} \, \in
\, \pm \, \mathrm{side} \, \mathrm{of} \, \widetilde{\Sigma}_{p}^{o}}{z^{
\prime} \to z}} \int_{\widetilde{\Sigma}_{p}^{o}} \dfrac{z^{\prime}f(s)}{s(s
\! - \! z^{\prime})} \, \dfrac{\md s}{2 \pi \mi};
\end{equation*}
furthermore, $\norm{C_{\pm}^{0}}_{\mathscr{N}_{2}(\widetilde{\Sigma}_{p}^{o})}
\! < \! \infty$.
\begin{bbbbb}
Let $\mathscr{R}^{o} \colon \mathbb{C} \setminus \widetilde{\Sigma}_{p}^{o}
\! \to \! \operatorname{SL}_{2}(\mathbb{C})$ solve the following, equivalent
{\rm RHP:} {\rm (i)} $\mathscr{R}^{o}(z)$ is holomorphic for $z \! \in \!
\mathbb{C} \setminus \widetilde{\Sigma}_{p}^{o};$ {\rm (ii)} $\mathscr{R}^{
o}_{\pm}(z) \! := \! \lim_{\underset{z^{\prime} \, \in \, \pm \, \mathrm{side}
\, \mathrm{of} \, \widetilde{\Sigma}_{p}^{o}}{z^{\prime} \to z}} \! \mathscr{
R}^{o}(z^{\prime})$ satisfy the boundary condition
\begin{equation*}
\mathscr{R}_{+}^{o}(z) \! = \! \mathscr{R}_{-}^{o}(z) \upsilon_{\mathscr{R}}^{
o}(z), \quad z \! \in \! \widetilde{\Sigma}_{p}^{o},
\end{equation*}
where $\upsilon_{\mathscr{R}}^{o}(z)$, for $z \! \in \! \widetilde{\Sigma}_{
p}^{o}$, is defined in Lemma~{\rm 4.8} and satisfies the asymptotic (as $n
\! \to \! \infty)$ estimates given in Proposition~{\rm 5.1;} {\rm (iii)}
$\mathscr{R}^{o}(z) \! =_{\underset{z \in \mathbb{C} \setminus \widetilde{
\Sigma}_{p}^{o}}{z \to 0}} \! \mathrm{I} \! + \! \mathcal{O}(z);$ and
{\rm (iv)} $\mathscr{R}^{o}(z) \! =_{\underset{z \in \mathbb{C} \setminus
\widetilde{\Sigma}_{p}^{o}}{z \to \infty}} \! \mathcal{O}(1)$. Then:
\begin{compactenum}
\item[{\rm (1)}] for $z \! \in \! (-\infty,b_{0}^{o} \! - \! \delta_{b_{0}}^{
o}) \cup (a_{N+1}^{o} \! + \! \delta_{a_{N+1}}^{o},+\infty) \! =: \! \Sigma_{
p}^{o,1}$,
\begin{gather*}
\norm{w_{+}^{\Sigma_{\mathscr{R}}^{o}}(\cdot)}_{\mathcal{L}^{q}_{\mathrm{M}_{2}
(\mathbb{C})}(\Sigma_{p}^{o,1})} \! \underset{n \to \infty}{=} \! \mathcal{O}
\! \left(\dfrac{f(n) \me^{-(n+\frac{1}{2})c}}{(n \! + \! \frac{1}{2})^{1/q}}
\right), \quad q \! = \! 1,2, \\
\norm{w_{+}^{\Sigma_{\mathscr{R}}^{o}}(\cdot)}_{\mathcal{L}^{\infty}_{\mathrm{
M}_{2}(\mathbb{C})}(\Sigma_{p}^{o,1})} \! \underset{n \to \infty}{=} \!
\mathcal{O} \! \left(f(n) \me^{-(n+\frac{1}{2})c} \right),
\end{gather*}
where $c \! > \! 0$ and $f(n) \! =_{n \to \infty} \! \mathcal{O}(1);$
\item[{\rm (2)}] for $z \! \in \! (a_{j}^{o} \! + \! \delta_{a_{j}}^{o},b_{
j}^{o} \! - \! \delta_{b_{j}}^{o}) \! =: \! \Sigma_{p,j}^{o,2} \! \subset
\! \cup_{l=1}^{N} \Sigma_{p,l}^{o,2} \! =: \Sigma_{p}^{o,2}$, $j \! = \! 1,
\dotsc,N$,
\begin{gather*}
\norm{w_{+}^{\Sigma_{\mathscr{R}}^{o}}(\cdot)}_{\mathcal{L}^{q}_{\mathrm{M}_{2}
(\mathbb{C})}(\Sigma_{p,j}^{o,2})} \! \underset{n \to \infty}{=} \! \mathcal{
O} \! \left(\dfrac{f_{j}(n) \me^{-(n+\frac{1}{2})c_{j}}}{(n \! + \! \frac{1}{
2})^{1/q}}\right), \quad q \! = \! 1,2, \\
\norm{w_{+}^{\Sigma_{\mathscr{R}}^{o}}(\cdot)}_{\mathcal{L}^{\infty}_{\mathrm{
M}_{2}(\mathbb{C})}(\Sigma_{p,j}^{o,2})} \! \underset{n \to \infty}{=} \!
\mathcal{O} \! \left(f_{j}(n) \me^{-(n+\frac{1}{2})c_{j}}\right),
\end{gather*}
where $c_{j} \! > \! 0$ and $f_{j}(n) \! =_{n \to \infty} \! \mathcal{O}(1);$
\item[{\rm (3)}] for $z \! \in \! \cup_{j=1}^{N+1}(J_{j}^{o,\smallfrown}
\setminus (\mathbb{C}_{+} \cap (\mathbb{U}_{\delta_{b_{j-1}}}^{o} \cup
\mathbb{U}_{\delta_{a_{j}}}^{o}))) \! =: \! \Sigma_{p}^{o,3}$,
\begin{gather*}
\norm{w_{+}^{\Sigma_{\mathscr{R}}^{o}}(\cdot)}_{\mathcal{L}^{q}_{\mathrm{M}_{2}
(\mathbb{C})}(\Sigma_{p}^{o,3})} \! \underset{n \to \infty}{=} \! \mathcal{O}
\! \left(\dfrac{f(n) \me^{-(n+\frac{1}{2})c}}{(n \! + \! \frac{1}{2})^{1/q}}
\right), \quad q \! = \! 1,2, \\
\norm{w_{+}^{\Sigma_{\mathscr{R}}^{o}}(\cdot)}_{\mathcal{L}^{\infty}_{\mathrm{
M}_{2}(\mathbb{C})}(\Sigma_{p}^{o,3})} \! \underset{n \to \infty}{=} \!
\mathcal{O} \! \left(f(n) \me^{-(n+\frac{1}{2})c}\right),
\end{gather*}
where $c \! > \! 0$ and $f(n) \! =_{n \to \infty} \! \mathcal{O}(1);$
\item[{\rm (4)}] for $z \! \in \! \cup_{j=1}^{N+1}(J_{j}^{o,\smallsmile}
\setminus (\mathbb{C}_{-} \cap (\mathbb{U}_{\delta_{b_{j-1}}}^{o} \cup
\mathbb{U}_{\delta_{a_{j}}}^{o}))) \! =: \! \Sigma_{p}^{o,4}$,
\begin{gather*}
\norm{w_{+}^{\Sigma_{\mathscr{R}}^{o}}(\cdot)}_{\mathcal{L}^{q}_{\mathrm{M}_{2}
(\mathbb{C})}(\Sigma_{p}^{o,4})} \! \underset{n \to \infty}{=} \! \mathcal{O}
\! \left(\dfrac{f(n) \me^{-(n+\frac{1}{2})c}}{(n \! + \! \frac{1}{2})^{1/q}}
\right), \quad q \! = \! 1,2, \\
\norm{w_{+}^{\Sigma_{\mathscr{R}}^{o}}(\cdot)}_{\mathcal{L}^{\infty}_{\mathrm{
M}_{2}(\mathbb{C})}(\Sigma_{p}^{o,4})} \! \underset{n \to \infty}{=} \!
\mathcal{O} \! \left(f(n) \me^{-(n+\frac{1}{2})c}\right),
\end{gather*}
where $c \! > \! 0$ and $f(n) \! =_{n \to \infty} \! \mathcal{O}(1);$ and
\item[{\rm (5)}] for $z \! \in \! \cup_{j=1}^{N+1}(\partial \mathbb{U}_{
\delta_{b_{j-1}}}^{o} \cup \partial \mathbb{U}_{\delta_{a_{j}}}^{o}) \! =: \!
\Sigma_{p}^{o,5}$,
\begin{equation*}
\norm{w_{+}^{\Sigma_{\mathscr{R}}^{o}}(\cdot)}_{\mathcal{L}^{q}_{\mathrm{M}_{2}
(\mathbb{C})}(\Sigma_{p}^{o,5})} \! \underset{n \to \infty}{=} \! \mathcal{O}
\! \left((n \! + \! \tfrac{1}{2})^{-1}f(n) \right), \quad q \! \in \! \lbrace
1,2,\infty \rbrace,
\end{equation*}
where $f(n) \! =_{n \to \infty} \! \mathcal{O}(1)$.
\end{compactenum}
Furthermore,
\begin{equation*}
\norm{C^{0}_{w^{\Sigma_{\mathscr{R}}^{o}}}}_{\mathscr{N}_{r}(\widetilde{
\Sigma}_{p}^{o})} \! \underset{n \to \infty}{\leqslant} \! \mathcal{O} \!
\left((n \! + \! \tfrac{1}{2})^{-1+\frac{1}{r}}f(n) \right), \quad r \! \in
\! \lbrace 2,\infty \rbrace,
\end{equation*}
where $f(n) \! =_{n \to \infty} \! \mathcal{O}(1);$ in particular, $(\id \! -
\! C^{0}_{w^{\Sigma_{\mathscr{R}}^{o}}})^{-1} \! \! \upharpoonright_{\mathcal{
L}^{2}_{\mathrm{M}_{2}(\mathbb{C})}(\widetilde{\Sigma}_{p}^{o})}$ exists, that
is,
\begin{equation*}
\norm{(\id \! - \! C^{0}_{w^{\Sigma_{\mathscr{R}}^{o}}})^{-1}}_{\mathscr{N}_{2}
(\widetilde{\Sigma}_{p}^{o})} \underset{n \to \infty}{=} \mathcal{O}(1),
\end{equation*}
and it can be expanded in a Neumann series.
\end{bbbbb}

\emph{Proof.} Without loss of generality, assume that $0 \! \in \! \Sigma^{o,
1}_{p}$ (cf. Proposition~5.1). Recall that $w_{+}^{\Sigma_{\mathscr{R}}^{o}}
(z) \! = \! \upsilon_{\mathscr{R}}^{o}(z) \! - \! \mathrm{I}$, $z \! \in \!
\widetilde{\Sigma}_{p}^{o}$. For $z \! \in \! \Sigma_{p}^{o,1}$, using the
asymptotic (as $n \! \to \! \infty)$ estimate for $\upsilon_{\mathscr{R}}^{
o}(z)$ given in item~(1) of Proposition~5.1, one gets that
\begin{align*}
\norm{w_{+}^{\Sigma_{\mathscr{R}}^{o}}(\cdot)}_{\mathcal{L}^{\infty}_{\mathrm{
M}_{2}(\mathbb{C})}(\Sigma_{p}^{o,1})} &:= \max_{i,j =1,2} \, \, \sup_{z \in
\Sigma_{p}^{o,1}} \vert (w_{+}^{\Sigma_{\mathscr{R}}^{o}}(z))_{ij} \vert \!
\underset{n \to \infty}{=} \! \mathcal{O} \! \left(f(n) \me^{-(n+\frac{1}{2})
c} \right), \\
\norm{w_{+}^{\Sigma_{\mathscr{R}}^{o}}(\cdot)}_{\mathcal{L}^{1}_{\mathrm{M}_{
2}(\mathbb{C})}(\Sigma_{p}^{o,1})} &:= \int_{\Sigma_{p}^{o,1}} \vert w_{+}^{
\Sigma_{\mathscr{R}}^{o}}(z) \vert \, \vert \md z \vert \! = \! \int_{(\Sigma_{
p}^{o,1} \setminus \mathbb{U}_{0}^{o}) \cup \mathbb{U}_{0}^{o}} \vert w_{+}^{
\Sigma_{\mathscr{R}}^{o}}(z) \vert \, \vert \md z \vert \\
&= \, \left(\int_{\Sigma_{p}^{o,1} \setminus \mathbb{U}_{0}^{o}} \! + \!
\int_{\mathbb{U}_{0}^{o}} \right) \! \! \left(\sum_{i,j=1}^{2}(w_{+}^{\Sigma_{
\mathscr{R}}^{o}}(z))_{ij} \overline{(w_{+}^{\Sigma_{\mathscr{R}}^{o}}(z))_{
ij}} \right)^{1/2} \, \vert \md z \vert \\
& \, \underset{n \to \infty}{=} \! \mathcal{O} \! \left((n \! + \! \tfrac{1}{
2})^{-1}f(n) \me^{-(n+\frac{1}{2})c} \right) \! + \! \mathcal{O} \! \left((n
\! + \! \tfrac{1}{2})^{-1}f(n) \me^{-(n+\frac{1}{2})c} \right) \\
& \, \underset{n \to \infty}{=} \! \mathcal{O} \! \left((n \! + \! \tfrac{1}{
2})^{-1}f(n) \me^{-(n+\frac{1}{2})c} \right)
\end{align*}
$(\vert \md z \vert$ denotes arc length), and
\begin{align*}
\norm{w_{+}^{\Sigma_{\mathscr{R}}^{o}}(\cdot)}_{\mathcal{L}^{2}_{\mathrm{M}_{
2}(\mathbb{C})}(\Sigma_{p}^{o,1})} &:= \left(\int_{\Sigma_{p}^{o,1}} \vert
w_{+}^{\Sigma_{\mathscr{R}}^{o}}(z) \vert^{2} \, \vert \md z \vert \right)^{
1/2} \! = \! \left(\int_{(\Sigma_{p}^{o,1} \setminus \mathbb{U}_{0}^{o}) \cup
\mathbb{U}_{0}^{o}} \vert w_{+}^{\Sigma_{\mathscr{R}}^{o}}(z) \vert^{2} \,
\vert \md z \vert \right)^{1/2} \\
&= \, \left(\! \left(\int_{\Sigma_{p}^{o,1} \setminus \mathbb{U}_{0}^{o}} \!
+ \! \int_{\mathbb{U}_{0}^{o}} \right) \! \! \left(\sum_{i,j=1}^{2}(w_{+}^{
\Sigma_{\mathscr{R}}^{o}}(z))_{ij} \overline{(w_{+}^{\Sigma_{\mathscr{R}}^{o}}
(z))_{ij}} \right) \, \vert \md z \vert \right)^{1/2} \\
& \, \underset{n \to \infty}{=} \! \left(\mathcal{O} \! \left((n \! + \!
\tfrac{1}{2})^{-1}f(n) \me^{-(n+\frac{1}{2})c} \right) \! + \! \mathcal{O}
\! \left((n \! + \! \tfrac{1}{2})^{-1}f(n) \me^{-(n+\frac{1}{2})c} \right)
\right)^{1/2} \\
& \, \underset{n \to \infty}{=} \! \mathcal{O} \! \left((n \! + \! \tfrac{1}{
2})^{-1/2}f(n) \me^{-(n+\frac{1}{2})c} \right),
\end{align*}
where $c \! > \! 0$ and $f(n) \! =_{n \to \infty} \! \mathcal{O}(1)$.

For $z \! \in \! (a_{j}^{o} \! + \! \delta_{a_{j}}^{o},b_{j}^{o} \! - \!
\delta_{b_{j}}^{o}) \! =: \! \Sigma_{p,j}^{o,2} \! \subset \! \cup_{l=1}^{N}
\Sigma_{p,l}^{o,2} \! =: \! \Sigma_{p}^{o,2}$, $j \! = \! 1,\dotsc,N$, using
the asymptotic (as $n \! \to \! \infty)$ estimate for $\upsilon_{\mathscr{R}
}^{o}(z)$ given in item~(2) of Proposition~5.1, one gets that
\begin{align*}
\norm{w_{+}^{\Sigma_{\mathscr{R}}^{o}}(\cdot)}_{\mathcal{L}^{\infty}_{\mathrm{
M}_{2}(\mathbb{C})}(\Sigma_{p,j}^{o,2})} &:= \max_{k,m =1,2} \, \, \sup_{z
\in \Sigma_{p,j}^{o,2}} \vert (w_{+}^{\Sigma_{\mathscr{R}}^{o}}(z))_{km} \vert
\! \underset{n \to \infty}{=} \! \mathcal{O} \! \left(f_{j}(n) \me^{-(n+\frac{
1}{2})c_{j}} \right), \\
\norm{w_{+}^{\Sigma_{\mathscr{R}}^{o}}(\cdot)}_{\mathcal{L}^{1}_{\mathrm{M}_{
2}(\mathbb{C})}(\Sigma_{p,j}^{o,2})} &:= \int_{\Sigma_{p,j}^{o,2}} \vert w_{
+}^{\Sigma_{\mathscr{R}}^{o}}(z) \vert \, \vert \md z \vert \! = \! \int_{
\Sigma_{p,j}^{o,2}} \! \left(\sum_{i,j=1}^{2}(w_{+}^{\Sigma_{\mathscr{R}}^{o}}
(z))_{ij} \overline{(w_{+}^{\Sigma_{\mathscr{R}}^{o}}(z))_{ij}} \right)^{1/2}
\, \vert \md z \vert \\
& \, \underset{n \to \infty}{=} \! \mathcal{O} \! \left((n \! + \! \tfrac{1}{
2})^{-1}f_{j}(n) \me^{-(n+\frac{1}{2})c_{j}} \right),
\end{align*}
and
\begin{align*}
\norm{w_{+}^{\Sigma_{\mathscr{R}}^{o}}(\cdot)}_{\mathcal{L}^{2}_{\mathrm{M}_{
2}(\mathbb{C})}(\Sigma_{p,j}^{o,2})} &:= \left(\int_{\Sigma_{p,j}^{o,2}} \vert
w_{+}^{\Sigma_{\mathscr{R}}^{o}}(z) \vert^{2} \, \vert \md z \vert \right)^{
1/2} \! = \! \left(\int_{\Sigma_{p,j}^{o,2}} \sum_{k,l=1}^{2}(w_{+}^{\Sigma_{
\mathscr{R}}^{o}}(z))_{kl} \overline{(w_{+}^{\Sigma_{\mathscr{R}}^{o}}(z))_{k
l}} \, \vert \md z \vert \right)^{1/2} \\
& \, \underset{n \to \infty}{=} \! \mathcal{O} \! \left((n \! + \! \tfrac{1}{
2})^{-1/2}f_{j}(n) \me^{-(n+\frac{1}{2})c_{j}} \right),
\end{align*}
where $c_{j} \! > \! 0$ and $f_{j}(n) \! =_{n \to \infty} \! \mathcal{O}(1)$,
$j \! = \! 1,\dotsc,N$.

For $z \! \in \! \cup_{j=1}^{N+1}(J_{j}^{o,\smallfrown} \setminus (\mathbb{
C}_{+} \cap (\mathbb{U}_{\delta_{b_{j-1}}}^{o} \cup \mathbb{U}_{\delta_{a_{
j}}}^{o}))) \! =: \! \Sigma_{p}^{o,3}$ $(\subset \! \widetilde{\Sigma}_{p}^{
o})$, using the asymptotic (as $n \! \to \! \infty)$ estimate for $\upsilon_{
\mathscr{R}}^{o}(z)$ given in item~(3) of Proposition~5.1, one gets that
\begin{align*}
\norm{w_{+}^{\Sigma_{\mathscr{R}}^{o}}(\cdot)}_{\mathcal{L}^{\infty}_{\mathrm{
M}_{2}(\mathbb{C})}(\Sigma_{p}^{o,3})} &:= \max_{i,j =1,2} \, \, \sup_{z \in
\Sigma_{p}^{o,3}} \vert (w_{+}^{\Sigma_{\mathscr{R}}^{o}}(z))_{ij} \vert \!
\underset{n \to \infty}{=} \! \mathcal{O} \! \left(f(n) \me^{-(n+\frac{1}{2})
c} \right), \\
\norm{w_{+}^{\Sigma_{\mathscr{R}}^{o}}(\cdot)}_{\mathcal{L}^{1}_{\mathrm{M}_{
2}(\mathbb{C})}(\Sigma_{p}^{o,3})} &:= \int_{\Sigma_{p}^{o,3}} \vert w_{+}^{
\Sigma_{\mathscr{R}}^{o}}(z) \vert \, \vert \md z \vert \! = \! \int_{\Sigma_{
p}^{o,3}} \! \left(\sum_{i,j=1}^{2}(w_{+}^{\Sigma_{\mathscr{R}}^{o}}(z))_{ij}
\overline{(w_{+}^{\Sigma_{\mathscr{R}}^{o}}(z))_{ij}} \right)^{1/2} \, \vert
\md z \vert \\
& \, \underset{n \to \infty}{=} \! \mathcal{O} \! \left((n \! + \! \tfrac{1}{
2})^{-1}f(n) \me^{-(n+\frac{1}{2})c} \right),
\end{align*}
and
\begin{align*}
\norm{w_{+}^{\Sigma_{\mathscr{R}}^{o}}(\cdot)}_{\mathcal{L}^{2}_{\mathrm{M}_{
2}(\mathbb{C})}(\Sigma_{p}^{o,3})} &:= \left(\int_{\Sigma_{p}^{o,3}} \vert w_{
+}^{\Sigma_{\mathscr{R}}^{o}}(z) \vert^{2} \, \vert \md z \vert \right)^{1/2}
\! = \! \left(\int_{\Sigma_{p}^{o,3}} \sum_{i,j=1}^{2}(w_{+}^{\Sigma_{
\mathscr{R}}^{o}}(z))_{ij} \overline{(w_{+}^{\Sigma_{\mathscr{R}}^{o}}(z))_{
ij}} \, \vert \md z \vert \right)^{1/2} \\
& \, \underset{n \to \infty}{=} \! \mathcal{O} \! \left((n \! + \! \tfrac{1}{
2})^{-1/2}f(n) \me^{-(n+\frac{1}{2})c} \right),
\end{align*}
where $c \! > \! 0$ and $f(n) \! =_{n \to \infty} \! \mathcal{O}(1)$: the
above analysis applies, \emph{mutatis mutandis}, for the analogous estimates
on $\Sigma_{p}^{o,4} \! := \! \cup_{j=1}^{N+1}(J_{j}^{o,\smallsmile} \setminus
(\mathbb{C}_{-} \cap (\mathbb{U}_{\delta_{b_{j-1}}}^{o} \cup \mathbb{U}_{
\delta_{a_{j}}}^{o}))) \! \ni \! z$.

For $z \! \in \! \cup_{j=1}^{N+1}(\partial \mathbb{U}_{\delta_{b_{j-1}}}^{o}
\cup \partial \mathbb{U}_{\delta_{a_{j}}}^{o}) \! =: \! \Sigma_{p}^{o,5}$
$(\subset \! \widetilde{\Sigma}_{p}^{o})$, using the $(2(N \! + \! 1))$
asymptotic (as $n \! \to \! \infty)$ estimates for $\upsilon_{\mathscr{
R}}^{o}(z)$ given in item~(5) of Proposition~5.1, one gets that
\begin{align*}
\norm{w_{+}^{\Sigma_{\mathscr{R}}^{o}}(\cdot)}_{\mathcal{L}^{\infty}_{\mathrm{
M}_{2}(\mathbb{C})}(\Sigma_{p}^{o,5})} &:= \max_{i,j =1,2} \, \, \sup_{z \in
\Sigma_{p}^{o,5}} \vert (w_{+}^{\Sigma_{\mathscr{R}}^{o}}(z))_{ij} \vert \!
\underset{n \to \infty}{=} \! \mathcal{O} \! \left((n \! + \! \tfrac{1}{2})^{
-1}f(n) \right), \\
\norm{w_{+}^{\Sigma_{\mathscr{R}}^{o}}(\cdot)}_{\mathcal{L}^{1}_{\mathrm{M}_{
2}(\mathbb{C})}(\Sigma_{p}^{o,5})} &:= \int_{\Sigma_{p}^{o,5}} \vert w_{+}^{
\Sigma_{\mathscr{R}}^{o}}(z) \vert \, \vert \md z \vert \! = \! \int_{\cup_{j=
1}^{N+1}(\partial \mathbb{U}_{\delta_{b_{j-1}}}^{o} \cup \partial \mathbb{U}_{
\delta_{a_{j}}}^{o})} \vert w_{+}^{\Sigma_{\mathscr{R}}^{o}}(z) \vert \, \vert
\md z \vert \\
&= \, \sum_{k=1}^{N+1} \! \left(\int_{\partial \mathbb{U}_{\delta_{b_{k-1}}}^{
o}} \! + \! \int_{\partial \mathbb{U}_{\delta_{a_{k}}}^{o}} \right) \! \left(
\sum_{i,j=1}^{2}(w_{+}^{\Sigma_{\mathscr{R}}^{o}}(z))_{ij} \overline{(w_{+}^{
\Sigma_{\mathscr{R}}^{o}}(z))_{ij}}\right)^{1/2} \, \vert \md z \vert,
\end{align*}
whence, (cf. Lemma~4.5) using the fact that the respective factors $(\gamma^{o}
(0))^{-1} \gamma^{o}(z) \! \pm \! \gamma^{o}(0)(\gamma^{o}(z))^{-1}$ and
$\boldsymbol{\theta}^{o}(\pm \boldsymbol{u}^{o}(z) \! - \! \tfrac{1}{2 \pi}(n
\! + \! \tfrac{1}{2}) \boldsymbol{\Omega}^{o} \! \pm \! \boldsymbol{d}_{o})$
are uniformly bounded (with respect to $z)$ in compact intervals outside open
intervals surrounding the end-points of the support of the `odd' equilibrium
measure, one arrives at
\begin{align*}
\norm{w_{+}^{\Sigma_{\mathscr{R}}^{o}}(\cdot)}_{\mathcal{L}^{1}_{\mathrm{M}_{
2}(\mathbb{C})}(\Sigma_{p}^{o,5})} &\underset{n \to \infty}{=} \dfrac{1}{(n \!
+ \! \frac{1}{2})} \sum_{k=1}^{N+1} \! \left(\int_{\partial \mathbb{U}_{
\delta_{b_{k-1}}}^{o}} \dfrac{\vert \star_{b_{k-1}}^{o}(z;n) \vert}{(z \! - \!
b_{k-1}^{o})^{3/2}} \, \vert \md z \vert \! + \! \int_{\partial \mathbb{U}_{
\delta_{a_{k}}}^{o}} \dfrac{\vert \star_{a_{k}}^{o}(z;n) \vert}{(z \! - \!
a_{k}^{o})^{3/2}} \, \vert \md z \vert \right) \\
&\underset{n \to \infty}{=} \, \mathcal{O} \! \left((n \! + \! \tfrac{1}{2})^{
-1}f(n) \right),
\end{align*}
and, similarly,
\begin{align*}
\norm{w_{+}^{\Sigma_{\mathscr{R}}^{o}}(\cdot)}_{\mathcal{L}^{2}_{\mathrm{M}_{
2}(\mathbb{C})}(\Sigma_{p}^{o,5})} &:= \left(\int_{\Sigma_{p}^{o,5}} \vert w_{
+}^{\Sigma_{\mathscr{R}}^{o}}(z) \vert^{2} \, \vert \md z \vert \right)^{1/2}
\! = \! \left(\int_{\cup_{j=1}^{N+1}(\partial \mathbb{U}_{\delta_{b_{j-1}}}^{
o} \cup \partial \mathbb{U}_{\delta_{a_{j}}}^{o})} \vert w_{+}^{\Sigma_{
\mathscr{R}}^{o}}(z) \vert^{2} \, \vert \md z \vert \right)^{1/2} \\
&= \, \left(\sum_{k=1}^{N+1} \! \left(\int_{\partial \mathbb{U}_{\delta_{b_{k
-1}}}^{o}} \! + \! \int_{\partial \mathbb{U}_{\delta_{a_{k}}}^{o}} \right) \!
\sum_{i,j=1}^{2}(w_{+}^{\Sigma_{\mathscr{R}}^{o}}(z))_{ij} \overline{(w_{+}^{
\Sigma_{\mathscr{R}}^{o}}(z))_{ij}} \, \vert \md z \vert \right)^{1/2} \\
&\underset{n \to \infty}{=} \left(\dfrac{1}{(n \! + \! \tfrac{1}{2})^{2}}
\sum_{k=1}^{N+1} \! \left(\int_{\partial \mathbb{U}_{\delta_{b_{k-1}}}^{o}} \!
\dfrac{\vert \star_{b_{k-1}}^{o}(z;n) \vert}{(z \! - \! b_{k-1}^{o})^{3}} \,
\vert \md z \vert \! + \! \int_{\partial \mathbb{U}_{\delta_{a_{k}}}^{o}} \!
\dfrac{\vert \star_{a_{k}}^{o}(z;n) \vert}{(z \! - \! a_{k}^{o})^{3}} \, \vert
\md z \vert \right) \right)^{1/2} \\
&\underset{n \to \infty}{=} \, \mathcal{O} \! \left((n \! + \! \tfrac{1}{2})^{
-1}f(n) \right),
\end{align*}
where $f(n) \! =_{n \to \infty} \! \mathcal{O}(1)$.

Recall that $C^{0}_{w^{\Sigma_{\mathscr{R}}^{o}}}f \! := \! C^{0}_{-}(fw^{
\Sigma_{\mathscr{R}}^{o}}_{+})$, where $(C^{0}_{-}g)(z) \! := \! \lim_{
\underset{z^{\prime} \in -\widetilde{\Sigma}_{p}^{o}}{z^{\prime} \to z}} \!
\int_{\widetilde{\Sigma}_{p}^{o}} \tfrac{z^{\prime}g(s)}{s(s-z^{\prime})} \,
\tfrac{\md s}{2 \pi \mi}$, with $-\widetilde{\Sigma}_{p}^{o}$ shorthand for
`the $-$ side of $\widetilde{\Sigma}_{p}^{o}$'. For the $\norm{C^{0}_{w^{
\Sigma_{\mathscr{R}}^{o}}}}_{\mathscr{N}_{\infty}(\widetilde{\Sigma}_{p}^{
o})}$ norm, one proceeds as follows:
\begin{align*}
\norm{C^{0}_{w^{\Sigma_{\mathscr{R}}^{o}}}g}_{\mathcal{L}^{\infty}_{\mathrm{
M}_{2}(\mathbb{C})}(\widetilde{\Sigma}_{p}^{o})} &:= \max_{j,l=1,2} \, \,
\sup_{z \in \widetilde{\Sigma}_{p}^{o}} \vert (C^{0}_{w^{\Sigma_{\mathscr{R}
}^{o}}}g)_{jl}(z) \vert \! = \! \max_{j,l=1,2} \, \, \sup_{z \in \widetilde{
\Sigma}_{p}^{o}} \! \left\vert \lim_{\underset{z^{\prime} \in -\widetilde{
\Sigma}_{p}^{o}}{z^{\prime} \to z}} \int_{\widetilde{\Sigma}_{p}^{o}} \dfrac{
z^{\prime}(g(s)w^{\Sigma_{\mathscr{R}}^{o}}_{+}(s))_{jl}}{s(s \! - \!
z^{\prime})} \, \dfrac{\md s}{2 \pi \mi} \right\vert \\
&\leqslant \, \norm{g(\cdot)}_{\mathcal{L}^{\infty}_{\mathrm{M}_{2}(\mathbb{
C})}(\widetilde{\Sigma}_{p}^{o})} \max_{j,l=1,2} \, \, \sup_{z \in \widetilde{
\Sigma}_{p}^{o}} \! \left\vert \lim_{\underset{z^{\prime} \in -\widetilde{
\Sigma}_{p}^{o}}{z^{\prime} \to z}} \int_{\widetilde{\Sigma}_{p}^{o}} \dfrac{
z^{\prime}(w^{\Sigma_{\mathscr{R}}^{o}}_{+}(s))_{jl}}{s(s \! - \! z^{\prime})
} \, \dfrac{\md s}{2 \pi \mi} \right\vert \\
&\leqslant \, \norm{g(\cdot)}_{\mathcal{L}^{\infty}_{\mathrm{M}_{2}(\mathbb{
C})}(\widetilde{\Sigma}_{p}^{o})} \max_{j,l=1,2} \, \, \sup_{z \in \widetilde{
\Sigma}_{p}^{o}} \! \left\vert \lim_{\underset{z^{\prime} \in -\widetilde{
\Sigma}_{p}^{o}}{z^{\prime} \to z}} \! \left(\int_{\Sigma_{p}^{o,1} \setminus
\mathbb{U}_{0}^{o}} \! + \! \int_{\mathbb{U}_{0}^{o}} \! + \! \sum_{k=1}^{N+
1} \int_{\Sigma_{p,k}^{o,2}} \! + \! \int_{\Sigma_{p}^{o,3}} \right. \right. \\
&+ \left. \left. \, \int_{\Sigma_{p}^{o,4}} \! + \! \sum_{k=1}^{N+1} \! \left(
\int_{\partial \mathbb{U}_{\delta_{b_{k-1}}}^{o}} \! + \! \int_{\partial
\mathbb{U}_{a_{k}}^{o}} \right) \right) \! \dfrac{z^{\prime}(w^{\Sigma_{
\mathscr{R}}^{o}}_{+}(s))_{jl}}{s(s \! - \! z^{\prime})} \, \dfrac{\md s}{2
\pi \mi} \right\vert \\
&\underset{n \to \infty}{\leqslant} \norm{g(\cdot)}_{\mathcal{L}^{\infty}_{
\mathrm{M}_{2}(\mathbb{C})}(\widetilde{\Sigma}_{p}^{o})} \max_{j,l=1,2} \, \,
\sup_{z \in \widetilde{\Sigma}_{p}^{o}} \! \left\vert \lim_{\underset{z^{
\prime} \in -\widetilde{\Sigma}_{p}^{o}}{z^{\prime} \to z}} \! \left(\int_{
\Sigma_{p}^{o,1} \setminus \mathbb{U}_{0}^{o}} \dfrac{(\mathcal{O}(\me^{-(n+
\frac{1}{2})c_{\infty} \vert s \vert}f_{\infty}(n)))_{jl}z^{\prime}}{s(s \!
- \! z^{\prime})} \, \dfrac{\md s}{2 \pi \mi} \right. \right. \\
&+ \left. \left. \, \int_{\mathbb{U}_{0}^{o}} \dfrac{(\mathcal{O}(\me^{-(n+
\frac{1}{2})c_{0} \vert s \vert^{-1}}f_{0}(n)))_{jl}z^{\prime}}{s(s \! - \!
z^{\prime})} \, \dfrac{\md s}{2 \pi \mi} \! + \! \sum_{k=1}^{N} \int_{\Sigma_{
p,k}^{o,2}} \dfrac{(\mathcal{O}(\me^{-(n+\frac{1}{2})c_{k}(s-a_{k}^{o})}f_{k}
(n)))_{jl}z^{\prime}}{s(s \! - \! z^{\prime})} \right. \right. \\
&\times \left. \left. \, \dfrac{\md s}{2 \pi \mi} \! + \! \int_{\Sigma_{p}^{o,
3}} \dfrac{(\mathcal{O}(\me^{-(n+\frac{1}{2}) \overset{\smallfrown}{c} \vert s
\vert} \overset{\smallfrown}{f}(n)))_{jl}z^{\prime}}{s(s \! - \! z^{\prime})}
\, \dfrac{\md s}{2 \pi \mi} \! + \! \int_{\Sigma_{p}^{o,4}} \dfrac{(\mathcal{
O}(\me^{-(n+\frac{1}{2}) \overset{\smallsmile}{c} \vert s \vert} \overset{
\smallsmile}{f}(n)))_{jl}z^{\prime}}{s(s \! - \! z^{\prime})} \, \dfrac{\md
s}{2 \pi \mi} \right. \right. \\
&+ \left. \left. \, \sum_{k=1}^{N+1} \! \left(\int_{\partial \mathbb{U}_{
\delta_{b_{k-1}}}^{o}} \! \left(\! \mathcal{O} \! \left(\dfrac{z^{\prime}
\overset{o}{\mathfrak{M}}^{\raise-1.0ex\hbox{$\scriptstyle \infty$}}(s)
\left(
\begin{smallmatrix}
\ast & \ast \\
\ast & \ast
\end{smallmatrix}
\right)
(\overset{o}{\mathfrak{M}}^{\raise-1.0ex\hbox{$\scriptstyle \infty$}}(s))^{-1}
}{(n \! + \! \frac{1}{2})s(s \! - \! z^{\prime})(s \! - \! b_{k-1}^{o})^{3/2}
G_{b_{k-1}}^{o}(s)} \right) \right)_{jl} \dfrac{\md s}{2 \pi \mi} \right.
\right. \right. \\
&+ \left. \left. \left. \, \int_{\partial \mathbb{U}_{\delta_{a_{k}}}^{o}} \!
\left(\! \mathcal{O} \! \left(\dfrac{z^{\prime}
\overset{o}{\mathfrak{M}}^{\raise-1.0ex\hbox{$\scriptstyle \infty$}}(s)
\left(
\begin{smallmatrix}
\ast & \ast \\
\ast & \ast
\end{smallmatrix}
\right)
(\overset{o}{\mathfrak{M}}^{\raise-1.0ex\hbox{$\scriptstyle \infty$}}(s))^{-1}
}{(n \! + \! \frac{1}{2})s(s \! - \! z^{\prime})(s \! - \! a_{k}^{o})^{3/2}G_{
a_{k}}^{o}(s)} \right) \right)_{jl} \dfrac{\md s}{2 \pi \mi} \right) \right)
\right\vert,
\end{align*}
whence, taking note of the partial fraction decomposition $\tfrac{z^{\prime}}{
s(s-z^{\prime})} \! = \! -\tfrac{1}{s}+\tfrac{1}{s-z^{\prime}}$, and (cf.
Lemma~4.5) using the fact that the respective factors $(\gamma^{o}(0))^{-1}
\gamma^{o}(z) \! \pm \! \gamma^{o}(0)(\gamma^{o}(z))^{-1}$ and $\boldsymbol{
\theta}^{o}(\pm \boldsymbol{u}^{o}(z) \! - \! \tfrac{1}{2 \pi}(n \! + \!
\tfrac{1}{2}) \boldsymbol{\Omega}^{o} \! \pm \! \boldsymbol{d}_{o})$ are
uniformly bounded (with respect to $z)$ in compact intervals outside open
intervals surrounding the end-points of the support of the `odd' equilibrium
measure, one arrives at, after a straightforward integration argument and an
application of the Maximum Length $(\mathrm{ML})$ Theorem,
\begin{align*}
\norm{C^{0}_{w^{\Sigma_{\mathscr{R}}^{o}}}g}_{\mathcal{L}^{\infty}_{\mathrm{
M}_{2}(\mathbb{C})}(\widetilde{\Sigma}_{p}^{o})} &\underset{n \to \infty}{
\leqslant} \norm{g(\cdot)}_{\mathcal{L}^{\infty}_{\mathrm{M}_{2}(\mathbb{C})}
(\widetilde{\Sigma}_{p}^{o})} \! \left(\! \mathcal{O} \! \left(\dfrac{f(n)
\me^{-(n+\frac{1}{2})c}}{(n \! + \! \frac{1}{2}) \min \{1,\operatorname{dist}
(z,\widetilde{\Sigma}_{p}^{o})\}} \right) \right. \\
+&\left. \, \mathcal{O} \! \left(\dfrac{f(n)}{(n \! + \! \frac{1}{2}) \min \{1,
\operatorname{dist}(z,\widetilde{\Sigma}_{p}^{o})\}} \right) \right) \!
\underset{n \to \infty}{=} \! \norm{g(\cdot)}_{\mathcal{L}^{\infty}_{\mathrm{
M}_{2}(\mathbb{C})}(\widetilde{\Sigma}_{p}^{o})} \mathcal{O} \! \left(\dfrac{
f(n)}{n \! + \! \frac{1}{2}} \right),
\end{align*}
where $\operatorname{dist}(z,\widetilde{\Sigma}_{p}^{o}) \! := \! \inf
\left\lbrace \mathstrut \vert z \! - \! r \vert; \, r \! \in \! \widetilde{
\Sigma}_{p}^{o}, \, z \! \in \! \mathbb{C} \setminus \widetilde{\Sigma}_{p}^{
o} \right\rbrace$ $(> \! 0)$, and $f(n) \! =_{n \to \infty} \! \mathcal{O}
(1)$, whence one obtains the asymptotic (as $n \to \infty)$ estimate for
$\norm{C^{0}_{w^{\Sigma_{\mathscr{R}}^{o}}}}_{\mathscr{N}_{\infty}(\widetilde{
\Sigma}_{p}^{o})}$ stated in the Proposition. Similarly, for $\norm{C^{0}_{w^{
\Sigma_{\mathscr{R}}^{o}}}}_{\mathscr{N}_{2}(\widetilde{\Sigma}_{p}^{o})}$, as
$n \! \to \! \infty$:
\begin{align*}
\norm{C^{0}_{w^{\Sigma_{\mathscr{R}}^{o}}}g}_{\mathcal{L}^{2}_{\mathrm{M}_{2}
(\mathbb{C})}(\widetilde{\Sigma}_{p}^{o})} &:= \left(\int_{\widetilde{\Sigma}_{
p}^{o}} \vert (C^{0}_{w^{\Sigma_{\mathscr{R}}^{o}}}g)(z) \vert^{2} \, \vert
\md z \vert \right)^{1/2} \! = \! \left(\int_{\widetilde{\Sigma}_{p}^{o}}
\sum_{j,l=1}^{2}(C^{0}_{w^{\Sigma_{\mathscr{R}}^{o}}}g)_{jl}(z) \overline{(
C^{0}_{w^{\Sigma_{\mathscr{R}}^{o}}}g)_{jl}(z)} \, \vert \md z \vert \right)^{
1/2} \\
&= \, \left(\int_{\widetilde{\Sigma}_{p}^{o}} \sum_{j,l=1}^{2} \! \left\vert
\lim_{\underset{z^{\prime} \in -\widetilde{\Sigma}_{p}^{o}}{z^{\prime} \to z}}
\int_{\widetilde{\Sigma}_{p}^{o}} \dfrac{z^{\prime}(g(s)w^{\Sigma_{\mathscr{
R}}^{o}}_{+}(s))_{jl}}{s(s \! - \! z^{\prime})} \, \dfrac{\md s}{2 \pi \mi}
\right\vert^{2} \vert \md z \vert \right)^{1/2} \\
&\leqslant \, \norm{g(\cdot)}_{\mathcal{L}^{2}_{\mathrm{M}_{2}(\mathbb{C})}
(\widetilde{\Sigma}_{p}^{o})} \! \left(\int_{\widetilde{\Sigma}_{p}^{o}}
\sum_{j,l=1}^{2} \! \left\vert \lim_{\underset{z^{\prime} \in -\widetilde{
\Sigma}_{p}^{o}}{z^{\prime} \to z}} \int_{\widetilde{\Sigma}_{p}^{o}} \dfrac{
z^{\prime}(w^{\Sigma_{\mathscr{R}}^{e}}_{+}(s))_{jl}}{s(s \! - \! z^{\prime})}
\, \dfrac{\md s}{2 \pi \mi} \right\vert^{2} \vert \md z \vert \right)^{1/2} \\
&\leqslant \, \norm{g(\cdot)}_{\mathcal{L}^{2}_{\mathrm{M}_{2}(\mathbb{C})}
(\widetilde{\Sigma}_{p}^{o})} \! \left(\int_{\widetilde{\Sigma}_{p}^{o}} \sum_{
j,l=1}^{2} \! \left\vert \lim_{\underset{z^{\prime} \in -\widetilde{\Sigma}_{
p}^{o}}{z^{\prime} \to z}} \! \left(\int_{\Sigma_{p}^{o,1} \setminus \mathbb{
U}_{0}^{o}} \! + \! \int_{\mathbb{U}_{0}^{o}} \! + \! \sum_{k=1}^{N+1} \int_{
\Sigma_{p,k}^{o,2}} \! + \! \int_{\Sigma_{p}^{o,3}} \right. \right. \right. \\
&+ \left. \left. \left. \, \int_{\Sigma_{p}^{o,4}} \! + \! \sum_{k=1}^{N+1} \!
\left(\int_{\partial \mathbb{U}_{\delta_{b_{k-1}}}^{o}} \! + \! \int_{\partial
\mathbb{U}_{a_{k}}^{o}} \right) \right) \! \dfrac{z^{\prime}(w^{\Sigma_{
\mathscr{R}}^{o}}_{+}(s))_{jl}}{s(s \! - \! z^{\prime})} \, \dfrac{\md s}{2
\pi \mi} \right\vert^{2} \vert \md z \vert \right)^{1/2} \\
&\underset{n \to \infty}{\leqslant} \norm{g(\cdot)}_{\mathcal{L}^{2}_{\mathrm{
M}_{2}(\mathbb{C})}(\widetilde{\Sigma}_{p}^{o})} \! \left(\int_{\widetilde{
\Sigma}_{p}^{o}} \sum_{j,l=1}^{2} \! \left\vert \lim_{\underset{z^{\prime} \in
-\widetilde{\Sigma}_{p}^{o}}{z^{\prime} \to z}} \! \left(\int_{\Sigma_{p}^{o,
1} \setminus \mathbb{U}_{0}^{o}} \dfrac{(\mathcal{O}(\me^{-(n+\frac{1}{2})c_{
\infty} \vert s \vert}f_{\infty}(n)))_{jl}z^{\prime}}{s(s \! - \! z^{\prime})}
\, \dfrac{\md s}{2 \pi \mi} \right. \right. \right. \\
&+ \left. \left. \left. \int_{\mathbb{U}_{0}^{o}} \dfrac{(\mathcal{O}(\me^{
-(n+\frac{1}{2})c_{0} \vert s \vert^{-1}}f_{0}(n)))_{jl}z^{\prime}}{s(s \! -
\! z^{\prime})} \, \dfrac{\md s}{2 \pi \mi} \! + \! \sum_{k=1}^{N} \int_{
\Sigma_{p,k}^{o,2}} \dfrac{(\mathcal{O}(\me^{-(n+\frac{1}{2})c_{k}(s-a_{k}^{o}
)}f_{k}(n)))_{jl}z^{\prime}}{s(s \! - \! z^{\prime})} \right. \right. \right.
\\
&\times \left. \left. \left. \, \dfrac{\md s}{2 \pi \mi} \! + \! \int_{\Sigma_{
p}^{o,3}} \dfrac{(\mathcal{O}(\me^{-(n+\frac{1}{2}) \overset{\smallfrown}{c}
\vert s \vert} \overset{\smallfrown}{f}(n)))_{jl}z^{\prime}}{s(s \! - \! z^{
\prime})} \, \dfrac{\md s}{2 \pi \mi} \! + \! \int_{\Sigma_{p}^{o,4}} \dfrac{
(\mathcal{O}(\me^{-(n+\frac{1}{2}) \overset{\smallsmile}{c} \vert s \vert}
\overset{\smallsmile}{f}(n)))_{jl}z^{\prime}}{s(s \! - \! z^{\prime})} \,
\dfrac{\md s}{2 \pi \mi} \right. \right. \right. \\
&+ \left. \left. \left. \, \sum_{k=1}^{N+1} \! \left(\int_{\partial \mathbb{
U}_{\delta_{b_{k-1}}}^{o}} \! \left(\! \mathcal{O} \! \left(\dfrac{z^{\prime}
\overset{o}{\mathfrak{M}}^{\raise-1.0ex\hbox{$\scriptstyle \infty$}}(s)
\left(
\begin{smallmatrix}
\ast & \ast \\
\ast & \ast
\end{smallmatrix}
\right)
(\overset{o}{\mathfrak{M}}^{\raise-1.0ex\hbox{$\scriptstyle \infty$}}(s))^{-1}
}{(n \! + \! \frac{1}{2})s(s \! - \! z^{\prime})(s \! - \! b_{k-1}^{o})^{3/2}
G_{b_{k-1}}^{o}(s)} \right) \right)_{jl} \dfrac{\md s}{2 \pi \mi} \right.
\right. \right. \right. \\
&+ \left. \left. \left. \left. \, \int_{\partial \mathbb{U}_{\delta_{a_{k}}}^{
o}} \! \left(\! \mathcal{O} \! \left(\dfrac{z^{\prime}
\overset{o}{\mathfrak{M}}^{\raise-1.0ex\hbox{$\scriptstyle \infty$}}(s)
\left(
\begin{smallmatrix}
\ast & \ast \\
\ast & \ast
\end{smallmatrix}
\right)
(\overset{o}{\mathfrak{M}}^{\raise-1.0ex\hbox{$\scriptstyle \infty$}}(s))^{-1}
}{(n \! + \! \frac{1}{2})s(s \! - \! z^{\prime})(s \! - \! a_{k}^{o})^{3/2}G_{
a_{k}}^{o}(s)} \right) \right)_{jl} \dfrac{\md s}{2 \pi \mi} \right) \right)
\right\vert^{2} \vert \md z \vert \right)^{1/2},
\end{align*}
whence, taking note of the partial fraction decomposition $\tfrac{z^{\prime}}{
s(s-z^{\prime})} \! = \! -\tfrac{1}{s}+\tfrac{1}{s-z^{\prime}}$, and (cf.
Lemma~4.5) using the fact that the respective factors $(\gamma^{o}(0))^{-1}
\gamma^{o}(z) \! \pm \! \gamma^{o}(0)(\gamma^{o}(z))^{-1}$ and $\boldsymbol{
\theta}^{o}(\pm \boldsymbol{u}^{o}(z) \! - \! \tfrac{1}{2 \pi}(n \! + \!
\tfrac{1}{2}) \boldsymbol{\Omega}^{o} \! \pm \! \boldsymbol{d}_{o})$ are
uniformly bounded (with respect to $z)$ in compact intervals outside open
intervals surrounding the end-points of the support of the `odd' equilibrium
measure, one arrives at, after a straightforward integration argument and an
application of the $\mathrm{ML}$ Theorem,
\begin{align*}
\norm{C^{0}_{w^{\Sigma_{\mathscr{R}}^{o}}}g}_{\mathcal{L}^{2}_{\mathrm{M}_{2}
(\mathbb{C})}(\widetilde{\Sigma}_{p}^{o})} &\underset{n \to \infty}{\leqslant}
\norm{g(\cdot)}_{\mathcal{L}^{2}_{\mathrm{M}_{2}(\mathbb{C})}(\widetilde{
\Sigma}_{p}^{o})} \! \left(\! \mathcal{O} \! \left(\dfrac{f(n) \me^{-(n+\frac{
1}{2})c}}{(n \! + \! \frac{1}{2})^{1/2} \min \{1,\operatorname{dist}(z,
\widetilde{\Sigma}_{p}^{o})\}} \right) \right. \\
+&\left. \, \mathcal{O} \! \left(\dfrac{f(n)}{(n \! + \! \frac{1}{2})^{1/2}
\min \{1,\operatorname{dist}(z,\widetilde{\Sigma}_{p}^{o})\}} \right) \right)
\! \underset{n \to \infty}{=} \! \norm{g(\cdot)}_{\mathcal{L}^{2}_{\mathrm{M}_{
2}(\mathbb{C})}(\widetilde{\Sigma}_{p}^{o})} \mathcal{O} \! \left(\dfrac{f(n)}{
\sqrt{\smash[b]{n \! + \! \frac{1}{2}}}} \right),
\end{align*}
where $f(n) \! =_{n \to \infty} \! \mathcal{O}(1)$, whence one obtains the
asymptotic (as $n \to \infty)$ inequality for $\norm{C^{0}_{w^{\Sigma_{
\mathscr{R}}^{o}}}}_{\mathscr{N}_{2}(\widetilde{\Sigma}_{p}^{o})}$ stated in
the Proposition. The above analysis establishes the fact that, as $n \! \to \!
\infty$, $C^{0}_{w^{\Sigma_{\mathscr{R}}^{o}}} \! \in \! \mathscr{N}_{2}
(\widetilde{\Sigma}_{p}^{o})$, with operator norm $\norm{C^{0}_{w^{\Sigma_{
\mathscr{R}}^{o}}}}_{\mathscr{N}_{2}(\widetilde{\Sigma}_{p}^{o})} \! =_{n \to
\infty} \! \mathcal{O}((n \! + \! 1/2)^{-1/2}f(n))$, where $f(n) \! =_{n \to
\infty} \! \mathcal{O}(1)$; due to a well-known result for bounded linear
operators in Hilbert space \cite{a92}, it follows, thus, that $(\id \! - \!
C^{0}_{w^{\Sigma_{\mathscr{R}}^{o}}})^{-1} \! \! \upharpoonright_{\mathcal{
L}^{2}_{\mathrm{M}_{2}(\mathbb{C})}(\widetilde{\Sigma}_{p}^{o})}$ exists, and
$(\id \! - \! C^{0}_{w^{\Sigma_{\mathscr{R}}^{o}}}) \! \! \upharpoonright_{
\mathcal{L}^{2}_{\mathrm{M}_{2}(\mathbb{C})}(\widetilde{\Sigma}_{p}^{o})}$ can
be inverted by a Neumann series (as $n \! \to \! \infty)$, with $\norm{(\id \!
- \! C^{0}_{w^{\Sigma_{\mathscr{R}}^{o}}})^{-1}}_{\mathscr{N}_{2}(\widetilde{
\Sigma}_{p}^{o})} \! \leqslant_{n \to \infty} \! (1 \! - \! \norm{C^{0}_{w^{
\Sigma_{\mathscr{R}}^{o}}}}_{\mathscr{N}_{2}(\widetilde{\Sigma}_{p}^{o})})^{-
1} \! =_{n \to \infty} \! \mathcal{O}(1)$. \hfill $\qed$
\begin{ccccc}
Set $\Sigma_{\circlearrowright}^{o} \! := \! \Sigma_{p}^{o,5}$ $(= \! \cup_{j
=1}^{N+1}(\partial \mathbb{U}_{\delta_{b_{j-1}}}^{o} \cup \partial \mathbb{U}_{
\delta_{a_{j}}}^{o}))$ and $\Sigma_{\scriptscriptstyle \blacksquare}^{o} \! :=
\! \widetilde{\Sigma}_{p}^{o} \setminus \Sigma_{\circlearrowright}^{o}$, and
let $\mathscr{R}^{o} \colon \mathbb{C} \setminus \widetilde{\Sigma}_{p}^{o}
\! \to \! \operatorname{SL}_{2}(\mathbb{C})$ solve the (equivalent) {\rm RHP}
$(\mathscr{R}^{o}(z),\upsilon_{\mathscr{R}}^{o}(z),\widetilde{\Sigma}_{p}^{o}
)$ formulated in Proposition~{\rm 5.2} with integral representation given by
Equation~{\rm (5.1)}. Let the asymptotic (as $n \! \to \! \infty)$ estimates
and bounds given in Propositions~{\rm 5.1} and~{\rm 5.2} be valid. Then,
uniformly for compact subsets of $\mathbb{C} \setminus \widetilde{\Sigma}_{
p}^{o} \! \ni \! z$,
\begin{equation*}
\mathscr{R}^{o}(z) \underset{\underset{z \in \mathbb{C} \setminus \widetilde{
\Sigma}^{o}_{p}}{n \to \infty}}{=} \mathrm{I} \! + \! \int_{\Sigma_{
\circlearrowright}^{o}} \dfrac{zw^{\Sigma_{\circlearrowright}^{o}}_{+}(s)
}{s(s \! - \! z)} \, \dfrac{\md s}{2 \pi \mi} \! + \! \mathcal{O} \! \left(
\dfrac{f(n)}{(n \! + \! \frac{1}{2})^{2} \min \{1,\operatorname{dist}
(z,\widetilde{\Sigma}_{p}^{o})\}} \right),
\end{equation*}
where $w_{+}^{\Sigma_{\circlearrowright}^{o}}(z) \! := \! w_{+}^{\Sigma_{
\mathscr{R}}^{o}}(z) \! \! \upharpoonright_{\Sigma_{\circlearrowright}^{o}}$,
and $(f(n))_{ij} \! =_{n \to \infty} \! \mathcal{O}(1)$, $i,j \! = \! 1,2$.
\end{ccccc}

\emph{Proof.} Define $\Sigma_{\circlearrowright}^{o}$ and $\Sigma^{o}_{
\scriptscriptstyle \blacksquare}$ as in the Lemma, and write $\widetilde{
\Sigma}_{p}^{o} \! = \! (\widetilde{\Sigma}_{p}^{o} \setminus \Sigma^{o}_{
\circlearrowright}) \cup \Sigma_{\circlearrowright}^{o} \! := \! \Sigma^{
o}_{\scriptscriptstyle \blacksquare} \cup \Sigma_{\circlearrowright}^{o}$
(with $\Sigma^{o}_{\scriptscriptstyle \blacksquare} \cap \Sigma_{
\circlearrowright}^{o} \! = \! \varnothing)$. Recall, {}from Equation~(5.1),
the integral representation for $\mathscr{R}^{o} \colon \mathbb{C} \setminus
\widetilde{\Sigma}_{p}^{o} \! \to \! \operatorname{SL}_{2}(\mathbb{C})$:
\begin{equation*}
\mathscr{R}^{o}(z) \! = \! \mathrm{I} \! + \! \int_{\widetilde{\Sigma}_{p}^{o}}
\dfrac{z \mu^{\Sigma_{\mathscr{R}}^{o}}(s)w_{+}^{\Sigma_{\mathscr{R}}^{o}}(s)}{
s(s \! - \! z)} \, \dfrac{\md s}{2 \pi \mi}, \quad z \! \in \! \mathbb{C}
\setminus \widetilde{\Sigma}_{p}^{o}.
\end{equation*}
Using the linearity property of the Cauchy integral operator $C^{0}_{w^{
\Sigma_{\mathscr{R}}^{o}}}$, one shows that $C^{0}_{w^{\Sigma_{\mathscr{R}}^{
o}}} \! = \! C^{0}_{w^{\Sigma_{\circlearrowright}^{o}}} \! + \! C^{0}_{w^{
\Sigma_{\scriptscriptstyle \blacksquare}^{o}}}$. Via a repeated application of
the second resolvent identity\footnote{For general operators $\mathscr{A}$ and
$\mathscr{B}$, if $(\id \! - \! \mathscr{A})^{-1}$ and $(\id \! - \! \mathscr{
B})^{-1}$ exist, then $(\id \! - \! \mathscr{B})^{-1} \! - \! (\id \! - \!
\mathscr{A})^{-1} \! = \! (\id \! - \! \mathscr{B})^{-1}(\mathscr{B} \! - \!
\mathscr{A})(\id \! - \! \mathscr{A})^{-1}$ \cite{a92}.}:
\begin{align*}
\mu^{\Sigma_{\mathscr{R}}^{o}}(z) &= \mathrm{I} \! + \! ((\id \! - \! C^{0}_{
w^{\Sigma_{\mathscr{R}}^{o}}})^{-1}C^{0}_{w^{\Sigma_{\mathscr{R}}^{o}}}
\mathrm{I})(z) \! = \! \mathrm{I} \! + \! ((\id \! - \! C^{0}_{w^{\Sigma_{
\circlearrowright}^{o}}} \! - \! C^{0}_{w^{\Sigma_{\scriptscriptstyle
\blacksquare}^{o}}})^{-1}(C^{0}_{w^{\Sigma_{\circlearrowright}^{o}}} \! + \!
C^{0}_{w^{\Sigma_{\scriptscriptstyle \blacksquare}^{o}}}) \mathrm{I})(z) \\
&=\mathrm{I} \! + \! ((\id \! - \! C^{0}_{w^{\Sigma_{\circlearrowright}^{o}}}
\! - \! C^{0}_{w^{\Sigma_{\scriptscriptstyle \blacksquare}^{o}}})^{-1}C^{0}_{
w^{\Sigma_{\circlearrowright}^{o}}} \mathrm{I})(z) \! + \! ((\id \! - \!
C^{0}_{w^{\Sigma_{\circlearrowright}^{o}}} \! - \! C^{0}_{w^{\Sigma_{
\scriptscriptstyle \blacksquare}^{o}}})^{-1}C^{0}_{w^{\Sigma_{
\scriptscriptstyle \blacksquare}^{o}}} \mathrm{I})(z) \\
&=\mathrm{I} \! + \! (((\id \! - \! C^{0}_{w^{\Sigma_{\circlearrowright}^{o}}})
(\id \! - \! (\id \! - \! C^{0}_{w^{\Sigma_{\circlearrowright}^{o}}})^{-1}C^{
0}_{w^{\Sigma_{\scriptscriptstyle \blacksquare}^{o}}}))^{-1}C^{0}_{w^{\Sigma_{
\circlearrowright}^{o}}} \mathrm{I})(z) \\
&+(((\id \! - \! C^{0}_{w^{\Sigma_{\scriptscriptstyle \blacksquare}^{o}}})
(\id \! - \! (\id \! - \! C^{0}_{w^{\Sigma_{\scriptscriptstyle \blacksquare}^{
o}}})^{-1}C^{0}_{w^{\Sigma_{\circlearrowright}^{o}}}))^{-1}
C^{0}_{w^{\Sigma_{\scriptscriptstyle \blacksquare}^{o}}} \mathrm{I})(z) \\
&=\mathrm{I} \! + \! ((\id \! - \! (\id \! - \! C^{0}_{w^{\Sigma_{
\circlearrowright}^{o}}})^{-1}C^{0}_{w^{\Sigma_{\scriptscriptstyle
\blacksquare}^{o}}})^{-1}(\id \! + \! (\id \! - \! C^{0}_{w^{\Sigma_{
\circlearrowright}^{o}}})^{-1}C^{0}_{w^{\Sigma_{\circlearrowright}^{o}}})
C^{0}_{w^{\Sigma_{\circlearrowright}^{o}}} \mathrm{I})(z) \\
&+((\id \! - \! (\id \! - \! C^{0}_{w^{\Sigma_{\scriptscriptstyle
\blacksquare}^{o}}})^{-1}C^{0}_{w^{\Sigma_{\circlearrowright}^{o}}})^{-1}
(\id \! + \! (\id \! - \! C^{0}_{w^{\Sigma_{\scriptscriptstyle \blacksquare}^{
o}}})^{-1}C^{0}_{w^{\Sigma_{\scriptscriptstyle \blacksquare}^{o}}})C^{0}_{w^{
\Sigma_{\scriptscriptstyle \blacksquare}^{o}}} \mathrm{I})(z) \\
&=\mathrm{I} \! + \! ((\id \! - \! (\id \! - \! C^{0}_{w^{\Sigma_{
\circlearrowright}^{o}}})^{-1}C^{0}_{w^{\Sigma_{\scriptscriptstyle
\blacksquare}^{o}}})^{-1}((\id \! - \! C^{0}_{w^{\Sigma_{\circlearrowright}^{
o}}})^{-1}C^{0}_{w^{\Sigma_{\circlearrowright}^{o}}})(C^{0}_{w^{\Sigma_{
\circlearrowright}^{o}}} \mathrm{I}))(z) \\
&+((\id \! - \! (\id \! - \! C^{0}_{w^{\Sigma_{\scriptscriptstyle
\blacksquare}^{o}}})^{-1}C^{0}_{w^{\Sigma_{\circlearrowright}^{o}}})^{-1}
((\id \! - \! C^{0}_{w^{\Sigma_{\scriptscriptstyle \blacksquare}^{o}}})^{
-1}C^{0}_{w^{\Sigma_{\scriptscriptstyle \blacksquare}^{o}}})(C^{0}_{w^{
\Sigma_{\scriptscriptstyle \blacksquare}^{o}}} \mathrm{I}))(z) \\
&+((\id \! - \! (\id \! - \! C^{0}_{w^{\Sigma_{\circlearrowright}^{o}}})^{
-1}C^{0}_{w^{\Sigma_{\scriptscriptstyle \blacksquare}^{o}}})^{-1}(C^{0}_{w^{
\Sigma_{\circlearrowright}^{o}}} \mathrm{I}))(z) \! + \! ((\id \! - \! (\id \!
- \! C^{0}_{w^{\Sigma_{\scriptscriptstyle \blacksquare}^{o}}})^{-1}C^{0}_{w^{
\Sigma_{\circlearrowright}^{o}}})^{-1}(C^{0}_{w^{\Sigma_{\scriptscriptstyle
\blacksquare}^{o}}} \mathrm{I}))(z) \\
&=\mathrm{I} \! + \! ((\id \! + \! (\id \! - \! (\id \! - \! C^{0}_{w^{
\Sigma_{\circlearrowright}^{o}}})^{-1}C^{0}_{w^{\Sigma_{\scriptscriptstyle
\blacksquare}^{o}}})^{-1}(\id \! - \! C^{0}_{w^{\Sigma_{\circlearrowright}^{
o}}})^{-1}C^{0}_{w^{\Sigma_{\scriptscriptstyle \blacksquare}^{o}}})(C^{0}_{
w^{\Sigma_{\circlearrowright}^{o}}} \mathrm{I}))(z) \\
&+((\id \! + \! (\id \! - \! (\id \! - \! C^{0}_{w^{\Sigma_{\scriptscriptstyle
\blacksquare}^{o}}})^{-1}C^{0}_{w^{\Sigma_{\circlearrowright}^{o}}})^{-1}(\id
\! - \! C^{0}_{w^{\Sigma_{\scriptscriptstyle \blacksquare}^{o}}})^{-1}C^{0}_{
w^{\Sigma_{\circlearrowright}^{o}}})(C^{0}_{w^{\Sigma_{\scriptscriptstyle
\blacksquare}^{o}}} \mathrm{I}))(z) \\
&+((\id \! + \! (\id \! - \! (\id \! - \! C^{0}_{w^{\Sigma_{
\circlearrowright}^{o}}})^{-1}C^{0}_{w^{\Sigma_{\scriptscriptstyle
\blacksquare}^{o}}})^{-1}(\id \! - \! C^{0}_{w^{\Sigma_{\circlearrowright}^{
o}}})^{-1}C^{0}_{w^{\Sigma_{\scriptscriptstyle \blacksquare}^{o}}})(\id \! -
\! C^{0}_{w^{\Sigma_{\circlearrowright}^{o}}})^{-1}C^{0}_{w^{\Sigma_{
\circlearrowright}^{o}}}(C^{0}_{w^{\Sigma_{\circlearrowright}^{o}}} \mathrm{
I}))(z) \\
&+((\id \! + \! (\id \! - \! (\id \! - \! C^{0}_{w^{\Sigma_{\scriptscriptstyle
\blacksquare}^{o}}})^{-1}C^{0}_{w^{\Sigma_{\circlearrowright}^{o}}})^{-1}(\id
\! - \! C^{0}_{w^{\Sigma_{\scriptscriptstyle \blacksquare}^{o}}})^{-1}C^{0}_{
w^{\Sigma_{\circlearrowright}^{o}}})(\id \! - \! C^{0}_{w^{\Sigma_{
\scriptscriptstyle \blacksquare}^{o}}})^{-1}C^{0}_{w^{\Sigma_{
\scriptscriptstyle \blacksquare}^{o}}}(C^{0}_{w^{\Sigma_{\scriptscriptstyle
\blacksquare}^{o}}} \mathrm{I}))(z) \\
&=\mathrm{I} \! + \! (C^{0}_{w^{\Sigma_{\circlearrowright}^{o}}} \mathrm{I})
(z) \! + \! (C^{0}_{w^{\Sigma_{\scriptscriptstyle \blacksquare}^{o}}} \mathrm{
I})(z) \! + \! ((\id \! - \! C^{0}_{w^{\Sigma_{\circlearrowright}^{o}}})^{-1}
C^{0}_{w^{\Sigma_{\circlearrowright}^{o}}}(C^{0}_{w^{\Sigma_{
\circlearrowright}^{o}}} \mathrm{I}))(z) \\
&+((\id \! - \! C^{0}_{w^{\Sigma_{\scriptscriptstyle \blacksquare}^{o}}})^{-1}
C^{0}_{w^{\Sigma_{\scriptscriptstyle \blacksquare}^{o}}}(C^{0}_{w^{\Sigma_{
\scriptscriptstyle \blacksquare}^{o}}} \mathrm{I}))(z) \! + \! ((\id \! - \!
(\id \! - \! C^{0}_{w^{\Sigma_{\circlearrowright}^{o}}})^{-1}C^{0}_{w^{\Sigma_{
\scriptscriptstyle \blacksquare}^{o}}})^{-1}(\id \! - \! C^{0}_{w^{\Sigma_{
\circlearrowright}^{o}}})^{-1}C^{0}_{w^{\Sigma_{\scriptscriptstyle
\blacksquare}^{o}}} \\
&\times (C^{0}_{w^{\Sigma_{\circlearrowright}^{o}}} \mathrm{I}))(z) \! + \!
((\id \! - \! (\id \! - \! C^{0}_{w^{\Sigma_{\scriptscriptstyle \blacksquare}^{
o}}})^{-1}C^{0}_{w^{\Sigma_{\circlearrowright}^{o}}})^{-1}(\id \! - \! C^{0}_{
w^{\Sigma_{\scriptscriptstyle \blacksquare}^{o}}})^{-1}C^{0}_{w^{\Sigma_{
\circlearrowright}^{o}}}(C^{0}_{w^{\Sigma_{\scriptscriptstyle \blacksquare}^{
o}}} \mathrm{I}))(z) \\
&+((\id \! - \! (\id \! - \! C^{0}_{w^{\Sigma_{\circlearrowright}^{o}}})^{-1}
C^{0}_{w^{\Sigma_{\scriptscriptstyle \blacksquare}^{o}}})^{-1}(\id \! - \! C^{
0}_{w^{\Sigma_{\circlearrowright}^{o}}})^{-1}C^{0}_{w^{\Sigma_{
\scriptscriptstyle \blacksquare}^{o}}}(\id \! - \! C^{0}_{w^{\Sigma_{
\circlearrowright}^{o}}})^{-1}C^{0}_{w^{\Sigma_{\circlearrowright}^{o}}}(C^{
0}_{w^{\Sigma_{\circlearrowright}^{o}}} \mathrm{I}))(z) \\
&+((\id \! - \! (\id \! - \! C^{0}_{w^{\Sigma_{\scriptscriptstyle
\blacksquare}^{o}}})^{-1}C^{0}_{w^{\Sigma_{\circlearrowright}^{o}}})^{-1}(\id
\! - \! C^{0}_{w^{\Sigma_{\scriptscriptstyle \blacksquare}^{o}}})^{-1}C^{0}_{
w^{\Sigma_{\circlearrowright}^{o}}}(\id \! - \! C^{0}_{w^{\Sigma_{
\scriptscriptstyle \blacksquare}^{o}}})^{-1}C^{0}_{w^{\Sigma_{
\scriptscriptstyle \blacksquare}^{o}}}(C^{0}_{w^{\Sigma_{\scriptscriptstyle
\blacksquare}^{o}}} \mathrm{I}))(z) \\
&=\mathrm{I} \! + \! (C^{0}_{w^{\Sigma_{\circlearrowright}^{o}}} \mathrm{I})
(z) \! + \! (C^{0}_{w^{\Sigma_{\scriptscriptstyle \blacksquare}^{o}}} \mathrm{
I})(z) \! + \! ((\id \! - \! C^{0}_{w^{\Sigma_{\circlearrowright}^{o}}})^{-1}
C^{0}_{w^{\Sigma_{\circlearrowright}^{o}}}(C^{0}_{w^{\Sigma_{
\circlearrowright}^{o}}} \mathrm{I}))(z) \! + \! ((\id \! - \! C^{0}_{w^{
\Sigma_{\scriptscriptstyle \blacksquare}^{o}}})^{-1}C^{0}_{w^{\Sigma_{
\scriptscriptstyle \blacksquare}^{o}}} \\
&\times \! (C^{0}_{w^{\Sigma_{\scriptscriptstyle \blacksquare}^{o}}} \mathrm{
I}))(z) \! + \! ((\id \! - \! (\id \! - \! C^{0}_{w^{\Sigma_{
\circlearrowright}^{o}}})^{-1}(\id \! - \! C^{0}_{w^{\Sigma_{
\scriptscriptstyle \blacksquare}^{o}}})^{-1}C^{0}_{w^{\Sigma_{
\scriptscriptstyle \blacksquare}^{o}}}C^{0}_{w^{\Sigma_{\circlearrowright}^{
o}}})^{-1}(\id \! - \! C^{0}_{w^{\Sigma_{\circlearrowright}^{o}}})^{-1}(\id
\! - \! C^{0}_{w^{\Sigma_{\scriptscriptstyle \blacksquare}^{o}}})^{-1} \\
&\times C^{0}_{w^{\Sigma_{\circlearrowright}^{o}}}(C^{0}_{w^{\Sigma_{
\scriptscriptstyle \blacksquare}^{o}}} \mathrm{I}))(z) \! + \! ((\id \! - \!
(\id \! - \! C^{0}_{w^{\Sigma_{\scriptscriptstyle \blacksquare}^{o}}})^{-1}
(\id \! - \! C^{0}_{w^{\Sigma_{\circlearrowright}^{o}}})^{-1}C^{0}_{w^{
\Sigma_{\circlearrowright}^{o}}}C^{0}_{w^{\Sigma_{\scriptscriptstyle
\blacksquare}^{o}}})^{-1}(\id \! - \! C^{0}_{w^{\Sigma_{\scriptscriptstyle
\blacksquare}^{o}}})^{-1} \\
&\times (\id \! - \! C^{0}_{w^{\Sigma_{\circlearrowright}^{o}}})^{-1}C^{0}_{
w^{\Sigma_{\scriptscriptstyle \blacksquare}^{o}}}(C^{0}_{w^{\Sigma_{
\circlearrowright}^{o}}} \mathrm{I}))(z) \! + \! ((\id \! - \! (\id \! - \!
C^{0}_{w^{\Sigma_{\circlearrowright}^{o}}})^{-1}(\id \! - \! C^{0}_{w^{
\Sigma_{\scriptscriptstyle \blacksquare}^{o}}})^{-1}C^{0}_{w^{\Sigma_{
\scriptscriptstyle \blacksquare}^{o}}}C^{0}_{w^{\Sigma_{\circlearrowright}^{
o}}})^{-1} \\
&\times (\id \! - \! C^{0}_{w^{\Sigma_{\circlearrowright}^{o}}})^{-1}(\id \!
- \! C^{0}_{w^{\Sigma_{\scriptscriptstyle \blacksquare}^{o}}})^{-1}C^{0}_{
w^{\Sigma_{\circlearrowright}^{o}}}(\id \! - \! C^{0}_{w^{\Sigma_{
\scriptscriptstyle \blacksquare}^{o}}})^{-1}C^{0}_{w^{\Sigma_{
\scriptscriptstyle \blacksquare}^{o}}}(C^{0}_{w^{\Sigma_{\scriptscriptstyle
\blacksquare}^{o}}} \mathrm{I}))(z) \! + \! ((\id \! - \! (\id \! - \! C^{0}_{
w^{\Sigma_{\scriptscriptstyle \blacksquare}^{o}}})^{-1} \\
&\times (\id \! - \! C^{0}_{w^{\Sigma_{\circlearrowright}^{o}}})^{-1}C^{0}_{
w^{\Sigma_{\circlearrowright}^{o}}}C^{0}_{w^{\Sigma_{\scriptscriptstyle
\blacksquare}^{o}}})^{-1}(\id \! - \! C^{0}_{w^{\Sigma_{\scriptscriptstyle
\blacksquare}^{o}}})^{-1}(\id \! - \! C^{0}_{w^{\Sigma_{\circlearrowright}^{
o}}})^{-1}C^{0}_{w^{\Sigma_{\scriptscriptstyle \blacksquare}^{o}}}(\id \!
- \! C^{0}_{w^{\Sigma_{\circlearrowright}^{o}}})^{-1}C^{0}_{w^{\Sigma_{
\circlearrowright}^{o}}} \\
&\times (C^{0}_{w^{\Sigma_{\circlearrowright}^{o}}} \mathrm{I}))(z);
\end{align*}
hence, recalling the integral representation for $\mathscr{R}^{o}(z)$ given
above, one arrives at, for $\mathbb{C} \setminus \widetilde{\Sigma}_{p}^{o}
\! \ni \! z$, upon using the partial fraction decomposition $\tfrac{z}{s(s-z)}
\! = \! -\tfrac{1}{s} \! + \! \tfrac{1}{s-z}$,
\begin{equation}
\mathscr{R}^{o}(z) \! - \! \mathrm{I} \! - \! \int_{\Sigma_{\circlearrowright}
^{o}} \dfrac{zw^{\Sigma_{\circlearrowright}^{o}}_{+}(s)}{s(s \! - \! z)} \,
\dfrac{\md s}{2 \pi \mi} \! = \! \int_{\Sigma_{\scriptscriptstyle
\blacksquare}^{o}}w_{+}^{\Sigma_{\scriptscriptstyle \blacksquare}^{o}}(s) \!
\left(\dfrac{1}{s \! - \! z} \! - \! \dfrac{1}{s} \right) \dfrac{\md s}{2 \pi
\mi} \! + \! \sum_{k=1}^{8}I_{k}^{o},
\end{equation}
where $w_{+}^{\Sigma_{\circlearrowright}^{o}}(z) \! := \! w_{+}^{\Sigma_{
\mathscr{R}}^{o}}(z) \! \! \upharpoonright_{\Sigma_{\circlearrowright}^{o}}$,
$w_{+}^{\Sigma_{\scriptscriptstyle \blacksquare}^{o}}(z) \! := \! w_{+}^{
\Sigma_{\mathscr{R}}^{o}}(z) \! \! \upharpoonright_{\Sigma_{\scriptscriptstyle
\blacksquare}^{o}}$,
\begin{align*}
I_{1}^{o} :=& \int_{\widetilde{\Sigma}_{p}^{o}} (C^{0}_{w^{\Sigma_{
\scriptscriptstyle \blacksquare}^{o}}} \mathrm{I})(s)w_{+}^{\Sigma_{\mathscr{
R}}^{o}}(s) \! \left(\dfrac{1}{s \! - \! z} \! - \! \dfrac{1}{s} \right)
\dfrac{\md s}{2 \pi \mi}, \qquad \quad I_{2}^{o} \! := \! \int_{\widetilde{
\Sigma}_{p}^{o}} (C^{0}_{w^{\Sigma_{\circlearrowright}^{o}}} \mathrm{I})(s)
w_{+}^{\Sigma_{\mathscr{R}}^{o}}(s) \! \left(\dfrac{1}{s \! - \! z} \! - \!
\dfrac{1}{s} \right) \dfrac{\md s}{2 \pi \mi}, \\
I_{3}^{o} :=& \int_{\widetilde{\Sigma}_{p}^{o}} ((\id \! - \! C^{0}_{w^{
\Sigma_{\scriptscriptstyle \blacksquare}^{o}}})^{-1}C^{0}_{w^{\Sigma_{
\scriptscriptstyle \blacksquare}^{o}}}(C^{0}_{w^{\Sigma_{\scriptscriptstyle
\blacksquare}^{o}}} \mathrm{I}))(s)w_{+}^{\Sigma_{\mathscr{R}}^{o}}(s) \!
\left(\dfrac{1}{s \! - \! z} \! - \! \dfrac{1}{s} \right) \dfrac{\md s}{2 \pi
\mi}, \\
I_{4}^{o} :=& \int_{\widetilde{\Sigma}_{p}^{o}} ((\id \! - \! C^{0}_{w^{
\Sigma_{\circlearrowright}^{o}}})^{-1}C^{0}_{w^{\Sigma_{\circlearrowright}^{
o}}}(C^{0}_{w^{\Sigma_{\circlearrowright}^{o}}} \mathrm{I}))(s)w_{+}^{\Sigma_{
\mathscr{R}}^{o}}(s) \! \left(\dfrac{1}{s \! - \! z} \! - \! \dfrac{1}{s}
\right) \dfrac{\md s}{2 \pi \mi}, \\
I_{5}^{o} :=& \int_{\widetilde{\Sigma}_{p}^{o}}((\id \! - \! (\id \! - \! C^{
0}_{w^{\Sigma_{\circlearrowright}^{o}}})^{-1}(\id \! - \! C^{0}_{w^{\Sigma_{
\scriptscriptstyle \blacksquare}^{o}}})^{-1}C^{0}_{w^{\Sigma_{
\scriptscriptstyle \blacksquare}^{o}}}C^{0}_{w^{\Sigma_{\circlearrowright}^{
o}}})^{-1}(\id \! - \! C^{0}_{w^{\Sigma_{\circlearrowright}^{o}}})^{-1}(\id \!
- \! C^{0}_{w^{\Sigma_{\scriptscriptstyle \blacksquare}^{o}}})^{-1} \\
\times& \, C^{0}_{w^{\Sigma_{\circlearrowright}^{o}}}(C^{0}_{w^{\Sigma_{
\scriptscriptstyle \blacksquare}^{o}}} \mathrm{I}))(s)w_{+}^{\Sigma_{\mathscr{
R}}^{o}}(s) \! \left(\dfrac{1}{s \! - \! z} \! - \! \dfrac{1}{s} \right)
\dfrac{\md s}{2 \pi \mi}, \\
I_{6}^{o} :=& \int_{\widetilde{\Sigma}_{p}^{o}}((\id \! - \! (\id \! - \! C^{
0}_{w^{\Sigma_{\scriptscriptstyle \blacksquare}^{o}}})^{-1}(\id \! - \! C^{
0}_{w^{\Sigma_{\circlearrowright}^{o}}})^{-1}C^{0}_{w^{\Sigma_{
\circlearrowright}^{o}}}C^{0}_{w^{\Sigma_{\scriptscriptstyle \blacksquare}^{
o}}})^{-1}(\id \! - \! C^{0}_{w^{\Sigma_{\scriptscriptstyle \blacksquare}^{
o}}})^{-1}(\id \! - \! C^{0}_{w^{\Sigma_{\circlearrowright}^{o}}})^{-1} \\
\times& \, C^{0}_{w^{\Sigma_{\scriptscriptstyle \blacksquare}^{o}}}(C^{0}_{
w^{\Sigma_{\circlearrowright}^{o}}} \mathrm{I}))(s)w_{+}^{\Sigma_{\mathscr{
R}}^{o}}(s) \! \left(\dfrac{1}{s \! - \! z} \! - \! \dfrac{1}{s} \right)
\dfrac{\md s}{2 \pi \mi}, \\
I_{7}^{o} :=& \int_{\widetilde{\Sigma}_{p}^{o}}((\id \! - \! (\id \! - \! C^{
0}_{w^{\Sigma_{\circlearrowright}^{o}}})^{-1}(\id \! - \! C^{0}_{w^{\Sigma_{
\scriptscriptstyle \blacksquare}^{o}}})^{-1}C^{0}_{w^{\Sigma_{
\scriptscriptstyle \blacksquare}^{o}}}C^{0}_{w^{\Sigma_{\circlearrowright}^{
o}}})^{-1}(\id \! - \! C^{0}_{w^{\Sigma_{\circlearrowright}^{o}}})^{-1}(\id \!
- \! C^{0}_{w^{\Sigma_{\scriptscriptstyle \blacksquare}^{o}}})^{-1} \\
\times& \, C^{0}_{w^{\Sigma_{\circlearrowright}^{o}}}(\id \! - \! C^{0}_{
w^{\Sigma_{\scriptscriptstyle \blacksquare}^{o}}})^{-1}C^{0}_{w^{\Sigma_{
\scriptscriptstyle \blacksquare}^{o}}}(C^{0}_{w^{\Sigma_{\scriptscriptstyle
\blacksquare}^{o}}} \mathrm{I}))(s)w_{+}^{\Sigma_{\mathscr{R}}^{o}}(s) \!
\left(\dfrac{1}{s \! - \! z} \! - \! \dfrac{1}{s} \right) \dfrac{\md s}{2 \pi
\mi}, \\
I_{8}^{o} :=& \int_{\widetilde{\Sigma}_{p}^{o}}((\id \! - \! (\id \! - \! C^{
0}_{w^{\Sigma_{\scriptscriptstyle \blacksquare}^{o}}})^{-1}(\id \! - \! C^{
0}_{w^{\Sigma_{\circlearrowright}^{o}}})^{-1}C^{0}_{w^{\Sigma_{
\circlearrowright}^{o}}}C^{0}_{w^{\Sigma_{\scriptscriptstyle \blacksquare}^{
o}}})^{-1}(\id \! - \! C^{0}_{w^{\Sigma_{\scriptscriptstyle \blacksquare}^{
o}}})^{-1}(\id \! - \! C^{0}_{w^{\Sigma_{\circlearrowright}^{o}}})^{-1} \\
\times& \, C^{0}_{w^{\Sigma_{\scriptscriptstyle \blacksquare}^{o}}}(\id \!
- \! C^{0}_{w^{\Sigma_{\circlearrowright}^{o}}})^{-1}C^{0}_{w^{\Sigma_{
\circlearrowright}^{o}}}(C^{0}_{w^{\Sigma_{\circlearrowright}^{o}}} \mathrm{
I}))(s)w_{+}^{\Sigma_{\mathscr{R}}^{o}}(s) \! \left(\dfrac{1}{s \! - \! z} \!
- \! \dfrac{1}{s} \right) \dfrac{\md s}{2 \pi \mi}.
\end{align*}
One now proceeds to estimate, as $n \! \to \! \infty$, and without loss of
generality, the respective terms on the right-hand side of Equation~(5.2)
corresponding to the (standard) Cauchy kernel, $\tfrac{1}{s-z} \, \tfrac{
\md s}{2 \pi \mi}$, using the estimates and bounds given in Propositions~5.1
and~5.2.
\begin{equation*}
\left\vert \int_{\Sigma_{\scriptscriptstyle \blacksquare}^{o}} \dfrac{w_{+}^{
\Sigma_{\scriptscriptstyle \blacksquare}^{o}}(s)}{s \! - \! z} \, \dfrac{\md
s}{2 \pi \mi} \right\vert \! \leqslant \! \int_{\Sigma_{\scriptscriptstyle
\blacksquare}^{o}} \dfrac{\vert w_{+}^{\Sigma_{\scriptscriptstyle
\blacksquare}^{o}}(s) \vert}{\vert s \! - \! z \vert} \, \dfrac{\vert \md s
\vert}{2 \pi} \! \leqslant \! \dfrac{\norm{w_{+}^{\Sigma_{\scriptscriptstyle
\blacksquare}^{o}}(\cdot)}_{\mathcal{L}^{1}_{\mathrm{M}_{2}(\mathbb{C})}
(\Sigma_{\scriptscriptstyle \blacksquare}^{o})}}{2 \pi \operatorname{dist}
(z,\Sigma_{\scriptscriptstyle \blacksquare}^{o})} \! \underset{n \to \infty}{
\leqslant} \! \mathcal{O} \! \left(\dfrac{f(n) \me^{-(n+\frac{1}{2})c}}{(n \!
+ \! \frac{1}{2}) \operatorname{dist}(z,\widetilde{\Sigma}_{p}^{o})} \right),
\end{equation*}
where, here and below, $(f(n) \! > \! 0$ and) $f(n) \! =_{n \to \infty} \!
\mathcal{O}(1)$, and $c \! > \! 0$. One estimates the `Cauchy part' of $I_{
1}^{o}$, denoted $I^{o,\mathcal{C}}_{1}$, as follows:
\begin{align*}
\vert I_{1}^{o,\mathcal{C}} \vert &\leqslant \int_{\widetilde{\Sigma}_{p}^{o}}
\dfrac{\vert (C^{0}_{w^{\Sigma_{\scriptscriptstyle \blacksquare}^{o}}} \mathrm{
I})(s) \vert \vert w_{+}^{\Sigma_{\mathscr{R}}^{o}}(s) \vert}{\vert s \! - \!
z \vert} \, \dfrac{\vert \md s \vert}{2 \pi} \! \leqslant \! \dfrac{\norm{(C^{
0}_{w^{\Sigma_{\scriptscriptstyle \blacksquare}^{o}}} \mathrm{I})(\cdot)}_{
\mathcal{L}^{2}_{\mathrm{M}_{2}(\mathbb{C})}(\widetilde{\Sigma}_{p}^{o})}
\norm{w_{+}^{\Sigma_{\mathscr{R}}^{o}}(\cdot)}_{\mathcal{L}^{2}_{\mathrm{M}_{
2}(\mathbb{C})}(\widetilde{\Sigma}_{p}^{o})}}{2 \pi \operatorname{dist}(z,
\widetilde{\Sigma}_{p}^{o})} \\
&\leqslant \dfrac{\operatorname{const.} \, \norm{w_{+}^{\Sigma_{
\scriptscriptstyle \blacksquare}^{o}}(\cdot)}_{\mathcal{L}^{2}_{\mathrm{M}_{2}
(\mathbb{C})}(\Sigma^{o}_{\scriptscriptstyle \blacksquare})}(\norm{w_{+}^{
\Sigma_{\circlearrowright}^{o}}(\cdot)}_{\mathcal{L}^{2}_{\mathrm{M}_{2}
(\mathbb{C})}(\Sigma^{o}_{\circlearrowright})} \! + \! \norm{w_{+}^{\Sigma_{
\scriptscriptstyle \blacksquare}^{o}}(\cdot)}_{\mathcal{L}^{2}_{\mathrm{M}_{
2}(\mathbb{C})}(\Sigma^{o}_{\scriptscriptstyle \blacksquare})})}{2 \pi
\operatorname{dist}(z,\widetilde{\Sigma}_{p}^{o})} \\
&\underset{n \to \infty}{\leqslant} \mathcal{O} \! \left(\dfrac{f(n) \me^{-
(n+\frac{1}{2})c}}{\sqrt{n \! + \! \frac{1}{2}} \, \operatorname{dist}(z,
\widetilde{\Sigma}_{p}^{o})} \right) \! \left(\! \mathcal{O} \! \left(\dfrac{
f(n)}{n \! + \! \frac{1}{2}} \right) \! + \! \mathcal{O} \! \left(\dfrac{f(n)
\me^{-(n+\frac{1}{2})c}}{\sqrt{n \! + \! \frac{1}{2}}} \right) \right) \\
&\underset{n \to \infty}{\leqslant} \! \mathcal{O} \! \left(\dfrac{f(n) \me^{
-(n+\frac{1}{2})c}}{(n \! + \! \frac{1}{2}) \operatorname{dist}(z,\widetilde{
\Sigma}_{p}^{o})} \right),
\end{align*}
where, here and below, $\operatorname{const.}$ denotes some positive,
$\mathcal{O}(1)$ constant; in going {}from the first to the second (resp.,
second to the third) line in the above asymptotic (as $n \! \to \! \infty)$
estimation for $I_{1}^{o,\mathcal{C}}$, one uses the fact that $\norm{(C^{
0}_{w^{\Sigma^{o}_{\scriptscriptstyle \blacksquare}}} \mathrm{I})(\cdot)}_{
\mathcal{L}^{2}_{\mathrm{M}_{2}(\mathbb{C})}(\Sigma^{o}_{\circlearrowright})}
\leqslant_{n \to \infty} \mathcal{O}(f(n)(n \! + \! 1/2)^{-1} \me^{-(n+\frac{
1}{2})c})$ and $\norm{(C^{0}_{w^{\Sigma^{o}_{\scriptscriptstyle \blacksquare}
}} \mathrm{I})(\cdot)}_{\mathcal{L}^{2}_{\mathrm{M}_{2}(\mathbb{C})}(\Sigma^{
o}_{\scriptscriptstyle \blacksquare})} \linebreak[4]
\leqslant_{n \to \infty} \! \mathcal{O}(f(n)(n \! + \! 1/2)^{-1/2} \me^{-(n+
\frac{1}{2})c})$ (resp., for $a,b \! > \! 0$, $\sqrt{\smash[b]{a^{2} \! + \!
b^{2}}} \! \leqslant \! \sqrt{\smash[b]{a^{2}}} \! + \! \sqrt{\smash[b]{b^{2}}
})$ (facts used repeatedly below). One estimates the Cauchy part of $I_{2}^{
o}$, denoted $I_{2}^{o,\mathcal{C}}$, as follows:
\begin{align*}
\vert I_{2}^{o,\mathcal{C}} \vert &\leqslant \int_{\widetilde{\Sigma}_{p}^{o}}
\dfrac{\vert (C^{0}_{w^{\Sigma_{\circlearrowright}^{o}}} \mathrm{I})(s) \vert
\vert w_{+}^{\Sigma_{\mathscr{R}}^{o}}(s) \vert}{\vert s \! - \! z \vert} \,
\dfrac{\vert \md s \vert}{2 \pi} \! \leqslant \! \dfrac{\norm{(C^{0}_{w^{
\Sigma_{\circlearrowright}^{o}}} \mathrm{I})(\cdot)}_{\mathcal{L}^{2}_{\mathrm{
M}_{2}(\mathbb{C})}(\widetilde{\Sigma}_{p}^{o})} \norm{w_{+}^{\Sigma_{\mathscr{
R}}^{o}}(\cdot)}_{\mathcal{L}^{2}_{\mathrm{M}_{2}(\mathbb{C})}(\widetilde{
\Sigma}_{p}^{o})}}{2 \pi \operatorname{dist}(z,\widetilde{\Sigma}_{p}^{o})} \\
&\leqslant \dfrac{\operatorname{const.} \, \norm{w_{+}^{\Sigma_{
\circlearrowright}^{o}}(\cdot)}_{\mathcal{L}^{2}_{\mathrm{M}_{2}(\mathbb{C})}
(\Sigma^{o}_{\circlearrowright})}(\norm{w_{+}^{\Sigma_{\circlearrowright}^{o}}
(\cdot)}_{\mathcal{L}^{2}_{\mathrm{M}_{2}(\mathbb{C})}(\Sigma^{o}_{
\circlearrowright})} \! + \! \norm{w_{+}^{\Sigma_{\scriptscriptstyle
\blacksquare}^{o}}(\cdot)}_{\mathcal{L}^{2}_{\mathrm{M}_{2}(\mathbb{C})}
(\Sigma^{o}_{\scriptscriptstyle \blacksquare})})}{2 \pi \operatorname{dist}
(z,\widetilde{\Sigma}_{p}^{o})} \\
&\underset{n \to \infty}{\leqslant} \mathcal{O} \! \left(\dfrac{f(n)}{(n \! +
\! \frac{1}{2}) \operatorname{dist}(z,\widetilde{\Sigma}_{p}^{o})} \right) \!
\left(\! \mathcal{O} \! \left(\dfrac{f(n)}{n \! + \! \frac{1}{2}} \right) \! +
\! \mathcal{O} \! \left(\dfrac{f(n) \me^{-(n+\frac{1}{2})c}}{\sqrt{n \! + \!
\frac{1}{2}}} \right) \right) \\
&\underset{n \to \infty}{\leqslant} \! \mathcal{O} \! \left(\dfrac{f(n)}{(n \!
+ \! \frac{1}{2})^{2} \operatorname{dist}(z,\widetilde{\Sigma}_{p}^{o})}
\right):
\end{align*}
in going {}from the second to the third line in the above asymptotic (as $n \!
\to \! \infty)$ estimation for $I_{2}^{o,\mathcal{C}}$, one uses the fact that
$\norm{(C^{0}_{w^{\Sigma^{o}_{\circlearrowright}}} \mathrm{I})(\cdot)}_{
\mathcal{L}^{2}_{\mathrm{M}_{2}(\mathbb{C})}(\Sigma^{o}_{\circlearrowright})}
\! \leqslant_{n \to \infty} \! \mathcal{O}(f(n)(n \! + \! 1/2)^{-1})$ and
$\norm{(C^{0}_{w^{\Sigma^{o}_{\circlearrowright}}} \mathrm{I})(\cdot)}_{
\mathcal{L}^{2}_{\mathrm{M}_{2}(\mathbb{C})}(\Sigma^{o}_{\scriptscriptstyle
\blacksquare})} \! \leqslant_{n \to \infty} \! \mathcal{O}(f(n)(n \! + \!
1/2)^{-1})$. One estimates the Cauchy part of $I_{3}^{o}$, denoted $I_{3}^{o,
\mathcal{C}}$, as follows:
\begin{align*}
\vert I_{3}^{o,\mathcal{C}} \vert &\leqslant \int_{\widetilde{\Sigma}_{p}^{o}}
\dfrac{\vert ((\id \! - \! C^{0}_{w^{\Sigma_{\scriptscriptstyle \blacksquare}^{
o}}})^{-1}C^{0}_{w^{\Sigma_{\scriptscriptstyle \blacksquare}^{o}}}(C^{0}_{w^{
\Sigma_{\scriptscriptstyle \blacksquare}^{o}}} \mathrm{I}))(s) \vert \vert
w_{+}^{\Sigma_{\mathscr{R}}^{o}}(s) \vert}{\vert s \! - \! z \vert} \, \dfrac{
\vert \md s \vert}{2 \pi} \\
&\leqslant \dfrac{\norm{((\id \! - \! C^{0}_{w^{\Sigma_{\scriptscriptstyle
\blacksquare}^{o}}})^{-1}C^{0}_{w^{\Sigma_{\scriptscriptstyle \blacksquare}^{
o}}}(C^{0}_{w^{\Sigma_{\scriptscriptstyle \blacksquare}^{o}}} \mathrm{I}))
(\cdot)}_{\mathcal{L}^{2}_{\mathrm{M}_{2}(\mathbb{C})}(\widetilde{\Sigma}_{
p}^{o})} \norm{w_{+}^{\Sigma_{\mathscr{R}}^{o}}(\cdot)}_{\mathcal{L}^{2}_{
\mathrm{M}_{2}(\mathbb{C})}(\widetilde{\Sigma}_{p}^{o})}}{2 \pi
\operatorname{dist}(z,\widetilde{\Sigma}_{p}^{o})} \\
&\leqslant \dfrac{\norm{(\id \! - \! C^{0}_{w^{\Sigma_{\scriptscriptstyle
\blacksquare}^{o}}})^{-1}}_{\mathscr{N}_{2}(\widetilde{\Sigma}_{p}^{o})}
\norm{C^{0}_{w^{\Sigma_{\scriptscriptstyle \blacksquare}^{o}}}}_{\mathscr{
N}_{2}(\widetilde{\Sigma}_{p}^{o})} \norm{(C^{0}_{w^{\Sigma_{
\scriptscriptstyle \blacksquare}^{o}}} \mathrm{I})(\cdot)}_{\mathcal{L}^{2}_{
\mathrm{M}_{2}(\mathbb{C})}(\widetilde{\Sigma}_{p}^{o})} \norm{w_{+}^{\Sigma_{
\mathscr{R}}^{o}}(\cdot)}_{\mathcal{L}^{2}_{\mathrm{M}_{2}(\mathbb{C})}
(\widetilde{\Sigma}_{p}^{o})}}{2 \pi \operatorname{dist}(z,\widetilde{\Sigma}_{
p}^{o})} \\
&\leqslant \dfrac{\operatorname{const.} \, \norm{(\id \! - \! C^{0}_{w^{
\Sigma_{\scriptscriptstyle \blacksquare}^{o}}})^{-1}}_{\mathscr{N}_{2}
(\widetilde{\Sigma}_{p}^{o})} \norm{C^{0}_{w^{\Sigma_{\scriptscriptstyle
\blacksquare}^{o}}}}_{\mathscr{N}_{2}(\widetilde{\Sigma}_{p}^{o})} \norm{w_{
+}^{\Sigma_{\scriptscriptstyle \blacksquare}^{o}}(\cdot)}_{\mathcal{L}^{2}_{
\mathrm{M}_{2}(\mathbb{C})}(\Sigma_{\scriptscriptstyle \blacksquare}^{o})}}{
2 \pi \operatorname{dist}(z,\widetilde{\Sigma}_{p}^{o})} \\
&\times \left(\norm{w_{+}^{\Sigma_{\circlearrowright}^{o}}(\cdot)}_{\mathcal{
L}^{2}_{\mathrm{M}_{2}(\mathbb{C})}(\Sigma^{o}_{\circlearrowright})} \! + \!
\norm{w_{+}^{\Sigma_{\scriptscriptstyle \blacksquare}^{o}}(\cdot)}_{\mathcal{
L}^{2}_{\mathrm{M}_{2}(\mathbb{C})}(\Sigma^{o}_{\scriptscriptstyle
\blacksquare})} \right);
\end{align*}
using the fact that (cf. Proposition~5.2) $\norm{(\id \! - \! C^{0}_{w^{
\Sigma_{\scriptscriptstyle \blacksquare}^{o}}})^{-1}}_{\mathscr{N}_{2}
(\widetilde{\Sigma}_{p}^{o})} \! =_{n \to \infty} \! \mathcal{O}(1)$ (via a
Neuman series inversion argument, since $\norm{C^{0}_{w^{\Sigma_{
\scriptscriptstyle \blacksquare}^{o}}}}_{\mathscr{N}_{2}(\widetilde{\Sigma}_{
p}^{o})} \! \leqslant_{n \to \infty} \! \mathcal{O}((n \! + \! 1/2)^{-1/2}f(n)
\me^{-(n+\frac{1}{2})c}))$, one gets that
\begin{align*}
\vert I_{3}^{o,\mathcal{C}} \vert \underset{n \to \infty}{\leqslant}& \,
\mathcal{O} \! \left(\dfrac{f(n) \me^{-(n+\frac{1}{2})c}}{(n \! + \! \frac{1}{
2}) \operatorname{dist}(z,\widetilde{\Sigma}_{p}^{o})} \right) \! \mathcal{O}
\! \left(\dfrac{f(n) \me^{-(n+\frac{1}{2})c}}{\sqrt{n \! + \! \frac{1}{2}}}
\right) \! \left(\! \mathcal{O} \! \left(\dfrac{f(n)}{n \! + \! \frac{1}{2}}
\right) \! + \! \mathcal{O} \! \left(\dfrac{f(n) \me^{-(n+\frac{1}{2})c}}{
\sqrt{n \! + \! \frac{1}{2}}} \right) \right) \\
\underset{n \to \infty}{\leqslant}& \, \mathcal{O} \! \left(\dfrac{f(n) \me^{-
(n+\frac{1}{2})c}}{(n \! + \! \frac{1}{2})^{2} \operatorname{dist}(z,
\widetilde{\Sigma}_{p}^{o})} \right).
\end{align*}
One estimates the Cauchy part of $I_{4}^{o}$, denoted $I_{4}^{o,\mathcal{C}}$,
as follows:
\begin{align*}
\vert I_{4}^{o,\mathcal{C}} \vert &\leqslant \int_{\widetilde{\Sigma}_{p}^{o}}
\dfrac{\vert ((\id \! - \! C^{0}_{w^{\Sigma_{\circlearrowright}^{o}}})^{-1}
C^{0}_{w^{\Sigma_{\circlearrowright}^{o}}}(C^{0}_{w^{\Sigma_{
\circlearrowright}^{o}}} \mathrm{I}))(s) \vert \vert w_{+}^{\Sigma_{\mathscr{
R}}^{o}}(s) \vert}{\vert s \! - \! z \vert} \, \dfrac{\vert \md s \vert}{2
\pi} \\
&\leqslant \dfrac{\norm{((\id \! - \! C^{0}_{w^{\Sigma_{\circlearrowright}^{
o}}})^{-1}C^{0}_{w^{\Sigma_{\circlearrowright}^{o}}}(C^{0}_{w^{\Sigma_{
\circlearrowright}^{o}}} \mathrm{I}))(\cdot)}_{\mathcal{L}^{2}_{\mathrm{M}_{2}
(\mathbb{C})}(\widetilde{\Sigma}_{p}^{o})} \norm{w_{+}^{\Sigma_{\mathscr{R}}^{
o}}(\cdot)}_{\mathcal{L}^{2}_{\mathrm{M}_{2}(\mathbb{C})}(\widetilde{\Sigma}_{
p}^{o})}}{2 \pi \operatorname{dist}(z,\widetilde{\Sigma}_{p}^{o})} \\
&\leqslant \dfrac{\norm{(\id \! - \! C^{0}_{w^{\Sigma_{\circlearrowright}^{
o}}})^{-1}}_{\mathscr{N}_{2}(\widetilde{\Sigma}_{p}^{o})} \norm{C^{0}_{w^{
\Sigma_{\circlearrowright}^{o}}}}_{\mathscr{N}_{2}(\widetilde{\Sigma}_{p}^{
o})} \norm{(C^{0}_{w^{\Sigma_{\circlearrowright}^{o}}} \mathrm{I})(\cdot)}_{
\mathcal{L}^{2}_{\mathrm{M}_{2}(\mathbb{C})}(\widetilde{\Sigma}_{p}^{o})}
\norm{w_{+}^{\Sigma_{\mathscr{R}}^{o}}(\cdot)}_{\mathcal{L}^{2}_{\mathrm{M}_{2}
(\mathbb{C})}(\widetilde{\Sigma}_{p}^{o})}}{2 \pi \operatorname{dist}(z,
\widetilde{\Sigma}_{p}^{o})} \\
&\leqslant \dfrac{\operatorname{const.} \, \norm{(\id \! - \! C^{0}_{w^{
\Sigma_{\circlearrowright}^{o}}})^{-1}}_{\mathscr{N}_{2}(\widetilde{\Sigma}_{
p}^{o})} \norm{C^{0}_{w^{\Sigma_{\circlearrowright}^{o}}}}_{\mathscr{N}_{2}
(\widetilde{\Sigma}_{p}^{o})} \norm{w_{+}^{\Sigma_{\circlearrowright}^{o}}
(\cdot)}_{\mathcal{L}^{2}_{\mathrm{M}_{2}(\mathbb{C})}(\Sigma_{
\circlearrowright}^{o})}}{2 \pi \operatorname{dist}(z,\widetilde{\Sigma}_{
p}^{o})} \\
&\times \left(\norm{w_{+}^{\Sigma_{\circlearrowright}^{o}}(\cdot)}_{\mathcal{
L}^{2}_{\mathrm{M}_{2}(\mathbb{C})}(\Sigma^{o}_{\circlearrowright})} \! + \!
\norm{w_{+}^{\Sigma_{\scriptscriptstyle \blacksquare}^{o}}(\cdot)}_{\mathcal{
L}^{2}_{\mathrm{M}_{2}(\mathbb{C})}(\Sigma^{o}_{\scriptscriptstyle
\blacksquare})} \right);
\end{align*}
using the fact that (cf. Proposition~5.2) $\norm{(\id \! - \! C^{0}_{w^{
\Sigma_{\circlearrowright}^{o}}})^{-1}}_{\mathscr{N}_{2}(\widetilde{\Sigma}_{
p}^{o})} \! =_{n \to \infty} \! \mathcal{O}(1)$ (via a Neuman series inversion
argument, since $\norm{C^{0}_{w^{\Sigma_{\circlearrowright}^{o}}}}_{\mathscr{
N}_{2}(\widetilde{\Sigma}_{p}^{o})} \! =_{n \to \infty} \! \mathcal{O}((n \! +
\! 1/2)^{-1}f(n)))$, one gets that
\begin{align*}
\vert I_{4}^{o,\mathcal{C}} \vert \underset{n \to \infty}{\leqslant}& \,
\mathcal{O} \! \left(\dfrac{f(n)}{(n \! + \! \frac{1}{2}) \operatorname{dist}
(z,\widetilde{\Sigma}_{p}^{o})} \right) \! \mathcal{O} \! \left(\dfrac{f(n)}{
n \! + \! \frac{1}{2}} \right) \! \left(\! \mathcal{O} \! \left(\dfrac{f(n)}{
n \! + \! \frac{1}{2}} \right) \! + \! \mathcal{O} \! \left(\dfrac{f(n) \me^{
-(n+\frac{1}{2})c}}{\sqrt{n \! + \! \frac{1}{2}}} \right) \right) \\
\underset{n \to \infty}{\leqslant}& \, \mathcal{O} \! \left(\dfrac{f(n)}{(n
\! + \! \frac{1}{2})^{3} \operatorname{dist}(z,\widetilde{\Sigma}_{p}^{o})}
\right).
\end{align*}
One estimates the Cauchy part of $I_{5}^{o}$, denoted $I_{5}^{o,\mathcal{C}}$,
as follows:
\begin{align*}
\vert I_{5}^{o,\mathcal{C}} \vert &\leqslant \int_{\widetilde{\Sigma}_{p}^{o}}
\vert ((\id \! - \! (\id \! - \! C^{0}_{w^{\Sigma_{\circlearrowright}^{o}}})^{
-1}(\id \! - \! C^{0}_{w^{\Sigma_{\scriptscriptstyle \blacksquare}^{o}}})^{-1}
C^{0}_{w^{\Sigma_{\scriptscriptstyle \blacksquare}^{o}}}C^{0}_{w^{\Sigma_{
\circlearrowright}^{o}}})^{-1}(\id \! - \! C^{0}_{w^{\Sigma_{
\circlearrowright}^{o}}})^{-1}(\id \! - \! C^{0}_{w^{\Sigma_{
\scriptscriptstyle \blacksquare}^{o}}})^{-1} \\
&\times \, \dfrac{C^{0}_{w^{\Sigma_{\circlearrowright}^{o}}}(C^{0}_{w^{
\Sigma_{\scriptscriptstyle \blacksquare}^{o}}} \mathrm{I}))(s) \vert \vert
w_{+}^{\Sigma_{\mathscr{R}}^{o}}(s) \vert}{\vert s \! - \! z \vert} \, \dfrac{
\vert \md s \vert}{2 \pi} \\
&\leqslant \, \vert \vert ((\id \! - \! (\id \! - \! C^{0}_{w^{\Sigma_{
\circlearrowright}^{o}}})^{-1}(\id \! - \! C^{0}_{w^{\Sigma_{
\scriptscriptstyle \blacksquare}^{o}}})^{-1}C^{0}_{w^{\Sigma_{
\scriptscriptstyle \blacksquare}^{o}}}C^{0}_{w^{\Sigma_{\circlearrowright}^{
o}}})^{-1}(\id \! - \! C^{0}_{w^{\Sigma_{\circlearrowright}^{o}}})^{-1}(\id \!
- \! C^{0}_{w^{\Sigma_{\scriptscriptstyle \blacksquare}^{o}}})^{-1} \\
&\times \, \dfrac{C^{0}_{w^{\Sigma_{\circlearrowright}^{o}}}(C^{0}_{w^{
\Sigma_{\scriptscriptstyle \blacksquare}^{o}}} \mathrm{I}))(\cdot) \vert
\vert_{\mathcal{L}^{2}_{\mathrm{M}_{2}(\mathbb{C})}(\widetilde{\Sigma}_{p}^{
o})} \norm{w_{+}^{\Sigma_{\mathscr{R}}^{o}}(\cdot)}_{\mathcal{L}^{2}_{\mathrm{
M}_{2}(\mathbb{C})}(\widetilde{\Sigma}_{p}^{o})}}{2 \pi \operatorname{dist}
(z,\widetilde{\Sigma}_{p}^{o})} \\
&\leqslant \, \norm{(\id \! - \! (\id \! - \! C^{0}_{w^{\Sigma_{
\circlearrowright}^{o}}})^{-1}(\id \! - \! C^{0}_{w^{\Sigma_{
\scriptscriptstyle \blacksquare}^{o}}})^{-1}C^{0}_{w^{\Sigma_{
\scriptscriptstyle \blacksquare}^{o}}}C^{0}_{w^{\Sigma_{\circlearrowright}^{
o}}})^{-1}}_{\mathscr{N}_{2}(\widetilde{\Sigma}_{p}^{o})} \norm{(\id \! - \!
C^{0}_{w^{\Sigma_{\circlearrowright}^{o}}})^{-1}}_{\mathscr{N}_{2}(\widetilde{
\Sigma}_{p}^{o})} \\
&\times \, \dfrac{\norm{(\id \! - \! C^{0}_{w^{\Sigma_{\scriptscriptstyle
\blacksquare}^{o}}})^{-1}}_{\mathscr{N}_{2}(\widetilde{\Sigma}_{p}^{o})}
\norm{C^{0}_{w^{\Sigma_{\circlearrowright}^{o}}}}_{\mathscr{N}_{2}
(\widetilde{\Sigma}_{p}^{o})} \norm{(C^{0}_{w^{\Sigma_{\scriptscriptstyle
\blacksquare}^{o}}} \mathrm{I})(\cdot)}_{\mathcal{L}^{2}_{\mathrm{M}_{2}
(\mathbb{C})}(\widetilde{\Sigma}_{p}^{o})} \norm{w_{+}^{\Sigma_{\mathscr{R}}^{
o}}(\cdot)}_{\mathcal{L}^{2}_{\mathrm{M}_{2}(\mathbb{C})}(\widetilde{\Sigma}_{
p}^{o})}}{2 \pi \operatorname{dist}(z,\widetilde{\Sigma}_{p}^{o})} \\
&\leqslant \, \norm{(\id \! - \! (\id \! - \! C^{0}_{w^{\Sigma_{
\circlearrowright}^{o}}})^{-1}(\id \! - \! C^{0}_{w^{\Sigma_{
\scriptscriptstyle \blacksquare}^{o}}})^{-1}C^{0}_{w^{\Sigma_{
\scriptscriptstyle \blacksquare}^{o}}}C^{0}_{w^{\Sigma_{\circlearrowright}^{
o}}})^{-1}}_{\mathscr{N}_{2}(\widetilde{\Sigma}_{p}^{o})} \norm{(\id \! -
\! C^{0}_{w^{\Sigma_{\circlearrowright}^{o}}})^{-1}}_{\mathscr{N}_{2}
(\widetilde{\Sigma}_{p}^{o})} \\
&\times \, \operatorname{const.} \, \norm{(\id \! - \! C^{0}_{w^{\Sigma_{
\scriptscriptstyle \blacksquare}^{o}}})^{-1}}_{\mathscr{N}_{2}(\widetilde{
\Sigma}_{p}^{o})} \norm{C^{0}_{w^{\Sigma_{\circlearrowright}^{o}}}}_{\mathscr{
N}_{2}(\widetilde{\Sigma}_{p}^{o})} \norm{w_{+}^{\Sigma_{\scriptscriptstyle
\blacksquare}^{o}}(\cdot)}_{\mathcal{L}^{2}_{\mathrm{M}_{2}(\mathbb{C})}
(\Sigma_{\scriptscriptstyle \blacksquare}^{o})} \\
&\times \dfrac{(\norm{w_{+}^{\Sigma_{\circlearrowright}^{o}}(\cdot)}_{\mathcal{
L}^{2}_{\mathrm{M}_{2}(\mathbb{C})}(\Sigma^{o}_{\circlearrowright})} \! + \!
\norm{w_{+}^{\Sigma_{\scriptscriptstyle \blacksquare}^{o}}(\cdot)}_{\mathcal{
L}^{2}_{\mathrm{M}_{2}(\mathbb{C})}(\Sigma^{o}_{\scriptscriptstyle
\blacksquare})})}{2 \pi \operatorname{dist}(z,\widetilde{\Sigma}_{p}^{o})};
\end{align*}
using the fact that (cf. Proposition~5.2) $\norm{(\id \! - \! (\id \! - \! C^{
0}_{w^{\Sigma_{\circlearrowright}^{o}}})^{-1}(\id \! - \! C^{0}_{w^{\Sigma_{
\scriptscriptstyle \blacksquare}^{o}}})^{-1}C^{0}_{w^{\Sigma_{
\scriptscriptstyle \blacksquare}^{o}}}C^{0}_{w^{\Sigma_{\circlearrowright}^{
o}}})^{-1}}_{\mathscr{N}_{2}(\widetilde{\Sigma}_{p}^{o})} \! =_{n \to \infty}
\! \mathcal{O}(1)$ (via a Neuman series inversion argument, since $\norm{C^{
0}_{w^{\Sigma_{\scriptscriptstyle \blacksquare}^{o}}}}_{\mathscr{N}_{2}
(\widetilde{\Sigma}_{p}^{o})} \! \leqslant_{n \to \infty} \! \mathcal{O}((n \!
+ \! \tfrac{1}{2})^{-1/2}f(n) \me^{-(n+\frac{1}{2})c})$ and $\norm{C^{0}_{w^{
\Sigma_{\circlearrowright}^{o}}}}_{\mathscr{N}_{2}(\widetilde{\Sigma}_{p}^{o}
)} \! =_{n \to \infty} \! \mathcal{O}((n \! + \! 1/2)^{-1}f(n)))$, one gets
that
\begin{align*}
\vert I_{5}^{o,\mathcal{C}} \vert \underset{n \to \infty}{\leqslant}& \,
\mathcal{O} \! \left(\dfrac{f(n)}{(n \! + \! \frac{1}{2}) \operatorname{dist}
(z,\widetilde{\Sigma}_{p}^{o})} \right) \! \mathcal{O} \! \left(\dfrac{f(n)
\me^{-(n+\frac{1}{2})c}}{\sqrt{n \! + \! \frac{1}{2}}} \right) \! \left(\!
\mathcal{O} \! \left(\dfrac{f(n)}{n \! + \! \frac{1}{2}} \right) \! + \!
\mathcal{O} \! \left(\dfrac{f(n) \me^{-(n+\frac{1}{2})c}}{\sqrt{n \! + \!
\frac{1}{2}}} \right) \right) \\
\underset{n \to \infty}{\leqslant}& \, \mathcal{O} \! \left(\dfrac{f(n) \me^{-
(n+\frac{1}{2})c}}{(n \! + \! \frac{1}{2})^{2} \operatorname{dist}(z,
\widetilde{\Sigma}_{p}^{o})} \right).
\end{align*}
One estimates the Cauchy part of $I_{6}^{o}$, denoted $I_{6}^{o,\mathcal{C}}$,
as follows:
\begin{align*}
\vert I_{6}^{o,\mathcal{C}} \vert &\leqslant \int_{\widetilde{\Sigma}_{p}^{o}}
\vert ((\id \! - \! (\id \! - \! C^{0}_{w^{\Sigma_{\scriptscriptstyle
\blacksquare}^{o}}})^{-1}(\id \! - \! C^{0}_{w^{\Sigma_{\circlearrowright}^{
o}}})^{-1}C^{0}_{w^{\Sigma_{\circlearrowright}^{o}}}C^{0}_{w^{\Sigma_{
\scriptscriptstyle \blacksquare}^{o}}})^{-1}(\id \! - \! C^{0}_{w^{\Sigma_{
\scriptscriptstyle \blacksquare}^{o}}})^{-1}(\id \! - \! C^{0}_{w^{\Sigma_{
\circlearrowright}^{o}}})^{-1} \\
&\times \, \dfrac{C^{0}_{w^{\Sigma_{\scriptscriptstyle \blacksquare}^{o}}}
(C^{0}_{w^{\Sigma_{\circlearrowright}^{o}}} \mathrm{I}))(s) \vert \vert
w_{+}^{\Sigma_{\mathscr{R}}^{o}}(s) \vert}{\vert s \! - \! z \vert} \, \dfrac{
\vert \md s \vert}{2 \pi} \\
&\leqslant \, \vert \vert ((\id \! - \! (\id \! - \! C^{0}_{w^{\Sigma_{
\scriptscriptstyle \blacksquare}^{o}}})^{-1}(\id \! - \! C^{0}_{w^{\Sigma_{
\circlearrowright}^{o}}})^{-1}C^{0}_{w^{\Sigma_{\circlearrowright}^{o}}}C^{
0}_{w^{\Sigma_{\scriptscriptstyle \blacksquare}^{o}}})^{-1}(\id \! - \! C^{
0}_{w^{\Sigma_{\scriptscriptstyle \blacksquare}^{o}}})^{-1}(\id \! - \! C^{
0}_{w^{\Sigma_{\circlearrowright}^{o}}})^{-1} \\
&\times \, \dfrac{C^{0}_{w^{\Sigma_{\scriptscriptstyle \blacksquare}^{o}}}
(C^{0}_{w^{\Sigma_{\circlearrowright}^{o}}} \mathrm{I}))(\cdot) \vert \vert_{
\mathcal{L}^{2}_{\mathrm{M}_{2}(\mathbb{C})}(\widetilde{\Sigma}_{p}^{o})}
\norm{w_{+}^{\Sigma_{\mathscr{R}}^{o}}(\cdot)}_{\mathcal{L}^{2}_{\mathrm{M}_{2}
(\mathbb{C})}(\widetilde{\Sigma}_{p}^{o})}}{2 \pi \operatorname{dist}(z,
\widetilde{\Sigma}_{p}^{o})} \\
&\leqslant \, \norm{(\id \! - \! (\id \! - \! C^{0}_{w^{\Sigma_{
\scriptscriptstyle \blacksquare}^{o}}})^{-1}(\id \! - \! C^{0}_{w^{\Sigma_{
\circlearrowright}^{o}}})^{-1}C^{0}_{w^{\Sigma_{\circlearrowright}^{o}}}C^{
0}_{w^{\Sigma_{\scriptscriptstyle \blacksquare}^{o}}})^{-1}}_{\mathscr{N}_{2}
(\widetilde{\Sigma}_{p}^{o})} \norm{(\id \! - \! C^{0}_{w^{\Sigma_{
\scriptscriptstyle \blacksquare}^{o}}})^{-1}}_{\mathscr{N}_{2}(\widetilde{
\Sigma}_{p}^{o})} \\
&\times \, \dfrac{\norm{(\id \! - \! C^{0}_{w^{\Sigma_{\circlearrowright}^{
o}}})^{-1}}_{\mathscr{N}_{2}(\widetilde{\Sigma}_{p}^{o})} \norm{C^{0}_{w^{
\Sigma_{\scriptscriptstyle \blacksquare}^{o}}}}_{\mathscr{N}_{2}(\widetilde{
\Sigma}_{p}^{o})} \norm{(C^{0}_{w^{\Sigma_{\circlearrowright}^{o}}} \mathrm{I})
(\cdot)}_{\mathcal{L}^{2}_{\mathrm{M}_{2}(\mathbb{C})}(\widetilde{\Sigma}_{p}^{
o})} \norm{w_{+}^{\Sigma_{\mathscr{R}}^{o}}(\cdot)}_{\mathcal{L}^{2}_{\mathrm{
M}_{2}(\mathbb{C})}(\widetilde{\Sigma}_{p}^{o})}}{2 \pi \operatorname{dist}(z,
\widetilde{\Sigma}_{p}^{o})} \\
&\leqslant \, \norm{(\id \! - \! (\id \! - \! C^{0}_{w^{\Sigma_{
\scriptscriptstyle \blacksquare}^{o}}})^{-1}(\id \! - \! C^{0}_{w^{\Sigma_{
\circlearrowright}^{o}}})^{-1}C^{0}_{w^{\Sigma_{\circlearrowright}^{o}}}C^{
0}_{w^{\Sigma_{\scriptscriptstyle \blacksquare}^{o}}})^{-1}}_{\mathscr{N}_{2}
(\widetilde{\Sigma}_{p}^{o})} \norm{(\id \! - \! C^{0}_{w^{\Sigma_{
\scriptscriptstyle \blacksquare}^{o}}})^{-1}}_{\mathscr{N}_{2}(\widetilde{
\Sigma}_{p}^{o})} \\
&\times \, \operatorname{const.} \, \norm{(\id \! - \! C^{0}_{w^{\Sigma_{
\circlearrowright}^{o}}})^{-1}}_{\mathscr{N}_{2}(\widetilde{\Sigma}_{p}^{o})}
\norm{C^{0}_{w^{\Sigma_{\scriptscriptstyle \blacksquare}^{o}}}}_{\mathscr{N}_{
2}(\widetilde{\Sigma}_{p}^{o})} \norm{w_{+}^{\Sigma_{\circlearrowright}^{o}}
(\cdot)}_{\mathcal{L}^{2}_{\mathrm{M}_{2}(\mathbb{C})}(\Sigma_{
\circlearrowright}^{o})} \\
&\times \, \dfrac{(\norm{w_{+}^{\Sigma_{\circlearrowright}^{o}}(\cdot)}_{
\mathcal{L}^{2}_{\mathrm{M}_{2}(\mathbb{C})}(\Sigma^{o}_{\circlearrowright})}
\! + \! \norm{w_{+}^{\Sigma_{\scriptscriptstyle \blacksquare}^{o}}(\cdot)}_{
\mathcal{L}^{2}_{\mathrm{M}_{2}(\mathbb{C})}(\Sigma^{o}_{\scriptscriptstyle
\blacksquare})})}{2 \pi \operatorname{dist}(z,\widetilde{\Sigma}_{p}^{o})};
\end{align*}
using the fact that (cf. Proposition~5.2) $\norm{(\id \! - \! (\id \! - \! C^{
0}_{w^{\Sigma_{\scriptscriptstyle \blacksquare}^{o}}})^{-1}(\id \! - \! C^{
0}_{w^{\Sigma_{\circlearrowright}^{o}}})^{-1}C^{0}_{w^{\Sigma_{
\circlearrowright}^{o}}}C^{0}_{w^{\Sigma_{\scriptscriptstyle \blacksquare}^{
o}}})^{-1}}_{\mathscr{N}_{2}(\widetilde{\Sigma}_{p}^{o})} \! =_{n \to \infty}
\! \mathcal{O}(1)$ (via a Neuman series inversion argument, since $\norm{C^{
0}_{w^{\Sigma_{\scriptscriptstyle \blacksquare}^{o}}}}_{\mathscr{N}_{2}
(\widetilde{\Sigma}_{p}^{o})} \! =_{n \to \infty} \! \mathcal{O}((n \! + \!
\tfrac{1}{2})^{-1/2}f(n) \me^{-(n+\frac{1}{2})c})$ and $\norm{C^{0}_{w^{
\Sigma_{\circlearrowright}^{o}}}}_{\mathscr{N}_{2}(\widetilde{\Sigma}_{p}^{o}
)} \! =_{n \to \infty} \! \mathcal{O}((n \! + \! 1/2)^{-1}f(n)))$, one gets
that
\begin{align*}
\vert I_{6}^{o,\mathcal{C}} \vert \underset{n \to \infty}{\leqslant}& \,
\mathcal{O} \! \left(\dfrac{f(n) \me^{-(n+\frac{1}{2})c}}{\sqrt{n \! + \!
\frac{1}{2}} \operatorname{dist}(z,\widetilde{\Sigma}_{p}^{o})} \right) \!
\mathcal{O} \! \left(\dfrac{f(n)}{n \! + \! \frac{1}{2}} \right) \! \left(\!
\mathcal{O} \! \left(\dfrac{f(n)}{n \! + \! \frac{1}{2}} \right) \! + \!
\mathcal{O} \! \left(\dfrac{f(n) \me^{-(n+\frac{1}{2})c}}{\sqrt{n \! + \!
\frac{1}{2}}} \right) \right) \\
\underset{n \to \infty}{\leqslant}& \, \mathcal{O} \! \left(\dfrac{f(n) \me^{-
(n+\frac{1}{2})c}}{(n \! + \! \frac{1}{2})^{2} \operatorname{dist}(z,
\widetilde{\Sigma}_{p}^{o})} \right).
\end{align*}
One estimates, succinctly, the Cauchy part of $I_{7}^{o}$, denoted $I_{7}^{o,
\mathcal{C}}$, as follows:
\begin{align*}
\vert I_{7}^{o,\mathcal{C}} \vert &\leqslant \int_{\widetilde{\Sigma}_{p}^{o}}
\vert ((\id \! - \! (\id \! - \! C^{0}_{w^{\Sigma_{\circlearrowright}^{o}}})^{
-1}(\id \! - \! C^{0}_{w^{\Sigma_{\scriptscriptstyle \blacksquare}^{o}}})^{-1}
C^{0}_{w^{\Sigma_{\scriptscriptstyle \blacksquare}^{o}}}C^{0}_{w^{\Sigma_{
\circlearrowright}^{o}}})^{-1}(\id \! - \! C^{0}_{w^{\Sigma_{
\circlearrowright}^{o}}})^{-1}(\id \! - \! C^{0}_{w^{\Sigma_{
\scriptscriptstyle \blacksquare}^{o}}})^{-1} \\
&\times \, \dfrac{C^{0}_{w^{\Sigma_{\circlearrowright}^{o}}}(\id \! - \! C^{
0}_{w^{\Sigma_{\scriptscriptstyle \blacksquare}^{o}}})^{-1}C^{0}_{w^{\Sigma_{
\scriptscriptstyle \blacksquare}^{o}}}(C^{0}_{w^{\Sigma_{\scriptscriptstyle
\blacksquare}^{o}}} \mathrm{I}))(s) \vert \vert w_{+}^{\Sigma_{\mathscr{R}}^{
o}}(s) \vert}{\vert s \! - \! z \vert} \, \dfrac{\vert \md s \vert}{2 \pi} \\
&\leqslant \, \dfrac{\norm{(\id \! - \! (\id \! - \! C^{0}_{w^{\Sigma_{
\circlearrowright}^{o}}})^{-1}(\id \! - \! C^{0}_{w^{\Sigma_{
\scriptscriptstyle \blacksquare}^{o}}})^{-1}C^{0}_{w^{\Sigma_{
\scriptscriptstyle \blacksquare}^{o}}}C^{0}_{w^{\Sigma_{\circlearrowright}^{
o}}})^{-1}}_{\mathscr{N}_{2}(\widetilde{\Sigma}_{p}^{o})} \norm{(\id \! - \!
C^{0}_{w^{\Sigma_{\circlearrowright}^{o}}})^{-1}}_{\mathscr{N}_{2}(\widetilde{
\Sigma}_{p}^{o})}}{2 \pi \operatorname{dist}(z,\widetilde{\Sigma}_{p}^{o})} \\
&\times \, \norm{(\id \! - \! C^{0}_{w^{\Sigma_{\scriptscriptstyle
\blacksquare}^{o}}})^{-1}}_{\mathscr{N}_{2}(\widetilde{\Sigma}_{p}^{o})}
\norm{C^{0}_{w^{\Sigma_{\circlearrowright}^{o}}}}_{\mathscr{N}_{2}(\widetilde{
\Sigma}_{p}^{o})} \norm{(\id \! - \! C^{0}_{w^{\Sigma_{\scriptscriptstyle
\blacksquare}^{o}}})^{-1}}_{\mathscr{N}_{2}(\widetilde{\Sigma}_{p}^{o})}
\norm{C^{0}_{w^{\Sigma_{\scriptscriptstyle \blacksquare}^{o}}}}_{\mathscr{N}_{
2}(\widetilde{\Sigma}_{p}^{o})} \\
&\times \, \operatorname{const.} \, \norm{w_{+}^{\Sigma_{\scriptscriptstyle
\blacksquare}^{o}}(\cdot)}_{\mathcal{L}^{2}_{\mathrm{M}_{2}(\mathbb{C})}
(\Sigma^{o}_{\scriptscriptstyle \blacksquare})} \! \left(\norm{w_{+}^{\Sigma_{
\scriptscriptstyle \blacksquare}^{o}}(\cdot)}_{\mathcal{L}^{2}_{\mathrm{M}_{2}
(\mathbb{C})}(\Sigma^{o}_{\scriptscriptstyle \blacksquare})} \! + \! \norm{w_{
+}^{\Sigma_{\circlearrowright}^{o}}(\cdot)}_{\mathcal{L}^{2}_{\mathrm{M}_{2}
(\mathbb{C})}(\Sigma^{o}_{\circlearrowright})} \right);
\end{align*}
using the fact that (established above) $\norm{(\id \! - \! (\id \! - \! C^{
0}_{w^{\Sigma_{\circlearrowright}^{o}}})^{-1}(\id \! - \! C^{0}_{w^{\Sigma_{
\scriptscriptstyle \blacksquare}^{o}}})^{-1}C^{0}_{w^{\Sigma_{
\scriptscriptstyle \blacksquare}^{o}}}C^{0}_{w^{\Sigma_{\circlearrowright}^{
o}}})^{-1}}_{\mathscr{N}_{2}(\widetilde{\Sigma}_{p}^{o})} \! =_{n \to \infty}
\! \mathcal{O}(1)$, one gets that
\begin{align*}
\vert I_{7}^{o,\mathcal{C}} \vert \underset{n \to \infty}{\leqslant}& \,
\mathcal{O} \! \left(\dfrac{f(n)}{(n \! + \! \frac{1}{2}) \operatorname{dist}
(z,\widetilde{\Sigma}_{p}^{o})} \right) \! \mathcal{O} \! \left(\dfrac{f(n)
\me^{-(n+\frac{1}{2})c}}{\sqrt{n \! + \! \frac{1}{2}}} \right) \! \mathcal{O}
\! \left(\dfrac{f(n) \me^{-(n+\frac{1}{2})c}}{\sqrt{n \! + \! \frac{1}{2}}}
\right) \! \left(\! \mathcal{O} \! \left(\dfrac{f(n)}{n \! + \! \frac{1}{2}}
\right) \right. \\
+&\left. \, \mathcal{O} \! \left(\dfrac{f(n) \me^{-(n+\frac{1}{2})c}}{\sqrt{n
\! + \! \frac{1}{2}}} \right) \right) \! \underset{n \to \infty}{\leqslant}
\mathcal{O} \! \left(\dfrac{f(n) \me^{-(n+\frac{1}{2})c}}{(n \! + \! \frac{
1}{2})^{3} \operatorname{dist}(z,\widetilde{\Sigma}_{p}^{o})} \right).
\end{align*}
One estimates, succinctly, the Cauchy part of $I_{8}^{o}$, denoted $I_{8}^{o,
\mathcal{C}}$, as follows:
\begin{align*}
\vert I_{8}^{o,\mathcal{C}} \vert &\leqslant \int_{\widetilde{\Sigma}_{p}^{o}}
\vert ((\id \! - \! (\id \! - \! C^{0}_{w^{\Sigma_{\scriptscriptstyle
\blacksquare}^{o}}})^{-1}(\id \! - \! C^{0}_{w^{\Sigma_{\circlearrowright}^{
o}}})^{-1}C^{0}_{w^{\Sigma_{\circlearrowright}^{o}}}C^{0}_{w^{\Sigma_{
\scriptscriptstyle \blacksquare}^{o}}})^{-1}(\id \! - \! C^{0}_{w^{\Sigma_{
\scriptscriptstyle \blacksquare}^{o}}})^{-1}(\id \! - \! C^{0}_{w^{\Sigma_{
\circlearrowright}^{o}}})^{-1} \\
&\times \, \dfrac{C^{0}_{w^{\Sigma_{\scriptscriptstyle \blacksquare}^{o}}}
(\id \! - \! C^{0}_{w^{\Sigma_{\circlearrowright}^{o}}})^{-1}C^{0}_{w^{
\Sigma_{\circlearrowright}^{o}}}(C^{0}_{w^{\Sigma_{\circlearrowright}^{o}}}
\mathrm{I}))(s) \vert \vert w_{+}^{\Sigma_{\mathscr{R}}^{o}}(s) \vert}{\vert
s \! - \! z \vert} \, \dfrac{\vert \md s \vert}{2 \pi} \\
&\leqslant \, \dfrac{\norm{(\id \! - \! (\id \! - \! C^{0}_{w^{\Sigma_{
\scriptscriptstyle \blacksquare}^{o}}})^{-1}(\id \! - \! C^{0}_{w^{\Sigma_{
\circlearrowright}^{o}}})^{-1}C^{0}_{w^{\Sigma_{\circlearrowright}^{o}}}C^{
0}_{w^{\Sigma_{\scriptscriptstyle \blacksquare}^{o}}})^{-1}}_{\mathscr{N}_{2}
(\widetilde{\Sigma}_{p}^{o})} \norm{(\id \! - \! C^{0}_{w^{\Sigma_{
\scriptscriptstyle \blacksquare}^{o}}})^{-1}}_{\mathscr{N}_{2}(\widetilde{
\Sigma}_{p}^{o})}}{2 \pi \operatorname{dist}(z,\widetilde{\Sigma}_{p}^{o})} \\
&\times \, \norm{(\id \! - \! C^{0}_{w^{\Sigma_{\circlearrowright}^{o}}})^{-
1}}_{\mathscr{N}_{2}(\widetilde{\Sigma}_{p}^{o})} \norm{C^{0}_{w^{\Sigma_{
\scriptscriptstyle \blacksquare}^{o}}}}_{\mathscr{N}_{2}(\widetilde{\Sigma}_{
p}^{o})} \norm{(\id \! - \! C^{0}_{w^{\Sigma_{\circlearrowright}^{o}}})^{-1}
}_{\mathscr{N}_{2}(\widetilde{\Sigma}_{p}^{o})} \norm{C^{0}_{w^{\Sigma_{
\circlearrowright}^{o}}}}_{\mathscr{N}_{2}(\widetilde{\Sigma}_{p}^{o})} \\
&\times \, \operatorname{const.} \, \norm{w_{+}^{\Sigma_{\circlearrowright}^{
o}}(\cdot)}_{\mathcal{L}^{2}_{\mathrm{M}_{2}(\mathbb{C})}(\Sigma^{o}_{
\circlearrowright})} \! \left(\norm{w_{+}^{\Sigma_{\scriptscriptstyle
\blacksquare}^{o}}(\cdot)}_{\mathcal{L}^{2}_{\mathrm{M}_{2}(\mathbb{C})}
(\Sigma^{o}_{\scriptscriptstyle \blacksquare})} \! + \! \norm{w_{+}^{\Sigma_{
\circlearrowright}^{o}}(\cdot)}_{\mathcal{L}^{2}_{\mathrm{M}_{2}(\mathbb{C})}
(\Sigma^{o}_{\circlearrowright})} \right);
\end{align*}
using the fact that (established above) $\norm{(\id \! - \! (\id \! - \! C^{
0}_{w^{\Sigma_{\scriptscriptstyle \blacksquare}^{o}}})^{-1}(\id \! - \! C^{
0}_{w^{\Sigma_{\circlearrowright}^{o}}})^{-1}C^{0}_{w^{\Sigma_{
\circlearrowright}^{o}}}C^{0}_{w^{\Sigma_{\scriptscriptstyle \blacksquare}^{
o}}})^{-1}}_{\mathscr{N}_{2}(\widetilde{\Sigma}_{p}^{o})} \! =_{n \to \infty}
\! \mathcal{O}(1)$, one gets that
\begin{align*}
\vert I_{8}^{o,\mathcal{C}} \vert &\underset{n \to \infty}{\leqslant} \mathcal{
O} \! \left(\dfrac{f(n) \me^{-(n+\frac{1}{2})c}}{\sqrt{n \! + \! \frac{1}{2}}
\, \operatorname{dist}(z,\widetilde{\Sigma}_{p}^{o})} \right) \! \mathcal{O}
\! \left(\dfrac{f(n)}{n \! + \! \frac{1}{2}} \right) \! \mathcal{O} \! \left(
\dfrac{f(n)}{n \! + \! \frac{1}{2}} \right) \! \left(\! \mathcal{O} \! \left(
\dfrac{f(n)}{n \! + \! \frac{1}{2}} \right) \! + \! \mathcal{O} \! \left(
\dfrac{f(n) \me^{-(n+\frac{1}{2})c}}{\sqrt{n \! + \! \frac{1}{2}}} \right)
\right) \\
&\underset{n \to \infty}{\leqslant} \mathcal{O} \! \left(\dfrac{f(n) \me^{-
(n+\frac{1}{2})c}}{(n \! + \! \frac{1}{2})^{3} \operatorname{dist}(z,
\widetilde{\Sigma}_{p}^{o})} \right).
\end{align*}
Analogously, estimating (as $n \! \to \! \infty)$ the non-Cauchy contributions 
(corresponding to the kernel $\tfrac{1}{s})$ of the terms on the right-hand 
side of Equation~(5.2) which, too, are $\mathcal{O}((n \! + \! 1/2)^{-2})$, 
and gathering all derived (upper) bounds, one arrives at the result stated in 
the Lemma. \hfill $\qed$
\begin{ccccc}
Let $\mathscr{R}^{o} \colon \mathbb{C} \setminus \widetilde{\Sigma}_{p}^{o} 
\! \to \! \operatorname{SL}_{2}(\mathbb{C})$ be the solution of the {\rm RHP} 
$(\mathscr{R}^{o}(z),\upsilon_{\mathscr{R}}^{o}(z),\widetilde{\Sigma}_{p}^{o}
)$ formulated in Proposition~{\rm 5.2} with the $n \! \to \! \infty$ integral 
representation given in Lemma~{\rm 5.2}. Then, uniformly for compact subsets 
of $\mathbb{C} \setminus \widetilde{\Sigma}_{p}^{o} \! \ni \! z$,
\begin{equation*}
\mathscr{R}^{o}(z) \underset{\underset{z \in \mathbb{C} \setminus \widetilde{
\Sigma}_{p}^{o}}{n \to \infty}}{=} \mathrm{I} \! + \! \dfrac{1}{(n \! + \!
\frac{1}{2})} \! \left(\mathscr{R}^{o}_{0}(z) \! - \! \widetilde{\mathscr{R}
}^{o}_{0}(z) \right) \! + \! \mathcal{O} \! \left(\dfrac{f(z;n)}{(n \! + \!
\frac{1}{2})^{2}} \right),
\end{equation*}
where $\mathscr{R}^{o}_{0}(z)$ is defined in Theorem~{\rm 2.3.1}, 
Equations~{\rm (2.23)--(2.57)}, $\widetilde{\mathscr{R}}^{o}_{0}(z)$ 
is defined in Theorem~{\rm 2.3.1}, Equations~{\rm (2.14)--(2.20)} 
and~{\rm (2.70)--(2.74)}, and $f(z;n)$, where the $n$-dependence arises 
due to the $n$-dependence of the associated Riemann theta functions, is a 
bounded (with respect to $z$ and $n)$, 
$\operatorname{GL}_{2}(\mathbb{C})$-valued function which is analytic (with 
respect to $z)$ for $z \! \in \! \mathbb{C} \setminus \widetilde{\Sigma}_{
p}^{o}$, and $(f(\pmb{\cdot};n))_{kl} \! =_{n \to \infty} \! \mathcal{O}
(1)$, $k,l \! = \! 1,2$.
\end{ccccc}
\begin{eeeee}
Note {}from the formulation of Lemma~5.3 above that (cf. Theorem~2.3.1, 
Equations (2.24)--(2.27)), for $j \! = \! 1,\dotsc,N \! + \! 1$, 
$\operatorname{tr}(\mathscr{A}^{o}(a_{j}^{o})) \! = \! \operatorname{tr}
(\mathscr{A}^{o}(b_{j-1}^{o})) \! = \! \operatorname{tr}(\mathscr{B}^{o}
(a_{j}^{o})) \! = \! \operatorname{tr}(\mathscr{B}^{o}(b_{j-1}^{o})) \! = \! 
0$. \hfill $\blacksquare$
\end{eeeee}

\emph{Proof.} Recall the integral representation for $\mathscr{R}^{o} \colon
\mathbb{C} \setminus \widetilde{\Sigma}_{p}^{o} \! \to \! \operatorname{SL}_{
2}(\mathbb{C})$ given in Lemma~5.2:
\begin{equation*}
\mathscr{R}^{o}(z) \underset{n \to \infty}{=} \mathrm{I} \! + \! \int_{
\Sigma^{o}_{\circlearrowright}} \dfrac{zw_{+}^{\Sigma^{o}_{\circlearrowright}}
(s)}{s(s \! - \! z)} \, \dfrac{\md s}{2 \pi \mi} \! + \! \mathcal{O} \! \left(
\! \dfrac{f(n)}{(n \! + \! \frac{1}{2})^{2} \min \{1,\operatorname{dist}(z,
\widetilde{\Sigma}_{p}^{o})\}} \right), \quad z \! \in \! \mathbb{C} \setminus
\widetilde{\Sigma}_{p}^{o},
\end{equation*}
where $\Sigma^{o}_{\circlearrowright} \! := \! \cup_{j=1}^{N+1}(\partial 
\mathbb{U}^{o}_{\delta_{b_{j-1}}} \cup \partial \mathbb{U}^{o}_{\delta_{a_{
j}}})$, and $(f(n))_{kl} \! =_{n \to \infty} \! \mathcal{O}(1)$, $k,l \! = 
\! 1,2$. Recalling that the radii of the open discs $\mathbb{U}^{o}_{\delta_{
b_{j-1}}},\mathbb{U}^{o}_{\delta_{a_{j}}}$, $j \! = \! 1,\dotsc,N \! + \! 
1$, are chosen, amongst other factors (cf. Lemmae~4.6 and~4.7), such that 
$\mathbb{U}^{o}_{\delta_{b_{j-1}}} \cap \mathbb{U}^{o}_{\delta_{a_{k}}} \!
= \! \varnothing$, $j,k \! = \! 1,\dotsc,N \! + \! 1$, it follows {}from 
the above integral representation that
\begin{equation*}
\mathscr{R}^{o}(z) \! \underset{n \to \infty}{=} \! \mathrm{I} \! - \! \sum_{
j=1}^{N+1} \! \left(\! \oint_{\partial \mathbb{U}^{o}_{\delta_{b_{j-1}}}}+
\oint_{\partial \mathbb{U}^{o}_{\delta_{a_{j}}}} \right) \! \dfrac{zw_{+}^{
\Sigma^{o}_{\circlearrowright}}(s)}{s(s \! - \! z)} \, \dfrac{\md s}{2 \pi
\mi} \! + \! \mathcal{O} \! \left(\! \dfrac{f(n)}{(n \! + \! \frac{1}{2})^{2}
\min \{1,\operatorname{dist}(z,\widetilde{\Sigma}_{p}^{o})\}} \right), \quad z
\! \in \! \mathbb{C} \setminus \widetilde{\Sigma}_{p}^{o},
\end{equation*}
where $\oint_{\partial \mathbb{U}^{o}_{\delta_{b_{j-1}}}},\oint_{\partial 
\mathbb{U}^{o}_{\delta_{a_{j}}}}$, $j \! = \! 1,\dotsc,N \! + \! 1$, are 
counter-clockwise-oriented, closed (contour) integrals (Figure~10) about 
the end-points of the support of the `odd' equilibrium measure, $\lbrace 
b_{j-1}^{o},a_{j}^{o} \rbrace_{j=1}^{N+1}$. Noting the partial fraction 
decomposition $\tfrac{z}{s(s-z)} \! = \! -\tfrac{1}{s} \! + \! \tfrac{1}{s
-z}$, the evaluation of these $4(N \! + \! 1)$ contour integrals requires 
the application of the Cauchy and Residue Theorems; and, since the evaluation 
of the respective integrals entails analogous calculations, consider, say, 
and without loss of generality, the evaluation of the integrals corresponding 
to the (standard) Cauchy kernel, $\tfrac{1}{s-z} \, \tfrac{\md s}{2 \pi \mi}$, 
about the right-most end-points $a_{j}^{o}$, $j \! = \! 1,\dotsc,N$, namely:
\begin{equation*}
\oint_{\partial \mathbb{U}^{o}_{\delta_{a_{j}}}} \dfrac{w_{+}^{\Sigma^{o}_{
\circlearrowright}}(s)}{s \! - \! z} \, \dfrac{\md s}{2 \pi \mi}, \quad j \!
= \! 1,\dotsc,N.
\end{equation*}
Recalling {}from Lemma~4.7 that $\xi_{a_{j}}^{o}(z) \! = \! (z \! - \! a_{j}^{
o})^{3/2}G_{a_{j}}^{o}(z)$, $z \! \in \! \mathbb{U}^{o}_{\delta_{a_{j}}}
\setminus (-\infty,a_{j}^{o})$, $j \! = \! 1,\dotsc,N$, it follows {}from
item~(5) of Proposition~5.1 that, since $w_{+}^{\Sigma^{o}_{\circlearrowright}
}(z) \! = \! \upsilon_{\mathscr{R}}^{o}(z) \! - \! \mathrm{I}$, for $j \! =
\! 1,\dotsc,N$,
\begin{align*}
w_{+}^{\Sigma^{o}_{\circlearrowright}}(z) \underset{\underset{z \in \mathbb{
C}_{\pm} \cap \partial \mathbb{U}_{\delta_{a_{j}}}^{o}}{n \to \infty}}{=}& \,
\dfrac{1}{(n \! + \! \frac{1}{2})(z \! - \! a_{j}^{o})^{3/2}G_{a_{j}}^{o}(z)}
\overset{o}{\mathfrak{M}}^{\raise-1.0ex\hbox{$\scriptstyle \infty$}}(z) \!
\begin{pmatrix}
\mp (s_{1}+t_{1}) & \pm \mi (s_{1}-t_{1}) \me^{\mi (n+\frac{1}{2}) \Omega_{
j}^{o}} \\
\pm \mi (s_{1}-t_{1}) \me^{-\mi (n+\frac{1}{2}) \Omega_{j}^{o}} & \pm (s_{
1}+t_{1})
\end{pmatrix} \\
\times& \,
(\overset{o}{\mathfrak{M}}^{\raise-1.0ex\hbox{$\scriptstyle \infty$}}(z))^{-
1} \! + \! \mathcal{O} \! \left(\dfrac{1}{(n \! + \! \frac{1}{2})^{2}(z \! -
\! a_{j}^{o})^{3}(G_{a_{j}}^{o}(z))^{2}}
\overset{o}{\mathfrak{M}}^{\raise-1.0ex\hbox{$\scriptstyle \infty$}}(z)f_{a_{
j}}^{o}(n)(\overset{o}{\mathfrak{M}}^{\raise-1.0ex\hbox{$\scriptstyle \infty$}}
(z))^{-1} \right),
\end{align*}
where $\overset{o}{\mathfrak{M}}^{\raise-1.0ex\hbox{$\scriptstyle \infty$}}
(z)$ and $\Omega_{j}^{o}$ are defined in Lemma~4.5, and $(f_{a_{j}}^{o}
(n))_{kl} \! =_{n \to \infty} \! \mathcal{O}(1)$, $k,l \! = \! 1,2$. A
matrix-multiplication argument shows that
$\overset{o}{\mathfrak{M}}^{\raise-1.0ex\hbox{$\scriptstyle \infty$}}(z) \!
\left(
\begin{smallmatrix}
\mp (s_{1}+t_{1}) & \pm \mi (s_{1}-t_{1}) \me^{\mi (n+\frac{1}{2}) \Omega_{
j}^{o}} \\
\pm \mi (s_{1}-t_{1}) \me^{-\mi (n+\frac{1}{2}) \Omega_{j}^{o}} & \pm (s_{1}
+t_{1})
\end{smallmatrix}
\right) \!
(\overset{o}{\mathfrak{M}}^{\raise-1.0ex\hbox{$\scriptstyle \infty$}}(z))^{-
1}$ is given by
\begin{equation*}
\begin{pmatrix}
\boxed{\begin{matrix}
\mp \frac{1}{4}(s_{1} \! + \! t_{1}) \! \left(\frac{(\gamma^{o}(z))^{2}+
(\gamma^{o}(0))^{2}}{\gamma^{o}(0) \gamma^{o}(z)} \right)^{2} \! \mathfrak{m}^{
o}_{11}(z) \mathfrak{m}^{o}_{22}(z) \\
\mp \frac{1}{4}(s_{1} \! + \! t_{1}) \! \left(\frac{(\gamma^{o}(z))^{2}-
(\gamma^{o}(0))^{2}}{\gamma^{o}(0) \gamma^{o}(z)} \right)^{2} \! \mathfrak{m}^{
o}_{12}(z) \mathfrak{m}^{o}_{21}(z) \\
\mp \frac{1}{4}(s_{1} \! - \! t_{1}) \! \left(\frac{(\gamma^{o}(z))^{4}-
(\gamma^{o}(0))^{4}}{(\gamma^{o}(0) \gamma^{o}(z))^{2}} \right) \! \frac{
\mathfrak{m}^{o}_{11}(z) \mathfrak{m}^{o}_{21}(z)}{\me^{-\mi (n+\frac{1}{2})
\Omega_{j}^{o}}} \\
\mp \frac{1}{4}(s_{1} \! - \! t_{1}) \! \left(\frac{(\gamma^{o}(z))^{4}-
(\gamma^{o}(0))^{4}}{(\gamma^{o}(0) \gamma^{o}(z))^{2}} \right) \! \frac{
\mathfrak{m}^{o}_{12}(z) \mathfrak{m}^{o}_{22}(z)}{\me^{\mi (n+\frac{1}{2})
\Omega_{j}^{o}}}
\end{matrix}} &
\boxed{\begin{matrix}
\pm \frac{\mi}{2}(s_{1} \! + \! t_{1}) \! \left(\frac{(\gamma^{o}(z))^{4}-
(\gamma^{o}(0))^{4}}{(\gamma^{o}(0) \gamma^{o}(z))^{2}} \right) \! \mathfrak{
m}^{o}_{11}(z) \mathfrak{m}^{o}_{12}(z) \\
\pm \frac{\mi}{4}(s_{1} \! - \! t_{1}) \! \left(\frac{(\gamma^{o}(z))^{2}+
(\gamma^{o}(0))^{2}}{\gamma^{o}(0) \gamma^{o}(z)} \right)^{2} \! \frac{
(\mathfrak{m}^{o}_{11}(z))^{2}}{\me^{-\mi (n+\frac{1}{2}) \Omega_{j}^{o}}} \\
\pm \frac{\mi}{4}(s_{1} \! - \! t_{1}) \! \left(\frac{(\gamma^{o}(z))^{2}-
(\gamma^{o}(0))^{2}}{\gamma^{o}(0) \gamma^{o}(z)} \right)^{2} \! \frac{
(\mathfrak{m}^{o}_{12}(z))^{2}}{\me^{\mi (n+\frac{1}{2}) \Omega_{j}^{o}}}
\end{matrix}} \\
\boxed{\begin{matrix}
\pm \frac{\mi}{2}(s_{1} \! + \! t_{1}) \! \left(\frac{(\gamma^{o}(z))^{4}-
(\gamma^{o}(0))^{4}}{(\gamma^{o}(0) \gamma^{o}(z))^{2}} \right) \! \mathfrak{
m}^{o}_{21}(z) \mathfrak{m}^{o}_{22}(z) \\
\pm \frac{\mi}{4}(s_{1} \! - \! t_{1}) \! \left(\frac{(\gamma^{o}(z))^{2}-
(\gamma^{o}(0))^{2}}{\gamma^{o}(0) \gamma^{o}(z)} \right)^{2} \! \frac{
(\mathfrak{m}^{o}_{21}(z))^{2}}{\me^{-\mi (n+\frac{1}{2}) \Omega_{j}^{o}}} \\
\pm \frac{\mi}{4}(s_{1} \! - \! t_{1}) \! \left(\frac{(\gamma^{o}(z))^{2}+
(\gamma^{o}(0))^{2}}{\gamma^{o}(0) \gamma^{o}(z)} \right)^{2} \! \frac{
(\mathfrak{m}^{o}_{22}(z))^{2}}{\me^{\mi (n+\frac{1}{2}) \Omega_{j}^{o}}}
\end{matrix}} &
\boxed{\begin{matrix}
\pm \frac{1}{4}(s_{1} \! + \! t_{1}) \! \left(\frac{(\gamma^{o}(z))^{2}+
(\gamma^{o}(0))^{2}}{\gamma^{o}(0) \gamma^{o}(z)} \right)^{2} \! \mathfrak{m}^{
o}_{11}(z) \mathfrak{m}^{o}_{22}(z) \\
\pm \frac{1}{4}(s_{1} \! + \! t_{1}) \! \left(\frac{(\gamma^{o}(z))^{2}-
(\gamma^{o}(0))^{2}}{\gamma^{o}(0) \gamma^{o}(z)} \right)^{2} \! \mathfrak{m}^{
o}_{12}(z) \mathfrak{m}^{o}_{21}(z) \\
\pm \frac{1}{4}(s_{1} \! - \! t_{1}) \! \left(\frac{(\gamma^{o}(z))^{4}-
(\gamma^{o}(0))^{4}}{(\gamma^{o}(0) \gamma^{o}(z))^{2}} \right) \! \frac{
\mathfrak{m}^{o}_{11}(z) \mathfrak{m}^{o}_{21}(z)}{\me^{-\mi (n+\frac{1}{2})
\Omega_{j}^{o}}} \\
\pm \frac{1}{4}(s_{1} \! - \! t_{1}) \! \left(\frac{(\gamma^{o}(z))^{4}-
(\gamma^{o}(0))^{4}}{(\gamma^{o}(0) \gamma^{o}(z))^{2}} \right) \! \frac{
\mathfrak{m}^{o}_{12}(z) \mathfrak{m}^{o}_{22}(z)}{\me^{\mi (n+\frac{1}{2})
\Omega_{j}^{o}}}
\end{matrix}}
\end{pmatrix},
\end{equation*}
where $s_{1}$ and $t_{1}$ are given in Theorem~2.3.1, Equations~(2.28), 
$\gamma^{o}(z)$ and $\gamma^{o}(0)$ are defined in Lemma~4.4, and 
$\mathfrak{m}^{o}_{kl}(z)$, $k,l \! = \! 1,2$, are defined in Theorem~2.3.1, 
Equations~(2.17)--(2.20). Recall that, for $j \! = \! 1,\dotsc,N$, $\omega_{
j}^{o} \! = \! \sum_{k=1}^{N}c_{jk}^{o}(R_{o}(z))^{-1/2}z^{N-k} \, \md z$, 
where $c_{jk}^{o}$, $j,k \! = \! 1,\dotsc,N$, are obtained {}from 
Equations~(O1) and~(O2), and (the multi-valued function) $(R_{o}(z))^{1/2}$ 
is defined in Theorem~2.3.1, Equation~(2.8). One shows that
\begin{equation*}
\omega_{m}^{o} \underset{\underset{j=1,\dotsc,N}{z \to a_{j}^{o}}}{=} \dfrac{
(\mathfrak{f}_{o}(a_{j}^{o}))^{-1}}{\sqrt{\smash[b]{z \! - \! a_{j}^{o}}}}
\! \left(\mathfrak{p}_{m}^{\natural}(a_{j}^{o}) \! + \! \mathfrak{q}_{m}^{
\natural}(a_{j}^{o})(z \! - \! a_{j}^{o}) \! + \! \mathfrak{r}_{m}^{\natural}
(a_{j}^{o})(z \! - \! a_{j}^{o})^{2} \! + \! \mathcal{O}((z \! - \! a_{j}^{
o})^{3}) \right) \md z, \quad m \! = \! 1,\dotsc,N,
\end{equation*}
where
\begin{gather*}
\mathfrak{f}_{o}(\xi) \! = \! (-1)^{N-j+1} \! \left((a_{N+1}^{o} \! - \! \xi)
(\xi \! - \! b_{0}^{o})(b_{j}^{o} \! - \! \xi) \prod_{k=1}^{j-1}(\xi \! - \! 
b_{k}^{o})(\xi \! - \! a_{k}^{o}) \prod_{l=j+1}^{N}(b_{l}^{o} \! - \! \xi)
(a_{l}^{o} \! - \! \xi) \right)^{1/2}, \\
\mathfrak{p}_{m}^{\natural}(\xi) \! = \! \sum_{k=1}^{N}c_{mk}^{o} \xi^{N-k}, 
\qquad \qquad \quad \mathfrak{q}_{m}^{\natural}(\xi) \! = \! \sum_{k=1}^{N}
c_{mk}^{o} \xi^{N-k-1} \! \left(N \! - \! k \! - \! \dfrac{\xi \mathfrak{
f}_{o}^{\prime}(\xi)}{\mathfrak{f}_{o}(\xi)} \right), \\
\mathfrak{r}_{m}^{\natural}(\xi) \! = \! \sum_{k=1}^{N}c_{mk}^{o} \xi^{N-k-2}
\! \left(\dfrac{(N \! - \! k)(N \! - \! k \! - \! 1)}{2} \! - \! \dfrac{(N 
\! - \! k) \xi \mathfrak{f}_{o}^{\prime}(\xi)}{\mathfrak{f}_{o}(\xi)} \! + \! 
\xi^{2} \! \left(\left(\dfrac{\mathfrak{f}_{o}^{\prime}(\xi)}{\mathfrak{f}_{
o}(\xi)} \right)^{2} \! - \! \dfrac{\mathfrak{f}_{o}^{\prime \prime}(\xi)}{2 
\mathfrak{f}_{o}(\xi)} \right) \right),
\end{gather*}
with $(-1)^{-N+j-1} \mathfrak{f}_{o}(a_{j}^{o}) \! > \! 0$,
\begin{align*}
\mathfrak{f}_{o}^{\prime}(\xi)=& \, \dfrac{1}{2} \mathfrak{f}_{o}(\xi) \! 
\left(\sum_{\substack{k=1\\k \not= j}}^{N} \left(\dfrac{1}{\xi \! - \! b_{
k}^{o}} \! + \! \dfrac{1}{\xi \! - \! a_{k}^{o}} \right) \! + \! \dfrac{1}{
\xi \! - \! b_{j}^{o}} \! + \! \dfrac{1}{\xi \! - \! a_{N+1}^{o}} \! + \! 
\dfrac{1}{\xi \! - \! b_{0}^{o}} \right), \\
\mathfrak{f}_{o}^{\prime \prime}(\xi)=& \, -\dfrac{1}{2} \mathfrak{f}_{o}
(\xi) \! \left(\sum_{\substack{k=1\\k \not= j}}^{N} \left(\dfrac{1}{(\xi \! 
- \! b_{k}^{o})^{2}} \! + \! \dfrac{1}{(\xi \! - \! a_{k}^{o})^{2}} \right) 
\! + \! \dfrac{1}{(\xi \! - \! b_{j}^{o})^{2}} \! + \! \dfrac{1}{(\xi \! - \! 
a_{N+1}^{o})^{2}} \! + \! \dfrac{1}{(\xi \! - \! b_{0}^{o})^{2}} \right) \\
+& \, \dfrac{1}{4} \mathfrak{f}_{o}(\xi) \! \left(\sum_{\substack{k=1\\k\not= 
j}}^{N} \left(\dfrac{1}{\xi \! - \! b_{k}^{o}} \! + \! \dfrac{1}{\xi \! - \! 
a_{k}^{o}} \right) \! + \! \dfrac{1}{\xi \! - \! b_{j}^{o}} \! + \! \dfrac{
1}{\xi \! - \! a_{N+1}^{o}} \! + \! \dfrac{1}{\xi \! - \! b_{0}^{o}} 
\right)^{2}.
\end{align*}
Recall (cf. Lemma~4.5), also, that $\bm{u}^{o} \! \equiv \! \int_{a_{N+1}^{o}
}^{z} \bm{\omega}^{o}$ $(\in \! \operatorname{Jac}(\mathcal{Y}_{o}))$, where 
$\equiv$ denotes congruence modulo the period lattice, with $\bm{\omega}^{o} 
\! := \! (\omega_{1}^{o},\omega_{2}^{o},\dotsc,\omega_{N}^{o})$; hence, via 
the above expansion (as $z \! \to \! a_{j}^{o}$, $j \! = \! 1,\dotsc,N)$ for 
$\omega_{m}^{o}$, $m \! = \! 1,\dotsc,N$, one arrives at
\begin{equation*}
\int_{a_{j}^{o}}^{z} \omega_{m}^{o} \underset{\underset{j=1,\dotsc,N}{z
\to a_{j}^{o}}}{\equiv} \dfrac{2 \mathfrak{p}_{m}^{\natural}(a_{j}^{o})}{
\mathfrak{f}_{o}(a_{j}^{o})}(z \! - \! a_{j}^{o})^{1/2} \! + \! \dfrac{2
\mathfrak{q}_{m}^{\natural}(a_{j}^{o})}{3 \mathfrak{f}_{o}(a_{j}^{o})}(z \!
- \! a_{j}^{o})^{3/2} \! + \! \dfrac{2 \mathfrak{r}_{m}^{\natural}(a_{j}^{
o})}{5 \mathfrak{f}_{o}(a_{j}^{o})}(z \! - \! a_{j}^{o})^{5/2} \! + \!
\mathcal{O}((z \! - \! a_{j}^{o})^{7/2}).
\end{equation*}
{}From the definition of $\mathfrak{m}^{o}_{kl}(z)$, $k,l \! = \! 1,2$, given
in Theorem~2.3.1, Equations~(2.17)--(2.20), the definition of the `odd'
Riemann theta function given by Equation~(2.1), and recalling that $\mathfrak{
m}^{o}_{kl}(z)$, $k,l \! = \! 1,2$, satisfy the jump relation (cf. Lemma~4.5)
$\mathfrak{m}^{o}_{+}(z) \! = \! \mathfrak{m}^{o}_{-}(z)(\exp (-\mi (n \! + \!
\tfrac{1}{2}) \Omega_{j}^{o}) \sigma_{-} \! + \! \exp (\mi (n \! + \! \tfrac{
1}{2}) \Omega_{j}^{o}) \sigma_{+})$, via the above asymptotic expansion (as $z
\! \to \! a_{j}^{o}$, $j \! = \! 1,\dotsc,N)$ for $\int_{a_{j}^{o}}^{z}
\omega_{m}^{o}$, $m \! = \! 1,\dotsc,N$, one arrives at
\begin{align*}
\mathfrak{m}^{o}_{11}(z) \underset{\underset{j=1,\dotsc,N}{z \to a_{j}^{o}}}{
=}& \, \varkappa_{1}^{o}(a_{j}^{o}) \! \left(1 \! + \! \mi \aleph^{1}_{1}(a_{
j}^{o})(z \! - \! a_{j}^{o})^{1/2} \! + \! \daleth^{1}_{1}(a_{j}^{o})(z \! -
\! a_{j}^{o}) \! + \! \mi \beth^{1}_{1}(a_{j}^{o})(z \! - \! a_{j}^{o})^{3/2}
\! + \! \gimel^{1}_{1}(a_{j}^{o})(z \! - \! a_{j}^{o})^{2} \right. \\
+&\left. \, \mathcal{O}((z \! - \! a_{j}^{o})^{5/2}) \right), \\
\mathfrak{m}^{o}_{12}(z) \underset{\underset{j=1,\dotsc,N}{z \to a_{j}^{o}}}{
=}& \, \varkappa_{1}^{o}(a_{j}^{o}) \! \left(1 \! - \! \mi \aleph^{-1}_{1}
(a_{j}^{o})(z \! - \! a_{j}^{o})^{1/2} \! + \! \daleth^{-1}_{1}(a_{j}^{o})(z
\! - \! a_{j}^{o}) \! - \! \mi \beth^{-1}_{1}(a_{j}^{o})(z \! - \! a_{j}^{o}
)^{3/2} \! + \! \gimel^{-1}_{1}(a_{j}^{o})(z \! - \! a_{j}^{o})^{2} \right. \\
+&\left. \, \mathcal{O}((z \! - \! a_{j}^{o})^{5/2}) \right) \! \exp \! \left(
\mi \! \left(n \! + \! \dfrac{1}{2} \right) \! \Omega_{j}^{o} \right), \\
\mathfrak{m}^{o}_{21}(z) \underset{\underset{j=1,\dotsc,N}{z \to a_{j}^{o}}}{
=}& \, \varkappa_{2}^{o}(a_{j}^{o}) \! \left(1 \! + \! \mi \aleph^{1}_{-1}(a_{
j}^{o})(z \! - \! a_{j}^{o})^{1/2} \! + \! \daleth^{1}_{-1}(a_{j}^{o})(z \! -
\! a_{j}^{o}) \! + \! \mi \beth^{1}_{-1}(a_{j}^{o})(z \! - \! a_{j}^{o})^{3/2}
\! + \! \gimel^{1}_{-1}(a_{j}^{o})(z \! - \! a_{j}^{o})^{2} \right. \\
+&\left. \, \mathcal{O}((z \! - \! a_{j}^{o})^{5/2}) \right), \\
\mathfrak{m}^{o}_{22}(z) \underset{\underset{j=1,\dotsc,N}{z \to a_{j}^{o}}}{
=}& \, \varkappa_{2}^{o}(a_{j}^{o}) \! \left(1 \! - \! \mi \aleph^{-1}_{-1}
(a_{j}^{o})(z \! - \! a_{j}^{o})^{1/2} \! + \! \daleth^{-1}_{-1}(a_{j}^{o})
(z \! - \! a_{j}^{o}) \! - \! \mi \beth^{-1}_{-1}(a_{j}^{o})(z \! - \! a_{
j}^{o})^{3/2} \! + \! \gimel^{-1}_{-1}(a_{j}^{o})(z \! - \! a_{j}^{o})^{2}
\right. \\
+&\left. \, \mathcal{O}((z \! - \! a_{j}^{o})^{5/2}) \right) \! \exp \! 
\left(\mi \! \left(n \! + \! \dfrac{1}{2} \right) \! \Omega_{j}^{o} \right),
\end{align*}
where, for $\varepsilon_{1},\varepsilon_{2} \! = \! \pm 1$,
\begin{align*}
\varkappa_{1}^{o}(\xi) =& \, \dfrac{1}{\mathbb{E}} \dfrac{\bm{\theta}^{o}
(\bm{u}^{o}_{+}(0) \! + \! \bm{d}_{o}) \bm{\theta}^{o}(\bm{u}^{o}_{+}(\xi) 
\! - \! \frac{1}{2 \pi}(n \! + \! \frac{1}{2}) \bm{\Omega}^{o} \! + \! \bm{
d}_{o})}{\bm{\theta}^{o}(\bm{u}^{o}_{+}(0) \! - \! \frac{1}{2 \pi}(n \! + 
\! \frac{1}{2}) \bm{\Omega}^{o} \! + \! \bm{d}_{o}) \bm{\theta}^{o}(\bm{u}^{
o}_{+}(\xi) \! + \! \bm{d}_{o})}, \\
\varkappa_{2}^{o}(\xi) =& \, \mathbb{E} \dfrac{\bm{\theta}^{o}(-\bm{u}^{o}_{
+}(0) \! - \! \bm{d}_{o}) \bm{\theta}^{o}(\bm{u}^{o}_{+}(\xi) \! - \! \frac{
1}{2 \pi}(n \! + \! \frac{1}{2}) \bm{\Omega}^{o} \! - \! \bm{d}_{o})}{\bm{
\theta}^{o}(-\bm{u}^{o}_{+}(0) \! - \! \frac{1}{2 \pi}(n \! + \! \frac{1}{2}) 
\bm{\Omega}^{o} \! - \! \bm{d}_{o}) \bm{\theta}^{o}(\bm{u}^{o}_{+}(\xi) \! 
- \! \bm{d}_{o})}, \\
\aleph^{\varepsilon_{1}}_{\varepsilon_{2}}(\xi) =& \, -\dfrac{\mathfrak{u}^{
o}(\varepsilon_{1},\varepsilon_{2},\bm{0};\xi)}{\bm{\theta}^{o}(\varepsilon_{
1} \bm{u}^{o}_{+}(\xi) \! + \! \varepsilon_{2} \bm{d}_{o})} \! + \! \dfrac{
\mathfrak{u}^{o}(\varepsilon_{1},\varepsilon_{2},\bm{\Omega}^{o};\xi)}{\bm{
\theta}^{o}(\varepsilon_{1} \bm{u}^{o}_{+}(\xi) \! - \! \frac{1}{2 \pi}(n \!
+ \! \frac{1}{2}) \bm{\Omega}^{o} \! + \! \varepsilon_{2} \bm{d}_{o})}, \\
\daleth^{\varepsilon_{1}}_{\varepsilon_{2}}(\xi) =& \, -\dfrac{\mathfrak{v}^{
o}(\varepsilon_{1},\varepsilon_{2},\bm{0};\xi)}{\bm{\theta}^{o}(\varepsilon_{
1} \bm{u}^{o}_{+}(\xi) \! + \! \varepsilon_{2} \bm{d}_{o})} \! + \! \dfrac{
\mathfrak{v}^{o}(\varepsilon_{1},\varepsilon_{2},\bm{\Omega}^{o};\xi)}{\bm{
\theta}^{o}(\varepsilon_{1} \bm{u}^{o}_{+}(\xi) \! - \! \frac{1}{2 \pi}(n \!
+ \! \frac{1}{2}) \bm{\Omega}^{o} \! + \! \varepsilon_{2} \bm{d}_{o})} \! - 
\! \left(\dfrac{\mathfrak{u}^{o}(\varepsilon_{1},\varepsilon_{2},\bm{0};
\xi)}{\bm{\theta}^{o}(\varepsilon_{1} \bm{u}^{o}_{+}(\xi) \! + \! 
\varepsilon_{2} \bm{d}_{o})} \right)^{2} \\
+& \, \dfrac{\mathfrak{u}^{o}(\varepsilon_{1},\varepsilon_{2},\bm{0};\xi) 
\mathfrak{u}^{o}(\varepsilon_{1},\varepsilon_{2},\bm{\Omega}^{o};\xi)}{\bm{
\theta}^{o}(\varepsilon_{1} \bm{u}^{o}_{+}(\xi) \! + \! \varepsilon_{2} 
\bm{d}_{o}) \bm{\theta}^{o}(\varepsilon_{1} \bm{u}^{o}_{+}(\xi) \! - \! 
\frac{1}{2 \pi}(n \! + \! \frac{1}{2}) \bm{\Omega}^{o} \! + \! 
\varepsilon_{2} \bm{d}_{o})}, \\
\beth^{\varepsilon_{1}}_{\varepsilon_{2}}(\xi) =& \, -\dfrac{\mathfrak{w}^{
o}(\varepsilon_{1},\varepsilon_{2},\bm{0};\xi)}{\bm{\theta}^{o}(\varepsilon_{
1} \bm{u}^{o}_{+}(\xi) \! + \! \varepsilon_{2} \bm{d}_{o})} \! + \! \dfrac{
\mathfrak{w}^{o}(\varepsilon_{1},\varepsilon_{2},\bm{\Omega}^{o};\xi)}{\bm{
\theta}^{o}(\varepsilon_{1} \bm{u}^{o}_{+}(\xi) \! - \! \frac{1}{2 \pi}(n 
\! + \! \frac{1}{2}) \bm{\Omega}^{o} \! + \! \varepsilon_{2} \bm{d}_{o})} \! 
+ \! \dfrac{2 \mathfrak{u}^{o}(\varepsilon_{1},\varepsilon_{2},\bm{0};\xi) 
\mathfrak{v}^{o}(\varepsilon_{1},\varepsilon_{2},\bm{0};\xi)}{(\bm{\theta}^{
o}(\varepsilon_{1} \bm{u}^{o}_{+}(\xi) \! + \! \varepsilon_{2} \bm{d}_{o}))^{
2}} \\
-& \, \dfrac{\mathfrak{v}^{o}(\varepsilon_{1},\varepsilon_{2},\bm{0};\xi) 
\mathfrak{u}^{o}(\varepsilon_{1},\varepsilon_{2},\bm{\Omega}^{o};\xi)}{\bm{
\theta}^{o}(\varepsilon_{1} \bm{u}^{o}_{+}(\xi) \! + \! \varepsilon_{2} \bm{
d}_{o}) \bm{\theta}^{o}(\varepsilon_{1} \bm{u}^{o}_{+}(\xi) \! - \! \frac{
1}{2 \pi}(n \! + \! \frac{1}{2}) \bm{\Omega}^{o} \! + \! \varepsilon_{2} 
\bm{d}_{o})} \! + \! \left(\dfrac{\mathfrak{u}^{o}(\varepsilon_{1},
\varepsilon_{2},\bm{0};\xi)}{\bm{\theta}^{o}(\varepsilon_{1} \bm{u}^{o}_{+}
(\xi) \! + \! \varepsilon_{2} \bm{d}_{o})} \right)^{3} \\
-& \, \dfrac{\mathfrak{u}^{o}(\varepsilon_{1},\varepsilon_{2},\bm{0};\xi) 
\mathfrak{v}^{o}(\varepsilon_{1},\varepsilon_{2},\bm{\Omega}^{o};\xi)}{\bm{
\theta}^{o}(\varepsilon_{1} \bm{u}^{o}_{+}(\xi) \! + \! \varepsilon_{2} \bm{
d}_{o}) \bm{\theta}^{o}(\varepsilon_{1} \bm{u}^{o}_{+}(\xi) \! - \! \frac{1}{
2 \pi}(n \! + \! \frac{1}{2}) \bm{\Omega}^{o} \! + \! \varepsilon_{2} \bm{
d}_{o})} \! - \! \dfrac{(\mathfrak{u}^{o}(\varepsilon_{1},\varepsilon_{2},
\bm{0};\xi))^{2}}{(\bm{\theta}^{o}(\varepsilon_{1} \bm{u}^{o}_{+}(\xi) \! + 
\! \varepsilon_{2}\bm{d}_{o}))^{2}} \\
\times& \, \dfrac{\mathfrak{u}^{o}(\varepsilon_{1},\varepsilon_{2},\bm{
\Omega}^{o};\xi)}{\bm{\theta}^{o}(\varepsilon_{1} \bm{u}^{o}_{+}(\xi) \! - \! 
\frac{1}{2 \pi}(n \! + \! \frac{1}{2}) \bm{\Omega}^{o} \! + \! \varepsilon_{
2} \bm{d}_{o})}, \\
\gimel^{\varepsilon_{1}}_{\varepsilon_{2}}(\xi) =& \, -\dfrac{\mathfrak{z}^{
o}(\varepsilon_{1},\varepsilon_{2},\bm{0};\xi)}{\bm{\theta}^{o}(\varepsilon_{
1} \bm{u}^{o}_{+}(\xi) \! + \! \varepsilon_{2} \bm{d}_{o})} \! + \! \dfrac{
\mathfrak{z}^{o}(\varepsilon_{1},\varepsilon_{2},\bm{\Omega}^{o};\xi)}{\bm{
\theta}^{o}(\varepsilon_{1} \bm{u}^{o}_{+}(\xi) \! - \! \frac{1}{2 \pi}(n 
\! + \! \frac{1}{2}) \bm{\Omega}^{o} \! + \! \varepsilon_{2} \bm{d}_{o})} \! 
+ \! \left(\dfrac{\mathfrak{v}^{o}(\varepsilon_{1},\varepsilon_{2},\bm{0};
\xi)}{\bm{\theta}^{o}(\varepsilon_{1} \bm{u}^{o}_{+}(\xi) \! + \! 
\varepsilon_{2} \bm{d}_{o})} \right)^{2} \\
-& \, \dfrac{\mathfrak{v}^{o}(\varepsilon_{1},\varepsilon_{2},\bm{0};\xi) 
\mathfrak{v}^{o}(\varepsilon_{1},\varepsilon_{2},\bm{\Omega}^{o};\xi)}{\bm{
\theta}^{o}(\varepsilon_{1}\bm{u}^{o}_{+}(\xi) \! + \! \varepsilon_{2} \bm{
d}_{o}) \bm{\theta}^{o}(\varepsilon_{1} \bm{u}^{o}_{+}(\xi) \! - \! \frac{
1}{2 \pi}(n \! + \! \frac{1}{2}) \bm{\Omega}^{o} \! + \! \varepsilon_{2} 
\bm{d}_{o})} \! - \! \dfrac{2 \mathfrak{u}^{o}(\varepsilon_{1},\varepsilon_{
2},\bm{0};\xi) \mathfrak{w}^{o}(\varepsilon_{1},\varepsilon_{2},\bm{0};\xi)
}{(\bm{\theta}^{o}(\varepsilon_{1} \bm{u}^{o}_{+}(\xi) \! + \! \varepsilon_{
2} \bm{d}_{o}))^{2}} \\
+& \, \dfrac{\mathfrak{w}^{o}(\varepsilon_{1},\varepsilon_{2},\bm{0};\xi) 
\mathfrak{u}^{o}(\varepsilon_{1},\varepsilon_{2},\bm{\Omega}^{o};\xi)}{\bm{
\theta}^{o}(\varepsilon_{1} \bm{u}^{o}_{+}(\xi) \! + \! \varepsilon_{2} \bm{
d}_{o}) \bm{\theta}^{o}(\varepsilon_{1} \bm{u}^{o}_{+}(\xi) \! - \! \frac{
1}{2 \pi}(n \! + \! \frac{1}{2}) \bm{\Omega}^{o} \! + \! \varepsilon_{2} 
\bm{d}_{o})} \! + \! \dfrac{3(\mathfrak{u}^{o}(\varepsilon_{1},\varepsilon_{
2},\bm{0};\xi))^{2} \mathfrak{v}^{o}(\varepsilon_{1},\varepsilon_{2},\bm{0};
\xi)}{(\bm{\theta}^{o}(\varepsilon_{1} \bm{u}^{o}_{+}(\xi) \! + \! 
\varepsilon_{2} \bm{d}_{o}))^{3}} \\
+& \, \dfrac{\mathfrak{u}^{o}(\varepsilon_{1},\varepsilon_{2},\bm{0};\xi) 
\mathfrak{w}^{o}(\varepsilon_{1},\varepsilon_{2},\bm{\Omega}^{o};\xi)}{\bm{
\theta}^{o}(\varepsilon_{1} \bm{u}^{o}_{+}(\xi) \! + \! \varepsilon_{2} \bm{
d}_{o}) \bm{\theta}^{o}(\varepsilon_{1} \bm{u}^{o}_{+}(\xi) \! - \! \frac{
1}{2 \pi}(n \! + \! \frac{1}{2}) \bm{\Omega}^{o} \! + \! \varepsilon_{2} 
\bm{d}_{o})} \! + \! \left(\dfrac{\mathfrak{u}^{o}(\varepsilon_{1},
\varepsilon_{2},\bm{0};\xi)}{\bm{\theta}^{o}(\varepsilon_{1} \bm{u}^{o}_{+}
(\xi) \! + \! \varepsilon_{2} \bm{d}_{o})} \right)^{4} \\
-& \, \dfrac{2 \mathfrak{u}^{o}(\varepsilon_{1},\varepsilon_{2},\bm{0};\xi) 
\mathfrak{v}^{o}(\varepsilon_{1},\varepsilon_{2},\bm{0};\xi) \mathfrak{u}^{
o}(\varepsilon_{1},\varepsilon_{2},\bm{\Omega}^{o};\xi)}{(\bm{\theta}^{o}
(\varepsilon_{1} \bm{u}^{o}_{+}(\xi) \! + \! \varepsilon_{2} \bm{d}_{o}))^{2} 
\bm{\theta}^{o}(\varepsilon_{1} \bm{u}^{o}_{+}(\xi) \! - \! \frac{1}{2 \pi}
(n \! + \! \frac{1}{2}) \bm{\Omega}^{o} \! + \! \varepsilon_{2} \bm{d}_{o})} 
\! - \! \dfrac{(\mathfrak{u}^{o}(\varepsilon_{1},\varepsilon_{2},\bm{0};\xi)
)^{2}}{(\bm{\theta}^{o}(\varepsilon_{1} \bm{u}^{o}_{+}(\xi) \! + \! 
\varepsilon_{2} \bm{d}_{o}))^{2}} \\
\times& \, \dfrac{\mathfrak{v}^{o}(\varepsilon_{1},\varepsilon_{2},\bm{
\Omega}^{o};\xi)}{\bm{\theta}^{o}(\varepsilon_{1} \bm{u}^{o}_{+}(\xi) \! - \! 
\frac{1}{2 \pi}(n \! + \! \frac{1}{2}) \bm{\Omega}^{o} \! + \! \varepsilon_{
2} \bm{d}_{o})} \! - \! \dfrac{(\mathfrak{u}^{o}(\varepsilon_{1},\varepsilon_{
2},\bm{0};\xi))^{3} \mathfrak{u}^{o}(\varepsilon_{1},\varepsilon_{2},\bm{
\Omega}^{o};\xi)}{(\bm{\theta}^{o}(\varepsilon_{1} \bm{u}^{o}_{+}(\xi) \! + 
\! \varepsilon_{2} \bm{d}_{o}))^{3} \bm{\theta}^{o}(\varepsilon_{1} \bm{u}^{
o}_{+}(\xi) \! - \! \frac{1}{2 \pi}(n \! + \! \frac{1}{2}) \bm{\Omega}^{o} \! 
+ \! \varepsilon_{2} \bm{d}_{o})},
\end{align*}
with $\bm{0} \! := \! (0,0,\dotsc,0)^{\operatorname{T}}$ $(\in \! \mathbb{
R}^{N})$,
\begin{gather*}
\mathfrak{u}^{o}(\varepsilon_{1},\varepsilon_{2},\bm{\Omega}^{o};\xi) \!
:= \! 2 \pi \overset{o}{\Lambda}^{\raise-1.0ex\hbox{$\scriptstyle 1$}}_{0}
(\varepsilon_{1},\varepsilon_{2},\bm{\Omega}^{o};\xi), \qquad \qquad 
\mathfrak{v}^{o}(\varepsilon_{1},\varepsilon_{2},\bm{\Omega}^{o};\xi) \! := 
\! -2 \pi^{2} \overset{o}{\Lambda}^{\raise-1.0ex\hbox{$\scriptstyle 2$}}_{0}
(\varepsilon_{1},\varepsilon_{2},\bm{\Omega}^{o};\xi), \\
\mathfrak{w}^{o}(\varepsilon_{1},\varepsilon_{2},\bm{\Omega}^{o};\xi) \! := \! 
2 \pi \! \left(\overset{o}{\Lambda}^{\raise-1.0ex\hbox{$\scriptstyle 0$}}_{1}
(\varepsilon_{1},\varepsilon_{2},\bm{\Omega}^{o};\xi) \! - \! \dfrac{2 \pi^{
2}}{3} \overset{o}{\Lambda}^{\raise-1.0ex\hbox{$\scriptstyle 3$}}_{0}
(\varepsilon_{1},\varepsilon_{2},\bm{\Omega}^{o};\xi) \right), \\
\mathfrak{z}^{o}(\varepsilon_{1},\varepsilon_{2},\bm{\Omega}^{o};\xi) \! := 
\! -(2 \pi)^{2} \! \left(\overset{o}{\Lambda}^{\raise-1.0ex\hbox{$\scriptstyle 
1$}}_{1}(\varepsilon_{1},\varepsilon_{2},\bm{\Omega}^{o};\xi) \! - \! \dfrac{
\pi^{2}}{6} \overset{o}{\Lambda}^{\raise-1.0ex\hbox{$\scriptstyle 4$}}_{0}
(\varepsilon_{1},\varepsilon_{2},\bm{\Omega}^{o};\xi) \right), \\
\overset{o}{\Lambda}^{\raise-1.0ex\hbox{$\scriptstyle j_{1}$}}_{j_{2}}
(\varepsilon_{1},\varepsilon_{2},\bm{\Omega}^{o};\xi) \! = \! \sum_{m \in
\mathbb{Z}^{N}}(\mathfrak{r}_{o}(\xi))^{j_{1}}(\mathfrak{s}_{o}(\xi))^{j_{
2}} \, \me^{2 \pi \mi (m,\varepsilon_{1} \bm{u}^{o}_{+}(\xi)-\frac{1}{2 \pi}
(n+\frac{1}{2}) \bm{\Omega}^{o}+ \varepsilon_{2} \bm{d}_{o})+ \pi \mi 
(m,\bm{\tau}^{o}m)}, \\
\mathfrak{r}_{o}(\xi) \! := \! \dfrac{2(m,\vec{\moo}^{o}_{1}(\xi))}{\mathfrak{
f}_{o}(\xi)}, \qquad \qquad \mathfrak{s}_{o}(\xi) \! := \! \dfrac{2(m,\vec{
\moo}^{o}_{2}(\xi))}{3 \mathfrak{f}_{o}(\xi)}, \\
\vec{\moo}^{o}_{1}(\xi) \! = \! \left(\mathfrak{p}_{1}^{\natural}(\xi),
\mathfrak{p}_{2}^{\natural}(\xi),\dotsc,\mathfrak{p}_{N}^{\natural}(\xi) 
\right), \qquad \qquad \vec{\moo}^{o}_{2}(\xi) \! = \! \left(\mathfrak{q}_{
1}^{\natural}(\xi),\mathfrak{q}_{2}^{\natural}(\xi),\dotsc,\mathfrak{q}_{
N}^{\natural}(\xi) \right).
\end{gather*}

Recall the definition of $\gamma^{o}(z)$ given in Lemma~4.4: a careful
analysis of the branch cuts shows that, for $j \! = \! 1,\dotsc,N$,
\begin{align*}
(\gamma^{o}(z))^{2} &\underset{z \in \mathbb{C}_{\pm}}{=} \pm \, \dfrac{\left(
\! (z \! - \! b_{j}^{o}) \prod_{\substack{k=1\\k \not= j}}^{N} \! \left(\!
\frac{z-b_{k}^{o}}{z-a_{k}^{o}} \right) \! \left(\! \frac{z-b_{0}^{o}}{z-a_{
N+1}^{o}} \right) \right)^{1/2}}{\sqrt{\smash[b]{z \! - \! a_{j}^{o}}}} \\
&\underset{\mathbb{C}_{\pm} \ni z \to a_{j}^{o}}{=} \pm \, \dfrac{\left(
Q_{0}^{o}(a_{j}^{o}) \! + \! Q_{1}^{o}(a_{j}^{o})(z \! - \! a_{j}^{o}) \!
+ \! \frac{1}{2}Q_{2}^{o}(a_{j}^{o})(z \! - \! a_{j}^{o})^{2} \! + \!
\mathcal{O}((z \! - \! a_{j}^{o})^{3}) \right)}{\sqrt{\smash[b]{z \! - \!
a_{j}^{o}}}},
\end{align*}
where $Q_{0}^{o}(a_{j}^{o}),Q_{1}^{o}(a_{j}^{o})$, $j \! = \! 1,\dotsc,N$,
are given in Theorem~2.3.1, Equations~(2.35) and~(2.36), and
\begin{align*}
Q_{2}^{o}(a_{j}^{o})=& -\dfrac{1}{2}Q_{0}^{o}(a_{j}^{o}) \! \left(\sum_{
\substack{k=1\\k \not= j}}^{N} \! \left(\dfrac{1}{(a_{j}^{o} \! - \! b_{k}^{
o})^{2}} \! - \! \dfrac{1}{(a_{j}^{o} \! - \! a_{k}^{o})^{2}} \right) \! + \!
\dfrac{1}{(a_{j}^{o} \! - \! b_{0}^{o})^{2}} \! - \! \dfrac{1}{(a_{j}^{o} \!
- \! a_{N+1}^{o})^{2}} \! + \! \dfrac{1}{(a_{j}^{o} \! - \! b_{j}^{o})^{2}}
\right) \\
+& \, \dfrac{1}{4}Q_{0}^{o}(a_{j}^{o}) \! \left(\sum_{\substack{k=1\\k \not=
j}}^{N} \! \left(\dfrac{1}{a_{j}^{o} \! - \! b_{k}^{o}} \! - \! \dfrac{1}{a_{
j}^{o} \! - \! a_{k}^{o}} \right) \! + \! \dfrac{1}{a_{j}^{o} \! - \! b_{0}^{
o}} \! - \! \dfrac{1}{a_{j}^{o} \! - \! a_{N+1}^{o}} \! + \! \dfrac{1}{a_{
j}^{o} \! - \! b_{j}^{o}} \right)^{2}, \quad j \! = \! 1,\dotsc,N.
\end{align*}
Recall the formula above for
$\overset{o}{\mathfrak{M}}^{\raise-1.0ex\hbox{$\scriptstyle \infty$}}(z) \!
\left(
\begin{smallmatrix}
\mp (s_{1}+t_{1}) & \pm \mi (s_{1}-t_{1}) \me^{\mi (n+\frac{1}{2}) \Omega_{j}^{
o}} \\
\pm \mi (s_{1}-t_{1}) \me^{-\mi (n+\frac{1}{2}) \Omega_{j}^{o}} & \pm (s_{1}+
t_{1})
\end{smallmatrix}
\right) \!
(\overset{o}{\mathfrak{M}}^{\raise-1.0ex\hbox{$\scriptstyle \infty$}}(z))^{-
1}$. Substituting the a\-b\-o\-v\-e expansions (as $z \! \to \! a_{j}^{o}$, 
$j \! = \! 1,\dotsc,N)$ into this formula, equating coefficients of like 
powers of $(z \! - \! a_{j}^{o})^{-p/2}(\widetilde{n}G_{a_{j}}^{o}(z))^{-1}$, 
$p \! \in \! \lbrace 4,3,2,1,0 \rbrace$, with $\widetilde{n} \! := \! n \! 
+ \! \tfrac{1}{2}$, and considering, say, the $(1 \, 1)$-element of the 
resulting (asymptotic) expansions, one arrives at (modulo a minus sign, this 
result is equally applicable to the $(2 \, 2)$-element, since $\operatorname{
tr}(w_{+}^{\Sigma^{o}_{\circlearrowright}}(z)) \! = \! 0)$, up to terms that 
are $\mathcal{O}((\widetilde{n}^{2}(z \! - \! a_{j}^{o})^{3}(G_{a_{j}}^{o}
(z))^{2})^{-1} 
\overset{o}{\mathfrak{M}}^{\raise-1.0ex\hbox{$\scriptstyle \infty$}}(z)f_{a_{
j}}^{o}(n)(\overset{o}{\mathfrak{M}}^{\raise-1.0ex\hbox{$\scriptstyle \infty$}
}(z))^{-1})$, upon setting, for economy of notation, $Q_{q}^{o}(a_{j}^{o})
\! =: \! Q_{q}$, $q \! = \! 0,1,2$, $\gamma^{o}(0) \! =: \! \gamma^{o}$,
$\varkappa_{1}^{o}(a_{j}^{o}) \! =: \! \varkappa_{1}^{o}$, $\varkappa_{2}^{o}
(a_{j}^{o}) \! =: \! \varkappa_{2}^{o}$, $\aleph^{\varepsilon_{1}}_{
\varepsilon_{2}}(a_{j}^{o}) \! =: \! \aleph^{\varepsilon_{1}}_{\varepsilon_{
2}}$, $\daleth^{\varepsilon_{1}}_{\varepsilon_{2}}(a_{j}^{o}) \! =: \!
\daleth^{\varepsilon_{1}}_{\varepsilon_{2}}$, $\beth^{\varepsilon_{1}}_{
\varepsilon_{2}}(a_{j}^{o}) \! =: \! \beth^{\varepsilon_{1}}_{\varepsilon_{
2}}$, and $\gimel^{\varepsilon_{1}}_{\varepsilon_{2}}(a_{j}^{o}) \! =: \!
\gimel^{\varepsilon_{1}}_{\varepsilon_{2}}$:
\begin{align*}
\mathcal{O} \! \left(\dfrac{(z \! - \! a_{j}^{o})^{-2} \, \me^{\mi 
\widetilde{n} \Omega_{j}^{o}}}{\widetilde{n}G_{a_{j}}^{o}(z)} \right):& 
-\dfrac{(s_{1} \! + \! t_{1}) \varkappa_{1}^{o} \varkappa_{2}^{o}Q_{0}}{
4(\gamma^{o})^{2}} \! - \! \dfrac{(s_{1} \! + \! t_{1}) \varkappa_{1}^{o} 
\varkappa_{2}^{o}Q_{0}}{4(\gamma^{o})^{2}} \! - \! \dfrac{(s_{1} \! - \! 
t_{1}) \varkappa_{1}^{o} \varkappa_{2}^{o}Q_{0}}{4(\gamma^{o})^{2}} \! - 
\! \dfrac{(s_{1} \! - \! t_{1}) \varkappa_{1}^{o} \varkappa_{2}^{o}Q_{0}}{4
(\gamma^{o})^{2}}; \\
\mathcal{O} \! \left(\dfrac{(z \! - \! a_{j}^{o})^{-3/2} \, \me^{\mi 
\widetilde{n} \Omega_{j}^{o}}}{\widetilde{n}G_{a_{j}}^{o}(z)} \right):& 
-\dfrac{\mi (s_{1} \! + \! t_{1}) \varkappa_{1}^{o} \varkappa_{2}^{o}Q_{0}
(\aleph_{1}^{1} \! - \! \aleph_{-1}^{-1})}{4(\gamma^{o})^{2}} \! - \! 
\dfrac{\mi (s_{1} \! + \! t_{1}) \varkappa_{1}^{o} \varkappa_{2}^{o}Q_{0}
(\aleph^{1}_{-1} \! - \! \aleph^{-1}_{1})}{4(\gamma^{o})^{2}} \\
&-\dfrac{\mi (s_{1} \! - \! t_{1}) \varkappa_{1}^{o} \varkappa_{2}^{o}Q_{0}
(\aleph^{1}_{-1} \! + \! \aleph^{1}_{1})}{4(\gamma^{o})^{2}} \! + \! 
\dfrac{\mi (s_{1} \! - \! t_{1}) \varkappa_{1}^{o} \varkappa_{2}^{o}Q_{0}
(\aleph^{-1}_{-1} \! + \! \aleph^{-1}_{1})}{4(\gamma^{o})^{2}} \\
&-\dfrac{(s_{1} \! + \! t_{1}) \varkappa_{1}^{o} \varkappa_{2}^{o}}{2} \!
+ \! \dfrac{(s_{1} \! + \! t_{1}) \varkappa_{1}^{o} \varkappa_{2}^{o}}{2};
\\
\mathcal{O} \! \left(\dfrac{(z \! - \! a_{j}^{o})^{-1} \, \me^{\mi 
\widetilde{n} \Omega_{j}^{o}}}{\widetilde{n}G_{a_{j}}^{o}(z)} \right):& 
-\dfrac{(s_{1} \! + \! t_{1}) \varkappa_{1}^{o} \varkappa_{2}^{o}}{4
(\gamma^{o})^{2}} \! \left(Q_{1} \! + \! Q_{0} \! \left(\daleth^{-1}_{-1} 
\! + \! \daleth^{1}_{1} \! + \! \aleph^{1}_{1} \aleph^{-1}_{-1} \right) 
\right) \! - \! \dfrac{(s_{1} \! + \! t_{1}) \varkappa_{1}^{o} \varkappa_{
2}^{o}(\gamma^{o})^{2}}{4Q_{0}} \\
&-\dfrac{(s_{1} \! + \! t_{1}) \varkappa_{1}^{o} \varkappa_{2}^{o}}{4
(\gamma^{o})^{2}} \! \left(Q_{1} \! + \! Q_{0} \! \left(\daleth^{1}_{-1} \! 
+ \! \daleth^{-1}_{1} \! + \! \aleph^{-1}_{1} \aleph^{1}_{-1} \right) \right) 
\! - \! \dfrac{(s_{1} \! + \! t_{1}) \varkappa_{1}^{o} \varkappa_{2}^{o}
(\gamma^{o})^{2}}{4Q_{0}} \\
&-\dfrac{(s_{1} \! - \! t_{1}) \varkappa_{1}^{o} \varkappa_{2}^{o}}{4
(\gamma^{o})^{2}} \! \left(Q_{1} \! + \! Q_{0} \! \left(\daleth^{1}_{-1} \! 
+ \! \daleth^{1}_{1} \! - \! \aleph^{1}_{1} \aleph^{1}_{-1} \right) \right) 
\! + \! \dfrac{(s_{1} \! - \! t_{1}) \varkappa_{1}^{o} \varkappa_{2}^{o}
(\gamma^{o})^{2}}{4Q_{0}} \\
&-\dfrac{(s_{1} \! - \! t_{1}) \varkappa_{1}^{o} \varkappa_{2}^{o}}{4
(\gamma^{o})^{2}} \! \left(Q_{1} \! + \! Q_{0} \! \left(\daleth^{-1}_{-1} \! 
+ \! \daleth^{-1}_{1} \! - \! \aleph^{-1}_{1} \aleph^{-1}_{-1} \right) \right) 
\! + \! \dfrac{(s_{1} \! - \! t_{1}) \varkappa_{1}^{o} \varkappa_{2}^{o}
(\gamma^{o})^{2}}{4Q_{0}} \\
&-\dfrac{\mi (s_{1} \! + \! t_{1}) \varkappa_{1}^{o} \varkappa_{2}^{o}}{2} 
\! \left(\aleph^{1}_{1} \! - \! \aleph^{-1}_{-1} \right) \! + \! \dfrac{\mi 
(s_{1} \! + \! t_{1}) \varkappa_{1}^{o} \varkappa_{2}^{o}}{2} \! \left(
\aleph^{1}_{-1} \! - \!\aleph^{-1}_{1} \right); \\
\mathcal{O} \! \left(\dfrac{(z \! - \! a_{j}^{o})^{-1/2} \, \me^{\mi 
\widetilde{n} \Omega_{j}^{o}}}{\widetilde{n}G_{a_{j}}^{o}(z)} \right):& 
-\dfrac{\mi (s_{1} \! + \! t_{1}) \varkappa_{1}^{o} \varkappa_{2}^{o}}{4
(\gamma^{o})^{2}} \! \left(Q_{1} \! \left(\aleph^{1}_{1} \! - \! \aleph^{-
1}_{-1} \right) \! + \! Q_{0} \! \left(\beth^{1}_{1} \! - \! \beth^{-1}_{-1} 
\! + \! \aleph^{1}_{1} \daleth^{-1}_{-1} \! - \! \aleph^{-1}_{-1} \daleth^{
1}_{1} \right) \right) \\
&-\dfrac{\mi (s_{1} \! + \! t_{1}) \varkappa_{1}^{o} \varkappa_{2}^{o}}{4
(\gamma^{o})^{2}} \! \left(Q_{1} \! \left(\aleph^{1}_{-1} \! - \! \aleph^{-
1}_{1} \right) \! + \! Q_{0} \! \left(\beth^{1}_{-1} \! - \! \beth^{-1}_{1} 
\! + \! \aleph^{1}_{-1} \daleth^{-1}_{1} \! - \! \aleph^{-1}_{1} \daleth^{
1}_{-1} \right) \right) \\
&-\dfrac{\mi (s_{1} \! - \! t_{1}) \varkappa_{1}^{o} \varkappa_{2}^{o}}{4
(\gamma^{o})^{2}} \! \left(Q_{1} \! \left(\aleph^{1}_{-1} \! + \! \aleph^{
1}_{1} \right) \! + \! Q_{0} \! \left(\beth^{1}_{-1} \! + \! \beth^{1}_{1} 
\! + \! \aleph^{1}_{1} \daleth^{1}_{-1} \! + \! \aleph^{1}_{-1} \daleth^{
1}_{1} \right) \right) \\
&+\dfrac{\mi (s_{1} \! - \! t_{1}) \varkappa_{1}^{o} \varkappa_{2}^{o}}{4
(\gamma^{o})^{2}} \! \left(Q_{1} \! \left(\aleph^{-1}_{-1} \! + \! \aleph^{-
1}_{1} \right) \! + \! Q_{0} \! \left(\beth^{-1}_{-1} \! + \! \beth^{-1}_{1} 
\! + \! \aleph^{-1}_{1} \daleth^{-1}_{-1} \! + \! \aleph^{-1}_{-1} \daleth^{-
1}_{1} \right) \right) \\
&-\dfrac{\mi (s_{1} \! + \! t_{1}) \varkappa_{1}^{o} \varkappa_{2}^{o}
(\gamma^{o})^{2}}{4Q_{0}} \! \left(\aleph^{1}_{1} \! - \! \aleph^{-1}_{-1} 
\right) \! - \! \dfrac{(s_{1} \! + \! t_{1}) \varkappa_{1}^{o} \varkappa_{
2}^{o}}{2} \! \left(\daleth^{-1}_{-1} \! + \! \daleth^{1}_{1} \! + \! 
\aleph^{1}_{1} \aleph^{-1}_{-1} \right) \\
&-\dfrac{\mi (s_{1} \! + \! t_{1}) \varkappa_{1}^{o} \varkappa_{2}^{o}
(\gamma^{o})^{2}}{4Q_{0}} \! \left(\aleph^{1}_{-1} \! - \! \aleph^{-1}_{1} 
\right) \! + \! \dfrac{(s_{1} \! + \! t_{1}) \varkappa_{1}^{o} \varkappa_{
2}^{o}}{2} \! \left(\daleth^{1}_{-1} \! + \! \daleth^{-1}_{1} \! + \! 
\aleph^{-1}_{1} \aleph^{1}_{-1} \right) \\
&+\dfrac{\mi (s_{1} \! - \! t_{1}) \varkappa_{1}^{o} \varkappa_{2}^{o}
(\gamma^{o})^{2}}{4Q_{0}} \! \left(\aleph^{1}_{-1} \! + \! \aleph^{1}_{1} 
\right) \! - \! \dfrac{\mi (s_{1} \! - \! t_{1}) \varkappa_{1}^{o} 
\varkappa_{2}^{o}(\gamma^{o})^{2}}{4Q_{0}} \! \left(\aleph^{-1}_{1} \! + \! 
\aleph^{-1}_{-1} \right); \\
\mathcal{O} \! \left(\dfrac{\me^{\mi \widetilde{n} \Omega_{j}^{o}}}{\widetilde{
n}G_{a_{j}}^{o}(z)} \right):& -\dfrac{(s_{1} \! + \! t_{1}) \varkappa_{1}^{o}
\varkappa_{2}^{o}}{4(\gamma^{o})^{2}} \! \left(Q_{0} \! \left(\gimel^{-1}_{-1}
\! + \! \gimel^{1}_{1} \! + \! \daleth^{1}_{1} \daleth^{-1}_{-1} \! + \! 
\aleph^{-1}_{-1} \beth^{1}_{1} \! + \! \aleph^{1}_{1} \beth^{-1}_{-1} \right) 
\! + \! Q_{1} \! \left(\daleth^{-1}_{-1} \! + \! \daleth^{1}_{1} \! + \! 
\aleph^{1}_{1} \aleph^{-1}_{-1} \right) \right. \\
&\left. + \, \dfrac{1}{2}Q_{2} \right) \! - \! \dfrac{(s_{1} \! + \! t_{1}) 
\varkappa_{1}^{o} \varkappa_{2}^{o}}{4(\gamma^{o})^{2}} \! \left(Q_{0} \! 
\left(\gimel^{1}_{-1} \! + \! \gimel^{-1}_{1} \! + \! \daleth^{-1}_{1} 
\daleth^{1}_{-1} \! + \! \aleph^{-1}_{1} \beth^{1}_{-1} \! + \! \aleph^{1}_{
-1} \beth^{-1}_{1} \right) \! + \! \dfrac{1}{2}Q_{2} \right. \\
&\left. + \, Q_{1} \! \left(\daleth^{1}_{-1} \! + \! \daleth^{-1}_{1} \! + \! 
\aleph^{-1}_{1} \aleph^{1}_{-1} \right) \right) \! - \! \dfrac{(s_{1} \! - \! 
t_{1}) \varkappa_{1}^{o} \varkappa_{2}^{o}}{4(\gamma^{o})^{2}} \! \left(Q_{0} 
\! \left(\gimel^{1}_{-1} \! + \! \gimel^{1}_{1} \! + \! \daleth^{1}_{1} 
\daleth^{1}_{-1} \! - \! \aleph^{1}_{1} \beth^{1}_{-1} \! - \! \aleph^{1}_{-
1} \beth^{1}_{1} \right) \right. \\
&\left. + \, \dfrac{1}{2}Q_{2} \! + \! Q_{1} \! \left(\daleth^{1}_{-1} \! + \! 
\daleth^{1}_{1} \! - \! \aleph^{1}_{1} \aleph^{1}_{-1} \right) \right) \! - \! 
\dfrac{(s_{1} \! - \! t_{1}) \varkappa_{1}^{o} \varkappa_{2}^{o}}{4(\gamma^{
o})^{2}} \! \left(Q_{0} \! \left(\gimel^{-1}_{-1} \! + \! \gimel^{-1}_{1} \! 
+ \! \daleth^{-1}_{1} \daleth^{-1}_{-1} \! - \! \aleph^{-1}_{1} \beth^{-1}_{-
1} \right. \right. \\
&\left. \left. - \, \aleph^{-1}_{-1} \beth^{-1}_{1} \right) \! + \! \dfrac{
1}{2}Q_{2} \! + \! Q_{1} \! \left(\daleth^{-1}_{-1} \! + \! \daleth^{-1}_{1} 
\! - \! \aleph^{-1}_{1} \aleph^{-1}_{-1} \right) \right) \! - \! \dfrac{(s_{
1} \! + \! t_{1}) \varkappa_{1}^{o} \varkappa_{2}^{o}(\gamma^{o})^{2}}{4Q_{
0}} \! \left(\daleth^{-1}_{-1} \! + \! \daleth^{1}_{1} \! + \! \aleph^{1}_{1} 
\aleph^{-1}_{-1} \right. \\
&\left. - \, Q_{1}(Q_{0})^{-1} \right) \! - \! \dfrac{(s_{1} \! + \! t_{1}) 
\varkappa_{1}^{o} \varkappa_{2}^{o}(\gamma^{o})^{2}}{4Q_{0}} \! \left(
\daleth^{1}_{-1} \! + \! \daleth^{-1}_{1} \! + \! \aleph^{-1}_{1} \aleph^{
1}_{-1} \! - \! Q_{1}(Q_{0})^{-1} \right) \! - \! \dfrac{\mi (s_{1} \! + \! 
t_{1}) \varkappa_{1}^{o} \varkappa_{2}^{o}}{2} \\
&\times \, \left(\beth^{1}_{1} \! - \! \beth^{-1}_{-1} \! + \! \aleph^{1}_{1} 
\daleth^{-1}_{-1} \! - \! \aleph^{-1}_{-1} \daleth^{1}_{1} \right) \! + \! 
\dfrac{\mi (s_{1} \! + \! t_{1}) \varkappa_{1}^{o} \varkappa_{2}^{o}}{2} \! 
\left(\beth^{1}_{-1} \! - \! \beth^{-1}_{1} \! + \! \aleph^{1}_{-1} \daleth^{
-1}_{1} \! - \! \aleph^{-1}_{1} \daleth^{1}_{-1} \right) \\
&+ \, \dfrac{(s_{1} \! - \! t_{1}) \varkappa_{1}^{o} \varkappa_{2}^{o}
(\gamma^{o})^{2}}{4Q_{0}} \! \left( \daleth^{1}_{-1} \! + \! \daleth^{1}_{1} 
\! - \! \aleph^{1}_{1} \aleph^{1}_{-1} \! - \! Q_{1}(Q_{0})^{-1} \right) 
\! + \! \dfrac{(s_{1} \! - \! t_{1}) \varkappa_{1}^{o} \varkappa_{2}^{o}
(\gamma^{o})^{2}}{4Q_{0}} \! \left(\daleth^{-1}_{-1} \! + \! \daleth^{-1}_{1} 
\right. \\
&\left. - \, \aleph^{-1}_{1} \aleph^{-1}_{-1} \! - \!Q_{1}(Q_{0})^{-1} 
\right).
\end{align*}
Repeating the above analysis, \emph{mutatis mutandis}, for the $(1 \, 2)$- 
and $(2 \, 1)$-elements, substituting $\aleph^{-1}_{1} \! = \! \aleph^{1}_{
1}$, $\aleph^{-1}_{-1} \! = \! \aleph^{1}_{-1}$, $\daleth^{-1}_{1} \! = \! 
\daleth^{1}_{1}$, $\daleth^{-1}_{-1} \! = \! \daleth^{1}_{-1}$, $\beth^{-1}_{
1} \! = \! \beth^{1}_{1}$, $\beth^{-1}_{-1} \! = \! \beth^{1}_{-1}$, 
$\gimel^{-1}_{1} \! = \! \gimel^{1}_{1}$, and $\gimel^{-1}_{-1} \! = \! 
\gimel^{1}_{-1}$ into the above (and resulting) `coefficient equations', and 
simplifying, one shows that: (i) the coefficients of the terms that are 
$\mathcal{O}((z \! - \! a_{j}^{o})^{-p/2}((n \! + \! \tfrac{1}{2})G_{a_{j}}^{
o}(z))^{-1} \exp (\mi (n \! + \! \tfrac{1}{2}) \Omega_{j}^{o}))$, $p \! = \! 
1,3$, are equal to zero; and (ii) recalling {}from Lemma~4.7 that, for $z \! 
\in \! \mathbb{U}^{o}_{\delta_{a_{j}}} \setminus (-\infty,a_{j}^{o})$, $j 
\! = \! 1,\dotsc,N$, $G_{a_{j}}^{o}(z) \! =_{z \to a_{j}^{o}} \! \widehat{
\alpha}_{0} \! + \! \widehat{\alpha}_{1}(z \! - \! a_{j}^{o}) \! + \! 
\widehat{\alpha}_{2}(z \! - \! a_{j}^{o})^{2} \! + \! \mathcal{O}((z \! - \! 
a_{j}^{o})^{3})$, where $\widehat{\alpha}_{0} \! = \! \widehat{\alpha}_{0}^{
o}(a_{j}^{o}) \! := \! \tfrac{4}{3}f(a_{j}^{o})$, $\widehat{\alpha}_{1} \! 
= \! \widehat{\alpha}_{1}^{o}(a_{j}^{o}) \! := \! \tfrac{4}{3}f^{\prime}
(a_{j}^{o})$, and $\widehat{\alpha}_{2} \! = \! \widehat{\alpha}_{2}^{o}
(a_{j}^{o}) \! := \! \tfrac{2}{7}f^{\prime \prime}(a_{j}^{o})$, with $f
(a_{j}^{o})$, $f^{\prime}(a_{j}^{o})$, and $f^{\prime \prime}(a_{j}^{o})$ 
given in Lemma~4.7, substituting the expansion for $G^{o}_{a_{j}}(z)$ (as 
$z \! \to \! a_{j}^{o}$, $j \! = \! 1,\dotsc,N)$ into the remaining non-zero 
coefficient equations, collecting coefficients of like powers of $(z \! - \! 
a_{j}^{o})^{-p}$, $p \! = \! 0,1,2$, and continuing, analytically, the 
resulting (rational) expressions to $\partial \mathbb{U}^{o}_{\delta_{a_{
j}}}$, $j \! = \! 1,\dotsc,N$, one arrives at, after a lengthy algebraic 
calculation and reinserting explicit $a_{j}^{o}$, $j \! = \! 1,\dotsc,
N$, dependencies,
\begin{align}
w_{+}^{\Sigma^{o}_{\circlearrowright}}(z) \underset{\underset{z \in \partial
\mathbb{U}^{o}_{\delta_{a_{j}}}}{n \to \infty}}{=}& \dfrac{1}{n \! + \! \frac{
1}{2}} \! \left(\dfrac{\mathscr{A}^{o}(a_{j}^{o})}{\widehat{\alpha}_{0}^{o}
(a_{j}^{o})(z \! - \! a_{j}^{o})^{2}} \! + \! \dfrac{(\mathscr{B}^{o}(a_{j}^{
o}) \widehat{\alpha}_{0}^{o}(a_{j}^{o}) \! - \! \mathscr{A}^{o}(a_{j}^{o})
\widehat{\alpha}_{1}^{o}(a_{j}^{o}))}{(\widehat{\alpha}_{0}^{o}(a_{j}^{o}))^{2}
(z \! - \! a_{j}^{o})} \right. \nonumber \\
&\left. + \, \dfrac{\left(\mathscr{A}^{o}(a_{j}^{o}) \widehat{\alpha}_{0}^{o}
(a_{j}^{o}) \! \left(\! \left(\frac{\widehat{\alpha}_{1}^{o}(a_{j}^{o})}{
\widehat{\alpha}_{0}^{o}(a_{j}^{o})} \right)^{2} \! - \! \frac{\widehat{
\alpha}_{2}^{o}(a_{j}^{o})}{\widehat{\alpha}_{0}^{o}(a_{j}^{o})} \right) \! -
\! \mathscr{B}^{o}(a_{j}^{o}) \widehat{\alpha}_{1}^{o}(a_{j}^{o}) \! + \!
\mathscr{C}^{o}(a_{j}^{o}) \widehat{\alpha}_{0}^{o}(a_{j}^{o}) \right)}{(
\widehat{\alpha}_{0}^{o}(a_{j}^{o}))^{2}} \right) \nonumber \\
&+ \, \mathcal{O} \! \left(\dfrac{\sum_{k=1}^{\infty}f_{k}^{o}(n)(z \! - \!
a_{j}^{o})^{k}}{n \! + \! \frac{1}{2}} \right) \! + \! \mathcal{O} \! \left(
\dfrac{\overset{o}{\mathfrak{M}}^{\raise-1.0ex\hbox{$\scriptstyle \infty$}}
(z)f^{o}_{a_{j}}(n)
(\overset{o}{\mathfrak{M}}^{\raise-1.0ex\hbox{$\scriptstyle \infty$}}
(z))^{-1}}{(n \! + \! \frac{1}{2})^{2}(z \! - \! a_{j}^{o})^{3}(G^{o}_{a_{j}}
(z))^{2}} \right),
\end{align}
where $\mathscr{A}^{o}(a_{j}^{o}),\mathscr{B}^{o}(a_{j}^{o})$, $j \! = \! 1,
\dotsc,N$, are defined in Theorem~2.3.1, Equations~(2.25), (2.27), (2.28), 
(2.35)--(2.45), (2.49), (2.56), and~(2.57),
\begin{equation*}
\dfrac{\mathscr{C}^{o}(a_{j}^{o})}{\me^{\mi (n+\frac{1}{2}) \Omega_{j}^{o}}} 
\! := \!
\begin{pmatrix}
\boxed{\begin{matrix} -\varkappa_{1}^{o}(a_{j}^{o}) \varkappa_{2}^{o}(a_{j}^{
o}) \! \left(s_{1} \! \left\{\frac{Q_{0}^{o}(a_{j}^{o})}{(\gamma^{o}(0))^{2}} 
\! \left[\gimel^{1}_{-1}(a_{j}^{o}) \right. \right. \right. \\
\left. \left. \left. + \, \gimel^{1}_{1}(a_{j}^{o}) \! + \! \daleth^{1}_{1}
(a_{j}^{o}) \daleth^{1}_{-1}(a_{j}^{o}) \right] \! + \! \frac{Q_{1}^{o}(a_{
j}^{o})}{(\gamma^{o}(0))^{2}} \right. \right. \\
\left. \left. \times \left[\daleth^{1}_{1}(a_{j}^{o}) \! + \! \daleth^{1}_{-1}
(a_{j}^{o}) \right] \! + \! \frac{1}{2} \frac{Q_{2}^{o}(a_{j}^{o})}{(\gamma^{
o}(0))^{2}} \right. \right. \\
\left. \left. + \, \frac{(\gamma^{o}(0))^{2}}{Q_{0}^{o}(a_{j}^{o})} \aleph^{
1}_{1}(a_{j}^{o}) \aleph^{1}_{-1}(a_{j}^{o}) \right\} \right. \\
\left. + \, t_{1} \! \left\{\frac{Q_{0}^{o}(a_{j}^{o})}{(\gamma^{o}(0))^{2}} 
\! \left[\aleph^{1}_{-1}(a_{j}^{o}) \beth^{1}_{1}(a_{j}^{o}) \right. \right. 
\right. \\
\left. \left. \left. + \, \aleph^{1}_{1}(a_{j}^{o}) \beth^{1}_{-1}(a_{j}^{o}) 
\right] \! + \! \frac{(\gamma^{o}(0))^{2}}{Q_{0}^{o}(a_{j}^{o})} \right. 
\right. \\
\left. \left. \times \left[\daleth^{1}_{-1}(a_{j}^{o}) \! + \! \daleth^{1}_{1}
(a_{j}^{o}) \! - \! \frac{Q_{1}^{o}(a_{j}^{o})}{Q_{0}^{o}(a_{j}^{o})} \right] 
\right. \right. \\
\left. \left. + \, \frac{Q_{1}^{o}(a_{j}^{o})}{(\gamma^{o}(0))^{2}} \aleph^{
1}_{1}(a_{j}^{o}) \aleph^{1}_{-1}(a_{j}^{o}) \right\} \right. \\
\left. + \, \mi (s_{1} \! + \! t_{1}) \! \left\{\beth^{1}_{1}(a_{j}^{o}) \! 
+ \! \aleph^{1}_{1}(a_{j}^{o}) \daleth^{1}_{-1}(a_{j}^{o}) \right. \right. \\
\left. \left. - \, \beth^{1}_{-1}(a_{j}^{o}) \! - \! \aleph^{1}_{-1}(a_{j}^{
o}) \daleth^{1}_{1}(a_{j}^{o}) \right\} \right)
\end{matrix}} & \boxed{\begin{matrix}
(\varkappa_{1}^{o}(a_{j}^{o}))^{2} \! \left(\mi s_{1} \! \left\{\frac{Q_{0}^{
o}(a_{j}^{o})}{(\gamma^{o}(0))^{2}} \! \left[2 \gimel^{1}_{1}(a_{j}^{o}) 
\right. \right. \right. \\
\left. \left. \left. + \, (\daleth^{1}_{1}(a_{j}^{o}))^{2} \right] \! + \! 
\frac{2Q_{1}^{o}(a_{j}^{o})}{(\gamma^{o}(0))^{2}} \daleth^{1}_{1}(a_{j}^{o}) 
\right. \right. \\
\left. \left. - \, \frac{(\gamma^{o}(0))^{2}}{Q_{0}^{o}(a_{j}^{o})}(\aleph^{
1}_{1}(a_{j}^{o}))^{2} \! + \! \frac{1}{2} \frac{Q_{2}^{o}(a_{j}^{o})}{
(\gamma^{o}(0))^{2}} \right\} \right. \\
\left. + \, \mi t_{1} \! \left\{\frac{2Q_{0}^{o}(a_{j}^{o})}{(\gamma^{o}(0)
)^{2}} \aleph^{1}_{1}(a_{j}^{o}) \beth^{1}_{1}(a_{j}^{o}) \right. \right. \\
\left. \left. + \, \frac{Q_{1}^{o}(a_{j}^{o})}{(\gamma^{o}(0))^{2}}(\aleph^{
1}_{1}(a_{j}^{o}))^{2} \! + \! \frac{(\gamma^{o}(0))^{2}}{Q_{0}^{o}(a_{j}^{
o})} \right. \right. \\
\left. \left. \times \left[\frac{Q_{1}^{o}(a_{j}^{o})}{Q_{0}^{o}(a_{j}^{
o})} \! - \! 2 \daleth^{1}_{1}(a_{j}^{o}) \right] \right\} \! - \! 2(s_{1} 
\! - \! t_{1}) \right. \\
\left. \times \left\{\beth^{1}_{1}(a_{j}^{o}) \! + \! \aleph^{1}_{1}(a_{j}^{
o}) \daleth^{1}_{1}(a_{j}^{o}) \right\} \right)
\end{matrix}} \\
\boxed{\begin{matrix} (\varkappa_{2}^{o}(a_{j}^{o}))^{2} \! \left(\mi s_{1} 
\! \left\{\frac{Q_{0}^{o}(a_{j}^{o})}{(\gamma^{o}(0))^{2}} \! \left[2 
\gimel^{1}_{-1}(a_{j}^{o}) \right. \right. \right. \\
\left. \left. \left. + \, (\daleth^{1}_{-1}(a_{j}^{o}))^{2} \right] \! + \! 
\frac{2Q_{1}^{o}(a_{j}^{o})}{(\gamma^{o}(0))^{2}} \daleth^{1}_{-1}(a_{j}^{o}) 
\right. \right. \\
\left. \left. - \, \frac{(\gamma^{o}(0))^{2}}{Q_{0}^{o}(a_{j}^{o})}(\aleph^{
1}_{-1}(a_{j}^{o}))^{2} \! + \! \frac{1}{2} \frac{Q_{2}^{o}(a_{j}^{o})}{
(\gamma^{o}(0))^{2}} \right\} \right. \\
\left. + \, \mi t_{1} \! \left\{\frac{2Q_{0}^{o}(a_{j}^{o})}{(\gamma^{o}(0)
)^{2}} \aleph^{1}_{-1}(a_{j}^{o}) \beth^{1}_{-1}(a_{j}^{o}) \right. \right. \\
\left. \left. + \, \frac{Q_{1}^{o}(a_{j}^{o})}{(\gamma^{o}(0))^{2}}(\aleph^{
1}_{-1}(a_{j}^{o}))^{2} \! + \! \frac{(\gamma^{o}(0))^{2}}{Q_{0}^{o}(a_{j}^{
o})} \right. \right. \\
\left. \left. \times \left[\frac{Q_{1}^{o}(a_{j}^{o})}{Q_{0}^{o}(a_{j}^{
o})} \! - \! 2 \daleth^{1}_{-1}(a_{j}^{o}) \right] \right\} \! + \! 2
(s_{1} \! - \! t_{1}) \right. \\
\left. \times \left\{\beth^{1}_{-1}(a_{j}^{o}) \! + \! \aleph^{1}_{-1}
(a_{j}^{o}) \daleth^{1}_{-1}(a_{j}^{o}) \right\} \right)
\end{matrix}} & 
\boxed{\begin{matrix} \varkappa_{1}^{o}(a_{j}^{o}) \varkappa_{2}^{o}(a_{j}^{
o}) \! \left(s_{1} \! \left\{\frac{Q_{0}^{o}(a_{j}^{o})}{(\gamma^{o}(0))^{2}} 
\! \left[\gimel^{1}_{-1}(a_{j}^{o}) \right. \right. \right. \\
\left. \left. \left. + \, \gimel^{1}_{1}(a_{j}^{o}) \! + \! \daleth^{1}_{1}
(a_{j}^{o}) \daleth^{1}_{-1}(a_{j}^{o}) \right] \! + \! \frac{Q_{1}^{o}(a_{
j}^{o})}{(\gamma^{o}(0))^{2}} \right. \right. \\
\left. \left. \times \left[\daleth^{1}_{1}(a_{j}^{o}) \! + \! \daleth^{1}_{-1}
(a_{j}^{o}) \right] \! + \! \frac{1}{2} \frac{Q_{2}^{o}(a_{j}^{o})}{(\gamma^{
o}(0))^{2}} \right. \right. \\
\left. \left. + \, \frac{(\gamma^{o}(0))^{2}}{Q_{0}^{o}(a_{j}^{o})} \aleph^{
1}_{1}(a_{j}^{o}) \aleph^{1}_{-1}(a_{j}^{o}) \right\} \right. \\
\left. + \, t_{1} \! \left\{\frac{Q_{0}^{o}(a_{j}^{o})}{(\gamma^{o}(0))^{2}} 
\! \left[\aleph^{1}_{-1}(a_{j}^{o}) \beth^{1}_{1}(a_{j}^{o}) \right. \right. 
\right. \\
\left. \left. \left. + \, \aleph^{1}_{1}(a_{j}^{o}) \beth^{1}_{-1}(a_{j}^{o}) 
\right] \! + \! \frac{(\gamma^{o}(0))^{2}}{Q_{0}^{o}(a_{j}^{o})} \right. 
\right. \\
\left. \left. \times \left[\daleth^{1}_{-1}(a_{j}^{o}) \! + \! \daleth^{1}_{1}
(a_{j}^{o}) \! - \! \frac{Q_{1}^{o}(a_{j}^{o})}{Q_{0}^{o}(a_{j}^{o})} \right] 
\right. \right. \\
\left. \left. + \, \frac{Q_{1}^{o}(a_{j}^{o})}{(\gamma^{o}(0))^{2}} \aleph^{
1}_{1}(a_{j}^{o}) \aleph^{1}_{-1}(a_{j}^{o}) \right\} \right. \\
\left. + \, \mi (s_{1} \! + \! t_{1}) \! \left\{\beth^{1}_{1}(a_{j}^{o}) \! 
+ \! \aleph^{1}_{1}(a_{j}^{o}) \daleth^{1}_{-1}(a_{j}^{o}) \right. \right. \\
\left. \left. - \, \beth^{1}_{-1}(a_{j}^{o}) \! - \! \aleph^{1}_{-1}(a_{j}^{
o}) \daleth^{1}_{1}(a_{j}^{o}) \right\} \right)
\end{matrix}}
\end{pmatrix},
\end{equation*}
(with $\operatorname{tr}(\mathscr{C}^{o}(a_{j}^{o})) \! = \! 0)$, and $(f^{
o}_{k}(n))_{ij} \! =_{n \to \infty} \! \mathcal{O}(1)$, $k \! \in \! \mathbb{
N}$, $i,j \! = \! 1,2$. (The expression for $\mathscr{C}^{o}(a_{j}^{o})$ is
necessary for obtaining asymptotics \textbf{at} the end-points $\lbrace a_{
j}^{o} \rbrace_{j=1}^{N}$, as well as for Remark~5.2 below.) Returning to the
counter-clockwise-oriented integrals $\oint_{\partial \mathbb{U}^{o}_{\delta_{
a_{j}}}} \tfrac{w_{+}^{\Sigma^{o}_{\circlearrowright}}(s)}{s-z} \, \tfrac{\md
s}{2 \pi \mi}$, $z \! \in \! \mathbb{C} \setminus \widetilde{\Sigma}_{p}^{o}$,
it follows, via the Residue and Cauchy Theorems, that, for $j \! = \! 1,\dotsc,
N$,
\begin{equation*}
\oint_{\partial \mathbb{U}^{o}_{\delta_{a_{j}}}} \dfrac{w_{+}^{\Sigma^{o}_{
\circlearrowright}}(s)}{s \! - \! z} \, \dfrac{\md s}{2 \pi \mi} \underset{
n \to \infty}{=}
\begin{cases}
-\dfrac{\widehat{\mathscr{A}}^{o}(a_{j}^{o})}{(n \! + \! \frac{1}{2})(z \! -
\! a_{j}^{o})^{2}} \! - \! \dfrac{\widehat{\mathscr{B}}^{o}(a_{j}^{o})}{(n \!
+ \! \frac{1}{2})(z \! - \! a_{j}^{o})} \! + \! \mathcal{O} \! \left(\dfrac{
\widehat{f}^{o}(z;n)}{(n \! + \! \frac{1}{2})^{2}} \right),
&\text{$z \! \in \! \mathbb{U}^{o,\ast}_{\delta_{a_{j}}}$,} \\
-\dfrac{\widehat{\mathscr{A}}^{o}(a_{j}^{o})}{(n \! + \! \frac{1}{2})(z \! -
\! a_{j}^{o})^{2}} \! - \! \dfrac{\widehat{\mathscr{B}}^{o}(a_{j}^{o})}{(n \!
+ \! \frac{1}{2})(z \! - \! a_{j}^{o})} \! + \! \dfrac{\mathscr{R}^{0}_{a_{
j}^{o}}(z)}{n \! + \! \frac{1}{2}} \! + \! \mathcal{O} \! \left(\dfrac{
\widehat{f}^{o}(z;n)}{(n \! + \! \frac{1}{2})^{2}} \right), &\text{$z \! \in
\! \mathbb{U}^{o}_{\delta_{a_{j}}}$,}
\end{cases}
\end{equation*}
where $\mathbb{U}^{o,\ast}_{\delta_{a_{j}}} \! := \! \mathbb{C} \setminus 
(\mathbb{U}^{o}_{\delta_{a_{j}}} \cup \partial \mathbb{U}^{o}_{\delta_{a_{j}}
})$, $\widehat{\mathscr{A}}^{o}(a_{j}^{o}) \! := \! (\widehat{\alpha}_{0}^{o}
(a_{j}^{o}))^{-1} \mathscr{A}^{o}(a_{j}^{o})$, $\widehat{\mathscr{B}}^{o}(a_{
j}^{o}) \! := \! (\widehat{\alpha}_{0}^{o}(a_{j}^{o}))^{-2}(\mathscr{B}^{o}
(a_{j}^{o}) \widehat{\alpha}_{0}^{o}(a_{j}^{o}) \! - \! \mathscr{A}^{o}(a_{
j}^{o}) \widehat{\alpha}_{1}^{o}(a_{j}^{o}))$, $\mathscr{R}^{o}_{a_{j}^{o}}
(z)$ is given in Theorem~2.3.1, Equations~(2.73) and~(2.74), and $\widehat{
f}^{o}(z;n)$, where the $n$-dependence arises due to the $n$-dependence of 
the associated Riemann theta functions, denote some bounded (with respect to 
both $z$ and $n)$, analytic (for $\mathbb{C} \setminus \widetilde{\Sigma}^{
o}_{p} \! \ni \! z)$, $\operatorname{GL}_{2}(\mathbb{C})$-valued functions 
for which $(\widehat{f}^{o}(z;n))_{kl} \! =_{\underset{z \in \mathbb{C} 
\setminus \widetilde{\Sigma}^{o}_{p}}{n \to \infty}} \! \mathcal{O}(1)$, $k,
l \! = \! 1,2$. Similarly, one shows that, for $j \! = \! 1,\dotsc,N$,
\begin{equation*}
\oint_{\partial \mathbb{U}^{o}_{\delta_{a_{j}}}}s^{-1}w_{+}^{\Sigma^{o}_{
\circlearrowright}}(s) \, \dfrac{\md s}{2 \pi \mi} \underset{n \to \infty}{=}
\dfrac{(\mathscr{B}^{o}(a_{j}^{o}) \widehat{\alpha}^{o}_{0}(a_{j}^{o}) \! - \!
\mathscr{A}^{o}(a_{j}^{o})(\widehat{\alpha}^{o}_{1}(a_{j}^{o}) \! + \! (a_{
j}^{o})^{-1} \widehat{\alpha}^{o}_{0}(a_{j}^{o})))}{(n \! + \! \frac{1}{2})
(\widehat{\alpha}^{o}_{0}(a_{j}^{o}))^{2}a_{j}^{o}} \! + \! \mathcal{O} \!
\left(\dfrac{f_{j}(n)}{(n \! + \! \frac{1}{2})^{2}} \right),
\end{equation*}
where $(f_{j}(n))_{kl} \! =_{n \to \infty} \! \mathcal{O}(1)$, $k,l \! = \! 1,
2$. Repeating the above analysis for the remaining end-points of the support
of the `odd' equilibrium measure, that is, $\lbrace b_{0}^{o},\dotsc,b_{N}^{
o},a_{N+1}^{o} \rbrace$, one arrives at the result stated in the Lemma.
\hfill $\qed$
\begin{eeeee}
A brisk perusing of the asymptotic (as $n \! \to \! \infty)$ result for 
$\mathscr{R}^{o}(z)$ stated in Lemma~5.3 seems to imply that, at first glance,
there are second-order poles at $\lbrace b_{j-1}^{o},a_{j}^{o} \rbrace_{j=1}^{
N+1}$; however, this is not the case. As the proof of Lemma~5.3 demonstrates 
(cf. the analysis leading up to Equations~(5.3)), Laurent series expansions 
about $\lbrace b_{j-1}^{o},a_{j}^{o} \rbrace_{j=1}^{N+1}$ show that, as $n \! 
\to \! \infty$, all expansions are, indeed, analytic; in particular: (i) for 
$z \! \in \! \mathbb{U}^{o}_{\delta_{a_{j}}}$, $j \! = \! 1,\dotsc,N \! + \! 
1$ (all contour integrals are counter-clockwise oriented),
\begin{align*}
\oint_{\partial \mathbb{U}^{o}_{\delta_{a_{j}}}} \dfrac{w_{+}^{\Sigma^{o}_{
\circlearrowright}}(s)}{s \! - \! z} \, \dfrac{\md s}{2 \pi \mi} \! \underset{
\underset{z \in \mathbb{U}^{o}_{\delta_{a_{j}}}}{n \to \infty}}{=}& \dfrac{
\left(\mathscr{A}^{o}(a_{j}^{o}) \widehat{\alpha}_{0}^{o}(a_{j}^{o}) \! \left(
\! \left(\frac{\widehat{\alpha}_{1}^{o}(a_{j}^{o})}{\widehat{\alpha}_{0}^{o}
(a_{j}^{o})} \right)^{2} \! - \! \frac{\widehat{\alpha}_{2}^{o}(a_{j}^{o})}{
\widehat{\alpha}_{0}^{o}(a_{j}^{o})} \right) \! - \! \mathscr{B}^{o}(a_{j}^{
o}) \widehat{\alpha}_{1}^{o}(a_{j}^{o}) \! + \! \mathscr{C}^{o}(a_{j}^{o})
\widehat{\alpha}_{0}^{o}(a_{j}^{o}) \right)}{(n \! + \! \frac{1}{2})(\widehat{
\alpha}_{0}^{o}(a_{j}^{o}))^{2}} \\
&+ \, \dfrac{1}{(n \! + \! \frac{1}{2})} \sum_{k=1}^{\infty}f_{k}^{a_{j}^{o}}
(n)(z \! - \! a_{j}^{o})^{k} \! + \! \mathcal{O} \! \left(\dfrac{\widehat{f}^{
o}(z;n)}{(n \! + \! \frac{1}{2})^{2}} \right),
\end{align*}
where, for $j \! = \! 1,\dotsc,N$, $\mathscr{A}^{o}(a_{j}^{o})$, $\mathscr{
B}^{o}(a_{j}^{o})$, $\mathscr{C}^{o}(a_{j}^{o})$, $\widehat{\alpha}_{0}^{o}
(a_{j}^{o})$, $\widehat{\alpha}_{1}^{o}(a_{j}^{o})$, and $\widehat{\alpha}_{
2}^{o}(a_{j}^{o})$ are given in (the proof of) Lemma~5.3, $\mathscr{A}^{o}
(a_{N+1}^{o})$, $\mathscr{B}^{o}(a_{N+1}^{o})$ are given in Theorem~2.3.1, 
Equations~(2.25), (2.27), (2.28), (2.31), (2.32), (2.37)--(2.45), (2.47), 
(2.52), and~(2.53), $\widehat{\alpha}_{0}^{o}(a_{N+1}^{o}) \! := \! \tfrac{
4}{3}f(a_{N+1}^{o})$, $\widehat{\alpha}_{1}^{o}(a_{N+1}^{o}) \! := \! 
\tfrac{4}{5}f^{\prime}(a_{N+1}^{o})$, and $\widehat{\alpha}_{2}^{o}(a_{N+1}^{
o}) \! := \! \tfrac{2}{7}f^{\prime \prime}(a_{N+1}^{o})$, with $f(a_{N+1}^{
o})$, $f^{\prime}(a_{N+1}^{o})$, and $f^{\prime \prime}(a_{N+1}^{o})$ given 
in Lemma~4.7, $\mathscr{C}^{o}(a_{N+1}^{o})$ is given by the same expression 
as $\mathscr{C}^{o}(a_{j}^{o})$ above subject to the modifications $\Omega_{
j}^{o} \! \to \! 0$, $a_{j}^{o} \! \to \! a_{N+1}^{o}$, $Q_{0}^{o}(a_{j}^{o}) 
\! \to \! Q_{0}^{o}(a_{N+1}^{o})$, $Q_{1}^{o}(a_{j}^{o}) \! \to \! Q_{1}^{o}
(a_{N+1}^{o})$, with $Q_{0}^{o}(a_{N+1}^{o})$, $Q_{1}^{o}(a_{N+1}^{o})$ given 
in Theorem~2.3.1, Equations~(2.31) and~(2.32), and $Q_{2}^{o}(a_{j}^{o}) \! 
\to \! Q_{2}^{o}(a_{N+1}^{o})$, where
\begin{align*}
Q_{2}^{o}(a_{N+1}^{o})&= -\dfrac{1}{2}Q_{0}^{o}(a_{N+1}^{o}) \! \left(\sum_{k
=1}^{N} \! \left(\dfrac{1}{(a_{N+1}^{o} \! - \! b_{k}^{o})^{2}} \! - \! \dfrac{
1}{(a_{N+1}^{o} \! - \! a_{k}^{o})^{2}} \right) \! + \! \dfrac{1}{(a_{N+1}^{o}
\! - \! b_{0}^{o})^{2}} \right) \\
&+ \, \dfrac{1}{4}Q_{0}^{o}(a_{N+1}^{o}) \! \left(\sum_{k=1}^{N} \! \left(
\dfrac{1}{a_{N+1}^{o} \! - \! b_{k}^{o}} \! - \! \dfrac{1}{a_{N+1}^{o} \! - \!
a_{k}^{o}} \right) \! + \! \dfrac{1}{a_{N+1}^{o} \! - \! b_{0}^{o}} \right)^{
2},
\end{align*}
$\widehat{f}^{o}(z;n)$ is characterised completely at the end of the proof 
of Lemma~5.3, $(f_{k}^{a_{j}^{o}}(n))_{l_{1}l_{2}} \! =_{n \to \infty} \! 
\mathcal{O}(1)$, $k \! \in \! \mathbb{N}$, $l_{1},l_{2} \! = \! 1,2$, and
\begin{align*}
\oint_{\partial \mathbb{U}^{o}_{\delta_{a_{N+1}}}}s^{-1}w_{+}^{\Sigma^{o}_{
\circlearrowright}}(s) \, \dfrac{\md s}{2 \pi \mi} \underset{n \to \infty}{=}&
\, \dfrac{(\mathscr{B}^{o}(a_{N+1}^{o}) \widehat{\alpha}^{o}_{0}(a_{N+1}^{o}
) \! - \! \mathscr{A}^{o}(a_{N+1}^{o})(\widehat{\alpha}^{o}_{1}(a_{N+1}^{o})
\! + \! (a_{N+1}^{o})^{-1} \widehat{\alpha}^{o}_{0}(a_{N+1}^{o})))}{(n \! +
\! \frac{1}{2})(\widehat{\alpha}^{o}_{0}(a_{N+1}^{o}))^{2}a_{N+1}^{o}} \\
&+ \, \mathcal{O} \! \left(\dfrac{f(n)}{(n \! + \! \frac{1}{2})^{2}} \right),
\end{align*}
where $(f(n))_{kl} \! =_{n \to \infty} \! \mathcal{O}(1)$, $k,l \! = \! 1,2$: 
and (ii) for $z \! \in \! \mathbb{U}^{o}_{\delta_{b_{j}}}$, $j \! = \! 0,
\dotsc,N$,
\begin{align*}
\oint_{\partial \mathbb{U}^{o}_{\delta_{b_{j}}}} \dfrac{w_{+}^{\Sigma^{o}_{
\circlearrowright}}(s)}{s \! - \! z} \, \dfrac{\md s}{2 \pi \mi} \! \underset{
\underset{z \in \mathbb{U}^{o}_{\delta_{b_{j}}}}{n \to \infty}}{=}& \dfrac{
\left(\mathscr{A}^{o}(b_{j}^{o}) \widehat{\alpha}_{0}^{o}(b_{j}^{o}) \! \left(
\! \left(\frac{\widehat{\alpha}_{1}^{o}(b_{j}^{o})}{\widehat{\alpha}_{0}^{o}
(b_{j}^{o})} \right)^{2} \! - \! \frac{\widehat{\alpha}_{2}^{o}(b_{j}^{o})}{
\widehat{\alpha}_{0}^{o}(b_{j}^{o})} \right) \! - \! \mathscr{B}^{o}(b_{j}^{o}
) \widehat{\alpha}_{1}^{o}(b_{j}^{o}) \! + \! \mathscr{C}^{o}(b_{j}^{o})
\widehat{\alpha}_{0}^{o}(b_{j}^{o}) \right)}{(n \! + \! \frac{1}{2})(
\widehat{\alpha}_{0}^{o}(b_{j}^{o}))^{2}} \\
&+ \, \dfrac{1}{(n \! + \! \frac{1}{2})} \sum_{k=1}^{\infty}f_{k}^{b_{j}^{o}}
(n)(z \! - \! b_{j}^{o})^{k} \! + \! \mathcal{O} \! \left(\dfrac{\widehat{f}^{
o}(z;n)}{(n \! + \! \frac{1}{2})^{2}} \right),
\end{align*}
where, for $j \! = \! 1,\dotsc,N \! + \! 1$, $\mathscr{A}^{o}(b_{j-1}^{o})$, 
$\mathscr{B}^{o}(b_{j-1}^{o})$ are given in Theorem~2.3.1, Equations~(2.24), 
(2.26), (2.28), (2.29), (2.30), (2.33), (2.34), (2.37)--(2.45), (2.46), 
(2.48), (2.50), (2.51), (2.54), and~(2.55), $\widehat{\alpha}_{0}^{o}(b_{j-
1}^{o}) \! := \! \tfrac{4}{3}f(b_{j-1}^{o})$, $\widehat{\alpha}_{1}^{o}(b_{j-
1}^{o}) \! := \! \tfrac{4}{5}f^{\prime}(b_{j-1}^{o})$, and $\widehat{\alpha}_{
2}^{o}(b_{j-1}^{o}) \! := \! \tfrac{2}{7}f^{\prime \prime}(b_{j-1}^{o})$, with 
$f(b_{j-1}^{o})$, $f^{\prime}(b_{j-1}^{o})$, and $f^{\prime \prime}(b_{j-1}^{
o})$ given in Lemma~4.6, for $j \! = \! 1,\dotsc,N$,
\begin{equation*}
\dfrac{\mathscr{C}^{o}(b_{j}^{o})}{\me^{\mi (n+\frac{1}{2}) \Omega_{j}^{o}}} 
\! := \! 
\begin{pmatrix}
\boxed{\begin{matrix} \varkappa_{1}^{o}(b_{j}^{o}) \varkappa_{2}^{o}(b_{j}^{
o}) \! \left(s_{1} \! \left\{-\frac{Q_{0}^{o}(b_{j}^{o})}{(\gamma^{o}(0))^{
2}} \aleph^{1}_{1}(b_{j}^{o}) \right. \right. \\
\left. \left. \times \, \aleph^{1}_{-1}(b_{j}^{o}) \! - \! \frac{(Q_{1}^{o}
(b_{j}^{o}) \gamma^{o}(0))^{2}}{(Q_{0}^{o}(b_{j}^{o}))^{3}} \! + \! \frac{
1}{2} \frac{Q_{2}^{o}(b_{j}^{o})(\gamma^{o}(0))^{2}}{(Q_{0}^{o}(b_{j}^{o}))^{
2}} \right. \right. \\
\left. \left. + \, \frac{Q_{1}^{o}(b_{j}^{o})(\gamma^{o}(0))^{2}}{(Q_{0}^{o}
(b_{j}^{o}))^{2}} \! \left[\daleth^{1}_{1}(b_{j}^{o}) \! + \! \daleth^{1}_{
-1}(b_{j}^{o}) \right] \! - \! \frac{(\gamma^{o}(0))^{2}}{Q_{0}^{o}(b_{j}^{
o})} \right. \right. \\
\left. \left. \times \left[\gimel^{1}_{1}(b_{j}^{o}) \! + \! \gimel^{1}_{-1}
(b_{j}^{o}) \! + \! \daleth^{1}_{1}(b_{j}^{o}) \daleth^{1}_{-1}(b_{j}^{o}) 
\right] \right\} \right. \\
\left. + \, t_{1} \! \left\{-\frac{Q_{1}^{o}(b_{j}^{o})}{(\gamma^{o}(0))^{2}} 
\! - \! \frac{Q_{0}^{o}(b_{j}^{o})}{(\gamma^{o}(0))^{2}} \! \left[\daleth^{
1}_{-1}(b_{j}^{o}) \! + \! \daleth^{1}_{1}(b_{j}^{o}) \right] \right. \right. 
\\
\left. \left. + \, \frac{Q_{1}^{o}(b_{j}^{o})(\gamma^{o}(0))^{2}}{(Q_{0}^{o}
(b_{j}^{o}))^{2}} \aleph^{1}_{1}(b_{j}^{o}) \aleph^{1}_{-1}(b_{j}^{o}) \! 
- \! \frac{(\gamma^{o}(0))^{2}}{Q_{0}^{o}(b_{j}^{o})} \right. \right. \\
\left. \left. \times \left[\aleph^{1}_{1}(b_{j}^{o}) \beth^{1}_{-1}(b_{j}^{o}) 
\! + \! \aleph^{1}_{-1}(b_{j}^{o}) \beth^{1}_{1}(b_{j}^{o}) \right] \right\} 
\right. \\
\left. + \, \mi (s_{1} \! + \! t_{1}) \! \left[\beth^{1}_{-1}(b_{j}^{o}) \! 
- \! \aleph^{1}_{1}(b_{j}^{o}) \daleth^{1}_{-1}(b_{j}^{o}) \right. \right. \\
\left. \left. + \, \aleph^{1}_{-1}(b_{j}^{o}) \daleth^{1}_{1}(b_{j}^{o}) \! 
- \! \beth^{1}_{1}(b_{j}^{o}) \right] \right)
\end{matrix}} & 
\boxed{\begin{matrix} (\varkappa_{1}^{o}(b_{j}^{o}))^{2} \! \left(\mi s_{1} 
\! \left\{\frac{Q_{0}^{o}(b_{j}^{o})}{(\gamma^{o}(0))^{2}}(\aleph^{1}_{1}
(b_{j}^{o}))^{2} \right. \right. \\
\left. \left. - \, \frac{(Q_{1}^{o}(b_{j}^{o}) \gamma^{o}(0))^{2}}{(Q_{0}^{o}
(b_{j}^{o}))^{3}} \! + \! \frac{1}{2} \frac{Q_{2}^{o}(b_{j}^{o})(\gamma^{o}
(0))^{2}}{(Q_{0}^{o}(b_{j}^{o}))^{2}} \right. \right. \\
\left. \left. + \, \frac{2Q_{1}^{o}(b_{j}^{o})(\gamma^{o}(0))^{2}}{(Q_{0}^{o}
(b_{j}^{o}))^{2}} \daleth^{1}_{1}(b_{j}^{o}) \! - \! \frac{(\gamma^{o}(0))^{
2}}{Q_{0}^{o}(b_{j}^{o})} \right. \right. \\
\left. \left. \times \left[2 \gimel^{1}_{1}(b_{j}^{o}) \! + \! (\daleth^{
1}_{1}(b_{j}^{o}))^{2} \right] \right\} \! + \! \mi t_{1} \right. \\
\left. \times \left\{\frac{2Q_{0}^{o}(b_{j}^{o})}{(\gamma^{o}(0))^{2}} 
\daleth^{1}_{1}(b_{j}^{o}) \! + \! \frac{Q_{1}^{o}(b_{j}^{o})}{(\gamma^{o}
(0))^{2}} \right. \right. \\
\left. \left. + \, \frac{Q_{1}^{o}(b_{j}^{o})(\gamma^{o}(0))^{2}}{(Q_{0}^{o}
(b_{j}^{o}))^{2}}(\aleph^{1}_{1}(b_{j}^{o}))^{2} \! - \! \frac{2(\gamma^{o}
(0))^{2}}{Q_{0}^{o}(b_{j}^{o})} \right. \right. \\
\left. \left. \times \, \aleph^{1}_{1}(b_{j}^{o}) \beth^{1}_{1}(b_{j}^{o}) 
\right\} \! + \! 2(s_{1} \! - \! t_{1}) \right. \\
\left. \times \left\{\beth^{1}_{1}(b_{j}^{o}) \! + \! \aleph^{1}_{1}(b_{j}^{
o}) \daleth^{1}_{1}(b_{j}^{o}) \right\} \right)
\end{matrix}} \\
\boxed{\begin{matrix} (\varkappa_{2}^{o}(b_{j}^{o}))^{2} \! \left(\mi s_{1} 
\! \left\{\frac{Q_{0}^{o}(b_{j}^{o})}{(\gamma^{o}(0))^{2}}(\aleph^{1}_{-1}
(b_{j}^{o}))^{2} \right. \right. \\
\left. \left. - \, \frac{(Q_{1}^{o}(b_{j}^{o}) \gamma^{o}(0))^{2}}{(Q_{0}^{o}
(b_{j}^{o}))^{3}} \! + \! \frac{1}{2} \frac{Q_{2}^{o}(b_{j}^{o})(\gamma^{o}
(0))^{2}}{(Q_{0}^{o}(b_{j}^{o}))^{2}} \right. \right. \\
\left. \left. + \, \frac{2Q_{1}^{o}(b_{j}^{o})(\gamma^{o}(0))^{2}}{(Q_{0}^{o}
(b_{j}^{o}))^{2}} \daleth^{1}_{-1}(b_{j}^{o}) \! - \! \frac{(\gamma^{o}(0))^{
2}}{Q_{0}^{o}(b_{j}^{o})} \right. \right. \\
\left. \left. \times \left[2 \gimel^{1}_{-1}(b_{j}^{o}) \! + \! (\daleth^{
1}_{-1}(b_{j}^{o}))^{2} \right] \right\} \! + \! \mi t_{1} \right. \\
\left. \times \left\{\frac{2Q_{0}^{o}(b_{j}^{o})}{(\gamma^{o}(0))^{2}} 
\daleth^{1}_{-1}(b_{j}^{o}) \! + \! \frac{Q_{1}^{o}(b_{j}^{o})}{(\gamma^{o}
(0))^{2}} \right. \right. \\
\left. \left. + \, \frac{Q_{1}^{o}(b_{j}^{o})(\gamma^{o}(0))^{2}}{(Q_{0}^{o}
(b_{j}^{o}))^{2}}(\aleph^{1}_{-1}(b_{j}^{o}))^{2} \! - \! \frac{2(\gamma^{o}
(0))^{2}}{Q_{0}^{o}(b_{j}^{o})} \right. \right. \\
\left. \left. \times \, \aleph^{1}_{-1}(b_{j}^{o}) \beth^{1}_{-1}(b_{j}^{o}) 
\right\} \! - \! 2(s_{1} \! - \! t_{1}) \right. \\
\left. \times \left\{\beth^{1}_{-1}(b_{j}^{o}) \! + \! \aleph^{1}_{-1}
(b_{j}^{o}) \daleth^{1}_{-1}(b_{j}^{o}) \right\} \right)
\end{matrix}} & 
\boxed{\begin{matrix} -\varkappa_{1}^{o}(b_{j}^{o}) \varkappa_{2}^{o}
(b_{j}^{o}) \! \left(s_{1} \! \left\{-\frac{Q_{0}^{o}(b_{j}^{o})}{(\gamma^{
o}(0))^{2}} \aleph^{1}_{1}(b_{j}^{o}) \right. \right. \\
\left. \left. \times \, \aleph^{1}_{-1}(b_{j}^{o}) \! - \! \frac{(Q_{1}^{o}
(b_{j}^{o}) \gamma^{o}(0))^{2}}{(Q_{0}^{o}(b_{j}^{o}))^{3}} \! + \! \frac{
1}{2} \frac{Q_{2}^{o}(b_{j}^{o})(\gamma^{o}(0))^{2}}{(Q_{0}^{o}(b_{j}^{o}
))^{2}} \right. \right. \\
\left. \left. + \, \frac{Q_{1}^{o}(b_{j}^{o})(\gamma^{o}(0))^{2}}{(Q_{0}^{o}
(b_{j}^{o}))^{2}} \! \left[\daleth^{1}_{1}(b_{j}^{o}) \! + \! \daleth^{1}_{
-1}(b_{j}^{o}) \right] \! - \! \frac{(\gamma^{o}(0))^{2}}{Q_{0}^{o}(b_{j}^{
o})} \right. \right. \\
\left. \left. \times \left[\gimel^{1}_{1}(b_{j}^{o}) \! + \! \gimel^{1}_{-1}
(b_{j}^{o}) \! + \! \daleth^{1}_{1}(b_{j}^{o}) \daleth^{1}_{-1}(b_{j}^{o}) 
\right] \right\} \right. \\
\left. + \, t_{1} \! \left\{-\frac{Q_{1}^{o}(b_{j}^{o})}{(\gamma^{o}(0))^{2}} 
\! - \! \frac{Q_{0}^{o}(b_{j}^{o})}{(\gamma^{o}(0))^{2}} \! \left[\daleth^{
1}_{-1}(b_{j}^{o}) \! + \! \daleth^{1}_{1}(b_{j}^{o}) \right] \right. \right. 
\\
\left. \left. + \, \frac{Q_{1}^{o}(b_{j}^{o})(\gamma^{o}(0))^{2}}{(Q_{0}^{o}
(b_{j}^{o}))^{2}} \aleph^{1}_{1}(b_{j}^{o}) \aleph^{1}_{-1}(b_{j}^{o}) \! 
- \! \frac{(\gamma^{o}(0))^{2}}{Q_{0}^{o}(b_{j}^{o})} \right. \right. \\
\left. \left. \times \left[\aleph^{1}_{1}(b_{j}^{o}) \beth^{1}_{-1}(b_{j}^{o}) 
\! + \! \aleph^{1}_{-1}(b_{j}^{o}) \beth^{1}_{1}(b_{j}^{o}) \right] \right\} 
\right. \\
\left. + \, \mi (s_{1} \! + \! t_{1}) \! \left[\beth^{1}_{-1}(b_{j}^{o}) \! 
- \! \aleph^{1}_{1}(b_{j}^{o}) \daleth^{1}_{-1}(b_{j}^{o}) \right. \right. \\
\left. \left. + \, \aleph^{1}_{-1}(b_{j}^{o}) \daleth^{1}_{1}(b_{j}^{o}) \! - 
\! \beth^{1}_{1}(b_{j}^{o}) \right] \right)
\end{matrix}}
\end{pmatrix}
\end{equation*}
(with $\operatorname{tr}(\mathscr{C}^{o}(b_{j}^{o})) \! = \! 0)$, where
$Q_{0}^{o}(b_{j}^{o})$, $Q_{1}^{o}(b_{j}^{o})$ are given in Theorem~2.3.1,
Equations~(2.33) and~(2.34),
\begin{align*}
Q_{2}^{o}(b_{j}^{o}) &= -\dfrac{1}{2}Q_{0}^{o}(b_{j}^{o}) \! \left(\sum_{
\substack{k=1\\k \not= j}}^{N} \! \left(\dfrac{1}{(b_{j}^{o} \! - \! b_{k}^{
o})^{2}} \! - \! \dfrac{1}{(b_{j}^{o} \! - \! a_{k}^{o})^{2}} \right) \! + \!
\dfrac{1}{(b_{j}^{o} \! - \! b_{0}^{o})^{2}} \! - \! \dfrac{1}{(b_{j}^{o} \!
- \! a_{N+1}^{o})^{2}} \! - \! \dfrac{1}{(b_{j}^{o} \! - \! a_{j}^{o})^{2}}
\right) \\
&+ \, \dfrac{1}{4}Q_{0}^{o}(b_{j}^{o}) \! \left(\sum_{\substack{k=1\\k \not=
j}}^{N} \! \left(\dfrac{1}{b_{j}^{o} \! - \! b_{k}^{o}} \! - \! \dfrac{1}{b_{
j}^{o} \! - \! a_{k}^{o}} \right) \! + \! \dfrac{1}{b_{j}^{o} \! - \! b_{0}^{
o}} \! - \! \dfrac{1}{b_{j}^{o} \! - \! a_{N+1}^{o}} \! - \! \dfrac{1}{b_{j}^{
o} \! - \! a_{j}^{o}} \right)^{2},
\end{align*}
$\mathscr{C}^{o}(b_{0}^{o})$ is given by the same expression as $\mathscr{C}^{
o}(b_{j}^{o})$ above subject to the modifications $\Omega_{j}^{o} \! \to \! 0$,
$b_{j}^{o} \! \to \! b_{0}^{o}$, $Q_{0}^{o}(b_{j}^{o}) \! \to \! Q_{0}^{o}(b_{
0}^{o})$, $Q_{1}^{o}(b_{j}^{o}) \! \to \! Q_{1}^{o}(b_{0}^{o})$, with $Q_{0}^{
o}(b_{0}^{o})$, $Q_{1}^{o}(b_{0}^{o})$ given in Theorem~2.3.1, Equations~(2.29)
and~(2.30), and $Q_{2}^{o}(b_{j}^{o}) \! \to \! Q_{2}^{o}(b_{0}^{o})$, where
\begin{align*}
Q_{2}^{o}(b_{0}^{o})&= -\dfrac{1}{2}Q_{0}^{o}(b_{0}^{o}) \! \left(\sum_{k=1}^{
N} \! \left(\dfrac{1}{(b_{0}^{o} \! - \! b_{k}^{o})^{2}} \! - \! \dfrac{1}{(
b_{0}^{o} \! - \! a_{k}^{o})^{2}} \right) \! - \! \dfrac{1}{(b_{0}^{o} \! - \!
a_{N+1}^{o})^{2}} \right) \\
&+ \, \dfrac{1}{4}Q_{0}^{o}(b_{0}^{o}) \! \left(\sum_{k=1}^{N} \! \left(
\dfrac{1}{b_{0}^{o} \! - \! b_{k}^{o}} \! - \! \dfrac{1}{b_{0}^{o} \! - \!
a_{k}^{o}} \right) \! - \! \dfrac{1}{b_{0}^{o} \! - \! a_{N+1}^{o}} \right)^{
2},
\end{align*}
$(f_{k}^{b_{j-1}^{o}}(n))_{l_{1}l_{2}} \! =_{n \to \infty} \! \mathcal{O}(1)$,
$j \! = \! 1,\dotsc,N \! + \! 1$, $k \! \in \! \mathbb{N}$, $l_{1},l_{2} \! =
\! 1,2$, and, for $j \! = \! 1,\dotsc,N \! + \! 1$,
\begin{align*}
\oint_{\partial \mathbb{U}^{o}_{\delta_{b_{j-1}}}}s^{-1}w_{+}^{\Sigma^{o}_{
\circlearrowright}}(s) \, \dfrac{\md s}{2 \pi \mi} \underset{n \to \infty}{=}&
\, \dfrac{(\mathscr{B}^{o}(b_{j-1}^{o}) \widehat{\alpha}^{o}_{0}(b_{j-1}^{o}
) \! - \! \mathscr{A}^{o}(b_{j-1}^{o})(\widehat{\alpha}^{o}_{1}(b_{j-1}^{o})
\! + \! (b_{j-1}^{o})^{-1} \widehat{\alpha}^{o}_{0}(b_{j-1}^{o})))}{(n \! +
\! \frac{1}{2})(\widehat{\alpha}^{o}_{0}(b_{j-1}^{o}))^{2}b_{j-1}^{o}} \\
&+ \, \mathcal{O} \! \left(\dfrac{f_{j}(n)}{(n \! + \! \frac{1}{2})^{2}}
\right),
\end{align*}
where $(f_{j}(n))_{kl} \! =_{n \to \infty} \! \mathcal{O}(1)$, $k,l \! = \!
1,2$. \hfill $\blacksquare$
\end{eeeee}

Re-tracing the finite sequence of RHP transformations (all of which are
invertible) and definitions, namely, $\mathscr{R}^{o}(z)$ (Lemmae~5.3 and~4.8)
and $\mathscr{S}^{o}_{p}(z)$ (Lemma~4.8) $\to$ $\mathcal{X}^{o}(z)$ (Lemmae
4.6 and~4.7) $\to$ $\overset{o}{m}^{\raise-1.0ex\hbox{$\scriptstyle \infty$}}
(z)$ (Lemma~4.5) $\to$
$\overset{o}{\mathscr{M}}^{\raise-1.0ex\hbox{$\scriptstyle \sharp$}}(z)$
(Lemma~4.2) $\to$
$\overset{o}{\mathscr{M}}^{\raise-1.0ex\hbox{$\scriptstyle \flat$}}(z)$
(Proposition~4.1) $\to$
$\overset{o}{\mathscr{M}}(z)$ (Lemma~4.1) $\to$ $\overset{o}{\mathrm{Y}}(z)$
(Lemma~3.4), the asymptotic (as $n \! \to \! \infty)$ solution of the original
\textbf{RHP2}, that is, $(\overset{o}{\mathrm{Y}}(z),\mathrm{I} \! + \! \exp
(-n \widetilde{V}(z)) \sigma_{+},\mathbb{R})$, in the various bounded and
unbounded regions (Figure~7), is given by:
\begin{enumerate}
\item[(1)] for $z \! \in \! \Upsilon^{o}_{1}$,
\begin{equation*}
\overset{o}{\mathrm{Y}}(z) \! = \! \me^{\frac{n \ell_{o}}{2} \operatorname{ad}
(\sigma_{3})} \mathscr{R}^{o}(z)
\overset{o}{m}^{\raise-1.0ex\hbox{$\scriptstyle \infty$}}(z) \mathbb{E}^{
\sigma_{3}} \me^{n(g^{o}(z)-\mathfrak{Q}^{+}_{\mathscr{A}}) \sigma_{3}},
\end{equation*}
and, for $z \! \in \! \Upsilon^{o}_{2}$,
\begin{equation*}
\overset{o}{\mathrm{Y}}(z) \! = \! \me^{\frac{n \ell_{o}}{2} \operatorname{ad}
(\sigma_{3})} \mathscr{R}^{o}(z)
\overset{o}{m}^{\raise-1.0ex\hbox{$\scriptstyle \infty$}}(z) \mathbb{E}^{-
\sigma_{3}} \me^{n(g^{o}(z)-\mathfrak{Q}^{-}_{\mathscr{A}}) \sigma_{3}},
\end{equation*}
where $g^{o}(z)$ and $\mathfrak{Q}^{\pm}_{\mathscr{A}}$, $\ell_{o}$, $\mathbb{
E}$, $\overset{o}{m}^{\raise-1.0ex\hbox{$\scriptstyle \infty$}}(z)$, and
$\mathscr{R}^{o}(z)$ are given in Lemma~3.4, Lemma~3.6, Proposition~4.1,
Lemma~4.5, and Lemma~5.3, respectively;
\item[(2)] for $z \! \in \! \Upsilon^{o}_{3}$,
\begin{equation*}
\overset{o}{\mathrm{Y}}(z) \! = \! \me^{\frac{n \ell_{o}}{2} \operatorname{ad}
(\sigma_{3})} \mathscr{R}^{o}(z)
\overset{o}{m}^{\raise-1.0ex\hbox{$\scriptstyle \infty$}}(z) \! \left(
\mathrm{I} \! + \! \me^{-4(n+\frac{1}{2}) \pi \mi \int_{z}^{a_{N+1}^{o}}
\psi_{V}^{o}(s) \, \md s} \, \sigma_{-} \right) \! \mathbb{E}^{\sigma_{3}}
\me^{n(g^{o}(z)-\mathfrak{Q}^{+}_{\mathscr{A}}) \sigma_{3}},
\end{equation*}
where $\psi_{V}^{o}(z)$ is given in Lemma~3.5, and, for $z \! \in \!
\Upsilon^{o}_{4}$,
\begin{equation*}
\overset{o}{\mathrm{Y}}(z) \! = \! \me^{\frac{n \ell_{o}}{2} \operatorname{ad}
(\sigma_{3})} \mathscr{R}^{o}(z)
\overset{o}{m}^{\raise-1.0ex\hbox{$\scriptstyle \infty$}}(z) \!
\left(\mathrm{I} \! - \! \me^{4(n+\frac{1}{2}) \pi \mi \int_{z}^{a_{N+1}^{o}}
\psi_{V}^{o}(s) \, \md s} \, \sigma_{-} \right) \! \mathbb{E}^{-\sigma_{3}}
\me^{n(g^{o}(z)-\mathfrak{Q}^{-}_{\mathscr{A}}) \sigma_{3}};
\end{equation*}
\item[(3)] for $z \! \in \! \Omega^{o,1}_{b_{j-1}} \cup \Omega^{o,1}_{a_{j}}$,
$j \! = \! 1,\dotsc,N \! + \! 1$,
\begin{equation*}
\overset{o}{\mathrm{Y}}(z) \! = \! \me^{\frac{n \ell_{o}}{2} \operatorname{ad}
(\sigma_{3})} \mathscr{R}^{o}(z) \mathcal{X}^{o}(z) \mathbb{E}^{\sigma_{3}}
\me^{n(g^{o}(z)-\mathfrak{Q}^{+}_{\mathscr{A}}) \sigma_{3}},
\end{equation*}
and, for $z \! \in \! \Omega^{o,4}_{b_{j-1}} \cup \Omega^{o,4}_{a_{j}}$, $j \!
= \! 1,\dotsc,N \! + \! 1$,
\begin{equation*}
\overset{o}{\mathrm{Y}}(z) \! = \! \me^{\frac{n \ell_{o}}{2} \operatorname{ad}
(\sigma_{3})} \mathscr{R}^{o}(z) \mathcal{X}^{o}(z) \mathbb{E}^{-\sigma_{3}}
\me^{n(g^{o}(z)-\mathfrak{Q}^{-}_{\mathscr{A}}) \sigma_{3}},
\end{equation*}
where, for $z \! \in \! \mathbb{U}^{o}_{\delta_{b_{j-1}}}$ $(\supset \Omega^{
o,1}_{b_{j-1}} \cup \Omega^{o,4}_{b_{j-1}})$, $\mathcal{X}^{o}(z)$ is given by
Lemma~4.6, and, for $z \! \in \! \mathbb{U}^{o}_{\delta_{a_{j}}}$ $(\supset
\Omega^{o,1}_{a_{j}} \cup \Omega^{o,4}_{a_{j}})$, $\mathcal{X}^{o}(z)$ is
given by Lemma~4.7; and
\item[(4)] for $z \! \in \! \Omega^{o,2}_{b_{j-1}} \cup \Omega^{o,2}_{a_{j}}$,
$j \! = \! 1,\dotsc,N \! + \! 1$,
\begin{equation*}
\overset{o}{\mathrm{Y}}(z) \! = \! \me^{\frac{n \ell_{o}}{2} \operatorname{ad}
(\sigma_{3})} \mathscr{R}^{o}(z) \mathcal{X}^{o}(z) \! \left(\mathrm{I} \! +
\! \me^{-4(n+\frac{1}{2}) \pi \mi \int_{z}^{a_{N+1}^{o}} \psi_{V}^{o}(s) \,
\md s} \, \sigma_{-} \right) \! \mathbb{E}^{\sigma_{3}} \me^{n(g^{o}(z)-
\mathfrak{Q}^{+}_{\mathscr{A}}) \sigma_{3}},
\end{equation*}
and, for $z \! \in \! \Omega^{o,3}_{b_{j-1}} \cup \Omega^{o,3}_{a_{j}}$, $j \!
= \! 1,\dotsc,N \! + \! 1$,
\begin{equation*}
\overset{o}{\mathrm{Y}}(z) \! = \! \me^{\frac{n \ell_{o}}{2} \operatorname{ad}
(\sigma_{3})} \mathscr{R}^{o}(z) \mathcal{X}^{o}(z) \! \left(\mathrm{I} \! -
\! \me^{4(n+\frac{1}{2}) \pi \mi \int_{z}^{a_{N+1}^{o}} \psi_{V}^{o}(s) \,
\md s} \, \sigma_{-} \right) \! \mathbb{E}^{-\sigma_{3}} \me^{n(g^{o}(z)-
\mathfrak{Q}^{-}_{\mathscr{A}}) \sigma_{3}}.
\end{equation*}
\end{enumerate}
Multiplying the respective matrices in items~(1)--(4) above and collecting
$(1 \, 1)$- and $(1 \, 2)$-elements, one arrives at, finally, the asymptotic
(as $n \! \to \! \infty)$ results for $z \boldsymbol{\pi}_{2n+1}(z)$ and
$\int_{\mathbb{R}} \tfrac{(s \boldsymbol{\pi}_{2n+1}(s)) \exp (-n \widetilde{
V}(s))}{s(s-z)} \, \tfrac{\md s}{2 \pi \mi}$ (in the entire complex plane)
stated in Theorem~2.3.1.

In order to obtain asymptotics (as $n \! \to \! \infty)$ for $\xi^{(2n+1)}_{-
n-1}$ $(= \! \norm{\boldsymbol{\pi}_{2n+1}(\pmb{\cdot})}_{\mathscr{L}}^{-1} 
\! = \! (H^{(-2n)}_{2n+1}/H^{(-2n-2)}_{2n+2})^{1/2})$ and $\phi_{2n+1}(z)$ $(=
\! \xi^{(2n+1)}_{-n-1} \boldsymbol{\pi}_{2n+1}(z))$ stated in Theorem~2.3.2,
small-$z$ asymptotics for $\overset{o}{\mathrm{Y}}(z)$ are necessary.
\begin{eeeee}
Since a tedious algebraic exercise shows that $\mathbb{C}_{-} \! \ni \! z \! 
\to \! 0$ asymptotics of 
$\overset{o}{m}^{\raise-1.0ex\hbox{$\scriptstyle \infty$}}(z)$ are obtained 
by multiplying $\mathbb{C}_{+} \! \ni \! z \! \to \! 0$ asymptotics of 
$\overset{o}{m}^{\raise-1.0ex\hbox{$\scriptstyle \infty$}}(z)$ on the right 
by $\exp (\mi (n \! + \! \tfrac{1}{2}) \Omega_{j}^{o} \sigma_{3})$ and using 
the relation $\mathbb{E}^{-\sigma_{3}} \exp (\mi (n \! + \! \tfrac{1}{2}) 
\Omega_{j}^{o} \sigma_{3}) \! = \! \mathbb{E}^{\sigma_{3}}$, only the 
asymptotic expansion as $\mathbb{C}_{+} \! \ni \! z \! \to \! 0$ of 
$\overset{o}{m}^{\raise-1.0ex\hbox{$\scriptstyle \infty$}}(z)$ is presented 
in Proposition~5.3 below. \hfill $\blacksquare$
\end{eeeee}
\begin{bbbbb}
Let $\mathscr{R}^{o} \colon \mathbb{C} \setminus \widetilde{\Sigma}_{p}^{o}
\! \to \! \operatorname{SL}_{2}(\mathbb{C})$ be the solution of the {\rm RHP}
$(\mathscr{R}^{o}(z),\upsilon_{\mathscr{R}}^{o}(z),\widetilde{\Sigma}_{p}^{o}
)$ formulated in Proposition~{\rm 5.2} with $n \! \to \! \infty$ asymptotics
given in Lemma~{\rm 5.3}. Then,
\begin{equation*}
\mathscr{R}^{o}(z) \underset{z \to 0}{=} \mathrm{I} \! + \! \mathscr{R}^{o,
0}_{1}(n)z \! + \! \mathscr{R}^{o,0}_{2}(n)z^{2} \! + \! \mathcal{O}(z^{3}),
\end{equation*}
where, for $k \! = \! 2,3$,
\begin{equation*}
\mathscr{R}^{o,0}_{k-1}(n) \! := \! \int_{\Sigma^{o}_{\circlearrowright}}
s^{-k}w_{+}^{\Sigma^{o}_{\circlearrowright}}(s) \, \dfrac{\md s}{2 \pi \mi}
\! = \! -\sum_{j=1}^{N+1} \, \sum_{q \in \{b_{j-1}^{o},a_{j}^{o}\}}
\operatorname{Res} \! \left(z^{-k}w^{\Sigma^{o}_{\circlearrowright}}_{+}(z);
q \right),
\end{equation*}
with, in particular,
\begin{align*}
\mathscr{R}^{o,0}_{k-1}(n) \underset{n \to \infty}{=}& \, \dfrac{1}{(n \! + \!
\frac{1}{2})} \sum_{j=1}^{N+1} \! \left(\dfrac{(\mathscr{A}^{o}(b_{j-1}^{o})
(\widehat{\alpha}_{1}^{o}(b_{j-1}^{o}) \! + \! k(b_{j-1}^{o})^{-1} \widehat{
\alpha}_{0}^{o}(b_{j-1}^{o})) \! - \! \mathscr{B}^{o}(b_{j-1}^{o}) \widehat{
\alpha}_{0}^{o}(b_{j-1}^{o}))}{(b_{j-1}^{o})^{k}(\widehat{\alpha}_{0}^{o}(b_{j
-1}^{o}))^{2}} \right. \\
+&\left. \, \dfrac{(\mathscr{A}^{o}(a_{j}^{o})(\widehat{\alpha}_{1}^{o}(a_{j}^{
o}) \! + \! k(a_{j}^{o})^{-1} \widehat{\alpha}_{0}^{o}(a_{j}^{o})) \! - \!
\mathscr{B}^{o}(a_{j}^{o}) \widehat{\alpha}_{0}^{o}(a_{j}^{o}))}{(a_{j}^{o})^{
k}(\widehat{\alpha}_{0}^{o}(a_{j}^{o}))^{2}} \right) \! + \! \mathcal{O} \!
\left(\dfrac{1}{(n \! + \! \frac{1}{2})^{2}} \right),
\end{align*}
and all parameters are defined in Lemma~{\rm 5.3}.

Let $\overset{o}{m}^{\raise-1.0ex\hbox{$\scriptstyle \infty$}} \colon \mathbb{
C} \setminus J_{o}^{\infty} \! \to \! \operatorname{SL}_{2}(\mathbb{C})$ solve
the {\rm RHP} $(\overset{o}{m}^{\raise-1.0ex\hbox{$\scriptstyle \infty$}}(z),
J_{o}^{\infty},
\overset{o}{\upsilon}^{\raise-1.0ex\hbox{$\scriptstyle \infty$}}(z))$
formulated in Lemma~{\rm 4.3} with (unique) solution given by
Lemma~{\rm 4.5}. For $\varepsilon_{1},\varepsilon_{2} \! = \! \pm 1$, set
\begin{gather*}
\theta^{o}_{0}(\varepsilon_{1},\varepsilon_{2},\boldsymbol{\Omega}^{o}) \! :=
\! \boldsymbol{\theta}^{o}(\varepsilon_{1} \boldsymbol{u}^{o}_{+}(0) \! - \!
\tfrac{1}{2 \pi}(n \! + \! \tfrac{1}{2}) \boldsymbol{\Omega}^{o} \! + \!
\varepsilon_{2} \boldsymbol{d}_{o}), \\
\alpha^{o}_{0}(\varepsilon_{1},\varepsilon_{2},\boldsymbol{\Omega}^{o}) \!
:= \! 2 \pi \mi \varepsilon_{1} \sum_{m \in \mathbb{Z}^{N}} (m,\widehat{
\boldsymbol{\alpha}}^{o}_{0}) \me^{2 \pi \mi (m,\varepsilon_{1} \boldsymbol{
u}^{o}_{+}(0)-\frac{1}{2 \pi}(n+\frac{1}{2}) \boldsymbol{\Omega}^{o}+
\varepsilon_{2} \boldsymbol{d}_{o})+\pi \mi (m,\tau^{o}m)},
\end{gather*}
where $\widehat{\boldsymbol{\alpha}}^{o}_{0} \! = \! (\widehat{\alpha}_{0,
1}^{o},\widehat{\alpha}^{o}_{0,2},\dotsc,\widehat{\alpha}^{o}_{0,N})$, with
$\widehat{\alpha}^{o}_{0,j} \! := \! (-1)^{\mathcal{N}_{+}}(\prod_{i=1}^{N+1}
\vert b_{i-1}^{o}a_{i}^{o} \vert)^{-1/2}c_{jN}^{o}$, $j \! = \! 1,\dotsc,N$,
where $\mathcal{N}_{+} \! \in \! \lbrace 0,\dotsc,N \! + \! 1 \rbrace$ is the
number of bands to the right of $z \! = \! 0$, and
\begin{gather*}
\beta^{o}_{0}(\varepsilon_{1},\varepsilon_{2},\boldsymbol{\Omega}^{o}) \! :=
\! 2 \pi \sum_{m \in \mathbb{Z}^{N}} \! \left(\mi \varepsilon_{1}(m,\widehat{
\boldsymbol{\beta}}^{o}_{0}) \! - \!  \pi (m,\widehat{\boldsymbol{\alpha}}^{
o}_{0})^{2} \right) \! \me^{2 \pi \mi (m,\varepsilon_{1} \boldsymbol{u}^{o}_{+}
(0)-\frac{1}{2 \pi}(n+\frac{1}{2}) \boldsymbol{\Omega}^{o}+\varepsilon_{2}
\boldsymbol{d}_{o})+ \pi \mi (m,\tau^{o}m)},
\end{gather*}
where $\widehat{\boldsymbol{\beta}}^{o}_{0} \! = \! (\widehat{\beta}^{o}_{0,
1},\widehat{\beta}^{o}_{0,2},\dotsc,\widehat{\beta}^{o}_{0,N})$, with
$\widehat{\beta}^{o}_{0,j} \! := \! \tfrac{1}{2}(-1)^{\mathcal{N}_{+}}(\prod_{
i=1}^{N+1} \vert b_{i-1}^{o}a_{i}^{o} \vert)^{-1/2}(c^{o}_{jN-1} \! + \!
\tfrac{1}{2}c_{jN}^{o} \sum_{k=1}^{N+1}((a_{k}^{o})^{-1} \! + \! (b_{k-1}^{
o})^{-1}))$, $j \! = \! 1,\dotsc,N$, where $c_{jN}^{o},c_{jN-1}^{o}$, $j \! 
= \! 1,\dotsc,N$, are obtained {}from Equations~{\rm (O1)} and~{\rm (O2)}. 
Then,
\begin{gather*}
\overset{o}{m}^{\raise-1.0ex\hbox{$\scriptstyle \infty$}}(z) \underset{z
\to 0}{=} \mathbb{E}^{-\sigma_{3}} \! + \!
\overset{o}{m}_{1}^{\raise-1.0ex\hbox{$\scriptstyle 0$}}z \! + \!
\overset{o}{m}_{2}^{\raise-1.0ex\hbox{$\scriptstyle 0$}}z^{2} \! + \!
\mathcal{O}(z^{3}),
\end{gather*}
where
\begin{align*}
(\overset{o}{m}_{1}^{\raise-1.0ex\hbox{$\scriptstyle 0$}})_{11} &= -\dfrac{
\boldsymbol{\theta}^{o}(\boldsymbol{u}^{o}_{+}(0) \! + \! \boldsymbol{d}_{o})
\mathbb{E}^{-1}}{\boldsymbol{\theta}^{o}(\boldsymbol{u}^{o}_{+}(0) \! - \!
\frac{1}{2 \pi}(n \! + \! \frac{1}{2}) \boldsymbol{\Omega}^{o} \! + \!
\boldsymbol{d}_{o})} \! \left(\dfrac{\theta^{o}_{0}(1,1,\boldsymbol{\Omega}^{o}
) \alpha^{o}_{0}(1,1,\vec{\pmb{0}}) \! - \! \alpha^{o}_{0}(1,1,\boldsymbol{
\Omega}^{o}) \theta^{o}_{0}(1,1,\vec{\pmb{0}})}{(\theta^{o}_{0}(1,1,\vec{
\pmb{0}}))^{2}} \right), \\
(\overset{o}{m}_{1}^{\raise-1.0ex\hbox{$\scriptstyle 0$}})_{12} &= \dfrac{
1}{4 \mi} \! \left(\sum_{k=1}^{N+1} \! \left(\dfrac{1}{b_{k-1}^{o}} \! - \!
\dfrac{1}{a_{k}^{o}} \right) \right) \! \dfrac{\boldsymbol{\theta}^{o}
(\boldsymbol{u}^{o}_{+}(0) \! + \! \boldsymbol{d}_{o}) \theta^{o}_{0}(-1,1,
\boldsymbol{\Omega}^{o}) \mathbb{E}^{-1}}{\boldsymbol{\theta}^{o}(\boldsymbol{
u}^{o}_{+}(0) \! - \! \frac{1}{2 \pi}(n \! + \! \frac{1}{2}) \boldsymbol{
\Omega}^{o} \! + \! \boldsymbol{d}_{o}) \theta^{o}_{0}(-1,1,\vec{\pmb{0}})}, \\
(\overset{o}{m}_{1}^{\raise-1.0ex\hbox{$\scriptstyle 0$}})_{21}&= -\dfrac{1}{4
\mi} \! \left(\sum_{k=1}^{N+1} \! \left(\dfrac{1}{b_{k-1}^{o}} \! - \! \dfrac{
1}{a_{k}^{o}} \right) \right) \! \dfrac{\boldsymbol{\theta}^{o}(\boldsymbol{
u}^{o}_{+}(0) \! + \! \boldsymbol{d}_{o}) \theta^{o}_{0}(1,-1,\boldsymbol{
\Omega}^{o}) \mathbb{E}}{\boldsymbol{\theta}^{o}(-\boldsymbol{u}^{o}_{+}(0) \!
- \! \frac{1}{2 \pi}(n \! + \! \frac{1}{2}) \boldsymbol{\Omega}^{o} \! - \!
\boldsymbol{d}_{o}) \theta^{o}_{0}(1,-1,\vec{\pmb{0}})}, \\
(\overset{o}{m}_{1}^{\raise-1.0ex\hbox{$\scriptstyle 0$}})_{22} &= -\left(
\dfrac{\theta^{o}_{0}(-1,-1,\boldsymbol{\Omega}^{o}) \alpha^{o}_{0}(-1,-1,
\vec{\pmb{0}}) \! - \! \alpha^{o}_{0}(-1,-1,\boldsymbol{\Omega}^{o}) \theta^{
o}_{0}(-1,-1,\vec{\pmb{0}})}{(\theta^{o}_{0}(-1,-1,\vec{\pmb{0}}))^{2}}
\right) \\
&\times \dfrac{\boldsymbol{\theta}^{o}(\boldsymbol{u}^{o}_{+}(0) \! + \!
\boldsymbol{d}_{o}) \mathbb{E}}{\boldsymbol{\theta}^{o}(-\boldsymbol{u}^{o}_{+}
(0) \! - \! \frac{1}{2 \pi}(n \! + \! \frac{1}{2}) \boldsymbol{\Omega}^{o} \!
- \! \boldsymbol{d}_{o})}, \\
(\overset{o}{m}_{2}^{\raise-1.0ex\hbox{$\scriptstyle 0$}})_{11} &= \left(
\theta^{o}_{0}(1,1,\boldsymbol{\Omega}^{o}) \! \left(-\beta^{0}_{0}(1,1,\vec{
\pmb{0}}) \theta^{o}_{0}(1,1,\vec{\pmb{0}}) \! + \! (\alpha^{o}_{0}(1,1,\vec{
\pmb{0}}))^{2} \right) \! - \! \alpha^{o}_{0}(1,1,\boldsymbol{\Omega}^{o})
\alpha^{o}_{0}(1,1,\vec{\pmb{0}}) \theta^{o}_{0}(1,1,\vec{\pmb{0}}) \right. \\
&\left. + \, \beta^{o}_{0}(1,1,\boldsymbol{\Omega}^{o})(\theta^{o}_{0}(1,1,
\vec{\pmb{0}}))^{2} \right) \! \dfrac{(\theta^{o}_{0}(1,1,\vec{\pmb{0}}))^{-3}
\boldsymbol{\theta}^{o}(\boldsymbol{u}^{o}_{+}(0) \! + \! \boldsymbol{d}_{o})
\mathbb{E}^{-1}}{\boldsymbol{\theta}^{o}(\boldsymbol{u}^{o}_{+}(0) \! - \!
\frac{1}{2 \pi}(n \! + \! \frac{1}{2}) \boldsymbol{\Omega}^{o} \! + \!
\boldsymbol{d}_{o})} \! + \! \dfrac{1}{32} \! \left(\sum_{k=1}^{N+1} \! \left(
\dfrac{1}{b_{k-1}^{o}} \! - \! \dfrac{1}{a_{k}^{o}} \right) \right)^{2} \!
\mathbb{E}^{-1}, \\
(\overset{o}{m}_{2}^{\raise-1.0ex\hbox{$\scriptstyle 0$}})_{12} &= -\dfrac{
\boldsymbol{\theta}^{o}(\boldsymbol{u}^{o}_{+}(0) \! + \! \boldsymbol{d}_{o})
\mathbb{E}^{-1}}{\boldsymbol{\theta}^{o}(\boldsymbol{u}^{o}_{+}(0) \! - \!
\frac{1}{2 \pi}(n \! + \! \frac{1}{2}) \boldsymbol{\Omega}^{o} \! + \!
\boldsymbol{d}_{o})} \! \left(\! \left(\dfrac{\theta^{o}_{0}(-1,1,\boldsymbol{
\Omega}^{o}) \alpha^{o}_{0}(-1,1,\vec{\pmb{0}}) \! - \! \alpha^{o}_{0}(-1,1,
\boldsymbol{\Omega}^{o}) \theta^{o}_{0}(-1,1,\vec{\pmb{0}})}{(\theta^{o}_{0}
(-1,1,\vec{\pmb{0}}))^{2}} \right) \right. \\
&\left. \times \dfrac{1}{4 \mi} \! \left(\sum_{k=1}^{N+1} \! \left(\dfrac{1}{
b_{k-1}^{o}} \! - \! \dfrac{1}{a_{k}^{o}} \right) \right) \! - \! \dfrac{1}{8
\mi} \! \left(\sum_{k=1}^{N+1} \! \left(\dfrac{1}{(b_{k-1}^{o})^{2}} \! - \!
\dfrac{1}{(a_{k}^{o})^{2}} \right) \right) \! \dfrac{\theta^{o}_{0}(-1,1,
\boldsymbol{\Omega}^{o})}{\theta^{o}_{0}(-1,1,\vec{\pmb{0}})} \right), \\
(\overset{o}{m}_{2}^{\raise-1.0ex\hbox{$\scriptstyle 0$}})_{21} &= \dfrac{
\boldsymbol{\theta}^{o}(\boldsymbol{u}^{o}_{+}(0) \! + \! \boldsymbol{d}_{o})
\mathbb{E}}{\boldsymbol{\theta}^{o}(-\boldsymbol{u}^{o}_{+}(0) \! - \! \frac{
1}{2 \pi}(n \! + \! \frac{1}{2}) \boldsymbol{\Omega}^{o} \! - \! \boldsymbol{
d}_{o})} \! \left(\! \left(\dfrac{\theta^{o}_{0}(1,-1,\boldsymbol{\Omega}^{o})
\alpha^{o}_{0}(1,-1,\vec{\pmb{0}}) \! - \! \alpha^{o}_{0}(1,-1,\boldsymbol{
\Omega}^{o}) \theta^{o}_{0}(1,-1,\vec{\pmb{0}})}{(\theta^{o}_{0}(1,-1,\vec{
\pmb{0}}))^{2}} \right) \right. \\
&\left. \times \dfrac{1}{4 \mi} \! \left(\sum_{k=1}^{N+1} \! \left(\dfrac{1}{
b_{k-1}^{o}} \! - \! \dfrac{1}{a_{k}^{o}} \right) \right) \! - \! \dfrac{1}{8
\mi} \! \left(\sum_{k=1}^{N+1} \! \left(\dfrac{1}{(b_{k-1}^{o})^{2}} \! - \!
\dfrac{1}{(a_{k}^{o})^{2}} \right) \right) \! \dfrac{\theta^{o}_{0}(1,-1,
\boldsymbol{\Omega}^{o})}{\theta^{o}_{0}(1,-1,\vec{\pmb{0}})} \right), \\
(\overset{o}{m}_{2}^{\raise-1.0ex\hbox{$\scriptstyle 0$}})_{22} &= \left(
\theta^{o}_{0}(-1,-1,\boldsymbol{\Omega}^{o}) \! \left(-\beta^{o}_{0}(-1,-1,
\vec{\pmb{0}}) \theta^{o}_{0}(-1,-1,\vec{\pmb{0}}) \! + \! (\alpha^{o}_{0}(-1,
-1,\vec{\pmb{0}}))^{2} \right) \! - \! \alpha^{o}_{0}(-1,-1,\boldsymbol{
\Omega}^{o}) \right. \\
&\left. \times \, \alpha^{o}_{0}(-1,-1,\vec{\pmb{0}}) \theta^{o}_{0}(-1,-1,
\vec{\pmb{0}}) \! + \! \beta^{o}_{0}(-1,-1,\boldsymbol{\Omega}^{o})(\theta^{
o}_{0}(-1,-1,\vec{\pmb{0}}))^{2} \right) \! (\theta^{o}_{0}(-1,-1,\vec{\pmb{
0}}))^{-3} \\
&\times \dfrac{\boldsymbol{\theta}^{o}(\boldsymbol{u}^{o}_{+}(0) \! + \!
\boldsymbol{d}_{o}) \mathbb{E}}{\boldsymbol{\theta}^{o}(-\boldsymbol{u}^{o}_{+}
(0) \! - \! \frac{1}{2 \pi}(n \! + \! \frac{1}{2}) \boldsymbol{\Omega}^{o} \!
- \! \boldsymbol{d}_{o})} \! + \! \dfrac{1}{32} \! \left(\sum_{k=1}^{N+1} \!
\left(\dfrac{1}{b_{k-1}^{o}} \! - \! \dfrac{1}{a_{k}^{o}} \right) \right)^{2}
\! \mathbb{E},
\end{align*}
with $(\star)_{ij}$, $i,j \! = \! 1,2$, denoting the $(i \, j)$-element of 
$\star$, $\mathbb{E}$ defined in Proposition~{\rm 4.1}, and $\vec{\pmb{0}} \! 
:= \! (0,0,\dotsc,0)^{\mathrm{T}}$ $(\in \! \mathbb{R}^{N})$.

Let $\overset{o}{\mathrm{Y}} \colon \mathbb{C} \setminus \mathbb{R} \! \to
\! \operatorname{SL}_{2}(\mathbb{C})$ be the solution of {\rm \pmb{RHP2}}.
Then,
\begin{equation*}
\overset{o}{\mathrm{Y}}(z)z^{n \sigma_{3}} \underset{z \to 0}{=} \mathrm{I}
\! + \! z \mathrm{Y}^{o,0}_{1} \! + \! z^{2} \mathrm{Y}^{o,0}_{2} \! + \!
\mathcal{O}(z^{3}),
\end{equation*}
where
\begin{align*}
(\mathrm{Y}^{o,0}_{1})_{11} &= -(2n \! + \! 1) \int_{J_{o}}s^{-1} \psi_{V}^{o}
(s) \, \md s \! + \!
(\overset{o}{m}_{1}^{\raise-1.0ex\hbox{$\scriptstyle 0$}})_{11} \mathbb{E} \!
+ \! (\mathscr{R}^{o,0}_{1}(n))_{11}, \\
(\mathrm{Y}^{o,0}_{1})_{12} &= \me^{n \ell_{o}} \! \left(
(\overset{o}{m}_{1}^{\raise-1.0ex\hbox{$\scriptstyle 0$}})_{12} \mathbb{E}^{-
1} \! + \! (\mathscr{R}^{o,0}_{1}(n))_{12} \right), \\
(\mathrm{Y}^{o,0}_{1})_{21} &= \me^{-n \ell_{o}} \! \left(
(\overset{o}{m}_{1}^{\raise-1.0ex\hbox{$\scriptstyle 0$}})_{21} \mathbb{E} \!
+ \! (\mathscr{R}^{o,0}_{1}(n))_{21} \right), \\
(\mathrm{Y}^{o,0}_{1})_{22} &= (2n \! + \! 1) \int_{J_{o}}s^{-1} \psi_{V}^{o}
(s) \, \md s \! + \!
(\overset{o}{m}_{1}^{\raise-1.0ex\hbox{$\scriptstyle 0$}})_{22} \mathbb{E}^{-
1} \! + \! (\mathscr{R}^{o,0}_{1}(n))_{22}, \\
(\mathrm{Y}^{o,0}_{2})_{11} =& \, \tfrac{1}{2}(2n \! + \! 1)^{2} \! \left(
\int_{J_{o}}s^{-1} \psi_{V}^{o}(s) \, \md s \right)^{2} \! - \! \tfrac{1}{2}
(2n \! + \! 1) \int_{J_{o}}s^{-2} \psi_{V}^{o}(s) \, \md s \! - \! (2n \! + \!
1) \left((\overset{o}{m}_{1}^{\raise-1.0ex\hbox{$\scriptstyle 0$}})_{11}
\mathbb{E} \! + \! (\mathscr{R}^{o,0}_{1}(n))_{11} \right) \\
\times& \, \int_{J_{o}}s^{-1} \psi_{V}^{o}(s) \, \md s \! + \! (\overset{o}{
m}_{2}^{\raise-1.0ex\hbox{$\scriptstyle 0$}})_{11} \mathbb{E} \! + \!
(\mathscr{R}^{o,0}_{2}(n))_{11} \! + \! \left((\mathscr{R}^{o,0}_{1}(n))_{11}
(\overset{o}{m}_{1}^{\raise-1.0ex\hbox{$\scriptstyle 0$}})_{11} \! + \!
(\mathscr{R}^{o,0}_{1}(n))_{12}
(\overset{o}{m}_{1}^{\raise-1.0ex\hbox{$\scriptstyle 0$}})_{21} \right) \!
\mathbb{E}, \\
(\mathrm{Y}^{o,0}_{2})_{12} &= \me^{n \ell_{o}} \! \left((2n \! + \! 1) \!
\left((\overset{o}{m}_{1}^{\raise-1.0ex\hbox{$\scriptstyle 0$}})_{12}
\mathbb{E}^{-1} \! + \! (\mathscr{R}^{o,0}_{1}(n))_{12} \right) \! \int_{J_{o}}
s^{-1} \psi_{V}^{o}(s) \, \md s \! + \!
(\overset{o}{m}_{2}^{\raise-1.0ex\hbox{$\scriptstyle 0$}})_{12} \mathbb{E}^{-
1} \! + \! (\mathscr{R}^{o,0}_{2}(n))_{12} \right. \\
&\left. + \, \left((\mathscr{R}^{o,0}_{1}(n))_{11}
(\overset{o}{m}_{1}^{\raise-1.0ex\hbox{$\scriptstyle 0$}})_{12} \! + \!
(\mathscr{R}^{o,0}_{1}(n))_{12}
(\overset{o}{m}_{1}^{\raise-1.0ex\hbox{$\scriptstyle 0$}})_{22} \right) \!
\mathbb{E}^{-1} \right), \\
(\mathrm{Y}^{o,0}_{2})_{21} &= \me^{-n \ell_{o}} \! \left(-(2n \! + \! 1) \!
\left((\overset{o}{m}_{1}^{\raise-1.0ex\hbox{$\scriptstyle 0$}})_{21} \mathbb{
E} \! + \! (\mathscr{R}^{o,0}_{1}(n))_{21} \right) \! \int_{J_{o}}s^{-1} \psi_{
V}^{o}(s) \, \md s \! + \!
(\overset{o}{m}_{2}^{\raise-1.0ex\hbox{$\scriptstyle 0$}})_{21} \mathbb{E} \!
+ \! (\mathscr{R}^{o,0}_{2}(n))_{21} \right. \\
&\left. + \, \left((\mathscr{R}^{o,0}_{1}(n))_{21}
(\overset{o}{m}_{1}^{\raise-1.0ex\hbox{$\scriptstyle 0$}})_{11} \! + \!
(\mathscr{R}^{o,0}_{1}(n))_{22}
(\overset{o}{m}_{1}^{\raise-1.0ex\hbox{$\scriptstyle 0$}})_{21} \right)
\mathbb{E} \right), \\
(\mathrm{Y}^{o,0}_{2})_{22} =& \, \tfrac{1}{2}(2n \! + \! 1)^{2} \! \left(
\int_{J_{o}}s^{-1} \psi_{V}^{o}(s) \, \md s \right)^{2} \! + \! \tfrac{1}{2}
(2n \! + \! 1) \int_{J_{o}}s^{-2} \psi_{V}^{o}(s) \, \md s \! + \! (2n \! + \!
1) \! \left((\overset{o}{m}_{1}^{\raise-1.0ex\hbox{$\scriptstyle 0$}})_{22}
\mathbb{E}^{-1} \! + \! (\mathscr{R}^{o,0}_{1}(n))_{22} \right) \\
\times& \, \int_{J_{o}}s^{-1} \psi_{V}^{o}(s) \, \md s \! + \!
(\overset{o}{m}_{2}^{\raise-1.0ex\hbox{$\scriptstyle 0$}})_{22} \mathbb{E}^{-
1} \! + \! (\mathscr{R}^{o,0}_{2}(n))_{22} \! + \! \left((\mathscr{R}^{o,0}_{1}
(n))_{21}(\overset{o}{m}_{1}^{\raise-1.0ex\hbox{$\scriptstyle 0$}})_{12} \! +
\! (\mathscr{R}^{o,0}_{1}(n))_{22}
(\overset{o}{m}_{1}^{\raise-1.0ex\hbox{$\scriptstyle 0$}})_{22} \right) \!
\mathbb{E}^{-1}.
\end{align*}
\end{bbbbb}

\emph{Proof.} Let $\mathscr{R}^{o} \colon \mathbb{C} \setminus \widetilde{
\Sigma}^{o}_{p} \! \to \! \operatorname{SL}_{2}(\mathbb{C})$ be the solution
of the RHP $(\mathscr{R}^{o}(z),\upsilon_{\mathscr{R}}^{o}(z),\widetilde{
\Sigma}^{o}_{p})$ formulated in Proposition~5.2 with $n \! \to \! \infty$
asymptotics given in Lemma~5.3. For $\vert z \vert \! \ll \! \min_{j=1,
\dotsc,N+1} \lbrace \vert b_{j-1}^{o} \! - \! a_{j}^{o} \vert \rbrace$, via
the expansion $\tfrac{1}{z-s} \! = \! -\sum_{k=0}^{l} \tfrac{z^{k}}{s^{k+1}}
\! + \! \tfrac{z^{l+1}}{s^{l+1}(z-s)}$, $l \! \in \! \mathbb{Z}_{0}^{+}$,
where $s \! \in \! \lbrace b_{j-1}^{o},a_{j}^{o} \rbrace$, $j \! = \! 1,
\dotsc,N \! + \! 1$, one obtains the asymptotics for $\mathscr{R}^{o}(z)$
stated in the Proposition.

Let $\overset{o}{m}^{\raise-1.0ex\hbox{$\scriptstyle \infty$}} \colon \mathbb{
C} \setminus J_{o}^{\infty} \! \to \! \operatorname{SL}_{2}(\mathbb{C})$ solve
the RHP $(\overset{o}{m}^{\raise-1.0ex\hbox{$\scriptstyle \infty$}}(z),J_{o}^{
\infty},\overset{o}{\upsilon}^{\raise-1.0ex\hbox{$\scriptstyle \infty$}}(z))$
formulated in Lemma~4.3 with (unique) solution given by Lemma~4.5. In order
to obtain small-$z$ asymptotics of
$\overset{o}{m}^{\raise-1.0ex\hbox{$\scriptstyle \infty$}}(z)$, one needs
small-$z$ asymptotics of $(\gamma^{o}(z))^{\pm 1}$ and $\tfrac{\boldsymbol{
\theta}^{o}(\varepsilon_{1} \boldsymbol{u}^{o}(z)-\frac{1}{2 \pi}(n+\frac{1}{
2}) \boldsymbol{\Omega}^{o}+\varepsilon_{2} \boldsymbol{d}_{o})}{\boldsymbol{
\theta}^{o}(\varepsilon_{1} \boldsymbol{u}^{o}(z)+\varepsilon_{2} \boldsymbol{
d}_{o})}$, $\varepsilon_{1},\varepsilon_{2} \! = \! \pm 1$. Consider, say, and
without loss of generality, $z \! \to \! 0$ asymptotics for $z \! \in \!
\mathbb{C}_{+}$ (designated $z \! \to \! 0^{+})$, where, by definition, $\sqrt{
\smash[b]{\star (z)}} \! := \! +\sqrt{\smash[b]{\star (z)}}$: equivalently,
one may consider $z \! \to \! 0$ asymptotics for $z \! \in \! \mathbb{C}_{-}$
(designated $z \! \to \! 0^{-})$; however, recalling that $\sqrt{\smash[b]{
\star (z)}} \! \upharpoonright_{\mathbb{C}_{+}} \! = \! -\sqrt{\smash[b]{\star
(z)}} \! \upharpoonright_{\mathbb{C}_{-}}$, one obtains (in either case, and
via the sheet-interchange index) the same $z \! \to \! 0$ asymptotics (for
$\overset{o}{m}^{\raise-1.0ex\hbox{$\scriptstyle \infty$}}(z))$. Recall the
expression for $\gamma^{o}(z)$ given in Lemma~4.4: for $\vert z \vert \! \ll
\! \min_{j=1,\dotsc,N+1} \lbrace \vert b_{j-1}^{o} \! - \! a_{j}^{o} \vert
\rbrace$, via the expansions $\tfrac{1}{z-s} \! = \! -\sum_{k=0}^{l} \tfrac{
z^{k}}{s^{k+1}} \! + \! \tfrac{z^{l+1}}{s^{l+1}(z-s)}$, $l \! \in \! \mathbb{
Z}_{0}^{+}$, and $\ln (s \! - \! z) \! =_{\vert z \vert \to 0} \! \ln (s) \!
- \! \sum_{k=1}^{\infty} \tfrac{1}{k}(\tfrac{z}{s})^{k}$, where $s \! \in \!
\lbrace b_{j-1}^{o},a_{j}^{o} \rbrace$, $j \! = \! 1,\dotsc,N \! + \! 1$, one
shows that, upon setting $\gamma^{o}_{0} \! := \! \gamma^{o}(0) \! = \!
(\prod_{k=1}^{N+1}b_{k-1}^{o}(a_{k}^{o})^{-1})^{1/4}$ $(> \! 0)$,
\begin{align*}
(\gamma^{o}_{0})^{\mp 1}(\gamma^{o}(z))^{\pm 1} \underset{z \to 0^{+}}{=}& 1
\! + \! z \! \left(\pm \dfrac{1}{4} \sum_{k=1}^{N+1} \! \left(\dfrac{1}{a_{
k}^{o}} \! - \! \dfrac{1}{b_{k-1}^{o}} \right) \right) \! + \! z^{2} \! \left(
\pm \dfrac{1}{8} \sum_{k=1}^{N+1} \! \left(\dfrac{1}{(a_{k}^{o})^{2}} \! - \!
\dfrac{1}{(b_{k-1}^{o})^{2}} \right) \right. \\
+&\left. \dfrac{1}{32} \! \left(\sum_{k=1}^{N+1} \! \left(\dfrac{1}{a_{k}^{o}}
\! - \! \dfrac{1}{b_{k-1}^{o}} \right) \right)^{2} \, \right) \! + \!
\mathcal{O}(z^{3}),
\end{align*}
whence
\begin{equation*}
\dfrac{1}{2}((\gamma^{o}_{0})^{-1} \gamma^{o}(z) \! + \! \gamma^{o}_{0}
(\gamma^{o}(z))^{-1}) \! \underset{z \to 0^{+}}{=} \! 1 \! + \! z^{2} \!
\left(\dfrac{1}{32} \! \left(\sum_{k=1}^{N+1} \! \left(\dfrac{1}{a_{k}^{o}} \!
- \! \dfrac{1}{b_{k-1}^{o}} \right) \right)^{2} \, \right) \! + \! \mathcal{O}
(z^{3}),
\end{equation*}
and
\begin{align*}
\dfrac{1}{2 \mi}((\gamma^{o}_{0})^{-1} \gamma^{o}(z) \! - \! \gamma^{o}_{0}
(\gamma^{o}(z))^{-1}) \underset{z \to 0^{+}}{=}& \, z \! \left(\dfrac{1}{4
\mi} \sum_{k=1}^{N+1} \! \left(\dfrac{1}{a_{k}^{o}} \! - \! \dfrac{1}{b_{k-
1}^{o}} \right) \right) \! + \! z^{2} \! \left(\dfrac{1}{8 \mi} \sum_{k=1}^{N
+1} \! \left(\dfrac{1}{(a_{k}^{o})^{2}} \! - \! \dfrac{1}{(b_{k-1}^{o})^{2}}
\right) \right) \\
+& \, \mathcal{O}(z^{3}).
\end{align*}
Recall {}from Lemma~4.5 that $\boldsymbol{u}^{o}(z) \! := \! \int_{a_{N+1}^{
o}}^{z} \boldsymbol{\omega}^{o}$ $(\in \operatorname{Jac}(\mathcal{Y}_{o})$, 
with $\mathcal{Y}_{o} \! := \! \lbrace \mathstrut (y,z); \, y^{2} \! = \! 
R_{o}(z) \rbrace)$, where $\boldsymbol{\omega}^{o}$, the associated normalised 
basis of holomorphic one-forms of $\mathcal{Y}_{o}$, is given by $\boldsymbol{
\omega}^{o} \! = \! (\omega_{1}^{o},\omega_{2}^{o},\dotsc,\omega_{N}^{o})$, 
with $\omega_{j}^{o} \! := \! \sum_{k=1}^{N}c_{jk}^{o}(\prod_{i=1}^{N+1}(z \! 
- \! b_{i-1}^{o})(z \! - \! a_{i}^{o}))^{-1/2}z^{N-k} \, \md z$, $j \! = \! 1,
\dotsc,N$, where $c_{jk}^{o}$, $j,k \! = \! 1,\dotsc,N$, are obtained {}from 
Equations~(O1) and~(O2). Writing
\begin{equation*}
\boldsymbol{u}^{o}(z) \! = \! \left(\int_{a_{N+1}^{o}}^{0^{+}} \! + \! \int_{
0^{+}}^{z} \right) \! \boldsymbol{\omega}^{o} \! = \! \boldsymbol{u}^{o}_{+}
(0) \! + \! \int_{0^{+}}^{z} \boldsymbol{\omega}^{o},
\end{equation*}
where $\boldsymbol{u}^{o}_{+}(0) \! := \! \int_{a_{N+1}^{o}}^{0^{+}}
\boldsymbol{\omega}^{o}$ (cf. Lemma~4.5), for $\vert z \vert \! \ll \! \min_{
j=1,\dotsc,N+1} \lbrace \vert b_{j-1}^{o} \! - \! a_{j}^{o} \vert \rbrace$,
via the expansions $\tfrac{1}{z-s} \! = \! -\sum_{k=0}^{l} \tfrac{z^{k}}{s^{k
+1}} \! + \! \tfrac{z^{l+1}}{s^{l+1}(z-s)}$, $l \! \in \! \mathbb{Z}_{0}^{+}$,
and $\ln (z \! - \! s) \! =_{\vert z \vert \to 0} \! \ln (s) \! - \! \sum_{k=
1}^{\infty} \tfrac{1}{k}(\tfrac{z}{s})^{k}$, where $s \! \in \! \lbrace b_{k-
1}^{o},a_{k}^{o} \rbrace$, $k \! = \! 1,\dotsc,N \! + \! 1$, one shows that,
for $j \! = \! 1,\dotsc,N$,
\begin{equation*}
\omega_{j}^{o} \! \underset{z \to 0^{+}}{=} \! (-1)^{\mathcal{N}_{+}} \!
\left(\prod_{i=1}^{N+1} \vert b_{i-1}^{o}a_{i}^{o} \vert \right)^{-1/2} \!
\left(c_{jN}^{o} \, \md z \! + \! \left(c_{jN-1}^{o} \! + \! \dfrac{c_{jN}^{
o}}{2} \sum_{k=1}^{N+1} \! \left(\dfrac{1}{b_{k-1}^{o}} \! + \! \dfrac{1}{
a_{k}^{o}} \right) \right) \! z \, \md z \! + \! \mathcal{O}(z^{2} \, \md z)
\right),
\end{equation*}
where $\mathcal{N}_{+} \! \in \! \lbrace 0,\dotsc,N \! + \! 1 \rbrace$ is
the number of bands to the right of $z \! = \! 0$, whence
\begin{align*}
\int_{0^{+}}^{z} \omega_{j}^{o} \underset{z \to 0^{+}}{=}& (-1)^{\mathcal{N}_{
+}} \! \left(\prod_{i=1}^{N+1} \vert b_{i-1}^{o}a_{i}^{o} \vert \right)^{-1/2}
\! \left(c_{jN}^{o}z \! + \! \dfrac{1}{2} \! \left(c_{jN-1}^{o} \! + \!
\dfrac{c_{jN}^{o}}{2} \sum_{k=1}^{N+1} \! \left(\dfrac{1}{b_{k-1}^{o}} \! +
\! \dfrac{1}{a_{k}^{o}} \right) \right) \! z^{2} \! + \! \mathcal{O}(z^{3})
\right) \\
=:& \, \widehat{\alpha}^{o}_{0,j}z \! + \! \widehat{\beta}^{o}_{0,j}z^{2} \!
+ \! \mathcal{O}(z^{3}).
\end{align*}
Defining $\theta^{o}_{0}(\varepsilon_{1},\varepsilon_{2},\boldsymbol{\Omega}^{
o})$, $\alpha^{o}_{0}(\varepsilon_{1},\varepsilon_{2},\boldsymbol{\Omega}^{o}
)$, and $\beta^{o}_{0}(\varepsilon_{1},\varepsilon_{2},\boldsymbol{\Omega}^{o}
)$, $\varepsilon_{1},\varepsilon_{2} \! = \! \pm 1$, as in the Proposition,
recalling that $\boldsymbol{\omega}^{o} \! = \! (\omega_{1}^{o},\omega_{2}^{
o},\dotsc,\omega_{N}^{o})$, and that the associated $N \! \times \! N$ Riemann
matrix of $\boldsymbol{\beta}^{o}$-periods, $\tau^{o} \! = \! (\tau^{o})_{i,j
=1,\dotsc,N} \! := \! (\oint_{\boldsymbol{\beta}^{o}_{j}} \omega^{o}_{i})_{i,
j=1,\dotsc,N}$, is non-degenerate, symmetric, and $-\mi \tau^{o}$ is positive
definite, via the above asymptotic (as $z \! \to \! 0^{+})$ expansion for
$\int_{0^{+}}^{z} \omega^{o}_{j}$, $j \! = \! 1,\dotsc,N$, one shows that
\begin{equation*}
\dfrac{\boldsymbol{\theta}^{o}(\varepsilon_{1} \boldsymbol{u}^{o}(z) \! - \!
\frac{1}{2 \pi}(n \! + \! \frac{1}{2}) \boldsymbol{\Omega}^{o} \! + \!
\varepsilon_{2} \boldsymbol{d}_{o})}{\boldsymbol{\theta}^{o}(\varepsilon_{1}
\boldsymbol{u}^{o}(z) \! + \! \varepsilon_{2} \boldsymbol{d}_{o})} \underset{
z \to 0^{+}}{=} \digamma_{0}^{o} \! + \! \digamma_{1}^{o}z \! + \! \digamma_{
2}^{o}z^{2} \! + \! \mathcal{O}(z^{3}),
\end{equation*}
where
\begin{align*}
\digamma_{0}^{o} :=& \, \dfrac{\theta^{o}_{0}(\varepsilon_{1},\varepsilon_{2},
\boldsymbol{\Omega}^{o})}{\theta^{o}_{0}(\varepsilon_{1},\varepsilon_{2},\vec{
\pmb{0}})}, \\
\digamma_{1}^{o} :=& \, \dfrac{\alpha^{o}_{0}(\varepsilon_{1},\varepsilon_{2},
\boldsymbol{\Omega}^{o}) \theta^{o}_{0}(\varepsilon_{1},\varepsilon_{2},\vec{
\pmb{0}}) \! - \! \theta^{o}_{0}(\varepsilon_{1},\varepsilon_{2},\boldsymbol{
\Omega}^{o}) \alpha^{o}_{0}(\varepsilon_{1},\varepsilon_{2},\vec{\pmb{0}})}{
(\theta^{o}_{0}(\varepsilon_{1},\varepsilon_{2},\vec{\pmb{0}}))^{2}}, \\
\digamma^{o}_{2} :=& \left(\theta^{o}_{0}(\varepsilon_{1},\varepsilon_{2},
\boldsymbol{\Omega}^{o}) \! \left((\alpha^{o}_{0}(\varepsilon_{1},\varepsilon_{
2},\vec{\pmb{0}}))^{2} \! - \! \beta^{o}_{0}(\varepsilon_{1},\varepsilon_{2},
\vec{\pmb{0}}) \theta^{o}_{0}(\varepsilon_{1},\varepsilon_{2},\vec{\pmb{0}})
\right) \! - \! \alpha^{o}_{0}(\varepsilon_{1},\varepsilon_{2},\boldsymbol{
\Omega}^{o}) \right. \\
\times&\left. \alpha^{o}_{0}(\varepsilon_{1},\varepsilon_{2},\vec{\pmb{0}})
\theta^{o}_{0}(\varepsilon_{1},\varepsilon_{2},\vec{\pmb{0}}) \! + \! \beta^{
o}_{0}(\varepsilon_{1},\varepsilon_{2},\boldsymbol{\Omega}^{o})(\theta^{o}_{0}
(\varepsilon_{1},\varepsilon_{2},\vec{\pmb{0}}))^{2} \right)(\theta^{o}_{0}
(\varepsilon_{1},\varepsilon_{2},\vec{\pmb{0}}))^{-3},
\end{align*}
with $\vec{\pmb{0}} \! := \! (0,0,\dotsc,0)^{\mathrm{T}}$ $(\in \! \mathbb{
R}^{N})$. Via the above asymptotic (as $z \! \to \! 0^{+})$ expansions for
$\tfrac{1}{2}((\gamma^{o}_{0})^{-1} \gamma^{o}(z) \! + \! \gamma^{o}_{0}
(\gamma^{o}(z))^{-1})$, $\tfrac{1}{2 \mi}((\gamma^{o}_{0})^{-1} \gamma^{o}
(z) \! - \! \gamma^{o}_{0}(\gamma^{o}(z))^{-1})$, and $\tfrac{\boldsymbol{
\theta}^{o}(\varepsilon_{1} \boldsymbol{u}^{o}(z)-\frac{1}{2 \pi}(n+\frac{
1}{2}) \boldsymbol{\Omega}^{o}+\varepsilon_{2} \boldsymbol{d}_{o})}{
\boldsymbol{\theta}^{o}(\varepsilon_{1} \boldsymbol{u}^{o}(z)+\varepsilon_{2}
\boldsymbol{d}_{o})}$, one arrives at, upon recalling the expression for
$\overset{o}{m}^{\raise-1.0ex\hbox{$\scriptstyle \infty$}}(z)$ given in
Lemma 4.5, the asymptotic expansion for
$\overset{o}{m}^{\raise-1.0ex\hbox{$\scriptstyle \infty$}}(z)$ stated in the
Proposition.

Let $\overset{o}{\operatorname{Y}} \colon \mathbb{C} \setminus \mathbb{R}
\! \to \! \operatorname{SL}_{2}(\mathbb{C})$ be the (unique) solution of
\textbf{RHP2}, that is, $(\overset{o}{\operatorname{Y}}(z),\mathrm{I} \! +
\! \exp (-n \widetilde{V}(z)) \sigma_{+},\mathbb{R})$. Recall, also, that, 
for $z \! \in \! \Upsilon^{o}_{1}$ (Figure~7),
\begin{equation*}
\overset{o}{\mathrm{Y}}(z) \! = \! \me^{\frac{n \ell_{o}}{2} \operatorname{ad}
(\sigma_{3})} \mathscr{R}^{o}(z)
\overset{o}{m}^{\raise-1.0ex\hbox{$\scriptstyle \infty$}}(z) \mathbb{E}^{
\sigma_{3}} \me^{n(g^{o}(z)-\mathfrak{Q}^{+}_{\mathscr{A}}) \sigma_{3}},
\end{equation*}
and, for $z \! \in \! \Upsilon^{o}_{2}$ (Figure~7),
\begin{equation*}
\overset{o}{\mathrm{Y}}(z) \! = \! \me^{\frac{n \ell_{o}}{2} \operatorname{ad}
(\sigma_{3})} \mathscr{R}^{o}(z)
\overset{o}{m}^{\raise-1.0ex\hbox{$\scriptstyle \infty$}}(z) \mathbb{E}^{-
\sigma_{3}} \me^{n(g^{o}(z)-\mathfrak{Q}^{-}_{\mathscr{A}}) \sigma_{3}}:
\end{equation*}
consider, say, and without loss of generality, small-$z$ asymptotics for 
$\overset{o}{\operatorname{Y}}(z)$ for $z \! \in \! \Upsilon_{1}^{o}$. 
Recalling {}from Lemma~3.4 that $g^{o}(z) \! := \! \int_{J_{o}} \ln ((z \! 
- \! s)^{2+\frac{1}{n}}(zs)^{-1}) \psi_{V}^{o}(s) \, \md s$, $z \! \in \! 
\mathbb{C} \setminus (-\infty,\max \lbrace 0,a_{N+1}^{o} \rbrace)$, for 
$\vert z \vert \! \ll \! \min_{j=1,\dotsc,N+1} \lbrace \vert b_{j-1}^{o} \! 
- \! a_{j}^{o} \vert \rbrace$, in particular, $\vert z/s \vert \! \ll \! 1$ 
with $s \! \in \! J_{o}$, and noting that $\int_{J_{o}} \psi_{V}^{o}(s) \, 
\md s \! = \! 1$ and $\int_{J_{o}}s^{-m} \psi_{V}^{o}(s) \, \md s \! < \! 
\infty$, $m \! \in \! \mathbb{N}$, via the expansions $\tfrac{1}{z-s} \! = \! 
-\sum_{k=0}^{l} \tfrac{z^{k}}{s^{k+1}} \! + \! \tfrac{z^{l+1}}{s^{l+1}(z-s)}$, 
$l \! \in \! \mathbb{Z}_{0}^{+}$, and $\ln (s \! - \! z) \! =_{\vert z \vert 
\to 0} \! \ln (s) \! - \! \sum_{k=1}^{\infty} \tfrac{1}{k}(\tfrac{z}{s})^{k}$, 
one shows that
\begin{align*}
g^{o}(z) \underset{\mathbb{C}_{\pm} \ni z \to 0}{=}& \, -\ln (z) \! + \! 
\left(1 \! + \! \dfrac{1}{n} \right) \! \int_{J_{o}} \ln (\lvert s \rvert) 
\psi_{V}^{o}(s) \, \md s \! - \! \mi \pi \int_{J_{o} \cap \mathbb{R}_{-}} 
\psi_{V}^{o}(s) \, \md s \! \pm \! \mi \pi \! \left(2 \! + \! \dfrac{1}{n} 
\right) \! \int_{J_{o} \cap \mathbb{R}_{+}} \psi_{V}^{o}(s) \, \md s \\
+& \, z \! \left(-\! \left(2 \! + \! \dfrac{1}{n} \right) \int_{J_{o}}s^{-1}
\psi_{V}^{o}(s) \, \md s \right) \! + \! z^{2} \! \left(-\dfrac{1}{2} \! 
\left(2 \! + \! \dfrac{1}{n} \right) \! \int_{J_{o}}s^{-2} \psi_{V}^{o}(s) 
\, \md s \right) \! + \! \mathcal{O}(z^{3}),
\end{align*}
where $\int_{J_{o} \cap \mathbb{R}_{\pm}} \psi_{V}^{o}(s) \, \md s$ are given
in Lemma~3.4. (Explicit expressions for $\int_{J_{o}}s^{-k} \psi_{V}^{o}(s) \,
\md s$, $k \! = \! 1,2$, are given in Remark~3.2.) Using the asymptotic (as
$z \! \to \! 0)$ expansions for $g^{o}(z)$, $\mathscr{R}^{o}(z)$, and
$\overset{o}{m}^{\raise-1.0ex\hbox{$\scriptstyle \infty$}}(z)$ derived above,
upon recalling the formula for $\overset{o}{\operatorname{Y}}(z)$, one arrives
at, after a matrix-multiplication argument, the asymptotic expansion for
$\overset{o}{\operatorname{Y}}(z)z^{n \sigma_{3}}$ stated in the Proposition.
\hfill $\qed$
\begin{bbbbb}
Let $\overset{o}{\operatorname{Y}} \colon \mathbb{C} \setminus \mathbb{R} \!
\to \! \operatorname{SL}_{2}(\mathbb{C})$ be the solution of {\rm \pmb{RHP2}}
with $z$ $(\in \! \mathbb{C} \setminus \mathbb{R})$ $\to \! 0$ asymptotics
given in Proposition~{\rm 5.3}. Then,
\begin{equation*}
\xi^{(2n+1)}_{-n-1} \! = \! \dfrac{1}{\norm{\boldsymbol{\pi}_{2n+1}(\pmb{
\cdot})}_{\mathscr{L}}} \! = \sqrt{\dfrac{H^{(-2n)}_{2n+1}}{H^{(-2n-2)}_{2
n+2}}} = \! \left(\dfrac{1}{2 \pi \mi (\operatorname{Y}^{o,0}_{1})_{12}}
\right)^{1/2} \quad (> \! 0),
\end{equation*}
where $(\operatorname{Y}^{o,0}_{1})_{12} \! = \! \me^{n \ell_{o}} \! \left(
(\overset{o}{m}_{1}^{\raise-1.0ex\hbox{$\scriptstyle 0$}})_{12} \mathbb{
E}^{-1} \! + \! (\mathscr{R}^{o,0}_{1}(n))_{12} \right)$, with
$(\overset{o}{m}_{1}^{\raise-1.0ex\hbox{$\scriptstyle 0$}})_{12}$ and
$(\mathscr{R}^{o,0}_{1}(n))_{12}$ given in Proposition~{\rm 5.3}.
\end{bbbbb}

\emph{Proof.} Recall {}from Lemma~2.2.2 that $z \boldsymbol{\pi}_{2n+1}(z) \!
:= \! (\overset{o}{\operatorname{Y}}(z))_{11}$ and $(\overset{o}{\operatorname{
Y}}(z))_{12} \! = \! z \int_{\mathbb{R}} \tfrac{(s \boldsymbol{\pi}_{2n+1}(s))
\me^{-n \widetilde{V}(s)}}{s(s-z)} \linebreak[4]
\cdot \tfrac{\md s}{2 \pi \mi}$. Using (for $\vert z/s \vert \! \ll \! 1)$ the
expansion $\tfrac{1}{z-s} \! = \! -\sum_{k=0}^{l} \tfrac{z^{k}}{s^{k+1}} \!
+ \! \tfrac{z^{l+1}}{s^{l+1}(z-s)}$, $l \! \in \! \mathbb{Z}_{0}^{+}$, and
recalling that $\langle \boldsymbol{\pi}_{2n+1},z^{j} \rangle_{\mathscr{L}}
\! = \! 0$, $j \! = \! -n,\dotsc,n$, and $\phi_{2n+1}(z) \! = \! \xi^{(2n+1)
}_{-n-1}(z) \boldsymbol{\pi}_{2n+1}(z)$, one proceeds as follows:
\begin{align*}
\left(\overset{o}{\operatorname{Y}}(z) \right)_{12} \underset{\underset{z \in
\mathbb{C} \setminus \mathbb{R}}{z \to 0}}{=}& \, z \int_{\mathbb{R}} \dfrac{
\boldsymbol{\pi}_{2n+1}(s)}{s} \! \left(1 \! + \! \dfrac{z}{s} \! + \! \cdots
\! + \! \dfrac{z^{n-1}}{s^{n-1}} \! + \! \dfrac{z^{n}}{s^{n}} \! + \! \cdots
\right) \! \me^{-n \widetilde{V}(s)} \, \dfrac{\md s}{2 \pi \mi} \\
\underset{\underset{z \in \mathbb{C} \setminus \mathbb{R}}{z \to 0}}{=}& \, z
\int_{\mathbb{R}} \boldsymbol{\pi}_{2n+1}(s) \! \left(\dfrac{z^{n}}{s^{n+1}}
\right) \! \me^{-n \widetilde{V}(s)} \, \dfrac{\md s}{2 \pi \mi} \! + \!
\mathcal{O} \! \left(z^{n+2} \right) \\
\underset{\underset{z \in \mathbb{C} \setminus \mathbb{R}}{z \to 0}}{=}&
\, \dfrac{z^{n+1}}{\xi^{(2n+1)}_{-n-1}} \int_{\mathbb{R}} \underbrace{\xi_{-
n-1}^{(2n+1)} \boldsymbol{\pi}_{2n+1}(s)}_{= \, \phi_{2n+1}(s)} \dfrac{\me^{-
n \widetilde{V}(s)}}{\xi_{-n-1}^{(2n+1)}} \! \underbrace{\left(\xi_{n}^{(2n+1)
}s^{n} \! + \! \cdots \! + \! \dfrac{\xi_{-n-1}^{(2n+1)}}{s^{n+1}} \right)}_{
= \, \phi_{2n+1}(s)} \dfrac{\md s}{2 \pi \mi} \\
+& \, \mathcal{O} \! \left(z^{n+2} \right) \\
\underset{\underset{z \in \mathbb{C} \setminus \mathbb{R}}{z \to 0}}{=}&
\, \dfrac{z^{n+1}}{2 \pi \mi (\xi_{-n-1}^{(2n+1)})^{2}} \underbrace{\int_{
\mathbb{R}} \phi_{2n+1}(s) \phi_{2n+1}(s) \me^{-n \widetilde{V}(s)} \, \md
s}_{= \, 1} + \, \mathcal{O} \! \left(z^{n+2} \right) \Rightarrow \\
\left(\overset{o}{\operatorname{Y}}(z)z^{n \sigma_{3}} \right)_{12} \underset{
\underset{z \in \mathbb{C} \setminus \mathbb{R}}{z \to 0}}{=}& \, z \! \left(
\dfrac{1}{2 \pi \mi (\xi_{-n-1}^{(2n+1)})^{2}} \right) \! + \! \mathcal{O}
(z^{2});
\end{align*}
but, noting {}from Proposition~5.3 that
\begin{equation*}
\left(\overset{o}{\operatorname{Y}}(z)z^{n \sigma_{3}} \right)_{12} \underset{
\underset{z \in \mathbb{C} \setminus \mathbb{R}}{z \to 0}}{=} z \! \left(
\operatorname{Y}^{o,0}_{1} \right)_{12} \! + \!\mathcal{O}(z^{2}),
\end{equation*}
upon equating the above two asymptotic expansions for $(\overset{o}{
\operatorname{Y}}(z)z^{n \sigma_{3}})_{12}$, one arrives at the result 
stated in the Proposition. \hfill $\qed$

Using the results of Propositions~5.3 and~5.4, one obtains $n \! \to \! 
\infty$ asymptotics for $\xi^{(2n+1)}_{-n-1}$ and $\phi_{2n+1}(z)$ (in 
the entire complex plane) stated in Theorem~2.3.2.

Large-$z$ asymptotics for $\overset{o}{\mathrm{Y}}(z)$ are given in the 
Appendix (see Lemma~A.1): these latter asymptotics are necessary for the 
results of \cite{a40}.

\vspace*{1.00cm}
\textbf{\Large Acknowledgements}

K.~T.-R.~McLaughlin was supported, in part, by National Science Foundation 
Grant Nos.~DMS--9970328 and~DMS--0200749. X.~Zhou was supported, in part, 
by National Science Foundation Grant No.~DMS--0300844.
\clearpage
\section*{Appendix: Large-$z$ Asymptotics for $\overset{o}{\mathrm{Y}}(z)$}
\setcounter{section}{1}
\setcounter{z0}{1}
Even though the results of Lemma~A.1 below, namely, large-$z$ asymptotics (as 
$(\mathbb{C} \setminus \mathbb{R} \! \ni)$ $z \! \to \! \infty)$ of $\overset{
o}{\operatorname{Y}}(z)$, are not necessary in order to prove Theorems~2.3.1 
and~2.3.2, they are essential for the results of \cite{a40}, related to 
asymptotics of the coefficients of the system of three- and five-term 
recurrence relations and the corresponding Laurent-Jacobi matrices (cf. 
Section 1). For the sake of completeness, therefore, and in order to eschew 
any duplication of the analysis of this paper, $(\mathbb{C} \setminus 
\mathbb{R} \! \ni)$ $z \! \to \! \infty$ asymptotics for $\overset{o}{
\operatorname{Y}}(z)$ are presented here.
\begin{ay}
Let $\mathscr{R}^{o} \colon \mathbb{C} \setminus \widetilde{\Sigma}_{p}^{o}
\! \to \! \operatorname{SL}_{2}(\mathbb{C})$ be the solution of the {\rm RHP}
$(\mathscr{R}^{o}(z),\upsilon_{\mathscr{R}}^{o}(z),\widetilde{\Sigma}_{p}^{o}
)$ formulated in Proposition~{\rm 5.2} with $n \! \to \! \infty$ asymptotics
given in Lemma~{\rm 5.3}. Then,
\begin{equation*}
\mathscr{R}^{o}(z) \underset{z \to \infty}{=} \mathrm{I} \! + \! \mathscr{
R}^{o,\infty}_{0}(n) \! + \! \dfrac{1}{z} \mathscr{R}^{o,\infty}_{1}(n) \!
+ \! \dfrac{1}{z^{2}} \mathscr{R}^{o,\infty}_{2}(n) \! + \! \mathcal{O}
(z^{-3}),
\end{equation*}
where, for $k \! = \! -1,0,1$,
\begin{equation*}
\mathscr{R}^{o,\infty}_{k+1}(n) \! := \! -\int_{\Sigma^{o}_{\circlearrowright}}
s^{k}w^{\Sigma^{o}_{\circlearrowright}}_{+}(s) \, \dfrac{\md s}{2 \pi \mi}
\! = \! \sum_{j=1}^{N+1} \, \sum_{q \in \{b_{j-1}^{o},a_{j}^{o}\}}
\operatorname{Res} \! \left(z^{k}w^{\Sigma^{o}_{\circlearrowright}}_{+}(z);
q \right),
\end{equation*}
with, in particular,
\begin{align*}
\mathscr{R}^{o,\infty}_{0}(n) \underset{n \to \infty}{=}& \, \dfrac{1}{(n \! +
\! \frac{1}{2})} \sum_{j=1}^{N+1} \! \left(\dfrac{(\mathscr{B}^{o}(b_{j-1}^{o}
) \widehat{\alpha}_{0}^{o}(b_{j-1}^{o}) \! - \! \mathscr{A}^{o}(b_{j-1}^{o})
(\widehat{\alpha}_{1}^{o}(b_{j-1}^{o}) \! + \! (b_{j-1}^{o})^{-1} \widehat{
\alpha}_{0}^{o}(b_{j-1}^{o})))}{b_{j-1}^{o}(\widehat{\alpha}_{0}^{o}(b_{j-1}^{
o}))^{2}} \right. \\
+&\left. \, \dfrac{(\mathscr{B}^{o}(a_{j}^{o}) \widehat{\alpha}_{0}^{o}(a_{j}^{
o}) \! - \! \mathscr{A}^{o}(a_{j}^{o})(\widehat{\alpha}_{1}^{o}(a_{j}^{o}) \!
+ \! (a_{j}^{o})^{-1} \widehat{\alpha}_{0}^{o}(a_{j}^{o})))}{a_{j}^{o}
(\widehat{\alpha}_{0}^{o}(a_{j}^{o}))^{2}} \right) \! + \! \mathcal{O} \!
\left(\dfrac{1}{(n \! + \! \frac{1}{2})^{2}} \right), \\
\mathscr{R}^{o,\infty}_{1}(n) \underset{n \to \infty}{=}& \, \dfrac{1}{(n \! +
\! \frac{1}{2})} \sum_{j=1}^{N+1} \! \left(\dfrac{(\mathscr{B}^{o}(b_{j-1}^{o}
) \widehat{\alpha}_{0}^{o}(b_{j-1}^{o}) \! - \! \mathscr{A}^{o}(b_{j-1}^{o})
\widehat{\alpha}_{1}^{o}(b_{j-1}^{o}))}{(\widehat{\alpha}_{0}^{o}(b_{j-1}^{o})
)^{2}} \right. \\
+&\left. \, \dfrac{(\mathscr{B}^{o}(a_{j}^{o}) \widehat{\alpha}_{0}^{o}(a_{j}^{
o}) \! - \! \mathscr{A}^{o}(a_{j}^{o}) \widehat{\alpha}_{1}^{o}(a_{j}^{o}))}{
(\widehat{\alpha}_{0}^{o}(a_{j}^{o}))^{2}} \right) \! + \! \mathcal{O} \!
\left(\dfrac{1}{(n \! + \! \frac{1}{2})^{2}} \right), \\
\mathscr{R}^{o,\infty}_{2}(n) \underset{n \to \infty}{=}& \, \dfrac{1}{(n \! +
\! \frac{1}{2})} \sum_{j=1}^{N+1} \! \left(\dfrac{(\mathscr{B}^{o}(b_{j-1}^{o}
) \widehat{\alpha}_{0}^{o}(b_{j-1}^{o})b_{j-1}^{o} \! - \! \mathscr{A}^{o}(b_{
j-1}^{o})(b_{j-1}^{o} \widehat{\alpha}_{1}^{o}(b_{j-1}^{o}) \! - \! \widehat{
\alpha}_{0}^{o}(b_{j-1}^{o})))}{(\widehat{\alpha}_{0}^{o}(b_{j-1}^{o}))^{2}}
\right. \\
+&\left. \, \dfrac{(\mathscr{B}^{o}(a_{j}^{o}) \widehat{\alpha}_{0}^{o}(a_{j}^{
o})a_{j}^{o} \! - \! \mathscr{A}^{o}(a_{j}^{o})(a_{j}^{o} \widehat{\alpha}_{
1}^{o}(a_{j}^{o}) \! - \! \widehat{\alpha}_{0}^{o}(a_{j}^{o})))}{(\widehat{
\alpha}_{0}^{o}(a_{j}^{o}))^{2}} \right) \! + \! \mathcal{O} \! \left(\dfrac{
1}{(n \! + \! \frac{1}{2})^{2}} \right),
\end{align*}
and all parameters are defined in Lemma~{\rm 5.3}.

Let $\overset{o}{m}^{\raise-1.0ex\hbox{$\scriptstyle \infty$}} \colon \mathbb{
C} \setminus J_{o}^{\infty} \! \to \! \operatorname{SL}_{2}(\mathbb{C})$ solve
the {\rm RHP} $(\overset{o}{m}^{\raise-1.0ex\hbox{$\scriptstyle \infty$}}(z),
J_{o}^{\infty},
\overset{o}{\upsilon}^{\raise-1.0ex\hbox{$\scriptstyle \infty$}}(z))$
formulated in Lemma~{\rm 4.3} with (unique) solution given by Lemma~{\rm 4.5}.
For $\varepsilon_{1},\varepsilon_{2} \! = \! \pm 1$, set
\begin{gather*}
\theta^{o}_{\infty}(\varepsilon_{1},\varepsilon_{2},\boldsymbol{\Omega}^{o})
\! := \! \boldsymbol{\theta}^{o}(\varepsilon_{1} \boldsymbol{u}^{o}_{+}
(\infty) \! - \! \tfrac{1}{2 \pi}(n \! + \! \tfrac{1}{2}) \boldsymbol{\Omega}^{
o} \! + \! \varepsilon_{2} \boldsymbol{d}_{o}),
\end{gather*}
where\footnote{For $\mathscr{P} \! := \! (y,z) \! \in \! \lbrace \mathstrut
(z_{1},z_{2}); \, z_{1}^{2} \! = \! R_{o}(z_{2}) \rbrace$, $\mathscr{P} \!
\to \! \infty^{\pm} \! \Leftrightarrow \! z \! \to \! \infty$, $y \! \sim \!
\pm z^{N+1}$.} $\boldsymbol{u}^{o}_{+}(\infty) \! = \! \int_{a_{N+1}^{o}}^{
\infty^{+}} \boldsymbol{\omega}^{o}$,
\begin{gather*}
\widetilde{\alpha}^{o}_{\infty}(\varepsilon_{1},\varepsilon_{2},\boldsymbol{
\Omega}^{o}) \! := \! 2 \pi \mi \varepsilon_{1} \sum_{m \in \mathbb{Z}^{N}}(m,
\widehat{\boldsymbol{\alpha}}^{o}_{\infty}) \me^{2 \pi \mi (m,\varepsilon_{1}
\boldsymbol{u}^{o}_{+}(\infty)-\frac{1}{2 \pi}(n+\frac{1}{2}) \boldsymbol{
\Omega}^{o}+\varepsilon_{2} \boldsymbol{d}_{o})+ \pi \mi (m,\tau^{o}m)},
\end{gather*}
where $\widehat{\boldsymbol{\alpha}}^{o}_{\infty} \! = \! (\widehat{\alpha}_{
\infty,1}^{o},\widehat{\alpha}^{o}_{\infty,2},\dotsc,\widehat{\alpha}^{o}_{
\infty,N})$, with $\widehat{\alpha}^{o}_{\infty,j} \! := \! c_{j1}^{o}$, $j \!
= \! 1,\dotsc,N$,
\begin{gather*}
\beta^{o}_{\infty}(\varepsilon_{1},\varepsilon_{2},\boldsymbol{\Omega}^{o}) \!
:= \! 2 \pi \sum_{m \in \mathbb{Z}^{N}} \! \left(\mi \varepsilon_{1}(m,
\widehat{\boldsymbol{\beta}}^{o}_{\infty}) \! + \!  \pi (m,\widehat{
\boldsymbol{\alpha}}^{o}_{\infty})^{2} \right) \! \me^{2 \pi \mi (m,
\varepsilon_{1} \boldsymbol{u}^{o}_{+}(\infty)-\frac{1}{2 \pi}(n+\frac{1}{2})
\boldsymbol{\Omega}^{o}+\varepsilon_{2} \boldsymbol{d}_{o})+ \pi \mi (m,\tau^{
o}m)},
\end{gather*}
where $\widehat{\boldsymbol{\beta}}^{o}_{\infty} \! = \! (\widehat{\beta}^{o}_{
\infty,1},\widehat{\beta}^{o}_{\infty,2},\dotsc,\widehat{\beta}^{o}_{\infty,N}
)$, with $\widehat{\beta}^{o}_{\infty,j} \! := \! \tfrac{1}{2}(c^{o}_{j2} \! +
\! \tfrac{1}{2}c_{j1}^{o} \sum_{i=1}^{N+1}(b_{i-1}^{o} \! + \! a_{i}^{o}))$,
$j \! = \! 1,\dotsc,N$, where $c_{j1}^{o},c_{j2}^{o}$, $j \! = \! 1,\dotsc,N$,
are obtained {}from Equations~{\rm (O1)} and~{\rm (O2)}. Set $\gamma^{o}_{0}
\! := \! \gamma^{o}(0) \! = \! (\prod_{k=1}^{N+1}b_{k-1}^{o}(a_{k}^{o})^{-1}
)^{1/4}$ $(> \! 0)$. Then,
\begin{gather*}
\overset{o}{m}^{\raise-1.0ex\hbox{$\scriptstyle \infty$}}(z) \underset{z \to
\infty}{=} \overset{o}{m}_{0}^{\raise-1.0ex\hbox{$\scriptstyle \infty$}} \! +
\! \dfrac{1}{z} \overset{o}{m}_{1}^{\raise-1.0ex\hbox{$\scriptstyle \infty$}}
\! + \! \dfrac{1}{z^{2}}
\overset{o}{m}_{2}^{\raise-1.0ex\hbox{$\scriptstyle \infty$}} \! + \!
\mathcal{O}(z^{-3}),
\end{gather*}
where
\begin{align*}
(\overset{o}{m}_{0}^{\raise-1.0ex\hbox{$\scriptstyle \infty$}})_{11} &= \dfrac{
\boldsymbol{\theta}^{o}(\boldsymbol{u}^{o}_{+}(0) \! + \! \boldsymbol{d}_{o})
\mathbb{E}^{-1}}{\boldsymbol{\theta}^{o}(\boldsymbol{u}^{o}_{+}(0) \! - \!
\frac{1}{2 \pi}(n \! + \! \frac{1}{2}) \boldsymbol{\Omega}^{o} \! + \!
\boldsymbol{d}_{o})} \! \left(\dfrac{\gamma^{o}_{0} \! + \! (\gamma^{o}_{0})^{-
1}}{2} \right) \! \dfrac{\theta^{o}_{\infty}(1,1,\boldsymbol{\Omega}^{o})}{
\theta^{o}_{\infty}(1,1,\vec{\pmb{0}})}, \\
(\overset{o}{m}_{0}^{\raise-1.0ex\hbox{$\scriptstyle \infty$}})_{12} &= \dfrac{
\boldsymbol{\theta}^{o}(\boldsymbol{u}^{o}_{+}(0) \! + \! \boldsymbol{d}_{o})
\mathbb{E}^{-1}}{\boldsymbol{\theta}^{o}(\boldsymbol{u}^{o}_{+}(0) \! - \!
\frac{1}{2 \pi}(n \! + \! \frac{1}{2}) \boldsymbol{\Omega}^{o} \! + \!
\boldsymbol{d}_{o})} \! \left(\dfrac{\gamma^{o}_{0} \! - \! (\gamma^{o}_{0})^{
-1}}{2 \mi} \right) \! \dfrac{\theta^{o}_{\infty}(-1,1,\boldsymbol{\Omega}^{o}
)}{\theta^{o}_{\infty}(-1,1,\vec{\pmb{0}})}, \\
(\overset{o}{m}_{0}^{\raise-1.0ex\hbox{$\scriptstyle \infty$}})_{21} &=
-\dfrac{\boldsymbol{\theta}^{o}(\boldsymbol{u}^{o}_{+}(0) \! + \! \boldsymbol{
d}_{o}) \mathbb{E}}{\boldsymbol{\theta}^{o}(-\boldsymbol{u}^{o}_{+}(0) \! -
\! \frac{1}{2 \pi}(n \! + \! \frac{1}{2}) \boldsymbol{\Omega}^{o} \! - \!
\boldsymbol{d}_{o})} \! \left(\dfrac{\gamma^{o}_{0} \! - \! (\gamma^{o}_{0})^{-
1}}{2 \mi} \right) \! \dfrac{\theta^{o}_{\infty}(1,-1,\boldsymbol{\Omega}^{o})
}{\theta^{o}_{\infty}(1,-1,\vec{\pmb{0}})}, \\
(\overset{o}{m}_{0}^{\raise-1.0ex\hbox{$\scriptstyle \infty$}})_{22} &= \dfrac{
\boldsymbol{\theta}^{o}(\boldsymbol{u}^{o}_{+}(0) \! + \! \boldsymbol{d}_{o})
\mathbb{E}}{\boldsymbol{\theta}^{o}(-\boldsymbol{u}^{o}_{+}(0) \! - \! \frac{
1}{2 \pi}(n \! + \! \frac{1}{2}) \boldsymbol{\Omega}^{o} \! - \! \boldsymbol{
d}_{o})} \! \left(\dfrac{\gamma^{o}_{0} \! + \! (\gamma^{o}_{0})^{-1}}{2}
\right) \! \dfrac{\theta^{o}_{\infty}(-1,-1,\boldsymbol{\Omega}^{o})}{\theta^{
o}_{\infty}(-1,-1,\vec{\pmb{0}})}, \\
(\overset{o}{m}_{1}^{\raise-1.0ex\hbox{$\scriptstyle \infty$}})_{11} &= \dfrac{
\boldsymbol{\theta}^{o}(\boldsymbol{u}^{o}_{+}(0) \! + \! \boldsymbol{d}_{o})
\mathbb{E}^{-1}}{\boldsymbol{\theta}^{o}(\boldsymbol{u}^{o}_{+}(0) \! - \!
\frac{1}{2 \pi}(n \! + \! \frac{1}{2}) \boldsymbol{\Omega}^{o} \! + \!
\boldsymbol{d}_{o})} \! \left(\! \left(\dfrac{\widetilde{\alpha}_{\infty}^{o}
(1,1,\vec{\pmb{0}}) \theta^{o}_{\infty}(1,1,\boldsymbol{\Omega}^{o}) \! - \!
\widetilde{\alpha}_{\infty}^{o}(1,1,\boldsymbol{\Omega}^{o}) \theta^{o}_{
\infty}(1,1,\vec{\pmb{0}})}{(\theta^{o}_{\infty}(1,1,\vec{\pmb{0}}))^{2}}
\right) \right. \\
&\left. \times \left(\dfrac{\gamma^{o}_{0} \! + \! (\gamma^{o}_{0})^{-1}}{2}
\right) \! - \! \left(\dfrac{\gamma^{o}_{0} \! - \! (\gamma^{o}_{0})^{-1}}{8}
\right) \! \left(\sum_{k=1}^{N+1} \! \left(a_{k}^{o} \! - \! b_{k-1}^{o}
\right) \right) \! \dfrac{\theta^{o}_{\infty}(1,1,\boldsymbol{\Omega}^{o})}{
\theta^{o}_{\infty}(1,1,\vec{\pmb{0}})} \right), \\
(\overset{o}{m}_{1}^{\raise-1.0ex\hbox{$\scriptstyle \infty$}})_{12} &= \dfrac{
\boldsymbol{\theta}^{o}(\boldsymbol{u}^{o}_{+}(0) \! + \! \boldsymbol{d}_{o})
\mathbb{E}^{-1}}{\boldsymbol{\theta}^{o}(\boldsymbol{u}^{o}_{+}(0) \! - \!
\frac{1}{2 \pi}(n \! + \! \frac{1}{2}) \boldsymbol{\Omega}^{o} \! + \!
\boldsymbol{d}_{o})} \! \left(\! \left(\dfrac{\widetilde{\alpha}_{\infty}^{o}
(-1,1,\vec{\pmb{0}}) \theta^{o}_{\infty}(-1,1,\boldsymbol{\Omega}^{o}) \! -
\! \widetilde{\alpha}_{\infty}^{o}(-1,1,\boldsymbol{\Omega}^{o}) \theta^{o}_{
\infty}(-1,1,\vec{\pmb{0}})}{(\theta^{o}_{\infty}(-1,1,\vec{\pmb{0}}))^{2}}
\right) \right. \\
&\left. \times \left(\dfrac{\gamma^{o}_{0} \! - \! (\gamma^{o}_{0})^{-1}}{2
\mi} \right) \! - \! \left(\dfrac{\gamma^{o}_{0} \! + \! (\gamma^{o}_{0})^{-1}
}{8 \mi} \right) \! \left(\sum_{k=1}^{N+1} \! \left(a_{k}^{o} \! - \! b_{k-1}^{
o} \right) \right) \! \dfrac{\theta^{o}_{\infty}(-1,1,\boldsymbol{\Omega}^{o})
}{\theta^{o}_{\infty}(-1,1,\vec{\pmb{0}})} \right), \\
(\overset{o}{m}_{1}^{\raise-1.0ex\hbox{$\scriptstyle \infty$}})_{21} &=
-\dfrac{\boldsymbol{\theta}^{o}(\boldsymbol{u}^{o}_{+}(0) \! + \! \boldsymbol{
d}_{o}) \mathbb{E}}{\boldsymbol{\theta}^{o}(-\boldsymbol{u}^{o}_{+}(0) \! -
\! \frac{1}{2 \pi}(n \! + \! \frac{1}{2}) \boldsymbol{\Omega}^{o} \! - \!
\boldsymbol{d}_{o})} \! \left(\! \left(\dfrac{\widetilde{\alpha}_{\infty}^{o}
(1,-1,\vec{\pmb{0}}) \theta^{o}_{\infty}(1,-1,\boldsymbol{\Omega}^{o}) \! -
\! \widetilde{\alpha}_{\infty}^{o}(1,-1,\boldsymbol{\Omega}^{o}) \theta^{o}_{
\infty}(1,-1,\vec{\pmb{0}})}{(\theta^{o}_{\infty}(1,-1,\vec{\pmb{0}}))^{2}}
\right) \right. \\
&\left. \times \left(\dfrac{\gamma^{o}_{0} \! - \! (\gamma^{o}_{0})^{-1}}{2
\mi} \right) \! - \! \left(\dfrac{\gamma^{o}_{0} \! + \! (\gamma^{o}_{0})^{-1}
}{8 \mi} \right) \! \left(\sum_{k=1}^{N+1} \! \left(a_{k}^{o} \! - \! b_{k-1}^{
o} \right) \right) \! \dfrac{\theta^{o}_{\infty}(1,-1,\boldsymbol{\Omega}^{o}
)}{\theta^{o}_{\infty}(1,-1,\vec{\pmb{0}})} \right), \\
(\overset{o}{m}_{1}^{\raise-1.0ex\hbox{$\scriptstyle \infty$}})_{22} &= \dfrac{
\boldsymbol{\theta}^{o}(\boldsymbol{u}^{o}_{+}(0) \! + \! \boldsymbol{d}_{o})
\mathbb{E}}{\boldsymbol{\theta}^{o}(-\boldsymbol{u}^{o}_{+}(0) \! - \! \frac{
1}{2 \pi}(n \! + \! \frac{1}{2}) \boldsymbol{\Omega}^{o} \! - \! \boldsymbol{
d}_{o})} \! \left(\! \left(\dfrac{\widetilde{\alpha}_{\infty}^{o}(-1,-1,\vec{
\pmb{0}}) \theta^{o}_{\infty}(-1,-1,\boldsymbol{\Omega}^{o}) \! - \!
\widetilde{\alpha}_{\infty}^{o}(-1,-1,\boldsymbol{\Omega}^{o}) \theta^{o}_{
\infty}(-1,-1,\vec{\pmb{0}})}{(\theta^{o}_{\infty}(-1,-1,\vec{\pmb{0}}))^{2}}
\right) \right. \\
&\left. \times \left(\dfrac{\gamma^{o}_{0} \! + \! (\gamma^{o}_{0})^{-1}}{2}
\right) \! - \! \left(\dfrac{\gamma^{o}_{0} \! - \! (\gamma^{o}_{0})^{-1}}{8}
\right) \! \left(\sum_{k=1}^{N+1} \! \left(a_{k}^{o} \! - \! b_{k-1}^{o}
\right) \right) \! \dfrac{\theta^{o}_{\infty}(-1,-1,\boldsymbol{\Omega}^{o})}{
\theta^{o}_{\infty}(-1,-1,\vec{\pmb{0}})} \right), \\
(\overset{o}{m}_{2}^{\raise-1.0ex\hbox{$\scriptstyle \infty$}})_{11} &= \dfrac{
\boldsymbol{\theta}^{o}(\boldsymbol{u}^{o}_{+}(0) \! + \! \boldsymbol{d}_{o})
\mathbb{E}^{-1}}{\boldsymbol{\theta}^{o}(\boldsymbol{u}^{o}_{+}(0) \! - \!
\frac{1}{2 \pi}(n \! + \! \frac{1}{2}) \boldsymbol{\Omega}^{o} \! + \!
\boldsymbol{d}_{o})} \! \left(\! \left(\theta^{o}_{\infty}(1,1,\boldsymbol{
\Omega}^{o}) \! \left((\widetilde{\alpha}_{\infty}^{o}(1,1,\vec{\pmb{0}}))^{2}
\! + \! \beta^{o}_{\infty}(1,1,\vec{\pmb{0}}) \theta^{o}_{\infty}(1,1,\vec{
\pmb{0}}) \right) \right. \right. \\
&\left. \left. - \, \widetilde{\alpha}_{\infty}^{o}(1,1,\boldsymbol{\Omega}^{
o}) \widetilde{\alpha}_{\infty}^{o}(1,1,\vec{\pmb{0}}) \theta^{o}_{\infty}(1,1,
\vec{\pmb{0}}) \! - \! \beta^{o}_{\infty}(1,1,\boldsymbol{\Omega}^{o})(\theta^{
o}_{\infty}(1,1,\vec{\pmb{0}}))^{2} \right) \! \tfrac{1}{(\theta^{o}_{\infty}
(1,1,\vec{\pmb{0}}))^{3}} \right. \\
&\left. \times \left(\dfrac{\gamma^{o}_{0} \! + \! (\gamma^{o}_{0})^{-1}}{2}
\right) \! + \! \left(\dfrac{\widetilde{\alpha}_{\infty}^{o}(1,1,\vec{\pmb{0}}
) \theta^{o}_{\infty}(1,1,\boldsymbol{\Omega}^{o}) \! - \! \widetilde{\alpha}_{
\infty}^{o}(1,1,\boldsymbol{\Omega}^{o}) \theta^{o}_{\infty}(1,1,\vec{\pmb{0}}
)}{(\theta^{o}_{\infty}(1,1,\vec{\pmb{0}}))^{2}} \right) \right. \\
&\left. \times \left(\sum_{k=1}^{N+1} \! \left(b_{k-1}^{o} \! - \! a_{k}^{o}
\right) \right) \! \left(\dfrac{\gamma^{o}_{0} \! - \! (\gamma^{o}_{0})^{-1}}{
8} \right) \! + \! \left(\! \left(\dfrac{\gamma^{o}_{0} \! - \! (\gamma^{o}_{
0})^{-1}}{16} \right) \! \sum_{k=1}^{N+1} \! \left((b_{k-1}^{o})^{2} \! - \!
(a_{k}^{o})^{2} \right) \right. \right. \\
&\left. \left. +\left(\dfrac{\gamma^{o}_{0} \! + \! (\gamma^{o}_{0})^{-1}}{64}
\right) \! \left(\sum_{k=1}^{N+1} \! \left(a_{k}^{o} \! - \! b_{k-1}^{o}
\right) \right)^{2} \right) \! \dfrac{\theta^{o}_{\infty}(1,1,\boldsymbol{
\Omega}^{o})}{\theta^{o}_{\infty}(1,1,\vec{\pmb{0}})} \right), \\
(\overset{o}{m}_{2}^{\raise-1.0ex\hbox{$\scriptstyle \infty$}})_{12} &= \dfrac{
\boldsymbol{\theta}^{o}(\boldsymbol{u}^{o}_{+}(0) \! + \! \boldsymbol{d}_{o})
\mathbb{E}^{-1}}{\boldsymbol{\theta}^{o}(\boldsymbol{u}^{o}_{+}(0) \! - \!
\frac{1}{2 \pi}(n \! + \! \frac{1}{2}) \boldsymbol{\Omega}^{o} \! + \!
\boldsymbol{d}_{o})} \! \left(\! \left(\theta^{o}_{\infty}(-1,1,\boldsymbol{
\Omega}^{o}) \! \left((\widetilde{\alpha}_{\infty}^{o}(-1,1,\vec{\pmb{0}}))^{
2} \! + \! \beta^{o}_{\infty}(-1,1,\vec{\pmb{0}}) \theta^{o}_{\infty}(-1,1,
\vec{\pmb{0}}) \right) \right. \right. \\
&\left. \left. - \, \widetilde{\alpha}_{\infty}^{o}(-1,1,\boldsymbol{\Omega}^{
o}) \widetilde{\alpha}_{\infty}^{o}(-1,1,\vec{\pmb{0}}) \theta^{o}_{\infty}(-1,
1,\vec{\pmb{0}}) \! - \! \beta^{o}_{\infty}(-1,1,\boldsymbol{\Omega}^{o})
(\theta^{o}_{\infty}(-1,1,\vec{\pmb{0}}))^{2} \right) \! \tfrac{1}{(\theta^{
o}_{\infty}(-1,1,\vec{\pmb{0}}))^{3}} \right. \\
&\left. \times \left(\dfrac{\gamma^{o}_{0} \! - \! (\gamma^{o}_{0})^{-1}}{2
\mi} \right) \! - \! \left(\dfrac{\widetilde{\alpha}_{\infty}^{o}(-1,1,\vec{
\pmb{0}}) \theta^{o}_{\infty}(-1,1,\boldsymbol{\Omega}^{o}) \! - \! \widetilde{
\alpha}_{\infty}^{o}(-1,1,\boldsymbol{\Omega}^{o}) \theta^{o}_{\infty}(-1,1,
\vec{\pmb{0}})}{(\theta^{o}_{\infty}(-1,1,\vec{\pmb{0}}))^{2}} \right) \right.
\\
&\left. \times \left(\sum_{k=1}^{N+1} \! \left(a_{k}^{o} \! - \! b_{k-1}^{o}
\right) \right) \! \left(\dfrac{\gamma^{o}_{0} \! + \! (\gamma^{o}_{0})^{-1}}{
8 \mi} \right) \! - \! \left(\! \left(\dfrac{\gamma^{o}_{0} \! + \! (\gamma^{
o}_{0})^{-1}}{16 \mi} \right) \! \sum_{k=1}^{N+1} \! \left((a_{k}^{o})^{2} \!
- \! (b_{k-1}^{o})^{2} \right) \right. \right. \\
&\left. \left. -\left(\dfrac{\gamma^{o}_{0} \! - \! (\gamma^{o}_{0})^{-1}}{64
\mi} \right) \! \left(\sum_{k=1}^{N+1} \! \left(a_{k}^{o} \! - \! b_{k-1}^{o}
\right) \right)^{2} \right) \! \dfrac{\theta^{o}_{\infty}(-1,1,\boldsymbol{
\Omega}^{o})}{\theta^{o}_{\infty}(-1,1,\vec{\pmb{0}})} \right), \\
(\overset{o}{m}_{2}^{\raise-1.0ex\hbox{$\scriptstyle \infty$}})_{21} &=
-\dfrac{\boldsymbol{\theta}^{o}(\boldsymbol{u}^{o}_{+}(0) \! + \! \boldsymbol{
d}_{o}) \mathbb{E}}{\boldsymbol{\theta}^{o}(-\boldsymbol{u}^{o}_{+}(0) \! -
\! \frac{1}{2 \pi}(n \! + \! \frac{1}{2}) \boldsymbol{\Omega}^{o} \! - \!
\boldsymbol{d}_{o})} \! \left(\! \left(\theta^{o}_{\infty}(1,-1,\boldsymbol{
\Omega}^{o}) \! \left((\widetilde{\alpha}_{\infty}^{o}(1,-1,\vec{\pmb{0}}))^{
2} \! + \! \beta^{o}_{\infty}(1,-1,\vec{\pmb{0}}) \theta^{o}_{\infty}(1,-1,
\vec{\pmb{0}}) \right) \right. \right. \\
&\left. \left. - \, \widetilde{\alpha}_{\infty}^{o}(1,-1,\boldsymbol{\Omega}^{
o}) \widetilde{\alpha}_{\infty}^{o}(1,-1,\vec{\pmb{0}}) \theta^{o}_{\infty}(1,
-1,\vec{\pmb{0}}) \! - \! \beta^{o}_{\infty}(1,-1,\boldsymbol{\Omega}^{o})
(\theta^{o}_{\infty}(1,-1,\vec{\pmb{0}}))^{2} \right) \! \tfrac{1}{(\theta^{
o}_{\infty}(1,-1,\vec{\pmb{0}}))^{3}} \right. \\
&\left. \times \left(\dfrac{\gamma^{o}_{0} \! - \! (\gamma^{o}_{0})^{-1}}{2
\mi} \right) \! - \! \left(\dfrac{\widetilde{\alpha}_{\infty}^{o}(1,-1,\vec{
\pmb{0}}) \theta^{o}_{\infty}(1,-1,\boldsymbol{\Omega}^{o}) \! - \! \widetilde{
\alpha}_{\infty}^{o}(1,-1,\boldsymbol{\Omega}^{o}) \theta^{o}_{\infty}(1,-1,
\vec{\pmb{0}})}{(\theta^{o}_{\infty}(1,-1,\vec{\pmb{0}}))^{2}} \right) \right.
\\
&\left. \times \left(\sum_{k=1}^{N+1} \! \left(a_{k}^{o} \! - \! b_{k-1}^{o}
\right) \right) \! \left(\dfrac{\gamma^{o}_{0} \! + \! (\gamma^{o}_{0})^{-1}}{
8 \mi} \right) \! - \! \left(\! \left(\dfrac{\gamma^{o}_{0} \! + \! (\gamma^{
o}_{0})^{-1}}{16 \mi} \right) \! \sum_{k=1}^{N+1} \! \left((a_{k}^{o})^{2} \!
- \! (b_{k-1}^{o})^{2} \right) \right. \right. \\
&\left. \left. -\left(\dfrac{\gamma^{o}_{0} \! - \! (\gamma^{o}_{0})^{-1}}{64
\mi} \right) \! \left(\sum_{k=1}^{N+1} \! \left(a_{k}^{o} \! - \! b_{k-1}^{o}
\right) \right)^{2} \right) \! \dfrac{\theta^{o}_{\infty}(1,-1,\boldsymbol{
\Omega}^{o})}{\theta^{o}_{\infty}(1,-1,\vec{\pmb{0}})} \right), \\
(\overset{o}{m}_{2}^{\raise-1.0ex\hbox{$\scriptstyle \infty$}})_{22} &= \dfrac{
\boldsymbol{\theta}^{o}(\boldsymbol{u}^{o}_{+}(0) \! + \! \boldsymbol{d}_{o})
\mathbb{E}}{\boldsymbol{\theta}^{o}(-\boldsymbol{u}^{o}_{+}(0) \! - \! \frac{
1}{2 \pi}(n \! + \! \frac{1}{2}) \boldsymbol{\Omega}^{o} \! - \! \boldsymbol{
d}_{o})} \! \left(\! \left(\theta^{o}_{\infty}(-1,-1,\boldsymbol{\Omega}^{o})
\! \left((\widetilde{\alpha}_{\infty}^{o}(-1,-1,\vec{\pmb{0}}))^{2} \! + \!
\beta^{o}_{\infty}(-1,-1,\vec{\pmb{0}}) \theta^{o}_{\infty}(-1,-1,\vec{\pmb{0}
}) \right) \right. \right. \\
&\left. \left. - \, \widetilde{\alpha}_{\infty}^{o}(-1,-1,\boldsymbol{\Omega}^{
o}) \widetilde{\alpha}_{\infty}^{o}(-1,-1,\vec{\pmb{0}}) \theta^{o}_{\infty}
(-1,-1,\vec{\pmb{0}}) \! - \! \beta^{o}_{\infty}(-1,-1,\boldsymbol{\Omega}^{o})
(\theta^{o}_{\infty}(-1,-1,\vec{\pmb{0}}))^{2} \right) \! \tfrac{1}{(\theta^{
o}_{\infty}(-1,-1,\vec{\pmb{0}}))^{3}} \right. \\
&\left. \times \left(\dfrac{\gamma^{o}_{0} \! + \! (\gamma^{o}_{0})^{-1}}{2}
\right) \! + \! \left(\dfrac{\widetilde{\alpha}_{\infty}^{o}(-1,-1,\vec{\pmb{
0}}) \theta^{o}_{\infty}(-1,-1,\boldsymbol{\Omega}^{o}) \! - \! \widetilde{
\alpha}_{\infty}^{o}(-1,-1,\boldsymbol{\Omega}^{o}) \theta^{o}_{\infty}(-1,-1,
\vec{\pmb{0}})}{(\theta^{o}_{\infty}(-1,-1,\vec{\pmb{0}}))^{2}} \right)
\right. \\
&\left. \times \left(\sum_{k=1}^{N+1} \! \left(b_{k-1}^{o} \! - \! a_{k}^{o}
\right) \right) \! \left(\dfrac{\gamma^{o}_{0} \! - \! (\gamma^{o}_{0})^{-1}}{
8} \right) \! + \! \left(\! \left(\dfrac{\gamma^{o}_{0} \! - \! (\gamma^{o}_{
0})^{-1}}{16} \right) \! \sum_{k=1}^{N+1} \! \left((b_{k-1}^{o})^{2} \! - \!
(a_{k}^{o})^{2} \right) \right. \right. \\
&\left. \left. +\left(\dfrac{\gamma^{o}_{0} \! + \! (\gamma^{o}_{0})^{-1}}{64}
\right) \! \left(\sum_{k=1}^{N+1} \! \left(a_{k}^{o} \! - \! b_{k-1}^{o}
\right) \right)^{2} \right) \! \dfrac{\theta^{o}_{\infty}(-1,-1,\boldsymbol{
\Omega}^{o})}{\theta^{o}_{\infty}(-1,-1,\vec{\pmb{0}})} \right),
\end{align*}
with $(\star)_{ij}$, $i,j \! = \! 1,2$, denoting the $(i \, j)$-element of
$\star$, $\mathbb{E}$ defined in Proposition~{\rm 4.1}, and $\vec{\pmb{0}}
\! := \! (0,0,\dotsc,0)^{\mathrm{T}}$ $(\in \! \mathbb{R}^{N})$. Set
\begin{align*}
(\widehat{Q}^{o}_{0})_{11} :=& \,
(\overset{o}{m}_{0}^{\raise-1.0ex\hbox{$\scriptstyle \infty$}})_{11} \! \left(
1 \! + \! (\mathscr{R}^{o,\infty}_{0}(n))_{11} \right) \! + \! (\mathscr{R}^{o,
\infty}_{0}(n))_{12}
(\overset{o}{m}_{0}^{\raise-1.0ex\hbox{$\scriptstyle \infty$}})_{21}, \\
(\widehat{Q}^{o}_{0})_{12} :=& \,
(\overset{o}{m}_{0}^{\raise-1.0ex\hbox{$\scriptstyle \infty$}})_{12} \! \left(
1 \! + \! (\mathscr{R}^{o,\infty}_{0}(n))_{11} \right) \! + \! (\mathscr{R}^{o,
\infty}_{0}(n))_{12}
(\overset{o}{m}_{0}^{\raise-1.0ex\hbox{$\scriptstyle \infty$}})_{22}, \\
(\widehat{Q}^{o}_{0})_{21} :=& \,
(\overset{o}{m}_{0}^{\raise-1.0ex\hbox{$\scriptstyle \infty$}})_{21} \! \left(
1 \! + \! (\mathscr{R}^{o,\infty}_{0}(n))_{22} \right) \! + \! (\mathscr{R}^{o,
\infty}_{0}(n))_{21}
(\overset{o}{m}_{0}^{\raise-1.0ex\hbox{$\scriptstyle \infty$}})_{11}, \\
(\widehat{Q}^{o}_{0})_{22} :=& \,
(\overset{o}{m}_{0}^{\raise-1.0ex\hbox{$\scriptstyle \infty$}})_{22} \! \left(
1 \! + \! (\mathscr{R}^{o,\infty}_{0}(n))_{22} \right) \! + \! (\mathscr{R}^{o,
\infty}_{0}(n))_{21}
(\overset{o}{m}_{0}^{\raise-1.0ex\hbox{$\scriptstyle \infty$}})_{12}, \\
(\widehat{Q}^{o}_{1})_{11} :=& \,
(\overset{o}{m}_{1}^{\raise-1.0ex\hbox{$\scriptstyle \infty$}})_{11} \! \left(
1 \! + \! (\mathscr{R}^{o,\infty}_{0}(n))_{11} \right) \! + \! (\mathscr{R}^{o,
\infty}_{0}(n))_{12}
(\overset{o}{m}_{1}^{\raise-1.0ex\hbox{$\scriptstyle \infty$}})_{21} \! + \!
(\mathscr{R}^{o,\infty}_{1}(n))_{11}
(\overset{o}{m}_{0}^{\raise-1.0ex\hbox{$\scriptstyle \infty$}})_{11} \\
+& \, (\mathscr{R}^{o,\infty}_{1}(n))_{12}
(\overset{o}{m}_{0}^{\raise-1.0ex\hbox{$\scriptstyle \infty$}})_{21}, \\
(\widehat{Q}^{o}_{1})_{12} :=& \,
(\overset{o}{m}_{1}^{\raise-1.0ex\hbox{$\scriptstyle \infty$}})_{12} \! \left(
1 \! + \! (\mathscr{R}^{o,\infty}_{0}(n))_{11} \right) \! + \! (\mathscr{R}^{o,
\infty}_{0}(n))_{12}
(\overset{o}{m}_{1}^{\raise-1.0ex\hbox{$\scriptstyle \infty$}})_{22} \! + \!
(\mathscr{R}^{o,\infty}_{1}(n))_{11}
(\overset{o}{m}_{0}^{\raise-1.0ex\hbox{$\scriptstyle \infty$}})_{12} \\
+& \, (\mathscr{R}^{o,\infty}_{1}(n))_{12}
(\overset{o}{m}_{0}^{\raise-1.0ex\hbox{$\scriptstyle \infty$}})_{22}, \\
(\widehat{Q}^{o}_{1})_{21} :=& \,
(\overset{o}{m}_{1}^{\raise-1.0ex\hbox{$\scriptstyle \infty$}})_{21} \! \left(
1 \! + \! (\mathscr{R}^{o,\infty}_{0}(n))_{22} \right) \! + \! (\mathscr{R}^{o,
\infty}_{0}(n))_{21}
(\overset{o}{m}_{1}^{\raise-1.0ex\hbox{$\scriptstyle \infty$}})_{11} \! + \!
(\mathscr{R}^{o,\infty}_{1}(n))_{21}
(\overset{o}{m}_{0}^{\raise-1.0ex\hbox{$\scriptstyle \infty$}})_{11} \\
+& \, (\mathscr{R}^{o,\infty}_{1}(n))_{22}
(\overset{o}{m}_{0}^{\raise-1.0ex\hbox{$\scriptstyle \infty$}})_{21}, \\
(\widehat{Q}^{o}_{1})_{22} :=& \,
(\overset{o}{m}_{1}^{\raise-1.0ex\hbox{$\scriptstyle \infty$}})_{22} \! \left(
1 \! + \! (\mathscr{R}^{o,\infty}_{0}(n))_{22} \right) \! + \! (\mathscr{R}^{o,
\infty}_{0}(n))_{21}
(\overset{o}{m}_{1}^{\raise-1.0ex\hbox{$\scriptstyle \infty$}})_{12} \! + \!
(\mathscr{R}^{o,\infty}_{1}(n))_{21}
(\overset{o}{m}_{0}^{\raise-1.0ex\hbox{$\scriptstyle \infty$}})_{12} \\
+& \, (\mathscr{R}^{o,\infty}_{1}(n))_{22}
(\overset{o}{m}_{0}^{\raise-1.0ex\hbox{$\scriptstyle \infty$}})_{22}, \\
(\widehat{Q}^{o}_{2})_{11} :=& \,
(\overset{o}{m}_{2}^{\raise-1.0ex\hbox{$\scriptstyle \infty$}})_{11} \! \left(
1 \! + \! (\mathscr{R}^{o,\infty}_{0}(n))_{11} \right) \! + \! (\mathscr{R}^{o,
\infty}_{0}(n))_{12}
(\overset{o}{m}_{2}^{\raise-1.0ex\hbox{$\scriptstyle \infty$}})_{21} \! + \!
(\mathscr{R}^{o,\infty}_{1}(n))_{11}
(\overset{o}{m}_{1}^{\raise-1.0ex\hbox{$\scriptstyle \infty$}})_{11} \\
+& \, (\mathscr{R}^{o,\infty}_{1}(n))_{12}
(\overset{o}{m}_{1}^{\raise-1.0ex\hbox{$\scriptstyle \infty$}})_{21} \! + \!
(\mathscr{R}^{o,\infty}_{2}(n))_{11}
(\overset{o}{m}_{0}^{\raise-1.0ex\hbox{$\scriptstyle \infty$}})_{11} \! + \!
(\mathscr{R}^{o,\infty}_{2}(n))_{12}
(\overset{o}{m}_{0}^{\raise-1.0ex\hbox{$\scriptstyle \infty$}})_{21}, \\
(\widehat{Q}^{o}_{2})_{12} :=& \,
(\overset{o}{m}_{2}^{\raise-1.0ex\hbox{$\scriptstyle \infty$}})_{12} \! \left(
1 \! + \! (\mathscr{R}^{o,\infty}_{0}(n))_{11} \right) \! + \! (\mathscr{R}^{o,
\infty}_{0}(n))_{12}
(\overset{o}{m}_{2}^{\raise-1.0ex\hbox{$\scriptstyle \infty$}})_{22} \! + \!
(\mathscr{R}^{o,\infty}_{1}(n))_{11}
(\overset{o}{m}_{1}^{\raise-1.0ex\hbox{$\scriptstyle \infty$}})_{12} \\
+& \, (\mathscr{R}^{o,\infty}_{1}(n))_{12}
(\overset{o}{m}_{1}^{\raise-1.0ex\hbox{$\scriptstyle \infty$}})_{22} \! + \!
(\mathscr{R}^{o,\infty}_{2}(n))_{11}
(\overset{o}{m}_{0}^{\raise-1.0ex\hbox{$\scriptstyle \infty$}})_{12} \! + \!
(\mathscr{R}^{o,\infty}_{2}(n))_{12}
(\overset{o}{m}_{0}^{\raise-1.0ex\hbox{$\scriptstyle \infty$}})_{22}, \\
(\widehat{Q}^{o}_{2})_{21} :=& \,
(\overset{o}{m}_{2}^{\raise-1.0ex\hbox{$\scriptstyle \infty$}})_{21} \! \left(
1 \! + \! (\mathscr{R}^{o,\infty}_{0}(n))_{22} \right) \! + \! (\mathscr{R}^{o,
\infty}_{0}(n))_{21}
(\overset{o}{m}_{2}^{\raise-1.0ex\hbox{$\scriptstyle \infty$}})_{11} \! + \!
(\mathscr{R}^{o,\infty}_{1}(n))_{21}
(\overset{o}{m}_{1}^{\raise-1.0ex\hbox{$\scriptstyle \infty$}})_{11} \\
+& \, (\mathscr{R}^{o,\infty}_{1}(n))_{22}
(\overset{o}{m}_{1}^{\raise-1.0ex\hbox{$\scriptstyle \infty$}})_{21} \! + \!
(\mathscr{R}^{o,\infty}_{2}(n))_{21}
(\overset{o}{m}_{0}^{\raise-1.0ex\hbox{$\scriptstyle \infty$}})_{11} \! + \!
(\mathscr{R}^{o,\infty}_{2}(n))_{22}
(\overset{o}{m}_{0}^{\raise-1.0ex\hbox{$\scriptstyle \infty$}})_{21}, \\
(\widehat{Q}^{o}_{2})_{22} :=& \,
(\overset{o}{m}_{2}^{\raise-1.0ex\hbox{$\scriptstyle \infty$}})_{22} \! \left(
1 \! + \! (\mathscr{R}^{o,\infty}_{0}(n))_{22} \right) \! + \! (\mathscr{R}^{o,
\infty}_{0}(n))_{21}
(\overset{o}{m}_{2}^{\raise-1.0ex\hbox{$\scriptstyle \infty$}})_{12} \! + \!
(\mathscr{R}^{o,\infty}_{1}(n))_{21}
(\overset{o}{m}_{1}^{\raise-1.0ex\hbox{$\scriptstyle \infty$}})_{12} \\
+& \, (\mathscr{R}^{o,\infty}_{1}(n))_{22}
(\overset{o}{m}_{1}^{\raise-1.0ex\hbox{$\scriptstyle \infty$}})_{22} \! + \!
(\mathscr{R}^{o,\infty}_{2}(n))_{21}
(\overset{o}{m}_{0}^{\raise-1.0ex\hbox{$\scriptstyle \infty$}})_{12} \! + \!
(\mathscr{R}^{o,\infty}_{2}(n))_{22}
(\overset{o}{m}_{0}^{\raise-1.0ex\hbox{$\scriptstyle \infty$}})_{22}.
\end{align*}

Let $\overset{o}{\mathrm{Y}} \colon \mathbb{C} \setminus \mathbb{R} \! \to \!
\operatorname{SL}_{2}(\mathbb{C})$ be the solution of {\rm \pmb{RHP2}}. Then,
\begin{equation*}
\overset{o}{\mathrm{Y}}(z)z^{-(n+1) \sigma_{3}} \underset{z \to \infty}{=}
\mathrm{Y}^{o,\infty}_{0} \! + \! \dfrac{1}{z} \mathrm{Y}^{o,\infty}_{1} \!
+ \! \dfrac{1}{z^{2}} \mathrm{Y}^{o,\infty}_{2} \! + \! \mathcal{O}(z^{-3}),
\end{equation*}
where
\begin{align*}
(\mathrm{Y}^{o,\infty}_{0})_{11} =& \, (\widehat{Q}^{o}_{0})_{11} \me^{-2(n+
\frac{1}{2}) \int_{J_{o}} \ln (\lvert s \rvert) \psi_{V}^{o}(s) \, \md s}, \\
(\mathrm{Y}^{o,\infty}_{0})_{12} =& \, (\widehat{Q}^{o}_{0})_{12} \me^{n(
\ell_{o}+2(1+\frac{1}{2n}) \int_{J_{o}} \ln (\lvert s \rvert) \psi_{V}^{o}(s)
\, \md s)}, \\
(\mathrm{Y}^{o,\infty}_{0})_{21} =& \, (\widehat{Q}^{o}_{0})_{21} \me^{-n(
\ell_{o}+2(1+\frac{1}{2n}) \int_{J_{o}} \ln (\lvert s \rvert) \psi_{V}^{o}(s)
\, \md s)}, \\
(\mathrm{Y}^{o,\infty}_{0})_{22} =& \, (\widehat{Q}^{o}_{0})_{22} \me^{2(n+
\frac{1}{2}) \int_{J_{o}} \ln (\lvert s \rvert) \psi_{V}^{o}(s) \, \md s}, \\
(\mathrm{Y}^{o,\infty}_{1})_{11} =& \, \left(\! (\widehat{Q}^{o}_{1})_{11} \!
- \! (2n \! + \! 1)(\widehat{Q}^{o}_{0})_{11} \int_{J_{o}}s \psi_{V}^{o}(s) \,
\md s \right) \! \me^{-2(n+\frac{1}{2}) \int_{J_{o}} \ln (\lvert s \rvert)
\psi_{V}^{o}(s) \, \md s}, \\
(\mathrm{Y}^{o,\infty}_{1})_{12} =& \, \left(\! (\widehat{Q}^{o}_{1})_{12} \!
+ \! (2n \! + \! 1)(\widehat{Q}^{o}_{0})_{12} \int_{J_{o}}s \psi_{V}^{o}(s) \,
\md s \right) \! \me^{n(\ell_{o}+2(1+\frac{1}{2n}) \int_{J_{o}} \ln (\lvert s
\rvert) \psi_{V}^{o}(s) \, \md s)}, \\
(\mathrm{Y}^{o,\infty}_{1})_{21} =& \, \left(\! (\widehat{Q}^{o}_{1})_{21} \!
- \! (2n \! + \! 1)(\widehat{Q}^{o}_{0})_{21} \int_{J_{o}}s \psi_{V}^{o}(s) \,
\md s \right) \! \me^{-n(\ell_{o}+2(1+\frac{1}{2n}) \int_{J_{o}} \ln (\lvert s
\rvert) \psi_{V}^{o}(s) \, \md s)}, \\
(\mathrm{Y}^{o,\infty}_{1})_{22} =& \, \left(\! (\widehat{Q}^{o}_{1})_{22} \!
+ \! (2n \! + \! 1)(\widehat{Q}^{o}_{0})_{22} \int_{J_{o}}s \psi_{V}^{o}(s) \,
\md s \right) \! \me^{2(n+\frac{1}{2}) \int_{J_{o}} \ln (\lvert s \rvert)
\psi_{V}^{o}(s) \, \md s}, \\
(\mathrm{Y}^{o,\infty}_{2})_{11} =& \left(\! (\widehat{Q}^{o}_{2})_{11} \! -
\! (2n \! + \! 1)(\widehat{Q}^{o}_{1})_{11} \int_{J_{o}}s \psi_{V}^{o}(s) \,
\md s \! + \! (\widehat{Q}^{o}_{0})_{11} \! \left(\tfrac{1}{2}(2n \! + \! 1)^{
2} \! \left(\int_{J_{o}}s \psi_{V}^{o}(s) \, \md s \right)^{2} \right. \right.
\\
-& \left. \left. \tfrac{1}{2}(2n \! + \! 1) \int_{J_{o}}s^{2} \psi_{V}^{o}(s)
\, \md s \right) \right) \! \me^{-2(n+\frac{1}{2}) \int_{J_{o}} \ln (\lvert s
\rvert) \psi_{V}^{o}(s) \, \md s}, \\
(\mathrm{Y}^{o,\infty}_{2})_{12} =& \left(\! (\widehat{Q}^{o}_{2})_{12} \! +
\! (2n \! + \! 1)(\widehat{Q}^{o}_{1})_{12} \int_{J_{o}}s \psi_{V}^{o}(s) \,
\md s \! + \! (\widehat{Q}^{o}_{0})_{12} \! \left(\tfrac{1}{2}(2n \! + \! 1)^{
2} \! \left(\int_{J_{o}}s \psi_{V}^{o}(s) \, \md s \right)^{2} \right. \right.
\\
+& \left. \left. \tfrac{1}{2}(2n \! + \! 1) \int_{J_{o}}s^{2} \psi_{V}^{o}(s)
\, \md s \right) \right) \! \me^{n(\ell_{o}+2(1+\frac{1}{2n}) \int_{J_{o}} \ln
(\lvert s \rvert) \psi_{V}^{o}(s) \, \md s)}, \\
(\mathrm{Y}^{o,\infty}_{2})_{21} =& \left(\! (\widehat{Q}^{o}_{2})_{21} \! -
\! (2n \! + \! 1)(\widehat{Q}^{o}_{1})_{21} \int_{J_{o}}s \psi_{V}^{o}(s) \,
\md s \! + \! (\widehat{Q}^{o}_{0})_{21} \! \left(\tfrac{1}{2}(2n \! + \! 1)^{
2} \! \left(\int_{J_{o}}s \psi_{V}^{o}(s) \, \md s \right)^{2} \right. \right.
\\
-& \left. \left. \tfrac{1}{2}(2n \! + \! 1) \int_{J_{o}}s^{2} \psi_{V}^{o}(s)
\, \md s \right) \right) \! \me^{-n(\ell_{o}+2(1+\frac{1}{2n}) \int_{J_{o}}
\ln (\lvert s \rvert) \psi_{V}^{o}(s) \, \md s)}, \\
(\mathrm{Y}^{o,\infty}_{2})_{22} =& \left(\! (\widehat{Q}^{o}_{2})_{22} \! +
\! (2n \! + \! 1)(\widehat{Q}^{o}_{1})_{22} \int_{J_{o}}s \psi_{V}^{o}(s) \,
\md s \! + \! (\widehat{Q}^{o}_{0})_{22} \! \left(\tfrac{1}{2}(2n \! + \! 1)^{
2} \! \left(\int_{J_{o}}s \psi_{V}^{o}(s) \, \md s \right)^{2} \right. \right.
\\
+& \left. \left. \tfrac{1}{2}(2n \! + \! 1) \int_{J_{o}}s^{2} \psi_{V}^{o}(s)
\, \md s \right) \right) \! \me^{2(n+\frac{1}{2}) \int_{J_{o}} \ln (\lvert s
\rvert) \psi_{V}^{o}(s) \, \md s}.
\end{align*}
\end{ay}

\emph{Proof.} Let $\mathscr{R}^{o} \colon \mathbb{C} \setminus \widetilde{
\Sigma}^{o}_{p} \! \to \! \operatorname{SL}_{2}(\mathbb{C})$ be the solution
of the RHP $(\mathscr{R}^{o}(z),\upsilon_{\mathscr{R}}^{o}(z),\widetilde{
\Sigma}^{o}_{p})$ formulated in Proposition~5.2 with $n \! \to \! \infty$
asymptotics given in Lemma~5.3. For $\vert z \vert \! \gg \! \max_{j=1,
\dotsc,N+1} \lbrace \vert b_{j-1}^{o} \! - \! a_{j}^{o} \vert \rbrace$, via
the expansion $\tfrac{1}{s-z} \! = \! -\sum_{k=0}^{l} \tfrac{s^{k}}{z^{k+1}}
\! + \! \tfrac{s^{l+1}}{z^{l+1}(s-z)}$, $l \! \in \! \mathbb{Z}_{0}^{+}$,
where $s \! \in \! \lbrace b_{j-1}^{o},a_{j}^{o} \rbrace$, $j \! = \! 1,
\dotsc,N \! + \! 1$, one obtains the asymptotics for $\mathscr{R}^{o}(z)$
stated in the Proposition.

Let $\overset{o}{m}^{\raise-1.0ex\hbox{$\scriptstyle \infty$}} \colon \mathbb{
C} \setminus J_{o}^{\infty} \! \to \! \operatorname{SL}_{2}(\mathbb{C})$ solve
the RHP $(\overset{o}{m}^{\raise-1.0ex\hbox{$\scriptstyle \infty$}}(z),J_{o}^{
\infty},\overset{o}{\upsilon}^{\raise-1.0ex\hbox{$\scriptstyle \infty$}}(z))$
formulated in Lemma~4.3 with (unique) solution given by Lemma~4.5. In order
to obtain large-$z$ asymptotics of
$\overset{o}{m}^{\raise-1.0ex\hbox{$\scriptstyle \infty$}}(z)$, one needs
large-$z$ asymptotics of $(\gamma^{o}(z))^{\pm 1}$ and $\tfrac{\boldsymbol{
\theta}^{o}(\varepsilon_{1} \boldsymbol{u}^{o}(z)-\frac{1}{2 \pi}(n+\frac{1}{
2}) \boldsymbol{\Omega}^{o}+\varepsilon_{2} \boldsymbol{d}_{o})}{\boldsymbol{
\theta}^{o}(\varepsilon_{1} \boldsymbol{u}^{o}(z)+\varepsilon_{2} \boldsymbol{
d}_{o})}$, $\varepsilon_{1},\varepsilon_{2} \! = \! \pm 1$. Consider, say, and
without loss of generality, $z \! \to \! \infty$ asymptotics for $z \! \in \!
\mathbb{C}_{+}$, that is, $z \! \to \! \infty^{+}$, where, by definition,
$\sqrt{\smash[b]{\star (z)}} \! := \! +\sqrt{\smash[b]{\star (z)}}$:
equivalently, one may consider $z \! \to \! \infty$ asymptotics for $z \! \in
\! \mathbb{C}_{-}$, that is, $z \! \to \! \infty^{-}$; however, recalling
that $\sqrt{\smash[b]{\star (z)}} \! \upharpoonright_{\mathbb{C}_{+}} \! = \!
-\sqrt{\smash[b]{\star (z)}} \! \upharpoonright_{\mathbb{C}_{-}}$, one obtains
(in either case, and via the sheet-interchange index) the same $z \! \to \!
\infty$ asymptotics (for
$\overset{o}{m}^{\raise-1.0ex\hbox{$\scriptstyle \infty$}}(z))$. Recall the
expression for $\gamma^{o}(z)$ given in Lemma~4.4: for $\vert z \vert \! \gg
\! \max_{j=1,\dotsc,N+1} \lbrace \vert b_{j-1}^{o} \! - \! a_{j}^{o} \vert
\rbrace$, via the expansions $\tfrac{1}{s-z} \! = \! -\sum_{k=0}^{l} \tfrac{
s^{k}}{z^{k+1}} \! + \! \tfrac{s^{l+1}}{z^{l+1}(s-z)}$, $l \! \in \! \mathbb{
Z}_{0}^{+}$, and $\ln (z \! - \! s) \! =_{\vert z \vert \to \infty} \! \ln (z)
\! - \! \sum_{k=1}^{\infty} \tfrac{1}{k}(\tfrac{s}{z})^{k}$, where $s \! \in
\! \lbrace b_{j-1}^{o},a_{j}^{o} \rbrace$, $j \! = \! 1,\dotsc,N \! + \! 1$,
one shows that, upon defining $\gamma^{o}(0)$ as in the Proposition,
\begin{align*}
(\gamma^{o}(z))^{\pm 1} \underset{z \to \infty^{+}}{=}& 1 \! + \! \dfrac{1}{z}
\! \left(\pm \dfrac{1}{4} \sum_{k=1}^{N+1}(a_{k}^{o} \! - \! b_{k-1}^{o})
\right) \! + \! \dfrac{1}{z^{2}} \! \left(\pm \dfrac{1}{8} \sum_{k=1}^{N+1}
\! \left((a_{k}^{o})^{2} \! - \! (b_{k-1}^{o})^{2} \right) \right. \\
+&\left. \dfrac{1}{32} \! \left(\sum_{k=1}^{N+1}(a_{k}^{o} \! - \! b_{k-1}^{
o}) \right)^{2} \, \right) \! + \! \mathcal{O} \! \left(\dfrac{1}{z^{3}}
\right),
\end{align*}
whence
\begin{align*}
\dfrac{1}{2}((\gamma^{o}_{0})^{-1} \gamma^{o}(z) \! + \! \gamma^{o}_{0}
(\gamma^{o}(z))^{-1}) \underset{z \to \infty^{+}}{=}& \, \dfrac{(\gamma^{o}_{
0} \! + \! (\gamma^{o}_{0})^{-1})}{2} \! + \! \dfrac{1}{z} \! \left(\! \left(
\dfrac{\gamma^{o}_{0} \! - \! (\gamma^{o}_{0})^{-1}}{8} \right) \sum_{k=1}^{N
+1} \! \left(b_{k-1}^{o} \! - \! a_{k}^{o} \right) \right) \\
+& \, \dfrac{1}{z^{2}} \! \left(\! \left(\dfrac{\gamma^{o}_{0} \! - \!
(\gamma^{o}_{0})^{-1}}{16} \right) \sum_{k=1}^{N+1} \! \left((b_{k-1}^{o})^{2}
\! - \! (a_{k}^{o})^{2} \right) \! + \! \left(\dfrac{\gamma^{o}_{0} \! + \!
(\gamma^{o}_{0})^{-1}}{64} \right) \right. \\
\times&\left. \left(\sum_{k=1}^{N+1} \! \left(b_{k-1}^{o} \! - \! a_{k}^{o}
\right) \right)^{2} \, \right) \! + \! \mathcal{O} \! \left(\dfrac{1}{z^{3}}
\right),
\end{align*}
and
\begin{align*}
\dfrac{1}{2 \mi}((\gamma^{o}_{0})^{-1} \gamma^{o}(z) \! - \! \gamma^{o}_{0}
(\gamma^{o}(z))^{-1}) \underset{z \to \infty^{+}}{=}& \, -\dfrac{(\gamma^{o}_{
0} \! - \! (\gamma^{o}_{0})^{-1})}{2 \mi} \! + \! \dfrac{1}{z} \! \left(\!
\left(\dfrac{\gamma^{o}_{0} \! + \! (\gamma^{o}_{0})^{-1}}{8 \mi} \right)
\sum_{k=1}^{N+1} \! \left(a_{k}^{o} \! - \! b_{k-1}^{o} \right) \right) \\
+& \, \dfrac{1}{z^{2}} \! \left(\! \left(\dfrac{\gamma^{o}_{0} \! + \!
(\gamma^{o}_{0})^{-1}}{16 \mi} \right) \sum_{k=1}^{N+1} \! \left((a_{k}^{o})^{
2} \! - \! (b_{k-1}^{o})^{2} \right) \! - \! \left(\dfrac{\gamma^{o}_{0} \! -
\! (\gamma^{o}_{0})^{-1}}{64 \mi} \right) \right. \\
\times&\left. \left(\sum_{k=1}^{N+1} \! \left(a_{k}^{o} \! - \! b_{k-1}^{o}
\right) \right)^{2} \, \right) \! + \! \mathcal{O} \! \left(\dfrac{1}{z^{3}}
\right).
\end{align*}
Recall {}from Lemma~4.5 that $\boldsymbol{u}^{o}(z) \! := \! \int_{a_{N+1}^{
o}}^{z} \boldsymbol{\omega}^{o}$ $(\in \operatorname{Jac}(\mathcal{Y}_{o})$,
with $\mathcal{Y}_{o} \! := \! \lbrace \mathstrut (y,z); \, y^{2} \! = \!
R_{o}(z) \rbrace)$, where $\boldsymbol{\omega}^{o}$, the associated normalised
basis of holomorphic one-forms of $\mathcal{Y}_{o}$, is given by $\boldsymbol{
\omega}^{o} \! = \! (\omega_{1}^{o},\omega_{2}^{o},\dotsc,\omega_{N}^{o})$,
with $\omega_{j}^{o} \! := \! \sum_{k=1}^{N}c_{jk}^{o}(\prod_{i=1}^{N+1}(z 
\! - \! b_{i-1}^{o})(z \! - \! a_{i}^{o}))^{-1/2}z^{N-k} \, \md z$, $j \! = 
\! 1,\dotsc,N$, where $c_{jk}^{o}$, $j,k \! = \! 1,\dotsc,N$, are obtained 
{}from Equations~(O1) and~(O2). Writing
\begin{equation*}
\boldsymbol{u}^{o}(z) \! = \! \left(\int_{a_{N+1}^{o}}^{\infty^{+}} \! + \!
\int_{\infty^{+}}^{z} \right) \! \boldsymbol{\omega}^{o} \! = \! \boldsymbol{
u}^{o}_{+}(\infty) \! + \! \int_{\infty^{+}}^{z} \boldsymbol{\omega}^{o},
\end{equation*}
where $\boldsymbol{u}^{o}_{+}(\infty) \! := \! \int_{a_{N+1}^{o}}^{\infty^{+}}
\boldsymbol{\omega}^{o}$, for $\vert z \vert \! \gg \! \max_{j=1,\dotsc,N+1}
\lbrace \vert b_{j-1}^{o} \! - \! a_{j}^{o} \vert \rbrace$, via the expansions
$\tfrac{1}{s-z} \! = \! -\sum_{k=0}^{l} \tfrac{s^{k}}{z^{k+1}} \! + \! \tfrac{
s^{l+1}}{z^{l+1}(s-z)}$, $l \! \in \! \mathbb{Z}_{0}^{+}$, and $\ln (z \! -
\! s) \! =_{\vert z \vert \to \infty} \! \ln (z) \! - \! \sum_{k=1}^{\infty}
\tfrac{1}{k}(\tfrac{s}{z})^{k}$, where $s \! \in \! \lbrace b_{k-1}^{o},a_{
k}^{o} \rbrace$, $k \! = \! 1,\dotsc,N \! + \! 1$, one shows that, for $j \!
= \! 1,\dotsc,N$,
\begin{equation*}
\omega_{j}^{o} \! \underset{z \to \infty^{+}}{=} \! \dfrac{c_{j1}^{o}}{z^{2}}
\, \md z \! + \! \dfrac{(c_{j2}^{o} \! + \! \frac{1}{2}c_{j1}^{o} \sum_{i=1}^{
N+1}(a_{i}^{o} \! + \! b_{i-1}^{o}))}{z^{3}} \, \md z \! + \! \mathcal{O} \!
\left(\dfrac{\md z}{z^{4}} \right),
\end{equation*}
whence
\begin{align*}
\int_{\infty^{+}}^{z} \omega_{j}^{o} \underset{z \to \infty^{+}}{=}& -\dfrac{
c_{j1}^{o}}{z} \! - \! \dfrac{\frac{1}{2}(c_{j2}^{o} \! + \! \frac{1}{2}c_{j
1}^{o} \sum_{i=1}^{N+1}(a_{i}^{o} \! + \! b_{i-1}^{o}))}{z^{2}} \! + \!
\mathcal{O} \! \left(\dfrac{1}{z^{3}} \right) \\
=:& -\dfrac{\widehat{\alpha}^{o}_{\infty,j}}{z} \! - \! \dfrac{\widehat{
\beta}^{o}_{\infty,j}}{z^{2}} \! + \! \mathcal{O} \! \left(\dfrac{1}{z^{3}}
\right).
\end{align*}
Defining $\theta^{o}_{\infty}(\varepsilon_{1},\varepsilon_{2},\boldsymbol{
\Omega}^{o})$, $\widetilde{\alpha}^{o}_{\infty}(\varepsilon_{1},\varepsilon_{
2},\boldsymbol{\Omega}^{o})$, and $\beta^{o}_{\infty}(\varepsilon_{1},
\varepsilon_{2},\boldsymbol{\Omega}^{o})$, $\varepsilon_{1},\varepsilon_{2} \!
= \! \pm 1$, as in the Proposition, recalling that $\boldsymbol{\omega}^{o}
\! = \! (\omega_{1}^{o},\omega_{2}^{o},\dotsc,\omega_{N}^{o})$, and that the
associated $N \! \times \! N$ Riemann matrix of
$\boldsymbol{\beta}^{o}$-periods, $\tau^{o} \! = \! (\tau^{o})_{i,j=1,\dotsc,
N} \! := \! (\oint_{\boldsymbol{\beta}^{o}_{j}} \omega^{o}_{i})_{i,j=1,\dotsc,
N}$, is non-degenerate, symmetric, and $-\mi \tau^{o}$ is positive definite,
via the above asymptotic (as $z \! \to \! \infty^{+})$ expansion for $\int_{
\infty^{+}}^{z} \omega^{o}_{j}$, $j \! = \! 1,\dotsc,N$, one shows that
\begin{equation*}
\dfrac{\boldsymbol{\theta}^{o}(\varepsilon_{1} \boldsymbol{u}^{o}(z) \! - \!
\frac{1}{2 \pi}(n \! + \! \frac{1}{2}) \boldsymbol{\Omega}^{o} \! + \!
\varepsilon_{2} \boldsymbol{d}_{o})}{\boldsymbol{\theta}^{o}(\varepsilon_{1}
\boldsymbol{u}^{o}(z) \! + \! \varepsilon_{2} \boldsymbol{d}_{o})} \underset{z
\to \infty^{+}}{=} \varTheta_{0}^{o} \! + \! \dfrac{1}{z} \varTheta_{1}^{o} \!
+ \! \dfrac{1}{z^{2}} \varTheta_{2}^{o} \! + \! \mathcal{O} \! \left(\dfrac{
1}{z^{3}} \right),
\end{equation*}
where
\begin{align*}
\varTheta_{0}^{o} :=& \, \dfrac{\theta^{o}_{\infty}(\varepsilon_{1},
\varepsilon_{2},\boldsymbol{\Omega}^{o})}{\theta^{o}_{\infty}(\varepsilon_{1},
\varepsilon_{2},\vec{\pmb{0}})}, \\
\varTheta_{1}^{o} :=& \, \dfrac{\theta^{o}_{\infty}(\varepsilon_{1},
\varepsilon_{2},\boldsymbol{\Omega}^{o}) \widetilde{\alpha}^{o}_{\infty}
(\varepsilon_{1},\varepsilon_{2},\vec{\pmb{0}}) \! - \! \widetilde{\alpha}^{
o}_{\infty}(\varepsilon_{1},\varepsilon_{2},\boldsymbol{\Omega}^{o}) \theta^{
o}_{\infty}(\varepsilon_{1},\varepsilon_{2},\vec{\pmb{0}})}{(\theta^{o}_{
\infty}(\varepsilon_{1},\varepsilon_{2},\vec{\pmb{0}}))^{2}}, \\
\varTheta^{o}_{2} :=& \left(\theta^{o}_{\infty}(\varepsilon_{1},\varepsilon_{
2},\boldsymbol{\Omega}^{o}) \! \left(\beta^{o}_{\infty}(\varepsilon_{1},
\varepsilon_{2},\vec{\pmb{0}}) \theta^{o}_{\infty}(\varepsilon_{1},
\varepsilon_{2},\vec{\pmb{0}}) \! + \! (\widetilde{\alpha}^{o}_{\infty}
(\varepsilon_{1},\varepsilon_{2},\vec{\pmb{0}}))^{2} \right) \! - \!
\widetilde{\alpha}^{o}_{\infty}(\varepsilon_{1},\varepsilon_{2},\boldsymbol{
\Omega}^{o}) \right. \\
\times&\left. \widetilde{\alpha}^{o}_{\infty}(\varepsilon_{1},\varepsilon_{2},
\vec{\pmb{0}}) \theta^{o}_{\infty}(\varepsilon_{1},\varepsilon_{2},\vec{\pmb{
0}}) \! - \! \beta^{o}_{\infty}(\varepsilon_{1},\varepsilon_{2},\boldsymbol{
\Omega}^{o})(\theta^{o}_{\infty}(\varepsilon_{1},\varepsilon_{2},\vec{\pmb{0}}
))^{2} \right)(\theta^{o}_{\infty}(\varepsilon_{1},\varepsilon_{2},\vec{\pmb{
0}}))^{-3},
\end{align*}
with $\vec{\pmb{0}} \! := \! (0,0,\dotsc,0)^{\mathrm{T}}$ $(\in \! \mathbb{
R}^{N})$. Via the above asymptotic (as $z \! \to \! \infty^{+})$ expansions
for $\tfrac{1}{2}((\gamma^{o}_{0})^{-1} \gamma^{o}(z) \! + \! \gamma^{o}_{0}
(\gamma^{o}(z))^{-1})$, $\tfrac{1}{2 \mi}((\gamma^{o}_{0})^{-1} \gamma^{o}
(z) \! - \! \gamma^{o}_{0}(\gamma^{o}(z))^{-1})$, and $\tfrac{\boldsymbol{
\theta}^{o}(\varepsilon_{1} \boldsymbol{u}^{o}(z)-\frac{1}{2 \pi}(n+\frac{
1}{2}) \boldsymbol{\Omega}^{o}+\varepsilon_{2} \boldsymbol{d}_{o})}{
\boldsymbol{\theta}^{o}(\varepsilon_{1} \boldsymbol{u}^{o}(z)+\varepsilon_{2}
\boldsymbol{d}_{o})}$, one arrives at, upon recalling the expression for
$\overset{o}{m}^{\raise-1.0ex\hbox{$\scriptstyle \infty$}}(z)$ given in
Lemma~4.5, the asymptotic expansion for
$\overset{o}{m}^{\raise-1.0ex\hbox{$\scriptstyle \infty$}}(z)$ stated in the
Proposition.

Let $\overset{o}{\operatorname{Y}} \colon \mathbb{C} \setminus \mathbb{R}
\! \to \! \operatorname{SL}_{2}(\mathbb{C})$ be the (unique) solution of
\textbf{RHP2}, that is, $(\overset{o}{\operatorname{Y}}(z),\mathrm{I} \! +
\! \exp (-n \widetilde{V}(z)) \sigma_{+},\mathbb{R})$. Recall, also, that,
for $z \! \in \! \Upsilon^{o}_{1}$,
\begin{equation*}
\overset{o}{\mathrm{Y}}(z) \! = \! \me^{\frac{n \ell_{o}}{2} \operatorname{ad}
(\sigma_{3})} \mathscr{R}^{o}(z)
\overset{o}{m}^{\raise-1.0ex\hbox{$\scriptstyle \infty$}}(z) \mathbb{E}^{
\sigma_{3}} \me^{n(g^{o}(z)-\mathfrak{Q}^{+}_{\mathscr{A}}) \sigma_{3}},
\end{equation*}
and, for $z \! \in \! \Upsilon^{o}_{2}$,
\begin{equation*}
\overset{o}{\mathrm{Y}}(z) \! = \! \me^{\frac{n \ell_{o}}{2} \operatorname{ad}
(\sigma_{3})} \mathscr{R}^{o}(z)
\overset{o}{m}^{\raise-1.0ex\hbox{$\scriptstyle \infty$}}(z) \mathbb{E}^{-
\sigma_{3}} \me^{n(g^{o}(z)-\mathfrak{Q}^{-}_{\mathscr{A}}) \sigma_{3}}:
\end{equation*}
consider, say, and without loss of generality, large-$z$ asymptotics for
$\overset{o}{\operatorname{Y}}(z)$ for $z \! \in \! \Upsilon_{1}^{o}$.
Recalling {}from Lemma~3.4 that $g^{o}(z) \! := \! \int_{J_{o}} \ln ((z \!
- \! s)^{2+\frac{1}{n}}(zs)^{-1}) \psi_{V}^{o}(s) \, \md s$, $z \! \in \!
\mathbb{C} \setminus (-\infty,\max \lbrace 0,a_{N+1}^{o} \rbrace)$, for $\vert
z \vert \! \gg \! \max_{j=1,\dotsc,N+1} \lbrace \vert b_{j-1}^{o} \! - \!
a_{j}^{o} \vert \rbrace$, in particular, $\vert s/z \vert \! \ll \! 1$ with
$s \! \in \! J_{o}$, and noting that $\int_{J_{o}} \psi_{V}^{o}(s) \, \md s
\! = \! 1$ and $\int_{J_{o}}s^{m} \psi_{V}^{o}(s) \, \md s \! < \! \infty$,
$m \! \in \! \mathbb{N}$, via the expansions $\tfrac{1}{s-z} \! = \! -\sum_{
k=0}^{l} \tfrac{s^{k}}{z^{k+1}} \! + \! \tfrac{s^{l+1}}{z^{l+1}(s-z)}$, $l
\! \in \! \mathbb{Z}_{0}^{+}$, and $\ln (z \! - \! s) \! =_{\vert z \vert \to
\infty} \! \ln (z) \! - \! \sum_{k=1}^{\infty} \tfrac{1}{k}(\tfrac{s}{z})^{
k}$, one shows that
\begin{align*}
g^{o}(z) \underset{z \to \infty}{=}& \, \left(1 \! + \! \dfrac{1}{n} \right)
\! \ln (z) \! - \! \int_{J_{o}} \ln (\lvert s \rvert) \psi_{V}^{o}(s) \, \md s
\! - \! \mi \pi \int_{J_{o} \cap \mathbb{R}_{-}} \psi_{V}^{o}(s) \, \md s \!
+ \! \dfrac{1}{z} \! \left(- \! \left(2 \! + \! \dfrac{1}{n} \right) \! \int_{
J_{o}}s \psi_{V}^{o}(s) \, \md s \right) \\
+& \, \dfrac{1}{z^{2}} \! \left(-\dfrac{1}{2} \! \left(2 \! + \! \dfrac{1}{n}
\right) \! \int_{J_{o}}s^{2} \psi_{V}^{o}(s) \, \md s \right) \! + \!
\mathcal{O}(z^{-3}),
\end{align*}
where $\int_{J_{o} \cap \mathbb{R}_{-}} \psi_{V}^{o}(s) \, \md s$ is given in
Lemma~3.4. (Explicit expressions for $\int_{J_{o}}s^{k} \psi_{V}^{o}(s) \, \md
s$, $k \! = \! 1,2$, are given in Remark~3.2.) Using the asymptotic (as $z \!
\to \! \infty)$ expansions for $g^{o}(z)$, $\mathscr{R}^{o}(z)$, and
$\overset{o}{m}^{\raise-1.0ex\hbox{$\scriptstyle \infty$}}(z)$ derived above,
upon recalling the formula for $\overset{o}{\operatorname{Y}}(z)$, one arrives
at, after a matrix-multiplication argument, the asymptotic expansion for
$\overset{o}{\operatorname{Y}}(z)z^{-(n+1) \sigma_{3}}$ stated in the
Proposition. \hfill $\qed$

\end{document}